\documentclass{report}
\usepackage{amssymb}
\usepackage{amsfonts}
\usepackage{svsing2e}
\usepackage{amsmath}

\newtheorem{theorem}{Theorem}

\newtheorem{corollary}[theorem]{Corollary}

\newtheorem{definition}[theorem]{Definition}
\newtheorem{example}[theorem]{Example}

\newtheorem{lemma}[theorem]{Lemma}

\newtheorem{proposition}[theorem]{Proposition}
\newtheorem{remark}[theorem]{Remark}

\newenvironment{proof}[1][Proof]{\noindent\textbf{#1.} }{\ \rule{0.5em}{0.5em}}
\input{tcilatex}
\begin{document}

\title{Inequalities for Functions of Selfadjoint Operators on Hilbert Spaces}
\author{Silvestru Sever Dragomir}
\date{January 2011}
\maketitle

\begin{abstract}
The main aim of this book is to present recent results concerning
inequalities for continuous functions of selfadjoint operators on complex
Hilbert spaces. It is intended for use by both researchers in various fields
of Linear Operator Theory and Mathematical Inequalities, domains which have
grown exponentially in the last decade, as well as by postgraduate students
and scientists applying inequalities in their specific areas.
\end{abstract}

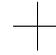
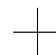
\tableofcontents

\chapter*{Preface}

\addcontentsline{toc}{chapter}{Preface}

\markboth{Preface}{Preface}

Linear Operator Theory in Hilbert spaces plays a central role in
contemporary mathematics with numerous applications for Partial Differential
Equations, in Approximation Theory, Optimization Theory, Numerical Analysis,
Probability Theory \& Statistics and other fields.

The main aim of this book is to present recent results concerning
inequalities for continuous functions of bounded selfadjoint operators on
complex Hilbert spaces.

The book is intended for use by both researchers in various fields of Linear
Operator Theory and Mathematical Inequalities, domains which have grown
exponentially in the last decade, as well as by postgraduate students and
scientists applying inequalities in their specific areas.

In the first chapter we recall some fundamental facts concerning bounded
selfadjoint operators on complex Hilbert spaces. The generalized Schwarz's
inequality for positive selfadjoint operators as well as some results for
the spectrum of this class of operators are presented. Then we introduce and
explore the fundamental results for polynomials in a linear operator,
continuous functions of selfadjoint operators as well as the step functions
of selfadjoint operators. By the use of these results we then introduce the
spectral decomposition of selfadjoint operators (the \textit{Spectral
Representation Theorem}) that will play a central role in the rest of the
book. This result is used as a key tool in obtaining various new
inequalities for continuous functions of selfadjoint operators, functions
which are of bounded variation, Lipschitzian, monotonic or absolutely
continuous. Another tool that will greatly simplify the error bounds
provided in the book is the \textit{Total Variation Schwarz's Inequality}
for which a simple proof is offered.

The chapter is concluded with some well known operator inequalities of
Jensen's type for convex and operator convex functions. Finally, some Gr\"{u}%
ss' type inequalities obtained in 1993 by Mond \& Pe\v{c}ari\'{c} are also
presented.

Jensen's type inequalities in their various settings ranging from discrete
to continuous case play an important role in different branches of Modern
Mathematics. A simple search in the \textit{MathSciNet} database of the
American Mathematical Society with the key words "jensen" and "inequality"
in the title reveals more than 300 items intimately devoted to this famous
result. However, the number of papers where this inequality is applied is a
lot larger and far more difficult to find.

In the second chapter we present some recent results obtained by the author
that deal with different aspects of this well research inequality than those
recently reported in the book \cite{I.FMPS}. They include but are not
restricted to the operator version of the Dragomir-Ionescu inequality,
Slater's type inequalities for operators and its inverses, Jensen's
inequality for twice differentiable functions whose second derivatives
satisfy some upper and lower bounds conditions, Jensen's type inequalities
for log-convex functions and for differentiable log-convex functions and
their applications to Ky Fan's inequality. Finally, some Hermite-Hadamard's
type inequalities for convex functions and Hermite-Hadamard's type
inequalities for operator convex functions are presented as well.

The third chapter is devoted to \v{C}eby\v{s}ev and Gr\"{u}ss' type
inequalities.

The \textit{\v{C}eby\v{s}ev}, or in a different spelling - \textit{%
Chebyshev, inequality} which compares the integral/discrete mean of the
product with the product of the integral/discrete means is famous in the
literature devoted to Mathematical Inequalities. It has been extended,
generalized, refined etc...by many authors during the last century. A simple
search utilizing either spellings and the key word "inequality" in the title
in the comprehensive \textit{MathSciNet} database produces more than 200
research articles devoted to this result.

The sister inequality due to Gr\"{u}ss which provides error bounds for the
magnitude of the difference between the integral mean of the product and the
product of the integral means has also attracted much interest since it has
been discovered in 1935 with more than 180 papers published, as a simple
search in the same database reveals. Far more publications have been devoted
to the applications of these inequalities and an accurate picture of the
impacted results in various fields of Modern Mathematics is difficult to
provide.

In this chapter, however, we present only some recent results due to the
author for the corresponding operator versions of these two famous
inequalities. Applications for particular functions of selfadjoint operators
such as the power, logarithmic and exponential functions are provided as
well.

The next chapter is devoted to the Ostrowski's type inequalities. They
provide sharp error estimates in approximating the value of a function by
its integral mean and can be utilized to obtain a priory error bounds for
different quadrature rules in approximating the Riemann integral by
different Riemann sums. They also shows, in general, that the mid-point rule
provides the best approximation in the class of all Riemann sums sampled in
the interior points of a given partition.

As revealed by a simple search in \textit{MathSciNet} with the key words
"Ostrowski" and "inequality" in the title, an exponential evolution of
research papers devoted to this result has been registered in the last
decade. There are now at least 280 papers that can be found by performing
the above search. Numerous extensions, generalizations in both the integral
and discrete case have been discovered. More general versions for $n$-time
differentiable functions, the corresponding versions on time scales, for
vector valued functions or multiple integrals have been established as well.
Numerous applications in Numerical Analysis, Probability Theory and other
fields have been also given.

In this chapter we present some recent results obtained by the author in
extending Ostrowski inequality in various directions for continuous
functions of selfadjoint operators in complex Hilbert spaces. Applications
for mid-point inequalities and some elementary functions of operators such
as the power function, the logarithmic and exponential functions are
provided as well.

From a complementary viewpoint to Ostrowski/mid-point inequalities,
trapezoidal type inequality provide a priory error bounds in approximating
the Riemann integral by a (generalized) trapezoidal formula.

Just like in the case of Ostrowski's inequality the development of these
kind of results have registered a sharp growth in the last decade with more
than 50 papers published, as one can easily asses this by performing a
search with the key word "trapezoid" and "inequality" in the title of the
papers reviewed by \textit{MathSciNet}.

Numerous extensions, generalizations in both the integral and discrete case
have been discovered. More general versions for $n$-time differentiable
functions, the corresponding versions on time scales, for vector valued
functions or multiple integrals have been established as well. Numerous
applications in Numerical Analysis, Probability Theory and other fields have
been also given.

In chapter five we present some recent results obtained by the author in
extending trapezoidal type inequality in various directions for continuous
functions of selfadjoint operators in complex Hilbert spaces. Applications
for some elementary functions of operators are provided as well.

In approximating $n$-time differentiable functions around a point, perhaps
the classical Taylor's expansion is one of the simplest and most convenient
and elegant methods that has been employed in the development of Mathematics
for the last three centuries.

In the sixth and last chapter of the book, we present some error bounds in
approximating $n$-time differentiable functions of selfadjoint operators by
the use of operator Taylor's type expansions around a point or two points
from its spectrum for which the remainder is known in an integral form. Some
applications for elementary functions including the exponential and
logarithmic functions are provided as well.

For the sake of completeness, all the results presented are completely
proved and the original references where they have been firstly obtained are
mentioned. The chapters are followed by the list of references used therein
and therefore are relatively independent and can be read separately.

\bigskip

\bigskip

The Author$^{\ast }$

\footnote[1]{$\ast $ This book is dedicated to my beloved children Sergiu \&
Camelia and granddaughter Sienna Clarisse.}

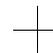
\chapter{Functions of Selfadjoint Operators in Hilbert Spaces}

\pagenumbering{arabic}

\section{Introduction}

In this introductory chapter we recall some fundamental facts concerning
bounded selfadjoint operators on complex Hilbert spaces. Since all the
operators considered in this book are supposed to be bounded, we no longer
mention this but understand it implicitly.

The generalized Schwarz's inequality for positive selfadjoint operators as
well as some results for the spectrum of this class of operators are
presented. Then we introduce and explore the fundamental results for
polynomials in a linear operator, continuous functions of selfadjoint
operators as well as the step functions of selfadjoint operators. By the use
of these results we then introduce the spectral decomposition of selfadjoint
operators (the \textit{Spectral Representation Theorem}) that will play a
central role in the rest of the book. This result is used as a key tool in
obtaining various new inequalities for continuous functions of selfadjoint
operators which are of bounded variation, Lipschitzian, monotonic or
absolutely continuous. Another tool that will greatly simplify the error
bounds provided in the book is the \textit{Total Variation Schwarz's
Inequality} for which a simple proof is offered.

The chapter is concluded with some well known operator inequalities of
Jensen's type for convex and operator convex functions. More results in this
spirit can be found in the recent book \cite{0.FMPS}.

Finally, some Gr\"{u}ss' type inequalities obtained in 1993 by Mond \& Pe%
\v{c}ari\'{c} are also presented. They are developed extensively in a
special chapter later in the book where some applications in relation with
classical power operator inequalities are provided as well.

\section{Bounded Selfadjoint Operators}

\subsection{Operator Order}

Let $(H;\left\langle .,.\right\rangle )$ be a Hilbert space over the complex
numbers field $\mathbb{C}.$

A bounded linear operator $A$ defined on $H$ is \textit{selfadjoint}, i.e., $%
A=A^{\ast }$ if and only if $\left\langle Ax,x\right\rangle \in \mathbb{R}$
for all $x\in H$ and if $A$ is selfadjoint, then%
\begin{equation}
\left\Vert A\right\Vert =\sup_{\left\Vert x\right\Vert =1}\left\vert
\left\langle Ax,x\right\rangle \right\vert =\sup_{\left\Vert x\right\Vert
=\left\Vert y\right\Vert =1}\left\vert \left\langle Ax,y\right\rangle
\right\vert .  \label{0.e.1.1}
\end{equation}

We assume in what follows that all operators are bounded on defined on the
whole Hilbert space $H.$ We denote by $\mathcal{B}\left( H\right) $ the
Banach algebra of all bounded linear operators defined on $H.$

\begin{definition}
\label{0.d.1.1}Let $A$ and $B$ be selfadjoint operators on $H.$ Then $A\leq
B $ ($A$ is less or equal to $B$) or, equivalently, $B\geq A$ if $%
\left\langle Ax,x\right\rangle \leq \left\langle Bx,x\right\rangle $ for all 
$x\in H.$ In particular, $A$ is called positive if $A\geq 0.$
\end{definition}

It is well known that for any operator $A\in \mathcal{B}\left( H\right) $
the composite operators $A^{\ast }A$ and $AA^{\ast }$ are positive
selfadjoint operators on $H$. However, the operators $A^{\ast }A$ and $%
AA^{\ast }$ are not comparable with each other in general.

The following result concerning the operator order holds (see for instance 
\cite[p. 220]{0.GH}):

\begin{theorem}
\label{0.t.1.1}Let $A,B,C\in \mathcal{B}\left( H\right) $ be selfadjoint
operators and let $\alpha ,\beta \in \mathbb{R}$. Then

\begin{enumerate}
\item $A\leq A;$

\item If $A\leq B$ and $B\leq C$, then $A\leq C;$

\item If $A\leq B$ and $B\leq A$, then $A=B;$

\item If $A\leq B$ and $\alpha \geq 0,$ then%
\begin{equation*}
A+C\leq B+C,\alpha A\leq \alpha B,-A\geq -B;
\end{equation*}

\item If $\alpha \leq \beta ,$ then $\alpha A\leq \beta A.$
\end{enumerate}
\end{theorem}

The following \textit{generalization of Schwarz's inequality} for positive
selfadjoint operators $A$ holds:

\begin{equation}
\left\vert \left\langle Ax,y\right\rangle \right\vert ^{2}\leq \left\langle
Ax,x\right\rangle \left\langle Ay,y\right\rangle  \label{0.GS}
\end{equation}%
for any $x,y\in H.$

The following inequality is of interest as well, see \cite[p. 221]{0.GH}

\begin{theorem}
\label{0.t.1.2}Let $A$ be a positive selfadjoint operator on $H.$ Then%
\begin{equation}
\left\Vert Ax\right\Vert ^{2}\leq \left\Vert A\right\Vert \left\langle
Ax,x\right\rangle  \label{0.e.1.2}
\end{equation}%
for any $x\in H.$
\end{theorem}

\begin{theorem}
\label{0.t.1.3}Let $A_{n},B\in \mathcal{B}\left( H\right) $ with $n\geq 1$
be selfadjoint operators with the property that%
\begin{equation*}
A_{1}\leq A_{2}\leq ...\leq A_{n}\leq ...\leq B.
\end{equation*}%
Then there exists a bounded selfadjoint operator $A$ defined on $H$ such that%
\begin{equation*}
A_{n}\leq A\leq B\text{ for all }n\geq 1
\end{equation*}%
and%
\begin{equation*}
\lim_{n\rightarrow \infty }A_{n}x=Ax\text{ for all }x\in H.
\end{equation*}
\end{theorem}

An analogous assertion holds if the sequence $\left\{ A_{n}\right\}
_{n=1}^{\infty }$ is decreasing and bounded below.

\begin{definition}
\label{0.d.1.2}We say that a sequence $\left\{ A_{n}\right\} _{n=1}^{\infty
}\subset \mathcal{B}\left( H\right) $ converges strongly to an operator $%
A\in \mathcal{B}\left( H\right) ,$ called the strong limit of the sequence $%
\left\{ A_{n}\right\} _{n=1}^{\infty }$ and we denote this by $\left(
s\right) \lim_{n\rightarrow \infty }A_{n}=A,$ if $\lim_{n\rightarrow \infty
}A_{n}x=Ax$ for all $x\in H.$
\end{definition}

The convergence in norm, i.e. lim$_{n\rightarrow \infty }\left\Vert
A_{n}-A\right\Vert =0$ will be called the \textit{"uniform convergence" }as
opposed to strong convergence. We denote lim$_{n\rightarrow \infty }A_{n}=A$
for the convergence in norm. From the inequality%
\begin{equation*}
\left\Vert A_{m}x-A_{n}x\right\Vert \leq \left\Vert A_{m}-A_{n}\right\Vert
\left\Vert x\right\Vert
\end{equation*}%
that holds for all $n,m$ and $x\in H$ it follows that uniform convergence of
the sequence $\left\{ A_{n}\right\} _{n=1}^{\infty }$ to $A$ implies strong
convergence of $\left\{ A_{n}\right\} _{n=1}^{\infty }$ to $A.$ However, the
converse of this assertion is false.

It is also possible to introduce yet another concept of \textit{"weak
convergence"} in $\mathcal{B}\left( H\right) $ by defining $\left( w\right)
\lim_{n\rightarrow \infty }A_{n}=A$ if and only if $\lim_{n\rightarrow
\infty }\left\langle A_{n}x,y\right\rangle =\left\langle Ax,y\right\rangle $
for all $x,y\in H.$

The following result holds (see \cite[p. 225]{0.GH}):

\begin{theorem}
\label{0.t.1.4}Let $A$ be a bounded selfadjoint operator on $H.$ Then%
\begin{eqnarray*}
\alpha _{1} &:&=\inf_{\left\Vert x\right\Vert =1}\left\langle
Ax,x\right\rangle =\max \left\{ \alpha \in \mathbb{R}\left\vert \alpha I\leq
A\right. \right\} ; \\
\alpha _{2} &:&=\sup_{\left\Vert x\right\Vert =1}\left\langle
Ax,x\right\rangle =\min \left\{ \alpha \in \mathbb{R}\left\vert A\leq \alpha
I\right. \right\} ;
\end{eqnarray*}%
and%
\begin{equation*}
\left\Vert A\right\Vert =\max \left\{ \left\vert \alpha _{1}\right\vert
,\left\vert \alpha _{2}\right\vert \right\} .
\end{equation*}%
Moreover, if $Sp\left( A\right) $ denotes the spectrum of $A,$ then $\alpha
_{1},\alpha _{2}\in Sp\left( A\right) $ and $Sp\left( A\right) \subset \left[
\alpha _{1},\alpha _{2}\right] .$
\end{theorem}

\begin{remark}
\label{0.r.1.1}We remark that, if $A,\alpha _{1},\alpha _{2}$ are as above,
then obviously%
\begin{eqnarray*}
\alpha _{1} &=&\min \left\{ \lambda \left\vert \lambda \in Sp\left( A\right)
\right. \right\} =:\min Sp\left( A\right) ; \\
\alpha _{2} &=&\max \left\{ \lambda \left\vert \lambda \in Sp\left( A\right)
\right. \right\} =:\max Sp\left( A\right) ; \\
\left\Vert A\right\Vert &=&\max \left\{ \left\vert \lambda \right\vert
\left\vert \lambda \in Sp\left( A\right) \right. \right\} .
\end{eqnarray*}%
We also observe that

\begin{enumerate}
\item $A$ is positive iff $\alpha _{1}\geq 0;$

\item $A$ is positive and invertible iff $\alpha _{1}>0;$

\item If $\alpha _{1}>0,$ then $A^{-1}$ is a positive selfadjoint operator
and $\min Sp\left( A^{-1}\right) =\alpha _{2}^{-1},\max Sp\left(
A^{-1}\right) =\alpha _{1}^{-1}.$
\end{enumerate}
\end{remark}

\section{Continuous Functions of Selfadjoint Operators}

\subsection{Polynomials in a Bounded Operator}

For two functions $\varphi ,\psi :\mathbb{C}\rightarrow \mathbb{C}$ we
adhere to the canonical notation:%
\begin{align*}
\left( \varphi +\psi \right) \left( s\right) & :=\varphi \left( s\right)
+\psi \left( s\right) , \\
\left( \lambda \varphi \right) \left( s\right) & :=\lambda \varphi \left(
s\right) , \\
\left( \varphi \psi \right) \left( s\right) & :=\varphi \left( s\right) \psi
\left( s\right)
\end{align*}%
for sum, scalar multiple and product of these functions. We denote by $\bar{%
\varphi}\left( s\right) $ the complex conjugate of $\varphi \left( s\right)
. $

As a first class of functions we consider the algebra $\mathcal{P}$ of all
polynomials in one variable with complex coefficients, namely%
\begin{equation*}
\mathcal{P}:=\left\{ \varphi \left( s\right) :=\sum_{k=0}^{n}\alpha
_{k}s^{k}\left\vert n\geq 0,\alpha _{k}\in \mathbb{C}\text{,}0\leq k\leq
n\right. \right\} .
\end{equation*}

\begin{theorem}
\label{0.t.1.5}Let $A\in \mathcal{B}\left( H\right) $ and for $\varphi
\left( s\right) :=\sum_{k=0}^{n}\alpha _{k}s^{k}\in \mathcal{P}$ define $%
\varphi \left( A\right) :=\sum_{k=0}^{n}\alpha _{k}A^{k}\in \mathcal{B}%
\left( H\right) \left( A^{0}=I\right) $ and $\bar{\varphi}\left( A\right)
:=\sum_{k=0}^{n}\bar{\alpha}_{k}\left( A^{\ast }\right) ^{k}\in \mathcal{B}%
\left( H\right) .$ Then the mapping $\varphi \left( s\right) \rightarrow
\varphi \left( A\right) $ has the following properties
\end{theorem}

\begin{enumerate}
\item[a)] $\left( \varphi +\psi \right) \left( A\right) =\varphi \left(
A\right) +\psi \left( A\right) ;$

\item[b)] $\left( \lambda \varphi \right) \left( A\right) =\lambda \varphi
\left( A\right) ;$

\item[c)] $\left( \varphi \psi \right) \left( A\right) =\varphi \left(
A\right) \psi \left( A\right) ;$

\item[d)] $\left[ \varphi \left( A\right) \right] ^{\ast }=\bar{\varphi}%
\left( A\right) .$
\end{enumerate}

Note that $\varphi \left( A\right) \psi \left( A\right) =\psi \left(
A\right) \varphi \left( A\right) $ and the constant polynomial $\varphi
\left( s\right) =\alpha _{0}$ is mapped into the operator.

Recall that, a mapping $a\rightarrow a^{\prime }$ of an algebra $\mathcal{U}$
into an algebra $\mathcal{U}^{\prime }$ is called a\textit{\ homomorphism}
if it has the properties

\begin{enumerate}
\item[a)] $\left( a+b\right) ^{\prime }=a^{\prime }+b^{\prime };$

\item[b)] $\left( \lambda \varphi \right) ^{\prime }=\lambda a^{\prime };$

\item[c)] $\left( ab\right) ^{\prime }=a^{\prime }b^{\prime }.$
\end{enumerate}

With this terminology, Theorem \ref{0.t.1.5} asserts that the mapping which
associates with any polynomial $\varphi \left( s\right) $ the operator $%
\varphi \left( A\right) $ is a homomorphism of $\mathcal{P}$ into $\mathcal{B%
}\left( H\right) $ satisfying the additional property d).

The following result provides a connection between the spectrum of $A$ and
the spectrum of the operator $\varphi \left( A\right) .$

\begin{theorem}
\label{0.t.1.6}If $A\in \mathcal{B}\left( H\right) $ and $\varphi \in 
\mathcal{P}$, then $Sp\left( \varphi \left( A\right) \right) =\varphi \left(
Sp\left( A\right) \right) .$
\end{theorem}

\begin{corollary}
\label{0.c.1.1}If $A\in \mathcal{B}\left( H\right) $ is selfadjoint and the
polynomial $\varphi \left( s\right) \in \mathcal{P}$ has real coefficients,
then $\varphi \left( A\right) $ is selfadjoint and 
\begin{equation}
\left\Vert \varphi \left( A\right) \right\Vert =\max \left\{ \left\vert
\varphi \left( \lambda \right) \right\vert ,\lambda \in Sp\left( A\right)
\right\} .  \label{0.e.1.3}
\end{equation}
\end{corollary}

\begin{remark}
\label{0.r.1.2} If $A\in \mathcal{B}\left( H\right) $ and $\varphi \in 
\mathcal{P}$, then

\begin{enumerate}
\item $\varphi \left( A\right) $ is invertible iff $\varphi \left( \lambda
\right) \neq 0$ for all $\lambda \in Sp\left( A\right) ;$

\item If $\varphi \left( A\right) $ is invertible, then $Sp\left( \varphi
\left( A\right) ^{-1}\right) =\left\{ \varphi \left( \lambda \right)
^{-1},\lambda \in Sp\left( A\right) \right\} .$
\end{enumerate}
\end{remark}

\subsection{Continuous Functions of Selfadjoint Operators}

Assume that $A$ is a bounded selfadjoint operator on the Hilbert space $H.$
If $\varphi $ is any function defined on $\mathbb{R}$ we define%
\begin{equation*}
\left\Vert \varphi \right\Vert _{A}=\sup \left\{ \left\vert \varphi \left(
\lambda \right) \right\vert ,\lambda \in Sp\left( A\right) \right\} .
\end{equation*}%
If $\varphi $ is continuous, in particular if $\varphi $ is a polynomial,
then the supremum is actually assumed for some points in $Sp\left( A\right) $
which is compact. Therefore the supremum may then be written as a maximum
and the formula (\ref{0.e.1.3}) can be written in the form $\left\Vert
\varphi \left( A\right) \right\Vert =\left\Vert \varphi \right\Vert _{A}.$

Consider $\mathcal{C}\left( \mathbb{R}\right) $ the algebra of all
continuous complex valued functions defined on $\mathbb{R}$. The following
fundamental result for continuous functional calculus holds, see for
instance \cite[p. 232]{0.GH}:

\begin{theorem}
\label{0.t.1.7}If $A$ is a bounded selfadjoint operator on the Hilbert space 
$H$ and $\varphi \in \mathcal{C}\left( \mathbb{R}\right) $, then there
exists a unique operator $\varphi \left( A\right) \in \mathcal{B}\left(
H\right) $ with the property that whenever $\left\{ \varphi _{n}\right\}
_{n=1}^{\infty }\subset \mathcal{P}$ such that lim$_{n\rightarrow \infty
}\left\Vert \varphi -\varphi _{n}\right\Vert _{A}=0,$ then $\varphi \left(
A\right) =\lim_{n\rightarrow \infty }\varphi _{n}\left( A\right) .$ The
mapping $\varphi \rightarrow \varphi \left( A\right) $ is a homomorphism of
the algebra $\mathcal{C}\left( \mathbb{R}\right) $ into $\mathcal{B}\left(
H\right) $ with the additional properties $\left[ \varphi \left( A\right) %
\right] ^{\ast }=\bar{\varphi}\left( A\right) $ and $\left\Vert \varphi
\left( A\right) \right\Vert \leq 2\left\Vert \varphi \right\Vert _{A}.$
Moreover, $\varphi \left( A\right) $ is a normal operator, i.e. $\left[
\varphi \left( A\right) \right] ^{\ast }\varphi \left( A\right) =\varphi
\left( A\right) \left[ \varphi \left( A\right) \right] ^{\ast }.$ If $%
\varphi $ is real-valued, then $\varphi \left( A\right) $ is selfadjoint.
\end{theorem}

As examples we notice that, if $A\in \mathcal{B}\left( H\right) $ is
selfadjoint and $\varphi \left( s\right) =e^{is},s\in \mathbb{R}$ then 
\begin{equation*}
e^{iA}=\sum_{k=0}^{\infty }\frac{1}{k!}\left( iA\right) ^{k}.
\end{equation*}%
Moreover, $e^{iA}$ is a unitary operator and its inverse is the operator%
\begin{equation*}
\left( e^{iA}\right) ^{\ast }=e^{-iA}=\sum_{k=0}^{\infty }\frac{1}{k!}\left(
-iA\right) ^{k}.
\end{equation*}

Now, if $\lambda \in \mathbb{C}\setminus \mathbb{R}$, $A\in \mathcal{B}%
\left( H\right) $ is selfadjoint and $\varphi \left( s\right) =\frac{1}{%
s-\lambda }\in \mathcal{C}\left( \mathbb{R}\right) ,$ then $\varphi \left(
A\right) =\left( A-\lambda I\right) ^{-1}.$

If the selfadjoint operator $A\in \mathcal{B}\left( H\right) $ and the
functions $\varphi ,\psi \in \mathcal{C}\left( \mathbb{R}\right) $ are
given, then we obtain the commutativity property $\varphi \left( A\right)
\psi \left( A\right) =\psi \left( A\right) \varphi \left( A\right) .$ This
property can be extended for another operator as follows, see for instance 
\cite[p. 235]{0.GH}:

\begin{theorem}
\label{0.t.1.8}Assume that $A\in \mathcal{B}\left( H\right) $ and the
function $\varphi \in \mathcal{C}\left( \mathbb{R}\right) $ are given. If $%
B\in \mathcal{B}\left( H\right) $ is such that $AB=BA,$ then $\varphi \left(
A\right) B=B\varphi \left( A\right) .$
\end{theorem}

The next result extends Theorem \ref{0.t.1.6} to the case of continuous
functions, see for instance \cite[p. 235]{0.GH}:

\begin{theorem}
\label{0.t.1.9}If $A$ is abounded selfadjoint operator on the Hilbert space $%
H$ and $\varphi $ is continuous, then $Sp\left( \varphi \left( A\right)
\right) =\varphi \left( Sp\left( A\right) \right) .$
\end{theorem}

As a consequence of this result we have:

\begin{corollary}
\label{0.c.1.2}With the assumptions in Theorem \ref{0.t.1.9} we have:

\begin{enumerate}
\item[a)] The operator $\varphi \left( A\right) $ is selfadjoint iff $%
\varphi \left( \lambda \right) \in \mathbb{R}$ for all $\lambda \in Sp\left(
A\right) ;$

\item[b)] The operator $\varphi \left( A\right) $ is unitary iff $\left\vert
\varphi \left( \lambda \right) \right\vert =1$ for all $\lambda \in Sp\left(
A\right) ;$

\item[c)] The operator $\varphi \left( A\right) $ is invertible iff $\varphi
\left( \lambda \right) \neq 0$ for all $\lambda \in Sp\left( A\right) ;$

\item[d)] If $\varphi \left( A\right) $ is selfadjoint, then $\left\Vert
\varphi \left( A\right) \right\Vert =\left\Vert \varphi \right\Vert _{A}.$
\end{enumerate}
\end{corollary}

In order to develop inequalities for functions of selfadjoint operators we
need the following result, see for instance \cite[p. 240]{0.GH}:

\begin{theorem}
\label{0.t.1.10}Let $A$ be a bounded selfadjoint operator on the Hilbert
space $H.$ The homomorphism $\varphi \rightarrow \varphi \left( A\right) $
of $\mathcal{C}\left( \mathbb{R}\right) $ into $\mathcal{B}\left( H\right) $
is order preserving, meaning that, if $\varphi ,\psi \in \mathcal{C}\left( 
\mathbb{R}\right) $ are real valued on $Sp\left( A\right) $ and $\varphi
\left( \lambda \right) \geq \psi \left( \lambda \right) $ for any $\lambda
\in Sp\left( A\right) ,$ then%
\begin{equation}
\varphi \left( A\right) \geq \psi \left( A\right) \text{ in the operator
order of }\mathcal{B}\left( H\right) .  \tag{P}  \label{P}
\end{equation}
\end{theorem}

The \textit{"square root"} of a positive bounded selfadjoint operator on $H$
can be defined as follows, see for instance \cite[p. 240]{0.GH}:

\begin{theorem}
\label{0.t.1.11}If the operator $A\in \mathcal{B}\left( H\right) $ is
selfadjoint and positive, then there exists a unique positive selfadjoint
operator $B:=\sqrt{A}\in \mathcal{B}\left( H\right) $ such that $B^{2}=A.$
If $A$ is invertible, then so is $B.$
\end{theorem}

If $A\in \mathcal{B}\left( H\right) ,$ then the operator $A^{\ast }A$ is
selfadjoint and positive. Define the \textit{"absolute value"} operator by $%
\left\vert A\right\vert :=\sqrt{A^{\ast }A}.$

Analogously to the familiar factorization of a complex number%
\begin{equation*}
\xi =\left\vert \xi \right\vert e^{i\arg \xi }
\end{equation*}%
a bounded normal operator on $H$ may be written as a commutative product of
a positive selfadjoint operator, representing its absolute value, and a
unitary operator, representing the factor of absolute value one.

In fact, the following more general result holds, see for instance \cite[p.
241]{0.GH}:

\begin{theorem}
\label{0.t.1.12}For every bounded linear operator $A$ on $H,$ there exists a
positive selfadjoint operator $B=\left\vert A\right\vert \in \mathcal{B}%
\left( H\right) $ and an isometric operator $C$ with the domain $\mathcal{D}%
_{C}=\overline{B\left( H\right) }$ and range $\mathcal{R}_{C}=C\left( 
\mathcal{D}_{C}\right) =\overline{A\left( H\right) }$ such that $A=CB.$
\end{theorem}

In particular, we have:

\begin{corollary}
\label{0.c.1.3}If the operator $A\in \mathcal{B}\left( H\right) $ is normal,
then there exists a positive selfadjoint operator $B=\left\vert A\right\vert
\in \mathcal{B}\left( H\right) $ and a unitary operator $C$ such that $%
A=BC=CB.$ Moreover, if $A$ is invertible, then $B$ and $C$ are uniquely
determined by these requirements.
\end{corollary}

\begin{remark}
\label{0.r.1.3}Now, suppose that $A=CB$ where $B\in \mathcal{B}\left(
H\right) $ is a positive selfadjoint operator and $C$ is an isometric
operator. Then

\begin{enumerate}
\item[a)] $B=\sqrt{A^{\ast }A};$ consequently $B$ is uniquely determined by
the stated requirements;

\item[b)] $C$ is uniquely determined by the stated requirements iff $A$ is
one-to-one.
\end{enumerate}
\end{remark}

\section{Step Functions of Selfadjoint Operators}

Let $A$ be a bonded selfadjoint operator on the Hilbert space $H.$ We intend
to extend the order preserving homomorphism $\varphi \rightarrow \varphi
\left( A\right) $ of the algebra $\mathcal{C}\left( \mathbb{R}\right) $ of
continuous functions $\varphi $ defined on $\mathbb{R}$ into $\mathcal{B}%
\left( H\right) ,$ restricted now to real-valued functions, to a larger
domain, namely an algebra of functions containing the "step functions" $%
\varphi _{\lambda },\lambda \in \mathbb{R}$, defined by%
\begin{equation*}
\varphi _{\lambda }\left( s\right) :=\left\{ 
\begin{array}{l}
1,\text{ for }-\infty <s\leq \lambda , \\ 
\\ 
0,\text{ for }\lambda <s<+\infty .%
\end{array}%
\right.
\end{equation*}%
Observe that $\overline{\varphi }_{\lambda }\left( s\right) =\varphi
_{\lambda }\left( s\right) $ and $\varphi _{\lambda }^{2}\left( s\right)
=\varphi _{\lambda }\left( s\right) $ which will imply that $\left[ \varphi
_{\lambda }\left( A\right) \right] ^{\ast }=\varphi _{\lambda }\left(
A\right) $ and $\left[ \varphi _{\lambda }\left( A\right) \right]
^{2}=\varphi _{\lambda }\left( A\right) ,$ i.e. $\varphi _{\lambda }\left(
A\right) $ will then be a projection. However, since the function $\varphi
_{\lambda }$ cannot be approximated uniformly by continuous functions on any
interval containing $\lambda ,$ then, in general, there is no way to define
an operator $\varphi _{\lambda }\left( A\right) $ as a uniform limit of
operators $\varphi _{\lambda ,n}\left( A\right) $ with $\varphi _{\lambda
,n}\in \mathcal{C}\left( \mathbb{R}\right) .$

The uniform limit of operators can be relaxed to the concept of strong limit
of operators (see Definition \ref{0.d.1.2}) in order to define the operator $%
\varphi _{\lambda }\left( A\right) .$ In order to do that, observe that the
function $\varphi _{\lambda }$ may be obtained as a pointwise limit of a
decreasing sequence of real-valued \ continuous functions $\varphi _{\lambda
,n}$ defined by%
\begin{equation*}
\varphi _{\lambda }\left( s\right) :=\left\{ 
\begin{array}{c}
1,\text{ for }-\infty <s\leq \lambda , \\ 
\\ 
1-n\left( s-\lambda \right) ,\text{ for }\lambda \leq s\leq \lambda +1/n \\ 
\\ 
0,\text{ for }\lambda <s<+\infty .%
\end{array}%
\right.
\end{equation*}%
By Theorem \ref{0.t.1.3} we observe that the sequence of corresponding
selfadjoint operators $\varphi _{\lambda ,n}\left( A\right) $ is
nondecreasing and bounded below by zero in the operator order of $\mathcal{B}%
\left( H\right) .$ It therefore converges strongly to some bounded
selfadjoint operator $\varphi _{\lambda }\left( A\right) $ on $H,$ see \cite[%
p. 244]{0.GH}.

To provide a formal presentation of the above, we need the following
definition.

\begin{definition}
\label{0.d.1.3}A real-valued function $\varphi $ on $\mathbb{R}$ is called
upper semi-continuous if it is a pointwise limit of a non-increasing
sequence of continuous real-valued functions on $\mathbb{R}$.
\end{definition}

We observe that it can be shown that a real-valued functions $\varphi $ on $%
\mathbb{R}$ is upper semi-continuous iff for every $s_{0}\in \mathbb{R}$ and
for every $\varepsilon >0$ there exists a $\delta >0$ such that 
\begin{equation*}
\varphi \left( s\right) <\varphi \left( s_{0}\right) +\varepsilon \text{ for
all }s\in \left( s_{0}-\delta ,s_{0}+\delta \right) .
\end{equation*}

We can introduce now the operator $\varphi \left( A\right) $ as follows, see
for instance \cite[p. 245]{0.GH}:

\begin{theorem}
\label{0.t.1.13}Let $A$ be a bonded selfadjoint operator on the Hilbert
space $H$ and let $\varphi $ be a nonnegative upper semi-continuous function
on $\mathbb{R}$. Then there exists a unique positive selfadjoint operator $%
\varphi \left( A\right) $ such that whenever $\left\{ \varphi _{n}\right\}
_{n=1}^{\infty }$ is any non-increasing sequence of non-negative functions
in $\mathcal{C}\left( \mathbb{R}\right) ,$ pointwise converging to $\varphi $
on $Sp\left( A\right) ,$ then $\varphi \left( A\right) =\left( s\right) \lim
\varphi _{n}\left( A\right) .$
\end{theorem}

If $\varphi $ is continuous, then the operator $\varphi \left( A\right) $
defined by Theorem \ref{0.t.1.7} coincides with the one defined by Theorem %
\ref{0.t.1.13}.

\begin{theorem}
\label{0.t.1.14}Let $A\in \mathcal{B}\left( H\right) $ be selfadjoint, let $%
\varphi $ and $\psi $ be non-negative upper semi-continuous functions on $%
\mathbb{R}$, and let $\alpha >0$ be given. Then the functions $\varphi +\psi
,$ $\alpha \varphi $ and $\varphi \psi $ are non-negative upper
semi-continuous and $\left( \varphi +\psi \right) \left( A\right) =\varphi
\left( A\right) +\psi \left( A\right) ,$ $\left( \alpha \varphi \right)
\left( A\right) =\alpha \varphi \left( A\right) $ and $\left( \varphi \psi
\right) \left( A\right) =\varphi \left( A\right) \psi \left( A\right) .$
Moreover, if $\varphi \left( s\right) \leq \psi \left( s\right) $ for all $%
s\in Sp\left( A\right) $ then $\varphi \left( A\right) \leq \psi \left(
A\right) .$
\end{theorem}

We enlarge the class of non-negative upper semi-continuous functions to an
algebra by defining $\mathcal{R}\left( \mathbb{R}\right) $ as the set of all
functions $\varphi =\varphi _{1}-\varphi _{2}$ where $\varphi _{1},\varphi
_{2}$ are nonnegative and upper semi-continuous functions defined on $%
\mathbb{R}$. It is easy to see that $\mathcal{R}\left( \mathbb{R}\right) $
endowed with pointwise sum, scalar multiple and product is an algebra.

The following result concerning functions of operators $\varphi \left(
A\right) $ with $\varphi \in \mathcal{R}\left( \mathbb{R}\right) $ can be
stated, see for instance \cite[p. 249-p. 250]{0.GH}:

\begin{theorem}
\label{0.t.1.15}Let $A\in \mathcal{B}\left( H\right) $ be selfadjoint and
let $\varphi \in \mathcal{R}\left( \mathbb{R}\right) .$ Then there exists a
unique selfadjoint operator $\varphi \left( A\right) \in \mathcal{B}\left(
H\right) $ such that if $\varphi =\varphi _{1}-\varphi _{2}$ where $\varphi
_{1},\varphi _{2}$ are nonnegative and upper semi-continuous functions
defined on $\mathbb{R}$, then $\varphi \left( A\right) =\varphi _{1}\left(
A\right) -\varphi _{2}\left( A\right) .$ The mapping $\varphi \rightarrow
\varphi \left( A\right) $ is a homomorphism of $\mathcal{R}\left( \mathbb{R}%
\right) $ into $\mathcal{B}\left( H\right) $ which is order preserving in
the following sense: if $\varphi ,\psi \in \mathcal{R}\left( \mathbb{R}%
\right) $ with the property that $\varphi \left( s\right) \leq \psi \left(
s\right) $ for any $s\in Sp\left( A\right) ,$ then $\varphi \left( A\right)
\leq \psi \left( A\right) .$ Moreover, if $B\in \mathcal{B}\left( H\right) $
satisfies the commutativity condition $AB=BA,$ then $\varphi \left( A\right)
B=B\varphi \left( A\right) .$
\end{theorem}

\section{The Spectral Decomposition of Selfadjoint Operators}

Let $A\in \mathcal{B}\left( H\right) $ be selfadjoint and let $\varphi
_{\lambda }$ defined for all $\lambda \in \mathbb{R}$ as follows%
\begin{equation*}
\varphi _{\lambda }\left( s\right) :=\left\{ 
\begin{array}{l}
1,\text{ for }-\infty <s\leq \lambda , \\ 
\\ 
0,\text{ for }\lambda <s<+\infty .%
\end{array}%
\right.
\end{equation*}%
Then for every $\lambda \in \mathbb{R}$ the operator 
\begin{equation}
E_{\lambda }:=\varphi _{\lambda }\left( A\right)  \label{Pr}
\end{equation}
is a projection which reduces $A.$

The properties of these projections are summed up in the following
fundamental result concerning the spectral decomposition of bounded
selfadjoint operators in Hilbert spaces, see for instance \cite[p. 256]{0.GH}

\begin{theorem}[Spectral Representation Theorem]
\label{0.t.1.16}Let $A$ be a bonded selfadjoint operator on the Hilbert
space $H$ and let $m=\min \left\{ \lambda \left\vert \lambda \in Sp\left(
A\right) \right. \right\} =:\min Sp\left( A\right) $ and $M=\max \left\{
\lambda \left\vert \lambda \in Sp\left( A\right) \right. \right\} =:\max
Sp\left( A\right) .$ Then there exists a family of projections $\left\{
E_{\lambda }\right\} _{\lambda \in \mathbb{R}}$, called the spectral family
of $A,$ with the following properties

\begin{enumerate}
\item[a)] $E_{\lambda }\leq E_{\lambda ^{\prime }}$ for $\lambda \leq
\lambda ^{\prime };$

\item[b)] $E_{m-0}=0,E_{M}=I$ and $E_{\lambda +0}=E_{\lambda }$ for all $%
\lambda \in \mathbb{R}$;

\item[c)] We have the representation%
\begin{equation}
A=\int_{m-0}^{M}\lambda dE_{\lambda }.  \label{0.e.1.4}
\end{equation}
\end{enumerate}

More generally, for every continuous complex-valued function $\varphi $
defined on $\mathbb{R}$ and for every $\varepsilon >0$ there exists a $%
\delta >0$ such that%
\begin{equation}
\left\Vert \varphi \left( A\right) -\sum_{k=1}^{n}\varphi \left( \lambda
_{k}^{\prime }\right) \left[ E_{\lambda _{k}}-E_{\lambda _{k-1}}\right]
\right\Vert \leq \varepsilon  \label{0.e.1.5}
\end{equation}%
whenever%
\begin{equation}
\left\{ 
\begin{array}{l}
\lambda _{0}<m=\lambda _{1}<...<\lambda _{n-1}<\lambda _{n}=M, \\ 
\\ 
\lambda _{k}-\lambda _{k-1}\leq \delta \text{ for }1\leq k\leq n, \\ 
\\ 
\lambda _{k}^{\prime }\in \left[ \lambda _{k-1},\lambda _{k}\right] \text{
for }1\leq k\leq n%
\end{array}%
\right.  \label{0.e.1.6}
\end{equation}%
this means that%
\begin{equation}
\varphi \left( A\right) =\int_{m-0}^{M}\varphi \left( \lambda \right)
dE_{\lambda },  \label{SAR}
\end{equation}%
where the integral is of Riemann-Stieltjes type.
\end{theorem}

\begin{corollary}
\label{0.c.1.5}With the assumptions of Theorem \ref{0.t.1.16} for $%
A,E_{\lambda }$ and $\varphi $ we have the representations%
\begin{equation}
\varphi \left( A\right) x=\int_{m-0}^{M}\varphi \left( \lambda \right)
dE_{\lambda }x\text{ \ for all }x\in H  \label{0.e.1.7}
\end{equation}%
and%
\begin{equation}
\left\langle \varphi \left( A\right) x,y\right\rangle =\int_{m-0}^{M}\varphi
\left( \lambda \right) d\left\langle E_{\lambda }x,y\text{ }\right\rangle 
\text{ \ for all }x,y\in H.  \label{0.e.1.8}
\end{equation}%
In particular,%
\begin{equation}
\left\langle \varphi \left( A\right) x,x\right\rangle =\int_{m-0}^{M}\varphi
\left( \lambda \right) d\left\langle E_{\lambda }x,x\text{ }\right\rangle 
\text{ \ for all }x\in H.  \label{0.e.1.9}
\end{equation}%
Moreover, we have the equality%
\begin{equation}
\left\Vert \varphi \left( A\right) x\right\Vert
^{2}=\int_{m-0}^{M}\left\vert \varphi \left( \lambda \right) \right\vert
^{2}d\left\Vert E_{\lambda }x\right\Vert ^{2}\text{ \ for all }x\in H.
\label{0.e.1.10}
\end{equation}
\end{corollary}

The next result shows that it is legitimate to talk about "\textit{the}"
spectral family of the bounded selfadjoint operator $A$ since it is uniquely
determined by the requirements a), b) and c) in Theorem \ref{0.t.1.16}, see
for instance \cite[p. 258]{0.GH}:

\begin{theorem}
\label{0.t.1.17}Let $A$ be a bonded selfadjoint operator on the Hilbert
space $H$ and let $m=\min Sp\left( A\right) $ and $M=\max Sp\left( A\right)
. $ If $\left\{ F_{\lambda }\right\} _{\lambda \in \mathbb{R}}$ is a family
of projections satisfying the requirements a), b) and c) in Theorem \ref%
{0.t.1.16}, then $F_{\lambda }=E_{\lambda }$ for all $\lambda \in \mathbb{R}$
where $E_{\lambda }$ is defined by (\ref{Pr}).
\end{theorem}

By the above two theorems, the spectral family $\left\{ E_{\lambda }\right\}
_{\lambda \in \mathbb{R}}$ uniquely determines and in turn is uniquely
determined by the bounded selfadjoint operator $A.$ The spectral family also
reflects in a direct way the properties of the operator $A$ as follows, see 
\cite[p. 263-p.266]{0.GH}

\begin{theorem}
\label{0.t.1.18}Let $\left\{ E_{\lambda }\right\} _{\lambda \in \mathbb{R}}$
be the spectral family of the bounded selfadjoint operator $A.$ If $B$ is a
bounded linear operator on $H$, then $AB=BA$ iff $E_{\lambda }B=BE_{\lambda
} $ for all $\lambda \in \mathbb{R}$. In particular $E_{\lambda
}A=AE_{\lambda }$ for all $\lambda \in \mathbb{R}$.
\end{theorem}

\begin{theorem}
\label{0.t.1.19}Let $\left\{ E_{\lambda }\right\} _{\lambda \in \mathbb{R}}$
be the spectral family of the bounded selfadjoint operator $A$ and $\mu \in 
\mathbb{R}$. Then

\begin{enumerate}
\item[a)] $\mu $ is a regular value of $A,$i.e., $A-\mu I$ is invertible iff
there exists a $\theta >0$ such that $E_{\mu -\theta }=E_{\mu +\theta };$

\item[b)] $\mu \in Sp\left( A\right) $ iff $E_{\mu -\theta }<E_{\mu +\theta
} $ for all $\theta >0;$

\item[c)] $\mu $ is an eigenvalue of $A$ iff $E_{\mu -0}<E_{\mu }.$
\end{enumerate}
\end{theorem}

The following result will play a key role in many results concerning
inequalities for bounded selfadjoint operators in Hilbert spaces. Since we
were not able to locate it in the literature, we will provide here a
complete proof:

\begin{theorem}[Total Variation Schwarz's Inequality]
\label{0.t.1.20}Let $\left\{ E_{\lambda }\right\} _{\lambda \in \mathbb{R}}$
be the spectral family of the bounded selfadjoint operator $A$ and let $%
m=\min Sp\left( A\right) $ and $M=\max Sp\left( A\right) .$ Then for any $%
x,y\in H$ the function $\lambda \rightarrow $ $\left\langle E_{\lambda
}x,y\right\rangle $ is of bounded variation on $\left[ m-s,M\right] ,$ for
any $s>0$ and we have the inequality%
\begin{equation}
\dbigvee\limits_{m-0}^{M}\left( \left\langle E_{\left( \cdot \right)
}x,y\right\rangle \right) \leq \left\Vert x\right\Vert \left\Vert
y\right\Vert ,  \tag{TVSI}  \label{TVSI}
\end{equation}%
where $\dbigvee\limits_{m-0}^{M}\left( \left\langle E_{\left( \cdot \right)
}x,y\right\rangle \right) $ denotes the limit $\lim_{s\rightarrow
0+}\dbigvee\limits_{m-s}^{M}\left( \left\langle E_{\left( \cdot \right)
}x,y\right\rangle \right) .$
\end{theorem}

\begin{proof}
If $P$ is a nonnegative selfadjoint operator on $H,$ i.e., $\left\langle
Px,x\right\rangle \geq 0$ for any $x\in H,$ then the following inequality is
a generalization of the Schwarz inequality in $H$%
\begin{equation}
\left\vert \left\langle Px,y\right\rangle \right\vert ^{2}\leq \left\langle
Px,x\right\rangle \left\langle Py,y\right\rangle ,  \label{0.e.1.11}
\end{equation}%
for any $x,y\in H.$

Now, if $d:m-s=t_{0}<t_{1}<...<t_{n-1}<t_{n}=M,$ where $s>0$ is an arbitrary
partition of the interval $\left[ m-s,M\right] ,$ then we have by Schwarz's
inequality for nonnegative operators (\ref{0.e.1.11}) that%
\begin{align}
& \dbigvee\limits_{m-s}^{M}\left( \left\langle E_{\left( \cdot \right)
}x,y\right\rangle \right)   \label{0.e.1.12} \\
& =\sup_{d}\left\{ \sum_{i=0}^{n-1}\left\vert \left\langle \left(
E_{t_{i+1}}-E_{t_{i}}\right) x,y\right\rangle \right\vert \right\}   \notag
\\
& \leq \sup_{d}\left\{ \sum_{i=0}^{n-1}\left[ \left\langle \left(
E_{t_{i+1}}-E_{t_{i}}\right) x,x\right\rangle ^{1/2}\left\langle \left(
E_{t_{i+1}}-E_{t_{i}}\right) y,y\right\rangle ^{1/2}\right] \right\} :=I. 
\notag
\end{align}%
By the Cauchy-Buniakovski-Schwarz inequality for sequences of real numbers
we also have that%
\begin{align}
I& \leq \sup_{d}\left\{ \left[ \sum_{i=0}^{n-1}\left\langle \left(
E_{t_{i+1}}-E_{t_{i}}\right) x,x\right\rangle \right] ^{1/2}\left[
\sum_{i=0}^{n-1}\left\langle \left( E_{t_{i+1}}-E_{t_{i}}\right)
y,y\right\rangle \right] ^{1/2}\right\}   \label{0.e.1.13} \\
& \leq \sup_{d}\left\{ \left[ \sum_{i=0}^{n-1}\left\langle \left(
E_{t_{i+1}}-E_{t_{i}}\right) x,x\right\rangle \right] ^{1/2}\sup_{d}\left[
\sum_{i=0}^{n-1}\left\langle \left( E_{t_{i+1}}-E_{t_{i}}\right)
y,y\right\rangle \right] ^{1/2}\right\}   \notag \\
& =\left[ \dbigvee\limits_{m-s}^{M}\left( \left\langle E_{\left( \cdot
\right) }x,x\right\rangle \right) \right] ^{1/2}\left[ \dbigvee%
\limits_{m-s}^{M}\left( \left\langle E_{\left( \cdot \right)
}y,y\right\rangle \right) \right] ^{1/2}  \notag \\
& =\left[ \left\Vert x\right\Vert ^{2}-\left\langle E_{m-s}x,x\right\rangle %
\right] ^{1/2}\left[ \left\Vert y\right\Vert ^{2}-\left\langle
E_{m-s}y,y\right\rangle \right] ^{1/2}  \notag
\end{align}%
for any $x,y\in H.$

On making use of (\ref{0.e.1.12}) and (\ref{0.e.1.13}) and letting $%
s\rightarrow 0+$ we deduce the desired result (\ref{TVSI}).
\end{proof}

\section{Jensen's Type Inequalities}

\subsection{Jensen's Inequality.}

The following result that provides an operator version for the \textit{%
Jensen inequality} is due to Mond \& Pe\v{c}ari\'{c} \cite{0.MP} (see also 
\cite[p. 5]{0.FMPS}):

\begin{theorem}[Mond- Pe\v{c}ari\'{c}, 1993, \protect\cite{0.MP}]
\label{0.t.2.1} Let $A$ be a selfadjoint operator on the Hilbert space $H$
and assume that $Sp\left( A\right) \subseteq \left[ m,M\right] $ for some
scalars $m,M$ with $m<M.$ If $f$ is a convex function on $\left[ m,M\right]
, $ then 
\begin{equation}
f\left( \left\langle Ax,x\right\rangle \right) \leq \left\langle f\left(
A\right) x,x\right\rangle  \tag{MP}  \label{0.MP}
\end{equation}%
for each $x\in H$ with $\left\Vert x\right\Vert =1.$
\end{theorem}

As a special case of Theorem \ref{0.t.2.1} we have the following \textit{H%
\"{o}lder-McCarthy inequality}:

\begin{theorem}[H\"{o}lder-McCarthy, 1967, \protect\cite{0.Mc}]
\label{0.t.2.2} Let $A$ be a selfadjoint positive operator on a Hilbert
space $H$. Then

(i) \ \ $\left\langle A^{r}x,x\right\rangle \geq \left\langle
Ax,x\right\rangle ^{r}$ for all $r>1$ and $x\in H$ with $\left\| x\right\|
=1;$

(ii) \ $\left\langle A^{r}x,x\right\rangle \leq \left\langle
Ax,x\right\rangle ^{r}$ for all $0<r<1$ and $x\in H$ with $\left\| x\right\|
=1;$

(iii) If $A$ is invertible, then $\left\langle A^{r}x,x\right\rangle \geq
\left\langle Ax,x\right\rangle ^{r}$ for all $r<0$ and $x\in H$ with $%
\left\| x\right\| =1.$
\end{theorem}

The following theorem is a multiple operator version of Theorem \ref{0.t.2.1}
(see for instance \cite[p. 5]{0.FMPS}):

\begin{theorem}[Furuta-Mi\'{c}i\'{c}-Pe\v{c}ari\'{c}-Seo, 2005, \protect\cite%
{0.FMPS}]
\label{0.t.2.3} Let $A_{j}$ be selfadjoint operators with $Sp\left(
A_{j}\right) \subseteq \left[ m,M\right] $, $j\in \left\{ 1,\dots ,n\right\} 
$ for some scalars $m<M$ and $x_{j}\in H,j\in \left\{ 1,\dots ,n\right\} $
with $\sum_{j=1}^{n}\left\Vert x_{j}\right\Vert ^{2}=1$. If $f$ is a convex
function on $\left[ m,M\right] $, then 
\begin{equation}
f\left( \sum_{j=1}^{n}\left\langle A_{j}x_{j},x_{j}\right\rangle \right)
\leq \sum_{j=1}^{n}\left\langle f\left( A_{j}\right)
x_{j},x_{j}\right\rangle .  \label{0.e.2.1}
\end{equation}
\end{theorem}

The following particular case is of interest.

\begin{corollary}
\label{0.c.2.1}Let $A_{j}$ be selfadjoint operators with $Sp\left(
A_{j}\right) \subseteq \left[ m,M\right] $, $j\in \left\{ 1,\dots ,n\right\} 
$ for some scalars $m<M.$ If $p_{j}\geq 0,$ $j\in \left\{ 1,\dots ,n\right\} 
$ with $\sum_{j=1}^{n}p_{j}=1,$ then 
\begin{equation}
f\left( \left\langle \sum_{j=1}^{n}p_{j}A_{j}x,x\right\rangle \right) \leq
\left\langle \sum_{j=1}^{n}p_{j}f\left( A_{j}\right) x,x\right\rangle ,
\label{0.e.2.2}
\end{equation}%
for any $x\in H$ with $\left\Vert x\right\Vert =1.$
\end{corollary}

\begin{proof}
Follows from Theorem \ref{0.t.2.3} by choosing $x_{j}=\sqrt{p_{j}}\cdot x,$ $%
j\in \left\{ 1,\dots ,n\right\} $ where $x\in H$ with $\left\Vert
x\right\Vert =1.$
\end{proof}

\begin{remark}
\label{0.r.2.1}The above inequality can be used to produce some norm
inequalities for the sum of positive operators in the case when the convex
function $f$ is nonnegative and monotonic nondecreasing on $\left[ 0,M\right]
.$ Namely, we have: 
\begin{equation}
f\left( \left\Vert \sum_{j=1}^{n}p_{j}A_{j}\right\Vert \right) \leq
\left\Vert \sum_{j=1}^{n}p_{j}f\left( A_{j}\right) \right\Vert .
\label{0.e.2.2.a}
\end{equation}%
The inequality (\ref{0.e.2.2.a}) reverses if the function is concave on $%
\left[ 0,M\right] $.

As particular cases we can state the following inequalities: 
\begin{equation}
\left\Vert \sum_{j=1}^{n}p_{j}A_{j}\right\Vert ^{p}\leq \left\Vert
\sum_{j=1}^{n}p_{j}A_{j}^{p}\right\Vert ,  \label{0.e.2.2.b}
\end{equation}%
for $p\geq 1$ and 
\begin{equation}
\left\Vert \sum_{j=1}^{n}p_{j}A_{j}\right\Vert ^{p}\geq \left\Vert
\sum_{j=1}^{n}p_{j}A_{j}^{p}\right\Vert  \label{0.e.2.2.c}
\end{equation}%
for $0<p<1.$

If $A_{j}$ are positive definite for each $j\in \left\{ 1,\dots ,n\right\} $
then (\ref{0.e.2.2.b}) also holds for $p<0.$

If one uses the inequality (\ref{0.e.2.2.a}) for the exponential function,
that one obtains the inequality 
\begin{equation}
\exp \left( \left\Vert \sum_{j=1}^{n}p_{j}A_{j}\right\Vert \right) \leq
\left\Vert \sum_{j=1}^{n}p_{j}\exp \left( A_{j}\right) \right\Vert ,
\label{0.e.2.2.d}
\end{equation}%
where $A_{j}$ are positive operators for each $j\in \left\{ 1,\dots
,n\right\} .$
\end{remark}

\subsection{Reverses of Jensen's Inequality}

In Section 2.4 of the monograph \cite{0.FMPS} there are numerous interesting
converses of the Jensen's type inequality (\ref{0.e.2.1}) from which we
would like to mention only two of the simplest.

The following result is an operator version of the well known Lah-Ribari\'{c}%
's reverse of the Jensen inequality for real functions of a real variable,
see for instance \cite{0.FMPS}:

\begin{theorem}
\label{0.t.2.5}Let $A_{j}$ be selfadjoint operators with $Sp\left(
A_{j}\right) \subseteq \left[ m,M\right] $, $j\in \left\{ 1,\dots ,n\right\} 
$ for some scalars $m<M$ and $x_{j}\in H,j\in \left\{ 1,\dots ,n\right\} $
with $\sum_{j=1}^{n}\left\Vert x_{j}\right\Vert ^{2}=1$. If $f$ is a
continuous convex function defined on $\left[ m,M\right] ,$ then%
\begin{align}
& \sum_{j=1}^{n}\left\langle f\left( A_{j}\right) x_{j},x_{j}\right\rangle
\label{I.e.1.4} \\
& \leq \frac{1}{M-m}\left[ f\left( M\right) \sum_{j=1}^{n}\left\langle
\left( A_{j}-mI\right) x_{j},x_{j}\right\rangle +f\left( m\right)
\sum_{j=1}^{n}\left\langle \left( MI-A_{j}\right) x_{j},x_{j}\right\rangle %
\right] .  \notag
\end{align}
\end{theorem}

\begin{theorem}[Mi\'{c}i\'{c}-Seo-Takahasi-Tominaga, 1999, \protect\cite%
{0.MSTT}]
\label{0.t.2.4}Let $A_{j}$ be selfadjoint operators with $Sp\left(
A_{j}\right) \subseteq \left[ m,M\right] $, $j\in \left\{ 1,\dots ,n\right\} 
$ for some scalars $m<M$ and $x_{j}\in H,j\in \left\{ 1,\dots ,n\right\} $
with $\sum_{j=1}^{n}\left\Vert x_{j}\right\Vert ^{2}=1$. If $f$ is a
strictly convex function twice differentiable on $\left[ m,M\right] $, then
for any positive real number $\alpha $ we have 
\begin{equation}
\sum_{j=1}^{n}\left\langle f\left( A_{j}\right) x_{j},x_{j}\right\rangle
\leq \alpha f\left( \sum_{j=1}^{n}\left\langle A_{j}x_{j},x_{j}\right\rangle
\right) +\beta ,  \label{I.e.1.3}
\end{equation}%
where 
\begin{gather*}
\beta =\mu _{f}t_{0}+\nu _{f}-\alpha f\left( t_{0}\right) , \\
\mu _{f}=\frac{f\left( M\right) -f\left( m\right) }{M-m},\quad \nu _{f}=%
\frac{Mf\left( m\right) -mf\left( M\right) }{M-m}
\end{gather*}%
and 
\begin{equation*}
t_{0}=\left\{ 
\begin{array}{ll}
f^{\prime -1}\left( \frac{\mu _{f}}{\alpha }\right) & \text{if }m<f^{\prime
-1}\left( \frac{\mu _{f}}{\alpha }\right) <M \\ 
&  \\ 
M & \text{if }M\leq f^{\prime -1}\left( \frac{\mu _{f}}{\alpha }\right) \\ 
&  \\ 
m & \text{if }f^{\prime -1}\left( \frac{\mu _{f}}{\alpha }\right) \leq m.%
\end{array}%
\right.
\end{equation*}
\end{theorem}

The case of equality was also analyzed, see \cite[p. 61]{0.FMPS} but will be
not stated in here.

\subsection{Operator Monotone and Operator Convex Functions}

We say that a real valued continuous function $f$ defined on an interval $I$
is said to be \textit{operator monotone} if it is monotone with respect to
the operator order, i.e. if $A$ and \ $B$ are bounded selfadjoint operators
with $A\leq B$ and $Sp\left( A\right) ,Sp\left( A\right) \subset I,$ then $%
f\left( A\right) \leq f\left( B\right) .$ The function is said to be \textit{%
operator convex} (\textit{operator concave}) if for any $A$,\ $B$ bounded
selfadjoint operators with $Sp\left( A\right) ,Sp\left( A\right) \subset I,$
we have%
\begin{equation}
f\left[ \left( 1-\lambda \right) A+\lambda B\right] \leq \left( \geq \right)
\left( 1-\lambda \right) f\left( A\right) +\lambda f\left( B\right)
\label{0.e.2.3}
\end{equation}%
for any $\lambda \in \left[ 0,1\right] .$

\begin{example}
\label{0.ex.1.1}The following examples are well know in the literature and
can be found for instance in \cite[p. 7-p. 9]{0.FMPS} where simple proofs
were also provided.

\begin{enumerate}
\item The affine function $f\left( t\right) =\alpha +\beta t$ is operator
monotone on every interval for all $\alpha \in \mathbb{R}$ and $\beta \geq
0. $ It is operator convex for all $\alpha ,\beta \in \mathbb{R}$;

\item If $f,g$ are operator monotone, and if $\alpha ,\beta \geq 0$ then the
linear combination $\alpha f+\beta g$ is also operator monotone. If the
functions $f_{n}$ are operator monotone and $f_{n}\left( t\right)
\rightarrow f\left( t\right) $ as $n\rightarrow \infty ,$ then $f$ is also
operator monotone;

\item The function $f\left( t\right) =t^{2}$ is operator convex on every
interval, however it is not operator monotone on $[0,\infty )$ even though
it is monotonic nondecreasing on this interval;

\item The function $f\left( t\right) =t^{3}$ is not operator convex on $%
[0,\infty )$ even though it is a convex function on this interval;

\item The function $f\left( t\right) =\frac{1}{t}$ is operator convex on $%
\left( 0,\infty \right) $ and $f\left( t\right) =-\frac{1}{t}$ is operator
monotone on $\left( 0,\infty \right) ;$

\item The function $f\left( t\right) =\ln t$ is operator monotone and
operator concave on $\left( 0,\infty \right) ;$

\item The entropy function $f\left( t\right) =-t\ln t$ is operator concave
on $\left( 0,\infty \right) ;$

\item The exponential function $f\left( t\right) =e^{t}$ is neither operator
convex nor operator monotone on any interval of $\mathbb{R}$.
\end{enumerate}
\end{example}

The following monotonicity property for the function $f\left( t\right)
=t^{r} $ with $r\in \left[ 0,1\right] $ is well known in the literature as
the \textit{L\"{o}wner-Heinz inequality} and was established essentially in
1934:

\begin{theorem}[L\"{o}wner-Heinz Inequality]
Let $A$ and $B$ be positive operators on a Hilbert space $H.$ If $A\geq
B\geq 0,$ then $A^{r}\geq B^{r}$ for all $r\in \left[ 0,1\right] .$
\end{theorem}

The following characterization of operator convexity holds, see \cite[p. 10]%
{0.FMPS}

\begin{theorem}[Jensen's Operator Inequality]
\label{0.t.2.6}Let $H$ and $K$ be Hilbert spaces. Let $f$ be a real valued
continuous function on an interval $J.$ Let $A$ and $A_{j}$ be selfadjoint
operators on $H$ with spectra contained in $J,$ for each $j=1,2,...,k.$ Then
the following conditions are mutually equivalent:

\begin{enumerate}
\item[(i)] $f$ is operator convex on $J;$

\item[(ii)] $f\left( C^{\ast }AC\right) \leq C^{\ast }f\left( A\right) C$
for every selfadjoint operator $A:H\rightarrow H$ and isometry $%
C:K\rightarrow H,i.e.,C^{\ast }C=1_{K};$

\item[(iii)] $f\left( C^{\ast }AC\right) \leq C^{\ast }f\left( A\right) C$
for every selfadjoint operator $A:H\rightarrow H$ and isometry $%
C:H\rightarrow H;$

\item[(iv)] $f\left( \sum_{j=1}^{k}C_{j}^{\ast }A_{j}C_{j}\right) \leq
\sum_{j=1}^{k}C_{j}^{\ast }f\left( A_{j}\right) C_{j}$ for every selfadjoint
operator $A_{j}:H\rightarrow H$ and bounded linear operators $%
C_{j}:K\rightarrow H,$ with $\sum_{j=1}^{k}C_{j}^{\ast }C_{j}=1_{K}\left(
j=1,...,k\right) ;$

\item[(v)] $f\left( \sum_{j=1}^{k}C_{j}^{\ast }A_{j}C_{j}\right) \leq
\sum_{j=1}^{k}C_{j}^{\ast }f\left( A_{j}\right) C_{j}$ for every selfadjoint
operator $A_{j}:H\rightarrow H$ and bounded linear operators $%
C_{j}:H\rightarrow H,$ with $\sum_{j=1}^{k}C_{j}^{\ast }C_{j}=1_{H}\left(
j=1,...,k\right) ;$

\item[(vi)] $f\left( \sum_{j=1}^{k}P_{j}A_{j}P_{j}\right) \leq
\sum_{j=1}^{k}P_{j}f\left( A_{j}\right) P_{j}$ for every selfadjoint
operator $A_{j}:H\rightarrow H$ and projection $P_{j}:H\rightarrow H,$ with $%
\sum_{j=1}^{k}P_{j}=1_{H}\left( j=1,...,k\right) .$
\end{enumerate}
\end{theorem}

The following well known result due to Hansen \& Pedersen also holds:

\begin{theorem}[Hansen-Pedersen-Jensen's Inequality]
Let $J$ be an interval containing $0$ and let $f$ be a real valued
continuous function defined on $J.$ Let $A$ and $A_{j}$ be selfadjoint
operators on $H$ with spectra contained in $J,$ for each $j=1,2,...,k.$ Then
the following conditions are mutually equivalent:

\begin{enumerate}
\item[(i)] $f$ is operator convex on $J$ and $f\left( 0\right) \leq 0;$

\item[(ii)] $f\left( C^{\ast }AC\right) \leq C^{\ast }f\left( A\right) C$
for every selfadjoint operator $A:H\rightarrow H$ and contraction $%
C:H\rightarrow H,i.e.,C^{\ast }C\leq 1_{H};$

\item[(iii)] $f\left( \sum_{j=1}^{k}C_{j}^{\ast }A_{j}C_{j}\right) \leq
\sum_{j=1}^{k}C_{j}^{\ast }f\left( A_{j}\right) C_{j}$ for every selfadjoint
operator $A_{j}:H\rightarrow H$ and bounded linear operators $%
C_{j}:H\rightarrow H,$ with $\sum_{j=1}^{k}C_{j}^{\ast }C_{j}\leq
1_{H}\left( j=1,...,k\right) ;$

\item[(iv)] $f\left( PAP\right) \leq Pf\left( A\right) P$ for every
selfadjoint operator $A:H\rightarrow H$ and projection $P.$
\end{enumerate}
\end{theorem}

The case of continuous and negative functions is as follows, \cite[p. 13]%
{0.FMPS}:

\begin{theorem}
\label{0.t.2.7}Let $f$ be continuous on $[0,\infty ).$ If $f\left( t\right)
\leq 0$ for all $t\in \lbrack 0,\infty ),$ then each of the conditions
(i)-(vi) from Theorem \ref{0.t.2.6} is equivalent with

\begin{enumerate}
\item[(vii)] $-f$ is an operator monotone function.
\end{enumerate}
\end{theorem}

\begin{corollary}
\label{0.c.2.2}Let $f$ be a real valued continuous function mapping the
positive half line $[0,\infty )$ into itself. Then $f$ is operator monotone
if and only if $f$ is operator concave.
\end{corollary}

The following result may be stated as well \cite[p. 14]{0.FMPS}:

\begin{theorem}
\label{0.t.2.8}Let $f$ be continuous on the interval $[0,r)$ with $r\leq
\infty .$ Then the following conditions are mutually equivalent:

\begin{enumerate}
\item[(i)] $f$ is operator convex and $f\left( 0\right) \leq 0;$

\item[(ii)] The function $t\mapsto \frac{f\left( t\right) }{t}$ is operator
monotone on $\left( 0,r\right) .$
\end{enumerate}
\end{theorem}

As a particular case of interest, we can state that \cite[p. 15]{0.FMPS}:

\begin{corollary}
\label{0.c.2.3}Let $f$ be continuous on $[0,\infty )$ and taking positive
values. The function $f$ is operator monotone if and only if the function $%
t\mapsto \frac{t}{f\left( t\right) }$ is operator monotone.
\end{corollary}

Finally we recall the following result as well \cite[p. 16]{0.FMPS}:

\begin{theorem}
\label{0.t.2.9}Let $f$ be a real valued continuous function on the interval $%
J=[\alpha ,\infty )$ and bounded below, i.e., there exists $m\in \mathbb{R}$
such that $m\leq f\left( t\right) $ for all $t\in J.$ Then the following
conditions are mutually equivalent:

\begin{enumerate}
\item[(i)] $f$ is operator concave on $J;$

\item[(ii)] $f$ is operator monotone on $J.$
\end{enumerate}
\end{theorem}

As a particular case of this result we note that, the function $f\left(
t\right) =t^{r}$ is operator monotone on $[0,\infty )$ if and only if $0\leq
r\leq 1.$ The function $f\left( t\right) =t^{r}$ is operator convex on $%
\left( 0,\infty \right) $ if either $1\leq r\leq 2$ or $-1\leq r\leq 0$ and
is operator concave on $\left( 0,\infty \right) $ if $0\leq r\leq 1.$

\section{Gr\"{u}ss' Type Inequalities}

The following operator version of the Gr\"{u}ss inequality was obtained by
Mond \& Pe\v{c}ari\'{c} in \cite{0.MP1}:

\begin{theorem}[Mond-Pe\v{c}ari\'{c}, 1993, \protect\cite{0.MP1}]
\label{0.t.3.1}Let $C_{j},$ $j\in \left\{ 1,\dots ,n\right\} $ be
selfadjoint operators on the Hilbert space $\left( H,\left\langle
.,.\right\rangle \right) $ and such that $m_{j}\cdot 1_{H}\leq C_{j}\leq
M_{j}\cdot 1_{H}$ for $j\in \left\{ 1,\dots ,n\right\} ,$ where $1_{H}$ is
the identity operator on $H.$ Further, let $g_{j},h_{j}:\left[ m_{j},Mj%
\right] \rightarrow \mathbb{R}$, $j\in \left\{ 1,\dots ,n\right\} $ be
functions such that 
\begin{equation}
\varphi \cdot 1_{H}\leq g_{j}\left( C_{j}\right) \leq \Phi \cdot 1_{H}\text{
\quad and \quad }\gamma \cdot 1_{H}\leq h_{j}\left( C_{j}\right) \leq \Gamma
\cdot 1_{H}  \label{0.e.3.1}
\end{equation}%
for each $j\in \left\{ 1,\dots ,n\right\} .$

If $x_{j}\in H,$ $j\in \left\{ 1,\dots ,n\right\} $ are such that $%
\sum_{j=1}^{n}\left\Vert x_{j}\right\Vert ^{2}=1,$ then 
\begin{align}
& \left\vert \sum_{j=1}^{n}\left\langle g_{j}\left( C_{j}\right) h_{j}\left(
C_{j}\right) x_{j},x_{j}\right\rangle -\sum_{j=1}^{n}\left\langle
g_{j}\left( C_{j}\right) x_{j},x_{j}\right\rangle \cdot
\sum_{j=1}^{n}\left\langle h_{j}\left( C_{j}\right) x_{j},x_{j}\right\rangle
\right\vert  \label{0.e.3.2} \\
& \leq \frac{1}{4}\left( \Phi -\varphi \right) \left( \Gamma -\gamma \right)
.  \notag
\end{align}
\end{theorem}

If $C_{j},j\in \left\{ 1,\dots ,n\right\} $ are selfadjoint operators such
that $Sp\left( C_{j}\right) \subseteq \left[ m,M\right] $ for $j\in \left\{
1,\dots ,n\right\} $ and for some scalars $m<M$ and if $g,h:\left[ m,M\right]
\longrightarrow \mathbb{R}$ are continuous then by the Mond-Pe\v{c}ari\'{c}
inequality we deduce the following version of the Gr\"{u}ss inequality for
operators 
\begin{align}
& \left\vert \sum_{j=1}^{n}\left\langle g\left( C_{j}\right) h\left(
C_{j}\right) x_{j},x_{j}\right\rangle -\sum_{j=1}^{n}\left\langle g\left(
C_{j}\right) x_{j},x_{j}\right\rangle \cdot \sum_{j=1}^{n}\left\langle
h\left( C_{j}\right) x_{j},x_{j}\right\rangle \right\vert  \label{0.e.3.3} \\
& \leq \frac{1}{4}\left( \Phi -\varphi \right) \left( \Gamma -\gamma \right)
,  \notag
\end{align}%
where $x_{j}\in H,$ $j\in \left\{ 1,\dots ,n\right\} $ are such that $%
\sum_{j=1}^{n}\left\Vert x_{j}\right\Vert ^{2}=1$ and $\varphi =\min_{t\in %
\left[ m,M\right] }g\left( t\right) ,$ $\Phi =\max_{t\in \left[ m,M\right]
}g\left( t\right) ,$ $\gamma =\min_{t\in \left[ m,M\right] }h\left( t\right) 
$ and $\Gamma =\max_{t\in \left[ m,M\right] }h\left( t\right) .$

In particular, if the selfadjoint operator $C$ satisfy the condition $%
Sp\left( C\right) \subseteq \left[ m,M\right] $ for some scalars $m<M$, then 
\begin{equation}
\left\vert \left\langle g\left( C\right) h\left( C\right) x,x\right\rangle
-\left\langle g\left( C\right) x,x\right\rangle \cdot \left\langle h\left(
C\right) x,x\right\rangle \right\vert \leq \frac{1}{4}\left( \Phi -\varphi
\right) \left( \Gamma -\gamma \right) ,  \label{0.e.3.4}
\end{equation}%
for any $x\in H$ with $\left\Vert x\right\Vert =1.$

\bigskip

\chapter{Inequalities for Convex Functions}

\section{Introduction}

Jensen's type inequalities in their various settings ranging from discrete
to continuous case play an important role in different branches of Modern
Mathematics. A simple search in the \textit{MathSciNet} database of the
American Mathematical Society with the key words "jensen" and "inequality"
in the title reveals more than 300 items intimately devoted to this famous
result. However, the number of papers where this inequality is applied is a
lot larger and far more difficult to find. It can be a good project in
itself for someone to write a monograph devoted to Jensen's inequality in
its different forms and its applications across Mathematics.

In the introductory chapter we have recalled a number of Jensen's type
inequalities for convex and operator convex functions of selfadjoint
operators in Hilbert spaces. In this chapter we present some recent results
obtained by the author that deal with different aspects of this well
research inequality than those recently reported in the book \cite{I.FMPS}.
They include but are not restricted to the operator version of the
Dragomir-Ionescu inequality, Slater's type inequalities for operators and
its inverses, Jensen's inequality for twice differentiable functions whose
second derivatives satisfy some upper and lower bounds conditions, Jensen's
type inequalities for log-convex functions and for differentiable log-convex
functions and their applications to Ky Fan's inequality.

Finally, some Hermite-Hadamard's type inequalities for convex functions and
Hermite-Hadamard's type inequalities for operator convex functions are
presented as well.

All the above results are exemplified for some classes of elementary
functions of interest such as the power function and the logarithmic
function.

\section{Reverses of the Jensen Inequality}

\subsection{An Operator Version of the Dragomir-Ionescu Inequality}

The following result holds:

\begin{theorem}[Dragomir, 2008, \protect\cite{I.SSD5}]
\label{I.t.2.1}Let $I$ be an interval and $f:I\rightarrow \mathbb{R}$ be a
convex and differentiable function on \r{I} (the interior of $I)$ whose
derivative $f^{\prime }$ is continuous on \r{I} $.$ If $A$ is a selfadjoint
operators on the Hilbert space $H$ with $Sp\left( A\right) \subseteq \left[
m,M\right] \subset $\r{I}$,$ then 
\begin{equation}
(0\leq )\left\langle f\left( A\right) x,x\right\rangle -f\left( \left\langle
Ax,x\right\rangle \right) \leq \left\langle f^{\prime }\left( A\right)
Ax,x\right\rangle -\left\langle Ax,x\right\rangle \cdot \left\langle
f^{\prime }\left( A\right) x,x\right\rangle  \label{I.e.2.1}
\end{equation}%
for any $x\in H$ with $\left\Vert x\right\Vert =1.$
\end{theorem}

\begin{proof}
Since $f$ is convex and differentiable, we have that 
\begin{equation*}
f\left( t\right) -f\left( s\right) \leq f^{\prime }\left( t\right) \cdot
\left( t-s\right)
\end{equation*}
for any $t,s\in \left[ m,M\right] .$

Now, if we chose in this inequality $s=\left\langle Ax,x\right\rangle \in %
\left[ m,M\right] $ for any $x\in H$ with $\left\Vert x\right\Vert =1$ since 
$Sp\left( A\right) \subseteq \left[ m,M\right] ,$ then we have 
\begin{equation}
f\left( t\right) -f\left( \left\langle Ax,x\right\rangle \right) \leq
f^{\prime }\left( t\right) \cdot \left( t-\left\langle Ax,x\right\rangle
\right)  \label{I.e.2.2}
\end{equation}%
for any $t\in \left[ m,M\right] $ any $x\in H$ with $\left\Vert x\right\Vert
=1.$

If we fix $x\in H$ with $\left\Vert x\right\Vert =1$ in (\ref{I.e.2.2}) and
apply the property (\ref{P}) then we get 
\begin{equation*}
\left\langle \left[ f\left( A\right) -f\left( \left\langle Ax,x\right\rangle
\right) 1_{H}\right] x,x\right\rangle \leq \left\langle f^{\prime }\left(
A\right) \cdot \left( A-\left\langle Ax,x\right\rangle 1_{H}\right)
x,x\right\rangle
\end{equation*}%
for each $x\in H$ with $\left\Vert x\right\Vert =1,$ which is clearly
equivalent to the desired inequality (\ref{I.e.2.1}).
\end{proof}

\begin{corollary}[Dragomir, 2008, \protect\cite{I.SSD5}]
\label{I.c.2.1}Assume that $f$ is as in the Theorem \ref{I.t.2.1}. If $A_{j}$
are selfadjoint operators with $Sp\left( A_{j}\right) \subseteq \left[ m,M%
\right] \subset $\r{I}, $j\in \left\{ 1,\dots ,n\right\} $ and $x_{j}\in
H,j\in \left\{ 1,\dots ,n\right\} $ with $\sum_{j=1}^{n}\left\Vert
x_{j}\right\Vert ^{2}=1$, then 
\begin{align}
(0& \leq )\sum_{j=1}^{n}\left\langle f\left( A_{j}\right)
x_{j},x_{j}\right\rangle -f\left( \sum_{j=1}^{n}\left\langle
A_{j}x_{j},x_{j}\right\rangle \right)  \label{I.e.2.3} \\
& \leq \sum_{j=1}^{n}\left\langle f^{\prime }\left( A_{j}\right)
A_{j}x_{j},x_{j}\right\rangle -\sum_{j=1}^{n}\left\langle
A_{j}x_{j},x_{j}\right\rangle \cdot \sum_{j=1}^{n}\left\langle f^{\prime
}\left( A_{j}\right) x_{j},x_{j}\right\rangle .  \notag
\end{align}
\end{corollary}

\begin{proof}
As in \cite[p. 6]{I.FMPS}, if we put 
\begin{equation*}
\widetilde{A}:=\left( 
\begin{array}{ccc}
A_{1} & \cdots & 0 \\ 
\vdots & \ddots & \vdots \\ 
0 & \cdots & A_{n}%
\end{array}%
\right) \qquad \text{and}\qquad \widetilde{x}=\left( 
\begin{array}{c}
x_{1} \\ 
\vdots \\ 
x_{n}%
\end{array}%
\right)
\end{equation*}

then we have $Sp\left( \widetilde{A}\right) \subseteq \left[ m,M\right] ,$ $%
\left\Vert \widetilde{x}\right\Vert =1,$%
\begin{equation*}
\left\langle f\left( \widetilde{A}\right) \widetilde{x},\widetilde{x}%
\right\rangle =\sum_{j=1}^{n}\left\langle f\left( A_{j}\right)
x_{j},x_{j}\right\rangle ,\qquad \left\langle \widetilde{A}\widetilde{x},%
\widetilde{x}\right\rangle =\sum_{j=1}^{n}\left\langle
A_{j}x_{j},x_{j}\right\rangle
\end{equation*}%
and so on$.$

Applying Theorem \ref{I.t.2.1} for $\widetilde{A}$ and $\widetilde{x}$ we
deduce the desired result (\ref{I.e.2.3}).
\end{proof}

\begin{corollary}[Dragomir, 2008, \protect\cite{I.SSD5}]
\label{I.c.2.2}Assume that $f$ is as in the Theorem \ref{I.t.2.1}. If $A_{j}$
are selfadjoint operators with $Sp\left( A_{j}\right) \subseteq \left[ m,M%
\right] \subset $\r{I}, $j\in \left\{ 1,\dots ,n\right\} $ and $p_{j}\geq 0,$
$j\in \left\{ 1,\dots ,n\right\} $ with $\sum_{j=1}^{n}p_{j}=1,$ then 
\begin{align}
(0& \leq )\left\langle \sum_{j=1}^{n}p_{j}f\left( A_{j}\right)
x,x\right\rangle -f\left( \left\langle
\sum_{j=1}^{n}p_{j}A_{j}x,x\right\rangle \right)  \label{I.e.2.4} \\
& \leq \left\langle \sum_{j=1}^{n}p_{j}f^{\prime }\left( A_{j}\right)
A_{j}x,x\right\rangle -\left\langle \sum_{j=1}^{n}p_{j}A_{j}x,x\right\rangle
\cdot \left\langle \sum_{j=1}^{n}p_{j}f^{\prime }\left( A_{j}\right)
x,x\right\rangle .  \notag
\end{align}%
for each $x\in H$ with $\left\Vert x\right\Vert =1.$
\end{corollary}

\begin{remark}
\label{I.r.2.1}The inequality (\ref{I.e.2.4}), in the scalar case, namely 
\begin{align}
(0& \leq )\sum_{j=1}^{n}p_{j}f\left( x_{j}\right) -f\left(
\sum_{j=1}^{n}p_{j}x_{j}\right)  \label{I.e.2.4.a} \\
& \leq \sum_{j=1}^{n}p_{j}f^{\prime }\left( x_{j}\right)
x_{j}-\sum_{j=1}^{n}p_{j}x_{j}\cdot \sum_{j=1}^{n}p_{j}f^{\prime }\left(
x_{j}\right) ,  \notag
\end{align}%
where $x_{j}\in $\r{I}, $j\in \left\{ 1,\dots ,n\right\} ,$ has been
obtained by the first time in 1994 by Dragomir \& Ionescu, see \cite{I.DI}.
\end{remark}

The following particular cases are of interest:

\begin{example}
\label{I.ex.2.1} \textbf{a.} Let $A$ be a positive definite operator on the
Hilbert space $H.$ Then we have the following inequality: 
\begin{equation}
(0\leq )\ln \left( \left\langle Ax,x\right\rangle \right) -\left\langle \ln
\left( A\right) x,x\right\rangle \leq \left\langle Ax,x\right\rangle \cdot
\left\langle A^{-1}x,x\right\rangle -1,  \label{I.e.2.5}
\end{equation}%
for each $x\in H$ with $\left\Vert x\right\Vert =1.$

\textbf{b.} If $A$ is a selfadjoint operator on $H$, then we have the
inequality: 
\begin{align}
(0& \leq )\left\langle \exp \left( A\right) x,x\right\rangle -\exp \left(
\left\langle Ax,x\right\rangle \right)  \label{I.e.2.6} \\
& \leq \left\langle A\exp \left( A\right) x,x\right\rangle -\left\langle
Ax,x\right\rangle \cdot \left\langle \exp \left( A\right) x,x\right\rangle ,
\notag
\end{align}%
for each $x\in H$ with $\left\Vert x\right\Vert =1.$

\textbf{c.} If $p\geq 1$ and $A$ is a positive operator on $H$, then 
\begin{equation}
(0\leq )\left\langle A^{p}x,x\right\rangle -\left\langle Ax,x\right\rangle
^{p}\leq p\left[ \left\langle A^{p}x,x\right\rangle -\left\langle
Ax,x\right\rangle \cdot \left\langle A^{p-1}x,x\right\rangle \right] ,
\label{I.e.2.7}
\end{equation}%
for each $x\in H$ with $\left\Vert x\right\Vert =1.$ If $A$ is positive
definite, then the inequality (\ref{I.e.2.7}) also holds for $p<0.$

If $0<p<1$ and $A$ is a positive definite operator then the reverse
inequality also holds 
\begin{equation}
\left\langle A^{p}x,x\right\rangle -\left\langle Ax,x\right\rangle ^{p}\geq
p \left[ \left\langle A^{p}x,x\right\rangle -\left\langle Ax,x\right\rangle
\cdot \left\langle A^{p-1}x,x\right\rangle \right] \geq 0,  \label{I.e.2.8}
\end{equation}%
for each $x\in H$ with $\left\Vert x\right\Vert =1.$
\end{example}

Similar results can be stated for sequences of operators, however the
details are omitted.

\subsection{Further Reverses}

In applications would be perhaps more useful to find upper bounds for the
quantity 
\begin{equation*}
\left\langle f\left( A\right) x,x\right\rangle -f\left( \left\langle
Ax,x\right\rangle \right) ,\text{\quad }x\in H\text{\quad with\quad }%
\left\Vert x\right\Vert =1,
\end{equation*}%
that are in terms of the spectrum margins $m,M$ and of the function $f$.

The following result may be stated:

\begin{theorem}[Dragomir, 2008, \protect\cite{I.SSD5}]
\label{I.t.3.1}Let $I$ be an interval and $f:I\rightarrow \mathbb{R}$ be a
convex and differentiable function on \r{I} (the interior of $I)$ whose
derivative $f^{\prime }$ is continuous on \r{I}$.$ If $A$ is a selfadjoint
operator on the Hilbert space $H$ with $Sp\left( A\right) \subseteq \left[
m,M\right] \subset $\r{I}$,$ then 
\begin{align}
(0& \leq )\left\langle f\left( A\right) x,x\right\rangle -f\left(
\left\langle Ax,x\right\rangle \right)  \label{I.e.3.1} \\
& \leq \left\{ 
\begin{array}{l}
\frac{1}{2}\cdot \left( M-m\right) \left[ \left\Vert f^{\prime }\left(
A\right) x\right\Vert ^{2}-\left\langle f^{\prime }\left( A\right)
x,x\right\rangle ^{2}\right] ^{1/2} \\ 
\\ 
\frac{1}{2}\cdot \left( f^{\prime }\left( M\right) -f^{\prime }\left(
m\right) \right) \left[ \left\Vert Ax\right\Vert ^{2}-\left\langle
Ax,x\right\rangle ^{2}\right] ^{1/2}%
\end{array}%
\right.  \notag \\
& \leq \frac{1}{4}\left( M-m\right) \left( f^{\prime }\left( M\right)
-f^{\prime }\left( m\right) \right) ,  \notag
\end{align}%
for any $x\in H$ with $\left\Vert x\right\Vert =1.$

We also have the inequality 
\begin{align}
(0& \leq )\left\langle f\left( A\right) x,x\right\rangle -f\left(
\left\langle Ax,x\right\rangle \right)  \label{I.e.3.1.0} \\
& \leq \frac{1}{4}\left( M-m\right) \left( f^{\prime }\left( M\right)
-f^{\prime }\left( m\right) \right)  \notag \\
& -\left\{ \hspace{-6pt}%
\begin{array}{l}
\left[ \left\langle Mx-Ax,Ax-mx\right\rangle \left\langle f^{\prime }\left(
M\right) x-f^{\prime }\left( A\right) x,f^{\prime }\left( A\right)
x-f^{\prime }\left( m\right) x\right\rangle \right] ^{\frac{1}{2}}, \\ 
\\ 
\left\vert \left\langle Ax,x\right\rangle -\frac{M+m}{2}\right\vert
\left\vert \left\langle f^{\prime }\left( A\right) x,x\right\rangle -\frac{%
f^{\prime }\left( M\right) +f^{\prime }\left( m\right) }{2}\right\vert%
\end{array}%
\right.  \notag \\
& \leq \frac{1}{4}\left( M-m\right) \left( f^{\prime }\left( M\right)
-f^{\prime }\left( m\right) \right) ,  \notag
\end{align}%
for any $x\in H$ with $\left\Vert x\right\Vert =1.$

Moreover, if $m>0$ and $f^{\prime }\left( m\right) >0,$ then we also have 
\begin{align}
(0& \leq )\left\langle f\left( A\right) x,x\right\rangle -f\left(
\left\langle Ax,x\right\rangle \right)  \label{I.e.3.1.1} \\
& \leq \left\{ 
\begin{array}{l}
\frac{1}{4}\cdot \frac{\left( M-m\right) \left( f^{\prime }\left( M\right)
-f^{\prime }\left( m\right) \right) }{\sqrt{Mmf^{\prime }\left( M\right)
f^{\prime }\left( m\right) }}\left\langle Ax,x\right\rangle \left\langle
f^{\prime }\left( A\right) x,x\right\rangle , \\[10pt] 
\left( \sqrt{M}-\sqrt{m}\right) \left( \sqrt{f^{\prime }\left( M\right) }-%
\sqrt{f^{\prime }\left( m\right) }\right) \left[ \left\langle
Ax,x\right\rangle \left\langle f^{\prime }\left( A\right) x,x\right\rangle %
\right] ^{\frac{1}{2}},%
\end{array}%
\right.  \notag
\end{align}%
for any $x\in H$ with $\left\Vert x\right\Vert =1.$
\end{theorem}

\begin{proof}
We use the following Gr\"{u}ss' type result we obtained in \cite{I.SSD1}:

Let $A$ be a selfadjoint operator on the Hilbert space $\left(
H;\left\langle .,.\right\rangle \right) $ and assume that $Sp\left( A\right)
\subseteq \left[ m,M\right] $ for some scalars $m<M.$ If $h$and $g$ are
continuous on $\left[ m,M\right] $ and $\gamma :=\min_{t\in \left[ m,M\right]
}h\left( t\right) $ and $\Gamma :=\max_{t\in \left[ m,M\right] }h\left(
t\right) ,$ then 
\begin{align}
\lefteqn{\left\vert \left\langle h\left( A\right) g\left( A\right)
x,x\right\rangle -\left\langle h\left( A\right) x,x\right\rangle \cdot
\left\langle g\left( A\right) x,x\right\rangle \right\vert }  \label{I.e.3.2}
\\
& \leq \frac{1}{2}\cdot \left( \Gamma -\gamma \right) \left[ \left\Vert
g\left( A\right) x\right\Vert ^{2}-\left\langle g\left( A\right)
x,x\right\rangle ^{2}\right] ^{1/2}  \notag \\
& \left( \leq \frac{1}{4}\left( \Gamma -\gamma \right) \left( \Delta -\delta
\right) \right)  \notag
\end{align}%
for each $x\in H$ with $\left\Vert x\right\Vert =1,$ where $\delta
:=\min_{t\in \left[ m,M\right] }g\left( t\right) $ and $\Delta :=\max_{t\in %
\left[ m,M\right] }g\left( t\right) .$

Therefore, we can state that 
\begin{align}
\lefteqn{\left\langle Af^{\prime }\left( A\right) x,x\right\rangle
-\left\langle Ax,x\right\rangle \cdot \left\langle f^{\prime }\left(
A\right) x,x\right\rangle }  \label{I.e.3.3} \\
& \leq \frac{1}{2}\cdot \left( M-m\right) \left[ \left\Vert f^{\prime
}\left( A\right) x\right\Vert ^{2}-\left\langle f^{\prime }\left( A\right)
x,x\right\rangle ^{2}\right] ^{1/2}  \notag \\
& \leq \frac{1}{4}\left( M-m\right) \left( f^{\prime }\left( M\right)
-f^{\prime }\left( m\right) \right)  \notag
\end{align}%
and 
\begin{align}
\lefteqn{\left\langle Af^{\prime }\left( A\right) x,x\right\rangle
-\left\langle Ax,x\right\rangle \cdot \left\langle f^{\prime }\left(
A\right) x,x\right\rangle }  \label{I.e.3.4} \\
& \leq \frac{1}{2}\cdot \left( f^{\prime }\left( M\right) -f^{\prime }\left(
m\right) \right) \left[ \left\Vert Ax\right\Vert ^{2}-\left\langle
Ax,x\right\rangle ^{2}\right] ^{1/2}  \notag \\
& \leq \frac{1}{4}\left( M-m\right) \left( f^{\prime }\left( M\right)
-f^{\prime }\left( m\right) \right)  \notag
\end{align}%
for each $x\in H$ with $\left\Vert x\right\Vert =1,$ which together with (%
\ref{I.e.2.1}) provide the desired result (\ref{I.e.3.1}).

On making use of the inequality obtained in \cite{I.SSD2}: 
\begin{align}
& \left\vert \left\langle h\left( A\right) g\left( A\right) x,x\right\rangle
-\left\langle h\left( A\right) x,x\right\rangle \left\langle g\left(
A\right) x,x\right\rangle \right\vert  \label{I.e.3.4.1} \\
& \leq \frac{1}{4}\cdot \left( \Gamma -\gamma \right) \left( \Delta -\delta
\right)  \notag \\
& -\left\{ 
\begin{array}{l}
\left[ \left\langle \Gamma x-h\left( A\right) x,f\left( A\right) x-\gamma
x\right\rangle \left\langle \Delta x-g\left( A\right) x,g\left( A\right)
x-\delta x\right\rangle \right] ^{\frac{1}{2}}, \\ 
\\ 
\left\vert \left\langle h\left( A\right) x,x\right\rangle -\frac{\Gamma
+\gamma }{2}\right\vert \left\vert \left\langle g\left( A\right)
x,x\right\rangle -\frac{\Delta +\delta }{2}\right\vert ,%
\end{array}%
\right.  \notag
\end{align}%
for each $x\in H$ with $\left\Vert x\right\Vert =1,$ we can state that 
\begin{align*}
& \left\langle Af^{\prime }\left( A\right) x,x\right\rangle -\left\langle
Ax,x\right\rangle \cdot \left\langle f^{\prime }\left( A\right)
x,x\right\rangle \\
& \leq \frac{1}{4}\left( M-m\right) \left( f^{\prime }\left( M\right)
-f^{\prime }\left( m\right) \right) \\
& -\left\{ 
\begin{array}{l}
\left[ \left\langle Mx-Ax,Ax-mx\right\rangle \left\langle f^{\prime }\left(
M\right) x-f^{\prime }\left( A\right) x,f^{\prime }\left( A\right)
x-f^{\prime }\left( m\right) x\right\rangle \right] ^{\frac{1}{2}}, \\ 
\\ 
\left\vert \left\langle Ax,x\right\rangle -\frac{M+m}{2}\right\vert
\left\vert \left\langle f^{\prime }\left( A\right) x,x\right\rangle -\frac{%
f^{\prime }\left( M\right) +f^{\prime }\left( m\right) }{2}\right\vert .%
\end{array}%
\right.
\end{align*}%
for each $x\in H$ with $\left\Vert x\right\Vert =1,$ which together with (%
\ref{I.e.2.1}) provide the desired result (\ref{I.e.3.1.0}).

Further, in order to prove the third inequality, we make use of the
following result of Gr\"{u}ss type obtained in \cite{I.SSD2}:

If $\gamma $ and $\delta $ are positive, then 
\begin{align}
& \left\vert \left\langle h\left( A\right) g\left( A\right) x,x\right\rangle
-\left\langle h\left( A\right) x,x\right\rangle \left\langle g\left(
A\right) x,x\right\rangle \right\vert  \label{I.e.3.5} \\
& \leq \left\{ 
\begin{array}{l}
\frac{1}{4}\cdot \frac{\left( \Gamma -\gamma \right) \left( \Delta -\delta
\right) }{\sqrt{\Gamma \gamma \Delta \delta }}\left\langle h\left( A\right)
x,x\right\rangle \left\langle g\left( A\right) x,x\right\rangle , \\ 
\\ 
\left( \sqrt{\Gamma }-\sqrt{\gamma }\right) \left( \sqrt{\Delta }-\sqrt{%
\delta }\right) \left[ \left\langle h\left( A\right) x,x\right\rangle
\left\langle g\left( A\right) x,x\right\rangle \right] ^{\frac{1}{2}}.%
\end{array}%
\right.  \notag
\end{align}%
for each $x\in H$ with $\left\Vert x\right\Vert =1.$

Now, on making use of (\ref{I.e.3.5}) we can state that 
\begin{align*}
& \left\langle Af^{\prime }\left( A\right) x,x\right\rangle -\left\langle
Ax,x\right\rangle \cdot \left\langle f^{\prime }\left( A\right)
x,x\right\rangle \\
& \leq \left\{ 
\begin{array}{l}
\frac{1}{4}\cdot \frac{\left( M-m\right) \left( f^{\prime }\left( M\right)
-f^{\prime }\left( m\right) \right) }{\sqrt{Mmf^{\prime }\left( M\right)
f^{\prime }\left( m\right) }}\left\langle Ax,x\right\rangle \left\langle
f^{\prime }\left( A\right) x,x\right\rangle , \\ 
\\ 
\left( \sqrt{M}-\sqrt{m}\right) \left( \sqrt{f^{\prime }\left( M\right) }-%
\sqrt{f^{\prime }\left( m\right) }\right) \left[ \left\langle
Ax,x\right\rangle \left\langle f^{\prime }\left( A\right) x,x\right\rangle %
\right] ^{\frac{1}{2}}.%
\end{array}%
\right.
\end{align*}%
for each $x\in H$ with $\left\Vert x\right\Vert =1,$ which together with (%
\ref{I.e.2.1}) provide the desired result (\ref{I.e.3.1.1}).
\end{proof}

\begin{corollary}[Dragomir, 2008, \protect\cite{I.SSD5}]
\label{I.c.3.1}Assume that $f$ is as in the Theorem \ref{I.t.3.1}. If $A_{j}$
are selfadjoint operators with $Sp\left( A_{j}\right) \subseteq \left[ m,M%
\right] \subset $\r{I}, $j\in \left\{ 1,\dots ,n\right\} $, then 
\begin{align}
(0& \leq )\sum_{j=1}^{n}\left\langle f\left( A_{j}\right)
x_{j},x_{j}\right\rangle -f\left( \sum_{j=1}^{n}\left\langle
A_{j}x_{j},x_{j}\right\rangle \right)  \label{I.e.3.6} \\
& \leq \left\{ \hspace{-6pt}%
\begin{array}{c}
\frac{1}{2}\cdot \left( M-m\right) \left[ \sum_{j=1}^{n}\left\Vert f^{\prime
}\left( A_{j}\right) x_{j}\right\Vert ^{2}-\left( \sum_{j=1}^{n}\left\langle
f^{\prime }\left( A_{j}\right) x_{j},x_{j}\right\rangle \right) ^{2}\right]
^{1/2}, \\ 
\\ 
\frac{1}{2}\cdot \left( f^{\prime }\left( M\right) -f^{\prime }\left(
m\right) \right) \left[ \sum_{j=1}^{n}\left\Vert A_{j}x_{j}\right\Vert
^{2}-\left( \sum_{j=1}^{n}\left\langle A_{j}x_{j},x_{j}\right\rangle \right)
^{2}\right] ^{1/2},%
\end{array}%
\right.  \notag \\
& \leq \frac{1}{4}\left( M-m\right) \left( f^{\prime }\left( M\right)
-f^{\prime }\left( m\right) \right) ,  \notag
\end{align}%
for any $x_{j}\in H,j\in \left\{ 1,\dots ,n\right\} $ with $%
\sum_{j=1}^{n}\left\Vert x_{j}\right\Vert ^{2}=1.$

We also have the inequality 
\begin{align}
(0& \leq )\sum_{j=1}^{n}\left\langle f\left( A_{j}\right)
x_{j},x_{j}\right\rangle -f\left( \sum_{j=1}^{n}\left\langle
A_{j}x_{j},x_{j}\right\rangle \right)  \label{I.e.3.7} \\
& \leq \frac{1}{4}\left( M-m\right) \left( f^{\prime }\left( M\right)
-f^{\prime }\left( m\right) \right)  \notag \\
& -\left\{ 
\begin{array}{l}
\left[ \sum\limits_{j=1}^{n}\left\langle
Mx_{j}-A_{j}x,A_{j}x_{j}-mx_{j}\right\rangle \right] ^{\frac{1}{2}} \\ 
\times \left[ \sum\limits_{j=1}^{n}\left\langle f^{\prime }\left( M\right)
x_{j}-f^{\prime }\left( A_{j}\right) x_{j},f^{\prime }\left( A_{j}\right)
x_{j}-f^{\prime }\left( m\right) x_{j}\right\rangle \right] ^{1/2}, \\ 
\\ 
\left\vert \sum\limits_{j=1}^{n}\left\langle A_{j}x_{j},x_{j}\right\rangle -%
\frac{M+m}{2}\right\vert \left\vert \sum\limits_{j=1}^{n}\left\langle
f^{\prime }\left( A_{j}\right) x_{j},x_{j}\right\rangle -\frac{f^{\prime
}\left( M\right) +f^{\prime }\left( m\right) }{2}\right\vert%
\end{array}%
\right.  \notag \\
& \leq \frac{1}{4}\left( M-m\right) \left( f^{\prime }\left( M\right)
-f^{\prime }\left( m\right) \right) ,  \notag
\end{align}%
for any $x_{j}\in H,j\in \left\{ 1,\dots ,n\right\} $ with $%
\sum_{j=1}^{n}\left\Vert x_{j}\right\Vert ^{2}=1.$

Moreover, if $m>0$ and $f^{\prime }\left( m\right) >0,$ then we also have 
\begin{align}
(0& \leq )\sum_{j=1}^{n}\left\langle f\left( A_{j}\right)
x_{j},x_{j}\right\rangle -f\left( \sum_{j=1}^{n}\left\langle
A_{j}x_{j},x_{j}\right\rangle \right)  \label{I.e.3.8} \\
& \leq \left\{ 
\begin{array}{l}
\frac{1}{4}\cdot \frac{\left( M-m\right) \left( f^{\prime }\left( M\right)
-f^{\prime }\left( m\right) \right) }{\sqrt{Mmf^{\prime }\left( M\right)
f^{\prime }\left( m\right) }}\sum_{j=1}^{n}\left\langle
A_{j}x_{j},x_{j}\right\rangle \sum_{j=1}^{n}\left\langle f^{\prime }\left(
A_{j}\right) x_{j},x_{j}\right\rangle , \\ 
\\ 
\left( \sqrt{M}-\sqrt{m}\right) \left( \sqrt{f^{\prime }\left( M\right) }-%
\sqrt{f^{\prime }\left( m\right) }\right) \\ 
\times \left[ \sum_{j=1}^{n}\left\langle A_{j}x_{j},x_{j}\right\rangle
\sum_{j=1}^{n}\left\langle f^{\prime }\left( A_{j}\right)
x_{j},x_{j}\right\rangle \right] ^{\frac{1}{2}},%
\end{array}%
\right.  \notag
\end{align}%
for any $x_{j}\in H,j\in \left\{ 1,\dots ,n\right\} $ with $%
\sum_{j=1}^{n}\left\Vert x_{j}\right\Vert ^{2}=1.$
\end{corollary}

The following corollary also holds:

\begin{corollary}[Dragomir, 2008, \protect\cite{I.SSD5}]
\label{I.c.3.2}Assume that $f$ is as in the Theorem \ref{I.t.2.1}. If $A_{j}$
are selfadjoint operators with $Sp\left( A_{j}\right) \subseteq \left[ m,M%
\right] \subset $\r{I}, $j\in \left\{ 1,\dots ,n\right\} $ and $p_{j}\geq 0,$
$j\in \left\{ 1,\dots ,n\right\} $ with $\sum_{j=1}^{n}p_{j}=1,$ then 
\begin{align}
(0& \leq )\left\langle \sum_{j=1}^{n}p_{j}f\left( A_{j}\right)
x,x\right\rangle -f\left( \left\langle
\sum_{j=1}^{n}p_{j}A_{j}x,x\right\rangle \right)  \label{I.e.3.9} \\
& \leq \left\{ 
\begin{array}{l}
\frac{1}{2}\cdot \left( M-m\right) \left[ \sum\limits_{j=1}^{n}p_{j}\left%
\Vert f^{\prime }\left( A_{j}\right) x\right\Vert ^{2}-\left\langle
\sum\limits_{j=1}^{n}p_{j}f^{\prime }\left( A_{j}\right) x,x\right\rangle
^{2}\right] ^{1/2}, \\ 
\\ 
\frac{1}{2}\cdot \left( f^{\prime }\left( M\right) -f^{\prime }\left(
m\right) \right) \left[ \sum\limits_{j=1}^{n}p_{j}\left\Vert
A_{j}x\right\Vert ^{2}-\left\langle
\sum\limits_{j=1}^{n}p_{j}A_{j}x,x\right\rangle ^{2}\right] ^{1/2},%
\end{array}%
\right.  \notag \\
& \leq \frac{1}{4}\left( M-m\right) \left( f^{\prime }\left( M\right)
-f^{\prime }\left( m\right) \right) ,  \notag
\end{align}%
for any $x\in H$ with $\left\Vert x\right\Vert =1.$

We also have the inequality 
\begin{align}
(0& \leq )\left\langle \sum_{j=1}^{n}p_{j}f\left( A_{j}\right)
x,x\right\rangle -f\left( \left\langle
\sum_{j=1}^{n}p_{j}A_{j}x,x\right\rangle \right)  \label{I.e.3.10} \\
& \leq \frac{1}{4}\left( M-m\right) \left( f^{\prime }\left( M\right)
-f^{\prime }\left( m\right) \right)  \notag \\
& -\left\{ 
\begin{array}{l}
\left[ \sum\limits_{j=1}^{n}p_{j}\left\langle
Mx-A_{j}x,A_{j}x-mx\right\rangle \right] ^{\frac{1}{2}} \\ 
\times \left[ \sum\limits_{j=1}^{n}p_{j}\left\langle f^{\prime }\left(
M\right) x-f^{\prime }\left( A_{j}\right) x,f^{\prime }\left( A_{j}\right)
x-f^{\prime }\left( m\right) x\right\rangle \right] ^{1/2}, \\ 
\\ 
\left\vert \left\langle \sum\limits_{j=1}^{n}p_{j}A_{j}x,x\right\rangle -%
\frac{M+m}{2}\right\vert \left\vert \left\langle
\sum\limits_{j=1}^{n}p_{j}f^{\prime }\left( A_{j}\right) x,x\right\rangle -%
\frac{f^{\prime }\left( M\right) +f^{\prime }\left( m\right) }{2}\right\vert%
\end{array}%
\right.  \notag \\
& \leq \frac{1}{4}\left( M-m\right) \left( f^{\prime }\left( M\right)
-f^{\prime }\left( m\right) \right) ,  \notag
\end{align}%
for any $x\in H$ with $\left\Vert x\right\Vert =1.$

Moreover, if $m>0$ and $f^{\prime }\left( m\right) >0,$ then we also have 
\begin{align}
(0& \leq )\left\langle \sum_{j=1}^{n}p_{j}f\left( A_{j}\right)
x,x\right\rangle -f\left( \left\langle
\sum_{j=1}^{n}p_{j}A_{j}x,x\right\rangle \right)  \label{I.e.3.11} \\
& \leq \left\{ 
\begin{array}{l}
\frac{1}{4}\cdot \frac{\left( M-m\right) \left( f^{\prime }\left( M\right)
-f^{\prime }\left( m\right) \right) }{\sqrt{Mmf^{\prime }\left( M\right)
f^{\prime }\left( m\right) }}\left\langle
\sum_{j=1}^{n}p_{j}A_{j}x,x\right\rangle \left\langle
\sum_{j=1}^{n}p_{j}f^{\prime }\left( A_{j}\right) x,x\right\rangle , \\ 
\\ 
\left( \sqrt{M}-\sqrt{m}\right) \left( \sqrt{f^{\prime }\left( M\right) }-%
\sqrt{f^{\prime }\left( m\right) }\right) \\ 
\times \left[ \left\langle \sum_{j=1}^{n}p_{j}A_{j}x,x\right\rangle
\left\langle \sum_{j=1}^{n}p_{j}f^{\prime }\left( A_{j}\right)
x,x\right\rangle \right] ^{\frac{1}{2}},%
\end{array}%
\right.  \notag
\end{align}%
for any $x\in H$ with $\left\Vert x\right\Vert =1.$
\end{corollary}

\begin{remark}
\label{I.r.3.1}Some of the inequalities in Corollary \ref{I.c.3.2} can be
used to produce reverse norm inequalities for the sum of positive operators
in the case when the convex function $f$ is nonnegative and monotonic
nondecreasing on $\left[ 0,M\right] .$

For instance, if we use the inequality (\ref{I.e.3.9}), then we have 
\begin{equation}
(0\leq )\left\Vert \sum_{j=1}^{n}p_{j}f\left( A_{j}\right) \right\Vert
-f\left( \left\Vert \sum_{j=1}^{n}p_{j}A_{j}\right\Vert \right) \leq \frac{1%
}{4}\left( M-m\right) \left( f^{\prime }\left( M\right) -f^{\prime }\left(
m\right) \right) .  \label{I.e.3.12}
\end{equation}%
Moreover, if we use the inequality (\ref{I.e.3.11}), then we obtain 
\begin{align}
(0& \leq )\left\Vert \sum_{j=1}^{n}p_{j}f\left( A_{j}\right) \right\Vert
-f\left( \left\Vert \sum_{j=1}^{n}p_{j}A_{j}\right\Vert \right)
\label{I.e.3.13} \\
& \leq \left\{ \hspace{-7pt}%
\begin{array}{l}
\frac{1}{4}\cdot \frac{\left( M-m\right) \left( f^{\prime }\left( M\right)
-f^{\prime }\left( m\right) \right) }{\sqrt{Mmf^{\prime }\left( M\right)
f^{\prime }\left( m\right) }}\left\Vert
\sum\limits_{j=1}^{n}p_{j}A_{j}\right\Vert \left\Vert
\sum\limits_{j=1}^{n}p_{j}f^{\prime }\left( A_{j}\right) \right\Vert , \\%
[10pt] 
\left( \sqrt{M}\hspace{-2pt}-\hspace{-2pt}\sqrt{m}\right) \hspace{-4pt}%
\left( \sqrt{f^{\prime }\left( M\right) }\hspace{-2pt}-\hspace{-2pt}\sqrt{%
f^{\prime }\left( m\right) }\right) \hspace{-4pt}\left[ \left\Vert
\sum\limits_{j=1}^{n}p_{j}A_{j}\right\Vert \left\Vert
\sum\limits_{j=1}^{n}p_{j}f^{\prime }\left( A_{j}\right) \right\Vert \right]
^{\frac{1}{2}}.%
\end{array}%
\right.  \notag
\end{align}
\end{remark}

\subsection{Some Particular Inequalities of Interest}

\textbf{1.} Consider the convex function $f:\left( 0,\infty \right)
\rightarrow \mathbb{R}$, $f\left( x\right) =-\ln x.$ On utilising the
inequality (\ref{I.e.3.1}), then for any positive definite operator $A$ on
the Hilbert space $H,$ we have the inequality 
\begin{align}
(0& \leq )\ln \left( \left\langle Ax,x\right\rangle \right) -\left\langle
\ln \left( A\right) x,x\right\rangle  \label{I.e.4.1} \\
& \leq \left\{ 
\begin{array}{l}
\frac{1}{2}\cdot \left( M-m\right) \left[ \left\Vert A^{-1}x\right\Vert
^{2}-\left\langle A^{-1}x,x\right\rangle ^{2}\right] ^{1/2} \\ 
\\ 
\frac{1}{2}\cdot \frac{M-m}{mM}\left[ \left\Vert Ax\right\Vert
^{2}-\left\langle Ax,x\right\rangle ^{2}\right] ^{1/2}%
\end{array}%
\right.  \notag \\
& \left( \leq \frac{1}{4}\cdot \frac{\left( M-m\right) ^{2}}{mM}\right) 
\notag
\end{align}%
for any $x\in H$ with $\left\Vert x\right\Vert =1.$

However, if we use the inequality (\ref{I.e.3.1.0}), then we have the
following result as well 
\begin{align}
(0& \leq )\ln \left( \left\langle Ax,x\right\rangle \right) -\left\langle
\ln \left( A\right) x,x\right\rangle  \label{I.e.4.2} \\
& \leq \frac{1}{4}\cdot \frac{\left( M-m\right) ^{2}}{mM}  \notag \\
& -\left\{ 
\begin{array}{l}
\left[ \left\langle Mx-Ax,Ax-mx\right\rangle \left\langle
M^{-1}x-A^{-1}x,A^{-1}x-m^{-1}x\right\rangle \right] ^{\frac{1}{2}}, \\ 
\\ 
\left\vert \left\langle Ax,x\right\rangle -\frac{M+m}{2}\right\vert
\left\vert \left\langle A^{-1}x,x\right\rangle -\frac{M+m}{2mM}\right\vert%
\end{array}%
\right.  \notag \\
& \left( \leq \frac{1}{4}\cdot \frac{\left( M-m\right) ^{2}}{mM}\right) 
\notag
\end{align}%
for any $x\in H$ with $\left\Vert x\right\Vert =1.$

\textbf{2. }Now consider the convex function $f:\left( 0,\infty \right)
\rightarrow \mathbb{R}$, $f\left( x\right) =x\ln x.$ On utilising the
inequality (\ref{I.e.3.1}), then for any positive definite operator $A$ on
the Hilbert space $H,$ we have the inequality 
\begin{align}
(0& \leq )\left\langle A\ln \left( A\right) x,x\right\rangle -\left\langle
Ax,x\right\rangle \ln \left( \left\langle Ax,x\right\rangle \right)
\label{I.e.4.3} \\
& \leq \left\{ 
\begin{array}{l}
\frac{1}{2}\cdot \left( M-m\right) \left[ \left\Vert \ln \left( eA\right)
x\right\Vert ^{2}-\left\langle \ln \left( eA\right) x,x\right\rangle ^{2}%
\right] ^{1/2} \\ 
\\ 
\ln \sqrt{\frac{M}{m}}\cdot \left[ \left\Vert Ax\right\Vert
^{2}-\left\langle Ax,x\right\rangle ^{2}\right] ^{1/2}%
\end{array}%
\right.  \notag \\
& \left( \leq \frac{1}{2}\left( M-m\right) \ln \sqrt{\frac{M}{m}}\right) 
\notag
\end{align}%
for any $x\in H$ with $\left\Vert x\right\Vert =1.$

If we now apply the inequality (\ref{I.e.3.1.0}), then we have the following
result as well 
\begin{align}
(0& \leq )\left\langle A\ln \left( A\right) x,x\right\rangle -\left\langle
Ax,x\right\rangle \ln \left( \left\langle Ax,x\right\rangle \right)
\label{I.e.4.4} \\
& \leq \frac{1}{2}\left( M-m\right) \ln \sqrt{\frac{M}{m}}  \notag \\
& -\left\{ \hspace{-4pt}%
\begin{array}{l}
\left[ \left\langle Mx\hspace{-2pt}-\hspace{-2pt}Ax,Ax\hspace{-2pt}-\hspace{%
-2pt}mx\right\rangle \left\langle \ln \left( M\right) x\hspace{-2pt}-\hspace{%
-2pt}\ln \left( A\right) x,\ln \left( A\right) x\hspace{-2pt}-\hspace{-2pt}%
\ln \left( m\right) x\right\rangle \right] ^{\frac{1}{2}}, \\ 
\\ 
\left\vert \left\langle Ax,x\right\rangle -\frac{M+m}{2}\right\vert
\left\vert \left\langle \ln \left( A\right) x,x\right\rangle -\ln \sqrt{mM}%
\right\vert%
\end{array}%
\right.  \notag \\
& \left( \leq \frac{1}{2}\left( M-m\right) \ln \sqrt{\frac{M}{m}}\right) 
\notag
\end{align}%
for any $x\in H$ with $\left\Vert x\right\Vert =1.$

Moreover, if we assume that $m>e^{-1},$ then, by utilising the inequality (%
\ref{I.e.3.1.1}) we can state the inequality 
\begin{align}
(0& \leq )\left\langle A\ln \left( A\right) x,x\right\rangle -\left\langle
Ax,x\right\rangle \ln \left( \left\langle Ax,x\right\rangle \right)
\label{I.e.4.5} \\
& \leq \left\{ 
\begin{array}{l}
\frac{1}{2}\cdot \frac{\left( M-m\right) \ln \sqrt{\frac{M}{m}}}{\sqrt{Mm\ln
\left( eM\right) \ln \left( em\right) }}\left\langle Ax,x\right\rangle
\left\langle \ln \left( eA\right) x,x\right\rangle , \\ 
\\ 
\left( \sqrt{M}-\sqrt{m}\right) \left( \sqrt{\ln \left( eM\right) }-\sqrt{%
\ln \left( em\right) }\right) \left[ \left\langle Ax,x\right\rangle
\left\langle \ln \left( eA\right) x,x\right\rangle \right] ^{\frac{1}{2}},%
\end{array}%
\right.  \notag
\end{align}%
for any $x\in H$ with $\left\Vert x\right\Vert =1.$

\textbf{3. }Consider now the following convex function $f:\mathbb{R}%
\rightarrow \left( 0,\infty \right) ,$ $f\left( x\right) =\exp \left( \alpha
x\right) $ with $\alpha >0.$ If we apply the inequalities (\ref{I.e.3.1}), (%
\ref{I.e.3.1.0}) and (\ref{I.e.3.1.1}) for $f\left( x\right) =\exp \left(
\alpha x\right) $ and for a selfadjoint operator $A,$ then we get the
following results 
\begin{align}
(0& \leq )\left\langle \exp \left( \alpha A\right) x,x\right\rangle -\exp
\left( \alpha \left\langle Ax,x\right\rangle \right)  \label{I.e.4.6} \\
& \leq \left\{ 
\begin{array}{l}
\frac{1}{2}\cdot \alpha \left( M-m\right) \left[ \left\Vert \exp \left(
\alpha A\right) x\right\Vert ^{2}-\left\langle \exp \left( \alpha A\right)
x,x\right\rangle ^{2}\right] ^{1/2} \\ 
\\ 
\frac{1}{2}\cdot \alpha \left( \exp \left( \alpha M\right) -\exp \left(
\alpha m\right) \right) \left[ \left\Vert Ax\right\Vert ^{2}-\left\langle
Ax,x\right\rangle ^{2}\right] ^{1/2}%
\end{array}%
\right.  \notag \\
& \left( \leq \frac{1}{4}\alpha \left( M-m\right) \left( \exp \left( \alpha
M\right) -\exp \left( \alpha m\right) \right) \right) ,  \notag
\end{align}%
and 
\begin{align}
(0& \leq )\left\langle \exp \left( \alpha A\right) x,x\right\rangle -\exp
\left( \alpha \left\langle Ax,x\right\rangle \right)  \label{I.e.4.7} \\
& \leq \frac{1}{4}\alpha \left( M-m\right) \left( \exp \left( \alpha
M\right) -\exp \left( \alpha m\right) \right)  \notag \\
& -\alpha \left\{ \hspace{-6pt}%
\begin{array}{l}
\left[ \left\langle Mx-Ax,Ax-mx\right\rangle \right] ^{1/2} \\ 
\quad \times \left[ \left\langle \exp \left( \alpha M\right) x-\exp \left(
\alpha A\right) x,\exp \left( \alpha A\right) x-\exp \left( \alpha m\right)
x\right\rangle \right] ^{\frac{1}{2}}, \\ 
\\ 
\left\vert \left\langle Ax,x\right\rangle -\frac{M+m}{2}\right\vert
\left\vert \left\langle \exp \left( \alpha A\right) x,x\right\rangle -\frac{%
\exp \left( \alpha M\right) +\exp \left( \alpha m\right) }{2}\right\vert%
\end{array}%
\right.  \notag \\
& \left( \leq \frac{1}{4}\alpha \left( M-m\right) \left( \exp \left( \alpha
M\right) -\exp \left( \alpha m\right) \right) \right)  \notag
\end{align}%
and 
\begin{align}
(0& \leq )\left\langle \exp \left( \alpha A\right) x,x\right\rangle -\exp
\left( \alpha \left\langle Ax,x\right\rangle \right)  \label{I.e.4.8} \\
& \leq \alpha \times \left\{ 
\begin{array}{l}
\frac{1}{4}\cdot \frac{\left( M-m\right) \left( \exp \left( \alpha M\right)
-\exp \left( \alpha m\right) \right) }{\sqrt{Mm}\exp \left[ \frac{\alpha
\left( M+m\right) }{2}\right] }\left\langle Ax,x\right\rangle \left\langle
\exp \left( \alpha A\right) x,x\right\rangle , \\ 
\\ 
\left( \sqrt{M}-\sqrt{m}\right) \left( \exp \left( \frac{\alpha M}{2}\right)
-\exp \left( \frac{\alpha m}{2}\right) \right) \\ 
\qquad \times \left[ \left\langle Ax,x\right\rangle \left\langle \exp \left(
\alpha A\right) x,x\right\rangle \right] ^{\frac{1}{2}}%
\end{array}%
\right.  \notag
\end{align}%
for any $x\in H$ with $\left\Vert x\right\Vert =1,$ respectively.

Now, consider the convex function $f:\mathbb{R}\rightarrow \left( 0,\infty
\right) ,$ $f\left( x\right) =\exp \left( -\beta x\right) $ with $\beta >0$.
If we apply the inequalities (\ref{I.e.3.1}) and (\ref{I.e.3.1.0}) for $%
f\left( x\right) =\exp \left( -\beta x\right) $ and for a selfadjoint
operator $A,$ then we get the following results 
\begin{align}
(0& \leq )\left\langle \exp \left( -\beta A\right) x,x\right\rangle -\exp
\left( -\beta \left\langle Ax,x\right\rangle \right)  \label{I.e.4.9} \\
& \leq \beta \times \left\{ 
\begin{array}{l}
\frac{1}{2}\cdot \left( M-m\right) \left[ \left\Vert \exp \left( -\beta
A\right) x\right\Vert ^{2}-\left\langle \exp \left( -\beta A\right)
x,x\right\rangle ^{2}\right] ^{1/2} \\ 
\\ 
\frac{1}{2}\cdot \left( \exp \left( -\beta m\right) -\exp \left( -\beta
M\right) \right) \left[ \left\Vert Ax\right\Vert ^{2}-\left\langle
Ax,x\right\rangle ^{2}\right] ^{1/2}%
\end{array}%
\right.  \notag \\
& \left( \leq \frac{1}{4}\beta \left( M-m\right) \left( \exp \left( -\beta
m\right) -\exp \left( -\beta M\right) \right) \right)  \notag
\end{align}%
and 
\begin{align}
(0& \leq )\left\langle \exp \left( -\beta A\right) x,x\right\rangle -\exp
\left( -\beta \left\langle Ax,x\right\rangle \right)  \label{I.e.4.10} \\
& \leq \frac{1}{4}\beta \left( M-m\right) \left( \exp \left( -\beta m\right)
-\exp \left( -\beta M\right) \right)  \notag \\
& -\beta \left\{ \hspace{-7pt}%
\begin{array}{l}
\left[ \left\langle Mx-Ax,Ax-mx\right\rangle \right] ^{1/2} \\ 
\times \left[ \left\langle \exp \left( -\beta M\right) x-\exp \left( -\beta
A\right) x,\exp \left( -\beta A\right) x-\exp \left( -\beta m\right)
x\right\rangle \right] ^{\frac{1}{2}}, \\ 
\\ 
\left\vert \left\langle Ax,x\right\rangle -\frac{M+m}{2}\right\vert
\left\vert \left\langle \exp \left( -\beta A\right) x,x\right\rangle -\frac{%
\exp \left( -\beta M\right) +\exp \left( -\beta m\right) }{2}\right\vert%
\end{array}%
\right.  \notag \\
& \left( \leq \frac{1}{4}\beta \left( M-m\right) \left( \exp \left( -\beta
m\right) -\exp \left( -\beta M\right) \right) \right)  \notag
\end{align}%
for any $x\in H$ with $\left\Vert x\right\Vert =1,$ respectively.

\textbf{4. }Finally, if we consider the convex function $f:\left[ 0,\infty
\right) \rightarrow \left[ 0,\infty \right) ,$ $f\left( x\right) =x^{p}$
with $p\geq 1,$ then on applying the inequalities (\ref{I.e.3.1}) and (\ref%
{I.e.3.1.0}) for the positive operator $A$ we have the inequalities 
\begin{align}
(0& \leq )\left\langle A^{p}x,x\right\rangle -\left\langle Ax,x\right\rangle
^{p}  \label{I.e.4.11} \\
& \leq p\times \left\{ 
\begin{array}{l}
\frac{1}{2}\cdot \left( M-m\right) \left[ \left\Vert A^{p-1}x\right\Vert
^{2}-\left\langle A^{p-1}x,x\right\rangle ^{2}\right] ^{1/2} \\ 
\\ 
\frac{1}{2}\cdot \left( M^{p-1}-m^{p-1}\right) \left[ \left\Vert
Ax\right\Vert ^{2}-\left\langle Ax,x\right\rangle ^{2}\right] ^{1/2}%
\end{array}%
\right.  \notag \\
& \left( \leq \frac{1}{4}p\left( M-m\right) \left( M^{p-1}-m^{p-1}\right)
\right)  \notag
\end{align}%
and 
\begin{align}
(0& \leq )\left\langle A^{p}x,x\right\rangle -\left\langle Ax,x\right\rangle
^{p}  \label{I.e.4.12} \\
& \leq \frac{1}{4}p\left( M-m\right) \left( M^{p-1}-m^{p-1}\right)  \notag \\
& \quad -p\left\{ \hspace{-6pt}%
\begin{array}{l}
\left[ \left\langle Mx-Ax,Ax-mx\right\rangle \left\langle
M^{p-1}x-A^{p-1}x,A^{p-1}x-m^{p-1}x\right\rangle \right] ^{\frac{1}{2}}, \\ 
\\ 
\left\vert \left\langle Ax,x\right\rangle -\frac{M+m}{2}\right\vert
\left\vert \left\langle A^{p-1}x,x\right\rangle -\frac{M^{p-1}+m^{p-1}}{2}%
\right\vert%
\end{array}%
\right.  \notag \\
& \left( \leq \frac{1}{4}p\left( M-m\right) \left( M^{p-1}-m^{p-1}\right)
\right)  \notag
\end{align}%
for any $x\in H$ with $\left\Vert x\right\Vert =1,$ respectively.

If the operator $A$ is positive definite $\left( m>0\right) $ then, by
utilising the inequality (\ref{I.e.3.1.1}), we have 
\begin{align}
(0& \leq )\left\langle A^{p}x,x\right\rangle -\left\langle Ax,x\right\rangle
^{p}  \label{I.e.4.13} \\
& \leq p\times \left\{ 
\begin{array}{l}
\frac{1}{4}\cdot \frac{\left( M-m\right) \left( M^{p-1}-m^{p-1}\right) }{%
M^{p/2}m^{p/2}}\left\langle Ax,x\right\rangle \left\langle
A^{p-1}x,x\right\rangle , \\ 
\\ 
\left( \sqrt{M}-\sqrt{m}\right) \left( M^{\left( p-1\right) /2}-m^{\left(
p-1\right) /2}\right) \left[ \left\langle Ax,x\right\rangle \left\langle
A^{p-1}x,x\right\rangle \right] ^{\frac{1}{2}},%
\end{array}%
\right.  \notag
\end{align}%
for any $x\in H$ with $\left\Vert x\right\Vert =1.$

Now, if we consider the convex function $f:\left[ 0,\infty \right)
\rightarrow \left[ 0,\infty \right) ,$ $f\left( x\right) =-x^{p}$ with $p\in
\left( 0,1\right) ,$ then from the inequalities (\ref{I.e.3.1}) and (\ref%
{I.e.3.1.0}) and for the positive definite operator $A$ we have the
inequalities

\begin{align}
(0& \leq )\left\langle Ax,x\right\rangle ^{p}-\left\langle
A^{p}x,x\right\rangle  \label{I.e.4.14} \\
& \leq p\times \left\{ 
\begin{array}{l}
\frac{1}{2}\cdot \left( M-m\right) \left[ \left\Vert A^{p-1}x\right\Vert
^{2}-\left\langle A^{p-1}x,x\right\rangle ^{2}\right] ^{1/2} \\ 
\\ 
\frac{1}{2}\cdot \left( m^{p-1}-M^{p-1}\right) \left[ \left\Vert
Ax\right\Vert ^{2}-\left\langle Ax,x\right\rangle ^{2}\right] ^{1/2}%
\end{array}%
\right.  \notag \\
& \left( \leq \frac{1}{4}p\left( M-m\right) \left( m^{p-1}-M^{p-1}\right)
\right)  \notag
\end{align}%
and 
\begin{align}
(0& \leq )\left\langle Ax,x\right\rangle ^{p}-\left\langle
A^{p}x,x\right\rangle  \label{I.e.4.15} \\
& \leq \frac{1}{4}p\left( M-m\right) \left( m^{p-1}-M^{p-1}\right)  \notag \\
& \quad -p\left\{ \hspace{-6pt}%
\begin{array}{l}
\left[ \left\langle Mx-Ax,Ax-mx\right\rangle \left\langle
M^{p-1}x-A^{p-1}x,A^{p-1}x-m^{p-1}x\right\rangle \right] ^{\frac{1}{2}}, \\ 
\\ 
\left\vert \left\langle Ax,x\right\rangle -\frac{M+m}{2}\right\vert
\left\vert \left\langle A^{p-1}x,x\right\rangle -\frac{M^{p-1}+m^{p-1}}{2}%
\right\vert%
\end{array}%
\right.  \notag \\
& \left( \leq \frac{1}{4}p\left( M-m\right) \left( m^{p-1}-M^{p-1}\right)
\right)  \notag
\end{align}%
for any $x\in H$ with $\left\Vert x\right\Vert =1,$ respectively.

Similar results may be stated for the convex function $f:\left( 0,\infty
\right) \rightarrow \left( 0,\infty \right) ,$ $f\left( x\right) =x^{p}$
with $p<0.$ However the details are left to the interested reader.

\section{Some Slater Type Inequalities}

\subsection{Slater Type Inequalities for Functions of Real Variables}

Suppose that $I$ is an interval of real numbers with interior \r{I} and $%
f:I\rightarrow \mathbb{R}$ is a convex function on $I$. Then $f$ is
continuous on \r{I} and has finite left and right derivatives at each point
of \r{I}. Moreover, if $x,y\in $\r{I} and $x<y,$ then $f_{-}^{\prime }\left(
x\right) \leq f_{+}^{\prime }\left( x\right) \leq f_{-}^{\prime }\left(
y\right) \leq f_{+}^{\prime }\left( y\right) $ which shows that both $%
f_{-}^{\prime }$ and $f_{+}^{\prime }$ are nondecreasing function on \r{I}.
It is also known that a convex function must be differentiable except for at
most countably many points.

For a convex function $f:I\rightarrow \mathbb{R}$, the subdifferential of $f$
denoted by $\partial f$ is the set of all functions $\varphi :I\rightarrow %
\left[ -\infty ,\infty \right] $ such that $\varphi \left( \text{\r{I}}%
\right) \subset \mathbb{R}$ and 
\begin{equation*}
f\left( x\right) \geq f\left( a\right) +\left( x-a\right) \varphi \left(
a\right) \text{\quad\ for any }x,a\in I.
\end{equation*}

It is also well known that if $f$ is convex on $I,$ then $\partial f$ is
nonempty, $f_{-}^{\prime }$, $f_{+}^{\prime }\in \partial f$ and if $\varphi
\in \partial f$, then 
\begin{equation*}
f_{-}^{\prime }\left( x\right) \leq \varphi \left( x\right) \leq
f_{+}^{\prime }\left( x\right) \text{\quad\ for any }x\in \text{\r{I}.}
\end{equation*}%
In particular, $\varphi $ is a nondecreasing function.

If $f$ is differentiable and convex on \r{I}, then $\partial f=\left\{
f^{\prime }\right\} .$

The following result is well known in the literature as \textit{the Slater
inequality:}

\begin{theorem}[Slater, 1981, \protect\cite{I.S}]
\label{I.stot.1.1.0}If $f:I\rightarrow \mathbb{R}$ is a nonincreasing
(nondecreasing) convex function, $x_{i}\in I,p_{i}\geq 0$ with $%
P_{n}:=\sum_{i=1}^{n}p_{i}>0$ and $\sum_{i=1}^{n}p_{i}\varphi \left(
x_{i}\right) \neq 0,$ where $\varphi \in \partial f,$ then 
\begin{equation}
\frac{1}{P_{n}}\sum_{i=1}^{n}p_{i}f\left( x_{i}\right) \leq f\left( \frac{%
\sum_{i=1}^{n}p_{i}x_{i}\varphi \left( x_{i}\right) }{\sum_{i=1}^{n}p_{i}%
\varphi \left( x_{i}\right) }\right) .  \label{I.stoe.1.1}
\end{equation}
\end{theorem}

As pointed out in \cite[p. 208]{I.SSDB}, the monotonicity assumption for the
derivative $\varphi $ can be replaced with the condition 
\begin{equation}
\frac{\sum_{i=1}^{n}p_{i}x_{i}\varphi \left( x_{i}\right) }{%
\sum_{i=1}^{n}p_{i}\varphi \left( x_{i}\right) }\in I,  \label{I.stoe.1.2}
\end{equation}%
which is more general and can hold for suitable points in $I$ and for not
necessarily monotonic functions.

\subsection{Some Slater Type Inequalities for Operators}

The following result holds:

\begin{theorem}[Dragomir, 2008, \protect\cite{I.SSD6}]
\label{I.stot.2.1}Let $I$ be an interval and $f:I\rightarrow \mathbb{R}$ be
a convex and differentiable function on \r{I} (the interior of $I)$ whose
derivative $f^{\prime }$ is continuous on \r{I}$.$ If $A$ is a selfadjoint
operator on the Hilbert space $H$ with $Sp\left( A\right) \subseteq \left[
m,M\right] \subset $\r{I} and $f^{\prime }\left( A\right) $ is a positive
definite operator on $H$ then 
\begin{align}
0& \leq f\left( \frac{\left\langle Af^{\prime }\left( A\right)
x,x\right\rangle }{\left\langle f^{\prime }\left( A\right) x,x\right\rangle }%
\right) -\left\langle f\left( A\right) x,x\right\rangle  \label{I.stoe.2.1}
\\
& \leq f^{\prime }\left( \frac{\left\langle Af^{\prime }\left( A\right)
x,x\right\rangle }{\left\langle f^{\prime }\left( A\right) x,x\right\rangle }%
\right) \left[ \frac{\left\langle Af^{\prime }\left( A\right)
x,x\right\rangle -\left\langle Ax,x\right\rangle \left\langle f^{\prime
}\left( A\right) x,x\right\rangle }{\left\langle f^{\prime }\left( A\right)
x,x\right\rangle }\right] ,  \notag
\end{align}%
for any $x\in H$ with $\left\Vert x\right\Vert =1.$
\end{theorem}

\begin{proof}
Since $f$ is convex and differentiable on \r{I}, then we have that 
\begin{equation}
f^{\prime }\left( s\right) \cdot \left( t-s\right) \leq f\left( t\right)
-f\left( s\right) \leq f^{\prime }\left( t\right) \cdot \left( t-s\right)
\label{I.stoe.2.1.a}
\end{equation}%
for any $t,s\in \left[ m,M\right] .$

Now, if we fix $t\in \left[ m,M\right] $ and apply the property (\ref{P})
for the operator $A,$ then for any $x\in H$ with $\left\Vert x\right\Vert =1$
we have 
\begin{align}
\left\langle f^{\prime }\left( A\right) \cdot \left( t\cdot 1_{H}-A\right)
x,x\right\rangle & \leq \left\langle \left[ f\left( t\right) \cdot
1_{H}-f\left( A\right) \right] x,x\right\rangle  \label{I.stoe.2.2} \\
& \leq \left\langle f^{\prime }\left( t\right) \cdot \left( t\cdot
1_{H}-A\right) x,x\right\rangle  \notag
\end{align}%
for any $t\in \left[ m,M\right] $ and any $x\in H$ with $\left\Vert
x\right\Vert =1.$

The inequality (\ref{I.stoe.2.2}) is equivalent with 
\begin{equation}
t\left\langle f^{\prime }\left( A\right) x,x\right\rangle -\left\langle
f^{\prime }\left( A\right) Ax,x\right\rangle \leq f\left( t\right)
-\left\langle f\left( A\right) x,x\right\rangle \leq f^{\prime }\left(
t\right) t-f^{\prime }\left( t\right) \left\langle Ax,x\right\rangle
\label{I.stoe.2.2.a}
\end{equation}%
for any $t\in \left[ m,M\right] $ any $x\in H$ with $\left\Vert x\right\Vert
=1.$

Now, since $A$ is selfadjoint with $mI\leq A\leq MI$ and $f^{\prime }\left(
A\right) $ is positive definite, then $mf^{\prime }\left( A\right) \leq
Af^{\prime }\left( A\right) \leq Mf^{\prime }\left( A\right) ,$ i.e., $%
m\left\langle f^{\prime }\left( A\right) x,x\right\rangle \leq \left\langle
Af^{\prime }\left( A\right) x,x\right\rangle \leq M\left\langle f^{\prime
}\left( A\right) x,x\right\rangle $ for any $x\in H$ with $\left\Vert
x\right\Vert =1,$ which shows that 
\begin{equation*}
t_{0}:=\frac{\left\langle Af^{\prime }\left( A\right) x,x\right\rangle }{%
\left\langle f^{\prime }\left( A\right) x,x\right\rangle }\in \left[ m,M%
\right] \text{ \quad for any }x\in H\text{ \quad with \quad }\left\Vert
x\right\Vert =1.
\end{equation*}%
Finally, if we put $t=t_{0}$ in the equation (\ref{I.stoe.2.2.a}), then we
get the desired result (\ref{I.stoe.2.1}).
\end{proof}

\begin{remark}
\label{I.stor.2.1}It is important to observe that, the condition that $%
f^{\prime }\left( A\right) $ is a positive definite operator on $H$ can be
replaced with the more general assumption that 
\begin{equation}
\frac{\left\langle Af^{\prime }\left( A\right) x,x\right\rangle }{%
\left\langle f^{\prime }\left( A\right) x,x\right\rangle }\in \text{\r{I}
\quad for any }x\in H\text{ \quad with\quad\ }\left\Vert x\right\Vert =1,
\label{I.stoe.2.3}
\end{equation}%
which may be easily verified for particular convex functions $f.$
\end{remark}

\begin{remark}
\label{I.stor.2.2}Now, if the functions is concave on \r{I} and the
condition (\ref{I.stoe.2.3}) holds, then we have the inequality 
\begin{align}
0& \leq \left\langle f\left( A\right) x,x\right\rangle -f\left( \frac{%
\left\langle Af^{\prime }\left( A\right) x,x\right\rangle }{\left\langle
f^{\prime }\left( A\right) x,x\right\rangle }\right)  \label{I.stoe.2.3.1} \\
& \leq f^{\prime }\left( \frac{\left\langle Af^{\prime }\left( A\right)
x,x\right\rangle }{\left\langle f^{\prime }\left( A\right) x,x\right\rangle }%
\right) \left[ \frac{\left\langle Ax,x\right\rangle \left\langle f^{\prime
}\left( A\right) x,x\right\rangle -\left\langle Af^{\prime }\left( A\right)
x,x\right\rangle }{\left\langle f^{\prime }\left( A\right) x,x\right\rangle }%
\right] ,  \notag
\end{align}%
for any $x\in H$ with $\left\Vert x\right\Vert =1.$
\end{remark}

The following examples are of interest:

\begin{example}
\label{I.stoEx.2.1} If $A$ is a positive definite operator on $H$, then 
\begin{equation}
(0\leq )\left\langle \ln Ax,x\right\rangle -\ln \left( \left\langle
A^{-1}x,x\right\rangle ^{-1}\right) \leq \left\langle Ax,x\right\rangle
\cdot \left\langle A^{-1}x,x\right\rangle -1,  \label{I.stoe.2.3.a}
\end{equation}%
for any $x\in H$ with $\left\Vert x\right\Vert =1.$
\end{example}

Indeed, we observe that if we consider the concave function $f:\left(
0,\infty \right) \rightarrow \mathbb{R}$, $f\left( t\right) =\ln t,$ then 
\begin{equation*}
\frac{\left\langle Af^{\prime }\left( A\right) x,x\right\rangle }{%
\left\langle f^{\prime }\left( A\right) x,x\right\rangle }=\frac{1}{%
\left\langle A^{-1}x,x\right\rangle }\in \left( 0,\infty \right) ,\text{
\quad for any }x\in H\text{ \quad with\quad\ }\left\Vert x\right\Vert =1
\end{equation*}%
and by the inequality (\ref{I.stoe.2.3.1}) we deduce the desired result (\ref%
{I.stoe.2.3.a}).

The following example concerning powers of operators is of interest as well:

\begin{example}
\label{I.stoEx.2.2} If $A$ is a positive definite operator on $H,$ then for
any $x\in H$ with $\left\Vert x\right\Vert =1$ we have 
\begin{align}
0& \leq \left\langle A^{p}x,x\right\rangle ^{p-1}-\left\langle
A^{p-1}x,x\right\rangle ^{p}  \label{I.stoe.2.4} \\
& \leq p\left\langle A^{p}x,x\right\rangle ^{p-2}\left[ \left\langle
A^{p}x,x\right\rangle -\left\langle Ax,x\right\rangle \left\langle
A^{p-1}x,x\right\rangle \right]  \notag
\end{align}%
for $p\geq 1,$%
\begin{align}
0& \leq \left\langle A^{p-1}x,x\right\rangle ^{p}-\left\langle
A^{p}x,x\right\rangle ^{p-1}  \label{I.stoe.2.5} \\
& \leq p\left\langle A^{p}x,x\right\rangle ^{p-2}\left[ \left\langle
Ax,x\right\rangle \left\langle A^{p-1}x,x\right\rangle -\left\langle
A^{p}x,x\right\rangle \right]  \notag
\end{align}%
for $0<p<1,$ and 
\begin{align}
0& \leq \left\langle A^{p}x,x\right\rangle ^{p-1}-\left\langle
A^{p-1}x,x\right\rangle ^{p}  \label{I.stoe.2.6} \\
& \leq \left( -p\right) \left\langle A^{p}x,x\right\rangle ^{p-2}\left[
\left\langle Ax,x\right\rangle \left\langle A^{p-1}x,x\right\rangle
-\left\langle A^{p}x,x\right\rangle \right]  \notag
\end{align}%
for $p<0.$
\end{example}

The proof follows from the inequalities (\ref{I.stoe.2.1}) and (\ref%
{I.stoe.2.3.1}) for the convex (concave) function $f\left( t\right)
=t^{p},p\in \left( -\infty ,0\right) \cup \lbrack 1,\infty )$ $\left( p\in
\left( 0,1\right) \right) $ by performing the required calculation. The
details are omitted.

\subsection{Further Reverses}

The following results that provide perhaps more useful upper bounds for the
nonnegative quantity 
\begin{equation*}
f\left( \frac{\left\langle Af^{\prime }\left( A\right) x,x\right\rangle }{%
\left\langle f^{\prime }\left( A\right) x,x\right\rangle }\right)
-\left\langle f\left( A\right) x,x\right\rangle \text{ \quad for }x\in H%
\text{\quad\ with \quad }\left\Vert x\right\Vert =1,
\end{equation*}%
can be stated:

\begin{theorem}[Dragomir, 2008, \protect\cite{I.SSD6}]
\label{I.stot.3.1}Let $I$ be an interval and $f:I\rightarrow \mathbb{R}$ be
a convex and differentiable function on \r{I} (the interior of $I)$ whose
derivative $f^{\prime }$ is continuous on \r{I}$.$ Assume that $A$ is a
selfadjoint operator on the Hilbert space $H$ with $Sp\left( A\right)
\subseteq \left[ m,M\right] \subset $\r{I} and $f^{\prime }\left( A\right) $
is a positive definite operator on $H.$ If we define 
\begin{equation*}
B\left( f^{\prime },A;x\right) :=\frac{1}{\left\langle f^{\prime }\left(
A\right) x,x\right\rangle }\cdot f^{\prime }\left( \frac{\left\langle
Af^{\prime }\left( A\right) x,x\right\rangle }{\left\langle f^{\prime
}\left( A\right) x,x\right\rangle }\right)
\end{equation*}%
then 
\begin{align}
(0& \leq )f\left( \frac{\left\langle Af^{\prime }\left( A\right)
x,x\right\rangle }{\left\langle f^{\prime }\left( A\right) x,x\right\rangle }%
\right) -\left\langle f\left( A\right) x,x\right\rangle  \label{I.stoe.3.1}
\\
& \leq B\left( f^{\prime },A;x\right) \times \left\{ 
\begin{array}{l}
\frac{1}{2}\cdot \left( M-m\right) \left[ \left\Vert f^{\prime }\left(
A\right) x\right\Vert ^{2}-\left\langle f^{\prime }\left( A\right)
x,x\right\rangle ^{2}\right] ^{1/2} \\[10pt] 
\frac{1}{2}\cdot \left( f^{\prime }\left( M\right) -f^{\prime }\left(
m\right) \right) \left[ \left\Vert Ax\right\Vert ^{2}-\left\langle
Ax,x\right\rangle ^{2}\right] ^{1/2}%
\end{array}%
\right.  \notag \\
& \leq \frac{1}{4}\left( M-m\right) \left( f^{\prime }\left( M\right)
-f^{\prime }\left( m\right) \right) B\left( f^{\prime },A;x\right)  \notag
\end{align}%
and 
\begin{align}
(0& \leq )f\left( \frac{\left\langle Af^{\prime }\left( A\right)
x,x\right\rangle }{\left\langle f^{\prime }\left( A\right) x,x\right\rangle }%
\right) -\left\langle f\left( A\right) x,x\right\rangle  \label{I.stoe.3.1.0}
\\
& \leq B\left( f^{\prime },A;x\right) \times \left[ \frac{1}{4}\left(
M-m\right) \left( f^{\prime }\left( M\right) -f^{\prime }\left( m\right)
\right) \right.  \notag \\
& \left. -\left\{ \hspace{-6pt}%
\begin{array}{l}
\left[ \left\langle Mx-Ax,Ax-mx\right\rangle \left\langle f^{\prime }\left(
M\right) x-f^{\prime }\left( A\right) x,f^{\prime }\left( A\right)
x-f^{\prime }\left( m\right) x\right\rangle \right] ^{\frac{1}{2}}, \\[10pt] 
\left\vert \left\langle Ax,x\right\rangle -\frac{M+m}{2}\right\vert
\left\vert \left\langle f^{\prime }\left( A\right) x,x\right\rangle -\frac{%
f^{\prime }\left( M\right) +f^{\prime }\left( m\right) }{2}\right\vert%
\end{array}%
\right. \hspace{-8pt}\right]  \notag \\
& \leq \frac{1}{4}\left( M-m\right) \left( f^{\prime }\left( M\right)
-f^{\prime }\left( m\right) \right) B\left( f^{\prime },A;x\right) ,  \notag
\end{align}%
for any $x\in H$ with $\left\Vert x\right\Vert =1,$ respectively.

Moreover, if $A$ is a positive definite operator, then 
\begin{align}
(0& \leq )f\left( \frac{\left\langle Af^{\prime }\left( A\right)
x,x\right\rangle }{\left\langle f^{\prime }\left( A\right) x,x\right\rangle }%
\right) -\left\langle f\left( A\right) x,x\right\rangle  \label{I.stoe.3.1.1}
\\
& \hspace{-3pt}\leq B\hspace{-2pt}\left( f^{\prime },A;x\right)  \notag \\
& \times \hspace{-3pt}\left\{ \hspace{-8pt}%
\begin{array}{l}
\frac{1}{4}\cdot \frac{\left( M-m\right) \left( f^{\prime }\left( M\right)
-f^{\prime }\left( m\right) \right) }{\sqrt{Mmf^{\prime }\left( M\right)
f^{\prime }\left( m\right) }}\left\langle Ax,x\right\rangle \left\langle
f^{\prime }\left( A\right) x,x\right\rangle , \\[10pt] 
\left( \hspace{-2pt}\sqrt{M}\hspace{-2pt}-\hspace{-2pt}\sqrt{m}\right) 
\hspace{-3pt}\left( \hspace{-2pt}\sqrt{f^{\prime }\left( M\right) }\hspace{%
-2pt}-\hspace{-2pt}\sqrt{f^{\prime }\left( m\right) }\right) \hspace{-2pt}%
\left[ \left\langle Ax,x\right\rangle \hspace{-2pt}\left\langle f^{\prime
}\left( A\right) x,x\right\rangle \right] ^{\frac{1}{2}},%
\end{array}%
\right.  \notag
\end{align}%
for any $x\in H$ with $\left\Vert x\right\Vert =1.$
\end{theorem}

\begin{proof}
We use the following Gr\"{u}ss' type result we obtained in \cite{I.SSD1}:

Let $A$ be a selfadjoint operator on the Hilbert space $\left(
H;\left\langle .,.\right\rangle \right) $ and assume that $Sp\left( A\right)
\subseteq \left[ m,M\right] $ for some scalars $m<M.$ If $h$ and $g$ are
continuous on $\left[ m,M\right] $ and $\gamma :=\min_{t\in \left[ m,M\right]
}h\left( t\right) $ and $\Gamma :=\max_{t\in \left[ m,M\right] }h\left(
t\right) ,$ then 
\begin{align}
& \left\vert \left\langle h\left( A\right) g\left( A\right) x,x\right\rangle
-\left\langle h\left( A\right) x,x\right\rangle \cdot \left\langle g\left(
A\right) x,x\right\rangle \right\vert  \label{I.stoe.3.2} \\
& \leq \frac{1}{2}\cdot \left( \Gamma -\gamma \right) \left[ \left\Vert
g\left( A\right) x\right\Vert ^{2}-\left\langle g\left( A\right)
x,x\right\rangle ^{2}\right] ^{1/2}  \notag \\
& \left( \leq \frac{1}{4}\left( \Gamma -\gamma \right) \left( \Delta -\delta
\right) \right) ,  \notag
\end{align}%
for each $x\in H$ with $\left\Vert x\right\Vert =1,$ where $\delta
:=\min_{t\in \left[ m,M\right] }g\left( t\right) $ and $\Delta :=\max_{t\in %
\left[ m,M\right] }g\left( t\right) .$

Therefore, we can state that 
\begin{align}
& \left\langle Af^{\prime }\left( A\right) x,x\right\rangle -\left\langle
Ax,x\right\rangle \cdot \left\langle f^{\prime }\left( A\right)
x,x\right\rangle  \label{I.stoe.3.3} \\
& \leq \frac{1}{2}\cdot \left( M-m\right) \left[ \left\Vert f^{\prime
}\left( A\right) x\right\Vert ^{2}-\left\langle f^{\prime }\left( A\right)
x,x\right\rangle ^{2}\right] ^{1/2}  \notag \\
& \leq \frac{1}{4}\left( M-m\right) \left( f^{\prime }\left( M\right)
-f^{\prime }\left( m\right) \right)  \notag
\end{align}%
and 
\begin{align}
& \left\langle Af^{\prime }\left( A\right) x,x\right\rangle -\left\langle
Ax,x\right\rangle \cdot \left\langle f^{\prime }\left( A\right)
x,x\right\rangle  \label{I.stoe.3.4} \\
& \leq \frac{1}{2}\cdot \left( f^{\prime }\left( M\right) -f^{\prime }\left(
m\right) \right) \left[ \left\Vert Ax\right\Vert ^{2}-\left\langle
Ax,x\right\rangle ^{2}\right] ^{1/2}  \notag \\
& \leq \frac{1}{4}\left( M-m\right) \left( f^{\prime }\left( M\right)
-f^{\prime }\left( m\right) \right) ,  \notag
\end{align}%
for each $x\in H$ with $\left\Vert x\right\Vert =1,$ which together with (%
\ref{I.stoe.2.1}) provide the desired result (\ref{I.stoe.3.1}).

On making use of the inequality obtained in \cite{I.SSD2} 
\begin{align}
& \left\vert \left\langle h\left( A\right) g\left( A\right) x,x\right\rangle
-\left\langle h\left( A\right) x,x\right\rangle \left\langle g\left(
A\right) x,x\right\rangle \right\vert  \label{I.stoe.3.4.1} \\
& \leq \frac{1}{4}\cdot \left( \Gamma -\gamma \right) \left( \Delta -\delta
\right)  \notag \\
& -\left\{ 
\begin{array}{l}
\left[ \left\langle \Gamma x-h\left( A\right) x,f\left( A\right) x-\gamma
x\right\rangle \left\langle \Delta x-g\left( A\right) x,g\left( A\right)
x-\delta x\right\rangle \right] ^{\frac{1}{2}}, \\ 
\\ 
\left\vert \left\langle h\left( A\right) x,x\right\rangle -\frac{\Gamma
+\gamma }{2}\right\vert \left\vert \left\langle g\left( A\right)
x,x\right\rangle -\frac{\Delta +\delta }{2}\right\vert ,%
\end{array}%
\right.  \notag
\end{align}%
for each $x\in H$ with $\left\Vert x\right\Vert =1,$ we can state that 
\begin{align*}
& \left\langle Af^{\prime }\left( A\right) x,x\right\rangle -\left\langle
Ax,x\right\rangle \cdot \left\langle f^{\prime }\left( A\right)
x,x\right\rangle \\
& \leq \frac{1}{4}\left( M-m\right) \left( f^{\prime }\left( M\right)
-f^{\prime }\left( m\right) \right) \\
& -\left\{ 
\begin{array}{l}
\left[ \left\langle Mx-Ax,Ax-mx\right\rangle \left\langle f^{\prime }\left(
M\right) x-f^{\prime }\left( A\right) x,f^{\prime }\left( A\right)
x-f^{\prime }\left( m\right) x\right\rangle \right] ^{\frac{1}{2}}, \\ 
\\ 
\left\vert \left\langle Ax,x\right\rangle -\frac{M+m}{2}\right\vert
\left\vert \left\langle f^{\prime }\left( A\right) x,x\right\rangle -\frac{%
f^{\prime }\left( M\right) +f^{\prime }\left( m\right) }{2}\right\vert ,%
\end{array}%
\right.
\end{align*}%
for each $x\in H$ with $\left\Vert x\right\Vert =1,$ which together with (%
\ref{I.stoe.2.1}) provide the desired result (\ref{I.stoe.3.1.0}).

Further, in order to prove the third inequality, we make use of the
following result of Gr\"{u}ss type obtained in \cite{I.SSD2}:

If $\gamma $ and $\delta $ are positive, then 
\begin{align}
& \left\vert \left\langle h\left( A\right) g\left( A\right) x,x\right\rangle
-\left\langle h\left( A\right) x,x\right\rangle \left\langle g\left(
A\right) x,x\right\rangle \right\vert  \label{I.stoe.3.5} \\
& \leq \left\{ 
\begin{array}{l}
\frac{1}{4}\cdot \frac{\left( \Gamma -\gamma \right) \left( \Delta -\delta
\right) }{\sqrt{\Gamma \gamma \Delta \delta }}\left\langle h\left( A\right)
x,x\right\rangle \left\langle g\left( A\right) x,x\right\rangle , \\ 
\\ 
\left( \sqrt{\Gamma }-\sqrt{\gamma }\right) \left( \sqrt{\Delta }-\sqrt{%
\delta }\right) \left[ \left\langle h\left( A\right) x,x\right\rangle
\left\langle g\left( A\right) x,x\right\rangle \right] ^{\frac{1}{2}},%
\end{array}%
\right.  \notag
\end{align}%
for each $x\in H$ with $\left\Vert x\right\Vert =1.$

Now, on making use of (\ref{I.stoe.3.5}) we can state that 
\begin{align}
& \left\langle Af^{\prime }\left( A\right) x,x\right\rangle -\left\langle
Ax,x\right\rangle \cdot \left\langle f^{\prime }\left( A\right)
x,x\right\rangle  \label{I.stoe.3.6} \\
& \leq \left\{ 
\begin{array}{l}
\frac{1}{4}\cdot \frac{\left( M-m\right) \left( f^{\prime }\left( M\right)
-f^{\prime }\left( m\right) \right) }{\sqrt{Mmf^{\prime }\left( M\right)
f^{\prime }\left( m\right) }}\left\langle Ax,x\right\rangle \left\langle
f^{\prime }\left( A\right) x,x\right\rangle , \\ 
\\ 
\left( \sqrt{M}-\sqrt{m}\right) \left( \sqrt{f^{\prime }\left( M\right) }-%
\sqrt{f^{\prime }\left( m\right) }\right) \left[ \left\langle
Ax,x\right\rangle \left\langle f^{\prime }\left( A\right) x,x\right\rangle %
\right] ^{\frac{1}{2}},%
\end{array}%
\right.  \notag
\end{align}%
for each $x\in H$ with $\left\Vert x\right\Vert =1,$ which together with (%
\ref{I.stoe.2.1}) provide the desired result (\ref{I.stoe.3.1.1}).
\end{proof}

\begin{remark}
\label{I.stor.3.1}We observe, from the first inequality in (\ref%
{I.stoe.3.1.1}), that 
\begin{equation*}
\left( 1\leq \right) \frac{\left\langle Af^{\prime }\left( A\right)
x,x\right\rangle }{\left\langle Ax,x\right\rangle \left\langle f^{\prime
}\left( A\right) x,x\right\rangle }\leq \frac{1}{4}\cdot \frac{\left(
M-m\right) \left( f^{\prime }\left( M\right) -f^{\prime }\left( m\right)
\right) }{\sqrt{Mmf^{\prime }\left( M\right) f^{\prime }\left( m\right) }}+1
\end{equation*}%
which implies that 
\begin{equation*}
f^{\prime }\left( \frac{\left\langle Af^{\prime }\left( A\right)
x,x\right\rangle }{\left\langle f^{\prime }\left( A\right) x,x\right\rangle }%
\right) \leq f^{\prime }\left( \left[ \frac{1}{4}\cdot \frac{\left(
M-m\right) \left( f^{\prime }\left( M\right) -f^{\prime }\left( m\right)
\right) }{\sqrt{Mmf^{\prime }\left( M\right) f^{\prime }\left( m\right) }}+1%
\right] \left\langle Ax,x\right\rangle \right) ,
\end{equation*}%
for each $x\in H$ with $\left\Vert x\right\Vert =1,$ since $f^{\prime }$ is
monotonic nondecreasing and $A$ is positive definite.

Now, the first inequality in (\ref{I.stoe.3.1.1}) implies the following
result 
\begin{align}
(0& \leq )f\left( \frac{\left\langle Af^{\prime }\left( A\right)
x,x\right\rangle }{\left\langle f^{\prime }\left( A\right) x,x\right\rangle }%
\right) -\left\langle f\left( A\right) x,x\right\rangle  \label{I.stoe.3.7}
\\
& \leq \frac{1}{4}\cdot \frac{\left( M-m\right) \left( f^{\prime }\left(
M\right) -f^{\prime }\left( m\right) \right) }{\sqrt{Mmf^{\prime }\left(
M\right) f^{\prime }\left( m\right) }}  \notag \\
& \times f^{\prime }\left( \left[ \frac{1}{4}\cdot \frac{\left( M-m\right)
\left( f^{\prime }\left( M\right) -f^{\prime }\left( m\right) \right) }{%
\sqrt{Mmf^{\prime }\left( M\right) f^{\prime }\left( m\right) }}+1\right]
\left\langle Ax,x\right\rangle \right) \left\langle Ax,x\right\rangle , 
\notag
\end{align}%
for each $x\in H$ with $\left\Vert x\right\Vert =1.$

From the second inequality in (\ref{I.stoe.3.1.1}) we also have 
\begin{align}
(0& \leq )f\left( \frac{\left\langle Af^{\prime }\left( A\right)
x,x\right\rangle }{\left\langle f^{\prime }\left( A\right) x,x\right\rangle }%
\right) -\left\langle f\left( A\right) x,x\right\rangle  \label{I.stoe.3.8}
\\
& \leq \left( \sqrt{M}-\sqrt{m}\right) \left( \sqrt{f^{\prime }\left(
M\right) }-\sqrt{f^{\prime }\left( m\right) }\right)  \notag \\
& \times f^{\prime }\left( \left[ \frac{1}{4}\cdot \frac{\left( M-m\right)
\left( f^{\prime }\left( M\right) -f^{\prime }\left( m\right) \right) }{%
\sqrt{Mmf^{\prime }\left( M\right) f^{\prime }\left( m\right) }}+1\right]
\left\langle Ax,x\right\rangle \right) \left[ \frac{\left\langle
Ax,x\right\rangle }{\left\langle f^{\prime }\left( A\right) x,x\right\rangle 
}\right] ^{\frac{1}{2}},  \notag
\end{align}%
for each $x\in H$ with $\left\Vert x\right\Vert =1.$
\end{remark}

\begin{remark}
\label{I.stor.3.2}If the condition that $f^{\prime }\left( A\right) $ is a
positive definite operator on $H$ from the Theorem \ref{I.stot.3.1} is
replaced by the condition (\ref{I.stoe.2.3}), then the inequalities (\ref%
{I.stoe.3.1}) and (\ref{I.stoe.3.2}) will still hold. Similar inequalities
for concave functions can be stated. However, the details are not provided
here.
\end{remark}

\subsection{Multivariate Versions}

The following result for sequences of operators can be stated.

\begin{theorem}[Dragomir, 2008, \protect\cite{I.SSD6}]
\label{I.stot.4.1}Let $I$ be an interval and $f:I\rightarrow \mathbb{R}$ be
a convex and differentiable function on \r{I} (the interior of $I)$ whose
derivative $f^{\prime }$ is continuous on \r{I}$.$ If $A_{j},j\in \left\{
1,\dots ,n\right\} $ are selfadjoint operators on the Hilbert space $H$ with 
$Sp\left( A_{j}\right) \subseteq \left[ m,M\right] \subset $\r{I} and 
\begin{equation}
\frac{\sum_{j=1}^{n}\left\langle A_{j}f^{\prime }\left( A_{j}\right)
x_{j},x_{j}\right\rangle }{\sum_{j=1}^{n}\left\langle f^{\prime }\left(
A_{j}\right) x_{j},x_{j}\right\rangle }\in \text{\r{I}}  \label{I.stoe.4.0}
\end{equation}%
for each $x_{j}\in H,j\in \left\{ 1,\dots ,n\right\} $ with $%
\sum_{j=1}^{n}\left\Vert x_{j}\right\Vert ^{2}=1,$ then 
\begin{align}
0& \leq f\left( \frac{\sum_{j=1}^{n}\left\langle A_{j}f^{\prime }\left(
A_{j}\right) x_{j},x_{j}\right\rangle }{\sum_{j=1}^{n}\left\langle f^{\prime
}\left( A_{j}\right) x_{j},x_{j}\right\rangle }\right)
-\sum_{j=1}^{n}\left\langle f\left( A_{j}\right) x_{j},x_{j}\right\rangle
\label{I.stoe.4.1} \\
& \leq f^{\prime }\left( \frac{\sum_{j=1}^{n}\left\langle A_{j}f^{\prime
}\left( A_{j}\right) x_{j},x_{j}\right\rangle }{\sum_{j=1}^{n}\left\langle
f^{\prime }\left( A_{j}\right) x_{j},x_{j}\right\rangle }\right)  \notag \\
& \times \left[ \frac{\sum_{j=1}^{n}\left\langle A_{j}f^{\prime }\left(
A_{j}\right) x_{j},x_{j}\right\rangle -\sum_{j=1}^{n}\left\langle
A_{j}x_{j},x_{j}\right\rangle \sum_{j=1}^{n}\left\langle f^{\prime }\left(
A_{j}\right) x_{j},x_{j}\right\rangle }{\sum_{j=1}^{n}\left\langle f^{\prime
}\left( A_{j}\right) x_{j},x_{j}\right\rangle }\right] ,  \notag
\end{align}%
for each $x_{j}\in H,j\in \left\{ 1,\dots ,n\right\} $ with $%
\sum_{j=1}^{n}\left\Vert x_{j}\right\Vert ^{2}=1.$
\end{theorem}

\begin{proof}
Follows from Theorem \ref{I.stot.2.1}. The details are omitted.
\end{proof}

The following particular case is of interest

\begin{corollary}[Dragomir, 2008, \protect\cite{I.SSD6}]
\label{I.stoc.4.1}Let $I$ be an interval and $f:I\rightarrow \mathbb{R}$ be
a convex and differentiable function on \r{I} (the interior of $I)$ whose
derivative $f^{\prime }$ is continuous on \r{I}$.$ If $A_{j},j\in \left\{
1,\dots ,n\right\} $ are selfadjoint operators on the Hilbert space $H$ with 
$Sp\left( A_{j}\right) \subseteq \left[ m,M\right] \subset $\r{I} and for $%
p_{j}\geq 0$ with $\sum_{j=1}^{n}p_{j}=1$ if we also assume that 
\begin{equation}
\frac{\left\langle \sum_{j=1}^{n}p_{j}A_{j}f^{\prime }\left( A_{j}\right)
x,x\right\rangle }{\left\langle \sum_{j=1}^{n}p_{j}f^{\prime }\left(
A_{j}\right) x,x\right\rangle }\in \text{\r{I}}  \label{I.stoe.4.2}
\end{equation}%
for each $x\in H$ with $\left\Vert x\right\Vert =1,$ then 
\begin{align}
0& \leq f\left( \frac{\left\langle \sum_{j=1}^{n}p_{j}A_{j}f^{\prime }\left(
A_{j}\right) x,x\right\rangle }{\left\langle \sum_{j=1}^{n}p_{j}f^{\prime
}\left( A_{j}\right) x,x\right\rangle }\right) -\left\langle
\sum_{j=1}^{n}p_{j}f\left( A_{j}\right) x,x\right\rangle  \label{I.stoe.4.3}
\\
& \leq f^{\prime }\left( \frac{\left\langle
\sum_{j=1}^{n}p_{j}A_{j}f^{\prime }\left( A_{j}\right) x,x\right\rangle }{%
\left\langle \sum_{j=1}^{n}p_{j}f^{\prime }\left( A_{j}\right)
x,x\right\rangle }\right)  \notag \\
& \times \hspace{-2pt}\left[ \hspace{-2pt}\frac{\left\langle
\sum_{j=1}^{n}p_{j}A_{j}f^{\prime }\left( A_{j}\right) x,x\right\rangle 
\hspace{-2pt}-\hspace{-2pt}\left\langle
\sum_{j=1}^{n}p_{j}A_{j}x,x\right\rangle \left\langle
\sum_{j=1}^{n}p_{j}f^{\prime }\left( A_{j}\right) x,x\right\rangle }{%
\left\langle \sum_{j=1}^{n}p_{j}f^{\prime }\left( A_{j}\right)
x,x\right\rangle }\right] ,  \notag
\end{align}%
for each $x\in H$ with $\left\Vert x\right\Vert =1.$
\end{corollary}

\begin{proof}
Follows from Theorem \ref{I.stot.4.1} on choosing $x_{j}=\sqrt{p_{j}}\cdot
x, $ $j\in \left\{ 1,\dots ,n\right\} ,$ where $p_{j}\geq 0,j\in \left\{
1,\dots ,n\right\} ,$ $\sum_{j=1}^{n}p_{j}=1$ and $x\in H,$ with $\left\Vert
x\right\Vert =1.$ The details are omitted.
\end{proof}

The following examples are interesting in themselves:

\begin{example}
\label{I.stoEx 4.1} If $A_{j}$, $j\in \left\{ 1,\dots ,n\right\} $ are
positive definite operators on $H$, then 
\begin{align}
(0& \leq )\sum_{j=1}^{n}\left\langle \ln A_{j}x_{j},x_{j}\right\rangle -\ln 
\left[ \left( \sum_{j=1}^{n}\left\langle A_{j}^{-1}x_{j},x_{j}\right\rangle
\right) ^{-1}\right]  \label{I.stoe.4.4} \\
& \leq \sum_{j=1}^{n}\left\langle A_{j}x_{j},x_{j}\right\rangle \cdot
\sum_{j=1}^{n}\left\langle A_{j}^{-1}x_{j},x_{j}\right\rangle -1,  \notag
\end{align}%
for each $x_{j}\in H,j\in \left\{ 1,\dots ,n\right\} $ with $%
\sum_{j=1}^{n}\left\Vert x_{j}\right\Vert ^{2}=1.$

If $p_{j}\geq 0,j\in \left\{ 1,\dots ,n\right\} $ with $%
\sum_{j=1}^{n}p_{j}=1,$ then we also have the inequality 
\begin{align}
(0& \leq )\left\langle \sum_{j=1}^{n}p_{j}\ln A_{j}x,x\right\rangle -\ln 
\left[ \left( \left\langle \sum_{j=1}^{n}p_{j}A_{j}^{-1}x,x\right\rangle
\right) ^{-1}\right]  \label{I.stoe.4.5} \\
& \leq \left\langle \sum_{j=1}^{n}p_{j}A_{j}x,x\right\rangle \cdot
\left\langle \sum_{j=1}^{n}p_{j}A_{j}^{-1}x,x\right\rangle -1,  \notag
\end{align}%
for each $x\in H$ with $\left\Vert x\right\Vert =1.$
\end{example}

\section{Other Inequalities for Convex Functions}

\subsection{Some Inequalities for Two Operators}

The following result holds:

\begin{theorem}[Dragomir, 2008, \protect\cite{I.SSD7}]
\label{I.cfot.2.1}Let $I$ be an interval and $f:I\rightarrow \mathbb{R}$ be
a convex and differentiable function on \r{I} (the interior of $I)$ whose
derivative $f^{\prime }$ is continuous on \r{I}$.$ If $A$ and $B$ are
selfadjoint operators on the Hilbert space $H$ with $Sp\left( A\right)
,Sp\left( B\right) \subseteq \left[ m,M\right] \subset $\r{I}$,$ then 
\begin{align}
\left\langle f^{\prime }\left( A\right) x,x\right\rangle \left\langle
By,y\right\rangle -\left\langle f^{\prime }\left( A\right) Ax,x\right\rangle
& \leq \left\langle f\left( B\right) y,y\right\rangle -\left\langle f\left(
A\right) x,x\right\rangle  \label{I.cfoe.2.1} \\
& \leq \left\langle f^{\prime }\left( B\right) By,y\right\rangle
-\left\langle Ax,x\right\rangle \left\langle f^{\prime }\left( B\right)
y,y\right\rangle  \notag
\end{align}%
for any $x,y\in H$ with $\left\Vert x\right\Vert =\left\Vert y\right\Vert
=1. $

In particular, we have 
\begin{align}
\left\langle f^{\prime }\left( A\right) x,x\right\rangle \left\langle
Ay,y\right\rangle -\left\langle f^{\prime }\left( A\right) Ax,x\right\rangle
& \leq \left\langle f\left( A\right) y,y\right\rangle -\left\langle f\left(
A\right) x,x\right\rangle  \label{I.cfoe.2.1.1} \\
& \leq \left\langle f^{\prime }\left( A\right) Ay,y\right\rangle
-\left\langle Ax,x\right\rangle \left\langle f^{\prime }\left( A\right)
y,y\right\rangle  \notag
\end{align}%
for any $x,y\in H$ with $\left\Vert x\right\Vert =\left\Vert y\right\Vert =1$
and 
\begin{align}
\left\langle f^{\prime }\left( A\right) x,x\right\rangle \left\langle
Bx,x\right\rangle -\left\langle f^{\prime }\left( A\right) Ax,x\right\rangle
& \leq \left\langle f\left( B\right) x,x\right\rangle -\left\langle f\left(
A\right) x,x\right\rangle  \label{I.cfoe.2.1.2} \\
& \leq \left\langle f^{\prime }\left( B\right) Bx,x\right\rangle
-\left\langle Ax,x\right\rangle \left\langle f^{\prime }\left( B\right)
x,x\right\rangle  \notag
\end{align}%
for any $x\in H$ with $\left\Vert x\right\Vert =1.$
\end{theorem}

\begin{proof}
Since $f$ is convex and differentiable on \r{I}, then we have that 
\begin{equation}
f^{\prime }\left( s\right) \cdot \left( t-s\right) \leq f\left( t\right)
-f\left( s\right) \leq f^{\prime }\left( t\right) \cdot \left( t-s\right)
\label{I.cfoe.2.1.a}
\end{equation}%
for any $t,s\in \left[ m,M\right] .$

Now, if we fix $t\in \left[ m,M\right] $ and apply the property (\ref{P})
for the operator $A,$ then for any $x\in H$ with $\left\Vert x\right\Vert =1$
we have 
\begin{align}
\left\langle f^{\prime }\left( A\right) \cdot \left( t\cdot 1_{H}-A\right)
x,x\right\rangle & \leq \left\langle \left[ f\left( t\right) \cdot
1_{H}-f\left( A\right) \right] x,x\right\rangle  \label{I.cfoe.2.2} \\
& \leq \left\langle f^{\prime }\left( t\right) \cdot \left( t\cdot
1_{H}-A\right) x,x\right\rangle  \notag
\end{align}%
for any $t\in \left[ m,M\right] $ and any $x\in H$ with $\left\Vert
x\right\Vert =1.$

The inequality (\ref{I.cfoe.2.2}) is equivalent with 
\begin{equation}
t\left\langle f^{\prime }\left( A\right) x,x\right\rangle -\left\langle
f^{\prime }\left( A\right) Ax,x\right\rangle \leq f\left( t\right)
-\left\langle f\left( A\right) x,x\right\rangle \leq f^{\prime }\left(
t\right) t-f^{\prime }\left( t\right) \left\langle Ax,x\right\rangle
\label{I.cfoe.2.2.a}
\end{equation}%
for any $t\in \left[ m,M\right] $ and any $x\in H$ with $\left\Vert
x\right\Vert =1.$

If we fix $x\in H$ with $\left\Vert x\right\Vert =1$ in (\ref{I.cfoe.2.2.a})
and apply the property (\ref{P}) for the operator $B,$ then we get 
\begin{align*}
\left\langle \left[ \left\langle f^{\prime }\left( A\right) x,x\right\rangle
B-\left\langle f^{\prime }\left( A\right) Ax,x\right\rangle 1_{H}\right]
y,y\right\rangle & \leq \left\langle \left[ f\left( B\right) -\left\langle
f\left( A\right) x,x\right\rangle 1_{H}\right] y,y\right\rangle \\
& \leq \left\langle \left[ f^{\prime }\left( B\right) B-\left\langle
Ax,x\right\rangle f^{\prime }\left( B\right) \right] y,y\right\rangle
\end{align*}%
for each $y\in H$ with $\left\Vert y\right\Vert =1,$ which is clearly
equivalent to the desired inequality (\ref{I.cfoe.2.1}).
\end{proof}

\begin{remark}
\label{I.cfor.2.1}If we fix $x\in H$ with $\left\Vert x\right\Vert =1$ and
choose $B=\left\langle Ax,x\right\rangle \cdot 1_{H},$ then we obtain from
the first inequality in (\ref{I.cfoe.2.1}) the reverse of the Mond-Pe\v{c}ari%
\'{c} inequality obtained by the author in \cite{I.SSD5}. The second
inequality will provide the Mond-Pe\v{c}ari\'{c} inequality for convex
functions whose derivatives are continuous.
\end{remark}

The following corollary is of interest:

\begin{corollary}
\label{I.cfoc.2.1}Let $I$ be an interval and $f:I\rightarrow \mathbb{R}$ be
a convex and differentiable function on \r{I} whose derivative $f^{\prime }$
is continuous on \r{I}$.$ Also, suppose that $A$ is a selfadjoint operator
on the Hilbert space $H$ with $Sp\left( A\right) \subseteq \left[ m,M\right]
\subset $\r{I}$.$ If $g$ is nonincreasing and continuous on $\left[ m,M%
\right] $ and 
\begin{equation}
f^{\prime }\left( A\right) \left[ g\left( A\right) -A\right] \geq 0\text{ }
\label{I.cfoe.2.3}
\end{equation}%
in the operator order of $B\left( H\right) ,$ then 
\begin{equation}
\left( f\circ g\right) \left( A\right) \geq f\left( A\right)
\label{I.cfoe.2.4}
\end{equation}%
in the operator order of $B\left( H\right) .$
\end{corollary}

\begin{proof}
If we apply the first inequality from (\ref{I.cfoe.2.1.2}) for $B=g\left(
A\right) $ we have 
\begin{equation}
\left\langle f^{\prime }\left( A\right) x,x\right\rangle \left\langle
g\left( A\right) x,x\right\rangle -\left\langle f^{\prime }\left( A\right)
Ax,x\right\rangle \leq \left\langle f\left( g\left( A\right) \right)
x,x\right\rangle -\left\langle f\left( A\right) x,x\right\rangle
\label{I.cfoe.2.5}
\end{equation}%
any $x\in H$ with $\left\Vert x\right\Vert =1.$

We use the following \v{C}eby\v{s}ev type inequality for functions of
operators established by the author in \cite{I.SSD4}:

Let $A$ be a selfadjoint operator with $Sp\left( A\right) \subseteq \left[
m,M\right] $ for some real numbers $m<M.$ If $h,g:\left[ m,M\right]
\longrightarrow \mathbb{R}$ are continuous and \textit{synchronous
(asynchronous) }on\textit{\ }$\left[ m,M\right] ,$ then 
\begin{equation}
\left\langle h\left( A\right) g\left( A\right) x,x\right\rangle \geq \left(
\leq \right) \left\langle h\left( A\right) x,x\right\rangle \cdot
\left\langle g\left( A\right) x,x\right\rangle  \label{I.cfoe.2.6}
\end{equation}%
for any $x\in H$ with $\left\Vert x\right\Vert =1.$

Now, since $f^{\prime }$ and $g$ are continuous and are asynchronous on $%
\left[ m,M\right] ,$ then by (\ref{I.cfoe.2.6}) we have the inequality 
\begin{equation}
\left\langle f^{\prime }\left( A\right) g\left( A\right) x,x\right\rangle
\leq \left\langle f^{\prime }\left( A\right) x,x\right\rangle \cdot
\left\langle g\left( A\right) x,x\right\rangle  \label{I.cfoe.2.7}
\end{equation}%
for any $x\in H$ with $\left\Vert x\right\Vert =1.$

Subtracting from both sides of (\ref{I.cfoe.2.7}) the quantity $\left\langle
f^{\prime }\left( A\right) Ax,x\right\rangle $ and taking into account, by (%
\ref{I.cfoe.2.3}), that $\left\langle f^{\prime }\left( A\right) \left[
g\left( A\right) -A\right] x,x\right\rangle \geq 0$ for any $x\in H$ with $%
\left\Vert x\right\Vert =1,$ we then\ have 
\begin{align*}
0& \leq \left\langle f^{\prime }\left( A\right) \left[ g\left( A\right) -A%
\right] x,x\right\rangle \\
& =\left\langle f^{\prime }\left( A\right) g\left( A\right) x,x\right\rangle
-\left\langle f^{\prime }\left( A\right) Ax,x\right\rangle \\
& \leq \left\langle f^{\prime }\left( A\right) x,x\right\rangle \cdot
\left\langle g\left( A\right) x,x\right\rangle -\left\langle f^{\prime
}\left( A\right) Ax,x\right\rangle
\end{align*}%
which together with (\ref{I.cfoe.2.5}) will produce the desired result (\ref%
{I.cfoe.2.4}).
\end{proof}

We provide now some particular inequalities of interest that can be derived
from Theorem \ref{I.cfot.2.1}:

\begin{example}
\label{I.cfoEx 2.1} \textbf{a.} Let $A,B$ two positive definite operators on 
$H.$ Then we have the inequalities 
\begin{equation}
1-\left\langle A^{-1}x,x\right\rangle \left\langle By,y\right\rangle \leq
\left\langle \ln Ax,x\right\rangle -\left\langle \ln By,y\right\rangle \leq
\left\langle Ax,x\right\rangle \left\langle B^{-1}y,y\right\rangle -1
\label{I.cfo2.7.1}
\end{equation}%
for any $x,y\in H$ with $\left\Vert x\right\Vert =\left\Vert y\right\Vert
=1. $

In particular, we have 
\begin{equation}
1-\left\langle A^{-1}x,x\right\rangle \left\langle Ay,y\right\rangle \leq
\left\langle \ln Ax,x\right\rangle -\left\langle \ln Ay,y\right\rangle \leq
\left\langle Ax,x\right\rangle \left\langle A^{-1}y,y\right\rangle -1
\label{I.cfo2.7.2}
\end{equation}%
for any $x,y\in H$ with $\left\Vert x\right\Vert =\left\Vert y\right\Vert =1$
and 
\begin{equation}
1-\left\langle A^{-1}x,x\right\rangle \left\langle Bx,x\right\rangle \leq
\left\langle \ln Ax,x\right\rangle -\left\langle \ln Bx,x\right\rangle \leq
\left\langle Ax,x\right\rangle \left\langle B^{-1}x,x\right\rangle -1
\label{I.cfo2.7.3}
\end{equation}%
for any $x\in H$ with $\left\Vert x\right\Vert =1.$

\textbf{b. }With the same assumption for $A$ and $B$ we have the
inequalities 
\begin{equation}
\left\langle By,y\right\rangle -\left\langle Ax,x\right\rangle \leq
\left\langle B\ln By,y\right\rangle -\left\langle \ln Ax,x\right\rangle
\left\langle By,y\right\rangle  \label{I.cfo2.7.4}
\end{equation}%
for any $x,y\in H$ with $\left\Vert x\right\Vert =\left\Vert y\right\Vert
=1. $

In particular, we have 
\begin{equation}
\left\langle Ay,y\right\rangle -\left\langle Ax,x\right\rangle \leq
\left\langle A\ln Ay,y\right\rangle -\left\langle \ln Ax,x\right\rangle
\left\langle Ay,y\right\rangle  \label{I.cfo2.7.6}
\end{equation}%
for any $x,y\in H$ with $\left\Vert x\right\Vert =\left\Vert y\right\Vert =1$
and 
\begin{equation}
\left\langle Bx,x\right\rangle -\left\langle Ax,x\right\rangle \leq
\left\langle B\ln Bx,x\right\rangle -\left\langle \ln Ax,x\right\rangle
\left\langle Bx,x\right\rangle  \label{I.cfo2.7.8}
\end{equation}%
for any $x\in H$ with $\left\Vert x\right\Vert =1.$
\end{example}

The proof of Example \textbf{a} follows from Theorem \ref{I.cfot.2.1} for
the convex function $f\left( x\right) =-\ln x$ while the proof of the second
example follows by the same theorem applied for the convex function $f\left(
x\right) =x\ln x$ and performing the required calculations. The details are
omitted.

The following result may be stated as well:

\begin{theorem}[Dragomir, 2008, \protect\cite{I.SSD7}]
\label{I.cfot.2.2} Let $I$ be an interval and $f:I\rightarrow \mathbb{R}$ be
a convex and differentiable function on \r{I} (the interior of $I)$ whose
derivative $f^{\prime }$ is continuous on \r{I}$.$ If $A$ and $B$ are
selfadjoint operators on the Hilbert space $H$ with $Sp\left( A\right)
,Sp\left( B\right) \subseteq \left[ m,M\right] \subset $\r{I}$,$ then 
\begin{align}
f^{\prime }\left( \left\langle Ax,x\right\rangle \right) \left( \left\langle
By,y\right\rangle -\left\langle Ax,x\right\rangle \right) & \leq
\left\langle f\left( B\right) y,y\right\rangle -f\left( \left\langle
Ax,x\right\rangle \right)  \label{I.cfoe.2.8} \\
& \leq \left\langle f^{\prime }\left( B\right) By,y\right\rangle
-\left\langle Ax,x\right\rangle \left\langle f^{\prime }\left( B\right)
y,y\right\rangle  \notag
\end{align}%
for any $x,y\in H$ with $\left\Vert x\right\Vert =\left\Vert y\right\Vert
=1. $

In particular, we have 
\begin{align}
f^{\prime }\left( \left\langle Ax,x\right\rangle \right) \left( \left\langle
Ay,y\right\rangle -\left\langle Ax,x\right\rangle \right) & \leq
\left\langle f\left( A\right) y,y\right\rangle -f\left( \left\langle
Ax,x\right\rangle \right)  \label{I.cfoe.2.9} \\
& \leq \left\langle f^{\prime }\left( A\right) Ay,y\right\rangle
-\left\langle Ax,x\right\rangle \left\langle f^{\prime }\left( A\right)
y,y\right\rangle  \notag
\end{align}%
for any $x,y\in H$ with $\left\Vert x\right\Vert =\left\Vert y\right\Vert =1$
and 
\begin{align}
f^{\prime }\left( \left\langle Ax,x\right\rangle \right) \left( \left\langle
Bx,x\right\rangle -\left\langle Ax,x\right\rangle \right) & \leq
\left\langle f\left( B\right) x,x\right\rangle -f\left( \left\langle
Ax,x\right\rangle \right)  \label{I.cfoe.2.10} \\
& \leq \left\langle f^{\prime }\left( B\right) Bx,x\right\rangle
-\left\langle Ax,x\right\rangle \left\langle f^{\prime }\left( B\right)
x,x\right\rangle  \notag
\end{align}%
for any $x\in H$ with $\left\Vert x\right\Vert =1.$
\end{theorem}

\begin{proof}
Since $f$ is convex and differentiable on \r{I}, then we have that 
\begin{equation}
f^{\prime }\left( s\right) \cdot \left( t-s\right) \leq f\left( t\right)
-f\left( s\right) \leq f^{\prime }\left( t\right) \cdot \left( t-s\right)
\label{I.cfoe.2.11}
\end{equation}%
for any $t,s\in \left[ m,M\right] .$

If we choose $s=\left\langle Ax,x\right\rangle \in \left[ m,M\right] ,$ with
a fix $x\in H$ with $\left\Vert x\right\Vert =1,$ then we have 
\begin{equation}
f^{\prime }\left( \left\langle Ax,x\right\rangle \right) \cdot \left(
t-\left\langle Ax,x\right\rangle \right) \leq f\left( t\right) -f\left(
\left\langle Ax,x\right\rangle \right) \leq f^{\prime }\left( t\right) \cdot
\left( t-\left\langle Ax,x\right\rangle \right)  \label{I.cfoe.2.12}
\end{equation}%
for any $t\in \left[ m,M\right] .$

Now, if we apply the property (\ref{P}) to the inequality (\ref{I.cfoe.2.12}%
) and the operator $B$, then we get 
\begin{align}
\left\langle f^{\prime }\left( \left\langle Ax,x\right\rangle \right) \cdot
\left( B-\left\langle Ax,x\right\rangle \cdot 1_{H}\right) y,y\right\rangle
& \leq \left\langle \left[ f\left( B\right) -f\left( \left\langle
Ax,x\right\rangle \right) \cdot 1_{H}\right] y,y\right\rangle
\label{I.cfoe.2.13} \\
& \leq \left\langle f^{\prime }\left( B\right) \cdot \left( B-\left\langle
Ax,x\right\rangle \cdot 1_{H}\right) y,y\right\rangle  \notag
\end{align}%
for any $x,y\in H$ with $\left\Vert x\right\Vert =\left\Vert y\right\Vert
=1, $ which is equivalent with the desired result (\ref{I.cfoe.2.8}).
\end{proof}

\begin{remark}
\label{I.cfor.2.2}We observe that if we choose $B=A$ in (\ref{I.cfoe.2.10})
or $y=x$ in (\ref{I.cfoe.2.9}) then we recapture the Mond-Pe\v{c}ari\'{c}
inequality and its reverse from (\ref{I.e.2.1}).
\end{remark}

The following particular case of interest follows from Theorem \ref%
{I.cfot.2.2}

\begin{corollary}[Dragomir, 2008, \protect\cite{I.SSD7}]
\label{I.cfoc.2.2}Assume that $f,A$ and $B$ are as in Theorem \ref%
{I.cfot.2.2}. If, either $f$ is increasing on $\left[ m,M\right] $ and $%
B\geq A$ in the operator order of $B\left( H\right) $ or $f$ is decreasing
and $B\leq A,$ then we have the Jensen's type inequality 
\begin{equation}
\left\langle f\left( B\right) x,x\right\rangle \geq f\left( \left\langle
Ax,x\right\rangle \right)  \label{I.cfoe.2.14}
\end{equation}%
for any $x\in H$ with $\left\Vert x\right\Vert =1.$
\end{corollary}

The proof is obvious by the first inequality in (\ref{I.cfoe.2.10}) and the
details are omitted.

We provide now some particular inequalities of interest that can be derived
from Theorem \ref{I.cfot.2.2}:

\begin{example}
\label{I.cfoEx 2.2} \textbf{a.} Let $A,B$ be two positive definite operators
on $H.$ Then we have the inequalities%
\begin{equation}
1-\left\langle Ax,x\right\rangle ^{-1}\left\langle By,y\right\rangle \leq
\ln \left( \left\langle Ax,x\right\rangle \right) -\left\langle \ln
By,y\right\rangle \leq \left\langle Ax,x\right\rangle \left\langle
B^{-1}y,y\right\rangle -1  \label{I.cfoe.2.15}
\end{equation}%
for any $x,y\in H$ with $\left\Vert x\right\Vert =\left\Vert y\right\Vert
=1. $

In particular, we have 
\begin{equation}
1-\left\langle Ax,x\right\rangle ^{-1}\left\langle Ay,y\right\rangle \leq
\ln \left( \left\langle Ax,x\right\rangle \right) -\left\langle \ln
Ay,y\right\rangle \leq \left\langle Ax,x\right\rangle \left\langle
A^{-1}y,y\right\rangle -1  \label{I.cfoe.2.16}
\end{equation}%
for any $x,y\in H$ with $\left\Vert x\right\Vert =\left\Vert y\right\Vert =1$
and 
\begin{equation}
1-\left\langle Ax,x\right\rangle ^{-1}\left\langle Bx,x\right\rangle \leq
\ln \left( \left\langle Ax,x\right\rangle \right) -\left\langle \ln
Bx,x\right\rangle \leq \left\langle Ax,x\right\rangle \left\langle
B^{-1}x,x\right\rangle -1  \label{I.cfoe.2.17}
\end{equation}%
for any $x\in H$ with $\left\Vert x\right\Vert =1.$

\textbf{b.} With the same assumption for $A$ and $B,$ we have the
inequalities 
\begin{equation}
\left\langle By,y\right\rangle -\left\langle Ax,x\right\rangle \leq
\left\langle B\ln By,y\right\rangle -\left\langle By,y\right\rangle \ln
\left( \left\langle Ax,x\right\rangle \right)  \label{I.cfoe.2.18}
\end{equation}%
for any $x,y\in H$ with $\left\Vert x\right\Vert =\left\Vert y\right\Vert
=1. $

In particular, we have 
\begin{equation}
\left\langle Ay,y\right\rangle -\left\langle Ax,x\right\rangle \leq
\left\langle A\ln Ay,y\right\rangle -\left\langle Ay,y\right\rangle \ln
\left( \left\langle Ax,x\right\rangle \right)  \label{I.cfoe.2.19}
\end{equation}%
for any $x,y\in H$ with $\left\Vert x\right\Vert =\left\Vert y\right\Vert =1$
and 
\begin{equation}
\left\langle Bx,x\right\rangle -\left\langle Ax,x\right\rangle \leq
\left\langle B\ln Bx,x\right\rangle -\left\langle Bx,x\right\rangle \ln
\left( \left\langle Ax,x\right\rangle \right)  \label{I.cfoe.2.20}
\end{equation}%
for any $x\in H$ with $\left\Vert x\right\Vert =1.$
\end{example}

\subsection{Inequalities for Two Sequences of Operators}

The following result may be stated:

\begin{theorem}[Dragomir, 2008, \protect\cite{I.SSD7}]
\label{I.cfot.3.1}Let $I$ be an interval and $f:I\rightarrow \mathbb{R}$ be
a convex and differentiable function on \r{I} (the interior of $I)$ whose
derivative $f^{\prime }$ is continuous on \r{I}$.$ If $A_{j}$ and $B_{j}$
are selfadjoint operators on the Hilbert space $H$ with $Sp\left(
A_{j}\right) ,Sp\left( B_{j}\right) \subseteq \left[ m,M\right] \subset $\r{I%
} for any $j\in \left\{ 1,\dots ,n\right\} ,$ then 
\begin{align}
& \sum_{j=1}^{n}\left\langle f^{\prime }\left( A_{j}\right)
x_{j},x_{j}\right\rangle \sum_{j=1}^{n}\left\langle
B_{j}y_{j},y_{j}\right\rangle -\sum_{j=1}^{n}\left\langle f^{\prime }\left(
A_{j}\right) A_{j}x_{j},x_{j}\right\rangle  \label{I.cfoe.3.1} \\
& \leq \sum_{j=1}^{n}\left\langle f\left( B_{j}\right)
y_{j},y_{j}\right\rangle -\sum_{j=1}^{n}\left\langle f\left( A_{j}\right)
x_{j},x_{j}\right\rangle  \notag \\
& \leq \sum_{j=1}^{n}\left\langle f^{\prime }\left( B_{j}\right)
B_{j}y_{j},y_{j}\right\rangle -\sum_{j=1}^{n}\left\langle
A_{j}x_{j},x_{j}\right\rangle \sum_{j=1}^{n}\left\langle f^{\prime }\left(
B_{j}\right) y_{j},y_{j}\right\rangle  \notag
\end{align}%
for any $x_{j},y_{j}\in H,$ $j\in \left\{ 1,\dots ,n\right\} $ with $%
\sum_{j=1}^{n}\left\Vert x_{j}\right\Vert ^{2}=\sum_{j=1}^{n}\left\Vert
y_{j}\right\Vert ^{2}=1.$

In particular, we have 
\begin{align}
& \sum_{j=1}^{n}\left\langle f^{\prime }\left( A_{j}\right)
x_{j},x_{j}\right\rangle \sum_{j=1}^{n}\left\langle
A_{j}y_{j},y_{j}\right\rangle -\sum_{j=1}^{n}\left\langle f^{\prime }\left(
A_{j}\right) A_{j}x_{j},x_{j}\right\rangle  \label{I.cfoe.3.2} \\
& \leq \sum_{j=1}^{n}\left\langle f\left( A_{j}\right)
y_{j},y_{j}\right\rangle -\sum_{j=1}^{n}\left\langle f\left( A_{j}\right)
x_{j},x_{j}\right\rangle  \notag \\
& \leq \sum_{j=1}^{n}\left\langle f^{\prime }\left( A_{j}\right)
A_{j}y_{j},y_{j}\right\rangle -\sum_{j=1}^{n}\left\langle
A_{j}x_{j},x_{j}\right\rangle \sum_{j=1}^{n}\left\langle f^{\prime }\left(
A_{j}\right) y_{j},y_{j}\right\rangle  \notag
\end{align}%
for any $x_{j},y_{j}\in H,$ $j\in \left\{ 1,\dots ,n\right\} $ with $%
\sum_{j=1}^{n}\left\Vert x_{j}\right\Vert ^{2}=\sum_{j=1}^{n}\left\Vert
y_{j}\right\Vert ^{2}=1$ and 
\begin{align}
& \sum_{j=1}^{n}\left\langle f^{\prime }\left( A_{j}\right)
x_{j},x_{j}\right\rangle \sum_{j=1}^{n}\left\langle
B_{j}x_{j},x_{j}\right\rangle -\sum_{j=1}^{n}\left\langle f^{\prime }\left(
A_{j}\right) A_{j}x_{j},x_{j}\right\rangle  \label{I.cfoe.3.3} \\
& \leq \sum_{j=1}^{n}\left\langle f\left( B_{j}\right)
x_{j},x_{j}\right\rangle -\sum_{j=1}^{n}\left\langle f\left( A_{j}\right)
x_{j},x_{j}\right\rangle  \notag \\
& \leq \sum_{j=1}^{n}\left\langle f^{\prime }\left( B_{j}\right)
B_{j}x_{j},x_{j}\right\rangle -\sum_{j=1}^{n}\left\langle
A_{j}x_{j},x_{j}\right\rangle \sum_{j=1}^{n}\left\langle f^{\prime }\left(
B_{j}\right) x_{j},x_{j}\right\rangle  \notag
\end{align}%
for any $x_{j}\in H,$ $j\in \left\{ 1,\dots ,n\right\} $ with $%
\sum_{j=1}^{n}\left\Vert x_{j}\right\Vert ^{2}=1.$
\end{theorem}

\begin{proof}
Follows from Theorem \ref{I.cfot.2.1} and the details are omitted.
\end{proof}

The following particular case may be of interest:

\begin{corollary}[Dragomir, 2008, \protect\cite{I.SSD7}]
\label{I.cfoc.3.1}Let $I$ be an interval and $f:I\rightarrow \mathbb{R}$ be
a convex and differentiable function on \r{I} (the interior of $I)$ whose
derivative $f^{\prime }$ is continuous on \r{I}$.$ If $A_{j}$ and $B_{j}$
are selfadjoint operators on the Hilbert space $H$ with $Sp\left(
A_{j}\right) ,Sp\left( B_{j}\right) \subseteq \left[ m,M\right] \subset $\r{I%
} for any $j\in \left\{ 1,\dots ,n\right\} ,$ then for any $p_{j},q_{j}\geq
0 $ with $\sum_{j=1}^{n}p_{j}=\sum_{j=1}^{n}q_{j}=1,$ we have the
inequalities 
\begin{align}
& \left\langle \sum_{j=1}^{n}p_{j}f^{\prime }\left( A_{j}\right)
x,x\right\rangle \left\langle \sum_{j=1}^{n}q_{j}B_{j}y,y\right\rangle
-\left\langle \sum_{j=1}^{n}p_{j}f^{\prime }\left( A_{j}\right)
A_{j}x,x\right\rangle  \label{I.cfoe.3.4} \\
& \leq \left\langle \sum_{j=1}^{n}q_{j}f\left( B_{j}\right) y,y\right\rangle
-\left\langle \sum_{j=1}^{n}p_{j}f\left( A_{j}\right) x,x\right\rangle 
\notag \\
& \leq \left\langle \sum_{j=1}^{n}q_{j}f^{\prime }\left( B_{j}\right)
B_{j}y,y\right\rangle -\left\langle \sum_{j=1}^{n}p_{j}A_{j}x,x\right\rangle
\left\langle \sum_{j=1}^{n}q_{j}f^{\prime }\left( B_{j}\right)
y,y\right\rangle  \notag
\end{align}%
for any $x,y\in H$ with $\left\Vert x\right\Vert =\left\Vert y\right\Vert
=1. $

In particular, we have 
\begin{align}
& \left\langle \sum_{j=1}^{n}p_{j}f^{\prime }\left( A_{j}\right)
x,x\right\rangle \left\langle \sum_{j=1}^{n}q_{j}A_{j}y,y\right\rangle
-\left\langle \sum_{j=1}^{n}p_{j}f^{\prime }\left( A_{j}\right)
A_{j}x,x\right\rangle  \label{I.cfoe.3.5} \\
& \leq \left\langle \sum_{j=1}^{n}q_{j}f\left( A_{j}\right) y,y\right\rangle
-\left\langle \sum_{j=1}^{n}p_{j}f\left( A_{j}\right) x,x\right\rangle 
\notag \\
& \leq \left\langle \sum_{j=1}^{n}q_{j}f^{\prime }\left( A_{j}\right)
B_{j}y,y\right\rangle -\left\langle \sum_{j=1}^{n}p_{j}A_{j}x,x\right\rangle
\left\langle \sum_{j=1}^{n}q_{j}f^{\prime }\left( A_{j}\right)
y,y\right\rangle  \notag
\end{align}%
for any $x,y\in H$ with $\left\Vert x\right\Vert =\left\Vert y\right\Vert =1$
and 
\begin{align}
& \left\langle \sum_{j=1}^{n}p_{j}f^{\prime }\left( A_{j}\right)
x,x\right\rangle \left\langle \sum_{j=1}^{n}p_{j}B_{j}x,x\right\rangle
-\left\langle \sum_{j=1}^{n}p_{j}f^{\prime }\left( A_{j}\right)
A_{j}x,x\right\rangle  \label{I.cfoe.3.6} \\
& \leq \left\langle \sum_{j=1}^{n}p_{j}f\left( B_{j}\right) x,x\right\rangle
-\left\langle \sum_{j=1}^{n}p_{j}f\left( A_{j}\right) x,x\right\rangle 
\notag \\
& \leq \left\langle \sum_{j=1}^{n}p_{j}f^{\prime }\left( B_{j}\right)
B_{j}x,x\right\rangle -\left\langle \sum_{j=1}^{n}p_{j}A_{j}x,x\right\rangle
\left\langle \sum_{j=1}^{n}p_{j}f^{\prime }\left( B_{j}\right)
x,x\right\rangle  \notag
\end{align}%
for any $x\in H$ with $\left\Vert x\right\Vert =1.$
\end{corollary}

\begin{proof}
Follows from Theorem \ref{I.cfot.3.1} on choosing $x_{j}=\sqrt{p_{j}}\cdot
x, $ $y_{j}=\sqrt{q_{j}}\cdot y,$ $j\in \left\{ 1,\dots ,n\right\} ,$ where $%
p_{j},q_{j}\geq 0,j\in \left\{ 1,\dots ,n\right\} ,$ $\sum_{j=1}^{n}p_{j}=%
\sum_{j=1}^{n}q_{j}=1$ and $x,y\in H,$ with $\left\Vert x\right\Vert
=\left\Vert y\right\Vert =1.$ The details are omitted.
\end{proof}

\begin{example}
\label{I.cfoEx 3.1}\textbf{a.} Let $A_{j},B_{j},$ $j\in \left\{ 1,\dots
,n\right\} ,$ be two sequences of positive definite operators on $H.$ Then
we have the inequalities 
\begin{align}
1-\sum_{j=1}^{n}\left\langle A_{j}^{-1}x_{j},x_{j}\right\rangle
\sum_{j=1}^{n}\left\langle B_{j}y_{j},y_{j}\right\rangle & \leq
\sum_{j=1}^{n}\left\langle \ln A_{j}x_{j},x_{j}\right\rangle
-\sum_{j=1}^{n}\left\langle \ln B_{j}y_{j},y_{j}\right\rangle
\label{I.cfoe.3.7} \\
& \leq \sum_{j=1}^{n}\left\langle A_{j}x_{j},x_{j}\right\rangle
\sum_{j=1}^{n}\left\langle B_{j}^{-1}y_{j},y_{j}\right\rangle -1  \notag
\end{align}%
for any $x_{j},y_{j}\in H,$ $j\in \left\{ 1,\dots ,n\right\} $ with $%
\sum_{j=1}^{n}\left\Vert x_{j}\right\Vert ^{2}=\sum_{j=1}^{n}\left\Vert
y_{j}\right\Vert ^{2}=1.$

\textbf{b. }With the same assumption for $A_{j}$ and $B_{j}$ we have the
inequalities 
\begin{align}
& \sum_{j=1}^{n}\left\langle B_{j}y_{j},y_{j}\right\rangle
-\sum_{j=1}^{n}\left\langle A_{j}x_{j},x_{j}\right\rangle  \label{I.cfoe.3.8}
\\
& \leq \sum_{j=1}^{n}\left\langle B_{j}\ln B_{j}y_{j},y_{j}\right\rangle
-\sum_{j=1}^{n}\left\langle \ln A_{j}x_{j},x_{j}\right\rangle
\sum_{j=1}^{n}\left\langle B_{j}y_{j},y_{j}\right\rangle  \notag
\end{align}%
for any $x_{j},y_{j}\in H,$ $j\in \left\{ 1,\dots ,n\right\} $ with $%
\sum_{j=1}^{n}\left\Vert x_{j}\right\Vert ^{2}=\sum_{j=1}^{n}\left\Vert
y_{j}\right\Vert ^{2}=1.$
\end{example}

Finally, we have

\begin{example}
\label{I.cfoEx 3.2}\textbf{a.} Let $A_{j},B_{j},$ $j\in \left\{ 1,\dots
,n\right\} ,$ be two sequences of positive definite operators on $H.$ Then
for any $p_{j},q_{j}\geq 0$ with $\sum_{j=1}^{n}p_{j}=\sum_{j=1}^{n}q_{j}=1,$
we have the inequalities 
\begin{align}
& 1-\left\langle \sum_{j=1}^{n}p_{j}A_{j}^{-1}x,x\right\rangle \left\langle
\sum_{j=1}^{n}q_{j}B_{j}y,y\right\rangle  \label{I.cfoe.3.9} \\
& \leq \left\langle \sum_{j=1}^{n}p_{j}\ln A_{j}x,x\right\rangle
-\left\langle \sum_{j=1}^{n}q_{j}\ln B_{j}y,y\right\rangle  \notag \\
& \leq \left\langle \sum_{j=1}^{n}p_{j}A_{j}x,x\right\rangle \left\langle
\sum_{j=1}^{n}q_{j}B_{j}^{-1}y,y\right\rangle -1  \notag
\end{align}%
for any $x,y\in H$ with $\left\Vert x\right\Vert =\left\Vert y\right\Vert
=1. $

\textbf{b. }With the same assumption for $A_{j}$, $B_{j},$ $p_{j}$ and $%
q_{j},$ we have the inequalities 
\begin{align}
& \left\langle \sum_{j=1}^{n}q_{j}B_{j}y,y\right\rangle -\left\langle
\sum_{j=1}^{n}p_{j}A_{j}x,x\right\rangle  \label{I.cfoe.3.10} \\
& \leq \left\langle \sum_{j=1}^{n}q_{j}B_{j}\ln B_{j}y,y\right\rangle
-\left\langle \sum_{j=1}^{n}p_{j}\ln A_{j}x,x\right\rangle \left\langle
\sum_{j=1}^{n}q_{j}B_{j}y,y\right\rangle  \notag
\end{align}%
for any $x,y\in H$ with $\left\Vert x\right\Vert =\left\Vert y\right\Vert
=1. $
\end{example}

\begin{remark}
\label{I.cfor.3.1}We observe that all the other inequalities for two
operators obtained in Subsection 3.1 can be extended for two sequences of
operators in a similar way. However, the details are left to the interested
reader.
\end{remark}

\section{Some Jensen Type Inequalities for Twice Differentiable Functions}

\subsection{Jensen's Inequality for Twice Differentiable Functions}

The following result may be stated:

\begin{theorem}[Dragomir, 2008, \protect\cite{I.SSD8}]
\label{I.jti2dt.3.1}Let $A$ be a positive definite operator on the Hilbert
space $H$ and assume that $Sp\left( A\right) \subseteq \left[ m,M\right] $
for some scalars $m,M$ with $0<$ $m<M.$ If $f$ is a twice differentiable
function on $\left( m,M\right) $ and for $p\in \left( -\infty ,0\right) \cup
\left( 1,\infty \right) $ we have for some $\gamma <\Gamma $ that 
\begin{equation}
\gamma \leq \frac{t^{2-p}}{p\left( p-1\right) }\cdot f^{\prime \prime
}\left( t\right) \leq \Gamma \text{\quad\ for any \quad }t\in \left(
m,M\right) ,  \label{I.jti2de.3.a}
\end{equation}%
then 
\begin{align}
\gamma \left( \left\langle A^{p}x,x\right\rangle -\left\langle
Ax,x\right\rangle ^{p}\right) & \leq \left\langle f\left( A\right)
x,x\right\rangle -f\left( \left\langle Ax,x\right\rangle \right)
\label{I.jti2de.3.1} \\
& \leq \Gamma \left( \left\langle A^{p}x,x\right\rangle -\left\langle
Ax,x\right\rangle ^{p}\right)  \notag
\end{align}%
for each $x\in H$ with $\left\Vert x\right\Vert =1.$

If 
\begin{equation}
\delta \leq \frac{t^{2-p}}{p\left( 1-p\right) }\cdot f^{\prime \prime
}\left( t\right) \leq \Delta \text{\quad\ for any\quad\ }t\in \left(
m,M\right)  \label{I.jti2de.3.b}
\end{equation}%
and for some $\delta <\Delta ,$ where $p\in \left( 0,1\right) $, then 
\begin{align}
\delta \left( \left\langle Ax,x\right\rangle ^{p}-\left\langle
A^{p}x,x\right\rangle \right) & \leq \left\langle f\left( A\right)
x,x\right\rangle -f\left( \left\langle Ax,x\right\rangle \right)
\label{I.jti2de.3.2} \\
& \leq \Delta \left( \left\langle Ax,x\right\rangle ^{p}-\left\langle
A^{p}x,x\right\rangle \right)  \notag
\end{align}%
for each $x\in H$ with $\left\Vert x\right\Vert =1.$
\end{theorem}

\begin{proof}
Consider the function $g_{\gamma ,p}:\left( m,M\right) \rightarrow \mathbb{R}
$ given by $g_{\gamma ,p}\left( t\right) =f\left( t\right) -\gamma t^{p}$
where $p\in \left( -\infty ,0\right) \cup \left( 1,\infty \right) .$ The
function $g_{\gamma ,p}$ is twice differentiable, 
\begin{equation*}
g_{\gamma ,p}^{\prime \prime }\left( t\right) =f^{\prime \prime }\left(
t\right) -\gamma p\left( p-1\right) t^{p-2}
\end{equation*}%
for any $t\in \left( m,M\right) $ and by (\ref{I.jti2de.3.a}) we deduce that 
$g_{\gamma ,p}$ is convex on $\left( m,M\right) .$ Now, applying the Mond \&
Pe\v{c}ari\'{c} inequality for $g_{\gamma ,p}$ we have 
\begin{align*}
0& \leq \left\langle \left( f\left( A\right) -\gamma A^{p}\right)
x,x\right\rangle -\left[ f\left( \left\langle Ax,x\right\rangle \right)
-\gamma \left\langle Ax,x\right\rangle ^{p}\right] \\
& =\left\langle f\left( A\right) x,x\right\rangle -f\left( \left\langle
Ax,x\right\rangle \right) -\gamma \left[ \left\langle A^{p}x,x\right\rangle
-\left\langle Ax,x\right\rangle ^{p}\right]
\end{align*}%
which is equivalent with the first inequality in (\ref{I.jti2de.3.1}).

By defining the function $g_{\Gamma ,p}:\left( m,M\right) \rightarrow 
\mathbb{R}$ given by $g_{\Gamma ,p}\left( t\right) =\Gamma t^{p}-f\left(
t\right) $ and applying the same argument we deduce the second part of (\ref%
{I.jti2de.3.1}).

The rest goes likewise and the details are omitted.
\end{proof}

\begin{remark}
\label{I.jti2dr.3.1}We observe that if $f$ is a twice differentiable
function on $\left( m,M\right) $ and $\varphi :=\inf_{t\in \left( m,M\right)
}f^{\prime \prime }\left( t\right) ,\Phi :=\sup_{t\in \left( m,M\right)
}f^{\prime \prime }\left( t\right) ,$ then by (\ref{I.jti2de.3.1}) we get
the inequality 
\begin{align}
\frac{1}{2}\varphi \left[ \left\langle A^{2}x,x\right\rangle -\left\langle
Ax,x\right\rangle ^{2}\right] & \leq \left\langle f\left( A\right)
x,x\right\rangle -f\left( \left\langle Ax,x\right\rangle \right)
\label{I.jti2de.3.3} \\
& \leq \frac{1}{2}\Phi \left[ \left\langle A^{2}x,x\right\rangle
-\left\langle Ax,x\right\rangle ^{2}\right]  \notag
\end{align}%
for each $x\in H$ with $\left\Vert x\right\Vert =1.$

We observe that the inequality (\ref{I.jti2de.3.3}) holds for selfadjoint
operators that are not necessarily positive.
\end{remark}

The following version for sequences of operators can be stated:

\begin{corollary}[Dragomir, 2008, \protect\cite{I.SSD7}]
\label{I.jti2dc.3.1}Let $A_{j}$ be positive definite operators with $%
Sp\left( A_{j}\right) \subseteq \left[ m,M\right] \subset \left( 0,\infty
\right) $ $\ j\in \left\{ 1,\dots ,n\right\} $. If $f$ is a twice
differentiable function on $\left( m,M\right) $ and for $p\in \left( -\infty
,0\right) \cup \left( 1,\infty \right) $ we have the condition (\ref%
{I.jti2de.3.a}), then 
\begin{align}
& \gamma \left[ \sum_{j=1}^{n}\left\langle A_{j}^{p}x_{j},x_{j}\right\rangle
-\left( \sum_{j=1}^{n}\left\langle A_{j}x_{j},x_{j}\right\rangle \right) ^{p}%
\right]  \label{I.jti2de.3.5} \\
& \leq \sum_{j=1}^{n}\left\langle f\left( A_{j}\right)
x_{j},x_{j}\right\rangle -f\left( \sum_{j=1}^{n}\left\langle
A_{j}x_{j},x_{j}\right\rangle \right)  \notag \\
& \leq \Gamma \left[ \sum_{j=1}^{n}\left\langle
A_{j}^{p}x_{j},x_{j}\right\rangle -\left( \sum_{j=1}^{n}\left\langle
A_{j}x_{j},x_{j}\right\rangle \right) ^{p}\right]  \notag
\end{align}%
for each $x_{j}\in H,j\in \left\{ 1,\dots ,n\right\} $ with $%
\sum_{j=1}^{n}\left\Vert x_{j}\right\Vert ^{2}=1$.

If we have the condition (\ref{I.jti2de.3.b}) for $p\in \left( 0,1\right) $,
then 
\begin{align}
& \delta \left[ \left( \sum_{j=1}^{n}\left\langle
A_{j}x_{j},x_{j}\right\rangle \right) ^{p}-\sum_{j=1}^{n}\left\langle
A_{j}^{p}x_{j},x_{j}\right\rangle \right]  \label{I.jti2de.3.6} \\
& \leq \sum_{j=1}^{n}\left\langle f\left( A_{j}\right)
x_{j},x_{j}\right\rangle -f\left( \sum_{j=1}^{n}\left\langle
A_{j}x_{j},x_{j}\right\rangle \right)  \notag \\
& \leq \Delta \left[ \left( \sum_{j=1}^{n}\left\langle
A_{j}x_{j},x_{j}\right\rangle \right) ^{p}-\sum_{j=1}^{n}\left\langle
A_{j}^{p}x_{j},x_{j}\right\rangle \right]  \notag
\end{align}%
for each $x_{j}\in H,j\in \left\{ 1,\dots ,n\right\} $ with $%
\sum_{j=1}^{n}\left\Vert x_{j}\right\Vert ^{2}=1$.
\end{corollary}

\begin{proof}
Follows from Theorem \ref{I.jti2dt.3.1}.
\end{proof}

\begin{corollary}[Dragomir, 2008, \protect\cite{I.SSD7}]
\label{I.jti2dc.3.2}Let $A_{j}$ be positive definite operators with $%
Sp\left( A_{j}\right) \subseteq \left[ m,M\right] \subset \left( 0,\infty
\right) $ $\ j\in \left\{ 1,\dots ,n\right\} $ and $p_{j}\geq 0,$ $j\in
\left\{ 1,\dots ,n\right\} $ with $\sum_{j=1}^{n}p_{j}=1$. If $f$ is a twice
differentiable function on $\left( m,M\right) $ and for $p\in \left( -\infty
,0\right) \cup \left( 1,\infty \right) $ we have the condition (\ref%
{I.jti2de.3.a}), then 
\begin{align}
& \gamma \left[ \left\langle \sum_{j=1}^{n}p_{j}A_{j}^{p}x,x\right\rangle
-\left\langle \sum_{j=1}^{n}p_{j}A_{j}x,x\right\rangle ^{p}\right]
\label{I.jti2de.3.7} \\
& \leq \left\langle \sum_{j=1}^{n}p_{j}f\left( A_{j}\right) x,x\right\rangle
-f\left( \left\langle \sum_{j=1}^{n}p_{j}A_{j}x,x\right\rangle \right) 
\notag \\
& \leq \Gamma \left[ \left\langle
\sum_{j=1}^{n}p_{j}A_{j}^{p}x,x\right\rangle -\left\langle
\sum_{j=1}^{n}p_{j}A_{j}x,x\right\rangle ^{p}\right]  \notag
\end{align}%
for each $x\in H$ with $\left\Vert x\right\Vert =1$.

If we have the condition (\ref{I.jti2de.3.b}) for $p\in \left( 0,1\right) $,
then 
\begin{align}
& \delta \left[ \left\langle \sum_{j=1}^{n}p_{j}A_{j}x,x\right\rangle
^{p}-\left\langle \sum_{j=1}^{n}p_{j}A_{j}^{p}x,x\right\rangle \right]
\label{I.jti2de.3.8} \\
& \leq \left\langle \sum_{j=1}^{n}p_{j}f\left( A_{j}\right) x,x\right\rangle
-f\left( \left\langle \sum_{j=1}^{n}p_{j}A_{j}x,x\right\rangle \right) 
\notag \\
& \leq \Delta \left[ \left\langle \sum_{j=1}^{n}p_{j}A_{j}x,x\right\rangle
^{p}-\left\langle \sum_{j=1}^{n}p_{j}A_{j}^{p}x,x\right\rangle \right] 
\notag
\end{align}%
for each $x\in H$ with $\left\Vert x\right\Vert =1$.
\end{corollary}

\begin{proof}
Follows from Corollary \ref{I.jti2dc.3.1} on choosing $x_{j}=\sqrt{p_{j}}%
\cdot x,$ $j\in \left\{ 1,\dots ,n\right\} ,$ where $p_{j}\geq 0,j\in
\left\{ 1,\dots ,n\right\} ,$ $\sum_{j=1}^{n}p_{j}=1$ and $x\in H,$ with $%
\left\Vert x\right\Vert =1.$ The details are omitted.
\end{proof}

\begin{remark}
\label{I.jti2dr.3.2}We observe that if $f$ is a twice differentiable
function on $\left( m,M\right) $ with $-\infty <m<M<\infty $, $Sp\left(
A_{j}\right) \subset \left[ m,M\right] ,$ $j\in \left\{ 1,\dots ,n\right\} $
and $\varphi :=\inf_{t\in \left( m,M\right) }f^{\prime \prime }\left(
t\right) ,\Phi :=\sup_{t\in \left( m,M\right) }f^{\prime \prime }\left(
t\right) ,$ then 
\begin{align}
& \varphi \left[ \sum_{j=1}^{n}\left\langle
A_{j}^{2}x_{j},x_{j}\right\rangle -\left( \sum_{j=1}^{n}\left\langle
A_{j}x_{j},x_{j}\right\rangle \right) ^{2}\right]  \label{I.jti2de.3.8.a} \\
& \leq \sum_{j=1}^{n}\left\langle f\left( A_{j}\right)
x_{j},x_{j}\right\rangle -f\left( \sum_{j=1}^{n}\left\langle
A_{j}x_{j},x_{j}\right\rangle \right)  \notag \\
& \leq \Phi \left[ \sum_{j=1}^{n}\left\langle
A_{j}^{2}x_{j},x_{j}\right\rangle -\left( \sum_{j=1}^{n}\left\langle
A_{j}x_{j},x_{j}\right\rangle \right) ^{2}\right]  \notag
\end{align}%
for each $x_{j}\in H,j\in \left\{ 1,\dots ,n\right\} $ with $%
\sum_{j=1}^{n}\left\Vert x_{j}\right\Vert ^{2}=1$.

Also, if $p_{j}\geq 0,$ $j\in \left\{ 1,\dots ,n\right\} $ with $%
\sum_{j=1}^{n}p_{j}=1,$ then 
\begin{align}
& \varphi \left[ \left\langle \sum_{j=1}^{n}p_{j}A_{j}^{2}x,x\right\rangle
-\left\langle \sum_{j=1}^{n}p_{j}A_{j}x,x\right\rangle ^{2}\right]
\label{I.jti2de.3.8.b} \\
& \leq \left\langle \sum_{j=1}^{n}p_{j}f\left( A_{j}\right) x,x\right\rangle
-f\left( \left\langle \sum_{j=1}^{n}p_{j}A_{j}x,x\right\rangle \right) 
\notag \\
& \leq \Phi \left[ \left\langle \sum_{j=1}^{n}p_{j}A_{j}^{2}x,x\right\rangle
-\left\langle \sum_{j=1}^{n}p_{j}A_{j}x,x\right\rangle ^{2}\right]  \notag
\end{align}
\end{remark}

The next result provides some inequalities for the function $f$ which
replace the cases $p=0$ and $p=1$ that were not allowed in Theorem \ref%
{I.jti2dt.3.1}:

\begin{theorem}[Dragomir, 2008, \protect\cite{I.SSD7}]
\label{I.jti2dt.3.2}Let $A$ be a positive definite operator on the Hilbert
space $H$ and assume that $Sp\left( A\right) \subseteq \left[ m,M\right] $
for some scalars $m,M$ with $0<$ $m<M.$ If $f$ is a twice differentiable
function on $\left( m,M\right) $ and we have for some $\gamma <\Gamma $ that 
\begin{equation}
\gamma \leq t^{2}\cdot f^{\prime \prime }\left( t\right) \leq \Gamma \text{%
\quad\ for any \quad }t\in \left( m,M\right) ,  \label{I.jti2de.3.9}
\end{equation}%
then 
\begin{align}
\gamma \left( \ln \left( \left\langle Ax,x\right\rangle \right)
-\left\langle \ln Ax,x\right\rangle \right) & \leq \left\langle f\left(
A\right) x,x\right\rangle -f\left( \left\langle Ax,x\right\rangle \right)
\label{I.jti2de.3.10} \\
& \leq \Gamma \left( \ln \left( \left\langle Ax,x\right\rangle \right)
-\left\langle \ln Ax,x\right\rangle \right)  \notag
\end{align}%
for each $x\in H$ with $\left\Vert x\right\Vert =1.$

If 
\begin{equation}
\delta \leq t\cdot f^{\prime \prime }\left( t\right) \leq \Delta \text{
\quad for any }\quad t\in \left( m,M\right)  \label{I.jti2de.3.11}
\end{equation}%
for some $\delta <\Delta $, then 
\begin{align}
\delta \left( \left\langle A\ln Ax,x\right\rangle -\left\langle
Ax,x\right\rangle \ln \left( \left\langle Ax,x\right\rangle \right) \right)
& \leq \left\langle f\left( A\right) x,x\right\rangle -f\left( \left\langle
Ax,x\right\rangle \right)  \label{I.jti2de.3.12} \\
& \leq \Delta \left( \left\langle A\ln Ax,x\right\rangle -\left\langle
Ax,x\right\rangle \ln \left( \left\langle Ax,x\right\rangle \right) \right) 
\notag
\end{align}%
for each $x\in H$ with $\left\Vert x\right\Vert =1.$
\end{theorem}

\begin{proof}
Consider the function $g_{\gamma ,0}:\left( m,M\right) \rightarrow \mathbb{R}
$ given by $g_{\gamma ,0}\left( t\right) =f\left( t\right) +\gamma \ln t.$
The function $g_{\gamma ,0}$ is twice differentiable, 
\begin{equation*}
g_{\gamma ,p}^{\prime \prime }\left( t\right) =f^{\prime \prime }\left(
t\right) -\gamma t^{-2}
\end{equation*}%
for any $t\in \left( m,M\right) $ and by (\ref{I.jti2de.3.9}) we deduce that 
$g_{\gamma ,0}$ is convex on $\left( m,M\right) .$ Now, applying the Mond \&
Pe\v{c}ari\'{c} inequality for $g_{\gamma ,0}$ we have 
\begin{align*}
0& \leq \left\langle \left( f\left( A\right) +\gamma \ln A\right)
x,x\right\rangle -\left[ f\left( \left\langle Ax,x\right\rangle \right)
+\gamma \ln \left( \left\langle Ax,x\right\rangle \right) \right] \\
& =\left\langle f\left( A\right) x,x\right\rangle -f\left( \left\langle
Ax,x\right\rangle \right) -\gamma \left[ \ln \left( \left\langle
Ax,x\right\rangle \right) -\left\langle \ln Ax,x\right\rangle \right]
\end{align*}%
which is equivalent with the first inequality in (\ref{I.jti2de.3.10}).

By defining the function $g_{\Gamma ,0}:\left( m,M\right) \rightarrow 
\mathbb{R}$ given by $g_{\Gamma ,0}\left( t\right) =-\Gamma \ln t-f\left(
t\right) $ and applying the same argument we deduce the second part of (\ref%
{I.jti2de.3.10}).

The rest goes likewise for the functions 
\begin{equation*}
g_{\delta ,1}\left( t\right) =f\left( t\right) -\delta t\ln t\quad \text{
and }\quad g_{\Delta ,0}\left( t\right) =\Delta t\ln t-f\left( t\right)
\end{equation*}%
and the details are omitted.
\end{proof}

\begin{corollary}[Dragomir, 2008, \protect\cite{I.SSD7}]
\label{I.jti2dc.3.3}Let $A_{j}$ be positive definite operators with $%
Sp\left( A_{j}\right) \subseteq \left[ m,M\right] \subset \left( 0,\infty
\right) $ $\ j\in \left\{ 1,\dots ,n\right\} $. If $f$ is a twice
differentiable function on $\left( m,M\right) $ and we have the condition (%
\ref{I.jti2de.3.9}), then 
\begin{align}
& \gamma \left( \ln \left( \sum_{j=1}^{n}\left\langle
A_{j}x_{j},x_{j}\right\rangle \right) -\sum_{j=1}^{n}\left\langle \ln
A_{j}x_{j},x_{j}\right\rangle \right)  \label{I.jti2de.3.13} \\
& \leq \sum_{j=1}^{n}\left\langle f\left( A_{j}\right)
x_{j},x_{j}\right\rangle -f\left( \sum_{j=1}^{n}\left\langle
A_{j}x_{j},x_{j}\right\rangle \right)  \notag \\
& \leq \Gamma \left( \ln \left( \sum_{j=1}^{n}\left\langle
A_{j}x_{j},x_{j}\right\rangle \right) -\sum_{j=1}^{n}\left\langle \ln
A_{j}x_{j},x_{j}\right\rangle \right)  \notag
\end{align}%
for each $x_{j}\in H,j\in \left\{ 1,\dots ,n\right\} $ with $%
\sum_{j=1}^{n}\left\Vert x_{j}\right\Vert ^{2}=1$.

If we have the condition (\ref{I.jti2de.3.11}), then 
\begin{align}
& \delta \left( \sum_{j=1}^{n}\left\langle A_{j}\ln
A_{j}x_{j},x_{j}\right\rangle -\sum_{j=1}^{n}\left\langle
A_{j}x_{j},x_{j}\right\rangle \ln \left( \sum_{j=1}^{n}\left\langle
A_{j}x_{j},x_{j}\right\rangle \right) \right)  \label{I.jti2de.3.14} \\
& \leq \sum_{j=1}^{n}\left\langle f\left( A_{j}\right)
x_{j},x_{j}\right\rangle -f\left( \sum_{j=1}^{n}\left\langle
A_{j}x_{j},x_{j}\right\rangle \right)  \notag \\
& \leq \Delta \left( \sum_{j=1}^{n}\left\langle A_{j}\ln
A_{j}x_{j},x_{j}\right\rangle -\sum_{j=1}^{n}\left\langle
A_{j}x_{j},x_{j}\right\rangle \ln \left( \sum_{j=1}^{n}\left\langle
A_{j}x_{j},x_{j}\right\rangle \right) \right)  \notag
\end{align}%
for each $x_{j}\in H,j\in \left\{ 1,\dots ,n\right\} $ with $%
\sum_{j=1}^{n}\left\Vert x_{j}\right\Vert ^{2}=1$.
\end{corollary}

The following particular case also holds:

\begin{corollary}[Dragomir, 2008, \protect\cite{I.SSD7}]
\label{I.jti2dc.3.4}Let $A_{j}$ be positive definite operators with $%
Sp\left( A_{j}\right) \subseteq \left[ m,M\right] \subset \left( 0,\infty
\right) $ $\ j\in \left\{ 1,\dots ,n\right\} $ and $p_{j}\geq 0,$ $j\in
\left\{ 1,\dots ,n\right\} $ with $\sum_{j=1}^{n}p_{j}=1$. If $f$ is a twice
differentiable function on $\left( m,M\right) $ and we have the condition (%
\ref{I.jti2de.3.9}), then 
\begin{align}
& \gamma \left( \ln \left( \left\langle
\sum_{j=1}^{n}p_{j}A_{j}x,x\right\rangle \right) -\left\langle
\sum_{j=1}^{n}p_{j}\ln A_{j}x,x\right\rangle \right)  \label{I.jti2de.3.15}
\\
& \leq \sum_{j=1}^{n}\left\langle f\left( A_{j}\right)
x_{j},x_{j}\right\rangle -f\left( \sum_{j=1}^{n}\left\langle
A_{j}x_{j},x_{j}\right\rangle \right)  \notag \\
& \leq \Gamma \left( \ln \left( \left\langle
\sum_{j=1}^{n}p_{j}A_{j}x,x\right\rangle \right) -\left\langle
\sum_{j=1}^{n}p_{j}\ln A_{j}x,x\right\rangle \right)  \notag
\end{align}%
for each $x\in H$ with $\left\Vert x\right\Vert =1$.

If we have the condition (\ref{I.jti2de.3.11}), then 
\begin{align}
& \delta \left( \left\langle \sum_{j=1}^{n}p_{j}A_{j}\ln
A_{j}x,x\right\rangle -\left\langle \sum_{j=1}^{n}p_{j}A_{j}x,x\right\rangle
\ln \left( \left\langle \sum_{j=1}^{n}p_{j}A_{j}x,x\right\rangle \right)
\right)  \label{I.jti2de.3.16} \\
& \leq \sum_{j=1}^{n}\left\langle f\left( A_{j}\right)
x_{j},x_{j}\right\rangle -f\left( \sum_{j=1}^{n}\left\langle
A_{j}x_{j},x_{j}\right\rangle \right)  \notag \\
& \leq \Delta \left( \left\langle \sum_{j=1}^{n}p_{j}A_{j}\ln
A_{j}x,x\right\rangle -\left\langle \sum_{j=1}^{n}p_{j}A_{j}x,x\right\rangle
\ln \left( \left\langle \sum_{j=1}^{n}p_{j}A_{j}x,x\right\rangle \right)
\right)  \notag
\end{align}%
for each $x\in H$ with $\left\Vert x\right\Vert =1$.
\end{corollary}

\subsection{Applications}

It is clear that the results from the previous section can be applied for
various particular functions which are twice differentiable and the second
derivatives satisfy the boundedness conditions from the statements of the
Theorems \ref{I.jti2dt.3.1}, \ref{I.jti2dt.3.2} and the Remark \ref%
{I.jti2dr.3.1}.

We point out here only some simple examples that are, in our opinion, of
large interest.

\textbf{1.} For a given $\alpha >0,$ consider the function $f\left( t\right)
=\exp \left( \alpha t\right) ,t\in \mathbb{R}$. We have $f^{\prime \prime
}\left( t\right) =\alpha ^{2}\exp \left( \alpha t\right) $ and for a
selfadjoint operator $A$ with $Sp\left( A\right) \subset \left[ m,M\right] $
(for some real numbers $m<M$) we also have 
\begin{equation*}
\varphi :=\inf_{t\in \left( m,M\right) }f^{\prime \prime }\left( t\right)
=\alpha ^{2}\exp \left( \alpha m\right) \text{ \quad and }\quad \Phi
:=\sup_{t\in \left( m,M\right) }f^{\prime \prime }\left( t\right) =\alpha
^{2}\exp \left( \alpha M\right) .
\end{equation*}%
Utilising the inequality (\ref{I.jti2de.3.3}) we get 
\begin{align}
\frac{1}{2}\alpha ^{2}\exp \left( \alpha m\right) \left[ \left\langle
A^{2}x,x\right\rangle \hspace{-2pt}-\hspace{-2pt}\left\langle
Ax,x\right\rangle ^{2}\right] & \hspace{-2pt}\leq \hspace{-2pt}\left\langle
\exp \left( \alpha A\right) x,x\right\rangle -\exp \left( \left\langle
\alpha Ax,x\right\rangle \right)  \label{I.jti2de.4.1} \\
& \hspace{-2pt}\leq \hspace{-2pt}\frac{1}{2}\alpha ^{2}\exp \left( \alpha
M\right) \left[ \left\langle A^{2}x,x\right\rangle \hspace{-2pt}-\hspace{-2pt%
}\left\langle Ax,x\right\rangle ^{2}\right] ,  \notag
\end{align}%
for each $x\in H$ with $\left\Vert x\right\Vert =1$.

Now, if $\beta >0$, then we also have 
\begin{align}
\frac{1}{2}\beta ^{2}\exp \hspace{-2pt}\left( -\beta M\right) \hspace{-2pt}%
\left[ \left\langle A^{2}x,x\right\rangle \hspace{-2pt}-\hspace{-2pt}%
\left\langle Ax,x\right\rangle ^{2}\right] & \hspace{-2pt}\leq \hspace{-2pt}%
\left\langle \exp \hspace{-2pt}\left( -\beta A\right) x,x\right\rangle 
\hspace{-2pt}-\hspace{-2pt}\exp \left( -\left\langle \beta Ax,x\right\rangle
\right)  \label{I.jti2de.4.2} \\
& \hspace{-2pt}\leq \hspace{-2pt}\frac{1}{2}\beta ^{2}\exp \hspace{-2pt}%
\left( \hspace{-2pt}-\hspace{-2pt}\beta m\right) \hspace{-2pt}\left[
\left\langle A^{2}x,x\right\rangle \hspace{-2pt}-\hspace{-2pt}\left\langle
Ax,x\right\rangle ^{2}\right] \hspace{-2pt},  \notag
\end{align}%
for each $x\in H$ with $\left\Vert x\right\Vert =1$.

\textbf{2. }Now, assume that $0<m<M$ and the operator $A$ satisfies the
condition $m\cdot 1_{H}\leq A\leq M\cdot 1_{H}.$ If we consider the function 
$f:\left( 0,\infty \right) \rightarrow \left( 0,\infty \right) $ defined by $%
f\left( t\right) =t^{p}$ with $p\in \left( -\infty ,0\right) \cup \left(
0,1\right) \cup \left( 1,\infty \right) .$ Then $f^{\prime \prime }\left(
t\right) =p\left( p-1\right) t^{p-2}$ and if we consider $\varphi
:=\inf_{t\in \left( m,M\right) }f^{\prime \prime }\left( t\right) $ and $%
\Phi :=\sup_{t\in \left( m,M\right) }f^{\prime \prime }\left( t\right) ,$
then we have 
\begin{equation*}
\varphi =p\left( p-1\right) m^{p-2},\Phi =p\left( p-1\right) M^{p-2}\quad 
\text{ for }p\in \lbrack 2,\infty ),
\end{equation*}%
\begin{equation*}
\varphi =p\left( p-1\right) M^{p-2},\Phi =p\left( p-1\right) m^{p-2}\quad 
\text{ for }p\in \left( 1,2\right) ,
\end{equation*}%
\begin{equation*}
\varphi =p\left( p-1\right) m^{p-2},\Phi =p\left( p-1\right) M^{p-2}\text{
\quad for }p\in \left( 0,1\right) ,
\end{equation*}%
and 
\begin{equation*}
\varphi =p\left( p-1\right) M^{p-2},\Phi =p\left( p-1\right) m^{p-2}\text{
\quad for }p\in \left( -\infty ,0\right) .
\end{equation*}%
Utilising the inequality (\ref{I.jti2de.3.3}) we then get the following
refinements an reverses of H\"{o}lder-McCarthy's inequalities: 
\begin{align}
& \frac{1}{2}p\left( p-1\right) m^{p-2}\left[ \left\langle
A^{2}x,x\right\rangle -\left\langle Ax,x\right\rangle ^{2}\right]
\label{I.jti2de.4.3} \\
& \leq \left\langle A^{p}x,x\right\rangle -\left\langle Ax,x\right\rangle
^{p}  \notag \\
& \leq \frac{1}{2}p\left( p-1\right) M^{p-2}\left[ \left\langle
A^{2}x,x\right\rangle -\left\langle Ax,x\right\rangle ^{2}\right] \quad 
\text{ for }p\in \lbrack 2,\infty ),  \notag
\end{align}%
\begin{align}
& \frac{1}{2}p\left( p-1\right) M^{p-2}\left[ \left\langle
A^{2}x,x\right\rangle -\left\langle Ax,x\right\rangle ^{2}\right]
\label{I.jti2de.4.4} \\
& \leq \left\langle A^{p}x,x\right\rangle -\left\langle Ax,x\right\rangle
^{p}  \notag \\
& \leq \frac{1}{2}p\left( p-1\right) m^{p-2}\left[ \left\langle
A^{2}x,x\right\rangle -\left\langle Ax,x\right\rangle ^{2}\right] \quad 
\text{ for }p\in \left( 1,2\right) ,  \notag
\end{align}%
\begin{align}
& \frac{1}{2}p\left( 1-p\right) M^{p-2}\left[ \left\langle
A^{2}x,x\right\rangle -\left\langle Ax,x\right\rangle ^{2}\right]
\label{I.jti2de.4.5} \\
& \leq \left\langle Ax,x\right\rangle ^{p}-\left\langle A^{p}x,x\right\rangle
\notag \\
& \leq \frac{1}{2}p\left( 1-p\right) m^{p-2}\left[ \left\langle
A^{2}x,x\right\rangle -\left\langle Ax,x\right\rangle ^{2}\right] \quad 
\text{ for }p\in \left( 0,1\right)  \notag
\end{align}%
and 
\begin{align}
& \frac{1}{2}p\left( p-1\right) M^{p-2}\left[ \left\langle
A^{2}x,x\right\rangle -\left\langle Ax,x\right\rangle ^{2}\right]
\label{I.jti2de.4.6} \\
& \leq \left\langle A^{p}x,x\right\rangle -\left\langle Ax,x\right\rangle
^{p}  \notag \\
& \leq \frac{1}{2}p\left( p-1\right) m^{p-2}\left[ \left\langle
A^{2}x,x\right\rangle -\left\langle Ax,x\right\rangle ^{2}\right] \text{
\quad for }p\in \left( -\infty ,0\right) ,  \notag
\end{align}%
for each $x\in H$ with $\left\Vert x\right\Vert =1$.

\textbf{3.} Now, if we consider the function $f:\left( 0,\infty \right)
\rightarrow \mathbb{R}$, $f\left( t\right) =-\ln t,$ then $f^{\prime \prime
}\left( t\right) =t^{-2}$ which gives that $\varphi =M^{-2}$ and $\Phi
=m^{-2}.$ Utilising the inequality (\ref{I.jti2de.3.3}) we then deduce the
bounds 
\begin{align}
& \frac{1}{2}M^{-2}\left[ \left\langle A^{2}x,x\right\rangle -\left\langle
Ax,x\right\rangle ^{2}\right]  \label{I.jti2de.4.7} \\
& \leq \ln \left( \left\langle Ax,x\right\rangle \right) -\left\langle \ln
Ax,x\right\rangle  \notag \\
& \leq \frac{1}{2}m^{-2}\left[ \left\langle A^{2}x,x\right\rangle
-\left\langle Ax,x\right\rangle ^{2}\right]  \notag
\end{align}%
for each $x\in H$ with $\left\Vert x\right\Vert =1$.

Moreover, if we consider the function $f:\left( 0,\infty \right) \rightarrow 
\mathbb{R}$, $f\left( t\right) =t\ln t,$ then $f^{\prime \prime }\left(
t\right) =t^{-1}$ which gives that $\varphi =M^{-1}$ and $\Phi =m^{-1}.$
Utilising the inequality (\ref{I.jti2de.3.3}) we then deduce the bounds 
\begin{align}
& \frac{1}{2}M^{-1}\left[ \left\langle A^{2}x,x\right\rangle -\left\langle
Ax,x\right\rangle ^{2}\right]  \label{I.jti2de.4.8} \\
& \leq \left\langle A\ln Ax,x\right\rangle -\left\langle Ax,x\right\rangle
\ln \left( \left\langle Ax,x\right\rangle \right)  \notag \\
& \leq \frac{1}{2}m^{-1}\left[ \left\langle A^{2}x,x\right\rangle
-\left\langle Ax,x\right\rangle ^{2}\right]  \notag
\end{align}%
for each $x\in H$ with $\left\Vert x\right\Vert =1.$

\begin{remark}
\label{I.jti2dr.4.1} Utilising Theorem \ref{I.jti2dt.3.1} for the particular
value of $p=-1$ we can state the inequality 
\begin{align}
\frac{1}{2}\psi \left( \left\langle A^{-1}x,x\right\rangle -\left\langle
Ax,x\right\rangle ^{-1}\right) & \leq \left\langle f\left( A\right)
x,x\right\rangle -f\left( \left\langle Ax,x\right\rangle \right)
\label{I.jti2de.4.9} \\
& \leq \frac{1}{2}\Psi \left( \left\langle A^{-1}x,x\right\rangle
-\left\langle Ax,x\right\rangle ^{-1}\right)  \notag
\end{align}%
for each $x\in H$ with $\left\Vert x\right\Vert =1,$ provided that $f$ is
twice differentiable on $\left( m,M\right) \subset \left( 0,\infty \right) $
and 
\begin{equation*}
\psi =\inf_{t\in \left( m,M\right) }t^{3}f^{\prime \prime }\left( t\right)
\quad \text{ while\quad\ }\Psi =\sup_{t\in \left( m,M\right) }t^{3}f^{\prime
\prime }\left( t\right)
\end{equation*}%
are assumed to be finite.

We observe that, by utilising the inequality (\ref{I.jti2de.4.9}) instead of
the inequality (\ref{I.jti2de.3.3}) we may obtain similar results in terms
of the quantity $\left\langle A^{-1}x,x\right\rangle -\left\langle
Ax,x\right\rangle ^{-1}$, $x\in H$ with $\left\Vert x\right\Vert =1.$
However the details are left to the interested reader.
\end{remark}

\section{Some Jensen's Type Inequalities for Log-convex Functions}

\subsection{Preliminary Results}

The following result that provides an operator version for the Jensen
inequality for convex functions is due to Mond and Pe\v{c}ari\'{c} \cite%
{I.MP} (see also \cite[p. 5]{I.FMPS}):

\begin{theorem}[Mond-Pe\v{c}ari\'{c}, 1993, \protect\cite{I.MP}]
\label{I.a.t.2.1} Let $A$ be a selfadjoint operator on the Hilbert space $H$
and assume that $Sp\left( A\right) \subseteq \left[ m,M\right] $ for some
scalars $m,M$ with $m<M.$ If $f$ is a convex function on $\left[ m,M\right]
, $ then 
\begin{equation}
f\left( \left\langle Ax,x\right\rangle \right) \leq \left\langle f\left(
A\right) x,x\right\rangle  \tag{MP}  \label{I.a.MPI}
\end{equation}%
for each $x\in H$ with $\left\Vert x\right\Vert =1.$
\end{theorem}

Taking into account the above result and its applications for various
concrete examples of convex functions, it is therefore natural to
investigate the corresponding results for the case of \textit{log-convex
functions}, namely functions $f:I\rightarrow \left( 0,\infty \right) $ for
which $\ln f$ is convex.

We observe that such functions satisfy the elementary inequality 
\begin{equation}
f\left( \left( 1-t\right) a+tb\right) \leq \left[ f\left( a\right) \right]
^{1-t}\left[ f\left( b\right) \right] ^{t}  \label{I.a.LC}
\end{equation}%
for any $a,b\in I$ and $t\in \left[ 0,1\right] .$ Also, due to the fact that
the weighted geometric mean is less than the weighted arithmetic mean, it
follows that any log-convex function is a convex functions. However,
obviously, there are functions that are convex but not log-convex.

As an immediate consequence of the Mond-Pe\v{c}ari\'{c} inequality above we
can provide the following result:

\begin{theorem}[Dragomir, 2010, \protect\cite{I.a.SSDLGC}]
\label{I.a.t.2.2}Let $A$ be a selfadjoint operator on the Hilbert space $H$
and assume that $Sp\left( A\right) \subseteq \left[ m,M\right] $ for some
scalars $m,M$ with $m<M.$ If $g:\left[ m,M\right] \rightarrow \left(
0,\infty \right) $ is log-convex, then 
\begin{equation}
g\left( \left\langle Ax,x\right\rangle \right) \leq \exp \left\langle \ln
g\left( A\right) x,x\right\rangle \leq \left\langle g\left( A\right)
x,x\right\rangle  \label{I.a.e.2.2}
\end{equation}%
for each $x\in H$ with $\left\Vert x\right\Vert =1.$
\end{theorem}

\begin{proof}
Consider the function $f:=\ln g,$ which is convex on $\left[ m,M\right] .$
Writing (\ref{I.a.MPI}) for $f$ we get $\ln \left[ g\left( \left\langle
Ax,x\right\rangle \right) \right] \leq \left\langle \ln g\left( A\right)
x,x\right\rangle ,$ for each $x\in H$ with $\left\Vert x\right\Vert =1,$
which, by taking the exponential, produces the first inequality in (\ref%
{I.a.e.2.2}).

If we also use (\ref{I.a.MPI}) for the exponential function, we get 
\begin{equation*}
\exp \left\langle \ln g\left( A\right) x,x\right\rangle \leq \left\langle
\exp \left[ \ln g\left( A\right) \right] x,x\right\rangle =\left\langle
g\left( A\right) x,x\right\rangle
\end{equation*}%
for each $x\in H$ with $\left\Vert x\right\Vert =1$ and the proof is
complete.
\end{proof}

The case of sequences of operators may be of interest and is embodied in the
following corollary:

\begin{corollary}[Dragomir, 2010, \protect\cite{I.a.SSDLGC}]
\label{I.a.c.2.1.a}Assume that $g$ is as in the Theorem \ref{I.a.t.2.2}. If $%
A_{j}$ are selfadjoint operators with $Sp\left( A_{j}\right) \subseteq \left[
m,M\right] $, $j\in \left\{ 1,...,n\right\} $ and $x_{j}\in H,j\in \left\{
1,...,n\right\} $ with $\sum_{j=1}^{n}\left\Vert x_{j}\right\Vert ^{2}=1$,
then 
\begin{equation}
g\left( \sum_{j=1}^{n}\left\langle A_{j}x_{j},x_{j}\right\rangle \right)
\leq \exp \left\langle \sum_{j=1}^{n}\ln g\left( A_{j}\right)
x_{j},x_{j}\right\rangle \leq \left\langle \sum_{j=1}^{n}g\left(
A_{j}\right) x_{j},x_{j}\right\rangle .  \label{I.a.e.2.2.a}
\end{equation}
\end{corollary}

\begin{proof}
Follows from Theorem \ref{I.a.t.2.2}and we omit the details.
\end{proof}

In particular we have:

\begin{corollary}[Dragomir, 2010, \protect\cite{I.a.SSDLGC}]
\label{I.a.c.2.1.b}Assume that $g$ is as in the Theorem \ref{I.a.t.2.2}. If $%
A_{j}$ are selfadjoint operators with $Sp\left( A_{j}\right) \subseteq \left[
m,M\right] \subset $\r{I}, $j\in \left\{ 1,...,n\right\} $ and $p_{j}\geq 0,$
$j\in \left\{ 1,...,n\right\} $ with $\sum_{j=1}^{n}p_{j}=1,$ then 
\begin{equation}
g\left( \left\langle \sum_{j=1}^{n}p_{j}A_{j}x,x\right\rangle \right) \leq
\left\langle \prod_{j=1}^{n}\left[ g\left( A_{j}\right) \right]
^{p_{j}}x,x\right\rangle \leq \left\langle \sum_{j=1}^{n}p_{j}g\left(
A_{j}\right) x,x\right\rangle  \label{I.a.e.2.2.b}
\end{equation}%
for each $x\in H$ with $\left\Vert x\right\Vert =1.$
\end{corollary}

\begin{proof}
Follows from Corollary \ref{I.a.c.2.1.a} by choosing $x_{j}=\sqrt{p_{j}}%
\cdot x,$ $j\in \left\{ 1,...,n\right\} $ where $x\in H$ with $\left\Vert
x\right\Vert =1.$
\end{proof}

It is also important to observe that, as a special case of Theorem \ref%
{I.a.t.2.1} we have the following important inequality in Operator Theory
that is well known as the H\"{o}lder-McCarthy inequality:

\begin{theorem}[H\"{o}lder-McCarthy, 1967, \protect\cite{I.Mc}]
\label{I.a.t.1.2} Let $A$ be a selfadjoint positive operator on a Hilbert
space $H$. Then

(i) \ \ $\left\langle A^{r}x,x\right\rangle \geq \left\langle
Ax,x\right\rangle ^{r}$ for all $r>1$ and $x\in H$ with $\left\| x\right\|
=1;$

(ii) \ $\left\langle A^{r}x,x\right\rangle \leq \left\langle
Ax,x\right\rangle ^{r}$ for all $0<r<1$ and $x\in H$ with $\left\| x\right\|
=1;$

(iii) If $A$ is invertible, then $\left\langle A^{-r}x,x\right\rangle \geq
\left\langle Ax,x\right\rangle ^{-r}$ for all $r>0$ and $x\in H$ with $%
\left\| x\right\| =1.$
\end{theorem}

Since the function $g\left( t\right) =t^{-r}$ for $r>0$ is log-convex, we
can improve the H\"{o}lder-McCarthy inequality as follows:

\begin{proposition}
\label{I.a.p.2.1}Let $A$ be a selfadjoint positive operator on a Hilbert
space $H.$ If $A$ is invertible, then 
\begin{equation}
\left\langle Ax,x\right\rangle ^{-r}\leq \exp \left\langle \ln \left(
A^{-r}\right) x,x\right\rangle \leq \left\langle A^{-r}x,x\right\rangle
\label{I.a.e.2.3}
\end{equation}%
for all $r>0$ and $x\in H$ with $\left\Vert x\right\Vert =1.$
\end{proposition}

The following reverse for the Mond-Pe\v{c}ari\'{c} inequality that
generalizes the scalar Lah-Ribari\'{c} inequality for convex functions is
well known, see for instance \cite[p. 57]{I.FMPS}:

\begin{theorem}
\label{I.a.t.2.3}Let $A$ be a selfadjoint operator on the Hilbert space $H$
and assume that $Sp\left( A\right) \subseteq \left[ m,M\right] $ for some
scalars $m,M$ with $m<M.$ If $f$ is a convex function on $\left[ m,M\right]
, $ then 
\begin{equation}
\left\langle f\left( A\right) x,x\right\rangle \leq \frac{M-\left\langle
Ax,x\right\rangle }{M-m}\cdot f\left( m\right) +\frac{\left\langle
Ax,x\right\rangle -m}{M-m}\cdot f\left( M\right)  \label{I.a.e.2.4}
\end{equation}%
for each $x\in H$ with $\left\Vert x\right\Vert =1.$
\end{theorem}

This result can be improved for log-convex functions as follows:

\begin{theorem}[Dragomir, 2010, \protect\cite{I.a.SSDLGC}]
\label{I.a.t.2.4}Let $A$ be a selfadjoint operator on the Hilbert space $H$
and assume that $Sp\left( A\right) \subseteq \left[ m,M\right] $ for some
scalars $m,M$ with $m<M.$ If $g:\left[ m,M\right] \rightarrow \left(
0,\infty \right) $ is log-convex, then 
\begin{align}
\left\langle g\left( A\right) x,x\right\rangle & \leq \left\langle \left[ %
\left[ g\left( m\right) \right] ^{\frac{M1_{H}-A}{M-m}}\left[ g\left(
M\right) \right] ^{\frac{A-m1_{H}}{M-m}}\right] x,x\right\rangle
\label{I.a.e.2.5} \\
& \leq \frac{M-\left\langle Ax,x\right\rangle }{M-m}\cdot g\left( m\right) +%
\frac{\left\langle Ax,x\right\rangle -m}{M-m\text{ }}\cdot g\left( M\right) 
\notag
\end{align}%
and%
\begin{align}
g\left( \left\langle Ax,x\right\rangle \right) & \leq \left[ g\left(
m\right) \right] ^{\frac{M-\left\langle Ax,x\right\rangle }{M-m}}\left[
g\left( M\right) \right] ^{\frac{\left\langle Ax,x\right\rangle -m}{M-m}}
\label{I.a.e.2.5.1} \\
& \leq \left\langle \left[ \left[ g\left( m\right) \right] ^{\frac{M1_{H}-A}{%
M-m}}\left[ g\left( M\right) \right] ^{\frac{A-m1_{H}}{M-m}}\right]
x,x\right\rangle  \notag
\end{align}%
for each $x\in H$ with $\left\Vert x\right\Vert =1.$
\end{theorem}

\begin{proof}
Observe that, by the log-convexity of $g,$ we have 
\begin{equation}
g\left( t\right) =g\left( \frac{M-t}{M-m}\cdot m+\frac{t-m}{M-m}\cdot
M\right) \leq \left[ g\left( m\right) \right] ^{\frac{M-t}{M-m}}\left[
g\left( M\right) \right] ^{\frac{t-m}{M-m}}  \label{I.a.e.2.5.2}
\end{equation}%
for any $t\in \left[ m,M\right] .$

Applying the property (\ref{P}) for the operator $A$, we have that 
\begin{equation*}
\left\langle g\left( A\right) x,x\right\rangle \leq \left\langle \Psi \left(
A\right) x,x\right\rangle
\end{equation*}%
for each $x\in H$ with $\left\Vert x\right\Vert =1,$ where $\Psi \left(
t\right) :=\left[ g\left( m\right) \right] ^{\frac{M-t}{M-m}}\left[ g\left(
M\right) \right] ^{\frac{t-m}{M-m}},$ $t\in \left[ m,M\right] .$ This proves
the first inequality in (\ref{I.a.e.2.5}).

Now, observe that, by the weighted arithmetic mean-geometric mean inequality
we have%
\begin{equation*}
\left[ g\left( m\right) \right] ^{\frac{M-t}{M-m}}\left[ g\left( M\right) %
\right] ^{\frac{t-m}{M-m}}\leq \frac{M-t}{M-m}\cdot g\left( m\right) +\frac{%
t-m}{M-m}\cdot g\left( M\right)
\end{equation*}
for any $t\in \left[ m,M\right] .$

Applying the property (\ref{P}) for the operator $A$ we deduce the second
inequality in (\ref{I.a.e.2.5}).

Further on, if we use the inequality (\ref{I.a.e.2.5.2}) for $t=\left\langle
Ax,x\right\rangle \in \left[ m,M\right] $ then we deduce the first part of (%
\ref{I.a.e.2.5.1}).

Now, observe that the function $\Psi $ introduced above can be rearranged to
read as%
\begin{equation*}
\Psi \left( t\right) =g\left( m\right) \left[ \frac{g\left( M\right) }{%
g\left( m\right) }\right] ^{\frac{t-m}{M-m}},t\in \left[ m,M\right]
\end{equation*}%
showing that $\Psi $ is a convex function on $\left[ m,M\right] .$

Applying Mond-Pe\v{c}ari\'{c}'s inequality for $\Psi $ we deduce the second
part of (\ref{I.a.e.2.5.1}) and the proof is complete.
\end{proof}

The case of sequences of operators is as follows:

\begin{corollary}[Dragomir, 2010, \protect\cite{I.a.SSDLGC}]
\label{I.a.c.2.2.a}Assume that $g$ is as in the Theorem \ref{I.a.t.2.2}. If $%
A_{j}$ are selfadjoint operators with $Sp\left( A_{j}\right) \subseteq \left[
m,M\right] $, $j\in \left\{ 1,...,n\right\} $ and $x_{j}\in H,j\in \left\{
1,...,n\right\} $ with $\sum_{j=1}^{n}\left\Vert x_{j}\right\Vert ^{2}=1$,
then 
\begin{align}
& \sum_{j=1}^{n}\left\langle g\left( A_{j}\right) x_{j},x_{j}\right\rangle
\label{I.a.e.2.5.a} \\
& \leq \left\langle \sum_{j=1}^{n}\left[ \left[ g\left( m\right) \right] ^{%
\frac{M1_{H}-A_{j}}{M-m}}\left[ g\left( M\right) \right] ^{\frac{A_{j}-m1_{H}%
}{M-m}}\right] x_{j},x_{j}\right\rangle  \notag \\
& \leq \frac{M-\sum_{j=1}^{n}\left\langle A_{j}x_{j},x_{j}\right\rangle }{M-m%
}\cdot g\left( m\right) +\frac{\sum_{j=1}^{n}\left\langle
A_{j}x_{j},x_{j}\right\rangle -m}{M-m}\cdot g\left( M\right)  \notag
\end{align}%
and%
\begin{align}
& g\left( \sum_{j=1}^{n}\left\langle A_{j}x_{j},x_{j}\right\rangle \right)
\label{I.a.e.2.5.a.1} \\
& \leq \left[ g\left( m\right) \right] ^{\frac{M-\sum_{j=1}^{n}\left\langle
A_{j}x_{j},x_{j}\right\rangle }{M-m}}\left[ g\left( M\right) \right] ^{\frac{%
\sum_{j=1}^{n}\left\langle A_{j}x_{j},x_{j}\right\rangle -m}{M-m}}  \notag \\
& \leq \left\langle \sum_{j=1}^{n}\left[ \left[ g\left( m\right) \right] ^{%
\frac{M1_{H}-A_{j}}{M-m}}\left[ g\left( M\right) \right] ^{\frac{A_{j}-m1_{H}%
}{M-m}}\right] x_{j},x_{j}\right\rangle .  \notag
\end{align}
\end{corollary}

In particular we have:

\begin{corollary}[Dragomir, 2010, \protect\cite{I.a.SSDLGC}]
\label{I.a.c.2.2.b}Assume that $g$ is as in the Theorem \ref{I.a.t.2.2}. If $%
A_{j}$ are selfadjoint operators with $Sp\left( A_{j}\right) \subseteq \left[
m,M\right] \subset $\r{I}, $j\in \left\{ 1,...,n\right\} $ and $p_{j}\geq 0,$
$j\in \left\{ 1,...,n\right\} $ with $\sum_{j=1}^{n}p_{j}=1,$ then 
\begin{align}
& \left\langle \sum_{j=1}^{n}p_{j}g\left( A_{j}\right) x,x\right\rangle
\label{I.a.e.2.5.b} \\
& \leq \left\langle \sum_{j=1}^{n}p_{j}\left[ g\left( m\right) \right] ^{%
\frac{M1_{H}-A_{j}}{M-m}}\left[ g\left( M\right) \right] ^{\frac{A_{j}-m1_{H}%
}{M-m}}x,x\right\rangle  \notag \\
& \leq \frac{M-\left\langle \sum_{j=1}^{n}p_{j}A_{j}x,x\right\rangle }{M-m}%
\cdot g\left( m\right) +\frac{\left\langle
\sum_{j=1}^{n}p_{j}A_{j}x,x\right\rangle -m}{M-m}\cdot g\left( M\right) 
\notag
\end{align}%
and%
\begin{align}
& g\left( \left\langle \sum_{j=1}^{n}p_{j}A_{j}x,x\right\rangle \right)
\label{I.a.e.2.5.c} \\
& \leq \left[ g\left( m\right) \right] ^{\frac{M-\left\langle
\sum_{j=1}^{n}p_{j}A_{j}x,x\right\rangle }{M-m}}\left[ g\left( M\right) %
\right] ^{\frac{\left\langle \sum_{j=1}^{n}p_{j}A_{j}x,x\right\rangle -m}{M-m%
}}  \notag \\
& \leq \left\langle \sum_{j=1}^{n}p_{j}\left[ g\left( m\right) \right] ^{%
\frac{M1_{H}-A_{j}}{M-m}}\left[ g\left( M\right) \right] ^{\frac{A_{j}-m1_{H}%
}{M-m}}x,x\right\rangle .  \notag
\end{align}
\end{corollary}

The above result from Theorem \ref{I.a.t.2.4} can be utilized to produce the
following reverse inequality for negative powers of operators:

\begin{proposition}
\label{I.a.p.2.2}Let $A$ be a selfadjoint positive operator on a Hilbert
space $H.$ If $A$ is invertible and $Sp\left( A\right) \subseteq \left[ m,M%
\right] \left( 0<m<M\right) ,$ then 
\begin{align}
\left\langle A^{-r}x,x\right\rangle & \leq \left\langle \left[ m^{\frac{%
M1_{H}-A}{M-m}}M^{\frac{A-m1_{H}}{M-m}}\right] ^{-r}x,x\right\rangle
\label{I.a.e.2.6} \\
& \leq \frac{M-\left\langle Ax,x\right\rangle }{M-m}\cdot m^{-r}+\frac{%
\left\langle Ax,x\right\rangle -m}{M-m}\cdot M^{-r}  \notag
\end{align}%
and%
\begin{align}
\left\langle Ax,x\right\rangle ^{-r}& \leq \left[ g\left( m\right) ^{\frac{%
M-\left\langle Ax,x\right\rangle }{M-m}}g\left( M\right) ^{\frac{%
\left\langle Ax,x\right\rangle -m}{M-m}}\right] ^{-r}  \label{I.a.e.2.6.1} \\
& \leq \left\langle \left[ m^{\frac{M1_{H}-A}{M-m}}M^{\frac{A-m1_{H}}{M-m}}%
\right] ^{-r}x,x\right\rangle  \notag
\end{align}%
for all $r>0$ and $x\in H$ with $\left\Vert x\right\Vert =1.$
\end{proposition}

\subsection{Jensen's Inequality for Differentiable Log-convex Functions}

The following result provides a reverse for the Jensen type inequality (\ref%
{I.a.MPI}):

\begin{theorem}[Dragomir, 2008, \protect\cite{I.SSD5}]
\label{I.a.t.3.1}Let $J$ be an interval and $f:J\rightarrow \mathbb{R}$ be a
convex and differentiable function on \r{J} (the interior of $J)$ whose
derivative $f^{\prime }$ is continuous on \r{J}$.$ If $A$ is a selfadjoint
operators on the Hilbert space $H$ with $Sp\left( A\right) \subseteq \left[
m,M\right] \subset \text{\r{J}},$ then 
\begin{equation}
\left( 0\leq \right) \left\langle f\left( A\right) x,x\right\rangle -f\left(
\left\langle Ax,x\right\rangle \right) \leq \left\langle f^{\prime }\left(
A\right) Ax,x\right\rangle -\left\langle Ax,x\right\rangle \cdot
\left\langle f^{\prime }\left( A\right) x,x\right\rangle  \label{I.a.e.3.1}
\end{equation}%
for any $x\in H$ with $\left\Vert x\right\Vert =1.$
\end{theorem}

The following result may be stated:

\begin{proposition}[Dragomir, 2010, \protect\cite{I.a.SSDLGC}]
\label{I.a.p.3.1}Let $J$ be an interval and $g:J\rightarrow \mathbb{R}$ be a
differentiable log-convex function on \r{J} whose derivative $g^{\prime }$
is continuous on \r{J}. If $A$ is a selfadjoint operator on the Hilbert
space $H$ with $Sp\left( A\right) \subseteq \left[ m,M\right] \subset \text{%
\r{J}},$ then 
\begin{align}
& \left( 1\leq \right) \frac{\exp \left\langle \ln g\left( A\right)
x,x\right\rangle }{g\left( \left\langle Ax,x\right\rangle \right) }
\label{I.a.e.3.2} \\
& \leq \exp \left[ \left\langle g^{\prime }\left( A\right) \left[ g\left(
A\right) \right] ^{-1}Ax,x\right\rangle -\left\langle Ax,x\right\rangle
\cdot \left\langle g^{\prime }\left( A\right) \left[ g\left( A\right) \right]
^{-1}x,x\right\rangle \right]  \notag
\end{align}%
for each $x\in H$ with $\left\Vert x\right\Vert =1.$
\end{proposition}

\begin{proof}
It follows by the inequality (\ref{I.a.e.3.1}) written for the convex
function $f=\ln g$ that 
\begin{align*}
\left\langle \ln g\left( A\right) x,x\right\rangle & \leq \ln g\left(
\left\langle Ax,x\right\rangle \right) \\
& +\left\langle g^{\prime }\left( A\right) \left[ g\left( A\right) \right]
^{-1}Ax,x\right\rangle -\left\langle Ax,x\right\rangle \cdot \left\langle
g^{\prime }\left( A\right) \left[ g\left( A\right) \right]
^{-1}x,x\right\rangle
\end{align*}%
for each $x\in H$ with $\left\Vert x\right\Vert =1.$

Now, taking the exponential and dividing by $g\left( \left\langle
Ax,x\right\rangle \right) >0$ for each $x\in H$ with $\left\Vert
x\right\Vert =1,$ we deduce the desired result (\ref{I.a.e.3.2}).
\end{proof}

\begin{corollary}[Dragomir, 2010, \protect\cite{I.a.SSDLGC}]
\label{I.a.c.3.1}Assume that $g$ is as in the Proposition \ref{I.a.p.3.1}
and $A_{j}$ are selfadjoint operators with $Sp\left( A_{j}\right) \subseteq %
\left[ m,M\right] \subset $\r{J}, $j\in \left\{ 1,...,n\right\} .$

If and $x_{j}\in H,j\in \left\{ 1,...,n\right\} $ with $\sum_{j=1}^{n}\left%
\Vert x_{j}\right\Vert ^{2}=1$, then 
\begin{align}
& \left( 1\leq \right) \frac{\exp \left\langle \sum_{j=1}^{n}\ln g\left(
A_{j}\right) x_{j},x_{j}\right\rangle }{g\left( \sum_{j=1}^{n}\left\langle
A_{j}x,x_{j}\right\rangle \right) }  \label{I.a.e.3.2.a} \\
& \leq \exp \left[ \left\langle \sum_{j=1}^{n}g^{\prime }\left( A_{j}\right) %
\left[ g\left( A_{j}\right) \right] ^{-1}A_{j}x_{j},x_{j}\right\rangle
\right.  \notag \\
& \left. -\sum_{j=1}^{n}\left\langle A_{j}x_{j},x_{j}\right\rangle \cdot
\sum_{j=1}^{n}\left\langle g^{\prime }\left( A_{j}\right) \left[ g\left(
A_{j}\right) \right] ^{-1}x_{j},x_{j}\right\rangle \right] .  \notag
\end{align}

If $p_{j}\geq 0,$ $j\in \left\{ 1,...,n\right\} $ with $%
\sum_{j=1}^{n}p_{j}=1,$ then 
\begin{align}
& \left( 1\leq \right) \frac{\left\langle \prod_{j=1}^{n}\left[ g\left(
A_{j}\right) \right] ^{p_{j}}x,x\right\rangle }{g\left( \left\langle
\sum_{j=1}^{n}p_{j}A_{j}x,x\right\rangle \right) }  \label{I.a.e.3.2.b} \\
& \leq \exp \left[ \left\langle \sum_{j=1}^{n}p_{j}g^{\prime }\left(
A_{j}\right) \left[ g\left( A_{j}\right) \right] ^{-1}A_{j}x,x\right\rangle
\right.  \notag \\
& \left. -\sum_{j=1}^{n}p_{j}\left\langle A_{j}x,x\right\rangle \cdot
\sum_{j=1}^{n}p_{j}\left\langle g^{\prime }\left( A_{j}\right) \left[
g\left( A_{j}\right) \right] ^{-1}x,x\right\rangle \right]  \notag
\end{align}%
for each $x\in H$ with $\left\Vert x\right\Vert =1.$
\end{corollary}

\begin{remark}
\label{I.a.r.3.1}Let $A$ be a selfadjoint positive operator on a Hilbert
space $H.$ If $A$ is invertible, then 
\begin{equation}
\left( 1\leq \right) \left\langle Ax,x\right\rangle ^{r}\exp \left\langle
\ln \left( A^{-r}\right) x,x\right\rangle \leq \exp \left[ r\left(
\left\langle Ax,x\right\rangle \cdot \left\langle A^{-1}x,x\right\rangle
-1\right) \right]  \label{I.a.e.3.2.c}
\end{equation}%
for all $r>0$ and $x\in H$ with $\left\Vert x\right\Vert =1.$
\end{remark}

The following result that provides both a refinement and a reverse of the
multiplicative version of Jensen's inequality can be stated as well:

\begin{theorem}[Dragomir, 2010, \protect\cite{I.a.SSDLGC}]
\label{I.a.t.3.2}Let $J$ be an interval and $g:J\rightarrow \mathbb{R}$ be a
log-convex differentiable function on \r{J} whose derivative $g^{\prime }$
is continuous on \r{J}. If $A$ is a selfadjoint operators on the Hilbert
space $H$ with $Sp\left( A\right) \subseteq \left[ m,M\right] \subset \text{%
\r{J}},$ then 
\begin{align}
1& \leq \left\langle \exp \left[ \frac{g^{\prime }\left( \left\langle
Ax,x\right\rangle \right) }{g\left( \left\langle Ax,x\right\rangle \right) }%
\left( A-\left\langle Ax,x\right\rangle 1_{H}\right) \right] x,x\right\rangle
\label{I.a.e.3.3} \\
& \leq \frac{\left\langle g\left( A\right) x,x\right\rangle }{g\left(
\left\langle Ax,x\right\rangle \right) }\leq \left\langle \exp \left[
g^{\prime }\left( A\right) \left[ g\left( A\right) \right] ^{-1}\left(
A-\left\langle Ax,x\right\rangle 1_{H}\right) \right] x,x\right\rangle 
\notag
\end{align}%
for each $x\in H$ with $\left\Vert x\right\Vert =1,$ where $1_{H}$ denotes
the identity operator on $H.$
\end{theorem}

\begin{proof}
It is well known that if $h:J\rightarrow \mathbb{R}$ is a convex
differentiable function on \r{J}, then the following \textit{gradient
inequality} holds 
\begin{equation*}
h\left( t\right) -h\left( s\right) \geq h^{\prime }\left( s\right) \left(
t-s\right)
\end{equation*}
for any $t,s\in $\r{J}.

Now, if we write this inequality for the convex function $h=\ln g,$ then we
get 
\begin{equation}
\ln g\left( t\right) -\ln g\left( s\right) \geq \frac{g^{\prime }\left(
s\right) }{g\left( s\right) }\left( t-s\right)  \label{I.a.e.3.4.0}
\end{equation}%
which is equivalent with 
\begin{equation}
g\left( t\right) \geq g\left( s\right) \exp \left[ \frac{g^{\prime }\left(
s\right) }{g\left( s\right) }\left( t-s\right) \right]  \label{I.a.e.3.4}
\end{equation}%
for any $t,s\in $\r{J}.

Further, if we take $s:=\left\langle Ax,x\right\rangle \in \left[ m,M\right]
\subset \text{\r{J}},$ for a fixed $x\in H$ with $\left\Vert x\right\Vert
=1, $ in the inequality (\ref{I.a.e.3.4}), then we get 
\begin{equation*}
g\left( t\right) \geq g\left( \left\langle Ax,x\right\rangle \right) \exp 
\left[ \frac{g^{\prime }\left( \left\langle Ax,x\right\rangle \right) }{%
g\left( \left\langle Ax,x\right\rangle \right) }\left( t-\left\langle
Ax,x\right\rangle \right) \right]
\end{equation*}%
for any $t\in $\r{J}.

Utilising the property (\ref{P}) for the operator $A$ and the Mond-Pe\v{c}ari%
\'{c} inequality for the exponential function, we can state the following
inequality that is of interest in itself as well: 
\begin{align}
\left\langle g\left( A\right) y,y\right\rangle & \geq g\left( \left\langle
Ax,x\right\rangle \right) \left\langle \exp \left[ \frac{g^{\prime }\left(
\left\langle Ax,x\right\rangle \right) }{g\left( \left\langle
Ax,x\right\rangle \right) }\left( A-\left\langle Ax,x\right\rangle
1_{H}\right) \right] y,y\right\rangle  \label{I.a.e.3.5} \\
& \geq g\left( \left\langle Ax,x\right\rangle \right) \exp \left[ \frac{%
g^{\prime }\left( \left\langle Ax,x\right\rangle \right) }{g\left(
\left\langle Ax,x\right\rangle \right) }\left( \left\langle
Ay,y\right\rangle -\left\langle Ax,x\right\rangle \right) \right]  \notag
\end{align}%
for each $x,y\in H$ with $\left\Vert x\right\Vert =\left\Vert y\right\Vert
=1.$

Further, if we put $y=x$ in (\ref{I.a.e.3.5}), then we deduce the first and
the second inequality in (\ref{I.a.e.3.3}).

Now, if we replace $s$ with $t$ in (\ref{I.a.e.3.4}) we can also write the
inequality 
\begin{equation*}
g\left( t\right) \exp \left[ \frac{g^{\prime }\left( t\right) }{g\left(
t\right) }\left( s-t\right) \right] \leq g\left( s\right)
\end{equation*}%
which is equivalent with 
\begin{equation}
g\left( t\right) \leq g\left( s\right) \exp \left[ \frac{g^{\prime }\left(
t\right) }{g\left( t\right) }\left( t-s\right) \right]  \label{I.a.e.3.7}
\end{equation}%
for any $t,s\in $\r{J}.

Further, if we take $s:=\left\langle Ax,x\right\rangle \in \left[ m,M\right]
\subset \text{\r{J}},$ for a fixed $x\in H$ with $\left\Vert x\right\Vert
=1, $ in the inequality (\ref{I.a.e.3.7}), then we get 
\begin{equation*}
g\left( t\right) \leq g\left( \left\langle Ax,x\right\rangle \right) \exp 
\left[ \frac{g^{\prime }\left( t\right) }{g\left( t\right) }\left(
t-\left\langle Ax,x\right\rangle \right) \right]
\end{equation*}%
for any $t\in $\r{J}.

Utilising the property (\ref{P}) for the operator $A$, then we can state the
following inequality that is of interest in itself as well: 
\begin{equation}
\left\langle g\left( A\right) y,y\right\rangle \leq g\left( \left\langle
Ax,x\right\rangle \right) \left\langle \exp \left[ g^{\prime }\left(
A\right) \left[ g\left( A\right) \right] ^{-1}\left( A-\left\langle
Ax,x\right\rangle 1_{H}\right) \right] y,y\right\rangle  \label{I.a.e.3.8}
\end{equation}%
for each $x,y\in H$ with $\left\Vert x\right\Vert =\left\Vert y\right\Vert
=1.$

Finally, if we put $y=x$ in (\ref{I.a.e.3.8}), then we deduce the last
inequality in (\ref{I.a.e.3.3}).
\end{proof}

The case of operator sequences is embodied in the following corollary:

\begin{corollary}[Dragomir, 2010, \protect\cite{I.a.SSDLGC}]
\label{I.a.c.3.2}Assume that $g$ is as in the Proposition \ref{I.a.p.3.1}
and $A_{j}$ are selfadjoint operators with $Sp\left( A_{j}\right) \subseteq %
\left[ m,M\right] \subset $\r{J}, $j\in \left\{ 1,...,n\right\} .$

If and $x_{j}\in H,j\in \left\{ 1,...,n\right\} $ with $\sum_{j=1}^{n}\left%
\Vert x_{j}\right\Vert ^{2}=1$, then 
\begin{align}
1& \leq \left\langle \sum_{j=1}^{n}\exp \left[ \frac{g^{\prime }\left(
\sum_{j=1}^{n}\left\langle A_{j}x_{j},x_{j}\right\rangle \right) }{g\left(
\sum_{j=1}^{n}\left\langle A_{j}x_{j},x_{j}\right\rangle \right) }\left(
A_{j}-\sum_{j=1}^{n}\left\langle A_{j}x_{j},x_{j}\right\rangle 1_{H}\right) %
\right] x_{j},x_{j}\right\rangle  \label{I.a.e.3.8.a} \\
& \leq \frac{\sum_{j=1}^{n}\left\langle g\left( A_{j}\right)
x_{j},x_{j}\right\rangle }{g\left( \sum_{j=1}^{n}\left\langle
A_{j}x_{j},x_{j}\right\rangle \right) }  \notag \\
& \leq \left\langle \sum_{j=1}^{n}\exp \left[ g^{\prime }\left( A_{j}\right) %
\left[ g\left( A_{j}\right) \right] ^{-1}\left(
A_{j}-\sum_{j=1}^{n}\left\langle A_{j}x_{j},x_{j}\right\rangle 1_{H}\right) %
\right] x_{j},x_{j}\right\rangle .  \notag
\end{align}

If $p_{j}\geq 0,$ $j\in \left\{ 1,...,n\right\} $ with $%
\sum_{j=1}^{n}p_{j}=1,$ then for each $x\in H$ with $\left\Vert x\right\Vert
=1$%
\begin{align}
1& \leq \left\langle \sum_{j=1}^{n}p_{j}\exp \left[ \frac{g^{\prime }\left(
\left\langle \sum_{j=1}^{n}p_{j}A_{j}x,x\right\rangle \right) }{g\left(
\left\langle \sum_{j=1}^{n}p_{j}A_{j}x,x\right\rangle \right) }\right.
\right.  \label{I.a.e.3.8.b} \\
& \left. \left. \times \left( A_{j}-\left\langle
\sum_{j=1}^{n}p_{j}A_{j}x,x\right\rangle 1_{H}\right) \right]
x,x\right\rangle  \notag \\
& \leq \frac{\left\langle \sum_{j=1}^{n}p_{j}g\left( A_{j}\right)
x,x\right\rangle }{g\left( \left\langle
\sum_{j=1}^{n}p_{j}A_{j}x,x\right\rangle \right) }  \notag \\
& \leq \left\langle \sum_{j=1}^{n}p_{j}\exp \left[ g^{\prime }\left(
A_{j}\right) \left[ g\left( A_{j}\right) \right] ^{-1}\left(
A_{j}-\left\langle \sum_{j=1}^{n}p_{j}A_{j}x,x\right\rangle 1_{H}\right) %
\right] x,x\right\rangle .  \notag
\end{align}
\end{corollary}

\begin{remark}
\label{I.a.r.3.2}Let $A$ be a selfadjoint positive operator on a Hilbert
space $H.$ If $A$ is invertible, then 
\begin{align}
1& \leq \left\langle \exp \left[ r\left( 1_{H}-\left\langle
Ax,x\right\rangle ^{-1}A\right) \right] x,x\right\rangle  \label{I.a.e.3.6}
\\
& \leq \left\langle A^{-r}x,x\right\rangle \left\langle Ax,x\right\rangle
^{r}\leq \left\langle \exp \left[ r\left( 1_{H}-\left\langle
Ax,x\right\rangle A^{-1}\right) \right] x,x\right\rangle  \notag
\end{align}%
for all $r>0$ and $x\in H$ with $\left\Vert x\right\Vert =1.$
\end{remark}

The following reverse inequality may be proven as well:

\begin{theorem}[Dragomir, 2010, \protect\cite{I.a.SSDLGC}]
\label{I.a.t.3.3}Let $J$ be an interval and $g:J\rightarrow \mathbb{R}$ be a
log-convex differentiable function on \r{J} whose derivative $g^{\prime }$
is continuous on \r{J}. If $A$ is a selfadjoint operators on the Hilbert
space $H$ with $Sp\left( A\right) \subseteq \left[ m,M\right] \subset \text{%
\r{J}},$ then 
\begin{align}
& \left( 1\leq \right) \frac{\left\langle \left[ g\left( M\right) \right] ^{%
\frac{A-m1_{H}}{M-m}}\left[ g\left( m\right) \right] ^{\frac{M1_{H}-A}{M-m}%
}x,x\right\rangle }{\left\langle g\left( A\right) x,x\right\rangle }
\label{I.a.e.3.6.a} \\
& \leq \frac{\left\langle g\left( A\right) \exp \left[ \frac{\left(
M1_{H}-A\right) \left( A-m1_{H}\right) }{M-m}\left( \frac{g^{\prime }\left(
M\right) }{g\left( M\right) }-\frac{g^{\prime }\left( m\right) }{g\left(
m\right) }\right) \right] x,x\right\rangle }{\left\langle g\left( A\right)
x,x\right\rangle }  \notag \\
& \leq \exp \left[ \frac{1}{4}\left( M-m\right) \left( \frac{g^{\prime
}\left( M\right) }{g\left( M\right) }-\frac{g^{\prime }\left( m\right) }{%
g\left( m\right) }\right) \right]  \notag
\end{align}%
for each $x\in H$ with $\left\Vert x\right\Vert =1.$
\end{theorem}

\begin{proof}
Utilising the inequality (\ref{I.a.e.3.4.0}) we have successively%
\begin{equation}
\frac{g\left( \left( 1-\lambda \right) t+\lambda s\right) }{g\left( s\right) 
}\geq \exp \left[ \left( 1-\lambda \right) \frac{g^{\prime }\left( s\right) 
}{g\left( s\right) }\left( t-s\right) \right]  \label{I.a.e.3.9}
\end{equation}%
and%
\begin{equation}
\frac{g\left( \left( 1-\lambda \right) t+\lambda s\right) }{g\left( t\right) 
}\geq \exp \left[ -\lambda \frac{g^{\prime }\left( t\right) }{g\left(
t\right) }\left( t-s\right) \right]  \label{I.a.e.3.10}
\end{equation}%
for any $t,s\in $\r{J} and any $\lambda \in \left[ 0,1\right] .$

Now, if we take the power $\lambda $ in the inequality (\ref{I.a.e.3.9}) and
the power $1-\lambda $ in (\ref{I.a.e.3.10}) and multiply the obtained
inequalities, we deduce%
\begin{align}
& \frac{\left[ g\left( t\right) \right] ^{1-\lambda }\left[ g\left( s\right) %
\right] ^{\lambda }}{g\left( \left( 1-\lambda \right) t+\lambda s\right) }
\label{I.a.e.3.11} \\
& \leq \exp \left[ \left( 1-\lambda \right) \lambda \left( \frac{g^{\prime
}\left( t\right) }{g\left( t\right) }-\frac{g^{\prime }\left( s\right) }{%
g\left( s\right) }\right) \left( t-s\right) \right]  \notag
\end{align}%
for any $t,s\in $\r{J} and any $\lambda \in \left[ 0,1\right] .$

Further on, if we choose in (\ref{I.a.e.3.11}) $t=M,s=m$ and $\lambda =\frac{%
M-u}{M-m},$ then, from (\ref{I.a.e.3.11}) we get the inequality%
\begin{align}
& \frac{\left[ g\left( M\right) \right] ^{\frac{u-m}{M-m}}\left[ g\left(
m\right) \right] ^{\frac{M-u}{M-m}}}{g\left( u\right) }  \label{I.a.e.3.12}
\\
& \leq \exp \left[ \frac{\left( M-u\right) \left( u-m\right) }{M-m}\left( 
\frac{g^{\prime }\left( M\right) }{g\left( M\right) }-\frac{g^{\prime
}\left( m\right) }{g\left( m\right) }\right) \right]  \notag
\end{align}%
which, together with the inequality%
\begin{equation*}
\frac{\left( M-u\right) \left( u-m\right) }{M-m}\leq \frac{1}{4}\left(
M-m\right)
\end{equation*}%
produce%
\begin{align}
& \left[ g\left( M\right) \right] ^{\frac{u-m}{M-m}}\left[ g\left( m\right) %
\right] ^{\frac{M-u}{M-m}}  \label{I.a.e.3.13} \\
& \leq g\left( u\right) \exp \left[ \frac{\left( M-u\right) \left(
u-m\right) }{M-m}\left( \frac{g^{\prime }\left( M\right) }{g\left( M\right) }%
-\frac{g^{\prime }\left( m\right) }{g\left( m\right) }\right) \right]  \notag
\\
& \leq g\left( u\right) \exp \left[ \frac{1}{4}\left( M-m\right) \left( 
\frac{g^{\prime }\left( M\right) }{g\left( M\right) }-\frac{g^{\prime
}\left( m\right) }{g\left( m\right) }\right) \right]  \notag
\end{align}%
for any $u\in \left[ m,M\right] .$

If we apply the property (\ref{P}) to the inequality (\ref{I.a.e.3.13}) and
for the operator $A$ we deduce the desired result.
\end{proof}

\begin{corollary}[Dragomir, 2010, \protect\cite{I.a.SSDLGC}]
\label{I.a.c.3.3}Assume that $g$ is as in the Theorem \ref{I.a.t.3.3} and $%
A_{j}$ are selfadjoint operators with $Sp\left( A_{j}\right) \subseteq \left[
m,M\right] \subset $\r{J}, $j\in \left\{ 1,...,n\right\} .$

If $x_{j}\in H,j\in \left\{ 1,...,n\right\} $ with $\sum_{j=1}^{n}\left\Vert
x_{j}\right\Vert ^{2}=1$, then 
\begin{align}
& \left( 1\leq \right) \frac{\sum_{j=1}^{n}\left\langle \left[ g\left(
M\right) \right] ^{\frac{A_{j}-m1_{H}}{M-m}}\left[ g\left( m\right) \right]
^{\frac{M1_{H}-A_{j}}{M-m}}x_{j},x_{j}\right\rangle }{\sum_{j=1}^{n}\left%
\langle g\left( A_{j}\right) x_{j},x_{j}\right\rangle }  \label{I.a.e.3.14}
\\
& \leq \frac{\sum_{j=1}^{n}\left\langle g\left( A_{j}\right) \exp \left[ 
\frac{\left( M1_{H}-A_{j}\right) \left( A_{j}-m1_{H}\right) }{M-m}\left( 
\frac{g^{\prime }\left( M\right) }{g\left( M\right) }-\frac{g^{\prime
}\left( m\right) }{g\left( m\right) }\right) \right] x_{j},x_{j}\right%
\rangle }{\sum_{j=1}^{n}\left\langle g\left( A_{j}\right)
x_{j},x_{j}\right\rangle }  \notag \\
& \leq \exp \left[ \frac{1}{4}\left( M-m\right) \left( \frac{g^{\prime
}\left( M\right) }{g\left( M\right) }-\frac{g^{\prime }\left( m\right) }{%
g\left( m\right) }\right) \right] .  \notag
\end{align}

If $p_{j}\geq 0,$ $j\in \left\{ 1,...,n\right\} $ with $%
\sum_{j=1}^{n}p_{j}=1,$ then for each $x\in H$ with $\left\Vert x\right\Vert
=1$%
\begin{align}
& \left( 1\leq \right) \frac{\left\langle \sum_{j=1}^{n}p_{j}\left[ g\left(
M\right) \right] ^{\frac{A_{j}-m1_{H}}{M-m}}\left[ g\left( m\right) \right]
^{\frac{M1_{H}-A_{j}}{M-m}}x,x\right\rangle }{\left\langle
\sum_{j=1}^{n}p_{j}g\left( A_{j}\right) x,x\right\rangle }
\label{I.a.e.3.15} \\
& \leq \frac{\left\langle \sum_{j=1}^{n}p_{j}g\left( A_{j}\right) \exp \left[
\frac{\left( M1_{H}-A_{j}\right) \left( A_{j}-m1_{H}\right) }{M-m}\left( 
\frac{g^{\prime }\left( M\right) }{g\left( M\right) }-\frac{g^{\prime
}\left( m\right) }{g\left( m\right) }\right) \right] x,x\right\rangle }{%
\left\langle \sum_{j=1}^{n}p_{j}g\left( A_{j}\right) x,x\right\rangle } 
\notag \\
& \leq \exp \left[ \frac{1}{4}\left( M-m\right) \left( \frac{g^{\prime
}\left( M\right) }{g\left( M\right) }-\frac{g^{\prime }\left( m\right) }{%
g\left( m\right) }\right) \right] .  \notag
\end{align}
\end{corollary}

\begin{remark}
\label{I.a.r.3.3}Let $A$ be a selfadjoint positive operator on a Hilbert
space $H.$ If $A$ is invertible and $Sp\left( A\right) \subseteq \left[ m,M%
\right] \left( 0<m<M\right) ,$ then 
\begin{align}
& \left( 1\leq \right) \frac{\left\langle \left[ g\left( M\right) \right] ^{%
\frac{r\left( m1_{H}-A\right) }{M-m}}\left[ g\left( m\right) \right] ^{\frac{%
r\left( A-M1_{H}\right) }{M-m}}x,x\right\rangle }{\left\langle
A^{-r}x,x\right\rangle }  \label{I.a.e.3.16} \\
& \leq \frac{\left\langle A^{-r}\exp \left[ \frac{r\left( M1_{H}-A\right)
\left( A-m1_{H}\right) }{Mm}\right] x,x\right\rangle }{\left\langle
A^{-r}x,x\right\rangle }\leq \exp \left[ \frac{1}{4}r\frac{\left( M-m\right)
^{2}}{mM}\right]  \notag
\end{align}
\end{remark}

\subsection{Applications for Ky Fan's Inequality}

Consider the function $g:\left( 0,1\right) \rightarrow \mathbb{R}$, $g\left(
t\right) =\left( \frac{1-t}{t}\right) ^{r},r>0.$ Observe that for the new
function $f:\left( 0,1\right) \rightarrow \mathbb{R}$, $f\left( t\right)
=\ln g\left( t\right) $ we have 
\begin{equation*}
f^{\prime }\left( t\right) =\frac{-r}{t\left( 1-t\right) }\text{ and }%
f^{\prime \prime }\left( t\right) =\frac{2r\left( \frac{1}{2}-t\right) }{%
t^{2}\left( 1-t\right) ^{2}}\text{ for }t\in \left( 0,1\right)
\end{equation*}
showing that the function $g$ is log-convex on the interval $\left( 0,\frac{1%
}{2}\right) .$

If $p_{i}>0$ for $i\in \left\{ 1,...,n\right\} $ with $\sum_{i=1}^{n}p_{i}=1$
and $t_{i}\in \left( 0,\frac{1}{2}\right) $ for $i\in \left\{
1,...,n\right\} ,$ then by applying the Jensen inequality for the convex
function $f$ (with $r=1$) on the interval $\left( 0,\frac{1}{2}\right) $ we
get 
\begin{equation}
\frac{\sum_{i=1}^{n}p_{i}t_{i}}{1-\sum_{i=1}^{n}p_{i}t_{i}}\geq
\prod_{i=1}^{n}\left( \frac{t_{i}}{1-t_{i}}\right) ^{p_{i}},
\label{I.a.e.4.1}
\end{equation}%
which is the weighted version of the celebrated \textit{Ky Fan's inequality}%
, see \cite[p. 3]{I.a.BB}.

This inequality is equivalent with 
\begin{equation*}
\prod_{i=1}^{n}\left( \frac{1-t_{i}}{t_{i}}\right) ^{p_{i}}\geq \frac{%
1-\sum_{i=1}^{n}p_{i}t_{i}}{\sum_{i=1}^{n}p_{i}t_{i}},
\end{equation*}
where $p_{i}>0$ for $i\in \left\{ 1,...,n\right\} $ with $%
\sum_{i=1}^{n}p_{i}=1$ and $t_{i}\in \left( 0,\frac{1}{2}\right) $ for $i\in
\left\{ 1,...,n\right\} .$

By the weighted arithmetic mean - geometric mean inequality we also have
that 
\begin{equation*}
\sum_{i=1}^{n}p_{i}\left( 1-t_{i}\right) t_{i}^{-1}\geq
\prod_{i=1}^{n}\left( \frac{1-t_{i}}{t_{i}}\right) ^{p_{i}}
\end{equation*}%
giving the double inequality 
\begin{equation}
\sum_{i=1}^{n}p_{i}\left( 1-t_{i}\right) t_{i}^{-1}\geq
\prod_{i=1}^{n}\left( \left( 1-t_{i}\right) t_{i}^{-1}\right) ^{p_{i}}\geq
\sum_{i=1}^{n}p_{i}\left( 1-t_{i}\right) \left(
\sum_{i=1}^{n}p_{i}t_{i}\right) ^{-1}.  \label{I.a.e.4.1.a}
\end{equation}

The following operator inequalities generalizing (\ref{I.a.e.4.1.a}) may be
stated:

\begin{proposition}
\label{I.a.p.4.1}Let $A$ be a selfadjoint positive operator on a Hilbert
space $H.$ If $A$ is invertible and $Sp\left( A\right) \subset \left( 0,%
\frac{1}{2}\right) ,$ then 
\begin{align}
\left\langle \left( A^{-1}\left( 1_{H}-A\right) \right) ^{r}x,x\right\rangle
& \geq \exp \left\langle \ln \left( A^{-1}\left( 1_{H}-A\right) \right)
^{r}x,x\right\rangle  \label{I.a.e.4.2} \\
& \geq \left( \left\langle \left( 1_{H}-A\right) x,x\right\rangle
\left\langle Ax,x\right\rangle ^{-1}\right) ^{r}  \notag
\end{align}%
for each $x\in H$ with $\left\Vert x\right\Vert =1$ and $r>0.$

In particular, 
\begin{align}
\left\langle A^{-1}\left( 1_{H}-A\right) x,x\right\rangle & \geq \exp
\left\langle \ln \left( A^{-1}\left( 1_{H}-A\right) \right) x,x\right\rangle
\label{I.a.e.4.3} \\
& \geq \left\langle \left( 1_{H}-A\right) x,x\right\rangle \left\langle
Ax,x\right\rangle ^{-1}  \notag
\end{align}%
for each $x\in H$ with $\left\Vert x\right\Vert =1$.
\end{proposition}

The proof follows by Theorem \ref{I.a.t.2.2} applied for the log-convex
function $g\left( t\right) =\left( \frac{1-t}{t}\right) ^{r},r>0,t\in \left(
0,\frac{1}{2}\right) .$

\begin{proposition}
\label{I.a.p.4.2}Let $A$ be a selfadjoint positive operator on a Hilbert
space $H.$ If $A$ is invertible and $Sp\left( A\right) \subseteq \left[ m,M%
\right] \subset \left( 0,\frac{1}{2}\right) ,$ then 
\begin{align}
& \left\langle \left( \left( 1_{H}-A\right) A^{-1}\right)
^{r}x,x\right\rangle  \label{I.a.e.4.4} \\
& \leq \left\langle \left[ \left( \frac{1-m}{m}\right) ^{\frac{r\left(
M1_{H}-A\right) }{M-m}}\left( \frac{1-M}{M}\right) ^{\frac{r\left(
A-m1_{H}\right) }{M-m}}\right] x,x\right\rangle  \notag \\
& \leq \frac{M-\left\langle Ax,x\right\rangle }{M-m}\cdot \left( \frac{1-m}{m%
}\right) ^{r}+\frac{\left\langle Ax,x\right\rangle -m}{M-m}\cdot \left( 
\frac{1-M}{M}\right) ^{r}  \notag
\end{align}%
and%
\begin{align}
& \left( \frac{1-\left\langle Ax,x\right\rangle }{\left\langle
Ax,x\right\rangle }\right) ^{r}  \label{I.a.e.4.4.1} \\
& \leq \left( \frac{1-m}{m}\right) ^{\frac{r\left( M-\left\langle
Ax,x\right\rangle \right) }{M-m}}\left( \frac{1-M}{M}\right) ^{\frac{r\left(
\left\langle Ax,x\right\rangle -m\right) }{M-m}}  \notag \\
& \leq \left\langle \left[ \left( \frac{1-m}{m}\right) ^{\frac{r\left(
M1_{H}-A\right) }{M-m}}\left( \frac{1-M}{M}\right) ^{\frac{r\left(
A-m1_{H}\right) }{M-m}}\right] x,x\right\rangle  \notag
\end{align}%
for each $x\in H$ with $\left\Vert x\right\Vert =1$ and $r>0.$
\end{proposition}

The proof follows by Theorem \ref{I.a.t.2.4} applied for the log-convex
function $g\left( t\right) =\left( \frac{1-t}{t}\right) ^{r},r>0,t\in \left(
0,\frac{1}{2}\right) .$

Finally we have:

\begin{proposition}
\label{I.a.p.4.3}Let $A$ be a selfadjoint positive operator on a Hilbert
space $H.$ If $A$ is invertible and $Sp\left( A\right) \subset \left( 0,%
\frac{1}{2}\right) ,$ then 
\begin{align}
& \left( 1\leq \right) \frac{\exp \left\langle \ln \left( \left(
1_{H}-A\right) A^{-1}\right) ^{r}x,x\right\rangle }{\left( \left(
1-\left\langle Ax,x\right\rangle \right) \left\langle Ax,x\right\rangle
^{-1}\right) ^{r}}  \label{I.a.e.4.5} \\
& \leq \exp \left[ r\left( \left\langle Ax,x\right\rangle \cdot \left\langle
A^{-1}\left( 1_{H}-A\right) ^{-1}x,x\right\rangle -\left\langle \left(
1_{H}-A\right) ^{-1}x,x\right\rangle \right) \right]  \notag
\end{align}%
and 
\begin{align}
1& \leq \left\langle \exp \left[ r\left( 1-\left\langle Ax,x\right\rangle
\right) ^{-1}\left( 1_{H}-\left\langle Ax,x\right\rangle ^{-1}A\right) %
\right] x,x\right\rangle  \label{I.a.e.4.6} \\
& \leq \frac{\left\langle \left( \left( 1_{H}-A\right) A^{-1}\right)
^{r}x,x\right\rangle }{\left( \left( 1-\left\langle Ax,x\right\rangle
\right) \left\langle Ax,x\right\rangle ^{-1}\right) ^{r}}  \notag \\
& \leq \left\langle \exp \left[ r\left( 1_{H}-A\right) ^{-1}\left(
\left\langle Ax,x\right\rangle A^{-1}-1_{H}\right) \right] x,x\right\rangle 
\notag
\end{align}

for each $x\in H$ with $\left\| x\right\| =1$ and $r>0.$
\end{proposition}

The proof follows by Proposition \ref{I.a.p.3.1} and Theorem \ref{I.a.t.3.2}
applied for the log-convex function $g\left( t\right) =\left( \frac{1-t}{t}%
\right) ^{r},r>0,t\in \left( 0,\frac{1}{2}\right) .$ The details are omitted.

\subsection{More Inequalities for Differentiable Log-convex Functions}

The following results providing companion inequalities for the Jensen
inequality for differentiable log-convex functions obtained above hold:

\begin{theorem}[Dragomir, 2010, \protect\cite{I.b.SSD5}]
\label{I.b.t.5.1} Let $A$ be a selfadjoint operator on the Hilbert space $H$
and assume that $Sp\left( A\right) \subseteq \left[ m,M\right] $ for some
scalars $m,M$ with $m<M.$ If $g:J\rightarrow \left( 0,\infty \right) $ is a
differentiable log-convex function with the derivative continuous on $%
\mathring{J}$ and $\left[ m,M\right] \subset \mathring{J}$, then 
\begin{align}
& \exp \left[ \frac{\left\langle g^{\prime }\left( A\right)
Ax,x\right\rangle }{\left\langle g\left( A\right) x,x\right\rangle }-\frac{%
\left\langle g\left( A\right) Ax,x\right\rangle }{\left\langle g\left(
A\right) x,x\right\rangle }\cdot \frac{\left\langle g^{\prime }\left(
A\right) x,x\right\rangle }{\left\langle g\left( A\right) x,x\right\rangle }%
\right]  \label{I.b.e.5.1} \\
& \geq \frac{\exp \left[ \frac{\left\langle g\left( A\right) \ln g\left(
A\right) x,x\right\rangle }{\left\langle g\left( A\right) x,x\right\rangle }%
\right] }{g\left( \frac{\left\langle g\left( A\right) Ax,x\right\rangle }{%
\left\langle g\left( A\right) x,x\right\rangle }\right) }\geq 1  \notag
\end{align}%
for each $x\in H$ with $\left\Vert x\right\Vert =1.$

If 
\begin{equation}
\frac{\left\langle g^{\prime }\left( A\right) Ax,x\right\rangle }{%
\left\langle g^{\prime }\left( A\right) x,x\right\rangle }\in \mathring{J}%
\text{ for each }x\in H\text{ with }\left\Vert x\right\Vert =1,\text{\ } 
\tag{C}  \label{I.b.C}
\end{equation}%
then 
\begin{align}
& \exp \left[ \frac{g^{\prime }\left( \frac{\left\langle g^{\prime }\left(
A\right) Ax,x\right\rangle }{\left\langle g^{\prime }\left( A\right)
x,x\right\rangle }\right) }{g\left( \frac{\left\langle g^{\prime }\left(
A\right) Ax,x\right\rangle }{\left\langle g^{\prime }\left( A\right)
x,x\right\rangle }\right) }\left( \frac{\left\langle g^{\prime }\left(
A\right) Ax,x\right\rangle }{\left\langle g^{\prime }\left( A\right)
x,x\right\rangle }-\frac{\left\langle Ag\left( A\right) x,x\right\rangle }{%
\left\langle g\left( A\right) x,x\right\rangle }\right) \right]
\label{I.b.e.5.2} \\
& \geq \frac{g\left( \frac{\left\langle g^{\prime }\left( A\right)
Ax,x\right\rangle }{\left\langle g^{\prime }\left( A\right) x,x\right\rangle 
}\right) }{\exp \left( \frac{\left\langle g\left( A\right) \ln g\left(
A\right) x,x\right\rangle }{\left\langle g\left( A\right) x,x\right\rangle }%
\right) }\geq 1,  \notag
\end{align}%
for each $x\in H$ with $\left\Vert x\right\Vert =1.$
\end{theorem}

\begin{proof}
By the gradient inequality for the convex function $\ln g$ we have%
\begin{equation}
\frac{g^{\prime }\left( t\right) }{g\left( t\right) }\left( t-s\right) \geq
\ln g\left( t\right) -\ln g\left( s\right) \geq \frac{g^{\prime }\left(
s\right) }{g\left( s\right) }\left( t-s\right)  \label{I.b.e.5.2.a}
\end{equation}%
for any $t,s\in \mathring{J}$, which by multiplication with $g\left(
t\right) >0$ is equivalent with%
\begin{equation}
g^{\prime }\left( t\right) \left( t-s\right) \geq g\left( t\right) \ln
g\left( t\right) -g\left( t\right) \ln g\left( s\right) \geq \frac{g^{\prime
}\left( s\right) }{g\left( s\right) }\left( tg\left( t\right) -sg\left(
t\right) \right)  \label{I.b.e.5.3}
\end{equation}%
for any $t,s\in \mathring{J}.$

Fix $s\in \mathring{J}$ and apply the property (\ref{P}) to get that%
\begin{align}
\left\langle g^{\prime }\left( A\right) Ax,x\right\rangle -s\left\langle
g^{\prime }\left( A\right) x,x\right\rangle & \geq \left\langle g\left(
A\right) \ln g\left( A\right) x,x\right\rangle -\left\langle g\left(
A\right) x,x\right\rangle \ln g\left( s\right)  \label{I.b.e.5.4} \\
& \geq \frac{g^{\prime }\left( s\right) }{g\left( s\right) }\left(
\left\langle Ag\left( A\right) x,x\right\rangle -s\left\langle g\left(
A\right) x,x\right\rangle \right)  \notag
\end{align}%
for any $x\in H$ with $\left\Vert x\right\Vert =1,$ which is an inequality
of interest in itself as well.

Since 
\begin{equation*}
\frac{\left\langle g\left( A\right) Ax,x\right\rangle }{\left\langle g\left(
A\right) x,x\right\rangle }\in \left[ m,M\right] \text{ for any }x\in H\text{
with }\left\Vert x\right\Vert =1
\end{equation*}%
then on choosing $s:=\frac{\left\langle g\left( A\right) Ax,x\right\rangle }{%
\left\langle g\left( A\right) x,x\right\rangle }$ in (\ref{I.b.e.5.4}) we get%
\begin{align*}
& \left\langle g^{\prime }\left( A\right) Ax,x\right\rangle -\frac{%
\left\langle g\left( A\right) Ax,x\right\rangle }{\left\langle g\left(
A\right) x,x\right\rangle }\left\langle g^{\prime }\left( A\right)
x,x\right\rangle \\
& \geq \left\langle g\left( A\right) \ln g\left( A\right) x,x\right\rangle
-\left\langle g\left( A\right) x,x\right\rangle \ln g\left( \frac{%
\left\langle g\left( A\right) Ax,x\right\rangle }{\left\langle g\left(
A\right) x,x\right\rangle }\right) \geq 0,
\end{align*}%
which, by division with $\left\langle g\left( A\right) x,x\right\rangle >0,$
produces%
\begin{align}
& \frac{\left\langle g^{\prime }\left( A\right) Ax,x\right\rangle }{%
\left\langle g\left( A\right) x,x\right\rangle }-\frac{\left\langle g\left(
A\right) Ax,x\right\rangle }{\left\langle g\left( A\right) x,x\right\rangle }%
\cdot \frac{\left\langle g^{\prime }\left( A\right) x,x\right\rangle }{%
\left\langle g\left( A\right) x,x\right\rangle }  \label{I.b.e.5.5} \\
& \geq \frac{\left\langle g\left( A\right) \ln g\left( A\right)
x,x\right\rangle }{\left\langle g\left( A\right) x,x\right\rangle }-\ln
g\left( \frac{\left\langle g\left( A\right) Ax,x\right\rangle }{\left\langle
g\left( A\right) x,x\right\rangle }\right) \geq 0  \notag
\end{align}%
for any $x\in H$ with $\left\Vert x\right\Vert =1.$

Taking the exponential in (\ref{I.b.e.5.5}) we deduce the desired inequality
(\ref{I.b.e.5.1}).

Now, assuming that the condition (\ref{I.b.C}) holds, then by choosing $s:=%
\frac{\left\langle g^{\prime }\left( A\right) Ax,x\right\rangle }{%
\left\langle g^{\prime }\left( A\right) x,x\right\rangle }$ in (\ref%
{I.b.e.5.4}) we get%
\begin{align*}
0& \geq \left\langle g\left( A\right) \ln g\left( A\right) x,x\right\rangle
-\left\langle g\left( A\right) x,x\right\rangle \ln g\left( \frac{%
\left\langle g^{\prime }\left( A\right) Ax,x\right\rangle }{\left\langle
g^{\prime }\left( A\right) x,x\right\rangle }\right) \\
& \geq \frac{g^{\prime }\left( \frac{\left\langle g^{\prime }\left( A\right)
Ax,x\right\rangle }{\left\langle g^{\prime }\left( A\right) x,x\right\rangle 
}\right) }{g\left( \frac{\left\langle g^{\prime }\left( A\right)
Ax,x\right\rangle }{\left\langle g^{\prime }\left( A\right) x,x\right\rangle 
}\right) }\left( \left\langle Ag\left( A\right) x,x\right\rangle -\frac{%
\left\langle g^{\prime }\left( A\right) Ax,x\right\rangle }{\left\langle
g^{\prime }\left( A\right) x,x\right\rangle }\left\langle g\left( A\right)
x,x\right\rangle \right)
\end{align*}%
which, by dividing with $\left\langle g\left( A\right) x,x\right\rangle >0$
and rearranging, is equivalent with%
\begin{align}
& \frac{g^{\prime }\left( \frac{\left\langle g^{\prime }\left( A\right)
Ax,x\right\rangle }{\left\langle g^{\prime }\left( A\right) x,x\right\rangle 
}\right) }{g\left( \frac{\left\langle g^{\prime }\left( A\right)
Ax,x\right\rangle }{\left\langle g^{\prime }\left( A\right) x,x\right\rangle 
}\right) }\left( \frac{\left\langle g^{\prime }\left( A\right)
Ax,x\right\rangle }{\left\langle g^{\prime }\left( A\right) x,x\right\rangle 
}-\frac{\left\langle Ag\left( A\right) x,x\right\rangle }{\left\langle
g\left( A\right) x,x\right\rangle }\right)  \label{I.b.e.5.6} \\
& \geq \ln g\left( \frac{\left\langle g^{\prime }\left( A\right)
Ax,x\right\rangle }{\left\langle g^{\prime }\left( A\right) x,x\right\rangle 
}\right) -\frac{\left\langle g\left( A\right) \ln g\left( A\right)
x,x\right\rangle }{\left\langle g\left( A\right) x,x\right\rangle }\geq 0 
\notag
\end{align}%
for any $x\in H$ with $\left\Vert x\right\Vert =1.$

Finally, on taking the exponential in (\ref{I.b.e.5.6}) we deduce the
desired inequality (\ref{I.b.e.5.2}).
\end{proof}

\begin{remark}
\label{I.b.r.5.1}We observe that a sufficient condition for (\ref{I.b.C}) to
hold is that either $g^{\prime }\left( A\right) $ or $-g^{\prime }\left(
A\right) $ is a positive definite operator on $H.$
\end{remark}

\begin{corollary}[Dragomir, 2010, \protect\cite{I.b.SSD5}]
\label{I.b.c.5.1}Assume that $A$ and $g$ are as in Theorem \ref{I.b.t.5.1}.
If the condition (\ref{I.b.C}) holds, then we have the double inequality%
\begin{equation}
\ln g\left( \frac{\left\langle g^{\prime }\left( A\right) Ax,x\right\rangle 
}{\left\langle g^{\prime }\left( A\right) x,x\right\rangle }\right) \geq 
\frac{\left\langle g\left( A\right) \ln g\left( A\right) x,x\right\rangle }{%
\left\langle g\left( A\right) x,x\right\rangle }\geq \ln g\left( \frac{%
\left\langle g\left( A\right) Ax,x\right\rangle }{\left\langle g\left(
A\right) x,x\right\rangle }\right) ,  \label{I.b.e.5.7}
\end{equation}%
for any $x\in H$ with $\left\Vert x\right\Vert =1.$
\end{corollary}

\begin{remark}
\label{I.b.r.5.2}Assume that $A$ is a positive definite operator on $H.$
Since for $r>0$ the function $g\left( t\right) =t^{-r}$ is log-convex on $%
\left( 0,\infty \right) $ and%
\begin{equation*}
\frac{\left\langle g^{\prime }\left( A\right) Ax,x\right\rangle }{%
\left\langle g^{\prime }\left( A\right) x,x\right\rangle }=\frac{%
\left\langle A^{-r}x,x\right\rangle }{\left\langle A^{-r-1}x,x\right\rangle }%
>0
\end{equation*}%
for any $x\in H$ with $\left\Vert x\right\Vert =1,$ then on applying the
inequality (\ref{I.b.e.5.7}) we deduce the following interesting result 
\begin{equation}
\ln \left( \frac{\left\langle A^{-r}x,x\right\rangle }{\left\langle
A^{-r-1}x,x\right\rangle }\right) \leq \frac{\left\langle A^{-r}\ln
Ax,x\right\rangle }{\left\langle A^{-r}x,x\right\rangle }\leq \ln \left( 
\frac{\left\langle A^{-r+1}x,x\right\rangle }{\left\langle
A^{-r}x,x\right\rangle }\right)  \label{I.b.e.5.8}
\end{equation}%
for any $x\in H$ with $\left\Vert x\right\Vert =1.$

The details of the proof are left to the interested reader.
\end{remark}

The case of sequences of operators is embodied in the following corollary:

\begin{corollary}[Dragomir, 2010, \protect\cite{I.b.SSD5}]
\label{I.b.c.5.2}Let $A_{j}$, $j\in \left\{ 1,...,n\right\} $ be selfadjoint
operators on the Hilbert space $H$ and assume that $Sp\left( A_{j}\right)
\subseteq \left[ m,M\right] $ for some scalars $m,M$ with $m<M$ and each $%
j\in \left\{ 1,...,n\right\} .$ If $g:J\rightarrow \left( 0,\infty \right) $
is a differentiable log-convex function with the derivative continuous on $%
\mathring{J}$ and $\left[ m,M\right] \subset \mathring{J}$, then 
\begin{align}
& \exp \left[ \frac{\sum_{j=1}^{n}\left\langle g^{\prime }\left(
A_{j}\right) A_{j}x_{j},x_{j}\right\rangle }{\sum_{j=1}^{n}\left\langle
g\left( A_{j}\right) x_{j},x_{j}\right\rangle }\right.  \label{I.b.e.5.8.1}
\\
& \left. -\frac{\sum_{j=1}^{n}\left\langle g\left( A_{j}\right)
A_{j}x_{j},x_{j}\right\rangle }{\sum_{j=1}^{n}\left\langle g\left(
A_{j}\right) x_{j},x_{j}\right\rangle }\cdot \frac{\sum_{j=1}^{n}\left%
\langle g^{\prime }\left( A_{j}\right) x_{j},x_{j}\right\rangle }{%
\sum_{j=1}^{n}\left\langle g\left( A_{j}\right) x_{j},x_{j}\right\rangle }%
\right]  \notag \\
& \geq \frac{\exp \left[ \frac{\sum_{j=1}^{n}\left\langle g\left(
A_{j}\right) \ln g\left( A_{j}\right) x_{j},x_{j}\right\rangle }{%
\sum_{j=1}^{n}\left\langle g\left( A_{j}\right) x_{j},x_{j}\right\rangle }%
\right] }{g\left( \frac{\sum_{j=1}^{n}\left\langle g\left( A_{j}\right)
A_{j}x_{j},x_{j}\right\rangle }{\sum_{j=1}^{n}\left\langle g\left(
A_{j}\right) x_{j},x_{j}\right\rangle }\right) }\geq 1  \notag
\end{align}%
for each $x_{j}\in H,$ $j\in \left\{ 1,...,n\right\} $ with $%
\sum_{j=1}^{n}\left\Vert x_{j}\right\Vert ^{2}=1.$

If 
\begin{equation}
\frac{\sum_{j=1}^{n}\left\langle g^{\prime }\left( A_{j}\right)
A_{j}x_{j},x_{j}\right\rangle }{\sum_{j=1}^{n}\left\langle g^{\prime }\left(
A_{j}\right) x_{j},x_{j}\right\rangle }\in \mathring{J}  \label{I.b.C1}
\end{equation}%
for each $x_{j}\in H,j\in \left\{ 1,...,n\right\} $ with $%
\sum_{j=1}^{n}\left\Vert x_{j}\right\Vert ^{2}=1,$ then 
\begin{align}
& \exp \left[ \frac{g^{\prime }\left( \frac{\sum_{j=1}^{n}\left\langle
g^{\prime }\left( A_{j}\right) A_{j}x_{j},x_{j}\right\rangle }{%
\sum_{j=1}^{n}\left\langle g^{\prime }\left( A_{j}\right)
x_{j},x_{j}\right\rangle }\right) }{g\left( \frac{\sum_{j=1}^{n}\left\langle
g^{\prime }\left( A_{j}\right) A_{j}x_{j},x_{j}\right\rangle }{%
\sum_{j=1}^{n}\left\langle g^{\prime }\left( A_{j}\right)
x_{j},x_{j}\right\rangle }\right) }\right.  \label{I.b.e.5.8.2} \\
& \left. \times \left( \frac{\sum_{j=1}^{n}\left\langle g^{\prime }\left(
A_{j}\right) A_{j}x_{j},x_{j}\right\rangle }{\sum_{j=1}^{n}\left\langle
g^{\prime }\left( A_{j}\right) x_{j},x_{j}\right\rangle }-\frac{%
\sum_{j=1}^{n}\left\langle A_{j}g\left( A_{j}\right)
x_{j},x_{j}\right\rangle }{\sum_{j=1}^{n}\left\langle g\left( A_{j}\right)
x_{j},x_{j}\right\rangle }\right) \right]  \notag \\
& \geq \frac{g\left( \frac{\sum_{j=1}^{n}\left\langle g^{\prime }\left(
A_{j}\right) A_{j}x_{j},x_{j}\right\rangle }{\sum_{j=1}^{n}\left\langle
g^{\prime }\left( A_{j}\right) x_{j},x_{j}\right\rangle }\right) }{\exp
\left( \frac{\sum_{j=1}^{n}\left\langle g\left( A_{j}\right) \ln g\left(
A_{j}\right) x_{j},x_{j}\right\rangle }{\sum_{j=1}^{n}\left\langle g\left(
A_{j}\right) x_{j},x_{j}\right\rangle }\right) }\geq 1,  \notag
\end{align}%
for each $x_{j}\in H,$ $j\in \left\{ 1,...,n\right\} $ with $%
\sum_{j=1}^{n}\left\Vert x_{j}\right\Vert ^{2}=1.$
\end{corollary}

The following particular case for sequences of operators also holds:

\begin{corollary}[Dragomir, 2010, \protect\cite{I.b.SSD5}]
\label{I.b.c.5.3}With the assumptions of Corollary \ref{I.b.c.5.2} and if $%
p_{j}\geq 0,$ $j\in \left\{ 1,...,n\right\} $ with $\sum_{j=1}^{n}p_{j}=1,$
then%
\begin{align}
& \exp \left[ \frac{\left\langle \sum_{j=1}^{n}p_{j}g^{\prime }\left(
A_{j}\right) A_{j}x,x\right\rangle }{\left\langle \sum_{j=1}^{n}p_{j}g\left(
A_{j}\right) x,x\right\rangle }\right.  \label{I.b.5.8.3} \\
& \left. -\frac{\left\langle \sum_{j=1}^{n}p_{j}g\left( A_{j}\right)
A_{j}x,x\right\rangle }{\left\langle \sum_{j=1}^{n}p_{j}g\left( A_{j}\right)
x,x\right\rangle }\cdot \frac{\left\langle \sum_{j=1}^{n}p_{j}g^{\prime
}\left( A_{j}\right) x,x\right\rangle }{\left\langle
\sum_{j=1}^{n}p_{j}g\left( A_{j}\right) x,x\right\rangle }\right]  \notag \\
& \geq \frac{\exp \left[ \frac{\left\langle \sum_{j=1}^{n}p_{j}g\left(
A_{j}\right) \ln g\left( A_{j}\right) x,x\right\rangle }{\left\langle
\sum_{j=1}^{n}p_{j}g\left( A_{j}\right) x,x\right\rangle }\right] }{g\left( 
\frac{\left\langle \sum_{j=1}^{n}p_{j}g\left( A_{j}\right)
A_{j}x,x\right\rangle }{\left\langle \sum_{j=1}^{n}p_{j}g\left( A_{j}\right)
x,x\right\rangle }\right) }\geq 1  \notag
\end{align}%
for each $x\in H,$ with $\left\Vert x\right\Vert =1.$

If 
\begin{equation}
\frac{\left\langle \sum_{j=1}^{n}p_{j}g^{\prime }\left( A_{j}\right)
A_{j}x,x\right\rangle }{\left\langle \sum_{j=1}^{n}p_{j}g^{\prime }\left(
A_{j}\right) x,x\right\rangle }\in \mathring{J}\text{ \ }  \label{I.b.5.8.4}
\end{equation}%
for each $x\in H,$ with $\left\Vert x\right\Vert =1,$ then 
\begin{align}
& \exp \left[ \frac{g^{\prime }\left( \frac{\left\langle
\sum_{j=1}^{n}p_{j}g^{\prime }\left( A_{j}\right) A_{j}x,x\right\rangle }{%
\left\langle \sum_{j=1}^{n}p_{j}g^{\prime }\left( A_{j}\right)
x,x\right\rangle }\right) }{g\left( \frac{\left\langle
\sum_{j=1}^{n}p_{j}g^{\prime }\left( A_{j}\right) A_{j}x,x\right\rangle }{%
\left\langle \sum_{j=1}^{n}p_{j}g^{\prime }\left( A_{j}\right)
x,x\right\rangle }\right) }\right.  \label{I.b.5.8.5} \\
& \left. \times \left( \frac{\left\langle \sum_{j=1}^{n}p_{j}g^{\prime
}\left( A_{j}\right) A_{j}x,x\right\rangle }{\left\langle
\sum_{j=1}^{n}p_{j}g^{\prime }\left( A_{j}\right) x,x\right\rangle }-\frac{%
\left\langle \sum_{j=1}^{n}p_{j}A_{j}g\left( A_{j}\right) x,x\right\rangle }{%
\left\langle \sum_{j=1}^{n}p_{j}g\left( A_{j}\right) x,x\right\rangle }%
\right) \right]  \notag \\
& \geq \frac{g\left( \frac{\left\langle \sum_{j=1}^{n}p_{j}g^{\prime }\left(
A_{j}\right) A_{j}x,x\right\rangle }{\left\langle
\sum_{j=1}^{n}p_{j}g^{\prime }\left( A_{j}\right) x,x\right\rangle }\right) 
}{\exp \left( \frac{\left\langle \sum_{j=1}^{n}p_{j}g\left( A_{j}\right) \ln
g\left( A_{j}\right) x,x\right\rangle }{\left\langle
\sum_{j=1}^{n}p_{j}g\left( A_{j}\right) x,x\right\rangle }\right) }\geq 1, 
\notag
\end{align}%
for each $x\in H,$ with $\left\Vert x\right\Vert =1.$
\end{corollary}

\begin{proof}
Follows from Corollary \ref{I.b.c.5.2} by choosing $x_{j}=\sqrt{p_{j}}\cdot
x,$ $j\in \left\{ 1,...,n\right\} $ where $x\in H$ with $\left\Vert
x\right\Vert =1.$
\end{proof}

The following result providing different inequalities also holds:

\begin{theorem}[Dragomir, 2010, \protect\cite{I.b.SSD5}]
\label{I.b.t.5.2}Let $A$ be a selfadjoint operator on the Hilbert space $H$
and assume that $Sp\left( A\right) \subseteq \left[ m,M\right] $ for some
scalars $m,M$ with $m<M.$ If $g:J\rightarrow \left( 0,\infty \right) $ is a
differentiable log-convex function with the derivative continuous on $%
\mathring{J}$ and $\left[ m,M\right] \subset \mathring{J}$, then 
\begin{align}
& \left\langle \exp \left[ g^{\prime }\left( A\right) \left( A-\frac{%
\left\langle g\left( A\right) Ax,x\right\rangle }{\left\langle g\left(
A\right) x,x\right\rangle }1_{H}\right) \right] x,x\right\rangle
\label{I.b.e.5.9} \\
& \geq \left\langle \left( \frac{g\left( A\right) }{g\left( \frac{%
\left\langle g\left( A\right) Ax,x\right\rangle }{\left\langle g\left(
A\right) x,x\right\rangle }\right) }\right) ^{g\left( A\right)
}x,x\right\rangle  \notag \\
& \geq \left\langle \exp \left[ \frac{g^{\prime }\left( \frac{\left\langle
g\left( A\right) Ax,x\right\rangle }{\left\langle g\left( A\right)
x,x\right\rangle }\right) }{g\left( \frac{\left\langle g\left( A\right)
Ax,x\right\rangle }{\left\langle g\left( A\right) x,x\right\rangle }\right) }%
\left( Ag\left( A\right) -\frac{\left\langle g\left( A\right)
Ax,x\right\rangle }{\left\langle g\left( A\right) x,x\right\rangle }g\left(
A\right) \right) \right] x,x\right\rangle \geq 1  \notag
\end{align}%
for each $x\in H$ with $\left\Vert x\right\Vert =1.$

If the condition (\ref{I.b.C}) from Theorem \ref{I.b.t.5.1} holds, then 
\begin{align}
& \left\langle \exp \left[ \frac{g^{\prime }\left( \frac{\left\langle
g^{\prime }\left( A\right) Ax,x\right\rangle }{\left\langle g^{\prime
}\left( A\right) x,x\right\rangle }\right) }{g\left( \frac{\left\langle
g^{\prime }\left( A\right) Ax,x\right\rangle }{\left\langle g^{\prime
}\left( A\right) x,x\right\rangle }\right) }\left( \frac{\left\langle
g^{\prime }\left( A\right) Ax,x\right\rangle }{\left\langle g^{\prime
}\left( A\right) x,x\right\rangle }g\left( A\right) -Ag\left( A\right)
\right) \right] x,x\right\rangle  \label{I.b.e.5.10} \\
& \geq \left\langle \left( g\left( \frac{\left\langle g^{\prime }\left(
A\right) Ax,x\right\rangle }{\left\langle g^{\prime }\left( A\right)
x,x\right\rangle }\right) \left[ g\left( A\right) \right] ^{-1}\right)
^{g\left( A\right) }x,x\right\rangle  \notag \\
& \geq \left\langle \exp \left[ g^{\prime }\left( A\right) \left( \frac{%
\left\langle g^{\prime }\left( A\right) Ax,x\right\rangle }{\left\langle
g^{\prime }\left( A\right) x,x\right\rangle }1_{H}-A\right) \right]
x,x\right\rangle \geq 1  \notag
\end{align}%
for each $x\in H$ with $\left\Vert x\right\Vert =1.$
\end{theorem}

\begin{proof}
By taking the exponential in (\ref{I.b.e.5.3}) we have the following
inequality%
\begin{equation}
\exp \left[ g^{\prime }\left( t\right) \left( t-s\right) \right] \geq \left( 
\frac{g\left( t\right) }{g\left( s\right) }\right) ^{g\left( t\right) }\geq
\exp \left[ \frac{g^{\prime }\left( s\right) }{g\left( s\right) }\left(
tg\left( t\right) -sg\left( t\right) \right) \right]  \label{I.b.e.5.11}
\end{equation}%
for any $t,s\in \mathring{J}.$

If we fix $s\in \mathring{J}$ and apply the property (\ref{P})\ to the
inequality (\ref{I.b.e.5.11}), we deduce 
\begin{align}
\left\langle \exp \left[ g^{\prime }\left( A\right) \left( A-s1_{H}\right) %
\right] x,x\right\rangle & \geq \left\langle \left( \frac{g\left( A\right) }{%
g\left( s\right) }\right) ^{g\left( A\right) }x,x\right\rangle
\label{I.b.e.5.12} \\
& \geq \left\langle \exp \left[ \frac{g^{\prime }\left( s\right) }{g\left(
s\right) }\left( Ag\left( A\right) -sg\left( A\right) \right) \right]
x,x\right\rangle  \notag
\end{align}%
for each $x\in H$ with $\left\Vert x\right\Vert =1,$ where $1_{H}$ is the
identity operator on $H.$

By Mond-Pe\v{c}ari\'{c}'s inequality applied for the convex function $\exp $
we also have%
\begin{align}
& \left\langle \exp \left[ \frac{g^{\prime }\left( s\right) }{g\left(
s\right) }\left( Ag\left( A\right) -sg\left( A\right) \right) \right]
x,x\right\rangle  \label{I.b.e.5.13} \\
& \geq \exp \left( \frac{g^{\prime }\left( s\right) }{g\left( s\right) }%
\left( \left\langle Ag\left( A\right) x,x\right\rangle -s\left\langle
g\left( A\right) x,x\right\rangle \right) \right)  \notag
\end{align}%
for each $s\in \mathring{J}$ and $x\in H$ with $\left\Vert x\right\Vert =1.$

Now, if we choose $s:=\frac{\left\langle g\left( A\right) Ax,x\right\rangle 
}{\left\langle g\left( A\right) x,x\right\rangle }\in \left[ m,M\right] $ in
(\ref{I.b.e.5.12}) and (\ref{I.b.e.5.13}) we deduce the desired result (\ref%
{I.b.e.5.9}).

Observe that, the inequality (\ref{I.b.e.5.11}) is equivalent with%
\begin{equation}
\exp \left[ \frac{g^{\prime }\left( s\right) }{g\left( s\right) }\left(
sg\left( t\right) -tg\left( t\right) \right) \right] \geq \left( \frac{%
g\left( s\right) }{g\left( t\right) }\right) ^{g\left( t\right) }\geq \exp %
\left[ g^{\prime }\left( t\right) \left( s-t\right) \right]
\label{I.b.e.5.14}
\end{equation}%
for any $t,s\in \mathring{J}.$

If we fix $s\in \mathring{J}$ and apply the property (\ref{P})\ to the
inequality (\ref{I.b.e.5.14}) we deduce%
\begin{align}
\left\langle \exp \left[ \frac{g^{\prime }\left( s\right) }{g\left( s\right) 
}\left( sg\left( A\right) -Ag\left( A\right) \right) \right]
x,x\right\rangle & \geq \left\langle \left( g\left( s\right) \left[ g\left(
A\right) \right] ^{-1}\right) ^{g\left( A\right) }x,x\right\rangle
\label{I.b.e.5.15} \\
& \geq \left\langle \exp \left[ g^{\prime }\left( A\right) \left(
s1_{H}-A\right) \right] x,x\right\rangle  \notag
\end{align}%
for each $x\in H$ with $\left\Vert x\right\Vert =1.$

By Mond-Pe\v{c}ari\'{c}'s inequality we also have%
\begin{equation}
\left\langle \exp \left[ g^{\prime }\left( A\right) \left( s1_{H}-A\right) %
\right] x,x\right\rangle \geq \exp \left[ s\left\langle g^{\prime }\left(
A\right) x,x\right\rangle -\left\langle g^{\prime }\left( A\right)
Ax,x\right\rangle \right]  \label{I.b.e.5.16}
\end{equation}%
for each $s\in \mathring{J}$ and $x\in H$ with $\left\Vert x\right\Vert =1.$

Taking into account that the condition (\ref{I.b.C}) is valid, then we can
choose in (\ref{I.b.e.5.15}) and (\ref{I.b.e.5.16}) $s:=\frac{\left\langle
g^{\prime }\left( A\right) Ax,x\right\rangle }{\left\langle g^{\prime
}\left( A\right) x,x\right\rangle }$ to get the desired result (\ref%
{I.b.e.5.10}).
\end{proof}

\begin{remark}
\label{r.3.4}If we apply, for instance, the inequality (\ref{I.b.e.5.9}) for
the log-convex function $g\left( t\right) =t^{-1},t>0,$ then, after simple
calculations, we get the inequality%
\begin{align}
\left\langle \exp \left( \frac{A^{-2}-\left\langle A^{-1}x,x\right\rangle
A^{-1}}{A^{-2}-\left\langle A^{-1}x,x\right\rangle }\right) x,x\right\rangle
& \geq \left\langle \left( \left\langle A^{-1}x,x\right\rangle A^{-1}\right)
^{A^{-1}}x,x\right\rangle  \label{I.b.e.5.17} \\
& \geq \left\langle \exp \left( \frac{A^{-1}-\left\langle
A^{-1}x,x\right\rangle 1_{H}}{\left\langle A^{-1}x,x\right\rangle ^{2}}%
\right) x,x\right\rangle  \notag \\
& \geq 1  \notag
\end{align}%
for each $x\in H$ with $\left\Vert x\right\Vert =1.$

Other similar results can be obtained from the inequality (\ref{I.b.e.5.10}%
), however the details are left to the interested reader.
\end{remark}

\subsection{A Reverse Inequality}

The following reverse inequality is also of interest:

\begin{theorem}[Dragomir, 2010, \protect\cite{I.b.SSD5}]
\label{I.b.t.3.3}Let $A$ be a selfadjoint operator on the Hilbert space $H$
and assume that $Sp\left( A\right) \subseteq \left[ m,M\right] $ for some
scalars $m,M$ with $m<M.$ If $g:J\rightarrow \left( 0,\infty \right) $ is a
differentiable log-convex function with the derivative continuous on $%
\mathring{J}$ and $\left[ m,M\right] \subset \mathring{J}$, then%
\begin{align}
& \left( 1\leq \right) \frac{\left[ g\left( m\right) \right] ^{\frac{%
M-\left\langle Ax,x\right\rangle }{M-m}}\left[ g\left( M\right) \right] ^{%
\frac{\left\langle Ax,x\right\rangle -m}{M-m}}}{\exp \left\langle \ln
g\left( A\right) x,x\right\rangle }  \label{I.b.e.3.9.0} \\
& \leq \exp \left[ \frac{\left\langle \left( M1_{H}-A\right) \left(
A-m1_{H}\right) x,x\right\rangle }{M-m}\left( \frac{g^{\prime }\left(
M\right) }{g\left( M\right) }-\frac{g^{\prime }\left( m\right) }{g\left(
m\right) }\right) \right]  \notag \\
& \leq \exp \left[ \frac{\left( M-\left\langle Ax,x\right\rangle \right)
\left( \left\langle Ax,x\right\rangle -m\right) }{M-m}\left( \frac{g^{\prime
}\left( M\right) }{g\left( M\right) }-\frac{g^{\prime }\left( m\right) }{%
g\left( m\right) }\right) \right]  \notag \\
& \leq \exp \left[ \frac{1}{4}\left( M-m\right) \left( \frac{g^{\prime
}\left( M\right) }{g\left( M\right) }-\frac{g^{\prime }\left( m\right) }{%
g\left( m\right) }\right) \right]  \notag
\end{align}%
for each $x\in H$ with $\left\Vert x\right\Vert =1.$
\end{theorem}

\begin{proof}
Utilising the inequality (\ref{I.b.e.5.2.a}) we have successively%
\begin{equation}
\ln g\left( \left( 1-\lambda \right) t+\lambda s\right) -\ln g\left(
s\right) \geq \left( 1-\lambda \right) \frac{g^{\prime }\left( s\right) }{%
g\left( s\right) }\left( t-s\right)  \label{I.b.e.3.9}
\end{equation}%
and%
\begin{equation}
\ln g\left( \left( 1-\lambda \right) t+\lambda s\right) -\ln g\left(
t\right) \geq -\lambda \frac{g^{\prime }\left( t\right) }{g\left( t\right) }%
\left( t-s\right)  \label{I.b.e.3.10}
\end{equation}%
for any $t,s\in $\r{J} and any $\lambda \in \left[ 0,1\right] .$

Now, if we multiply (\ref{I.b.e.3.9}) by $\lambda $ and (\ref{I.b.e.3.10})
by $1-\lambda $ and sum the obtained inequalities, we deduce%
\begin{align}
& \left( 1-\lambda \right) \ln g\left( t\right) +\lambda \ln g\left(
s\right) -\ln g\left( \left( 1-\lambda \right) t+\lambda s\right)
\label{I.b.e.3.11} \\
& \leq \left( 1-\lambda \right) \lambda \left[ \left( \frac{g^{\prime
}\left( t\right) }{g\left( t\right) }-\frac{g^{\prime }\left( s\right) }{%
g\left( s\right) }\right) \left( t-s\right) \right]  \notag
\end{align}%
for any $t,s\in $\r{J} and any $\lambda \in \left[ 0,1\right] .$

Now, if we choose $\lambda :=\frac{M-u}{M-m},s:=m$ and $t:=M$ in (\ref%
{I.b.e.3.11}) then we get the inequality%
\begin{align}
& \frac{u-m}{M-m}\ln g\left( M\right) +\frac{M-u}{M-m}\ln g\left( m\right)
-\ln g\left( u\right)  \label{I.b.e.3.12} \\
& \leq \left[ \frac{\left( M-u\right) \left( u-m\right) }{M-m}\left( \frac{%
g^{\prime }\left( M\right) }{g\left( M\right) }-\frac{g^{\prime }\left(
m\right) }{g\left( m\right) }\right) \right]  \notag
\end{align}%
for any $u\in \left[ m,M\right] .$

If we use the property (\ref{P}) for the operator $A$ we get%
\begin{align}
& \frac{\left\langle Ax,x\right\rangle -m}{M-m}\ln g\left( M\right) +\frac{%
M-\left\langle Ax,x\right\rangle }{M-m}\ln g\left( m\right) -\left\langle
\ln g\left( A\right) x,x\right\rangle  \label{I.b.e.3.13} \\
& \leq \left[ \frac{\left\langle \left( M1_{H}-A\right) \left(
A-m1_{H}\right) x,x\right\rangle }{M-m}\left( \frac{g^{\prime }\left(
M\right) }{g\left( M\right) }-\frac{g^{\prime }\left( m\right) }{g\left(
m\right) }\right) \right]  \notag
\end{align}%
for each $x\in H$ with $\left\Vert x\right\Vert =1.$

Taking the exponential in (\ref{I.b.e.3.13}) we deduce the first inequality
in (\ref{I.b.e.3.9.0}).

Now, consider the function $h:\left[ m,M\right] \rightarrow \mathbb{R}$, $%
h\left( t\right) =\left( M-t\right) \left( t-m\right) .$ This function is
concave in $\left[ m,M\right] $ and by Mond-Pe\v{c}ari\'{c}'s inequality we
have%
\begin{equation*}
\left\langle \left( M1_{H}-A\right) \left( A-m1_{H}\right) x,x\right\rangle
\leq \left( M-\left\langle Ax,x\right\rangle \right) \left( \left\langle
Ax,x\right\rangle -m\right)
\end{equation*}%
for each $x\in H$ with $\left\Vert x\right\Vert =1,$ which proves the second
inequality in (\ref{I.b.e.3.9.0}).

For the last inequality, we observe that%
\begin{equation*}
\left( M-\left\langle Ax,x\right\rangle \right) \left( \left\langle
Ax,x\right\rangle -m\right) \leq \frac{1}{4}\left( M-m\right) ^{2},
\end{equation*}%
and the proof is complete.
\end{proof}

\begin{corollary}[Dragomir, 2010, \protect\cite{I.b.SSD5}]
Assume that $g$ is as in Theorem \ref{I.b.t.3.3} and $A_{j}$ are selfadjoint
operators with $Sp\left( A_{j}\right) \subseteq \left[ m,M\right] \subset $%
\r{J}, $j\in \left\{ 1,...,n\right\} .$

If and $x_{j}\in H,j\in \left\{ 1,...,n\right\} $ with $\sum_{j=1}^{n}\left%
\Vert x_{j}\right\Vert ^{2}=1$, then%
\begin{align}
& \left( 1\leq \right) \frac{\left[ g\left( m\right) \right] ^{\frac{%
M-\sum_{j=1}^{n}\left\langle A_{j}x_{j},x_{j}\right\rangle }{M-m}}\left[
g\left( M\right) \right] ^{\frac{\sum_{j=1}^{n}\left\langle
A_{j}x_{j},x_{j}\right\rangle -m}{M-m}}}{\exp \left(
\sum_{j}^{n}\left\langle \ln g\left( A_{j}\right) x_{j},x_{j}\right\rangle
\right) }  \label{I.b.e.3.14} \\
& \leq \exp \left[ \frac{\sum_{j=1}^{n}\left\langle \left(
M1_{H}-A_{j}\right) \left( A_{j}-m1_{H}\right) x_{j},x_{j}\right\rangle }{M-m%
}\left( \frac{g^{\prime }\left( M\right) }{g\left( M\right) }-\frac{%
g^{\prime }\left( m\right) }{g\left( m\right) }\right) \right]  \notag \\
& \leq \exp \left[ \frac{\left( M-\sum_{j=1}^{n}\left\langle
A_{j}x_{j},x_{j}\right\rangle \right) \left( \sum_{j=1}^{n}\left\langle
A_{j}x_{j},x_{j}\right\rangle -m\right) }{M-m}\left( \frac{g^{\prime }\left(
M\right) }{g\left( M\right) }-\frac{g^{\prime }\left( m\right) }{g\left(
m\right) }\right) \right]  \notag \\
& \leq \exp \left[ \frac{1}{4}\left( M-m\right) \left( \frac{g^{\prime
}\left( M\right) }{g\left( M\right) }-\frac{g^{\prime }\left( m\right) }{%
g\left( m\right) }\right) \right] .  \notag
\end{align}

If $p_{j}\geq 0,$ $j\in \left\{ 1,...,n\right\} $ with $%
\sum_{j=1}^{n}p_{j}=1,$ then 
\begin{align}
& \left( 1\leq \right) \frac{\left[ g\left( m\right) \right] ^{\frac{%
M-\left\langle \sum_{j=1}^{n}p_{j}A_{j}x,x\right\rangle }{M-m}}\left[
g\left( M\right) \right] ^{\frac{\left\langle
\sum_{j=1}^{n}p_{j}A_{j}x,x\right\rangle -m}{M-m}}}{\left\langle
\prod_{j=1}^{n}\left[ g\left( A_{j}\right) \right] ^{p_{j}}x,x\right\rangle }
\label{I.b.e.3.15} \\
& \leq \exp \left[ \frac{\sum_{j=1}^{n}p_{j}\left\langle \left(
M1_{H}-A_{j}\right) \left( A_{j}-m1_{H}\right) x_{j},x_{j}\right\rangle }{M-m%
}\left( \frac{g^{\prime }\left( M\right) }{g\left( M\right) }-\frac{%
g^{\prime }\left( m\right) }{g\left( m\right) }\right) \right]  \notag \\
& \leq \exp \left[ \frac{\left( M-\left\langle
\sum_{j=1}^{n}p_{j}A_{j}x,x\right\rangle \right) \left( \left\langle
\sum_{j=1}^{n}p_{j}A_{j}x,x\right\rangle -m\right) }{M-m}\left( \frac{%
g^{\prime }\left( M\right) }{g\left( M\right) }-\frac{g^{\prime }\left(
m\right) }{g\left( m\right) }\right) \right]  \notag \\
& \leq \exp \left[ \frac{1}{4}\left( M-m\right) \left( \frac{g^{\prime
}\left( M\right) }{g\left( M\right) }-\frac{g^{\prime }\left( m\right) }{%
g\left( m\right) }\right) \right]  \notag
\end{align}%
for each $x\in H$ with $\left\Vert x\right\Vert =1.$
\end{corollary}

\begin{remark}
\label{I.b.r.3.3}Let $A$ be a selfadjoint positive operator on a Hilbert
space $H.$ If $A$ is invertible, then 
\begin{align}
\left( 1\leq \right) \frac{m^{\frac{\left\langle Ax,x\right\rangle -M}{M-m}%
}M^{\frac{m-\left\langle Ax,x\right\rangle }{M-m}}}{\exp \left\langle \ln
A^{-1}x,x\right\rangle }& \leq \exp \left[ \frac{\left\langle \left(
M1_{H}-A\right) \left( A-m1_{H}\right) x,x\right\rangle }{Mm}\right]
\label{I.b.e.3.16} \\
& \leq \exp \left[ \frac{\left( M-\left\langle Ax,x\right\rangle \right)
\left( \left\langle Ax,x\right\rangle -m\right) }{Mm}\right]  \notag \\
& \leq \exp \left[ \frac{1}{4}\frac{\left( M-m\right) ^{2}}{mM}\right] 
\notag
\end{align}%
for all $x\in H$ with $\left\Vert x\right\Vert =1.$
\end{remark}

\section{Hermite-Hadamard's Type Inequalities}

\subsection{Scalar Case}

If $f:I\rightarrow \mathbb{R}$ is a convex function on the interval $I,$
then for any $a,b\in I$ with $a\neq b$ we have the following double
inequality%
\begin{equation}
f\left( \frac{a+b}{2}\right) \leq \frac{1}{b-a}\int_{a}^{b}f\left( t\right)
dt\leq \frac{f\left( a\right) +f\left( b\right) }{2}.  \tag{HH}
\label{I.c.HH}
\end{equation}%
This remarkable result is well known in the literature as the \textit{%
Hermite-Hadamard inequality \cite{ML}}.

For various generalizations, extensions, reverses and related inequalities,
see \cite{I.c.AGP}, \cite{I.c.AZ}, \cite{I.c.F}, \cite{I.c.G}, \cite{I.c.LT}%
, \cite{I.c.L}, \cite{I.c.Mk}, \cite{ML}\ the monograph \cite{I.c.DP} and
the references therein.

\subsection{Some Inequalities for Convex Functions}

The following inequality related to the Mond-Pe\v{c}ari\'{c} result also
holds:

\begin{theorem}[Dragomir, 2010, \protect\cite{I.c.SSD5}]
\label{I.c.t.2.1}Let $A$ be a selfadjoint operator on the Hilbert space $H$
and assume that $Sp\left( A\right) \subseteq \left[ m,M\right] $ for some
scalars $m,M$ with $m<M.$

If $f$ is a convex function on $\left[ m,M\right] ,$ then 
\begin{align}
\frac{f\left( m\right) +f\left( M\right) }{2}& \geq \left\langle \frac{%
f\left( A\right) +f\left( \left( m+M\right) 1_{H}-A\right) }{2}%
x,x\right\rangle  \label{I.c.e.2.1} \\
& \geq \frac{f\left( \left\langle Ax,x\right\rangle \right) +f\left(
m+M-\left\langle Ax,x\right\rangle \right) }{2}  \notag \\
& \geq f\left( \frac{m+M}{2}\right)  \notag
\end{align}%
for each $x\in H$ with $\left\Vert x\right\Vert =1.$

In addition, if $x\in H$ with $\left\Vert x\right\Vert =1$ and $\left\langle
Ax,x\right\rangle \neq \frac{m+M}{2},$ then also%
\begin{align}
& \frac{f\left( \left\langle Ax,x\right\rangle \right) +f\left(
m+M-\left\langle Ax,x\right\rangle \right) }{2}  \label{I.c.e.2.1.a} \\
& \geq \frac{2}{\frac{m+M}{2}-\left\langle Ax,x\right\rangle }%
\int_{\left\langle Ax,x\right\rangle }^{m+M-\left\langle Ax,x\right\rangle
}f\left( u\right) du\geq f\left( \frac{m+M}{2}\right) .  \notag
\end{align}
\end{theorem}

\begin{proof}
Since $f$ is convex on $\left[ m,M\right] $ then for each $u\in \left[ m,M%
\right] $ we have the inequalities%
\begin{equation}
\frac{M-u}{M-m}f\left( m\right) +\frac{u-m}{M-m}f\left( M\right) \geq
f\left( \frac{M-u}{M-m}m+\frac{u-m}{M-m}M\right) =f\left( u\right)
\label{I.c.e.2.2}
\end{equation}%
and%
\begin{align}
\frac{M-u}{M-m}f\left( M\right) +\frac{u-m}{M-m}f\left( m\right) & \geq
f\left( \frac{M-u}{M-m}M+\frac{u-m}{M-m}m\right)  \label{I.c.e.2.3} \\
& =f\left( M+m-u\right) .  \notag
\end{align}%
If we add these two inequalities we get%
\begin{equation*}
f\left( m\right) +f\left( M\right) \geq f\left( u\right) +f\left(
M+m-u\right)
\end{equation*}%
for any $u\in \left[ m,M\right] ,$ which, by the property (\ref{P}) applied
for the operator $A,$ produces the first inequality in (\ref{I.c.e.2.1}).

By the Mond-Pe\v{c}ari\'{c} inequality we have 
\begin{equation*}
\left\langle f\left( \left( m+M\right) 1_{H}-A\right) x,x\right\rangle \geq
f\left( m+M-\left\langle Ax,x\right\rangle \right) ,
\end{equation*}%
which together with the same inequality produces the second inequality in (%
\ref{I.c.e.2.1}).

The third part follows by the convexity of $f.$

In order to prove (\ref{I.c.e.2.1.a}), we use the Hermite-Hadamard
inequality (\ref{I.c.HH}) for the convex functions $f$ and the choices $%
a=\left\langle Ax,x\right\rangle $ and $b=m+M-\left\langle Ax,x\right\rangle
.$

The proof is complete.
\end{proof}

\begin{remark}
\label{I.c.r2.1}We observe that, from the inequality (\ref{I.c.e.2.1}) we
have the following inequality in the operator order of $B\left( H\right) $%
\begin{equation}
\left[ \frac{f\left( m\right) +f\left( M\right) }{2}\right] 1_{H}\geq \frac{%
f\left( A\right) +f\left( \left( m+M\right) 1_{H}-A\right) }{2}\geq f\left( 
\frac{m+M}{2}\right) 1_{H},  \label{I.c.e.2.5.0}
\end{equation}%
where $f$ is a convex function on $\left[ m,M\right] $ and $A$ a selfadjoint
operator on the Hilbert space $H$ with $Sp\left( A\right) \subseteq \left[
m,M\right] $ for some scalars $m,M$ with $m<M.$
\end{remark}

The case of log-convex functions may be of interest for applications and
therefore is stated in:

\begin{corollary}[Dragomir, 2010, \protect\cite{I.c.SSD5}]
\label{I.c.c.2.1}If $g$ is a log-convex function on $\left[ m,M\right] ,$
then%
\begin{align}
\sqrt{g\left( m\right) g\left( M\right) }& \geq \exp \left\langle \ln \left[
g\left( A\right) g\left( \left( m+M\right) 1_{H}-A\right) \right]
^{1/2}x,x\right\rangle  \label{I.c.e.2.5} \\
& \geq \sqrt{g\left( \left\langle Ax,x\right\rangle \right) g\left(
m+M-\left\langle Ax,x\right\rangle \right) }  \notag \\
& \geq g\left( \frac{m+M}{2}\right)  \notag
\end{align}%
for each $x\in H$ with $\left\Vert x\right\Vert =1.$

In addition, if $x\in H$ with $\left\Vert x\right\Vert =1$ and $\left\langle
Ax,x\right\rangle \neq \frac{m+M}{2},$ then also%
\begin{align}
& \sqrt{g\left( \left\langle Ax,x\right\rangle \right) g\left(
m+M-\left\langle Ax,x\right\rangle \right) }  \label{I.c.e.2.6} \\
& \geq \exp \left[ \frac{2}{\frac{m+M}{2}-\left\langle Ax,x\right\rangle }%
\int_{\left\langle Ax,x\right\rangle }^{m+M-\left\langle Ax,x\right\rangle
}\ln g\left( u\right) du\right]  \notag \\
& \geq g\left( \frac{m+M}{2}\right) .  \notag
\end{align}
\end{corollary}

The following result also holds

\begin{theorem}[Dragomir, 2010, \protect\cite{I.c.SSD5}]
\label{I.c.t.2.2}Let $A$ and $B$ selfadjoint operators on the Hilbert space $%
H$ and assume that $Sp\left( A\right) ,Sp\left( B\right) \subseteq \left[ m,M%
\right] $ for some scalars $m,M$ with $m<M.$

If $f$ is a convex function on $\left[ m,M\right] ,$ then%
\begin{align}
& f\left( \left\langle \frac{A+B}{2}x,x\right\rangle \right)
\label{I.c.e.2.7} \\
& \leq \frac{1}{2}\left[ f\left( \left( 1-t\right) \left\langle
Ax,x\right\rangle +t\left\langle Bx,x\right\rangle \right) +f\left(
t\left\langle Ax,x\right\rangle +\left( 1-t\right) \left\langle
Bx,x\right\rangle \right) \right]  \notag \\
& \leq \left\langle \frac{1}{2}\left[ f\left( \left( 1-t\right) A+tB\right)
+f\left( tA+\left( 1-t\right) B\right) \right] x,x\right\rangle  \notag \\
& \leq \frac{M-\left\langle \frac{A+B}{2}x,x\right\rangle }{M-m}f\left(
m\right) +\frac{\left\langle \frac{A+B}{2}x,x\right\rangle -m}{M-m}f\left(
M\right)  \notag
\end{align}%
for any $t\in \left[ 0,1\right] $ and each $x\in H$ with $\left\Vert
x\right\Vert =1.$

Moreover, we have the Hermite-Hadamard's type inequalities:%
\begin{align}
& f\left( \left\langle \frac{A+B}{2}x,x\right\rangle \right)
\label{I.c.e.2.8} \\
& \leq \int_{0}^{1}f\left( \left( 1-t\right) \left\langle Ax,x\right\rangle
+t\left\langle Bx,x\right\rangle \right) dt  \notag \\
& \leq \left\langle \left[ \int_{0}^{1}f\left( \left( 1-t\right) A+tB\right)
dt\right] x,x\right\rangle  \notag \\
& \leq \frac{M-\left\langle \frac{A+B}{2}x,x\right\rangle }{M-m}f\left(
m\right) +\frac{\left\langle \frac{A+B}{2}x,x\right\rangle -m}{M-m}f\left(
M\right)  \notag
\end{align}%
each $x\in H$ with $\left\Vert x\right\Vert =1.$

In addition, if we assume that $B-A$ is a positive definite operator, then%
\begin{align}
& f\left( \left\langle \frac{A+B}{2}x,x\right\rangle \right) \left\langle
\left( B-A\right) x,x\right\rangle  \label{I.c.e.2.8.a} \\
& \leq \int_{\left\langle Ax,x\right\rangle }^{\left\langle
Bx,x\right\rangle }f\left( u\right) du\leq \left\langle \left( B-A\right)
x,x\right\rangle \left\langle \left[ \int_{0}^{1}f\left( \left( 1-t\right)
A+tB\right) dt\right] x,x\right\rangle  \notag \\
& \leq \left\langle \left( B-A\right) x,x\right\rangle \left[ \frac{%
M-\left\langle \frac{A+B}{2}x,x\right\rangle }{M-m}f\left( m\right) +\frac{%
\left\langle \frac{A+B}{2}x,x\right\rangle -m}{M-m}f\left( M\right) \right] .
\notag
\end{align}
\end{theorem}

\begin{proof}
It is obvious that for any $t\in \left[ 0,1\right] $ we have $Sp\left(
\left( 1-t\right) A+tB\right) ,Sp\left( tA+\left( 1-t\right) B\right)
\subseteq \left[ m,M\right] .$

On making use of the Mond-Pe\v{c}ari\'{c} inequality we have%
\begin{equation}
f\left( \left( 1-t\right) \left\langle Ax,x\right\rangle +t\left\langle
Bx,x\right\rangle \right) \leq \left\langle f\left( \left( 1-t\right)
A+tB\right) x,x\right\rangle  \label{I.c.e.2.9}
\end{equation}%
and%
\begin{equation}
f\left( t\left\langle Ax,x\right\rangle +\left( 1-t\right) \left\langle
Bx,x\right\rangle \right) \leq \left\langle f\left( tA+\left( 1-t\right)
B\right) x,x\right\rangle  \label{I.c.e.2.10}
\end{equation}%
for any $t\in \left[ 0,1\right] $ and each $x\in H$ with $\left\Vert
x\right\Vert =1.$

Adding (\ref{I.c.e.2.9}) with (\ref{I.c.e.2.10}) and utilising the convexity
of $f$ we deduce the first two inequalities in (\ref{I.c.e.2.7}).

By the Lah-Ribari\'{c} inequality (\ref{I.a.e.2.4}) we also have%
\begin{align}
\left\langle f\left( \left( 1-t\right) A+tB\right) x,x\right\rangle & \leq 
\frac{M-\left( 1-t\right) \left\langle Ax,x\right\rangle -t\left\langle
Bx,x\right\rangle }{M-m}\cdot f\left( m\right)  \label{I.c.e.2.11} \\
& +\frac{\left( 1-t\right) \left\langle Ax,x\right\rangle +t\left\langle
Bx,x\right\rangle -m}{M-m}\cdot f\left( M\right)  \notag
\end{align}%
and%
\begin{align}
\left\langle f\left( tA+\left( 1-t\right) B\right) x,x\right\rangle & \leq 
\frac{M-t\left\langle Ax,x\right\rangle -\left( 1-t\right) \left\langle
Bx,x\right\rangle }{M-m}\cdot f\left( m\right)  \label{I.c.e.2.12} \\
& +\frac{t\left\langle Ax,x\right\rangle +\left( 1-t\right) \left\langle
Bx,x\right\rangle -m}{M-m}\cdot f\left( M\right)  \notag
\end{align}%
for any $t\in \left[ 0,1\right] $ and each $x\in H$ with $\left\Vert
x\right\Vert =1.$

Now, if we add the inequalities (\ref{I.c.e.2.11}) with (\ref{I.c.e.2.12})
and divide by two, we deduce the last part in (\ref{I.c.e.2.7}).

Integrating the inequality over $t\in \left[ 0,1\right] $, utilising the
continuity property of the inner product and the properties of the integral
of operator-valued functions we have%
\begin{align}
& f\left( \left\langle \frac{A+B}{2}x,x\right\rangle \right)
\label{I.c.e.2.13} \\
& \leq \frac{1}{2}\left[ \int_{0}^{1}f\left( \left( 1-t\right) \left\langle
Ax,x\right\rangle +t\left\langle Bx,x\right\rangle \right)
dt+\int_{0}^{1}f\left( t\left\langle Ax,x\right\rangle +\left( 1-t\right)
\left\langle Bx,x\right\rangle \right) dt\right]  \notag \\
& \leq \left\langle \frac{1}{2}\left[ \int_{0}^{1}f\left( \left( 1-t\right)
A+tB\right) dt+\int_{0}^{1}f\left( tA+\left( 1-t\right) B\right) dt\right]
x,x\right\rangle  \notag \\
& \leq \frac{M-\left\langle \frac{A+B}{2}x,x\right\rangle }{M-m}f\left(
m\right) +\frac{\left\langle \frac{A+B}{2}x,x\right\rangle -m}{M-m}f\left(
M\right) .  \notag
\end{align}%
Since%
\begin{equation*}
\int_{0}^{1}f\left( \left( 1-t\right) \left\langle Ax,x\right\rangle
+t\left\langle Bx,x\right\rangle \right) dt=\int_{0}^{1}f\left(
t\left\langle Ax,x\right\rangle +\left( 1-t\right) \left\langle
Bx,x\right\rangle \right) dt
\end{equation*}%
and%
\begin{equation*}
\int_{0}^{1}f\left( \left( 1-t\right) A+tB\right) dt=\int_{0}^{1}f\left(
tA+\left( 1-t\right) B\right) dt
\end{equation*}%
then, by (\ref{I.c.e.2.13}), we deduce the inequality (\ref{I.c.e.2.8}).

The inequality (\ref{I.c.e.2.8.a}) follows from (\ref{I.c.e.2.8}) by
observing that for $\left\langle Bx,x\right\rangle >\left\langle
Ax,x\right\rangle $ we have%
\begin{equation*}
\int_{0}^{1}f\left( \left( 1-t\right) \left\langle Ax,x\right\rangle
+t\left\langle Bx,x\right\rangle \right) dt=\frac{1}{\left\langle
Bx,x\right\rangle -\left\langle Ax,x\right\rangle }\int_{\left\langle
Ax,x\right\rangle }^{\left\langle Bx,x\right\rangle }f\left( u\right) du
\end{equation*}%
for each $x\in H$ with $\left\Vert x\right\Vert =1.$
\end{proof}

\begin{remark}
\label{I.c.r2.2}We observe that, from the inequalities (\ref{I.c.e.2.7}) and
(\ref{I.c.e.2.8}) we have the following inequalities in the operator order
of $B\left( H\right) $%
\begin{align}
& \frac{1}{2}\left[ f\left( \left( 1-t\right) A+tB\right) +f\left( tA+\left(
1-t\right) B\right) \right]  \label{I.c.e.2.13.0} \\
& \leq f\left( m\right) \frac{M1_{H}-\frac{A+B}{2}}{M-m}+f\left( M\right) 
\frac{\frac{A+B}{2}-m1_{H}}{M-m},  \notag
\end{align}%
where $f$ is a convex function on $\left[ m,M\right] $ and $A,B$ are
selfadjoint operator on the Hilbert space $H$ with $Sp\left( A\right)
,Sp\left( B\right) \subseteq \left[ m,M\right] $ for some scalars $m,M$ with 
$m<M.$
\end{remark}

The case of log-convex functions is as follows:

\begin{corollary}[Dragomir, 2010, \protect\cite{I.c.SSD5}]
\label{I.c.c.2.3}If $g$ is a log-convex function on $\left[ m,M\right] ,$
then%
\begin{align}
& g\left( \left\langle \frac{A+B}{2}x,x\right\rangle \right)
\label{I.c.e.2.13.a} \\
& \leq \sqrt{g\left( \left( 1-t\right) \left\langle Ax,x\right\rangle
+t\left\langle Bx,x\right\rangle \right) g\left( t\left\langle
Ax,x\right\rangle +\left( 1-t\right) \left\langle Bx,x\right\rangle \right) }
\notag \\
& \leq \exp \left\langle \frac{1}{2}\left[ \ln g\left( \left( 1-t\right)
A+tB\right) +\ln g\left( tA+\left( 1-t\right) B\right) \right]
x,x\right\rangle  \notag \\
& \leq g\left( m\right) ^{\frac{M-\left\langle \frac{A+B}{2}x,x\right\rangle 
}{M-m}}g\left( M\right) ^{\frac{\left\langle \frac{A+B}{2}x,x\right\rangle -m%
}{M-m}}  \notag
\end{align}%
for any $t\in \left[ 0,1\right] $ and each $x\in H$ with $\left\Vert
x\right\Vert =1.$

Moreover, we have the Hermite-Hadamard's type inequalities:%
\begin{align}
& g\left( \left\langle \frac{A+B}{2}x,x\right\rangle \right)
\label{I.c.e.2.13.b} \\
& \leq \exp \left[ \int_{0}^{1}\ln g\left( \left( 1-t\right) \left\langle
Ax,x\right\rangle +t\left\langle Bx,x\right\rangle \right) dt\right]  \notag
\\
& \leq \exp \left\langle \left[ \int_{0}^{1}\ln g\left( \left( 1-t\right)
A+tB\right) dt\right] x,x\right\rangle  \notag \\
& \leq g\left( m\right) ^{\frac{M-\left\langle \frac{A+B}{2}x,x\right\rangle 
}{M-m}}g\left( M\right) ^{\frac{\left\langle \frac{A+B}{2}x,x\right\rangle -m%
}{M-m}}  \notag
\end{align}%
for each $x\in H$ with $\left\Vert x\right\Vert =1.$

In addition, if we assume that $B-A$ is a positive definite operator, then%
\begin{align}
& g\left( \left\langle \frac{A+B}{2}x,x\right\rangle \right) ^{\left\langle
\left( B-A\right) x,x\right\rangle }  \label{I.c.e.2.13.c} \\
& \leq \exp \left[ \int_{\left\langle Ax,x\right\rangle }^{\left\langle
Bx,x\right\rangle }\ln g\left( u\right) du\right]  \notag \\
& \leq \exp \left[ \left\langle \left( B-A\right) x,x\right\rangle
\left\langle \left[ \int_{0}^{1}\ln g\left( \left( 1-t\right) A+tB\right) dt%
\right] x,x\right\rangle \right]  \notag \\
& \leq \left[ g\left( m\right) ^{\frac{M-\left\langle \frac{A+B}{2}%
x,x\right\rangle }{M-m}}g\left( M\right) ^{\frac{\left\langle \frac{A+B}{2}%
x,x\right\rangle -m}{M-m}}\right] ^{\left\langle \left( B-A\right)
x,x\right\rangle }  \notag
\end{align}%
for each $x\in H$ with $\left\Vert x\right\Vert =1.$
\end{corollary}

From a different perspective we have the following result as well:

\begin{theorem}[Dragomir, 2010, \protect\cite{I.c.SSD5}]
\label{I.c.t.2.3}Let $A$ and $B$ selfadjoint operators on the Hilbert space $%
H$ and assume that $Sp\left( A\right) ,Sp\left( B\right) \subseteq \left[ m,M%
\right] $ for some scalars $m,M$ with $m<M.$ If $f$ is a convex function on $%
\left[ m,M\right] ,$ then%
\begin{align}
& f\left( \frac{\left\langle Ax,x\right\rangle +\left\langle
By,y\right\rangle }{2}\right)  \label{I.c.e.2.14} \\
& \leq \int_{0}^{1}f\left( \left( 1-t\right) \left\langle Ax,x\right\rangle
+t\left\langle By,y\right\rangle \right) dt  \notag \\
& \leq \left\langle \left[ \int_{0}^{1}f\left( \left( 1-t\right)
A+t\left\langle By,y\right\rangle 1_{H}\right) dt\right] x,x\right\rangle 
\notag \\
& \leq \frac{1}{2}\left[ \left\langle f\left( A\right) x,x\right\rangle
+f\left( \left\langle By,y\right\rangle \right) \right]  \notag \\
& \leq \frac{1}{2}\left[ \left\langle f\left( A\right) x,x\right\rangle
+\left\langle f\left( B\right) y,y\right\rangle \right]  \notag
\end{align}%
and%
\begin{align}
f\left( \frac{\left\langle Ax,x\right\rangle +\left\langle By,y\right\rangle 
}{2}\right) & \leq \left\langle f\left( \frac{A+\left\langle
By,y\right\rangle 1_{H}}{2}\right) x,x\right\rangle  \label{I.c.e.2.14.a} \\
& \leq \left\langle \left[ \int_{0}^{1}f\left( \left( 1-t\right)
A+t\left\langle By,y\right\rangle 1_{H}\right) dt\right] x,x\right\rangle 
\notag
\end{align}%
for each $x,y\in H$ with $\left\Vert x\right\Vert =\left\Vert y\right\Vert
=1.$
\end{theorem}

\begin{proof}
For a convex function $f$ and any $u,v\in \left[ m,M\right] $ and $t\in %
\left[ 0,1\right] $ we have the double inequality:%
\begin{align}
f\left( \frac{u+v}{2}\right) & \leq \frac{1}{2}\left[ f\left( \left(
1-t\right) u+tv\right) +f\left( tu+\left( 1-t\right) v\right) \right]
\label{I.c.e.2.15} \\
& \leq \frac{1}{2}\left[ f\left( u\right) +f\left( v\right) \right] .  \notag
\end{align}%
Utilising the second inequality in (\ref{I.c.e.2.15}) we have%
\begin{align}
& \frac{1}{2}\left[ f\left( \left( 1-t\right) u+t\left\langle
By,y\right\rangle \right) +f\left( tu+\left( 1-t\right) \left\langle
By,y\right\rangle \right) \right]  \label{I.c.e.2.16} \\
& \leq \frac{1}{2}\left[ f\left( u\right) +f\left( \left\langle
By,y\right\rangle \right) \right]  \notag
\end{align}%
for any $u\in \left[ m,M\right] $, $t\in \left[ 0,1\right] $ and $y\in H$
with $\left\Vert y\right\Vert =1.$

Now, on applying the property (\ref{P}) to the inequality (\ref{I.c.e.2.16})
for the operator $A$ we have 
\begin{align}
& \frac{1}{2}\left[ \left\langle f\left( \left( 1-t\right) A+t\left\langle
By,y\right\rangle \right) x,x\right\rangle +\left\langle f\left( tA+\left(
1-t\right) \left\langle By,y\right\rangle \right) x,x\right\rangle \right]
\label{I.c.e.2.17} \\
& \leq \frac{1}{2}\left[ \left\langle f\left( A\right) x,x\right\rangle
+f\left( \left\langle By,y\right\rangle \right) \right]  \notag
\end{align}%
for any $t\in \left[ 0,1\right] $ and $x,y\in H$ with $\left\Vert
x\right\Vert =\left\Vert y\right\Vert =1.$

On applying the Mond-Pe\v{c}ari\'{c} inequality we also have%
\begin{align}
& \frac{1}{2}\left[ f\left( \left( 1-t\right) \left\langle Ax,x\right\rangle
+t\left\langle By,y\right\rangle \right) +f\left( t\left\langle
Ax,x\right\rangle +\left( 1-t\right) \left\langle By,y\right\rangle \right) %
\right]  \label{I.c.e.2.18} \\
& \leq \frac{1}{2}\left[ \left\langle f\left( \left( 1-t\right)
A+t\left\langle By,y\right\rangle 1_{H}\right) x,x\right\rangle
+\left\langle f\left( tA+\left( 1-t\right) \left\langle By,y\right\rangle
1_{H}\right) x,x\right\rangle \right]  \notag
\end{align}%
for any $t\in \left[ 0,1\right] $ and $x,y\in H$ with $\left\Vert
x\right\Vert =\left\Vert y\right\Vert =1.$

Now, integrating over $t$ on $\left[ 0,1\right] $ the inequalities (\ref%
{I.c.e.2.17}) and (\ref{I.c.e.2.18}) and taking into account that%
\begin{align*}
& \int_{0}^{1}\left\langle f\left( \left( 1-t\right) A+t\left\langle
By,y\right\rangle 1_{H}\right) x,x\right\rangle dt \\
& =\int_{0}^{1}\left\langle f\left( tA+\left( 1-t\right) \left\langle
By,y\right\rangle 1_{H}\right) x,x\right\rangle dt \\
& =\left\langle \left[ \int_{0}^{1}f\left( \left( 1-t\right) A+t\left\langle
By,y\right\rangle 1_{H}\right) dt\right] x,x\right\rangle
\end{align*}%
and%
\begin{equation*}
\int_{0}^{1}f\left( \left( 1-t\right) \left\langle Ax,x\right\rangle
+t\left\langle By,y\right\rangle \right) dt=\int_{0}^{1}f\left(
t\left\langle Ax,x\right\rangle +\left( 1-t\right) \left\langle
By,y\right\rangle \right) dt,
\end{equation*}%
we obtain the second and the third inequality in (\ref{I.c.e.2.14}).

Further, on applying the Jensen integral inequality for the convex function $%
f$ we also have%
\begin{align*}
& \int_{0}^{1}f\left( \left( 1-t\right) \left\langle Ax,x\right\rangle
+t\left\langle By,y\right\rangle \right) dt \\
& \geq f\left( \int_{0}^{1}\left[ \left( 1-t\right) \left\langle
Ax,x\right\rangle +t\left\langle By,y\right\rangle \right] dt\right) \\
& =f\left( \frac{\left\langle Ax,x\right\rangle +\left\langle
By,y\right\rangle }{2}\right)
\end{align*}%
for each $x,y\in H$ with $\left\Vert x\right\Vert =\left\Vert y\right\Vert
=1 $, proving the first part of (\ref{I.c.e.2.14}).

Now, on utilising the first part of (\ref{I.c.e.2.15}) we can also state that%
\begin{equation}
f\left( \frac{u+\left\langle By,y\right\rangle }{2}\right) \leq \frac{1}{2}%
\left[ f\left( \left( 1-t\right) u+t\left\langle By,y\right\rangle \right)
+f\left( tu+\left( 1-t\right) \left\langle By,y\right\rangle \right) \right]
\label{I.c.e.2.19}
\end{equation}%
for any $u\in \left[ m,M\right] $, $t\in \left[ 0,1\right] $ and $y\in H$
with $\left\Vert y\right\Vert =1.$

Further, on applying the property (\ref{P}) to the inequality (\ref%
{I.c.e.2.19}) and for the operator $A$ we get 
\begin{align*}
& \left\langle f\left( \frac{A+\left\langle By,y\right\rangle 1_{H}}{2}%
\right) x,x\right\rangle \\
& \leq \frac{1}{2}\left[ \left\langle f\left( \left( 1-t\right)
A+t\left\langle By,y\right\rangle 1_{H}\right) x,x\right\rangle
+\left\langle f\left( tA+\left( 1-t\right) \left\langle By,y\right\rangle
1_{H}\right) x,x\right\rangle \right]
\end{align*}%
for each $x,y\in H$ with $\left\Vert x\right\Vert =\left\Vert y\right\Vert
=1,$ which, by integration over $t$ in $\left[ 0,1\right] $ produces the
second inequality in (\ref{I.c.e.2.14.a}). The first inequality is obvious.
\end{proof}

\begin{remark}
\label{I.c.r2.3}It is important to remark that, from the inequalities (\ref%
{I.c.e.2.14}) and (\ref{I.c.e.2.14.a}) we have the following
Hermite-Hadamard's type results in the operator order of $B\left( H\right) $
and for the convex function $f:\left[ m,M\right] \rightarrow \mathbb{R}$%
\begin{align}
f\left( \frac{A+\left\langle By,y\right\rangle 1_{H}}{2}\right) & \leq
\int_{0}^{1}f\left( \left( 1-t\right) A+t\left\langle By,y\right\rangle
1_{H}\right) dt  \label{I.c.e.2.20.0} \\
& \leq \frac{1}{2}\left[ f\left( A\right) +f\left( \left\langle
By,y\right\rangle \right) 1_{H}\right]  \notag
\end{align}%
for any $y\in H$ with $\left\Vert y\right\Vert =1$ and any selfadjoint
operators $A,B$ with spectra in $\left[ m,M\right] .$

In particular, we have from (\ref{I.c.e.2.20.0}) 
\begin{align}
f\left( \frac{A+\left\langle Ay,y\right\rangle 1_{H}}{2}\right) & \leq
\int_{0}^{1}f\left( \left( 1-t\right) A+t\left\langle Ay,y\right\rangle
1_{H}\right) dt  \label{I.c.e.2.20.a} \\
& \leq \frac{1}{2}\left[ f\left( A\right) +f\left( \left\langle
Ay,y\right\rangle \right) 1_{H}\right]  \notag
\end{align}%
for any $y\in H$ with $\left\Vert y\right\Vert =1$ and%
\begin{equation}
f\left( \frac{A+s1_{H}}{2}\right) \leq \int_{0}^{1}f\left( \left( 1-t\right)
A+ts1_{H}\right) dt\leq \frac{1}{2}\left[ f\left( A\right) +f\left( s\right)
1_{H}\right]  \label{I.c.e.2.20.b}
\end{equation}%
for any $s\in \left[ m,M\right] .$
\end{remark}

As a particular case of the above theorem we have the following refinement
of the Mond-Pe\v{c}ari\'{c} inequality:

\begin{corollary}[Dragomir, 2010, \protect\cite{I.c.SSD5}]
\label{I.c.c.2.2}Let $A$ be a selfadjoint operator on the Hilbert space $H$
and assume that $Sp\left( A\right) \subseteq \left[ m,M\right] $ for some
scalars $m,M$ with $m<M.$ If $f$ is a convex function on $\left[ m,M\right]
, $ then%
\begin{align}
f\left( \left\langle Ax,x\right\rangle \right) & \leq \left\langle f\left( 
\frac{A+\left\langle Ax,x\right\rangle 1_{H}}{2}\right) x,x\right\rangle
\label{I.c.e.2.20} \\
& \leq \left\langle \left[ \int_{0}^{1}f\left( \left( 1-t\right)
A+t\left\langle Ax,x\right\rangle 1_{H}\right) dt\right] x,x\right\rangle 
\notag \\
& \leq \frac{1}{2}\left[ \left\langle f\left( A\right) x,x\right\rangle
+f\left( \left\langle Ax,x\right\rangle \right) \right] \leq \left\langle
f\left( A\right) x,x\right\rangle .  \notag
\end{align}
\end{corollary}

Finally, the case of log-convex functions is as follows:

\begin{corollary}[Dragomir, 2010, \protect\cite{I.c.SSD5}]
\label{I.c.c.2.4}If $g$ is a log-convex function on $\left[ m,M\right] ,$
then%
\begin{align}
& g\left( \frac{\left\langle Ax,x\right\rangle +\left\langle
By,y\right\rangle }{2}\right)  \label{I.c.e.2.21} \\
& \leq \exp \left[ \int_{0}^{1}\ln g\left( \left( 1-t\right) \left\langle
Ax,x\right\rangle +t\left\langle By,y\right\rangle \right) dt\right]  \notag
\\
& \leq \exp \left\langle \left[ \int_{0}^{1}\ln g\left( \left( 1-t\right)
A+t\left\langle By,y\right\rangle 1_{H}\right) dt\right] x,x\right\rangle 
\notag \\
& \leq \exp \left[ \frac{1}{2}\left[ \left\langle \ln g\left( A\right)
x,x\right\rangle +\ln g\left( \left\langle By,y\right\rangle \right) \right] %
\right]  \notag \\
& \leq \exp \left[ \frac{1}{2}\left[ \left\langle \ln g\left( A\right)
x,x\right\rangle +\left\langle \ln g\left( B\right) y,y\right\rangle \right] %
\right]  \notag
\end{align}%
and%
\begin{align}
g\left( \frac{\left\langle Ax,x\right\rangle +\left\langle By,y\right\rangle 
}{2}\right) & \leq \exp \left\langle \ln g\left( \frac{A+\left\langle
By,y\right\rangle 1_{H}}{2}\right) x,x\right\rangle  \label{I.c.e.2.22} \\
& \leq \exp \left\langle \left[ \int_{0}^{1}\ln g\left( \left( 1-t\right)
A+t\left\langle By,y\right\rangle 1_{H}\right) dt\right] x,x\right\rangle 
\notag
\end{align}%
and%
\begin{align}
g\left( \left\langle Ax,x\right\rangle \right) & \leq \exp \left\langle \ln
g\left( \frac{A+\left\langle Ax,x\right\rangle 1_{H}}{2}\right)
x,x\right\rangle  \label{I.c.e.2.23} \\
& \leq \exp \left\langle \left[ \int_{0}^{1}\ln g\left( \left( 1-t\right)
A+t\left\langle Ax,x\right\rangle 1_{H}\right) dt\right] x,x\right\rangle 
\notag \\
& \leq \exp \left[ \frac{1}{2}\left[ \left\langle \ln g\left( A\right)
x,x\right\rangle +\ln g\left( \left\langle Ax,x\right\rangle \right) \right] %
\right] \leq \exp \left\langle \ln g\left( A\right) x,x\right\rangle  \notag
\end{align}%
respectively, for each $x\in H$ with $\left\Vert x\right\Vert =1$ and $A,B$
selfadjoint operators with spectra in $\left[ m,M\right] .$
\end{corollary}

It is obvious that all the above inequalities can be applied for particular
convex or log-convex functions of interest. However, we will restrict
ourselves to only a few examples that are connected with famous results such
as the H\"{o}lder-McCarthy inequality or the Ky Fan inequality.

\subsection{Applications for H\"{o}lder-McCarthy's Inequality}

We can improve the H\"{o}lder-McCarthy's inequality above as follows:

\begin{proposition}
\label{I.c.p.3.1}Let $A$ be a selfadjoint positive operator on a Hilbert
space $H$.

If $r>1,$ then%
\begin{align}
\left\langle Ax,x\right\rangle ^{r}& \leq \left\langle \left( \frac{%
A+\left\langle Ax,x\right\rangle 1_{H}}{2}\right) ^{r}x,x\right\rangle
\label{I.c.e.3.1} \\
& \leq \left\langle \left[ \int_{0}^{1}\left( \left( 1-t\right)
A+t\left\langle Ax,x\right\rangle 1_{H}\right) ^{r}dt\right] x,x\right\rangle
\notag \\
& \leq \frac{1}{2}\left[ \left\langle A^{r}x,x\right\rangle +\left\langle
Ax,x\right\rangle ^{r}\right] \leq \left\langle A^{r}x,x\right\rangle  \notag
\end{align}%
for any $x\in H$ with $\left\Vert x\right\Vert =1.$

If $0<r<1,$ then the inequalities reverse in (\ref{I.c.e.3.1}).

If $A$ is invertible and $r>0,$ then 
\begin{align}
\left\langle Ax,x\right\rangle ^{-r}& \leq \left\langle \left( \frac{%
A+\left\langle Ax,x\right\rangle 1_{H}}{2}\right) ^{-r}x,x\right\rangle
\label{I.c.e.3.2} \\
& \leq \left\langle \left[ \int_{0}^{1}\left( \left( 1-t\right)
A+t\left\langle Ax,x\right\rangle 1_{H}\right) ^{-r}dt\right]
x,x\right\rangle  \notag \\
& \leq \frac{1}{2}\left[ \left\langle A^{-r}x,x\right\rangle +\left\langle
Ax,x\right\rangle ^{-r}\right] \leq \left\langle A^{-r}x,x\right\rangle 
\notag
\end{align}%
for any $x\in H$ with $\left\Vert x\right\Vert =1.$
\end{proposition}

Follows from the inequality (\ref{I.c.e.2.21}) applied for the power
function.

Since the function $g\left( t\right) =t^{-r}$ for $r>0$ is log-convex, then
by utilising the inequality (\ref{I.c.e.2.23}) we can improve the H\"{o}%
lder-McCarthy inequality as follows:

\begin{proposition}
\label{I.c.p.3.2}Let $A$ be a selfadjoint positive operator on a Hilbert
space $H.$ If $A$ is invertible, then 
\begin{align}
\left\langle Ax,x\right\rangle ^{-r}& \leq \exp \left\langle \ln \left( 
\frac{A+\left\langle Ax,x\right\rangle 1_{H}}{2}\right) ^{-r}x,x\right\rangle
\label{e.3.3} \\
& \leq \exp \left\langle \left[ \int_{0}^{1}\ln \left( \left( 1-t\right)
A+t\left\langle Ax,x\right\rangle 1_{H}\right) ^{-r}dt\right]
x,x\right\rangle  \notag \\
& \leq \exp \left[ \frac{1}{2}\left[ \left\langle \ln A^{-r}x,x\right\rangle
+\ln \left\langle Ax,x\right\rangle ^{-r}\right] \right] \leq \exp
\left\langle \ln A^{-r}x,x\right\rangle  \notag
\end{align}%
for all $r>0$ and $x\in H$ with $\left\Vert x\right\Vert =1.$
\end{proposition}

Now, from a different perspective, we can state the following operator power
inequalities:

\begin{proposition}
\label{I.c.p.3.3}Let $A$ be a selfadjoint operator with $Sp\left( A\right)
\subseteq \left[ m,M\right] \subset \lbrack 0,\infty ),$ then 
\begin{align}
\frac{m^{r}+M^{r}}{2}& \geq \left\langle \frac{A^{r}+\left( \left(
m+M\right) 1_{H}-A\right) ^{r}}{2}x,x\right\rangle  \label{I.c.e.3.4} \\
& \geq \frac{\left\langle Ax,x\right\rangle ^{r}+\left( m+M-\left\langle
Ax,x\right\rangle \right) ^{r}}{2}\geq \left( \frac{m+M}{2}\right) ^{r} 
\notag
\end{align}%
for each $x\in H$ with $\left\Vert x\right\Vert =1$ and $r>1.$

If $0<r<1$ then the inequalities reverse in (\ref{I.c.e.3.4}).

If $A$ is positive definite and $r>0,$ then%
\begin{align}
\frac{m^{-r}+M^{-r}}{2}& \geq \left\langle \frac{A^{-r}+\left( \left(
m+M\right) 1_{H}-A\right) ^{-r}}{2}x,x\right\rangle  \label{I.c.e.3.5} \\
& \geq \frac{\left\langle Ax,x\right\rangle ^{-r}+\left( m+M-\left\langle
Ax,x\right\rangle \right) ^{-r}}{2}\geq \left( \frac{m+M}{2}\right) ^{-r} 
\notag
\end{align}%
for each $x\in H$ with $\left\Vert x\right\Vert =1.$
\end{proposition}

The proof follows by the inequality (\ref{I.c.e.2.1}).

Finally we have:

\begin{proposition}
\label{I.c.p.3.4}Assume that $A$ and $B$ are selfadjoint operators with
spectra in $\left[ m,M\right] \subset \lbrack 0,\infty )$ and $x\in H$ with $%
\left\Vert x\right\Vert =1$ and such that $\left\langle Ax,x\right\rangle
\neq \left\langle Bx,x\right\rangle .$

If $r>1$ or $r\in \left( \infty ,-1\right) \cup \left( -1,0\right) $ then we
have%
\begin{align}
\left\langle \left( \frac{A+B}{2}\right) x,x\right\rangle ^{r}& \leq \frac{1%
}{r+1}\cdot \frac{\left\langle Ax,x\right\rangle ^{r+1}-\left\langle
Bx,x\right\rangle ^{r+1}}{\left\langle Ax,x\right\rangle -\left\langle
Bx,x\right\rangle }  \label{I.c.e.3.6} \\
& \leq \left\langle \left[ \int_{0}^{1}\left( \left( 1-t\right) A+tB\right)
^{r}dt\right] x,x\right\rangle  \notag \\
& \leq \frac{M-\left\langle \frac{A+B}{2}x,x\right\rangle }{M-m}m^{r}+\frac{%
\left\langle \frac{A+B}{2}x,x\right\rangle -m}{M-m}M^{r}.  \notag
\end{align}

If $0<r<1,$ then the inequalities reverse in (\ref{I.c.e.3.6}).

If $A$ and $B$ are positive definite, then%
\begin{align}
\left\langle \left( \frac{A+B}{2}\right) x,x\right\rangle ^{-1}& \leq \frac{%
\ln \left\langle Bx,x\right\rangle -\ln \left\langle Ax,x\right\rangle }{%
\left\langle Bx,x\right\rangle -\left\langle Ax,x\right\rangle }
\label{I.c.e.3.7} \\
& \leq \left\langle \left[ \int_{0}^{1}\left( \left( 1-t\right) A+tB\right)
^{-1}dt\right] x,x\right\rangle  \notag \\
& \leq \frac{M-\left\langle \frac{A+B}{2}x,x\right\rangle }{\left(
M-m\right) m}+\frac{\left\langle \frac{A+B}{2}x,x\right\rangle -m}{\left(
M-m\right) M}.  \notag
\end{align}
\end{proposition}

\subsection{Applications for Ky Fan's Inequality}

The following results related to the Ky Fan inequality may be stated as well:

\begin{proposition}
\label{I.c.p.4.1}Let $A$ be a selfadjoint positive operator on a Hilbert
space $H.$ If $A$ is invertible and $Sp\left( A\right) \subset \left( 0,%
\frac{1}{2}\right) ,$ then%
\begin{align}
& \left( \left\langle \left( 1_{H}-A\right) x,x\right\rangle \left\langle
Ax,x\right\rangle ^{-1}\right) ^{r}  \label{I.c.e.4.2} \\
& \leq \exp \left\langle \ln \left( \left[ 1_{H}-A+\left\langle \left(
1_{H}-A\right) x,x\right\rangle 1_{H}\right] \left( A+\left\langle
Ax,x\right\rangle 1_{H}\right) ^{-1}\right) ^{r}x,x\right\rangle  \notag \\
& \leq \left\langle \exp \left[ \int_{0}^{1}\left[ \ln \left( \left(
1-t\right) \left( 1_{H}-A\right) +t\left\langle \left( 1_{H}-A\right)
x,x\right\rangle 1_{H}\right) \right. \right. \right.  \notag \\
& \left. \left. \times \left. \left( \left( 1-t\right) A+t\left\langle
Ax,x\right\rangle 1_{H}\right) ^{-1}\right] ^{r}dt\right] x,x\right\rangle 
\notag \\
& \leq \exp \left[ \frac{1}{2}\left[ \left\langle \ln \left[ \left(
1_{H}-A\right) A^{-1}\right] ^{r}x,x\right\rangle +\ln \left( \left\langle
\left( 1_{H}-A\right) x,x\right\rangle \left\langle Ax,x\right\rangle
^{-1}\right) ^{r}\right] \right]  \notag \\
& \leq \exp \left\langle \ln \left[ \left( 1_{H}-A\right) A^{-1}\right]
^{r}x,x\right\rangle  \notag
\end{align}%
for any $x\in H$ with $\left\Vert x\right\Vert =1.$
\end{proposition}

It follows from the inequality (\ref{I.c.e.2.23}) applied for the log-convex
function $g:\left( 0,1\right) \rightarrow \mathbb{R}$, $g\left( t\right)
=\left( \frac{1-t}{t}\right) ^{r},r>0.$

\begin{proposition}
\label{I.c.p.4.2}Assume that $A$ is a selfadjoint operator with $Sp\left(
A\right) \subset \left( 0,\frac{1}{2}\right) $ and $s\in \left( 0,\frac{1}{2}%
\right) .$ Then we have the following inequality in the operator order of $%
B\left( H\right) $:%
\begin{align}
& \ln \left[ \left[ \left( 2-s\right) 1_{H}-A\right] \left( A+s1_{H}\right)
^{-1}\right]  \label{I.c.e.4.3} \\
& \leq \int_{0}^{1}\ln \left( \left[ \left( 1-ts\right) 1_{H}-\left(
1-t\right) A\right] \left( \left( 1-t\right) A+ts1_{H}\right) ^{-1}\right) dt
\notag \\
& \leq \frac{1}{2}\left( \ln \left[ \left( 1_{H}-A\right) A^{-1}\right]
^{r}+\ln \left( \frac{1-s}{s}\right) ^{r}1_{H}\right) .  \notag
\end{align}
\end{proposition}

If follows from the inequality (\ref{I.c.e.2.20.b}) applied for the
log-convex function $g:\left( 0,1\right) \rightarrow \mathbb{R}$, $g\left(
t\right) =\left( \frac{1-t}{t}\right) ^{r},r>0.$

\section{Hermite-Hadamard's Type Inequalities for Operator Convex Functions}

\subsection{Introduction}

The following inequality holds for any convex function $f$ defined on $%
\mathbb{R}$ 
\begin{equation}
(b-a)f\left( \frac{a+b}{2}\right) <\int_{a}^{b}f(x)dx<(b-a)\frac{f(a)+f(b)}{2%
},\quad a,b\in \mathbb{R}.  \label{I.d.H}
\end{equation}%
It was firstly discovered by Ch. Hermite in 1881 in the journal \textit{%
Mathesis} (see \cite{ML}). But this result was nowhere mentioned in the
mathematical literature and was not widely known as Hermite's result \cite%
{I.d.PPT}.

E.F. Beckenbach, a leading expert on the history and the theory of convex
functions, wrote that this inequality was proven by J. Hadamard in 1893 \cite%
{I.a.BB}. In 1974, D.S. Mitrinovi\'{c} found Hermite's note in \textit{%
Mathesis} \cite{ML}. Since (\ref{I.d.H}) was known as Hadamard's inequality,
the inequality is now commonly referred as the Hermite-Hadamard inequality 
\cite{I.d.PPT}.

Let $X$ be a vector space, $x,y\in X,\ x\neq y$. Define the segment 
\begin{equation*}
\lbrack x,y]:=\{(1-t)x+ty,\ t\in \lbrack 0,1]\}.
\end{equation*}
We consider the function $f:[x,y]\rightarrow \mathbb{R}$ and the associated
function 
\begin{equation*}
g(x,y):[0,1]\rightarrow \mathbb{R},\ g(x,y)(t):=f[(1-t)x+ty],\ t\in \lbrack
0,1].
\end{equation*}
Note that $f$ is convex on $[x,y]$ if and only if $g(x,y)$ is convex on $%
[0,1]$.

For any convex function defined on a segment $[x.y]\subset X$, we have the
Hermite-Hadamard integral inequality (see \cite[p. 2]{I.d.SSD2.a}) 
\begin{equation}
f\left( \frac{x+y}{2}\right) \leq \int_{0}^{1}f[(1-t)x+ty]dt\leq \frac{%
f(x)+f(y)}{2},  \label{I.d.HH}
\end{equation}%
which can be derived from the classical Hermite-Hadamard inequality (\ref%
{I.d.H}) for the convex function $g(x,y):[0,1]\rightarrow \mathbb{R}$.

Since $f(x)=\Vert x\Vert ^{p}$ ($x\in X$ and $1\leq p<\infty $) is a convex
function, we have the following norm inequality from (\ref{I.d.HH}) (see 
\cite[p. 106]{I.d.PD}) 
\begin{equation}
\left\Vert \frac{x+y}{2}\right\Vert ^{p}\leq \int_{0}^{1}\Vert
(1-t)x+ty\Vert ^{p}dt\leq \frac{\Vert x\Vert ^{p}+\Vert y\Vert ^{p}}{2},
\label{I.d.HH1}
\end{equation}%
for any $x,y\in X$.

Motivated by the above results we investigate in this paper the operator
version of the Hermite-Hadamard inequality for operator convex functions.
The operator quasilinearity of some associated functionals are also provided.

A real valued continuous function $f$ on an interval $I$ is said to be 
\textit{operator convex (operator concave)} if 
\begin{equation}
f\left( \left( 1-\lambda \right) A+\lambda B\right) \leq \left( \geq \right)
\left( 1-\lambda \right) f\left( A\right) +\lambda f\left( B\right)  \tag{OC}
\label{I.d.OC}
\end{equation}%
in the operator order, for all $\lambda \in \left[ 0,1\right] $ and for
every selfadjoint operator $A$ and $B$ on a Hilbert space $H$ whose spectra
are contained in $I.$ Notice that a function $f$ is operator concave if $-f$
is operator convex.

A real valued continuous function $f$ on an interval $I$ is said to be 
\textit{operator monotone} if it is monotone with respect to the operator
order, i.e., $A\leq B$ with $Sp\left( A\right) ,Sp\left( B\right) \subset I$
imply $f\left( A\right) \leq f\left( B\right) .$

For some fundamental results on operator convex (operator concave) and
operator monotone functions, see \cite{I.FMPS} and the references therein.

As examples of such functions, we note that $f\left( t\right) =t^{r}$ is
operator monotone on $[0,\infty )$ if and only if $0\leq r\leq 1.$ The
function $f\left( t\right) =t^{r}$ is operator convex on $(0,\infty )$ if
either $1\leq r\leq 2$ or $-1\leq r\leq 0$ and is operator concave on $%
\left( 0,\infty \right) $ if $0\leq r\leq 1.$ The logarithmic function $%
f\left( t\right) =\ln t$ is operator monotone and operator concave on $%
(0,\infty ).$ The entropy function $f\left( t\right) =-t\ln t$ is operator
concave on $(0,\infty ).$ The exponential function$f\left( t\right) =e^{t}$
is neither operator convex nor operator monotone.

\subsection{Some Hermite-Hadamard's Type Inequalities}

We start with the following result:

\begin{theorem}[Dragomir, 2010, \protect\cite{I.d.SSD4}]
\label{I.d.t.2.1}Let $f:I\rightarrow \mathbb{R}$ be an operator convex
function on the interval $I.$ Then for any selfadjoint operators $A$ and $B$
with spectra in $I$ we have the inequality 
\begin{align}
& \left( f\left( \frac{A+B}{2}\right) \leq \right) \frac{1}{2}\left[ f\left( 
\frac{3A+B}{4}\right) +f\left( \frac{A+3B}{4}\right) \right]
\label{I.d.e.2.1} \\
& \leq \int_{0}^{1}f\left( \left( 1-t\right) A+tB\right) dt  \notag \\
& \leq \frac{1}{2}\left[ f\left( \frac{A+B}{2}\right) +\frac{f\left(
A\right) +f\left( B\right) }{2}\right] \left( \leq \frac{f\left( A\right)
+f\left( B\right) }{2}\right) .  \notag
\end{align}
\end{theorem}

\begin{proof}
First of all, since the function $f$ is continuos, the operator valued
integral $\int_{0}^{1}f\left( \left( 1-t\right) A+tB\right) dt$ exists for
any selfadjoint operators $A$ and $B$ with spectra in $I.$

We give here two proofs, the first using only the definition of operator
convex functions and the second using the classical Hermite-Hadamard
inequality for real valued functions.

1. By the definition of operator convex functions we have the double
inequality: 
\begin{align}
f\left( \frac{C+D}{2}\right) & \leq \frac{1}{2}\left[ f\left( \left(
1-t\right) C+tD\right) +f\left( \left( 1-t\right) D+tC\right) \right]
\label{I.d.e.2.2} \\
& \leq \frac{1}{2}\left[ f\left( C\right) +f\left( D\right) \right]  \notag
\end{align}%
for any $t\in \left[ 0,1\right] $ and any selfadjoint operators $C$ and $D$
with the spectra in $I.$

Integrating the inequality (\ref{I.d.e.2.2}) over $t\in \left[ 0,1\right] $
and taking into account that 
\begin{equation*}
\int_{0}^{1}f\left( \left( 1-t\right) C+tD\right) dt=\int_{0}^{1}f\left(
\left( 1-t\right) D+tC\right) dt
\end{equation*}%
then we deduce the Hermite-Hadamard inequality for operator convex functions 
\begin{equation}
f\left( \frac{C+D}{2}\right) \leq \int_{0}^{1}f\left( \left( 1-t\right)
C+tD\right) dt\leq \frac{1}{2}\left[ f\left( C\right) +f\left( D\right) %
\right]  \tag{HHO}  \label{I.d.HHO}
\end{equation}%
that holds for any selfadjoint operators $C$ and $D$ with the spectra in $I.$

Now, on making use of the change of variable $u=2t$ we have 
\begin{equation*}
\int_{0}^{1/2}f\left( \left( 1-t\right) A+tB\right) dt=\frac{1}{2}%
\int_{0}^{1}f\left( \left( 1-u\right) A+u\frac{A+B}{2}\right) du
\end{equation*}
and by the change of variable $u=2t-1$ we have 
\begin{equation*}
\int_{1/2}^{1}f\left( \left( 1-t\right) A+tB\right) dt=\frac{1}{2}%
\int_{0}^{1}f\left( \left( 1-u\right) \frac{A+B}{2}+uB\right) du.
\end{equation*}

Utilising the Hermite-Hadamard inequality (\ref{I.d.HHO}) we can write 
\begin{align*}
f\left( \frac{3A+B}{4}\right) & \leq \int_{0}^{1}f\left( \left( 1-u\right)
A+u\frac{A+B}{2}\right) du \\
& \leq \frac{1}{2}\left[ f\left( A\right) +f\left( \frac{A+B}{2}\right) %
\right]
\end{align*}%
and 
\begin{align*}
f\left( \frac{A+3B}{4}\right) & \leq \int_{0}^{1}f\left( \left( 1-u\right) 
\frac{A+B}{2}+uB\right) du \\
& \leq \frac{1}{2}\left[ f\left( A\right) +f\left( \frac{A+B}{2}\right) %
\right] ,
\end{align*}%
which by summation and division by two produces the desired result (\ref%
{I.d.e.2.1}).

2. Consider now $x\in H,\left\| x\right\| =1$ and two selfadjoint operators $%
A$ and $B$ with spectra in $I$. Define the real-valued function $\varphi
_{x,A,B}:\left[ 0,1\right] \rightarrow \mathbb{R}$ given by $\varphi
_{x,A,B}\left( t\right) =\left\langle f\left( \left( 1-t\right) A+tB\right)
x,x\right\rangle .$

Since $f$ is operator convex, then for any $t_{1},t_{2}\in \left[ 0,1\right] 
$ and $\alpha ,\beta \geq 0$ with $\alpha +\beta =1$ we have 
\begin{align*}
& \varphi _{x,A,B}\left( \alpha t_{1}+\beta t_{2}\right) \\
& =\left\langle f\left( \left( 1-\left( \alpha t_{1}+\beta t_{2}\right)
\right) A+\left( \alpha t_{1}+\beta t_{2}\right) B\right) x,x\right\rangle \\
& =\left\langle f\left( \alpha \left[ \left( 1-t_{1}\right) A+t_{1}B\right]
+\beta \left[ \left( 1-t_{2}\right) A+t_{2}B\right] \right) x,x\right\rangle
\\
& \leq \alpha \left\langle f\left( \left[ \left( 1-t_{1}\right) A+t_{1}B%
\right] \right) x,x\right\rangle +\beta \left\langle f\left( \left[ \left(
1-t_{2}\right) A+t_{2}B\right] \right) x,x\right\rangle \\
& =\alpha \varphi _{x,A,B}\left( t_{1}\right) +\beta \varphi _{x,A,B}\left(
t_{2}\right)
\end{align*}%
showing that $\varphi _{x,A,B}$ is a convex function on $\left[ 0,1\right] .$

Now we use the Hermite-Hadamard inequality for real-valued convex functions 
\begin{equation*}
g\left( \frac{a+b}{2}\right) \leq \frac{1}{b-a}\int_{a}^{b}g\left( s\right)
ds\leq \frac{g\left( a\right) +g\left( b\right) }{2}
\end{equation*}%
to get that 
\begin{equation*}
\varphi _{x,A,B}\left( \frac{1}{4}\right) \leq 2\int_{0}^{1/2}\varphi
_{x,A,B}\left( t\right) dt\leq \frac{\varphi _{x,A,B}\left( 0\right)
+\varphi _{x,A,B}\left( \frac{1}{2}\right) }{2}
\end{equation*}%
and 
\begin{equation*}
\varphi _{x,A,B}\left( \frac{3}{4}\right) \leq 2\int_{1/2}^{1}\varphi
_{x,A,B}\left( t\right) dt\leq \frac{\varphi _{x,A,B}\left( \frac{1}{2}%
\right) +\varphi _{x,A,B}\left( 1\right) }{2}
\end{equation*}%
which by summation and division by two produces 
\begin{align}
& \frac{1}{2}\left\langle \left[ f\left( \frac{3A+B}{4}\right) +f\left( 
\frac{A+3B}{4}\right) \right] x,x\right\rangle  \label{I.d.e.2.3} \\
& \leq \int_{0}^{1}\left\langle f\left( \left( 1-t\right) A+tB\right)
x,x\right\rangle dt  \notag \\
& \leq \frac{1}{2}\left\langle \left[ f\left( \frac{A+B}{2}\right) +\frac{%
f\left( A\right) +f\left( B\right) }{2}\right] x,x\right\rangle .  \notag
\end{align}%
Finally, since by the continuity of the function $f$ we have 
\begin{equation*}
\int_{0}^{1}\left\langle f\left( \left( 1-t\right) A+tB\right)
x,x\right\rangle dt=\left\langle \int_{0}^{1}f\left( \left( 1-t\right)
A+tB\right) dtx,x\right\rangle
\end{equation*}%
for any $x\in H,\left\Vert x\right\Vert =1$ and any two selfadjoint
operators $A$ and $B$ with spectra in $I,$ we deduce from (\ref{I.d.e.2.3})
the desired result (\ref{I.d.e.2.1}).
\end{proof}

A simple consequence of the above theorem is that the integral is closer to
the left bound than to the right, namely we can state:

\begin{corollary}[Dragomir, 2010, \protect\cite{I.d.SSD4}]
\label{I.d.c2.1.}With the assumptions in Theorem \ref{I.d.t.2.1} we have the
inequality 
\begin{align}
& \left( 0\leq \right) \int_{0}^{1}f\left( \left( 1-t\right) A+tB\right)
dt-f\left( \frac{A+B}{2}\right)  \label{I.d.e.2.4} \\
& \leq \frac{f\left( A\right) +f\left( B\right) }{2}-\int_{0}^{1}f\left(
\left( 1-t\right) A+tB\right) dt.  \notag
\end{align}
\end{corollary}

\begin{remark}
\label{I.d.r.2.1}Utilising different examples of operator convex or concave
functions, we can provide inequalities of interest.

If $r\in \left[ -1,0\right] \cup \left[ 1,2\right] $ then we have the
inequalities for powers of operators 
\begin{align}
& \left( \left( \frac{A+B}{2}\right) ^{r}\leq \right) \frac{1}{2}\left[
\left( \frac{3A+B}{4}\right) ^{r}+\left( \frac{A+3B}{4}\right) ^{r}\right]
\label{I.d.e.2.5} \\
& \leq \int_{0}^{1}\left( \left( 1-t\right) A+tB\right) ^{r}dt  \notag \\
& \leq \frac{1}{2}\left[ \left( \frac{A+B}{2}\right) ^{r}+\frac{A^{r}+B^{r}}{%
2}\right] \left( \leq \frac{A^{r}+B^{r}}{2}\right)  \notag
\end{align}%
for any two selfadjoint operators $A$ and $B$ with spectra in $\left(
0,\infty \right) .$

If $r\in \left( 0,1\right) $ the inequalities in (\ref{I.d.e.2.5}) hold with 
$"\geq "$ instead of $"\leq ".$

We also have the following inequalities for logarithm 
\begin{align}
& \left( \ln \left( \frac{A+B}{2}\right) \geq \right) \frac{1}{2}\left[ \ln
\left( \frac{3A+B}{4}\right) +\ln \left( \frac{A+3B}{4}\right) \right]
\label{I.d.e.2.6} \\
& \geq \int_{0}^{1}\ln \left( \left( 1-t\right) A+tB\right) dt  \notag \\
& \geq \frac{1}{2}\left[ \ln \left( \frac{A+B}{2}\right) +\frac{\ln \left(
A\right) +\ln \left( B\right) }{2}\right] \left( \geq \frac{\ln \left(
A\right) +\ln \left( B\right) }{2}\right)  \notag
\end{align}%
for any two selfadjoint operators $A$ and $B$ with spectra in $\left(
0,\infty \right) .$
\end{remark}

\subsection{Some Operator Quasilinearity Properties}

Consider an operator convex function $f:I\subset \mathbb{R}\rightarrow 
\mathbb{R}$ defined on the interval $I$ and two distinct selfadjoint
operators $A,B$ with the spectra in $I$. We denote by $\left[ A,B\right] $
the closed operator segment defined by the family of operators $\left\{
\left( 1-t\right) A+tB\text{, }t\in \left[ 0,1\right] \right\} .$ We also
define the operator-valued \ functional 
\begin{equation}
\Delta _{f}\left( A,B;t\right) :=\left( 1-t\right) f\left( A\right)
+tf\left( B\right) -f\left( \left( 1-t\right) A+tB\right) \geq 0
\label{I.d.eq.3.1}
\end{equation}%
in the operator order, for any $t\in \left[ 0,1\right] .$

The following result concerning an operator quasilinearity property for the
functional $\Delta _{f}\left( \cdot ,\cdot ;t\right) $ may be stated:

\begin{theorem}[Dragomir, 2010, \protect\cite{I.d.SSD4}]
\label{I.d.t.3.1}Let $f:I\subset \mathbb{R}\rightarrow \mathbb{R}$ be an
operator convex function on the interval $I$. Then for each $A,B$ two
distinct selfadjoint operators $A,B$ with the spectra in $I$ and $C\in $ $%
\left[ A,B\right] $ we have 
\begin{equation}
\left( 0\leq \right) \Delta _{f}\left( A,C;t\right) +\Delta _{f}\left(
C,B;t\right) \leq \Delta _{f}\left( A,B;t\right)  \label{I.d.eq.3.2}
\end{equation}%
for each $t\in \left[ 0,1\right] ,$ i.e., the functional $\Delta _{f}\left(
\cdot ,\cdot ;t\right) $ is operator superadditive as a function of interval.

If $\left[ C,D\right] \subset \left[ A,B\right] ,$ then 
\begin{equation}
\left( 0\leq \right) \Delta _{f}\left( C,D;t\right) \leq \Delta _{f}\left(
A,B;t\right)  \label{I.d.eq.3.3}
\end{equation}%
for each $t\in \left[ 0,1\right] ,$ i.e., the functional $\Delta _{f}\left(
\cdot ,\cdot ;t\right) $ is operator nondecreasing as a function of interval.
\end{theorem}

\begin{proof}
Let $C=\left( 1-s\right) A+sB$ with $s\in \left( 0,1\right) .$ For $t\in
\left( 0,1\right) $ we have 
\begin{align*}
\Delta _{f}\left( C,B;t\right) & =\left( 1-t\right) f\left( \left(
1-s\right) A+sB\right) +tf\left( B\right) \\
& -f\left( \left( 1-t\right) \left[ \left( 1-s\right) A+sB\right] +tB\right)
\end{align*}%
and 
\begin{align*}
\Delta _{f}\left( A,C;t\right) & =\left( 1-t\right) f\left( A\right)
+tf\left( \left( 1-s\right) A+sB\right) \\
& -f\left( \left( 1-t\right) A+t\left[ \left( 1-s\right) A+sB\right] \right)
\end{align*}%
giving that 
\begin{align}
& \Delta _{f}\left( A,C;t\right) +\Delta _{f}\left( C,B;t\right) -\Delta
_{f}\left( A,B;t\right)  \label{I.d.eq.3.4} \\
& =f\left( \left( 1-s\right) A+sB\right) +f\left( \left( 1-t\right)
A+tB\right)  \notag \\
& -f\left( \left( 1-t\right) \left( 1-s\right) A+\left[ \left( 1-t\right) s+t%
\right] B\right) -f\left( \left( 1-ts\right) A+tsB\right) .  \notag
\end{align}

Now, for a convex function $\varphi :I\subset \mathbb{R}\rightarrow \mathbb{R%
}$, where $I$ is an interval, and any real numbers $t_{1},t_{2},s_{1}$ and $%
s_{2}$ from $I$ and with the properties that $t_{1}\leq s_{1}$ and $%
t_{2}\leq s_{2}$ we have that 
\begin{equation}
\frac{\varphi \left( t_{1}\right) -\varphi \left( t_{2}\right) }{t_{1}-t_{2}}%
\leq \frac{\varphi \left( s_{1}\right) -\varphi \left( s_{2}\right) }{%
s_{1}-s_{2}}.  \label{I.d.eq.3.5}
\end{equation}%
Indeed, since $\varphi $ is convex on $I$ then for any $a\in I$ the function 
$\psi :I\backslash \left\{ a\right\} \rightarrow \mathbb{R}$ 
\begin{equation*}
\psi \left( t\right) :=\frac{\varphi \left( t\right) -\varphi \left(
a\right) }{t-a}
\end{equation*}%
is monotonic nondecreasing where is defined. Utilising this property
repeatedly we have 
\begin{eqnarray*}
\frac{\varphi \left( t_{1}\right) -\varphi \left( t_{2}\right) }{t_{1}-t_{2}}
&\leq &\frac{\varphi \left( s_{1}\right) -\varphi \left( t_{2}\right) }{%
s_{1}-t_{2}}=\frac{\varphi \left( t_{2}\right) -\varphi \left( s_{1}\right) 
}{t_{2}-s_{1}} \\
&\leq &\frac{\varphi \left( s_{2}\right) -\varphi \left( s_{1}\right) }{%
s_{2}-s_{1}}=\frac{\varphi \left( s_{1}\right) -\varphi \left( s_{2}\right) 
}{s_{1}-s_{2}}
\end{eqnarray*}%
which proves the inequality (\ref{I.d.eq.3.5}).

For a vector $x\in H$, with $\left\Vert x\right\Vert =1,$ consider the
function $\varphi _{x}:\left[ 0,1\right] \rightarrow \mathbb{R}$ given by $%
\varphi _{x}\left( t\right) :=\left\langle f\left( \left( 1-t\right)
A+tB\right) x,x\right\rangle .$ Since $f$ is operator convex on $I$ it
follows that $\varphi _{x}$ is convex on $\left[ 0,1\right] .$ Now, if we
consider, for given $t,s\in \left( 0,1\right) ,$ 
\begin{equation*}
t_{1}:=ts<s=:s_{1}\text{ and }t_{2}:=t<t+\left( 1-t\right) s=:s_{2},
\end{equation*}%
then we have $\varphi _{x}\left( t_{1}\right) =\left\langle f\left( \left(
1-ts\right) A+tsB\right) x,x\right\rangle $ and $\varphi _{x}\left(
t_{2}\right) =\left\langle f\left( \left( 1-t\right) A+tB\right)
x,x\right\rangle $ giving that 
\begin{equation*}
\frac{\varphi _{x}\left( t_{1}\right) -\varphi _{x}\left( t_{2}\right) }{%
t_{1}-t_{2}}=\left\langle \left[ \frac{f\left( \left( 1-ts\right)
A+tsB\right) -f\left( \left( 1-t\right) A+tB\right) }{t\left( s-1\right) }%
\right] x,x\right\rangle .
\end{equation*}%
Also $\varphi _{x}\left( s_{1}\right) =\left\langle f\left( \left(
1-s\right) A+sB\right) x,x\right\rangle $ and $\varphi _{x}\left(
s_{2}\right) =\left\langle f\left( \left( 1-t\right) \left( 1-s\right) A+%
\left[ \left( 1-t\right) s+t\right] B\right) x,x\right\rangle $ giving that 
\begin{eqnarray*}
&&\frac{\varphi _{x}\left( s_{1}\right) -\varphi _{x}\left( s_{2}\right) }{%
s_{1}-s_{2}} \\
&=&\left\langle \frac{f\left( \left( 1-s\right) A+sB\right) -f\left( \left(
1-t\right) \left( 1-s\right) A+\left[ \left( 1-t\right) s+t\right] B\right) 
}{t\left( s-1\right) }x,x\right\rangle .
\end{eqnarray*}%
Utilising the inequality (\ref{I.d.eq.3.5}) and multiplying with $t\left(
s-1\right) <0$ we deduce the following inequality in the operator order 
\begin{align}
& f\left( \left( 1-ts\right) A+tsB\right) -f\left( \left( 1-t\right)
A+tB\right)  \label{I.d.eq.3.6} \\
& \geq f\left( \left( 1-s\right) A+sB\right) -f\left( \left( 1-t\right)
\left( 1-s\right) A+\left[ \left( 1-t\right) s+t\right] B\right) .  \notag
\end{align}%
Finally, by (\ref{I.d.eq.3.4}) and (\ref{I.d.eq.3.6}) we get the desired
result (\ref{I.d.eq.3.2}).

Applying repeatedly the superadditivity property we have for $\left[ C,D%
\right] \subset \left[ A,B\right] $ that 
\begin{equation*}
\Delta _{f}\left( A,C;t\right) +\Delta _{f}\left( C,D;t\right) +\Delta
_{f}\left( D,B;t\right) \leq \Delta _{f}\left( A,B;t\right)
\end{equation*}%
giving that 
\begin{equation*}
0\leq \Delta _{f}\left( A,C;t\right) +\Delta _{f}\left( D,B;t\right) \leq
\Delta _{f}\left( A,B;t\right) -\Delta _{f}\left( C,D;t\right)
\end{equation*}%
which proves (\ref{I.d.eq.3.3}).
\end{proof}

For $t=\frac{1}{2}$ we consider the functional 
\begin{equation*}
\Delta _{f}\left( A,B\right) :=\Delta _{f}\left( A,B;\frac{1}{2}\right) =%
\frac{f\left( A\right) +f\left( B\right) }{2}-f\left( \frac{A+B}{2}\right) ,
\end{equation*}
which obviously inherits the superadditivity and monotonicity properties of
the functional $\Delta _{f}\left( \cdot ,\cdot ;t\right) .$ We are able then
to state the following

\begin{corollary}[Dragomir, 2010, \protect\cite{I.d.SSD4}]
\label{I.d.c.3.1}Let $f:I\subset \mathbb{R}\rightarrow \mathbb{R}$ be an
operator convex function on the interval $I$. Then for each $A,B$ two
distinct selfadjoint operators $A,B$ with the spectra in $I$ we have the
following bounds\ in the operator order 
\begin{equation}
\inf_{C\in \left[ A,B\right] }\left[ f\left( \frac{A+C}{2}\right) +f\left( 
\frac{C+B}{2}\right) -f\left( C\right) \right] =f\left( \frac{A+B}{2}\right)
\label{I.d.eq.3.7}
\end{equation}%
and 
\begin{equation}
\sup_{C,D\in \left[ A,B\right] }\left[ \frac{f\left( C\right) +f\left(
D\right) }{2}-f\left( \frac{C+D}{2}\right) \right] =\frac{f\left( A\right)
+f\left( B\right) }{2}-f\left( \frac{A+B}{2}\right) .  \label{I.d.eq.3.8}
\end{equation}
\end{corollary}

\begin{proof}
By the superadditivity of the functional $\Delta _{f}\left( \cdot ,\cdot
\right) $ we have for each $C\in \left[ A,B\right] $ that 
\begin{align*}
& \frac{f\left( A\right) +f\left( B\right) }{2}-f\left( \frac{A+B}{2}\right)
\\
& \geq \frac{f\left( A\right) +f\left( C\right) }{2}-f\left( \frac{A+C}{2}%
\right) +\frac{f\left( C\right) +f\left( B\right) }{2}-f\left( \frac{C+B}{2}%
\right)
\end{align*}%
which is equivalent with 
\begin{equation}
f\left( \frac{A+C}{2}\right) +f\left( \frac{C+B}{2}\right) -f\left( C\right)
\geq f\left( \frac{A+B}{2}\right) .  \label{I.d.eq.3.9}
\end{equation}%
Since the equality case in (\ref{I.d.eq.3.9}) is realized for either $C=A$
or $C=B$ we get the desired bound (\ref{I.d.eq.3.7}).

The bound (\ref{I.d.eq.3.8}) is obvious by the monotonicity of the
functional $\Delta _{f}\left( \cdot ,\cdot \right) $ as a function of
interval.
\end{proof}

Consider now the following functional 
\begin{equation*}
\Gamma _{f}\left( A,B;t\right) :=f\left( A\right) +f\left( B\right) -f\left(
\left( 1-t\right) A+tB\right) -f\left( \left( 1-t\right) B+tA\right) ,
\end{equation*}
where, as above, $f:C\subset X\rightarrow \mathbb{R}$ is a convex function
on the convex set $C$ and $A,B\in C$ while $t\in \left[ 0,1\right] .$

We notice that 
\begin{equation*}
\Gamma _{f}\left( A,B;t\right) =\Gamma _{f}\left( B,A;t\right) =\Gamma
_{f}\left( A,B;1-t\right)
\end{equation*}
and 
\begin{equation*}
\Gamma _{f}\left( A,B;t\right) =\Delta _{f}\left( A,B;t\right) +\Delta
_{f}\left( A,B;1-t\right) \geq 0
\end{equation*}
for any $A,B\in C$ and $t\in \left[ 0,1\right] .$

Therefore, we can state the following result as well

\begin{corollary}[Dragomir, 2010, \protect\cite{I.d.SSD4}]
\label{I.d.c.3.2}Let $f:I\subset \mathbb{R}\rightarrow \mathbb{R}$ be an
operator convex function on the interval $I$. Then for each $A,B$ two
distinct selfadjoint operators $A,B$ with the spectra in $I,$ the functional 
$\Gamma _{f}\left( \cdot ,\cdot ;t\right) $ is operator superadditive and
operator nondecreasing as a function of interval.
\end{corollary}

In particular, if $C\in \left[ A,B\right] $ then we have the inequality 
\begin{align}
& \frac{1}{2}\left[ f\left( \left( 1-t\right) A+tB\right) +f\left( \left(
1-t\right) B+tA\right) \right]  \label{I.d.eq.3.9.a} \\
& \leq \frac{1}{2}\left[ f\left( \left( 1-t\right) A+tC\right) +f\left(
\left( 1-t\right) C+tA\right) \right]  \notag \\
& +\frac{1}{2}\left[ f\left( \left( 1-t\right) C+tB\right) +f\left( \left(
1-t\right) B+tC\right) \right] -f\left( C\right) .  \notag
\end{align}

Also, if $C,D\in \left[ A,B\right] $ then we have the inequality 
\begin{align}
& f\left( A\right) +f\left( B\right) -f\left( \left( 1-t\right) A+tB\right)
-f\left( \left( 1-t\right) B+tA\right)  \label{I.d.eq.3.10} \\
& \geq f\left( C\right) +f\left( D\right) -f\left( \left( 1-t\right)
C+tD\right) -f\left( \left( 1-t\right) C+tD\right)  \notag
\end{align}%
for any $t\in \left[ 0,1\right] .$

Perhaps the most interesting functional we can consider is the following
one: 
\begin{equation}
\Theta _{f}\left( A,B\right) =\frac{f\left( A\right) +f\left( B\right) }{2}%
-\int_{0}^{1}f\left( \left( 1-t\right) A+tB\right) dt.  \label{I.d.eq.3.11}
\end{equation}%
Notice that, by the second Hermite-Hadamard inequality for operator convex
functions we have that $\Theta _{f}\left( A,B\right) \geq 0$ in the operator
order.

We also observe that 
\begin{equation}
\Theta _{f}\left( A,B\right) =\int_{0}^{1}\Delta _{f}\left( A,B;t\right)
dt=\int_{0}^{1}\Delta _{f}\left( A,B;1-t\right) dt.  \label{I.d.eq.3.12}
\end{equation}%
Utilising this representation, we can state the following result as well:

\begin{corollary}[Dragomir, 2010, \protect\cite{I.d.SSD4}]
\label{I.d.c.2.3}Let $f:I\subset \mathbb{R}\rightarrow \mathbb{R}$ be an
operator convex function on the interval $I$. Then for each $A,B$ two
distinct selfadjoint operators $A,B$ with the spectra in $I,$ the functional 
$\Theta _{f}\left( \cdot ,\cdot \right) $ is operator superadditive and
operator nondecreasing as a function of interval. Moreover, we have the
bounds in the operator order 
\begin{align}
& \inf_{C\in \left[ A,B\right] }\left[ \int_{0}^{1}\left[ f\left( \left(
1-t\right) A+tC\right) +f\left( \left( 1-t\right) C+tB\right) \right]
dt-f\left( C\right) \right]  \label{I.d.eq.3.13} \\
& =\int_{0}^{1}f\left( \left( 1-t\right) A+tB\right) dt  \notag
\end{align}%
and 
\begin{align}
& \sup_{C,D\in \left[ A,B\right] }\left[ \frac{f\left( C\right) +f\left(
D\right) }{2}-\int_{0}^{1}f\left( \left( 1-t\right) C+tD\right) dt\right]
\label{I.d.eq.3.14} \\
& =\frac{f\left( A\right) +f\left( B\right) }{2}-\int_{0}^{1}f\left( \left(
1-t\right) A+tB\right) dt.  \notag
\end{align}
\end{corollary}

\begin{remark}
\label{I.d.r.3.1}The above inequalities can be applied to various concrete
operator convex function of interest.

If we choose for instance the inequality (\ref{I.d.eq.3.14}), then we get
the following bounds in the operator order%
\begin{align}
& \sup_{C,D\in \left[ A,B\right] }\left[ \frac{C^{r}+D^{r}}{2}%
-\int_{0}^{1}\left( \left( 1-t\right) C+tD\right) ^{r}dt\right]
\label{I.d.eq.3.15} \\
& =\frac{A^{r}+B^{r}}{2}-\int_{0}^{1}\left( \left( 1-t\right) A+tB\right)
^{r}dt,  \notag
\end{align}%
where $r\in \left[ -1,0\right] \cup \left[ 1,2\right] $ and $A,B$ are
selfadjoint operators with spectra in $\left( 0,\infty \right) .$

If $r\in \left( 0,1\right) $ then%
\begin{align}
& \sup_{C,D\in \left[ A,B\right] }\left[ \int_{0}^{1}\left( \left(
1-t\right) C+tD\right) ^{r}dt-\frac{C^{r}+D^{r}}{2}\right]
\label{I.d.eq.3.16} \\
& =\int_{0}^{1}\left( \left( 1-t\right) A+tB\right) ^{r}dt-\frac{A^{r}+B^{r}%
}{2},  \notag
\end{align}%
and $A,B$ are selfadjoint operators with spectra in $\left( 0,\infty \right)
.$

We also have the operator bound for the logarithm%
\begin{align}
& \sup_{C,D\in \left[ A,B\right] }\left[ \int_{0}^{1}\ln \left( \left(
1-t\right) C+tD\right) dt-\frac{\ln \left( C\right) +\ln \left( D\right) }{2}%
\right]  \label{I.d.eq.3.17} \\
& =\int_{0}^{1}\ln \left( \left( 1-t\right) A+tB\right) dt-\frac{\ln \left(
A\right) +\ln \left( B\right) }{2},  \notag
\end{align}%
where $A,B$ are selfadjoint operators with spectra in $\left( 0,\infty
\right) .$
\end{remark}

\bigskip

\chapter{Inequalities for the \v{C}eby\v{s}ev Functional}

\section{Introduction}

The \textit{\v{C}eby\v{s}ev}, or in a different spelling, \textit{Chebyshev
inequality} which compares the integral/discrete mean of the product with
the product of the integral/discrete means is famous in the literature
devoted to Mathematical Inequalities. It has been extended, generalised,
refined etc...by many authors during the last century. A simple search
utilising either spellings and the key word "inequality" in the title in the
comprehensive \textit{MathSciNet} database of the \textit{American
Mathematical Society} produces more than 200 research articles devoted to
this result.

The sister result due to Gr\"{u}ss which provides error bounds for the
magnitude of the difference between the integral mean of the product and the
product of the integral means has also attracted much interest since it has
been discovered in 1935 with more than 180 papers published, as a simple
search in the same database reveals. Far more publications have been devoted
to the applications of these inequalities and an accurate picture of the
impacted results in various fields of Modern Mathematics is difficult to
provide.

In this chapter, however, we present only some recent results due to the
author for the corresponding operator versions of these two famous
inequalities. Applications for particular functions of selfadjoint operators
such as the power, logarithmic and exponential functions are provided as
well.

\section{\v{C}eby\v{s}ev's Inequality}

\subsection{\v{C}eby\v{s}ev's Inequality for Real Numbers}

First of all, let us recall a number of classical results for sequences of
real numbers concerning the celebrated \textit{\v{C}eby\v{s}ev inequality.}

Consider the real sequences $\left( n-\text{tuples}\right) $ $\mathbf{a}%
=\left( a_{1},\dots ,a_{n}\right) ,$ $\mathbf{b}=\left( b_{1},\dots
,b_{n}\right) $ and the nonnegative sequence $\mathbf{p}=\left( p_{1},\dots
,p_{n}\right) $ with $P_{n}:=\sum_{i=1}^{n}p_{i}>0.$ Define the \textit{%
weighted \v{C}eby\v{s}ev's functional} 
\begin{equation}
T_{n}\left( \mathbf{p};\mathbf{a},\mathbf{b}\right) :=\frac{1}{P_{n}}%
\sum_{i=1}^{n}p_{i}a_{i}b_{i}-\frac{1}{P_{n}}\sum_{i=1}^{n}p_{i}a_{i}\cdot 
\frac{1}{P_{n}}\sum_{i=1}^{n}p_{i}b_{i}.  \label{II.CFOa.1}
\end{equation}

In 1882 -- 1883, \v{C}eby\v{s}ev \cite{II.CE1} and \cite{II.CE2} proved that
if $\mathbf{a}$ and $\mathbf{b}$ are \textit{monotonic} in the same
(opposite) sense, then 
\begin{equation}
T_{n}\left( \mathbf{p};\mathbf{a},\mathbf{b}\right) \geq \left( \leq \right)
0.  \label{II.CFOa.2}
\end{equation}

In the special case $\mathbf{p}=\mathbf{a}\geq \mathbf{0}$, it appears that
the inequality (\ref{II.CFOa.2}) has been obtained by Laplace long before 
\v{C}eby\v{s}ev (see for example \cite[p. 240]{II.11b}).

The inequality (\ref{II.CFOa.2}) was mentioned by Hardy, Littlewood and P%
\'{o}lya in their book \cite{II.HLP} in 1934 in the more general setting of
synchronous sequences, i.e., if $\mathbf{a},$ $\mathbf{b}$ are \textit{%
synchronous} (\textit{asynchronous}), this means that 
\begin{equation}
\left( a_{i}-a_{j}\right) \left( b_{i}-b_{j}\right) \geq \left( \leq \right)
0\text{ for any }i,j\in \left\{ 1,\dots ,n\right\} ,  \label{II.CFOa.3}
\end{equation}%
then (\ref{II.CFOa.2}) holds true as well.

A relaxation of the synchronicity condition was provided by M. Biernacki in
1951, \cite{II.BIE}, which showed that, if $\mathbf{a},$ $\mathbf{b}$ are 
\textit{monotonic in mean} in the same sense, i.e., for $P_{k}:=%
\sum_{i=1}^{k}p_{i},$ $k=1,\dots ,n-1;$%
\begin{equation}
\frac{1}{P_{k}}\sum_{i=1}^{k}p_{i}a_{i}\leq \left( \geq \right) \frac{1}{%
P_{k+1}}\sum_{i=1}^{k+1}p_{i}a_{i},\;\;k\in \left\{ 1,\dots ,n-1\right\}
\label{II.CFOa.4}
\end{equation}%
and 
\begin{equation}
\frac{1}{P_{k}}\sum_{i=1}^{k}p_{i}b_{i}\leq \left( \geq \right) \frac{1}{%
P_{k+1}}\sum_{i=1}^{k+1}p_{i}b_{i},\;\;k\in \left\{ 1,\dots ,n-1\right\} ,
\label{II.CFOa.5}
\end{equation}%
then (\ref{II.CFOa.2}) holds with \textquotedblleft $\ \geq $ \
\textquotedblright . If if $\mathbf{a},$ $\mathbf{b}$ are monotonic in mean
in the opposite sense then (\ref{II.CFOa.2}) holds with \textquotedblleft $\
\leq $ \ \textquotedblright .

If one would like to drop the assumption of nonnegativity for the components
of $\mathbf{p},$ then one may state the following inequality obtained by
Mitrinovi\'{c} and Pe\v{c}ari\'{c} in 1991, \cite{II.MP}: If $0\leq
P_{i}\leq P_{n}$ for each $i\in \left\{ 1,\dots ,n-1\right\} ,$ then 
\begin{equation}
T_{n}\left( \mathbf{p};\mathbf{a},\mathbf{b}\right) \geq 0
\label{II.CFOa.11}
\end{equation}%
provided $\mathbf{a}$ and $\mathbf{b}$ are sequences with the same
monotonicity.

If $\mathbf{a}$ and $\mathbf{b}$ are monotonic in the opposite sense, the
sign of the inequality (\ref{II.CFOa.11}) reverses.

Similar integral inequalities may be stated, however we do not present them
here.

For other recent results on the \v{C}eby\v{s}ev inequality in either
discrete or integral form see \cite{II.BG}, \cite{II.cfo6b}, \cite{II.CFO8b}%
, \cite{II.SSD1}, \cite{II.DM}, \cite{II.CFO5b}, \cite{II.11b}, \cite{II.MP0}%
, \cite{II.MV}, \cite{II.cofPa2}, \cite{II.cofPa1}, \cite{II.Pe}, and the
references therein.

The main aim of the present section is to provide operator versions for the 
\v{C}eby\v{s}ev inequality in different settings. Related results and some
particular cases of interest are also given.

\subsection{A Version of the \v{C}eby\v{s}ev Inequality for One Operator}

We say that the functions $f,g:\left[ a,b\right] \longrightarrow \mathbb{R}$
are \textit{synchronous (asynchronous)} on the interval $\left[ a,b\right] $
if they satisfy the following condition:%
\begin{equation*}
\left( f\left( t\right) -f\left( s\right) \right) \left( g\left( t\right)
-g\left( s\right) \right) \geq \left( \leq \right) 0\text{ \quad for each }%
t,s\in \left[ a,b\right] .
\end{equation*}

It is obvious that, if $f,g$ are monotonic and have the same monotonicity on
the interval $\left[ a,b\right] ,$ then they are synchronous on $\left[ a,b%
\right] $ while if they have opposite monotonicity, they are asynchronous.

For some extensions of the discrete \textit{\v{C}eby\v{s}ev inequality }for%
\textit{\ synchronous (asynchronous) }sequences of vectors in an inner
product space, see \cite{II.DS} and \cite{II.DPS}.

The following result provides an inequality of \v{C}eby\v{s}ev type for
functions of selfadjoint operators.

\begin{theorem}[Dragomir, 2008, \protect\cite{II.SSD}]
\label{II.CFOt.2.1}Let $A$ be a selfadjoint operator with $Sp\left( A\right)
\subseteq \left[ m,M\right] $ for some real numbers $m<M.$ If $f,g:\left[ m,M%
\right] \longrightarrow \mathbb{R}$ are continuous and \textit{synchronous
(asynchronous) on }$\left[ m,M\right] ,$ \textit{then} 
\begin{equation}
\left\langle f\left( A\right) g\left( A\right) x,x\right\rangle \geq \left(
\leq \right) \left\langle f\left( A\right) x,x\right\rangle \cdot
\left\langle g\left( A\right) x,x\right\rangle  \label{II.CFOe.2.1}
\end{equation}%
for any $x\in H$ with $\left\Vert x\right\Vert =1.$
\end{theorem}

\begin{proof}
We consider only the case of synchronous functions. In this case we have
then 
\begin{equation}
f\left( t\right) g\left( t\right) +f\left( s\right) g\left( s\right) \geq
f\left( t\right) g\left( s\right) +f\left( s\right) g\left( t\right)
\label{II.CFOe.2.2}
\end{equation}%
for each $t,s\in \left[ a,b\right] .$

If we fix $s\in \left[ a,b\right] $ and apply the property (\ref{P}) for the
inequality (\ref{II.CFOe.2.2}) then we have for each $x\in H$ with $%
\left\Vert x\right\Vert =1$ that 
\begin{equation*}
\left\langle \left( f\left( A\right) g\left( A\right) +f\left( s\right)
g\left( s\right) 1_{H}\right) x,x\right\rangle \geq \left\langle \left(
g\left( s\right) f\left( A\right) +f\left( s\right) g\left( A\right) \right)
x,x\right\rangle ,
\end{equation*}%
which is clearly equivalent with 
\begin{equation}
\left\langle f\left( A\right) g\left( A\right) x,x\right\rangle +f\left(
s\right) g\left( s\right) \geq g\left( s\right) \left\langle f\left(
A\right) x,x\right\rangle +f\left( s\right) \left\langle g\left( A\right)
x,x\right\rangle  \label{II.CFOe.2.3}
\end{equation}%
for each $s\in \left[ a,b\right] .$

Now, if we apply again the property (\ref{P}) for the inequality (\ref%
{II.CFOe.2.3}), then \ we have for any $y\in H$ with $\left\Vert
y\right\Vert =1$ that 
\begin{align*}
& \left\langle \left( \left\langle f\left( A\right) g\left( A\right)
x,x\right\rangle 1_{H}+f\left( A\right) g\left( A\right) \right)
y,y\right\rangle \\
& \geq \left\langle \left( \left\langle f\left( A\right) x,x\right\rangle
g\left( A\right) +\left\langle g\left( A\right) x,x\right\rangle f\left(
A\right) \right) y,y\right\rangle ,
\end{align*}%
which is clearly equivalent with 
\begin{align}
& \left\langle f\left( A\right) g\left( A\right) x,x\right\rangle
+\left\langle f\left( A\right) g\left( A\right) y,y\right\rangle
\label{II.CFOe.2.4} \\
& \geq \left\langle f\left( A\right) x,x\right\rangle \left\langle g\left(
A\right) y,y\right\rangle +\left\langle f\left( A\right) y,y\right\rangle
\left\langle g\left( A\right) x,x\right\rangle  \notag
\end{align}%
for each $x,y\in H$ with $\left\Vert x\right\Vert =\left\Vert y\right\Vert
=1.$ This is an inequality of interest in itself.

Finally, on making $y=x$ in (\ref{II.CFOe.2.4}) we deduce the desired result
(\ref{II.CFOe.2.1}).
\end{proof}

Some particular cases are of interest for applications. In the first
instance we consider the case of power functions.

\begin{example}
\label{II.CFOc.2.3}Assume that $A$ is a positive operator on the Hilbert
space $H$ and $p,q>0.$ Then for each $x\in H$ with $\left\Vert x\right\Vert
=1$ we have the inequality 
\begin{equation}
\left\langle A^{p+q}x,x\right\rangle \geq \left\langle A^{p}x,x\right\rangle
\cdot \left\langle A^{q}x,x\right\rangle .  \label{II.CFOe.2.6}
\end{equation}%
If $A$ is positive definite then the inequality (\ref{II.CFOe.2.6}) also
holds for $p,q<0.$

If $A$ is positive definite and either $p>0,q<0$ or $p<0,q>0$, then the
reverse inequality holds in (\ref{II.CFOe.2.6}).
\end{example}

Another case of interest for applications is the exponential function.

\begin{example}
\label{II.CFOc.2.4}Assume that $A$ is a selfadjoint operator on $H.$ If $%
\alpha ,\beta >0$ or $\alpha ,\beta <0,$ then 
\begin{equation}
\left\langle \exp \left[ \left( \alpha +\beta \right) A\right]
x,x\right\rangle \geq \left\langle \exp \left( \alpha A\right)
x,x\right\rangle \cdot \left\langle \exp \left( \beta A\right)
x,x\right\rangle  \label{II.CFOe.2.7}
\end{equation}%
for each $x\in H$ with $\left\Vert x\right\Vert =1.$

If either $\alpha >0,\beta <0$ or $\alpha <0,\beta >0,$ then the reverse
inequality holds in (\ref{II.CFOe.2.7}).
\end{example}

The following particular cases may be of interest as well:

\begin{example}
\textbf{a.} Assume that $A$ is positive definite and $p>0.$ Then 
\begin{equation}
\left\langle A^{p}\log Ax,x\right\rangle \geq \left\langle
A^{p}x,x\right\rangle \cdot \left\langle \log Ax,x\right\rangle
\label{II.CFOe.2.8}
\end{equation}%
for each $x\in H$ with $\left\Vert x\right\Vert =1.$ If $p<0$ then the
reverse inequality holds in (\ref{II.CFOe.2.8}).

\textbf{b. }Assume that $A$ is positive definite and $Sp\left( A\right)
\subset \left( 0,1\right) .$ If $r,s>0$ or $r,s<0$ then 
\begin{align}
& \left\langle \left( 1_{H}-A^{r}\right) ^{-1}\left( 1_{H}-A^{s}\right)
^{-1}x,x\right\rangle  \label{II.CFOe.2.9} \\
& \geq \left\langle \left( 1_{H}-A^{r}\right) ^{-1}x,x\right\rangle \cdot
\left\langle \left( 1_{H}-A^{s}\right) ^{-1}x,x\right\rangle  \notag
\end{align}%
for each $x\in H$ with $\left\Vert x\right\Vert =1.$

If either $r>0,s<0$ or $r<0,s>0,$ then the reverse inequality holds in (\ref%
{II.CFOe.2.9}).
\end{example}

\begin{remark}
\label{II.CFOr.2.2}We observe, from the proof of the above theorem that, if $%
A$ and $B$ are selfadjoint operators and $Sp\left( A\right) ,Sp\left(
B\right) \subseteq \left[ m,M\right] ,$ then for any continuous synchronous
(asynchronous) functions $f,g:\left[ m,M\right] \longrightarrow \mathbb{R}$
we have the more general result 
\begin{align}
& \left\langle f\left( A\right) g\left( A\right) x,x\right\rangle
+\left\langle f\left( B\right) g\left( B\right) y,y\right\rangle
\label{II.CFOe.2.10} \\
& \geq \left( \leq \right) \left\langle f\left( A\right) x,x\right\rangle
\left\langle g\left( B\right) y,y\right\rangle +\left\langle f\left(
B\right) y,y\right\rangle \left\langle g\left( A\right) x,x\right\rangle 
\notag
\end{align}%
for each $x,y\in H$ with $\left\Vert x\right\Vert =\left\Vert y\right\Vert
=1.$

If $f:\left[ m,M\right] \longrightarrow \left( 0,\infty \right) $ is
continuous then the functions $f^{p},f^{q}$ are synchronous in the case when 
$p,q>0$ or $p,q<0$ and asynchronous when either $p>0,q<0$ or $p<0,q>0.$ In
this situation if $A$ and $B$ are positive definite operators then we have
the inequality 
\begin{align}
& \left\langle f^{p+q}\left( A\right) x,x\right\rangle +\left\langle
f^{p+q}\left( B\right) y,y\right\rangle  \label{II.CFOce.2.11} \\
& \geq \left\langle f^{p}\left( A\right) x,x\right\rangle \left\langle
f^{q}\left( B\right) y,y\right\rangle +\left\langle f^{p}\left( B\right)
y,y\right\rangle \left\langle f^{q}\left( A\right) x,x\right\rangle  \notag
\end{align}%
for each $x,y\in H$ with $\left\Vert x\right\Vert =\left\Vert y\right\Vert
=1 $ where either $p,q>0$ or $p,q<0.$ If $p>0,q<0$ or $p<0,q>0$ then the
reverse inequality also holds in (\ref{II.CFOce.2.11}).

As particular cases, we should observe that for $p=q=1$ and $f\left(
t\right) =t,$ we get from (\ref{II.CFOce.2.11}) the inequality 
\begin{equation}
\left\langle A^{2}x,x\right\rangle +\left\langle B^{2}y,y\right\rangle \geq
2\cdot \left\langle Ax,x\right\rangle \left\langle By,y\right\rangle
\label{II.CFOe.2.12}
\end{equation}%
for each $x,y\in H$ with $\left\Vert x\right\Vert =\left\Vert y\right\Vert
=1 $.

For $p=1$ and $q=-1$ we have from (\ref{II.CFOce.2.11}) 
\begin{equation}
\left\langle Ax,x\right\rangle \left\langle B^{-1}y,y\right\rangle
+\left\langle By,y\right\rangle \left\langle A^{-1}x,x\right\rangle \leq 2
\label{II.CFOe.2.13}
\end{equation}%
for each $x,y\in H$ with $\left\Vert x\right\Vert =\left\Vert y\right\Vert
=1.$
\end{remark}

\subsection{A Version of the \v{C}eby\v{s}ev Inequality for $n$ Operators}

The following multiple operator version of Theorem \ref{II.CFOt.2.1} holds:

\begin{theorem}[Dragomir, 2008, \protect\cite{II.SSD}]
\label{II.CFOt.3.1}Let $A_{j}$ be selfadjoint operators with $Sp\left(
A_{j}\right) \subseteq \left[ m,M\right] $ for $j\in \left\{ 1,\dots
,n\right\} $ and for some scalars $m<M.$ If $f,g:\left[ m,M\right]
\longrightarrow \mathbb{R}$ are continuous and \textit{synchronous
(asynchronous) on }$\left[ m,M\right] ,$ \textit{then} 
\begin{equation}
\sum_{j=1}^{n}\left\langle f\left( A_{j}\right) g\left( A_{j}\right)
x_{j},x_{j}\right\rangle \geq \left( \leq \right) \sum_{j=1}^{n}\left\langle
f\left( A_{j}\right) x_{j},x_{j}\right\rangle \cdot
\sum_{j=1}^{n}\left\langle g\left( A_{j}\right) x_{j},x_{j}\right\rangle ,
\label{II.CFOe.3.1}
\end{equation}%
for each $x_{j}\in H,j\in \left\{ 1,\dots ,n\right\} $ with $%
\sum_{j=1}^{n}\left\Vert x_{j}\right\Vert ^{2}=1.$
\end{theorem}

\begin{proof}
As in \cite[p. 6]{II.FMPS}, if we put 
\begin{equation*}
\widetilde{A}:=\left( 
\begin{array}{ccc}
A_{1} & \cdots & 0 \\ 
\vdots & \ddots & \vdots \\ 
0 & \cdots & A_{n}%
\end{array}%
\right) \text{ \quad and\quad\ }\widetilde{x}=\left( 
\begin{array}{c}
x_{1} \\ 
\vdots \\ 
x_{n}%
\end{array}%
\right)
\end{equation*}

then we have $Sp\left( \widetilde{A}\right) \subseteq \left[ m,M\right] ,$ $%
\left\Vert \widetilde{x}\right\Vert =1,$%
\begin{equation*}
\left\langle f\left( \widetilde{A}\right) g\left( \widetilde{A}\right) 
\widetilde{x},\widetilde{x}\right\rangle =\sum_{j=1}^{n}\left\langle f\left(
A_{j}\right) g\left( A_{j}\right) x_{j},x_{j}\right\rangle ,
\end{equation*}%
\begin{equation*}
\left\langle f\left( \widetilde{A}\right) \widetilde{x},\widetilde{x}%
\right\rangle =\sum_{j=1}^{n}\left\langle f\left( A_{j}\right)
x_{j},x_{j}\right\rangle \text{ \quad and\quad\ }\left\langle g\left( 
\widetilde{A}\right) \widetilde{x},\widetilde{x}\right\rangle
=\sum_{j=1}^{n}\left\langle g\left( A_{j}\right) x_{j},x_{j}\right\rangle .
\end{equation*}%
Applying Theorem \ref{II.CFOt.2.1} for $\widetilde{A}$ and $\widetilde{x}$
we deduce the desired result (\ref{II.CFOe.3.1}).
\end{proof}

The following particular cases may be of interest for applications.

\begin{example}
\label{II.CFOc.3.2}Assume that $A_{j},j\in \left\{ 1,\dots ,n\right\} $ are
positive operators on the Hilbert space $H$ and $p,q>0.$ Then for each $%
x_{j}\in H,j\in \left\{ 1,\dots ,n\right\} $ with $\sum_{j=1}^{n}\left\Vert
x_{j}\right\Vert ^{2}=1$ we have the inequality 
\begin{equation}
\left\langle \sum_{j=1}^{n}A_{j}^{p+q}x_{j},x_{j}\right\rangle \geq
\sum_{j=1}^{n}\left\langle A_{j}^{p}x_{j},x_{j}\right\rangle \cdot
\sum_{j=1}^{n}\left\langle A_{j}^{q}x_{j},x_{j}\right\rangle .
\label{II.CFOe.3.2}
\end{equation}%
If $A_{j}$ are positive definite then the inequality (\ref{II.CFOe.3.2})
also holds for $p,q<0.$

If $A_{j}$ are positive definite and either $p>0,q<0$ or $p<0,q>0$, then the
reverse inequality holds in (\ref{II.CFOe.3.2}).
\end{example}

Another case of interest for applications is the exponential function.

\begin{example}
\label{II.CFOc.3.3}Assume that $A_{j},j\in \left\{ 1,\dots ,n\right\} $ are
selfadjoint operators on $H.$ If $\alpha ,\beta >0$ or $\alpha ,\beta <0,$
then 
\begin{align}
& \left\langle \sum_{j=1}^{n}\exp \left[ \left( \alpha +\beta \right) A_{j}%
\right] x_{j},x_{j}\right\rangle  \label{II.CFOe.3.3} \\
& \geq \sum_{j=1}^{n}\left\langle \exp \left( \alpha A_{j}\right)
x_{j},x_{j}\right\rangle \cdot \sum_{j=1}^{n}\left\langle \exp \left( \beta
A_{j}\right) x_{j},x_{j}\right\rangle  \notag
\end{align}%
for each $x_{j}\in H,j\in \left\{ 1,\dots ,n\right\} $ with $%
\sum_{j=1}^{n}\left\Vert x_{j}\right\Vert ^{2}=1.$

If either $\alpha >0,\beta <0$ or $\alpha <0,\beta >0,$ then the reverse
inequality holds in (\ref{II.CFOe.3.3}).
\end{example}

The following particular cases may be of interest as well:

\begin{example}
\textbf{a.} Assume that $A_{j},j\in \left\{ 1,\dots ,n\right\} $ are
positive definite operators and $p>0.$ Then 
\begin{equation}
\left\langle \sum_{j=1}^{n}A_{j}^{p}\log A_{j}x_{j},x_{j}\right\rangle \geq
\sum_{j=1}^{n}\left\langle A_{j}^{p}x_{j},x_{j}\right\rangle \cdot
\sum_{j=1}^{n}\left\langle \log A_{j}x_{j},x_{j}\right\rangle
\label{II.CFOe.3.4}
\end{equation}%
for each $x_{j}\in H,j\in \left\{ 1,\dots ,n\right\} $ with $%
\sum_{j=1}^{n}\left\Vert x_{j}\right\Vert ^{2}=1.$ If $p<0$ then the reverse
inequality holds in (\ref{II.CFOe.3.4}).

\textbf{b. }If $A_{j}$ are positive definite and $Sp\left( A_{j}\right)
\subset \left( 0,1\right) $ for $j\in \left\{ 1,\dots ,n\right\} $ then for $%
r,s>0$ or $r,s<0$ we have the inequality 
\begin{align}
& \left\langle \sum_{j=1}^{n}\left( 1_{H}-A_{j}^{r}\right) ^{-1}\left(
1_{H}-A_{j}^{s}\right) ^{-1}x_{j},x_{j}\right\rangle  \label{II.CFOe.3.5} \\
& \geq \sum_{j=1}^{n}\left\langle \left( 1_{H}-A_{j}^{r}\right)
^{-1}x_{j},x_{j}\right\rangle \cdot \sum_{j=1}^{n}\left\langle \left(
1_{H}-A_{j}^{s}\right) ^{-1}x_{j},x_{j}\right\rangle  \notag
\end{align}%
for each $x_{j}\in H,j\in \left\{ 1,\dots ,n\right\} $ with $%
\sum_{j=1}^{n}\left\Vert x_{j}\right\Vert ^{2}=1.$

If either $r>0,s<0$ or $r<0,s>0,$ then the reverse inequality holds in (\ref%
{II.CFOe.3.5}).
\end{example}

\subsection{Another Version of the \v{C}eby\v{s}ev Inequality for $n$
Operators}

The following different version of the \v{C}eby\v{s}ev inequality for a
sequence of operators also holds:

\begin{theorem}[Dragomir, 2008, \protect\cite{II.SSD}]
\label{II.CFOt.4.1}Let $A_{j}$ be selfadjoint operators with $Sp\left(
A_{j}\right) \subseteq \left[ m,M\right] $ for $j\in \left\{ 1,\dots
,n\right\} $ and for some scalars $m<M.$ If $f,g:\left[ m,M\right]
\longrightarrow \mathbb{R}$ are continuous and \textit{synchronous
(asynchronous) on }$\left[ m,M\right] ,$ \textit{then} 
\begin{align}
& \left\langle \sum_{j=1}^{n}p_{j}f\left( A_{j}\right) g\left( A_{j}\right)
x,x\right\rangle  \label{II.CFOe.4.1} \\
& \geq \left( \leq \right) \left\langle \sum_{j=1}^{n}p_{j}f\left(
A_{j}\right) x,x\right\rangle \cdot \left\langle \sum_{j=1}^{n}p_{j}g\left(
A_{j}\right) x,x\right\rangle ,  \notag
\end{align}%
for any $p_{j}\geq 0,j\in \left\{ 1,\dots ,n\right\} $ with $%
\sum_{j=1}^{n}p_{j}=1$ and $x\in H$ with $\left\Vert x\right\Vert =1.$

In particular 
\begin{align}
& \left\langle \frac{1}{n}\sum_{j=1}^{n}f\left( A_{j}\right) g\left(
A_{j}\right) x,x\right\rangle  \label{II.CFOe.4.2} \\
& \geq \left( \leq \right) \left\langle \frac{1}{n}\sum_{j=1}^{n}f\left(
A_{j}\right) x,x\right\rangle \cdot \left\langle \frac{1}{n}%
\sum_{j=1}^{n}g\left( A_{j}\right) x,x\right\rangle ,  \notag
\end{align}%
for each $x\in H$ with $\left\Vert x\right\Vert =1.$
\end{theorem}

\begin{proof}
We provide here two proofs. The first is based on the inequality (\ref%
{II.CFOe.2.10}) and generates as a by-product a more general result. The
second is derived from Theorem \ref{II.CFOt.3.1}.

\textbf{1.} If we make use of the inequality (\ref{II.CFOe.2.10}), then we
can write 
\begin{align}
& \left\langle f\left( A_{j}\right) g\left( A_{j}\right) x,x\right\rangle
+\left\langle f\left( B_{k}\right) g\left( B_{k}\right) y,y\right\rangle
\label{II.CFOe.4.3} \\
& \geq \left( \leq \right) \left\langle f\left( A_{j}\right)
x,x\right\rangle \left\langle g\left( B_{k}\right) y,y\right\rangle
+\left\langle f\left( B_{k}\right) y,y\right\rangle \left\langle g\left(
A_{j}\right) x,x\right\rangle ,  \notag
\end{align}%
which holds for any $A_{j}$ and $B_{k}$ selfadjoint operators with $Sp\left(
A_{j}\right) ,Sp\left( B_{k}\right) \subseteq \left[ m,M\right] ,$ $j,k\in
\left\{ 1,\dots ,n\right\} $ and for each $x,y\in H$ with $\left\Vert
x\right\Vert =\left\Vert y\right\Vert =1.$

Now, if $p_{j}\geq 0,q_{k}\geq 0,j,k\in \left\{ 1,\dots ,n\right\} $ and $%
\sum_{j=1}^{n}p_{j}=\sum_{k=1}^{n}q_{k}=1$ then, by multiplying (\ref%
{II.CFOe.4.3}) with $p_{j}\geq 0,q_{k}\geq 0$ and summing over $j$ and $k$
from $1$ to $n$ we deduce the following inequality that is of interest in
its own right: 
\begin{align}
& \left\langle \sum_{j=1}^{n}p_{j}f\left( A_{j}\right) g\left( A_{j}\right)
x,x\right\rangle +\left\langle \sum_{k=1}^{n}q_{k}f\left( B_{k}\right)
g\left( B_{k}\right) y,y\right\rangle  \label{II.CFOe.4.4} \\
& \geq \left( \leq \right) \left\langle \sum_{j=1}^{n}p_{j}f\left(
A_{j}\right) x,x\right\rangle \left\langle \sum_{k=1}^{n}q_{k}g\left(
B_{k}\right) y,y\right\rangle  \notag \\
& +\left\langle \sum_{k=1}^{n}q_{k}f\left( B_{k}\right) y,y\right\rangle
\left\langle \sum_{j=1}^{n}p_{j}g\left( A_{j}\right) x,x\right\rangle  \notag
\end{align}%
for each $x,y\in H$ with $\left\Vert x\right\Vert =\left\Vert y\right\Vert
=1.$

Finally, the choice $B_{k}=A_{k},q_{k}=p_{k}$ and $y=x$ in (\ref{II.CFOe.4.4}%
) produces the desired result (\ref{II.CFOe.4.1}).

\textbf{2.} In we choose in Theorem \ref{II.CFOt.3.1} $x_{j}=\sqrt{p_{j}}%
\cdot x,$ $j\in \left\{ 1,\dots ,n\right\} ,$ where $p_{j}\geq 0,j\in
\left\{ 1,\dots ,n\right\} ,$ $\sum_{j=1}^{n}p_{j}=1$ and $x\in H,$ with $%
\left\Vert x\right\Vert =1$ then a simple calculation shows that the
inequality (\ref{II.CFOe.3.1}) becomes (\ref{II.CFOe.4.1}). The details are
omitted.
\end{proof}

\begin{remark}
\label{II.CFOr.4.2}We remark that the case $n=1$ in (\ref{II.CFOe.4.1})
produces the inequality (\ref{II.CFOe.2.1}).
\end{remark}

The following particular cases are of interest:

\begin{example}
\label{II.CFOc.4.3}Assume that $A_{j},j\in \left\{ 1,\dots ,n\right\} $ are
positive operators on the Hilbert space $H,$ $p_{j}\geq 0,j\in \left\{
1,\dots ,n\right\} $ with $\sum_{j=1}^{n}p_{j}=1$ and $p,q>0.$ Then for each 
$x\in H$ with $\left\Vert x\right\Vert =1$ we have the inequality 
\begin{equation}
\left\langle \sum_{j=1}^{n}p_{j}A_{j}^{p+q}x,x\right\rangle \geq
\left\langle \sum_{j=1}^{n}p_{j}A_{j}^{p}x,x\right\rangle \cdot \left\langle
\sum_{j=1}^{n}p_{j}A_{j}^{q}x,x\right\rangle .  \label{II.CFOe.4.5}
\end{equation}%
If $A_{j},j\in \left\{ 1,\dots ,n\right\} $ are positive definite then the
inequality (\ref{II.CFOe.4.5}) also holds for $p,q<0.$

If $A_{j},j\in \left\{ 1,\dots ,n\right\} $ are positive definite and either 
$p>0,q<0$ or $p<0,q>0$, then the reverse inequality holds in (\ref%
{II.CFOe.4.5}).
\end{example}

Another case of interest for applications is the exponential function.

\begin{example}
\label{II.CFOc.4.4}Assume that $A_{j},j\in \left\{ 1,\dots ,n\right\} $ are
selfadjoint operators on $H$ and $p_{j}\geq 0,j\in \left\{ 1,\dots
,n\right\} $ with $\sum_{j=1}^{n}p_{j}=1.$ If $\alpha ,\beta >0$ or $\alpha
,\beta <0,$ then 
\begin{align}
& \left\langle \sum_{j=1}^{n}p_{j}\exp \left[ \left( \alpha +\beta \right)
A_{j}\right] x,x\right\rangle  \label{II.CFOe.4.6} \\
& \geq \left\langle \sum_{j=1}^{n}p_{j}\exp \left( \alpha A_{j}\right)
x,x\right\rangle \cdot \left\langle \sum_{j=1}^{n}p_{j}\exp \left( \beta
A_{j}\right) x,x\right\rangle  \notag
\end{align}%
for each $x\in H$ with $\left\Vert x\right\Vert =1.$

If either $\alpha >0,\beta <0$ or $\alpha <0,\beta >0,$ then the reverse
inequality holds in (\ref{II.CFOe.4.6}).
\end{example}

The following particular cases may be of interest as well:

\begin{example}
\textbf{a.} Assume that $A_{j},j\in \left\{ 1,\dots ,n\right\} $ are
positive definite operators on the Hilbert space $H,$ $p_{j}\geq 0,j\in
\left\{ 1,\dots ,n\right\} $ with $\sum_{j=1}^{n}p_{j}=1$ and $p>0.$ Then 
\begin{equation}
\left\langle \sum_{j=1}^{n}p_{j}A_{j}^{p}\log A_{j}x,x\right\rangle \geq
\left\langle \sum_{j=1}^{n}p_{j}A_{j}^{p}x,x\right\rangle \cdot \left\langle
\sum_{j=1}^{n}p_{j}\log A_{j}x,x\right\rangle .  \label{II.CFOe.4.7}
\end{equation}%
If $p<0$ then the reverse inequality holds in (\ref{II.CFOe.4.7}).

\textbf{b. }Assume that $A_{j},j\in \left\{ 1,\dots ,n\right\} $ are
positive definite operators on the Hilbert space $H,Sp\left( A_{j}\right)
\subset \left( 0,1\right) $ and $p_{j}\geq 0,j\in \left\{ 1,\dots ,n\right\} 
$ with $\sum_{j=1}^{n}p_{j}=1$. If $r,s>0$ or $r,s<0$ then 
\begin{align}
& \left\langle \sum_{j=1}^{n}p_{j}\left( 1_{H}-A_{j}^{r}\right) ^{-1}\left(
1_{H}-A_{j}^{s}\right) ^{-1}x,x\right\rangle  \label{II.CFOe.4.8} \\
& \geq \left\langle \sum_{j=1}^{n}p_{j}\left( 1_{H}-A_{j}^{r}\right)
^{-1}x,x\right\rangle \cdot \left\langle \sum_{j=1}^{n}p_{j}\left(
1_{H}-A_{j}^{s}\right) ^{-1}x,x\right\rangle  \notag
\end{align}%
for each $x\in H$ with $\left\Vert x\right\Vert =1.$

If either $r>0,s<0$ or $r<0,s>0,$ then the reverse inequality holds in (\ref%
{II.CFOe.4.8}).
\end{example}

We remark that the following operator norm inequality can be stated as well:

\begin{corollary}
\label{II.CFOc.4.5}Let $A_{j}$ be selfadjoint operators with $Sp\left(
A_{j}\right) \subseteq \left[ m,M\right] $ for $j\in \left\{ 1,\dots
,n\right\} $ and for some scalars $m<M.$ If $f,g:\left[ m,M\right]
\longrightarrow \mathbb{R}$ are continuous, \textit{asynchronous on }$\left[
m,M\right] $ and for $p_{j}\geq 0,j\in \left\{ 1,\dots ,n\right\} $ with $%
\sum_{j=1}^{n}p_{j}=1$ the operator $\sum_{j=1}^{n}p_{j}f\left( A_{j}\right)
g\left( A_{j}\right) $ is positive, \textit{then} 
\begin{equation}
\left\Vert \sum_{j=1}^{n}p_{j}f\left( A_{j}\right) g\left( A_{j}\right)
\right\Vert \leq \left\Vert \sum_{j=1}^{n}p_{j}f\left( A_{j}\right)
\right\Vert \cdot \left\Vert \sum_{j=1}^{n}p_{j}g\left( A_{j}\right)
\right\Vert .  \label{II.CFOe.4.9}
\end{equation}
\end{corollary}

\begin{proof}
We have from (\ref{II.CFOe.4.1}) that%
\begin{equation*}
0\leq \left\langle \sum_{j=1}^{n}p_{j}f\left( A_{j}\right) g\left(
A_{j}\right) x,x\right\rangle \leq \left\langle \sum_{j=1}^{n}p_{j}f\left(
A_{j}\right) x,x\right\rangle \cdot \left\langle \sum_{j=1}^{n}p_{j}g\left(
A_{j}\right) x,x\right\rangle
\end{equation*}%
for each $x\in H$ with $\left\Vert x\right\Vert =1.$ Taking the supremum in
this inequality over $x\in H$ with $\left\Vert x\right\Vert =1$ we deduce
the desired result (\ref{II.CFOe.4.9}).
\end{proof}

The above Corollary \ref{II.CFOc.4.5} provides some interesting norm
inequalities for sums of positive operators as follows:

\begin{example}
\textbf{a.} If $A_{j},j\in \left\{ 1,\dots ,n\right\} $ are positive
definite and either $p>0,q<0$ or $p<0,q>0$, then for $p_{j}\geq 0,j\in
\left\{ 1,\dots ,n\right\} $ with $\sum_{j=1}^{n}p_{j}=1$ we have the norm
inequality: 
\begin{equation}
\left\Vert \sum_{j=1}^{n}p_{j}A_{j}^{p+q}\right\Vert \leq \left\Vert
\sum_{j=1}^{n}p_{j}A_{j}^{p}\right\Vert \cdot \left\Vert
\sum_{j=1}^{n}p_{j}A_{j}^{q}\right\Vert .  \label{II.CFOe.4.10}
\end{equation}%
In particular 
\begin{equation}
1\leq \left\Vert \sum_{j=1}^{n}p_{j}A_{j}^{r}\right\Vert \cdot \left\Vert
\sum_{j=1}^{n}p_{j}A_{j}^{-r}\right\Vert  \label{II.CFOe.4.11}
\end{equation}%
for any $r>0.$

\textbf{b. }Assume that $A_{j},j\in \left\{ 1,\dots ,n\right\} $ are
selfadjoint operators on $H$ and $p_{j}\geq 0,j\in \left\{ 1,\dots
,n\right\} $ with $\sum_{j=1}^{n}p_{j}=1.$ If either $\alpha >0,\beta <0$ or 
$\alpha <0,\beta >0,$ then%
\begin{equation}
\left\Vert \sum_{j=1}^{n}p_{j}\exp \left[ \left( \alpha +\beta \right) A_{j}%
\right] \right\Vert \leq \left\Vert \sum_{j=1}^{n}p_{j}\exp \left( \alpha
A_{j}\right) \right\Vert \cdot \left\Vert \sum_{j=1}^{n}p_{j}\exp \left(
\beta A_{j}\right) \right\Vert .  \label{II.CFOe.4.12}
\end{equation}

In particular%
\begin{equation*}
1\leq \left\Vert \sum_{j=1}^{n}p_{j}\exp \left( \gamma A_{j}\right)
\right\Vert \cdot \left\Vert \sum_{j=1}^{n}p_{j}\exp \left( -\gamma
A_{j}\right) \right\Vert .
\end{equation*}%
for any $\gamma >0.$
\end{example}

\subsection{Related Results for One Operator}

The following result that is related to the \v{C}eby\v{s}ev inequality may
be stated:

\begin{theorem}[Dragomir, 2008, \protect\cite{II.SSD}]
\label{II.CFOt.5.1}Let $A$ be a selfadjoint operator with $Sp\left( A\right)
\subseteq \left[ m,M\right] $ for some real numbers $m<M.$ If $f,g:\left[ m,M%
\right] \longrightarrow \mathbb{R}$ are continuous and \textit{synchronous
on }$\left[ m,M\right] ,$ \textit{then} 
\begin{align}
& \left\langle f\left( A\right) g\left( A\right) x,x\right\rangle
-\left\langle f\left( A\right) x,x\right\rangle \cdot \left\langle g\left(
A\right) x,x\right\rangle  \label{II.CFOe.5.1} \\
& \geq \left[ \left\langle f\left( A\right) x,x\right\rangle -f\left(
\left\langle Ax,x\right\rangle \right) \right] \cdot \left[ g\left(
\left\langle Ax,x\right\rangle \right) -\left\langle g\left( A\right)
x,x\right\rangle \right]  \notag
\end{align}%
for any $x\in H$ with $\left\Vert x\right\Vert =1.$

If $f,g$ are asynchronous, then 
\begin{align}
& \left\langle f\left( A\right) x,x\right\rangle \cdot \left\langle g\left(
A\right) x,x\right\rangle -\left\langle f\left( A\right) g\left( A\right)
x,x\right\rangle  \label{II.CFOe.5.1.a} \\
& \geq \left[ \left\langle f\left( A\right) x,x\right\rangle -f\left(
\left\langle Ax,x\right\rangle \right) \right] \cdot \left[ \left\langle
g\left( A\right) x,x\right\rangle -g\left( \left\langle Ax,x\right\rangle
\right) \right]  \notag
\end{align}%
for any $x\in H$ with $\left\Vert x\right\Vert =1.$
\end{theorem}

\begin{proof}
Since $f,g$ are synchronous and $m\leq \left\langle Ax,x\right\rangle \leq M$
for any $x\in H$ with $\left\Vert x\right\Vert =1,$ then we have 
\begin{equation}
\left[ f\left( t\right) -f\left( \left\langle Ax,x\right\rangle \right) %
\right] \left[ g\left( t\right) -g\left( \left\langle Ax,x\right\rangle
\right) \right] \geq 0  \label{II.CFOe.5.2}
\end{equation}%
for any $t\in \left[ a,b\right] $ and $x\in H$ with $\left\Vert x\right\Vert
=1.$

On utilising the property (\ref{P}) for the inequality (\ref{II.CFOe.5.2})
we have that 
\begin{equation}
\left\langle \left[ f\left( B\right) -f\left( \left\langle Ax,x\right\rangle
\right) \right] \left[ g\left( B\right) -g\left( \left\langle
Ax,x\right\rangle \right) \right] y,y\right\rangle \geq 0
\label{II.CFOe.5.3}
\end{equation}%
for any $B$ a bounded linear operator with $Sp\left( B\right) \subseteq %
\left[ m,M\right] $ and $y\in H$ with $\left\Vert y\right\Vert =1.$

Since 
\begin{align}
& \left\langle \left[ f\left( B\right) -f\left( \left\langle
Ax,x\right\rangle \right) \right] \left[ g\left( B\right) -g\left(
\left\langle Ax,x\right\rangle \right) \right] y,y\right\rangle
\label{II.CFOe.5.4} \\
& =\left\langle f\left( B\right) g\left( B\right) y,y\right\rangle +f\left(
\left\langle Ax,x\right\rangle \right) g\left( \left\langle
Ax,x\right\rangle \right)  \notag \\
& -\left\langle f\left( B\right) y,y\right\rangle g\left( \left\langle
Ax,x\right\rangle \right) -f\left( \left\langle Ax,x\right\rangle \right)
\left\langle g\left( B\right) y,y\right\rangle ,  \notag
\end{align}%
then from (\ref{II.CFOe.5.3}) we get 
\begin{align*}
& \left\langle f\left( B\right) g\left( B\right) y,y\right\rangle +f\left(
\left\langle Ax,x\right\rangle \right) g\left( \left\langle
Ax,x\right\rangle \right) \\
& \geq \left\langle f\left( B\right) y,y\right\rangle g\left( \left\langle
Ax,x\right\rangle \right) +f\left( \left\langle Ax,x\right\rangle \right)
\left\langle g\left( B\right) y,y\right\rangle
\end{align*}%
which is clearly equivalent with 
\begin{align}
& \left\langle f\left( B\right) g\left( B\right) y,y\right\rangle
-\left\langle f\left( A\right) y,y\right\rangle \cdot \left\langle g\left(
A\right) y,y\right\rangle  \label{II.CFOe.5.5} \\
& \geq \left[ \left\langle f\left( B\right) y,y\right\rangle -f\left(
\left\langle Ax,x\right\rangle \right) \right] \cdot \left[ g\left(
\left\langle Ax,x\right\rangle \right) -\left\langle g\left( B\right)
y,y\right\rangle \right]  \notag
\end{align}%
for each $x,y\in H$ with $\left\Vert x\right\Vert =\left\Vert y\right\Vert
=1.$ This inequality is of interest in its own right.

Now, if we choose $B=A$ and $y=x$ in (\ref{II.CFOe.5.5}), then we deduce the
desired result (\ref{II.CFOe.5.1}).
\end{proof}

The following result which improves the \v{C}eby\v{s}ev inequality may be
stated:

\begin{corollary}[Dragomir, 2008, \protect\cite{II.SSD}]
\label{II.CFOc.5.1}Let $A$ be a selfadjoint operator with $Sp\left( A\right)
\subseteq \left[ m,M\right] $ for some real numbers $m<M.$ If $f,g:\left[ m,M%
\right] \longrightarrow \mathbb{R}$ are continuous, \textit{synchronous and
one is convex while the other is concave on }$\left[ m,M\right] ,$ \textit{%
then} 
\begin{align}
& \left\langle f\left( A\right) g\left( A\right) x,x\right\rangle
-\left\langle f\left( A\right) x,x\right\rangle \cdot \left\langle g\left(
A\right) x,x\right\rangle  \label{II.CFOe.5.6} \\
& \geq \left[ \left\langle f\left( A\right) x,x\right\rangle -f\left(
\left\langle Ax,x\right\rangle \right) \right] \cdot \left[ g\left(
\left\langle Ax,x\right\rangle \right) -\left\langle g\left( A\right)
x,x\right\rangle \right] \geq 0  \notag
\end{align}%
for any $x\in H$ with $\left\Vert x\right\Vert =1.$

If $f,g$ are asynchronous and either both of them are convex or both of them
concave on $\left[ m,M\right] $, then 
\begin{align}
& \left\langle f\left( A\right) x,x\right\rangle \cdot \left\langle g\left(
A\right) x,x\right\rangle -\left\langle f\left( A\right) g\left( A\right)
x,x\right\rangle  \label{II.CFOe.5.7} \\
& \geq \left[ \left\langle f\left( A\right) x,x\right\rangle -f\left(
\left\langle Ax,x\right\rangle \right) \right] \cdot \left[ \left\langle
g\left( A\right) x,x\right\rangle -g\left( \left\langle Ax,x\right\rangle
\right) \right] \geq 0  \notag
\end{align}%
for any $x\in H$ with $\left\Vert x\right\Vert =1.$
\end{corollary}

\begin{proof}
The second inequality follows by making use of the result due to Mond \& Pe%
\v{c}ari\'{c}, see \cite{II.CFOMP1}, \cite{II.MP2} or \cite[p. 5]{II.FMPS}: 
\begin{equation}
\left\langle h\left( A\right) x,x\right\rangle \geq \left( \leq \right)
h\left( \left\langle Ax,x\right\rangle \right)  \tag{MP}  \label{II.CFOMP}
\end{equation}%
for any $x\in H$ with $\left\Vert x\right\Vert =1$ provided that $A$ is a
selfadjoint operator with $Sp\left( A\right) \subseteq \left[ m,M\right] $
for some real numbers $m<M$ and $h$ is convex (concave) on the given
interval $\left[ m,M\right] .$
\end{proof}

The above Corollary \ref{II.CFOc.5.1} offers the possibility to improve some
of the results established before for power function as follows:

\begin{example}
\textbf{a.} Assume that $A$ is a positive operator on the Hilbert space $H.$
If $p\in \left( 0,1\right) $ and $q\in \left( 1,\infty \right) ,$ then for
each $x\in H$ with $\left\Vert x\right\Vert =1$ we have the inequality 
\begin{align}
& \left\langle A^{p+q}x,x\right\rangle -\left\langle A^{p}x,x\right\rangle
\cdot \left\langle A^{q}x,x\right\rangle  \label{II.CFOe.5.8} \\
& \geq \left[ \left\langle A^{q}x,x\right\rangle -\left\langle
Ax,x\right\rangle ^{q}\right] \left[ \left\langle Ax,x\right\rangle
^{p}-\left\langle A^{p}x,x\right\rangle \right] \geq 0.  \notag
\end{align}%
If $A$ is positive definite and $p>1,q<0,$ then 
\begin{align}
& \left\langle A^{p}x,x\right\rangle \cdot \left\langle
A^{q}x,x\right\rangle -\left\langle A^{p+q}x,x\right\rangle
\label{II.CFOe.5.9} \\
& \geq \left[ \left\langle A^{q}x,x\right\rangle -\left\langle
Ax,x\right\rangle ^{q}\right] \left[ \left\langle A^{p}x,x\right\rangle
-\left\langle Ax,x\right\rangle ^{p}\right] \geq 0  \notag
\end{align}%
for each $x\in H$ with $\left\Vert x\right\Vert =1.$

\textbf{b. }Assume that $A$ is positive definite and $p>1.$ Then 
\begin{align}
& \left\langle A^{p}\log Ax,x\right\rangle -\left\langle
A^{p}x,x\right\rangle \cdot \left\langle \log Ax,x\right\rangle
\label{II.CFOe.5.10} \\
& \geq \left[ \left\langle A^{p}x,x\right\rangle -\left\langle
Ax,x\right\rangle ^{p}\right] \left[ \log \left\langle Ax,x\right\rangle
-\left\langle \log Ax,x\right\rangle \right] \geq 0  \notag
\end{align}%
for each $x\in H$ with $\left\Vert x\right\Vert =1.$
\end{example}

\subsection{Related Results for $n$ Operators}

We can state now the following generalisation of Theorem \ref{II.CFOt.5.1}
for $n$ operators:

\begin{theorem}[Dragomir, 2008, \protect\cite{II.SSD}]
\label{II.CFOt.6.1}Let $A_{j}$ be selfadjoint operators with $Sp\left(
A_{j}\right) \subseteq \left[ m,M\right] $ for $j\in \left\{ 1,\dots
,n\right\} $ and for some scalars $m<M.$

(i) \ If $f,g:\left[ m,M\right] \longrightarrow \mathbb{R}$ are continuous
and \textit{synchronous on }$\left[ m,M\right] ,$ \textit{then} 
\begin{align}
& \sum_{j=1}^{n}\left\langle f\left( A_{j}\right) g\left( A_{j}\right)
x_{j},x_{j}\right\rangle -\sum_{j=1}^{n}\left\langle f\left( A_{j}\right)
x_{j},x_{j}\right\rangle \cdot \sum_{j=1}^{n}\left\langle g\left(
A_{j}\right) x_{j},x_{j}\right\rangle  \label{II.CFOe.6.1} \\
& \geq \left[ \sum_{j=1}^{n}\left\langle f\left( A_{j}\right)
x_{j},x_{j}\right\rangle -f\left( \sum_{j=1}^{n}\left\langle
A_{j}x_{j},x_{j}\right\rangle \right) \right]  \notag \\
& \times \left[ g\left( \sum_{j=1}^{n}\left\langle
A_{j}x_{j},x_{j}\right\rangle \right) -\sum_{j=1}^{n}\left\langle g\left(
A_{j}\right) x_{j},x_{j}\right\rangle \right]  \notag
\end{align}%
for each $x_{j}\in H,j\in \left\{ 1,\dots ,n\right\} $ with $%
\sum_{j=1}^{n}\left\Vert x_{j}\right\Vert ^{2}=1.$ Moreover, if \textit{one
function is convex while the other is concave on }$\left[ m,M\right] ,$ 
\textit{then the right hand side of (\ref{II.CFOe.6.1}) is nonnegative. }

(ii) If $f,g$ are a\textit{synchronous on }$\left[ m,M\right] ,$ then 
\begin{align}
& \sum_{j=1}^{n}\left\langle f\left( A_{j}\right) x_{j},x_{j}\right\rangle
\cdot \sum_{j=1}^{n}\left\langle g\left( A_{j}\right)
x_{j},x_{j}\right\rangle -\sum_{j=1}^{n}\left\langle f\left( A_{j}\right)
g\left( A_{j}\right) x_{j},x_{j}\right\rangle  \label{II.CFOe.6.2} \\
& \geq \left[ \sum_{j=1}^{n}\left\langle f\left( A_{j}\right)
x_{j},x_{j}\right\rangle -f\left( \sum_{j=1}^{n}\left\langle
A_{j}x_{j},x_{j}\right\rangle \right) \right]  \notag \\
& \times \left[ \sum_{j=1}^{n}\left\langle g\left( A_{j}\right)
x_{j},x_{j}\right\rangle -g\left( \sum_{j=1}^{n}\left\langle
A_{j}x_{j},x_{j}\right\rangle \right) \right]  \notag
\end{align}%
for each $x_{j}\in H,j\in \left\{ 1,\dots ,n\right\} $ with $%
\sum_{j=1}^{n}\left\Vert x_{j}\right\Vert ^{2}=1.$ Moreover, if either both
of them are convex or both of them are concave on $\left[ m,M\right] $, then
the right hand side of (\ref{II.CFOe.6.2}) is nonnegative as well.
\end{theorem}

\begin{proof}
The argument is similar to the one from the proof of Theorem \ref%
{II.CFOt.3.1} on utilising the results from one operator obtained in Theorem %
\ref{II.CFOt.5.1}.

The nonnegativity of the right hand sides of the inequalities (\ref%
{II.CFOe.6.1}) and (\ref{II.CFOe.6.2}) follows by the use of the Jensen's
type result from \cite[p. 5]{II.FMPS} 
\begin{equation}
\sum_{j=1}^{n}\left\langle h\left( A_{j}\right) x_{j},x_{j}\right\rangle
\geq \left( \leq \right) h\left( \sum_{j=1}^{n}\left\langle
A_{j}x_{j},x_{j}\right\rangle \right)  \label{II.CFOe.6.2.a}
\end{equation}%
for each $x_{j}\in H,j\in \left\{ 1,\dots ,n\right\} $ with $%
\sum_{j=1}^{n}\left\Vert x_{j}\right\Vert ^{2}=1,$ which holds provided that 
$A_{j}$ are selfadjoint operators with $Sp\left( A_{j}\right) \subseteq %
\left[ m,M\right] $ for $j\in \left\{ 1,\dots ,n\right\} $ and for some
scalars $m<M$ and $h$ is convex (concave) on $\left[ m,M\right] .$

The details are omitted.
\end{proof}

\begin{example}
\textbf{a.} Assume that $A_{j},j\in \left\{ 1,\dots ,n\right\} $ are
positive operators on the Hilbert space $H.$ If $p\in \left( 0,1\right) $
and $q\in \left( 1,\infty \right) ,$ then for each $x_{j}\in H,j\in \left\{
1,\dots ,n\right\} $ with $\sum_{j=1}^{n}\left\Vert x_{j}\right\Vert ^{2}=1$
we have the inequality 
\begin{align}
& \sum_{j=1}^{n}\left\langle A_{j}^{p+q}x_{j},x_{j}\right\rangle
-\sum_{j=1}^{n}\left\langle A_{j}^{p}x_{j},x_{j}\right\rangle \cdot
\sum_{j=1}^{n}\left\langle A_{j}^{q}x_{j},x_{j}\right\rangle
\label{II.CFOe.6.3} \\
& \geq \left[ \sum_{j=1}^{n}\left\langle A_{j}^{q}x_{j},x_{j}\right\rangle
-\left( \sum_{j=1}^{n}\left\langle A_{j}x_{j},x_{j}\right\rangle \right) ^{q}%
\right]  \notag \\
& \times \left[ \left( \sum_{j=1}^{n}\left\langle
A_{j}x_{j},x_{j}\right\rangle \right) ^{p}-\sum_{j=1}^{n}\left\langle
A_{j}^{p}x_{j},x_{j}\right\rangle \right]  \notag \\
& \geq 0.  \notag
\end{align}%
If $A_{j}$ are positive definite and $p>1,q<0,$ then 
\begin{align}
& \sum_{j=1}^{n}\left\langle A_{j}^{p}x_{j},x_{j}\right\rangle \cdot
\sum_{j=1}^{n}\left\langle A_{j}^{q}x_{j},x_{j}\right\rangle
-\sum_{j=1}^{n}\left\langle A_{j}^{p+q}x_{j},x_{j}\right\rangle
\label{CFOe.6.4} \\
& \geq \left[ \sum_{j=1}^{n}\left\langle A_{j}^{q}x_{j},x_{j}\right\rangle
-\left( \sum_{j=1}^{n}\left\langle A_{j}x_{j},x_{j}\right\rangle \right) ^{q}%
\right]  \notag \\
& \times \left[ \sum_{j=1}^{n}\left\langle A_{j}^{p}x_{j},x_{j}\right\rangle
-\left( \sum_{j=1}^{n}\left\langle A_{j}x_{j},x_{j}\right\rangle \right) ^{p}%
\right]  \notag \\
& \geq 0  \notag
\end{align}%
for each $x_{j}\in H,j\in \left\{ 1,\dots ,n\right\} $ with $%
\sum_{j=1}^{n}\left\Vert x_{j}\right\Vert ^{2}=1.$

\textbf{b. }Assume that $A_{j}$ are positive definite and $p>1.$ Then 
\begin{align}
& \sum_{j=1}^{n}\left\langle A_{j}^{p}\log Ax_{j},x_{j}\right\rangle
-\sum_{j=1}^{n}\left\langle A_{j}^{p}x_{j},x_{j}\right\rangle \cdot
\sum_{j=1}^{n}\left\langle \log A_{j}x_{j},x_{j}\right\rangle
\label{II.CFOe.6.5} \\
& \geq \left[ \sum_{j=1}^{n}\left\langle A_{j}^{p}x_{j},x_{j}\right\rangle
-\left( \sum_{j=1}^{n}\left\langle A_{j}x_{j},x_{j}\right\rangle \right) ^{p}%
\right]  \notag \\
& \times \left[ \sum_{j=1}^{n}\log \left\langle
A_{j}x_{j},x_{j}\right\rangle -\log \left( \sum_{j=1}^{n}\left\langle
A_{j}x_{j},x_{j}\right\rangle \right) \right]  \notag \\
& \geq 0  \notag
\end{align}%
for each $x_{j}\in H,j\in \left\{ 1,\dots ,n\right\} $ with $%
\sum_{j=1}^{n}\left\Vert x_{j}\right\Vert ^{2}=1.$
\end{example}

The following result may be stated as well:

\begin{theorem}[Dragomir, 2008, \protect\cite{II.SSD}]
\label{II.CFOt.6.3.}Let $A_{j}$ be selfadjoint operators with $Sp\left(
A_{j}\right) \subseteq \left[ m,M\right] $ for $j\in \left\{ 1,\dots
,n\right\} $ and for some scalars $m<M.$

(i) \ If $f,g:\left[ m,M\right] \longrightarrow \mathbb{R}$ are continuous
and \textit{synchronous on }$\left[ m,M\right] ,$ \textit{then} 
\begin{align}
& \left\langle \sum_{j=1}^{n}p_{j}f\left( A_{j}\right) g\left( A_{j}\right)
x,x\right\rangle -\left\langle \sum_{j=1}^{n}p_{j}f\left( A_{j}\right)
x,x\right\rangle \cdot \left\langle \sum_{j=1}^{n}p_{j}g\left( A_{j}\right)
x,x\right\rangle  \label{II.CFOe.6.8} \\
& \geq \left[ f\left( \left\langle \sum_{j=1}^{n}p_{j}A_{j}x,x\right\rangle
\right) -\left\langle \sum_{j=1}^{n}p_{j}f\left( A_{j}\right)
x,x\right\rangle \right]  \notag \\
& \times \left[ \left\langle \sum_{j=1}^{n}p_{j}g\left( A_{j}\right)
x,x\right\rangle -g\left( \left\langle
\sum_{j=1}^{n}p_{j}A_{j}x,x\right\rangle \right) \right]  \notag
\end{align}%
for any $p_{j}\geq 0,j\in \left\{ 1,\dots ,n\right\} $ with $%
\sum_{j=1}^{n}p_{j}=1$ and $x\in H$ with $\left\Vert x\right\Vert =1.$
Moreover, if \textit{one is convex while the other is concave on }$\left[ m,M%
\right] ,$ \textit{then the right hand side of (\ref{II.CFOe.6.8}) is
nonnegative.}

(ii) If $f,g$ are a\textit{synchronous on }$\left[ m,M\right] ,$ then 
\begin{align}
& \left\langle \sum_{j=1}^{n}p_{j}f\left( A_{j}\right) x,x\right\rangle
\cdot \left\langle \sum_{j=1}^{n}p_{j}g\left( A_{j}\right) x,x\right\rangle
-\left\langle \sum_{j=1}^{n}p_{j}f\left( A_{j}\right) g\left( A_{j}\right)
x,x\right\rangle  \label{II.CFOe.6.9} \\
& \geq \left[ \left\langle \sum_{j=1}^{n}p_{j}f\left( A_{j}\right)
x,x\right\rangle -f\left( \left\langle
\sum_{j=1}^{n}p_{j}A_{j}x,x\right\rangle \right) \right]  \notag \\
& \times \left[ \left\langle \sum_{j=1}^{n}p_{j}g\left( A_{j}\right)
x,x\right\rangle -g\left( \left\langle
\sum_{j=1}^{n}p_{j}A_{j}x,x\right\rangle \right) \right]  \notag
\end{align}%
for any $p_{j}\geq 0,j\in \left\{ 1,\dots ,n\right\} $ with $%
\sum_{j=1}^{n}p_{j}=1$ and $x\in H$ with $\left\Vert x\right\Vert =1.$
Moreover, if either both of them are convex or both of them are concave on $%
\left[ m,M\right] $, then the right hand side of (\ref{II.CFOe.6.9}) is
nonnegative as well.
\end{theorem}

\begin{proof}
Follows from Theorem \ref{II.CFOt.6.1} on choosing $x_{j}=\sqrt{p_{j}}\cdot
x,$ $j\in \left\{ 1,\dots ,n\right\} ,$ where $p_{j}\geq 0,j\in \left\{
1,\dots ,n\right\} ,$ $\sum_{j=1}^{n}p_{j}=1$ and $x\in H,$ with $\left\Vert
x\right\Vert =1.$

Also, the positivity of the right hand term in (\ref{II.CFOe.6.8}) follows
by the Jensen's type inequality from the inequality (\ref{II.CFOe.6.2.a})
for the same choices, namely $x_{j}=\sqrt{p_{j}}\cdot x,$ $j\in \left\{
1,\dots ,n\right\} ,$ where $p_{j}\geq 0,j\in \left\{ 1,\dots ,n\right\} ,$ $%
\sum_{j=1}^{n}p_{j}=1$ and $x\in H,$ with $\left\Vert x\right\Vert =1$. The
details are omitted.
\end{proof}

Finally, we can list some particular inequalities that may be of interest
for applications. They improve some result obtained above:

\begin{example}
\textbf{a.} Assume that $A_{j},j\in \left\{ 1,\dots ,n\right\} $ are
positive operators on the Hilbert space $H$ and $p_{j}\geq 0,j\in \left\{
1,\dots ,n\right\} $ with $\sum_{j=1}^{n}p_{j}=1.$ If $p\in \left(
0,1\right) $ and $q\in \left( 1,\infty \right) ,$ then for each $x\in H$
with $\left\Vert x\right\Vert =1$ we have the inequality 
\begin{align}
& \left\langle \sum_{j=1}^{n}p_{j}A_{j}^{p+q}x,x\right\rangle -\left\langle
\sum_{j=1}^{n}p_{j}A_{j}^{p}x,x\right\rangle \cdot \left\langle
\sum_{j=1}^{n}p_{j}A_{j}^{q}x,x\right\rangle  \label{II.CFOe.6.13} \\
& \geq \left[ \left\langle \sum_{j=1}^{n}p_{j}A_{j}^{q}x,x\right\rangle
-\left\langle \sum_{j=1}^{n}p_{j}A_{j}x,x\right\rangle ^{q}\right]  \notag \\
& \times \left[ \left\langle \sum_{j=1}^{n}p_{j}A_{j}x,x\right\rangle
^{p}-\left\langle \sum_{j=1}^{n}p_{j}A_{j}^{p}x,x\right\rangle \right] 
\notag \\
& \geq 0.  \notag
\end{align}%
If $A_{j},j\in \left\{ 1,\dots ,n\right\} $ are positive definite and $%
p>1,q<0,$ then 
\begin{align}
& \left\langle \sum_{j=1}^{n}p_{j}A_{j}^{p}x,x\right\rangle \cdot
\left\langle \sum_{j=1}^{n}p_{j}A_{j}^{q}x,x\right\rangle -\left\langle
\sum_{j=1}^{n}p_{j}A_{j}^{p+q}x,x\right\rangle  \label{II.CFOe.6.14} \\
& \geq \left[ \left\langle \sum_{j=1}^{n}p_{j}A_{j}^{q}x,x\right\rangle
-\left\langle \sum_{j=1}^{n}p_{j}A_{j}x,x\right\rangle ^{q}\right]  \notag \\
& \times \left[ \left\langle \sum_{j=1}^{n}p_{j}A_{j}^{p}x,x\right\rangle
-\left\langle \sum_{j=1}^{n}p_{j}A_{j}x,x\right\rangle ^{p}\right]  \notag \\
& \geq 0  \notag
\end{align}%
for each $x\in H$ with $\left\Vert x\right\Vert =1.$

\textbf{b. }Assume that $A_{j}$, $j\in \left\{ 1,\dots ,n\right\} $ are
positive definite and $p>1.$ Then 
\begin{align}
& \left\langle \sum_{j=1}^{n}p_{j}A_{j}^{p}\log A_{j}x,x\right\rangle
-\left\langle \sum_{j=1}^{n}p_{j}A_{j}^{p}x,x\right\rangle \cdot
\left\langle \sum_{j=1}^{n}p_{j}\log A_{j}x,x\right\rangle
\label{II.CFOe.6.15} \\
& \geq \left[ \left\langle \sum_{j=1}^{n}p_{j}A_{j}^{p}x,x\right\rangle
-\left\langle \sum_{j=1}^{n}p_{j}A_{j}x,x\right\rangle ^{p}\right]  \notag \\
& \times \left[ \log \left\langle \sum_{j=1}^{n}p_{j}A_{j}x,x\right\rangle
-\left\langle \sum_{j=1}^{n}p_{j}\log A_{j}x,x\right\rangle \right]  \notag
\\
& \geq 0  \notag
\end{align}%
for each $x\in H$ with $\left\Vert x\right\Vert =1.$
\end{example}

\section{Gr\"{u}ss Inequality}

\subsection{Some Elementary Inequalities of Gr\"{u}ss Type}

In 1935, G. Gr\"{u}ss \cite{II.9b} proved the following integral inequality
which gives an approximation of the integral of the product in terms of the
product of the integrals as follows: 
\begin{align}
& \left\vert \frac{1}{b-a}\int_{a}^{b}f\left( x\right) g\left( x\right) dx-%
\frac{1}{b-a}\int_{a}^{b}f\left( x\right) dx\cdot \frac{1}{b-a}%
\int_{a}^{b}g\left( x\right) dx\right\vert  \label{II.1.1} \\
& \leq \frac{1}{4}\left( \Phi -\phi \right) \left( \Gamma -\gamma \right) , 
\notag
\end{align}%
where $f$, $g:\left[ a,b\right] \rightarrow \mathbb{R}$ are integrable on $%
\left[ a,b\right] $ and satisfy the condition 
\begin{equation}
\phi \leq f\left( x\right) \leq \Phi \text{,\qquad\ }\gamma \leq g\left(
x\right) \leq \Gamma  \label{II.1.2}
\end{equation}%
for each $x\in \left[ a,b\right] ,$ where $\phi ,\Phi ,\gamma ,\Gamma $ are
given real constants.

Moreover, the constant $\frac{1}{4}$ is sharp in the sense that it cannot be
replaced by a smaller one.

In 1950, M. Biernacki, H. Pidek and C. Ryll-Nardjewski \cite[Chapter X]%
{II.11b} established the following discrete version of Gr\"{u}ss' inequality:

Let $a=\left( a_{1},\dots ,a_{n}\right) ,\;b=\left( b_{1},\dots
,b_{n}\right) $ be two $n-$tuples of real numbers such that $r\leq a_{i}\leq
R$ and $s\leq b_{i}\leq S$ for $i=1,\dots ,n.$ Then one has 
\begin{equation}
\left\vert \frac{1}{n}\sum_{i=1}^{n}a_{i}b_{i}-\frac{1}{n}%
\sum_{i=1}^{n}a_{i}\cdot \frac{1}{n}\sum_{i=1}^{n}b_{i}\right\vert \leq 
\frac{1}{n}\left[ \frac{n}{2}\right] \left( 1-\frac{1}{n}\left[ \frac{n}{2}%
\right] \right) \left( R-r\right) \left( S-s\right) ,  \label{II.1.3}
\end{equation}%
where $\left[ x\right] $ denotes the integer part of $x,\;x\in \mathbb{R}.$

For a simple proof of (\ref{II.1.1}) as well as for some other integral
inequalities of Gr\"{u}ss type, see Chapter X of the recent book \cite%
{II.11b}. For other related results see the papers \cite{II.A1}-\cite%
{II.a.A3}, \cite{II.C3}-\cite{II.C1}, \cite{II.CD1}-\cite{II.CD3}, \cite%
{II.1b}-\cite{II.6b}, \cite{II.10b}, \cite{II.Pa1}, \cite{II.ZC} and the
references therein.

\subsection{An Inequality of Gr\"{u}ss' Type for One Operator}

The following result may be stated:

\begin{theorem}[Dragomir, 2008, \protect\cite{II.SSDG}]
\label{II.t.2.1}Let $A$ be a selfadjoint operator on the Hilbert space $%
\left( H;\left\langle .,.\right\rangle \right) $ and assume that $Sp\left(
A\right) \subseteq \left[ m,M\right] $ for some scalars $m<M.$ If $f$ and $g$
are continuous on $\left[ m,M\right] $ and $\gamma :=\min_{t\in \left[ m,M%
\right] }f\left( t\right) $ and $\Gamma :=\max_{t\in \left[ m,M\right]
}f\left( t\right) $ then 
\begin{align}
& \left\vert \left\langle f\left( A\right) g\left( A\right) y,y\right\rangle
-\left\langle f\left( A\right) y,y\right\rangle \cdot \left\langle g\left(
A\right) x,x\right\rangle \right.  \label{II.e.2.1.a} \\
& \left. -\frac{\gamma +\Gamma }{2}\left[ \left\langle g\left( A\right)
y,y\right\rangle -\left\langle g\left( A\right) x,x\right\rangle \right]
\right\vert  \notag \\
& \leq \frac{1}{2}\cdot \left( \Gamma -\gamma \right) \left[ \left\Vert
g\left( A\right) y\right\Vert ^{2}+\left\langle g\left( A\right)
x,x\right\rangle ^{2}-2\left\langle g\left( A\right) x,x\right\rangle
\left\langle g\left( A\right) y,y\right\rangle \right] ^{1/2}  \notag
\end{align}%
for any $x,y\in H$ with $\left\Vert x\right\Vert =\left\Vert y\right\Vert
=1. $
\end{theorem}

\begin{proof}
First of all, observe that, for each $\lambda \in \mathbb{R}$ and $x,y\in H,$
$\left\Vert x\right\Vert =\left\Vert y\right\Vert =1$ we have the identity 
\begin{align}
& \left\langle \left( f\left( A\right) -\lambda \cdot 1_{H}\right) \left(
g\left( A\right) -\left\langle g\left( A\right) x,x\right\rangle \cdot
1_{H}\right) y,y\right\rangle  \label{II.e.2.2} \\
& =\left\langle f\left( A\right) g\left( A\right) y,y\right\rangle -\lambda
\cdot \left[ \left\langle g\left( A\right) y,y\right\rangle -\left\langle
g\left( A\right) x,x\right\rangle \right]  \notag \\
& -\left\langle g\left( A\right) x,x\right\rangle \left\langle f\left(
A\right) y,y\right\rangle .  \notag
\end{align}

Taking the modulus in (\ref{II.e.2.2}) we have 
\begin{align}
& \left\vert \left\langle f\left( A\right) g\left( A\right) y,y\right\rangle
-\lambda \cdot \left[ \left\langle g\left( A\right) y,y\right\rangle
-\left\langle g\left( A\right) x,x\right\rangle \right] \right.
\label{II.e.2.3.a} \\
& \qquad \qquad \left. -\left\langle g\left( A\right) x,x\right\rangle
\left\langle f\left( A\right) y,y\right\rangle \right\vert  \notag \\
& =\left\vert \left\langle \left( g\left( A\right) -\left\langle g\left(
A\right) x,x\right\rangle \cdot 1_{H}\right) y,\left( f\left( A\right)
-\lambda \cdot 1_{H}\right) y\right\rangle \right\vert  \notag \\
& \leq \left\Vert g\left( A\right) y-\left\langle g\left( A\right)
x,x\right\rangle y\right\Vert \left\Vert f\left( A\right) y-\lambda
y\right\Vert  \notag \\
& =\left[ \left\Vert g\left( A\right) y\right\Vert ^{2}+\left\langle g\left(
A\right) x,x\right\rangle ^{2}-2\left\langle g\left( A\right)
x,x\right\rangle \left\langle g\left( A\right) y,y\right\rangle \right]
^{1/2}  \notag \\
& \qquad \qquad \times \left\Vert f\left( A\right) y-\lambda y\right\Vert 
\notag \\
& \leq \left[ \left\Vert g\left( A\right) y\right\Vert ^{2}+\left\langle
g\left( A\right) x,x\right\rangle ^{2}-2\left\langle g\left( A\right)
x,x\right\rangle \left\langle g\left( A\right) y,y\right\rangle \right]
^{1/2}  \notag \\
& \qquad \qquad \times \left\Vert f\left( A\right) -\lambda \cdot
1_{H}\right\Vert  \notag
\end{align}%
for any $x,y\in H,$ $\left\Vert x\right\Vert =\left\Vert y\right\Vert =1.$

Now, since $\gamma =\min_{t\in \left[ m,M\right] }f\left( t\right) $ and $%
\Gamma =\max_{t\in \left[ m,M\right] }f\left( t\right) ,$ then by the
property (\ref{P}) we have that $\gamma \leq \left\langle f\left( A\right)
y,y\right\rangle \leq \Gamma $ for each $y\in H$ with $\left\Vert
y\right\Vert =1$ which is clearly equivalent with 
\begin{equation*}
\left\vert \left\langle f\left( A\right) y,y\right\rangle -\frac{\gamma
+\Gamma }{2}\left\Vert y\right\Vert ^{2}\right\vert \leq \frac{1}{2}\left(
\Gamma -\gamma \right)
\end{equation*}%
or with 
\begin{equation*}
\left\vert \left\langle \left( f\left( A\right) -\frac{\gamma +\Gamma }{2}%
1_{H}\right) y,y\right\rangle \right\vert \leq \frac{1}{2}\left( \Gamma
-\gamma \right)
\end{equation*}%
for each $y\in H$ with $\left\Vert y\right\Vert =1.$

Taking the supremum in this inequality we get 
\begin{equation*}
\left\Vert f\left( A\right) -\frac{\gamma +\Gamma }{2}\cdot 1_{H}\right\Vert
\leq \frac{1}{2}\left( \Gamma -\gamma \right) ,
\end{equation*}%
which together with the inequality (\ref{II.e.2.3.a}) applied for $\lambda =%
\frac{\gamma +\Gamma }{2}$ produces the desired result (\ref{II.e.2.1.a}).
\end{proof}

As a particular case of interest we can derive from the above theorem the
following result of Gr\"{u}ss' type:

\begin{corollary}[Dragomir, 2008, \protect\cite{II.SSDG}]
\label{II.c.2.2}With the assumptions in Theorem \ref{II.t.2.1} we have 
\begin{align}
& \left\vert \left\langle f\left( A\right) g\left( A\right) x,x\right\rangle
-\left\langle f\left( A\right) x,x\right\rangle \cdot \left\langle g\left(
A\right) x,x\right\rangle \right\vert  \label{II.e.2.4} \\
& \leq \frac{1}{2}\cdot \left( \Gamma -\gamma \right) \left[ \left\Vert
g\left( A\right) x\right\Vert ^{2}-\left\langle g\left( A\right)
x,x\right\rangle ^{2}\right] ^{1/2}\left( \leq \frac{1}{4}\left( \Gamma
-\gamma \right) \left( \Delta -\delta \right) \right)  \notag
\end{align}%
for each $x\in H$ with $\left\Vert x\right\Vert =1,$ where $\delta
:=\min_{t\in \left[ m,M\right] }g\left( t\right) $ and $\Delta :=\max_{t\in %
\left[ m,M\right] }g\left( t\right) .$
\end{corollary}

\begin{proof}
The first inequality follows from (\ref{II.e.2.1.a}) by putting $y=x.$

Now, if we write the first inequality in (\ref{II.e.2.4}) for $f=g$ we get 
\begin{align*}
0& \leq \left\Vert g\left( A\right) x\right\Vert ^{2}-\left\langle g\left(
A\right) x,x\right\rangle ^{2}=\left\langle g^{2}\left( A\right)
x,x\right\rangle -\left\langle g\left( A\right) x,x\right\rangle ^{2} \\
& \leq \frac{1}{2}\left( \Delta -\delta \right) \left[ \left\Vert g\left(
A\right) x\right\Vert ^{2}-\left\langle g\left( A\right) x,x\right\rangle
^{2}\right] ^{1/2}
\end{align*}%
which implies that 
\begin{equation*}
\left[ \left\Vert g\left( A\right) x\right\Vert ^{2}-\left\langle g\left(
A\right) x,x\right\rangle ^{2}\right] ^{1/2}\leq \frac{1}{2}\left( \Delta
-\delta \right)
\end{equation*}%
for each $x\in H$ with $\left\Vert x\right\Vert =1.$

This together with the first part of (\ref{II.e.2.4}) proves the desired
bound.
\end{proof}

The following particular cases that hold for power function are of interest:

\begin{example}
Let $A$ be a selfadjoint operator with $Sp\left( A\right) \subseteq \left[
m,M\right] $ for some scalars $m<M.$

If \ $A$ is positive $\left( m\geq 0\right) $ and $p,q>0,$ then 
\begin{align}
(0& \leq )\left\langle A^{p+q}x,x\right\rangle -\left\langle
A^{p}x,x\right\rangle \cdot \left\langle A^{q}x,x\right\rangle
\label{II.e.2.5} \\
& \leq \frac{1}{2}\cdot \left( M^{p}-m^{p}\right) \left[ \left\Vert
A^{q}x\right\Vert ^{2}-\left\langle A^{q}x,x\right\rangle ^{2}\right] ^{1/2}
\notag \\
& \left[ \leq \frac{1}{4}\cdot \left( M^{p}-m^{p}\right) \left(
M^{q}-m^{q}\right) \right]  \notag
\end{align}%
for each $x\in H$ with $\left\Vert x\right\Vert =1.$

If \ $A$ is positive definite $\left( m>0\right) $ and $p,q<0,$ then 
\begin{align}
(0& \leq )\left\langle A^{p+q}x,x\right\rangle -\left\langle
A^{p}x,x\right\rangle \cdot \left\langle A^{q}x,x\right\rangle
\label{II.e.2.6} \\
& \leq \frac{1}{2}\cdot \frac{M^{-p}-m^{-p}}{M^{-p}m^{-p}}\left[ \left\Vert
A^{q}x\right\Vert ^{2}-\left\langle A^{q}x,x\right\rangle ^{2}\right] ^{1/2}
\notag \\
& \left[ \leq \frac{1}{4}\cdot \frac{M^{-p}-m^{-p}}{M^{-p}m^{-p}}\frac{%
M^{-q}-m^{-q}}{M^{-q}m^{-q}}\right]  \notag
\end{align}%
for each $x\in H$ with $\left\Vert x\right\Vert =1.$

If \ $A$ is positive definite $\left( m>0\right) $ and $p<0,$ $q>0$ then 
\begin{align}
(0& \leq )\left\langle A^{p}x,x\right\rangle \cdot \left\langle
A^{q}x,x\right\rangle -\left\langle A^{p+q}x,x\right\rangle  \label{II.e.2.7}
\\
& \leq \frac{1}{2}\cdot \frac{M^{-p}-m^{-p}}{M^{-p}m^{-p}}\left[ \left\Vert
A^{q}x\right\Vert ^{2}-\left\langle A^{q}x,x\right\rangle ^{2}\right] ^{1/2}
\notag \\
& \left[ \leq \frac{1}{4}\cdot \frac{M^{-p}-m^{-p}}{M^{-p}m^{-p}}\left(
M^{q}-m^{q}\right) \right]  \notag
\end{align}%
for each $x\in H$ with $\left\Vert x\right\Vert =1.$

If \ $A$ is positive definite $\left( m>0\right) $ and $p>0,$ $q<0$ then 
\begin{align}
(0& \leq )\left\langle A^{p}x,x\right\rangle \cdot \left\langle
A^{q}x,x\right\rangle -\left\langle A^{p+q}x,x\right\rangle  \label{II.e.2.8}
\\
& \leq \frac{1}{2}\cdot \left( M^{p}-m^{p}\right) \left[ \left\Vert
A^{q}x\right\Vert ^{2}-\left\langle A^{q}x,x\right\rangle ^{2}\right] ^{1/2}
\notag \\
& \left[ \leq \frac{1}{4}\cdot \left( M^{p}-m^{p}\right) \frac{M^{-q}-m^{-q}%
}{M^{-q}m^{-q}}\right]  \notag
\end{align}%
for each $x\in H$ with $\left\Vert x\right\Vert =1.$
\end{example}

We notice that the positivity of the quantities in the left hand side of the
above inequalities (\ref{II.e.2.5})-(\ref{II.e.2.8}) follows from the
Theorem \ref{II.CFOt.2.1}.

The following particular cases when one function is a power while the second
is the logarithm are of interest as well:

\begin{example}
Let $A$ be a positive definite operator with $Sp\left( A\right) \subseteq %
\left[ m,M\right] $ for some scalars $0<m<M.$

If $p>0$ then 
\begin{align}
(0& \leq )\left\langle A^{p}\ln Ax,x\right\rangle -\left\langle
A^{p}x,x\right\rangle \cdot \left\langle \ln Ax,x\right\rangle
\label{II.e.2.9} \\
& \leq \left\{ 
\begin{array}{l}
\frac{1}{2}\cdot \left( M^{p}-m^{p}\right) \left[ \left\Vert \ln
Ax\right\Vert ^{2}-\left\langle \ln Ax,x\right\rangle ^{2}\right] ^{1/2} \\ 
\\ 
\ln \sqrt{\frac{M}{m}}\cdot \left[ \left\Vert A^{p}x\right\Vert
^{2}-\left\langle A^{p}x,x\right\rangle ^{2}\right] ^{1/2}%
\end{array}%
\right.  \notag \\
& \left[ \leq \frac{1}{2}\cdot \left( M^{p}-m^{p}\right) \ln \sqrt{\frac{M}{m%
}}\right]  \notag
\end{align}%
for each $x\in H$ with $\left\Vert x\right\Vert =1.$

If $p<0$ then 
\begin{align}
(0& \leq )\left\langle A^{p}x,x\right\rangle \cdot \left\langle \ln
Ax,x\right\rangle -\left\langle A^{p}\ln Ax,x\right\rangle  \label{II.e.2.10}
\\
& \leq \left\{ 
\begin{array}{l}
\frac{1}{2}\cdot \frac{M^{-p}-m^{-p}}{M^{-p}m^{-p}}\left[ \left\Vert \ln
Ax\right\Vert ^{2}-\left\langle \ln Ax,x\right\rangle ^{2}\right] ^{1/2} \\ 
\\ 
\ln \sqrt{\frac{M}{m}}\cdot \left[ \left\Vert A^{p}x\right\Vert
^{2}-\left\langle A^{p}x,x\right\rangle ^{2}\right] ^{1/2}%
\end{array}%
\right.  \notag \\
& \left[ \leq \frac{1}{2}\cdot \frac{M^{-p}-m^{-p}}{M^{-p}m^{-p}}\ln \sqrt{%
\frac{M}{m}}\right]  \notag
\end{align}%
for each $x\in H$ with $\left\Vert x\right\Vert =1.$
\end{example}

\subsection{An Inequality of Gr\"{u}ss' Type for $n$ Operators}

The following multiple operator version of Theorem \ref{II.t.2.1} holds:

\begin{theorem}[Dragomir, 2008, \protect\cite{II.SSDG}]
\label{II.t.3.1}Let $A_{j}$ be selfadjoint operators with $Sp\left(
A_{j}\right) \subseteq \left[ m,M\right] $ for $j\in \left\{ 1,\dots
,n\right\} $ and for some scalars $m<M.$ If $f,g:\left[ m,M\right]
\longrightarrow \mathbb{R}$ are continuous and $\gamma :=\min_{t\in \left[
m,M\right] }f\left( t\right) $ and $\Gamma :=\max_{t\in \left[ m,M\right]
}f\left( t\right) $ then 
\begin{align}
& \left\vert \sum_{j=1}^{n}\left\langle f\left( A_{j}\right) g\left(
A_{j}\right) y_{j},y_{j}\right\rangle -\sum_{j=1}^{n}\left\langle f\left(
A_{j}\right) y_{j},y_{j}\right\rangle \cdot \sum_{j=1}^{n}\left\langle
g\left( A_{j}\right) x_{j},x_{j}\right\rangle \right.  \label{II.e.3.1} \\
& \left. -\frac{\gamma +\Gamma }{2}\left[ \sum_{j=1}^{n}\left\langle g\left(
A_{j}\right) y_{j},y_{j}\right\rangle -\sum_{j=1}^{n}\left\langle g\left(
A_{j}\right) x_{j},x_{j}\right\rangle \right] \right\vert  \notag \\
& \leq \frac{1}{2}\left( \Gamma -\gamma \right) \left[ \sum_{j=1}^{n}\left%
\Vert g\left( A_{j}\right) y_{j}\right\Vert ^{2}+\left(
\sum_{j=1}^{n}\left\langle g\left( A_{j}\right) x_{j},x_{j}\right\rangle
\right) ^{2}\right.  \notag \\
& -2\left. \sum_{j=1}^{n}\left\langle g\left( A_{j}\right)
x_{j},x_{j}\right\rangle \sum_{j=1}^{n}\left\langle g\left( A_{j}\right)
y_{j},y_{j}\right\rangle \right] ^{\frac{1}{2}}  \notag
\end{align}%
for each $x_{j},y_{j}\in H,j\in \left\{ 1,\dots ,n\right\} $ with $%
\sum_{j=1}^{n}\left\Vert x_{j}\right\Vert ^{2}=\sum_{j=1}^{n}\left\Vert
y_{j}\right\Vert ^{2}=1.$
\end{theorem}

\begin{proof}
Follows from Theorem \ref{II.t.2.1}.
\end{proof}

The following particular case provides a refinement of the Mond-Pe\v{c}ari%
\'{c} result.

\begin{corollary}[Dragomir, 2008, \protect\cite{II.SSDG}]
\label{II.c.3.2}With the assumptions of Theorem \ref{II.t.3.1} we have 
\begin{align}
& \left\vert \sum_{j=1}^{n}\left\langle f\left( A_{j}\right) g\left(
A_{j}\right) x_{j},x_{j}\right\rangle -\sum_{j=1}^{n}\left\langle f\left(
A_{j}\right) x_{j},x_{j}\right\rangle \cdot \sum_{j=1}^{n}\left\langle
g\left( A_{j}\right) x_{j},x_{j}\right\rangle \right\vert  \label{II.e.3.2}
\\
& \leq \frac{1}{2}\cdot \left( \Gamma -\gamma \right) \left[
\sum_{j=1}^{n}\left\Vert g\left( A_{j}\right) x_{j}\right\Vert ^{2}-\left(
\sum_{j=1}^{n}\left\langle g\left( A_{j}\right) x_{j},x_{j}\right\rangle
\right) ^{2}\right] ^{1/2}  \notag \\
& \left( \leq \frac{1}{4}\left( \Gamma -\gamma \right) \left( \Delta -\delta
\right) \right)  \notag
\end{align}%
for each $x_{j}\in H,j\in \left\{ 1,\dots ,n\right\} $ with $%
\sum_{j=1}^{n}\left\Vert x_{j}\right\Vert ^{2}=1$ where $\delta :=\min_{t\in %
\left[ m,M\right] }g\left( t\right) $ and $\Delta :=\max_{t\in \left[ m,M%
\right] }g\left( t\right) .$
\end{corollary}

\begin{example}
Let $A_{j},$ $j\in \left\{ 1,\dots ,n\right\} $ be a selfadjoint operators
with $Sp\left( A_{j}\right) \subseteq \left[ m,M\right] ,j\in \left\{
1,\dots ,n\right\} $ for some scalars $m<M.$

If \ $A_{j}$ are positive $\left( m\geq 0\right) $ and $p,q>0,$ then 
\begin{align}
(0& \leq )\sum_{j=1}^{n}\left\langle A_{j}^{p+q}x_{j},x_{j}\right\rangle
-\sum_{j=1}^{n}\left\langle A_{j}^{p}x_{j},x_{j}\right\rangle \cdot
\sum_{j=1}^{n}\left\langle A_{j}^{q}x_{j},x_{j}\right\rangle
\label{II.e.3.3} \\
& \leq \frac{1}{2}\cdot \left( M^{p}-m^{p}\right) \left[ \sum_{j=1}^{n}\left%
\Vert A_{j}^{q}x_{j}\right\Vert ^{2}-\left( \sum_{j=1}^{n}\left\langle
A_{j}^{q}x_{j},x_{j}\right\rangle \right) ^{2}\right] ^{1/2}  \notag \\
& \left[ \leq \frac{1}{4}\cdot \left( M^{p}-m^{p}\right) \left(
M^{q}-m^{q}\right) \right]  \notag
\end{align}%
for each $x_{j}\in H,j\in \left\{ 1,\dots ,n\right\} $ with $%
\sum_{j=1}^{n}\left\Vert x_{j}\right\Vert ^{2}=1.$

If \ $A_{j}$ are positive definite $\left( m>0\right) $ and $p,q<0,$ then 
\begin{align}
(0& \leq )\sum_{j=1}^{n}\left\langle A_{j}^{p+q}x_{j},x_{j}\right\rangle
-\sum_{j=1}^{n}\left\langle A_{j}^{p}x_{j},x_{j}\right\rangle \cdot
\sum_{j=1}^{n}\left\langle A_{j}^{q}x_{j},x_{j}\right\rangle
\label{II.e.3.4} \\
& \leq \frac{1}{2}\cdot \frac{M^{-p}-m^{-p}}{M^{-p}m^{-p}}\left[
\sum_{j=1}^{n}\left\Vert A_{j}^{q}x_{j}\right\Vert ^{2}-\left(
\sum_{j=1}^{n}\left\langle A_{j}^{q}x_{j},x_{j}\right\rangle \right) ^{2}%
\right] ^{1/2}  \notag \\
& \left[ \leq \frac{1}{4}\cdot \frac{M^{-p}-m^{-p}}{M^{-p}m^{-p}}\frac{%
M^{-q}-m^{-q}}{M^{-q}m^{-q}}\right]  \notag
\end{align}%
for each $x_{j}\in H,j\in \left\{ 1,\dots ,n\right\} $ with $%
\sum_{j=1}^{n}\left\Vert x_{j}\right\Vert ^{2}=1.$

If \ $A_{j}$ are positive definite $\left( m>0\right) $ and $p<0,$ $q>0$
then 
\begin{align}
(0& \leq )\sum_{j=1}^{n}\left\langle A_{j}^{p}x_{j},x_{j}\right\rangle \cdot
\sum_{j=1}^{n}\left\langle A_{j}^{q}x_{j},x_{j}\right\rangle
-\sum_{j=1}^{n}\left\langle A_{j}^{p+q}x_{j},x_{j}\right\rangle
\label{II.e.3.5} \\
& \leq \frac{1}{2}\cdot \frac{M^{-p}-m^{-p}}{M^{-p}m^{-p}}\left[
\sum_{j=1}^{n}\left\Vert A_{j}^{q}x_{j}\right\Vert ^{2}-\left(
\sum_{j=1}^{n}\left\langle A_{j}^{q}x_{j},x_{j}\right\rangle \right) ^{2}%
\right] ^{1/2}  \notag \\
& \left[ \leq \frac{1}{4}\cdot \frac{M^{-p}-m^{-p}}{M^{-p}m^{-p}}\left(
M^{q}-m^{q}\right) \right]  \notag
\end{align}%
for each $x_{j}\in H,j\in \left\{ 1,\dots ,n\right\} $ with $%
\sum_{j=1}^{n}\left\Vert x_{j}\right\Vert ^{2}=1.$

If \ $A_{j}$ are positive definite $\left( m>0\right) $ and $p>0,$ $q<0$
then 
\begin{align}
(0& \leq )\sum_{j=1}^{n}\left\langle A_{j}^{p}x_{j},x_{j}\right\rangle \cdot
\sum_{j=1}^{n}\left\langle A_{j}^{q}x_{j},x_{j}\right\rangle
-\sum_{j=1}^{n}\left\langle A_{j}^{p+q}x_{j},x_{j}\right\rangle
\label{II.e.3.6} \\
& \leq \frac{1}{2}\cdot \left( M^{p}-m^{p}\right) \left[ \sum_{j=1}^{n}\left%
\Vert A_{j}^{q}x_{j}\right\Vert ^{2}-\left( \sum_{j=1}^{n}\left\langle
A_{j}^{q}x_{j},x_{j}\right\rangle \right) ^{2}\right] ^{1/2}  \notag \\
& \left[ \leq \frac{1}{4}\cdot \left( M^{p}-m^{p}\right) \frac{M^{-q}-m^{-q}%
}{M^{-q}m^{-q}}\right]  \notag
\end{align}%
for each $x_{j}\in H,j\in \left\{ 1,\dots ,n\right\} $ with $%
\sum_{j=1}^{n}\left\Vert x_{j}\right\Vert ^{2}=1.$
\end{example}

We notice that the positivity of the quantities in the left hand side of the
above inequalities (\ref{II.e.3.3})-(\ref{II.e.3.6}) follows from the
Theorem \ref{II.CFOt.2.1}.

The following particular cases when one function is a power while the second
is the logarithm are of interest as well:

\begin{example}
Let $A_{j}$ be positive definite operators with $Sp\left( A_{j}\right)
\subseteq \left[ m,M\right] ,$ $j\in \left\{ 1,\dots ,n\right\} $ for some
scalars $0<m<M.$

If $p>0$ then 
\begin{align}
(0& \leq )\sum_{j=1}^{n}\left\langle A_{j}^{p}\ln
A_{j}x_{j},x_{j}\right\rangle -\sum_{j=1}^{n}\left\langle
A_{j}^{p}x_{j},x_{j}\right\rangle \cdot \sum_{j=1}^{n}\left\langle \ln
A_{j}x_{j},x_{j}\right\rangle  \label{II.e.3.7} \\
& \leq \left\{ 
\begin{array}{l}
\frac{1}{2}\cdot \left( M^{p}-m^{p}\right) \left[ \sum_{j=1}^{n}\left\Vert
\ln A_{j}x_{j}\right\Vert ^{2}-\left( \sum_{j=1}^{n}\left\langle \ln
A_{j}x_{j},x_{j}\right\rangle \right) ^{2}\right] ^{1/2} \\ 
\\ 
\ln \sqrt{\frac{M}{m}}\cdot \left[ \sum_{j=1}^{n}\left\Vert
A_{j}^{p}x_{j}\right\Vert ^{2}-\left( \sum_{j=1}^{n}\left\langle
A_{j}^{p}x_{j},x_{j}\right\rangle \right) ^{2}\right] ^{1/2}%
\end{array}%
\right.  \notag \\
& \left[ \leq \frac{1}{2}\cdot \left( M^{p}-m^{p}\right) \ln \sqrt{\frac{M}{m%
}}\right]  \notag
\end{align}%
for each $x_{j}\in H,j\in \left\{ 1,\dots ,n\right\} $ with $%
\sum_{j=1}^{n}\left\Vert x_{j}\right\Vert ^{2}=1.$

If $p<0$ then 
\begin{align}
(0& \leq )\sum_{j=1}^{n}\left\langle A_{j}^{p}x_{j},x_{j}\right\rangle \cdot
\sum_{j=1}^{n}\left\langle \ln A_{j}x_{j},x_{j}\right\rangle
-\sum_{j=1}^{n}\left\langle A_{j}^{p}\ln A_{j}x_{j},x_{j}\right\rangle
\label{II.e.3.8} \\
& \leq \left\{ 
\begin{array}{l}
\frac{1}{2}\frac{M^{-p}-m^{-p}}{M^{-p}m^{-p}}\left[ \sum_{j=1}^{n}\left\Vert
\ln A_{j}x_{j}\right\Vert ^{2}-\left( \sum_{j=1}^{n}\left\langle \ln
A_{j}x_{j},x_{j}\right\rangle \right) ^{2}\right] ^{1/2} \\ 
\\ 
\ln \sqrt{\frac{M}{m}}\cdot \left[ \sum_{j=1}^{n}\left\Vert
A_{j}^{p}x_{j}\right\Vert ^{2}-\left( \sum_{j=1}^{n}\left\langle
A_{j}^{p}x_{j},x_{j}\right\rangle \right) ^{2}\right] ^{1/2}%
\end{array}%
\right.  \notag \\
& \left[ \leq \frac{1}{2}\cdot \frac{M^{-p}-m^{-p}}{M^{-p}m^{-p}}\ln \sqrt{%
\frac{M}{m}}\right]  \notag
\end{align}%
for each $x_{j}\in H,j\in \left\{ 1,\dots ,n\right\} $ with $%
\sum_{j=1}^{n}\left\Vert x_{j}\right\Vert ^{2}=1.$
\end{example}

\subsection{Another Inequality of Gr\"{u}ss' Type for $n$ Operators}

The following different result for $n$ operators can be stated as well:

\begin{theorem}[Dragomir, 2008, \protect\cite{II.SSDG}]
\label{II.t.4.1}Let $A_{j}$ be selfadjoint operators with $Sp\left(
A_{j}\right) \subseteq \left[ m,M\right] $ for $j\in \left\{ 1,\dots
,n\right\} $ and for some scalars $m<M.$ If $f$ and $g$ are continuous on $%
\left[ m,M\right] $ and $\gamma :=\min_{t\in \left[ m,M\right] }f\left(
t\right) $ and $\Gamma :=\max_{t\in \left[ m,M\right] }f\left( t\right) $
then for any $p_{j}\geq 0,j\in \left\{ 1,\dots ,n\right\} $ with $%
\sum_{j=1}^{n}p_{j}=1$ we have%
\begin{align}
& \left\vert \left\langle \sum_{k=1}^{n}p_{k}f\left( A_{k}\right) g\left(
A_{k}\right) y,y\right\rangle \right.  \label{II.e.4.1} \\
& -\frac{\gamma +\Gamma }{2}\cdot \left[ \left\langle
\sum_{k=1}^{n}p_{k}g\left( A_{k}\right) y,y\right\rangle -\left\langle
\sum_{j=1}^{n}p_{j}g\left( A_{j}\right) x,x\right\rangle \right]  \notag \\
& \left. -\left\langle \sum_{k=1}^{n}p_{k}f\left( A_{k}\right)
y,y\right\rangle \cdot \left\langle \sum_{j=1}^{n}p_{j}g\left( A_{j}\right)
x,x\right\rangle \right\vert  \notag \\
& \leq \frac{\Gamma -\gamma }{2}\left[ \sum_{k=1}^{n}p_{k}\left\Vert g\left(
A_{k}\right) y\right\Vert ^{2}-2\left\langle \sum_{k=1}^{n}p_{k}g\left(
A_{k}\right) y,y\right\rangle \left\langle \sum_{j=1}^{n}p_{j}g\left(
A_{j}\right) x,x\right\rangle \right.  \notag \\
& \left. +\left\langle \sum_{j=1}^{n}p_{j}g\left( A_{j}\right)
x,x\right\rangle ^{2}\right] ^{1/2},  \notag
\end{align}%
for each $x,y\in H$ with $\left\Vert x\right\Vert =\left\Vert y\right\Vert
=1.$
\end{theorem}

\begin{proof}
Follows from Theorem \ref{II.t.3.1} on choosing $x_{j}=\sqrt{p_{j}}\cdot x,$ 
$y_{j}=\sqrt{p_{j}}\cdot y,$ $j\in \left\{ 1,\dots ,n\right\} ,$ where $%
p_{j}\geq 0,j\in \left\{ 1,\dots ,n\right\} ,$ $\sum_{j=1}^{n}p_{j}=1$ and $%
x,y\in H,$ with $\left\Vert x\right\Vert =\left\Vert y\right\Vert =1.$ The
details are omitted.
\end{proof}

\begin{remark}
The case $n=1$ (therefore $p=1$) in (\ref{II.e.4.1}) provides the result
from Theorem \ref{II.t.2.1}.
\end{remark}

As a particular case of interest we can derive from the above theorem the
following result of Gr\"{u}ss' type:

\begin{corollary}[Dragomir, 2008, \protect\cite{II.SSDG}]
\label{II.c.4.2}With the assumptions of Theorem \ref{II.t.4.1} we have 
\begin{align}
& \left\vert \left\langle \sum_{k=1}^{n}p_{k}f\left( A_{k}\right) g\left(
A_{k}\right) x,x\right\rangle -\left\langle \sum_{k=1}^{n}p_{k}f\left(
A_{k}\right) x,x\right\rangle \cdot \left\langle \sum_{k=1}^{n}p_{k}g\left(
A_{k}\right) x,x\right\rangle \right\vert  \label{II.e.4.4} \\
& \leq \frac{\Gamma -\gamma }{2}\left( \sum_{k=1}^{n}p_{k}\left\Vert g\left(
A_{k}\right) x\right\Vert ^{2}-\left\langle \sum_{k=1}^{n}p_{k}g\left(
A_{k}\right) x,x\right\rangle ^{2}\right) ^{1/2}  \notag \\
& \left[ \leq \frac{1}{4}\cdot \left( \Gamma -\gamma \right) \left( \Delta
-\delta \right) \right]  \notag
\end{align}%
for each $x\in H$ with $\left\Vert x\right\Vert =1,$ where $\delta
:=\min_{t\in \left[ m,M\right] }g\left( t\right) $ and $\Delta :=\max_{t\in %
\left[ m,M\right] }g\left( t\right) .$
\end{corollary}

\begin{proof}
It is similar with the proof from Corollary \ref{II.c.2.2} and the details
are omitted.
\end{proof}

The following particular cases that hold for power function are of interest:

\begin{example}
Let $A_{j},$ $j\in \left\{ 1,\dots ,n\right\} $ be a selfadjoint operators
with $Sp\left( A_{j}\right) \subseteq \left[ m,M\right] ,j\in \left\{
1,\dots ,n\right\} $ for some scalars $m<M$ and $p_{j}\geq 0,j\in \left\{
1,\dots ,n\right\} $ with $\sum_{j=1}^{n}p_{j}=1.$

If \ $A_{j},$ $j\in \left\{ 1,\dots ,n\right\} $ are positive $\left( m\geq
0\right) $ and $p,q>0,$ then 
\begin{align}
(0& \leq )\left\langle \sum_{k=1}^{n}p_{k}A_{k}^{p+q}x,x\right\rangle
-\left\langle \sum_{k=1}^{n}p_{k}A_{k}^{p}x,x\right\rangle \cdot
\left\langle \sum_{k=1}^{n}p_{k}A_{k}^{q}x,x\right\rangle  \label{II.e.4.5}
\\
& \leq \frac{1}{2}\cdot \left( M^{p}-m^{p}\right) \left[ \sum_{k=1}^{n}p_{k}%
\left\Vert A_{k}^{q}x\right\Vert ^{2}-\left\langle
\sum_{k=1}^{n}p_{k}A_{k}^{q}x,x\right\rangle ^{2}\right] ^{1/2}  \notag \\
& \left[ \leq \frac{1}{4}\cdot \left( M^{p}-m^{p}\right) \left(
M^{q}-m^{q}\right) \right]  \notag
\end{align}%
for each $x\in H$ with $\left\Vert x\right\Vert =1.$

If \ \ $A_{j},$ $j\in \left\{ 1,\dots ,n\right\} $ are positive definite $%
\left( m>0\right) $ and $p,q<0,$ then 
\begin{align}
(0& \leq )\left\langle \sum_{k=1}^{n}p_{k}A_{k}^{p+q}x,x\right\rangle
-\left\langle \sum_{k=1}^{n}p_{k}A_{k}^{p}x,x\right\rangle \cdot
\left\langle \sum_{k=1}^{n}p_{k}A_{k}^{q}x,x\right\rangle  \label{II.e.4.6}
\\
& \leq \frac{1}{2}\cdot \frac{M^{-p}-m^{-p}}{M^{-p}m^{-p}}\left[
\sum_{k=1}^{n}p_{k}\left\Vert A_{k}^{q}x\right\Vert ^{2}-\left\langle
\sum_{k=1}^{n}p_{k}A_{k}^{q}x,x\right\rangle ^{2}\right] ^{1/2}  \notag \\
& \left[ \leq \frac{1}{4}\cdot \frac{M^{-p}-m^{-p}}{M^{-p}m^{-p}}\frac{%
M^{-q}-m^{-q}}{M^{-q}m^{-q}}\right]  \notag
\end{align}%
for each $x\in H$ with $\left\Vert x\right\Vert =1.$

If \ $A_{j},$ $j\in \left\{ 1,\dots ,n\right\} $ are positive definite $%
\left( m>0\right) $ and $p<0,$ $q>0$ then 
\begin{align}
(0& \leq )\left\langle \sum_{k=1}^{n}p_{k}A_{k}^{p}x,x\right\rangle \cdot
\left\langle \sum_{k=1}^{n}p_{k}A_{k}^{q}x,x\right\rangle -\left\langle
\sum_{k=1}^{n}p_{k}A_{k}^{p+q}x,x\right\rangle  \label{II.e.4.7} \\
& \leq \frac{1}{2}\cdot \frac{M^{-p}-m^{-p}}{M^{-p}m^{-p}}\left[
\sum_{k=1}^{n}p_{k}\left\Vert A_{k}^{q}x\right\Vert ^{2}-\left\langle
\sum_{k=1}^{n}p_{k}A_{k}^{q}x,x\right\rangle ^{2}\right] ^{1/2}  \notag \\
& \left[ \leq \frac{1}{4}\cdot \frac{M^{-p}-m^{-p}}{M^{-p}m^{-p}}\left(
M^{q}-m^{q}\right) \right]  \notag
\end{align}%
for each $x\in H$ with $\left\Vert x\right\Vert =1.$

If \ $A_{j},$ $j\in \left\{ 1,\dots ,n\right\} $ are positive definite $%
\left( m>0\right) $ and $p>0,$ $q<0$ then 
\begin{align}
(0& \leq )\left\langle \sum_{k=1}^{n}p_{k}A_{k}^{p}x,x\right\rangle \cdot
\left\langle \sum_{k=1}^{n}p_{k}A_{k}^{q}x,x\right\rangle -\left\langle
\sum_{k=1}^{n}p_{k}A_{k}^{p+q}x,x\right\rangle  \label{II.e.4.8} \\
& \leq \frac{1}{2}\cdot \left( M^{p}-m^{p}\right) \left[ \sum_{k=1}^{n}p_{k}%
\left\Vert A_{k}^{q}x\right\Vert ^{2}-\left\langle
\sum_{k=1}^{n}p_{k}A_{k}^{q}x,x\right\rangle ^{2}\right] ^{1/2}  \notag \\
& \left[ \leq \frac{1}{4}\cdot \left( M^{p}-m^{p}\right) \frac{M^{-q}-m^{-q}%
}{M^{-q}m^{-q}}\right]  \notag
\end{align}%
for each $x\in H$ with $\left\Vert x\right\Vert =1.$
\end{example}

We notice that the positivity of the quantities in the left hand side of the
above inequalities (\ref{II.e.4.5})-(\ref{II.e.4.8}) follows from the
Theorem \ref{II.CFOt.2.1}.

The following particular cases when one function is a power while the second
is the logarithm are of interest as well:

\begin{example}
Let $A_{j},$ $j\in \left\{ 1,\dots ,n\right\} $ be positive definite
operators with $Sp\left( A_{j}\right) \subseteq \left[ m,M\right] ,$ $j\in
\left\{ 1,\dots ,n\right\} $ for some scalars $0<m<M$ and $p_{j}\geq 0,j\in
\left\{ 1,\dots ,n\right\} $ with $\sum_{j=1}^{n}p_{j}=1.$

If $p>0$ then 
\begin{align}
(0& \leq )\left\langle \sum_{k=1}^{n}p_{k}A_{k}^{p}\ln A_{k}x,x\right\rangle
-\left\langle \sum_{k=1}^{n}p_{k}A_{k}^{p}x,x\right\rangle \cdot
\left\langle \sum_{k=1}^{n}p_{k}\ln A_{k}x,x\right\rangle  \label{II.e.4.9}
\\
& \leq \left\{ 
\begin{array}{l}
\frac{1}{2}\cdot \left( M^{p}-m^{p}\right) \cdot \left[ \sum_{k=1}^{n}p_{k}%
\left\Vert \ln A_{k}x\right\Vert ^{2}-\left\langle \sum_{k=1}^{n}p_{k}\ln
A_{k}x,x\right\rangle ^{2}\right] ^{1/2} \\ 
\\ 
\ln \sqrt{\frac{M}{m}}\cdot \left[ \sum_{k=1}^{n}p_{k}\left\Vert
A_{k}^{p}x\right\Vert ^{2}-\left\langle
\sum_{k=1}^{n}p_{k}A_{k}^{p}x,x\right\rangle ^{2}\right] ^{1/2}%
\end{array}%
\right.  \notag \\
& \left[ \leq \frac{1}{2}\cdot \left( M^{p}-m^{p}\right) \ln \sqrt{\frac{M}{m%
}}\right]  \notag
\end{align}%
for each $x\in H$ with $\left\Vert x\right\Vert =1.$

If $p<0$ then 
\begin{align}
(0& \leq )\left\langle \sum_{k=1}^{n}p_{k}A_{k}^{p}x,x\right\rangle \cdot
\left\langle \sum_{k=1}^{n}p_{k}\ln A_{k}x,x\right\rangle -\left\langle
\sum_{k=1}^{n}p_{k}A_{k}^{p}\ln A_{k}x,x\right\rangle  \label{II.e.4.10} \\
& \leq \left\{ 
\begin{array}{l}
\frac{1}{2}\cdot \frac{M^{-p}-m^{-p}}{M^{-p}m^{-p}}\left[
\sum_{k=1}^{n}p_{k}\left\Vert \ln A_{k}x\right\Vert ^{2}-\left\langle
\sum_{k=1}^{n}p_{k}\ln A_{k}x,x\right\rangle ^{2}\right] ^{1/2} \\ 
\\ 
\ln \sqrt{\frac{M}{m}}\cdot \left[ \sum_{k=1}^{n}p_{k}\left\Vert
A_{k}^{p}x\right\Vert ^{2}-\left\langle
\sum_{k=1}^{n}p_{k}A_{k}^{p}x,x\right\rangle ^{2}\right] ^{1/2}%
\end{array}%
\right.  \notag \\
& \left[ \leq \frac{1}{2}\cdot \frac{M^{-p}-m^{-p}}{M^{-p}m^{-p}}\ln \sqrt{%
\frac{M}{m}}\right]  \notag
\end{align}%
for each $x\in H$ with $\left\Vert x\right\Vert =1.$
\end{example}

The following norm inequalities may be stated as well:

\begin{corollary}[Dragomir, 2008, \protect\cite{II.SSDG}]
\label{II.c.4.5}Let $A_{j}$ be selfadjoint operators with $Sp\left(
A_{j}\right) \subseteq \left[ m,M\right] $ for $j\in \left\{ 1,\dots
,n\right\} $ and for some scalars $m<M.$ If $f,g:\left[ m,M\right]
\longrightarrow \mathbb{R}$ are continuous, then for each $p_{j}\geq 0,j\in
\left\{ 1,\dots ,n\right\} $ with $\sum_{j=1}^{n}p_{j}=1$ we have the norm
inequality:%
\begin{equation}
\left\Vert \sum_{j=1}^{n}p_{j}f\left( A_{j}\right) g\left( A_{j}\right)
\right\Vert \leq \left\Vert \sum_{j=1}^{n}p_{j}f\left( A_{j}\right)
\right\Vert \cdot \left\Vert \sum_{j=1}^{n}p_{j}g\left( A_{j}\right)
\right\Vert +\frac{1}{4}\left( \Gamma -\gamma \right) \left( \Delta -\delta
\right) ,  \label{II.e.4.11}
\end{equation}%
where $\gamma :=\min_{t\in \left[ m,M\right] }f\left( t\right) $, $\Gamma
:=\max_{t\in \left[ m,M\right] }f\left( t\right) ,$ $\delta :=\min_{t\in %
\left[ m,M\right] }g\left( t\right) $ and $\Delta :=\max_{t\in \left[ m,M%
\right] }g\left( t\right) .$
\end{corollary}

\begin{proof}
Utilising the inequality (\ref{II.e.4.4}) we deduce the inequality%
\begin{align*}
\left\vert \left\langle \sum_{k=1}^{n}p_{k}f\left( A_{k}\right) g\left(
A_{k}\right) x,x\right\rangle \right\vert & \leq \left\vert \left\langle
\sum_{k=1}^{n}p_{k}f\left( A_{k}\right) x,x\right\rangle \right\vert \cdot
\left\vert \left\langle \sum_{k=1}^{n}p_{k}g\left( A_{k}\right)
x,x\right\rangle \right\vert \\
& +\frac{1}{4}\left( \Gamma -\gamma \right) \left( \Delta -\delta \right)
\end{align*}%
for each $x\in H$ with $\left\Vert x\right\Vert =1.$ Taking the supremum
over $\left\Vert x\right\Vert =1$ we deduce the desired inequality (\ref%
{II.e.4.11}).
\end{proof}

\begin{example}
\textbf{a. }Let $A_{j},$ $j\in \left\{ 1,\dots ,n\right\} $ be a selfadjoint
operators with $Sp\left( A_{j}\right) \subseteq \left[ m,M\right] ,j\in
\left\{ 1,\dots ,n\right\} $ for some scalars $m<M$ and $p_{j}\geq 0,j\in
\left\{ 1,\dots ,n\right\} $ with $\sum_{j=1}^{n}p_{j}=1.$

If \ $A_{j},$ $j\in \left\{ 1,\dots ,n\right\} $ are positive $\left( m\geq
0\right) $ and $p,q>0,$ then 
\begin{equation}
\left\Vert \sum_{k=1}^{n}p_{k}A_{k}^{p+q}\right\Vert \leq \left\Vert
\sum_{k=1}^{n}p_{k}A_{k}^{p}\right\Vert \cdot \left\Vert
\sum_{k=1}^{n}p_{k}A_{k}^{q}\right\Vert +\frac{1}{4}\cdot \left(
M^{p}-m^{p}\right) \left( M^{q}-m^{q}\right) .  \label{II.e.4.12}
\end{equation}

If \ $A_{j},$ $j\in \left\{ 1,\dots ,n\right\} $ are positive definite $%
\left( m>0\right) $ and $p,q<0,$ then 
\begin{equation}
\left\Vert \sum_{k=1}^{n}p_{k}A_{k}^{p+q}\right\Vert \leq \left\Vert
\sum_{k=1}^{n}p_{k}A_{k}^{p}\right\Vert \cdot \left\Vert
\sum_{k=1}^{n}p_{k}A_{k}^{q}\right\Vert +\frac{1}{4}\cdot \frac{M^{-p}-m^{-p}%
}{M^{-p}m^{-p}}\frac{M^{-q}-m^{-q}}{M^{-q}m^{-q}}.  \label{II.e.4.13}
\end{equation}

\textbf{b.} Let $A_{j},$ $j\in \left\{ 1,\dots ,n\right\} $ be positive
definite operators with $Sp\left( A_{j}\right) \subseteq \left[ m,M\right] ,$
$j\in \left\{ 1,\dots ,n\right\} $ for some scalars $0<m<M$ and $p_{j}\geq
0,j\in \left\{ 1,\dots ,n\right\} $ with $\sum_{j=1}^{n}p_{j}=1.$

If $p>0$ then 
\begin{equation}
\left\Vert \sum_{k=1}^{n}p_{k}A_{k}^{p}\ln A_{k}\right\Vert \leq \left\Vert
\sum_{k=1}^{n}p_{k}A_{k}^{p}\right\Vert \cdot \left\Vert
\sum_{k=1}^{n}p_{k}\ln A_{k}\right\Vert +\frac{1}{2}\cdot \left(
M^{p}-m^{p}\right) \ln \sqrt{\frac{M}{m}}.  \label{II.e.4.14}
\end{equation}
\end{example}

\section{More Inequalities of Gr\"{u}ss Type}

\subsection{Some Vectorial Gr\"{u}ss' Type Inequalities}

The following lemmas, that are of interest in their own right, collect some
Gr\"{u}ss type inequalities for vectors in inner product spaces obtained
earlier by the author:

\begin{lemma}[Dragomir, 2003 \& 2004, \protect\cite{II.MGOSSD1}, 
\protect\cite{II.MGOSSD2}]
\label{II.MGOl.3.2}Let $\left( H,\left\langle \cdot ,\cdot \right\rangle
\right) $ be an inner product space over the real or complex number field $%
\mathbb{K}$, $u,v,e\in H,$ $\left\Vert e\right\Vert =1,$ and $\alpha ,\beta
,\gamma ,\delta \in \mathbb{K}$ such that 
\begin{equation}
\func{Re}\left\langle \beta e-u,u-\alpha e\right\rangle \geq 0,\qquad \func{%
Re}\left\langle \delta e-v,v-\gamma e\right\rangle \geq 0
\label{II.MGOe.3.3}
\end{equation}%
or equivalently, 
\begin{equation}
\left\Vert u-\frac{\alpha +\beta }{2}e\right\Vert \leq \frac{1}{2}\left\vert
\beta -\alpha \right\vert ,\qquad \left\Vert v-\frac{\gamma +\delta }{2}%
e\right\Vert \leq \frac{1}{2}\left\vert \delta -\gamma \right\vert .
\label{II.MGOe.3.4}
\end{equation}%
Then 
\begin{align}
& \left\vert \left\langle u,v\right\rangle -\left\langle u,e\right\rangle
\left\langle e,v\right\rangle \right\vert  \label{II.MGOe.3.5} \\
& \leq \frac{1}{4}\cdot \left\vert \beta -\alpha \right\vert \left\vert
\delta -\gamma \right\vert  \notag \\
& -\left\{ 
\begin{array}{l}
\left[ \func{Re}\left\langle \beta e-u,u-\alpha e\right\rangle \func{Re}%
\left\langle \delta e-v,v-\gamma e\right\rangle \right] ^{\frac{1}{2}}, \\ 
\\ 
\left\vert \left\langle u,e\right\rangle -\frac{\alpha +\beta }{2}%
\right\vert \left\vert \left\langle v,e\right\rangle -\frac{\gamma +\delta }{%
2}\right\vert .%
\end{array}%
\right.  \notag
\end{align}
\end{lemma}

The first inequality has been obtained in \cite{II.MGOSSD1} (see also \cite[%
p. 44]{II.SSDM}) while the second result was established in \cite{II.MGOSSD2}
(see also \cite[p. 90]{II.SSDM}). They provide refinements of the earlier
result from \cite{II.2b} where only the first part of the bound, i.e., $%
\frac{1}{4}\left\vert \beta -\alpha \right\vert \left\vert \delta -\gamma
\right\vert $ has been given. Notice that, as pointed out in \cite%
{II.MGOSSD2}, the upper bounds for the Gr\"{u}ss functional incorporated in (%
\ref{II.MGOe.3.5}) cannot be compared in general, meaning that one is better
than the other depending on appropriate choices of the vectors and scalars
involved.

Another result of this type is the following one:

\begin{lemma}[Dragomir, 2004 \& 2006, \protect\cite{II.MGOSSD4}, 
\protect\cite{II.MGOSSD5}]
\label{II.MGOl.3.3}With the assumptions in Lemma \ref{II.MGOl.3.2} and if $%
\func{Re}\left( \beta \overline{\alpha }\right) >0,\func{Re}\left( \delta 
\overline{\gamma }\right) >0$ then 
\begin{align}
& \left\vert \left\langle u,v\right\rangle -\left\langle u,e\right\rangle
\left\langle e,v\right\rangle \right\vert  \label{II.MGOe.3.6} \\
& \leq \left\{ 
\begin{array}{l}
\frac{1}{4}\cdot \frac{\left\vert \beta -\alpha \right\vert \left\vert
\delta -\gamma \right\vert }{\left[ \func{Re}\left( \beta \overline{\alpha }%
\right) \func{Re}\left( \delta \overline{\gamma }\right) \right] ^{\frac{1}{2%
}}}\left\vert \left\langle u,e\right\rangle \left\langle e,v\right\rangle
\right\vert , \\ 
\\ 
\left[ \left( \left\vert \alpha +\beta \right\vert -2\left[ \func{Re}\left(
\beta \overline{\alpha }\right) \right] ^{\frac{1}{2}}\right) \left(
\left\vert \delta +\gamma \right\vert -2\left[ \func{Re}\left( \delta 
\overline{\gamma }\right) \right] ^{\frac{1}{2}}\right) \right] ^{\frac{1}{2}%
} \\ 
\times \left[ \left\vert \left\langle u,e\right\rangle \left\langle
e,v\right\rangle \right\vert \right] ^{\frac{1}{2}}.%
\end{array}%
\right.  \notag
\end{align}
\end{lemma}

The first inequality has been established in \cite{II.MGOSSD4} (see \cite[p.
62]{II.SSDM}) while the second one can be obtained in a canonical manner
from the reverse of the Schwarz inequality given in \cite{II.MGOSSD5}. The
details are omitted.

Finally, another inequality of Gr\"{u}ss type that has been obtained in \cite%
{II.MGOSSD6} (see also \cite[p. 65]{II.SSDM}) can be stated as:

\begin{lemma}[Dragomir, 2004, \protect\cite{II.MGOSSD6}]
\label{II.MGOl.3.4}With the assumptions in Lemma \ref{II.MGOl.3.2} and if $%
\beta \neq -\alpha ,$ $\delta \neq -\gamma $ then 
\begin{align}
& \left\vert \left\langle u,v\right\rangle -\left\langle u,e\right\rangle
\left\langle e,v\right\rangle \right\vert  \label{II.MGOe.3.7} \\
& \leq \frac{1}{4}\cdot \frac{\left\vert \beta -\alpha \right\vert
\left\vert \delta -\gamma \right\vert }{\left[ \left\vert \beta +\alpha
\right\vert \left\vert \delta +\gamma \right\vert \right] ^{\frac{1}{2}}}%
\left[ \left( \left\Vert u\right\Vert +\left\vert \left\langle
u,e\right\rangle \right\vert \right) \left( \left\Vert v\right\Vert
+\left\vert \left\langle v,e\right\rangle \right\vert \right) \right] ^{%
\frac{1}{2}}.  \notag
\end{align}
\end{lemma}

\subsection{Some Inequalities of Gr\"{u}ss' Type for One Operator}

The following results incorporates some new inequalities of Gr\"{u}ss' type
for two functions of a selfadjoint operator.

\begin{theorem}[Dragomir, 2008, \protect\cite{II.SSDG1}]
\label{II.MGOt.4.1}Let $A$ be a selfadjoint operator on the Hilbert space $%
\left( H;\left\langle .,.\right\rangle \right) $ and assume that $Sp\left(
A\right) \subseteq \left[ m,M\right] $ for some scalars $m<M.$ If $f$ and $g$
are continuous on $\left[ m,M\right] $ and $\gamma :=\min_{t\in \left[ m,M%
\right] }f\left( t\right) $, $\Gamma :=\max_{t\in \left[ m,M\right] }f\left(
t\right) $, $\delta :=\min_{t\in \left[ m,M\right] }g\left( t\right) $ and $%
\Delta :=\max_{t\in \left[ m,M\right] }g\left( t\right) $ then 
\begin{align}
& \left\vert \left\langle f\left( A\right) g\left( A\right) x,x\right\rangle
-\left\langle f\left( A\right) x,x\right\rangle \left\langle g\left(
A\right) x,x\right\rangle \right\vert  \label{II.MGOe.4.1} \\
& \leq \frac{1}{4}\left( \Gamma -\gamma \right) \left( \Delta -\delta \right)
\notag \\
& -\left\{ 
\begin{array}{l}
\left[ \left\langle \Gamma x-f\left( A\right) x,f\left( A\right) x-\gamma
x\right\rangle \left\langle \Delta x-g\left( A\right) x,g\left( A\right)
x-\delta x\right\rangle \right] ^{\frac{1}{2}}, \\ 
\\ 
\left\vert \left\langle f\left( A\right) x,x\right\rangle -\frac{\Gamma
+\gamma }{2}\right\vert \left\vert \left\langle g\left( A\right)
x,x\right\rangle -\frac{\Delta +\delta }{2}\right\vert ,%
\end{array}%
\right.  \notag
\end{align}%
for each $x\in H$ with $\left\Vert x\right\Vert =1.$

Moreover if $\gamma $ and $\delta $ are positive, then we also have 
\begin{align}
& \left\vert \left\langle f\left( A\right) g\left( A\right) x,x\right\rangle
-\left\langle f\left( A\right) x,x\right\rangle \left\langle g\left(
A\right) x,x\right\rangle \right\vert  \label{II.MGOe.4.2} \\
& \leq \left\{ 
\begin{array}{l}
\frac{1}{4}\cdot \frac{\left( \Gamma -\gamma \right) \left( \Delta -\delta
\right) }{\sqrt{\Gamma \gamma \Delta \delta }}\left\langle f\left( A\right)
x,x\right\rangle \left\langle g\left( A\right) x,x\right\rangle , \\ 
\\ 
\left( \sqrt{\Gamma }-\sqrt{\gamma }\right) \left( \sqrt{\Delta }-\sqrt{%
\delta }\right) \left[ \left\langle f\left( A\right) x,x\right\rangle
\left\langle g\left( A\right) x,x\right\rangle \right] ^{\frac{1}{2}},%
\end{array}%
\right.  \notag
\end{align}

while for $\Gamma +\gamma ,\Delta +\delta \neq 0$ we have 
\begin{align}
& \left\vert \left\langle f\left( A\right) g\left( A\right) x,x\right\rangle
-\left\langle f\left( A\right) x,x\right\rangle \left\langle g\left(
A\right) x,x\right\rangle \right\vert  \label{II.MGOe.4.3} \\
& \leq \frac{1}{4}\cdot \frac{\left( \Gamma -\gamma \right) \left( \Delta
-\delta \right) }{\left[ \left\vert \Gamma +\gamma \right\vert \left\vert
\Delta +\delta \right\vert \right] ^{\frac{1}{2}}}  \notag \\
& \times \left[ \left( \left\Vert f\left( A\right) x\right\Vert +\left\vert
\left\langle f\left( A\right) x,x\right\rangle \right\vert \right) \left(
\left\Vert g\left( A\right) x\right\Vert +\left\vert \left\langle g\left(
A\right) x,x\right\rangle \right\vert \right) \right] ^{\frac{1}{2}}  \notag
\end{align}%
for each $x\in H$ with $\left\Vert x\right\Vert =1.$
\end{theorem}

\begin{proof}
Since $\gamma :=\min_{t\in \left[ m,M\right] }f\left( t\right) $, $\Gamma
:=\max_{t\in \left[ m,M\right] }f\left( t\right) $, $\delta :=\min_{t\in %
\left[ m,M\right] }g\left( t\right) $ and $\Delta :=\max_{t\in \left[ m,M%
\right] }g\left( t\right) ,$ the by the property (\ref{P}) we have that 
\begin{equation*}
\gamma \cdot 1_{H}\leq f\left( A\right) \leq \Gamma \cdot 1_{H}\text{ \quad
and\quad\ }\delta \cdot 1_{H}\leq g\left( A\right) \leq \Delta \cdot 1_{H}%
\text{ }
\end{equation*}%
in the operator order, which imply that 
\begin{eqnarray}
\left[ f\left( A\right) -\gamma \cdot 1\right] \left[ \Gamma \cdot
1_{H}-f\left( A\right) \right] &\geq &0\text{ \ and }  \label{II.MGOe.4.4} \\
\text{\ }\left[ \Delta \cdot 1_{H}-g\left( A\right) \right] \left[ g\left(
A\right) -\delta \cdot 1_{H}\right] &\geq &0  \notag
\end{eqnarray}%
in the operator order.

We then have from (\ref{II.MGOe.4.4}) 
\begin{equation*}
\left\langle \left[ f\left( A\right) -\gamma \cdot 1\right] \left[ \Gamma
\cdot 1_{H}-f\left( A\right) \right] x,x\right\rangle \geq 0
\end{equation*}%
and 
\begin{equation*}
\left\langle \left[ \Delta \cdot 1_{H}-g\left( A\right) \right] \left[
g\left( A\right) -\delta \cdot 1_{H}\right] x,x\right\rangle \geq 0,
\end{equation*}%
for each $x\in H$ with $\left\Vert x\right\Vert =1,$ which, by the fact that
the involved operators are selfadjoint, are equivalent with the inequalities 
\begin{equation}
\left\langle \Gamma x-f\left( A\right) x,f\left( A\right) x-\gamma
x\right\rangle \geq 0\text{ \ and \ \ }\left\langle \Delta x-g\left(
A\right) x,g\left( A\right) x-\delta x\right\rangle \geq 0,\text{ }
\label{II.MGOe.4.5}
\end{equation}%
for each $x\in H$ with $\left\Vert x\right\Vert =1.$

Now, if we apply Lemma \ref{II.MGOl.3.2} for $u=f\left( A\right) x,$ $%
v=g\left( A\right) x$, $e=x,$ and the real scalars $\Gamma ,\gamma ,\Delta $
and $\delta $ defined in the statement of the theorem, then we can state the
inequality 
\begin{align}
& \left\vert \left\langle f\left( A\right) x,g\left( A\right) x\right\rangle
-\left\langle f\left( A\right) x,x\right\rangle \left\langle x,g\left(
A\right) x\right\rangle \right\vert  \label{II.MGOe.4.6} \\
& \leq \frac{1}{4}\cdot \left( \Gamma -\gamma \right) \left( \Delta -\delta
\right)  \notag \\
& -\left\{ 
\begin{array}{l}
\left[ \func{Re}\left\langle \Gamma x-f\left( A\right) x,f\left( A\right)
x-\gamma x\right\rangle \func{Re}\left\langle \Delta x-g\left( A\right)
x,g\left( A\right) x-\delta x\right\rangle \right] ^{\frac{1}{2}}, \\ 
\\ 
\left\vert \left\langle f\left( A\right) x,x\right\rangle -\frac{\Gamma
+\gamma }{2}\right\vert \left\vert \left\langle g\left( A\right)
x,x\right\rangle -\frac{\Delta +\delta }{2}\right\vert ,%
\end{array}%
\right.  \notag
\end{align}%
for each $x\in H$ with $\left\Vert x\right\Vert =1,$ which is clearly
equivalent with the inequality (\ref{II.MGOe.4.1}).

The inequalities (\ref{II.MGOe.4.2}) and (\ref{II.MGOe.4.3}) follow by Lemma %
\ref{II.MGOl.3.3} and Lemma \ref{II.MGOl.3.4} respectively and the details
are omitted.
\end{proof}

\begin{remark}
\label{II.MGOr.4.2}The first inequality in (\ref{II.MGOe.4.2}) can be
written in a more convenient way as 
\begin{equation}
\left\vert \frac{\left\langle f\left( A\right) g\left( A\right)
x,x\right\rangle }{\left\langle f\left( A\right) x,x\right\rangle
\left\langle g\left( A\right) x,x\right\rangle }-1\right\vert \leq \frac{1}{4%
}\cdot \frac{\left( \Gamma -\gamma \right) \left( \Delta -\delta \right) }{%
\sqrt{\Gamma \gamma \Delta \delta }}  \label{II.MGOe.4.7}
\end{equation}%
for each $x\in H$ with $\left\Vert x\right\Vert =1,$ while the second
inequality has the following equivalent form 
\begin{align}
& \left\vert \frac{\left\langle f\left( A\right) g\left( A\right)
x,x\right\rangle }{\left[ \left\langle f\left( A\right) x,x\right\rangle
\left\langle g\left( A\right) x,x\right\rangle \right] ^{1/2}}-\left[
\left\langle f\left( A\right) x,x\right\rangle \left\langle g\left( A\right)
x,x\right\rangle \right] ^{1/2}\right\vert  \label{II.MGOe.4.8} \\
& \leq \left( \sqrt{\Gamma }-\sqrt{\gamma }\right) \left( \sqrt{\Delta }-%
\sqrt{\delta }\right)  \notag
\end{align}%
for each $x\in H$ with $\left\Vert x\right\Vert =1.$

We know, from \cite{II.SSD} that if $f,g$ are synchronous (asynchronous)
functions on the interval $\left[ m,M\right] ,$ i.e., we recall that 
\begin{equation*}
\left[ f\left( t\right) -f\left( s\right) \right] \left[ g\left( t\right)
-g\left( s\right) \right] \left( \geq \right) \leq 0\text{ \quad for each }%
t,s\in \left[ m,M\right] ,
\end{equation*}%
then we have the inequality 
\begin{equation}
\left\langle f\left( A\right) g\left( A\right) x,x\right\rangle \geq \left(
\leq \right) \left\langle f\left( A\right) x,x\right\rangle \left\langle
g\left( A\right) x,x\right\rangle  \label{II.MGOe.4.9}
\end{equation}%
for each $x\in H$ with $\left\Vert x\right\Vert =1,$provided $f,g$ are
continuous on $\left[ m,M\right] $ and $A$ is a selfadjoint operator with $%
Sp\left( A\right) \subseteq \left[ m,M\right] $.

Therefore, if $f,g$ are synchronous then we have from (\ref{II.MGOe.4.7})
and from (\ref{II.MGOe.4.8}) the following results: 
\begin{equation}
0\leq \frac{\left\langle f\left( A\right) g\left( A\right) x,x\right\rangle 
}{\left\langle f\left( A\right) x,x\right\rangle \left\langle g\left(
A\right) x,x\right\rangle }-1\leq \frac{1}{4}\cdot \frac{\left( \Gamma
-\gamma \right) \left( \Delta -\delta \right) }{\sqrt{\Gamma \gamma \Delta
\delta }}  \label{II.MGOe.4.10}
\end{equation}%
and 
\begin{align}
0& \leq \frac{\left\langle f\left( A\right) g\left( A\right)
x,x\right\rangle }{\left[ \left\langle f\left( A\right) x,x\right\rangle
\left\langle g\left( A\right) x,x\right\rangle \right] ^{1/2}}-\left[
\left\langle f\left( A\right) x,x\right\rangle \left\langle g\left( A\right)
x,x\right\rangle \right] ^{1/2}  \label{II.MGOe.4.11} \\
& \leq \left( \sqrt{\Gamma }-\sqrt{\gamma }\right) \left( \sqrt{\Delta }-%
\sqrt{\delta }\right)  \notag
\end{align}%
for each $x\in H$ with $\left\Vert x\right\Vert =1,$ respectively.

If $f,g$ are asynchronous then 
\begin{equation}
0\leq 1-\frac{\left\langle f\left( A\right) g\left( A\right)
x,x\right\rangle }{\left\langle f\left( A\right) x,x\right\rangle
\left\langle g\left( A\right) x,x\right\rangle }\leq \frac{1}{4}\cdot \frac{%
\left( \Gamma -\gamma \right) \left( \Delta -\delta \right) }{\sqrt{\Gamma
\gamma \Delta \delta }}  \label{II.MGOe.4.12}
\end{equation}%
and 
\begin{align}
0& \leq \left[ \left\langle f\left( A\right) x,x\right\rangle \left\langle
g\left( A\right) x,x\right\rangle \right] ^{1/2}-\frac{\left\langle f\left(
A\right) g\left( A\right) x,x\right\rangle }{\left[ \left\langle f\left(
A\right) x,x\right\rangle \left\langle g\left( A\right) x,x\right\rangle %
\right] ^{1/2}}  \label{II.MGOe.4.13} \\
& \leq \left( \sqrt{\Gamma }-\sqrt{\gamma }\right) \left( \sqrt{\Delta }-%
\sqrt{\delta }\right)  \notag
\end{align}%
for each $x\in H$ with $\left\Vert x\right\Vert =1,$ respectively.
\end{remark}

It is obvious that all the inequalities from Theorem \ref{II.MGOt.4.1} can
be used to obtain reverse inequalities of Gr\"{u}ss' type for various
particular instances of operator functions, see for instance \cite{II.SSDG}.
However we give here only a few provided by the inequalities (\ref%
{II.MGOe.4.10}) and (\ref{II.MGOe.4.11}) above.

\begin{example}
\label{II.MGOex.1}Let $A$ be a selfadjoint operator with $Sp\left( A\right)
\subseteq \left[ m,M\right] $ for some scalars $m<M.$

If \ $A$ is positive $\left( m\geq 0\right) $ and $p,q>0,$ then 
\begin{equation}
0\leq \frac{\left\langle A^{p+q}x,x\right\rangle }{\left\langle
A^{p}x,x\right\rangle \cdot \left\langle A^{q}x,x\right\rangle }-1\leq \frac{%
1}{4}\cdot \frac{\left( M^{p}-m^{p}\right) \left( M^{q}-m^{q}\right) }{M^{%
\frac{p+q}{2}}m^{\frac{p+q}{2}}}  \label{II.MGOe.4.14}
\end{equation}%
and 
\begin{align}
0& \leq \frac{\left\langle A^{p+q}x,x\right\rangle }{\left[ \left\langle
A^{p}x,x\right\rangle \cdot \left\langle A^{q}x,x\right\rangle \right] ^{1/2}%
}-\left[ \left\langle A^{p}x,x\right\rangle \cdot \left\langle
A^{q}x,x\right\rangle \right] ^{1/2}  \label{II.MGOe.4.15} \\
& \leq \left( M^{\frac{p}{2}}-m^{\frac{p}{2}}\right) \left( M^{\frac{q}{2}%
}-m^{\frac{q}{2}}\right)  \notag
\end{align}%
for each $x\in H$ with $\left\Vert x\right\Vert =1.$

If $A$ is positive definite $\left( m>0\right) $ and $p,q<0,$ then 
\begin{equation}
0\leq \frac{\left\langle A^{p+q}x,x\right\rangle }{\left\langle
A^{p}x,x\right\rangle \cdot \left\langle A^{q}x,x\right\rangle }-1\leq \frac{%
1}{4}\cdot \frac{\left( M^{-p}-m^{-p}\right) \left( M^{-q}-m^{-q}\right) }{%
M^{-\frac{p+q}{2}}m^{-\frac{p+q}{2}}}  \label{II.MGOe.4.16}
\end{equation}%
and 
\begin{align}
0& \leq \frac{\left\langle A^{p+q}x,x\right\rangle }{\left[ \left\langle
A^{p}x,x\right\rangle \cdot \left\langle A^{q}x,x\right\rangle \right] ^{1/2}%
}-\left[ \left\langle A^{p}x,x\right\rangle \cdot \left\langle
A^{q}x,x\right\rangle \right] ^{1/2}  \label{II.MGOe.4.17} \\
& \leq \frac{\left( M^{-\frac{p}{2}}-m^{-\frac{p}{2}}\right) \left( M^{-%
\frac{q}{2}}-m^{-\frac{q}{2}}\right) }{M^{-\frac{p+q}{2}}m^{-\frac{p+q}{2}}}
\notag
\end{align}%
for each $x\in H$ with $\left\Vert x\right\Vert =1.$
\end{example}

Similar inequalities may be stated for either $p>0,q<0$ or $p<0,q>0.$ The
details are omitted.

\begin{example}
\label{II.MGOex.2}Let $A$ be a positive definite operator with $Sp\left(
A\right) \subseteq \left[ m,M\right] $ for some scalars $1<m<M.$ If $p>0$
then 
\begin{equation}
0\leq \frac{\left\langle A^{p}\ln Ax,x\right\rangle }{\left\langle
A^{p}x,x\right\rangle \cdot \left\langle \ln Ax,x\right\rangle }-1\leq \frac{%
1}{4}\cdot \frac{\left( M^{p}-m^{p}\right) \ln \frac{M}{m}}{M^{\frac{p}{2}%
}m^{\frac{p}{2}}\sqrt{\ln M\cdot \ln m}}  \label{II.MGOe.4.18}
\end{equation}%
and 
\begin{align}
0& \leq \frac{\left\langle A^{p}\ln Ax,x\right\rangle }{\left[ \left\langle
A^{p}x,x\right\rangle \cdot \left\langle \ln Ax,x\right\rangle \right] ^{1/2}%
}-\left[ \left\langle A^{p}x,x\right\rangle \cdot \left\langle \ln
Ax,x\right\rangle \right] ^{1/2}  \label{II.MGOe.4.19} \\
& \leq \left( M^{\frac{p}{2}}-m^{\frac{p}{2}}\right) \left[ \sqrt{\ln M}-%
\sqrt{\ln m}\right] ,  \notag
\end{align}%
for each $x\in H$ with $\left\Vert x\right\Vert =1.$
\end{example}

\subsection{Some Inequalities of Gr\"{u}ss' Type for $n$ Operators}

The following extension for sequences of operators can be stated:

\begin{theorem}[Dragomir, 2008, \protect\cite{II.SSDG1}]
\label{II.MGOt.5.1} Let $A_{j}$ be selfadjoint operators with $Sp\left(
A_{j}\right) \subseteq \left[ m,M\right] $ for $j\in \left\{ 1,\dots
,n\right\} $ and for some scalars $m<M.$ If $f$ and $g$ are continuous on $%
\left[ m,M\right] $ and $\gamma :=\min_{t\in \left[ m,M\right] }f\left(
t\right) $, $\Gamma :=\max_{t\in \left[ m,M\right] }f\left( t\right) $, $%
\delta :=\min_{t\in \left[ m,M\right] }g\left( t\right) $ and $\Delta
:=\max_{t\in \left[ m,M\right] }g\left( t\right) $ then 
\begin{align}
& \left\vert \sum_{j=1}^{n}\left\langle f\left( A_{j}\right) g\left(
A_{j}\right) x_{j},x_{j}\right\rangle -\sum_{j=1}^{n}\left\langle f\left(
A_{j}\right) x_{j},x_{j}\right\rangle \cdot \sum_{j=1}^{n}\left\langle
g\left( A_{j}\right) x_{j},x_{j}\right\rangle \right\vert
\label{II.MGOe.5.1} \\
& \leq \frac{1}{4}\cdot \left( \Gamma -\gamma \right) \left( \Delta -\delta
\right)  \notag \\
& -\left\{ 
\begin{array}{l}
\left[ \sum\limits_{j=1}^{n}\left\langle \Gamma x_{j}-f\left( A_{j}\right)
x_{j},f\left( A_{j}\right) x_{j}-\gamma x_{j}\right\rangle \right. \\ 
\qquad \qquad \times \left. \sum\limits_{j=1}^{n}\left\langle \Delta
x_{j}-g\left( A_{j}\right) x_{j},g\left( A_{j}\right) x-\delta
x_{j}\right\rangle \right] ^{\frac{1}{2}}, \\ 
\left\vert \sum\limits_{j=1}^{n}\left\langle f\left( A_{j}\right)
x_{j},x_{j}\right\rangle -\frac{\Gamma +\gamma }{2}\right\vert \left\vert
\sum\limits_{j=1}^{n}\left\langle g\left( A_{j}\right)
x_{j},x_{j}\right\rangle -\frac{\Delta +\delta }{2}\right\vert ,%
\end{array}%
\right.  \notag
\end{align}%
for each $x_{j}\in H,j\in \left\{ 1,\dots ,n\right\} $ with $%
\sum_{j=1}^{n}\left\Vert x_{j}\right\Vert ^{2}=1.$

Moreover if $\gamma $ and $\delta $ are positive, then we also have 
\begin{align}
& \left\vert \sum_{j=1}^{n}\left\langle f\left( A_{j}\right) g\left(
A_{j}\right) x_{j},x_{j}\right\rangle -\sum_{j=1}^{n}\left\langle f\left(
A_{j}\right) x_{j},x_{j}\right\rangle \cdot \sum_{j=1}^{n}\left\langle
g\left( A_{j}\right) x_{j},x_{j}\right\rangle \right\vert
\label{II.MGOe.5.2} \\
& \leq \left\{ 
\begin{array}{l}
\frac{1}{4}\cdot \frac{\left( \Gamma -\gamma \right) \left( \Delta -\delta
\right) }{\sqrt{\Gamma \gamma \Delta \delta }}\sum\limits_{j=1}^{n}\left%
\langle f\left( A_{j}\right) x_{j},x_{j}\right\rangle \cdot
\sum\limits_{j=1}^{n}\left\langle g\left( A_{j}\right)
x_{j},x_{j}\right\rangle , \\ 
\\ 
\left( \sqrt{\Gamma }-\sqrt{\gamma }\right) \left( \sqrt{\Delta }-\sqrt{%
\delta }\right) \\ 
\times \left[ \sum\limits_{j=1}^{n}\left\langle f\left( A_{j}\right)
x_{j},x_{j}\right\rangle \cdot \sum\limits_{j=1}^{n}\left\langle g\left(
A_{j}\right) x_{j},x_{j}\right\rangle \right] ^{\frac{1}{2}},%
\end{array}%
\right.  \notag
\end{align}

while for $\Gamma +\gamma ,\Delta +\delta \neq 0$ we have 
\begin{align}
& \left\vert \sum_{j=1}^{n}\left\langle f\left( A_{j}\right) g\left(
A_{j}\right) x_{j},x_{j}\right\rangle -\sum_{j=1}^{n}\left\langle f\left(
A_{j}\right) x_{j},x_{j}\right\rangle \cdot \sum_{j=1}^{n}\left\langle
g\left( A_{j}\right) x_{j},x_{j}\right\rangle \right\vert
\label{II.MGOe.5.3} \\
& \leq \frac{1}{4}\cdot \frac{\left( \Gamma -\gamma \right) \left( \Delta
-\delta \right) }{\left[ \left\vert \Gamma +\gamma \right\vert \left\vert
\Delta +\delta \right\vert \right] ^{\frac{1}{2}}}  \notag \\
& \times \left\lceil \left( \left( \sum_{j=1}^{n}\left\Vert f\left(
A_{j}\right) x_{j}\right\Vert ^{2}\right) ^{1/2}+\left\vert
\sum_{j=1}^{n}\left\langle f\left( A_{j}\right) x_{j},x_{j}\right\rangle
\right\vert \right) \right.  \notag \\
& \left. \times \left( \left( \sum_{j=1}^{n}\left\Vert g\left( A_{j}\right)
x_{j}\right\Vert ^{2}\right) ^{1/2}+\left\vert \sum_{j=1}^{n}\left\langle
g\left( A_{j}\right) x_{j},x_{j}\right\rangle \right\vert \right) \right]
^{1/2},  \notag
\end{align}%
for each $x_{j}\in H,j\in \left\{ 1,\dots ,n\right\} $ with $%
\sum_{j=1}^{n}\left\Vert x_{j}\right\Vert ^{2}=1.$
\end{theorem}

\begin{proof}
Follows from Theorem \ref{II.MGOt.4.1}. The details are omitted.
\end{proof}

\begin{remark}
\label{II.MGOr.5.2}The first inequality in (\ref{II.MGOe.5.2}) can be
written in a more convenient way as 
\begin{equation}
\left\vert \frac{\sum_{j=1}^{n}\left\langle f\left( A_{j}\right) g\left(
A_{j}\right) x_{j},x_{j}\right\rangle }{\sum_{j=1}^{n}\left\langle f\left(
A_{j}\right) x_{j},x_{j}\right\rangle \cdot \sum_{j=1}^{n}\left\langle
g\left( A_{j}\right) x_{j},x_{j}\right\rangle }-1\right\vert \leq \frac{1}{4}%
\cdot \frac{\left( \Gamma -\gamma \right) \left( \Delta -\delta \right) }{%
\sqrt{\Gamma \gamma \Delta \delta }}  \label{II.MGOe.5.4}
\end{equation}%
for each $x_{j}\in H,j\in \left\{ 1,\dots ,n\right\} $ with $%
\sum_{j=1}^{n}\left\Vert x_{j}\right\Vert ^{2}=1,$ while the second
inequality has the following equivalent form 
\begin{align}
& \left\vert \frac{\sum_{j=1}^{n}\left\langle f\left( A_{j}\right) g\left(
A_{j}\right) x_{j},x_{j}\right\rangle }{\left[ \sum_{j=1}^{n}\left\langle
f\left( A_{j}\right) x_{j},x_{j}\right\rangle \cdot
\sum_{j=1}^{n}\left\langle g\left( A_{j}\right) x_{j},x_{j}\right\rangle %
\right] ^{1/2}}\right.  \label{II.MGOe.5.5} \\
& \left. -\left[ \sum_{j=1}^{n}\left\langle f\left( A_{j}\right)
x_{j},x_{j}\right\rangle \cdot \sum_{j=1}^{n}\left\langle g\left(
A_{j}\right) x_{j},x_{j}\right\rangle \right] ^{1/2}\right\vert  \notag \\
& \leq \left( \sqrt{\Gamma }-\sqrt{\gamma }\right) \left( \sqrt{\Delta }-%
\sqrt{\delta }\right)  \notag
\end{align}%
for each $x_{j}\in H,j\in \left\{ 1,\dots ,n\right\} $ with $%
\sum_{j=1}^{n}\left\Vert x_{j}\right\Vert ^{2}=1.$

We know, from \cite{II.SSD} that if $f,g$ are synchronous (asynchronous)
functions on the interval $\left[ m,M\right] ,$ then we have the inequality 
\begin{equation}
\sum_{j=1}^{n}\left\langle f\left( A_{j}\right) g\left( A_{j}\right)
x_{j},x_{j}\right\rangle \geq \left( \leq \right) \sum_{j=1}^{n}\left\langle
f\left( A_{j}\right) x_{j},x_{j}\right\rangle \cdot
\sum_{j=1}^{n}\left\langle g\left( A_{j}\right) x_{j},x_{j}\right\rangle
\label{II.MGOe.5.6}
\end{equation}%
for each $x_{j}\in H,j\in \left\{ 1,\dots ,n\right\} $ with $%
\sum_{j=1}^{n}\left\Vert x_{j}\right\Vert ^{2}=1,$provided $f,g$ are
continuous on $\left[ m,M\right] $ and $A_{j}$ are selfadjoint operators
with $Sp\left( A_{j}\right) \subseteq \left[ m,M\right] $, $j\in \left\{
1,\dots ,n\right\} .$

Therefore, if $f,g$ are synchronous then we have from (\ref{II.MGOe.5.4})
and from (\ref{II.MGOe.5.5}) the following results: 
\begin{align}
0& \leq \frac{\sum_{j=1}^{n}\left\langle f\left( A_{j}\right) g\left(
A_{j}\right) x_{j},x_{j}\right\rangle }{\sum_{j=1}^{n}\left\langle f\left(
A_{j}\right) x_{j},x_{j}\right\rangle \cdot \sum_{j=1}^{n}\left\langle
g\left( A_{j}\right) x_{j},x_{j}\right\rangle }-1  \label{II.MGOe.5.7} \\
& \leq \frac{1}{4}\cdot \frac{\left( \Gamma -\gamma \right) \left( \Delta
-\delta \right) }{\sqrt{\Gamma \gamma \Delta \delta }}  \notag
\end{align}%
and 
\begin{align}
0& \leq \frac{\sum_{j=1}^{n}\left\langle f\left( A_{j}\right) g\left(
A_{j}\right) x_{j},x_{j}\right\rangle }{\left[ \sum_{j=1}^{n}\left\langle
f\left( A_{j}\right) x_{j},x_{j}\right\rangle \cdot
\sum_{j=1}^{n}\left\langle g\left( A_{j}\right) x_{j},x_{j}\right\rangle %
\right] ^{1/2}}  \label{II.MGOe.5.8} \\
& -\left[ \sum_{j=1}^{n}\left\langle f\left( A_{j}\right)
x_{j},x_{j}\right\rangle \cdot \sum_{j=1}^{n}\left\langle g\left(
A_{j}\right) x_{j},x_{j}\right\rangle \right] ^{1/2}  \notag \\
& \leq \left( \sqrt{\Gamma }-\sqrt{\gamma }\right) \left( \sqrt{\Delta }-%
\sqrt{\delta }\right)  \notag
\end{align}%
for each $x_{j}\in H,j\in \left\{ 1,\dots ,n\right\} $ with $%
\sum_{j=1}^{n}\left\Vert x_{j}\right\Vert ^{2}=1,$ respectively.

If $f,g$ are asynchronous then 
\begin{align}
0& \leq 1-\frac{\sum_{j=1}^{n}\left\langle f\left( A_{j}\right) g\left(
A_{j}\right) x_{j},x_{j}\right\rangle }{\sum_{j=1}^{n}\left\langle f\left(
A_{j}\right) x_{j},x_{j}\right\rangle \cdot \sum_{j=1}^{n}\left\langle
g\left( A_{j}\right) x_{j},x_{j}\right\rangle }  \label{II.MGOe.5.9} \\
& \leq \frac{1}{4}\cdot \frac{\left( \Gamma -\gamma \right) \left( \Delta
-\delta \right) }{\sqrt{\Gamma \gamma \Delta \delta }}  \notag
\end{align}%
and 
\begin{align}
0& \leq \left[ \sum_{j=1}^{n}\left\langle f\left( A_{j}\right)
x_{j},x_{j}\right\rangle \cdot \sum_{j=1}^{n}\left\langle g\left(
A_{j}\right) x_{j},x_{j}\right\rangle \right] ^{1/2}  \label{II.MGOe.5.10} \\
& -\frac{\sum_{j=1}^{n}\left\langle f\left( A_{j}\right) g\left(
A_{j}\right) x_{j},x_{j}\right\rangle }{\left[ \sum_{j=1}^{n}\left\langle
f\left( A_{j}\right) x_{j},x_{j}\right\rangle \cdot
\sum_{j=1}^{n}\left\langle g\left( A_{j}\right) x_{j},x_{j}\right\rangle %
\right] ^{1/2}}  \notag \\
& \leq \left( \sqrt{\Gamma }-\sqrt{\gamma }\right) \left( \sqrt{\Delta }-%
\sqrt{\delta }\right)  \notag
\end{align}%
for each $x_{j}\in H,j\in \left\{ 1,\dots ,n\right\} $ with $%
\sum_{j=1}^{n}\left\Vert x_{j}\right\Vert ^{2}=1,$ respectively.
\end{remark}

It is obvious that all the inequalities from Theorem \ref{II.MGOt.5.1} can
be used to obtain reverse inequalities of Gr\"{u}ss' type for various
particular instances of operator functions, see for instance \cite{II.SSDG}.
However we give here only a few provided by the inequalities (\ref%
{II.MGOe.5.7}) and (\ref{II.MGOe.5.8}) above.

\begin{example}
Let $A_{j}$ $j\in \left\{ 1,\dots ,n\right\} $ be selfadjoint operators with 
$Sp\left( A_{j}\right) \subseteq \left[ m,M\right] ,$ $j\in \left\{ 1,\dots
,n\right\} $ for some scalars $m<M.$

If $A_{j}$ are positive $\left( m\geq 0\right) $ and $p,q>0,$ then 
\begin{align}
0& \leq \frac{\sum_{j=1}^{n}\left\langle A_{j}^{p+q}x_{j},x_{j}\right\rangle 
}{\sum_{j=1}^{n}\left\langle A_{j}^{p}x_{j},x_{j}\right\rangle \cdot
\sum_{j=1}^{n}\left\langle A_{j}^{q}x_{j},x_{j}\right\rangle }-1
\label{II.MGOe.5.11} \\
& \leq \frac{1}{4}\cdot \frac{\left( M^{p}-m^{p}\right) \left(
M^{q}-m^{q}\right) }{M^{\frac{p+q}{2}}m^{\frac{p+q}{2}}}  \notag
\end{align}%
and 
\begin{align}
0& \leq \frac{\sum_{j=1}^{n}\left\langle A_{j}^{p+q}x_{j},x_{j}\right\rangle 
}{\left[ \sum_{j=1}^{n}\left\langle A_{j}^{p}x_{j},x_{j}\right\rangle \cdot
\sum_{j=1}^{n}\left\langle A_{j}^{q}x_{j},x_{j}\right\rangle \right] ^{1/2}}
\label{II.MGOe.5.12} \\
& -\left[ \sum_{j=1}^{n}\left\langle A_{j}^{p}x_{j},x_{j}\right\rangle \cdot
\sum_{j=1}^{n}\left\langle A_{j}^{q}x_{j},x_{j}\right\rangle \right] ^{1/2} 
\notag \\
& \leq \left( M^{\frac{p}{2}}-m^{\frac{p}{2}}\right) \left( M^{\frac{q}{2}%
}-m^{\frac{q}{2}}\right)  \notag
\end{align}%
for each $x_{j}\in H,j\in \left\{ 1,\dots ,n\right\} $ with $%
\sum_{j=1}^{n}\left\Vert x_{j}\right\Vert ^{2}=1.$

If\ $A$ is positive definite $\left( m>0\right) $ and $p,q<0,$ then 
\begin{align}
0& \leq \frac{\sum_{j=1}^{n}\left\langle A_{j}^{p+q}x_{j},x_{j}\right\rangle 
}{\sum_{j=1}^{n}\left\langle A_{j}^{p}x_{j},x_{j}\right\rangle \cdot
\sum_{j=1}^{n}\left\langle A_{j}^{q}x_{j},x_{j}\right\rangle }-1
\label{II.MGOe.5.13} \\
& \leq \frac{1}{4}\cdot \frac{\left( M^{-p}-m^{-p}\right) \left(
M^{-q}-m^{-q}\right) }{M^{-\frac{p+q}{2}}m^{-\frac{p+q}{2}}}  \notag
\end{align}%
and 
\begin{align}
0& \leq \left[ \sum_{j=1}^{n}\left\langle A_{j}^{p}x_{j},x_{j}\right\rangle
\cdot \sum_{j=1}^{n}\left\langle A_{j}^{q}x_{j},x_{j}\right\rangle \right]
^{1/2}  \label{II.MGOe.5.14} \\
& -\frac{\sum_{j=1}^{n}\left\langle A_{j}^{p+q}x,x\right\rangle }{\left[
\sum_{j=1}^{n}\left\langle A_{j}^{p}x_{j},x_{j}\right\rangle \cdot
\sum_{j=1}^{n}\left\langle A_{j}^{q}x_{j},x_{j}\right\rangle \right] ^{1/2}}
\notag \\
& \leq \frac{\left( M^{-\frac{p}{2}}-m^{-\frac{p}{2}}\right) \left( M^{-%
\frac{q}{2}}-m^{-\frac{q}{2}}\right) }{M^{-\frac{p+q}{2}}m^{-\frac{p+q}{2}}}
\notag
\end{align}%
for each $x_{j}\in H,j\in \left\{ 1,\dots ,n\right\} $ with $%
\sum_{j=1}^{n}\left\Vert x_{j}\right\Vert ^{2}=1.$
\end{example}

Similar inequalities may be stated for either $p>0,q<0$ or $p<0,q>0.$ The
details are omitted.

\begin{example}
Let $A$ be a positive definite operator with $Sp\left( A\right) \subseteq %
\left[ m,M\right] $ for some scalars $1<m<M.$ If $p>0$ then 
\begin{align}
0& \leq \frac{\sum_{j=1}^{n}\left\langle A_{j}^{p}\ln
A_{j}x_{j},x_{j}\right\rangle }{\sum_{j=1}^{n}\left\langle
A_{j}^{p}x_{j},x_{j}\right\rangle \cdot \sum_{j=1}^{n}\left\langle \ln
A_{j}x_{j},x_{j}\right\rangle }-1  \label{II.MGOe.5.15} \\
& \leq \frac{1}{4}\cdot \frac{\left( M^{p}-m^{p}\right) \ln \frac{M}{m}}{M^{%
\frac{p}{2}}m^{\frac{p}{2}}\sqrt{\ln M\cdot \ln m}}  \notag
\end{align}%
and 
\begin{align}
0& \leq \frac{\sum_{j=1}^{n}\left\langle A_{j}^{p}\ln
A_{j}x_{j},x_{j}\right\rangle }{\left[ \sum_{j=1}^{n}\left\langle
A_{j}^{p}x_{j},x_{j}\right\rangle \cdot \sum_{j=1}^{n}\left\langle \ln
A_{j}x_{j},x_{j}\right\rangle \right] ^{1/2}}  \label{II.MGOe.5.16} \\
& -\left[ \sum_{j=1}^{n}\left\langle A_{j}^{p}x_{j},x_{j}\right\rangle \cdot
\sum_{j=1}^{n}\left\langle \ln A_{j}x_{j},x_{j}\right\rangle \right] ^{1/2} 
\notag \\
& \leq \left( M^{\frac{p}{2}}-m^{\frac{p}{2}}\right) \left[ \sqrt{\ln M}-%
\sqrt{\ln m}\right] ,  \notag
\end{align}%
for each $x_{j}\in H,j\in \left\{ 1,\dots ,n\right\} $ with $%
\sum_{j=1}^{n}\left\Vert x_{j}\right\Vert ^{2}=1.$
\end{example}

Similar inequalities may be stated for $p<0.$ The details are omitted.

The following result for $n$ operators can be stated as well:

\begin{corollary}
\label{II.MGOc.5.1}Let $A_{j}$ be selfadjoint operators with $Sp\left(
A_{j}\right) \subseteq \left[ m,M\right] $ for $j\in \left\{ 1,\dots
,n\right\} $ and for some scalars $m<M.$ If $f$ and $g$ are continuous on $%
\left[ m,M\right] $ and $\gamma :=\min_{t\in \left[ m,M\right] }f\left(
t\right) $, $\Gamma :=\max_{t\in \left[ m,M\right] }f\left( t\right) $, $%
\delta :=\min_{t\in \left[ m,M\right] }g\left( t\right) $ and $\Delta
:=\max_{t\in \left[ m,M\right] }g\left( t\right) $ then for any $p_{j}\geq
0,j\in \left\{ 1,\dots ,n\right\} $ with $\sum_{j=1}^{n}p_{j}=1$ we have 
\begin{align}
& \left\vert \left\langle \sum_{j=1}^{n}p_{j}f\left( A_{j}\right) g\left(
A_{j}\right) x,x\right\rangle -\left\langle \sum_{j=1}^{n}p_{j}f\left(
A_{j}\right) x,x\right\rangle \cdot \left\langle \sum_{j=1}^{n}p_{j}g\left(
A_{j}\right) x,x\right\rangle \right\vert  \label{II.MGOe.5.17} \\
& \leq \frac{1}{4}\left( \Gamma -\gamma \right) \left( \Delta -\delta \right)
\notag \\
& -\left\{ 
\begin{array}{l}
\left[ \sum\limits_{j=1}^{n}p_{j}\left\langle \Gamma x-f\left( A_{j}\right)
x,f\left( A_{j}\right) x-\gamma x\right\rangle \right. \\ 
\qquad \qquad \times \left. \sum\limits_{j=1}^{n}p_{j}\left\langle \Delta
x-g\left( A_{j}\right) x,g\left( A_{j}\right) x-\delta x\right\rangle \right]
^{\frac{1}{2}}, \\ 
\\ 
\left\vert \left\langle \sum\limits_{j=1}^{n}p_{j}f\left( A_{j}\right)
x,x\right\rangle -\frac{\Gamma +\gamma }{2}\right\vert \left\vert
\left\langle \sum\limits_{j=1}^{n}p_{j}g\left( A_{j}\right) x,x\right\rangle
-\frac{\Delta +\delta }{2}\right\vert ,%
\end{array}%
\right.  \notag
\end{align}%
for each $x\in H,$ with $\left\Vert x\right\Vert ^{2}=1.$

Moreover if $\gamma $ and $\delta $ are positive, then we also have 
\begin{align}
& \left\vert \left\langle \sum_{j=1}^{n}p_{j}f\left( A_{j}\right) g\left(
A_{j}\right) x,x\right\rangle -\left\langle \sum_{j=1}^{n}p_{j}f\left(
A_{j}\right) x,x\right\rangle \cdot \left\langle \sum_{j=1}^{n}p_{j}g\left(
A_{j}\right) x,x\right\rangle \right\vert  \label{II.MGOe.5.18} \\
& \leq \left\{ 
\begin{array}{l}
\frac{1}{4}\cdot \frac{\left( \Gamma -\gamma \right) \left( \Delta -\delta
\right) }{\sqrt{\Gamma \gamma \Delta \delta }}\left\langle
\sum\limits_{j=1}^{n}p_{j}f\left( A_{j}\right) x,x\right\rangle \cdot
\left\langle \sum\limits_{j=1}^{n}p_{j}g\left( A_{j}\right) x,x\right\rangle
, \\ 
\\ 
\left( \sqrt{\Gamma }-\sqrt{\gamma }\right) \left( \sqrt{\Delta }-\sqrt{%
\delta }\right) \\ 
\times \left[ \left\langle \sum\limits_{j=1}^{n}p_{j}f\left( A_{j}\right)
x,x\right\rangle \cdot \left\langle \sum\limits_{j=1}^{n}p_{j}g\left(
A_{j}\right) x,x\right\rangle \right] ^{\frac{1}{2}}.%
\end{array}%
\right.  \notag
\end{align}%
while for $\Gamma +\gamma ,\Delta +\delta \neq 0$ we have 
\begin{align}
& \left\vert \left\langle \sum_{j=1}^{n}p_{j}f\left( A_{j}\right) g\left(
A_{j}\right) x,x\right\rangle -\left\langle \sum_{j=1}^{n}p_{j}f\left(
A_{j}\right) x,x\right\rangle \cdot \left\langle \sum_{j=1}^{n}p_{j}g\left(
A_{j}\right) x,x\right\rangle \right\vert  \label{II.MGOe.5.19} \\
& \leq \frac{1}{4}\cdot \frac{\left( \Gamma -\gamma \right) \left( \Delta
-\delta \right) }{\left[ \left\vert \Gamma +\gamma \right\vert \left\vert
\Delta +\delta \right\vert \right] ^{\frac{1}{2}}}  \notag \\
& \times \left\lceil \left( \left( \sum_{j=1}^{n}p_{j}\left\Vert f\left(
A_{j}\right) x\right\Vert ^{2}\right) ^{1/2}+\left\vert \left\langle
\sum_{j=1}^{n}p_{j}f\left( A_{j}\right) x,x\right\rangle \right\vert \right)
\right.  \notag \\
& \left. \times \left( \left( \sum_{j=1}^{n}p_{j}\left\Vert g\left(
A_{j}\right) x\right\Vert ^{2}\right) ^{1/2}+\left\vert \left\langle
\sum_{j=1}^{n}p_{j}g\left( A_{j}\right) x,x\right\rangle \right\vert \right) 
\right] ^{1/2}  \notag
\end{align}%
for each $x\in H,$ with $\left\Vert x\right\Vert ^{2}=1.$
\end{corollary}

\begin{proof}
Follows from Theorem \ref{II.MGOt.5.1} on choosing $x_{j}=\sqrt{p_{j}}\cdot
x,$ $j\in \left\{ 1,\dots ,n\right\} ,$ where $p_{j}\geq 0,j\in \left\{
1,\dots ,n\right\} ,$ $\sum_{j=1}^{n}p_{j}=1$ and $x\in H,$ with $\left\Vert
x\right\Vert =1.$ The details are omitted.
\end{proof}

\begin{remark}
The first inequality in (\ref{II.MGOe.5.18}) can be written in a more
convenient way as 
\begin{equation}
\left\vert \frac{\left\langle \sum_{j=1}^{n}p_{j}f\left( A_{j}\right)
g\left( A_{j}\right) x,x\right\rangle }{\left\langle
\sum_{j=1}^{n}p_{j}f\left( A_{j}\right) x,x\right\rangle \cdot \left\langle
\sum_{j=1}^{n}p_{j}g\left( A_{j}\right) x,x\right\rangle }-1\right\vert \leq 
\frac{1}{4}\cdot \frac{\left( \Gamma -\gamma \right) \left( \Delta -\delta
\right) }{\sqrt{\Gamma \gamma \Delta \delta }}  \label{II.MGOe.5.20}
\end{equation}%
for each $x\in H,$ with $\left\Vert x\right\Vert ^{2}=1,$ while the second
inequality has the following equivalent form 
\begin{align}
& \left\vert \frac{\left\langle \sum_{j=1}^{n}p_{j}f\left( A_{j}\right)
g\left( A_{j}\right) x,x\right\rangle }{\left[ \left\langle
\sum_{j=1}^{n}p_{j}f\left( A_{j}\right) x,x\right\rangle \cdot \left\langle
\sum_{j=1}^{n}p_{j}g\left( A_{j}\right) x,x\right\rangle \right] ^{1/2}}%
\right.  \label{II.MGOe.5.21} \\
& \left. -\left[ \left\langle \sum_{j=1}^{n}p_{j}f\left( A_{j}\right)
x,x\right\rangle \cdot \left\langle \sum_{j=1}^{n}p_{j}g\left( A_{j}\right)
x,x\right\rangle \right] ^{1/2}\right\vert  \notag \\
& \leq \left( \sqrt{\Gamma }-\sqrt{\gamma }\right) \left( \sqrt{\Delta }-%
\sqrt{\delta }\right)  \notag
\end{align}%
for each $x\in H,$ with $\left\Vert x\right\Vert ^{2}=1.$

We know, from \cite{II.SSD} that if $f,g$ are synchronous (asynchronous)
functions on the interval $\left[ m,M\right] ,$ then we have the inequality 
\begin{equation}
\left\langle \sum_{j=1}^{n}p_{j}f\left( A_{j}\right) g\left( A_{j}\right)
x,x\right\rangle \geq \left( \leq \right) \left\langle
\sum_{j=1}^{n}p_{j}f\left( A_{j}\right) x,x\right\rangle \cdot \left\langle
\sum_{j=1}^{n}p_{j}g\left( A_{j}\right) x,x\right\rangle
\label{II.MGOe.5.22}
\end{equation}%
for each $x\in H,$ with $\left\Vert x\right\Vert ^{2}=1,$ provided $f,g$ are
continuous on $\left[ m,M\right] $ and $A_{j}$ are selfadjoint operators
with $Sp\left( A_{j}\right) \subseteq \left[ m,M\right] $, $j\in \left\{
1,\dots ,n\right\} .$

Therefore, if $f,g$ are synchronous then we have from (\ref{II.MGOe.5.20})
and from (\ref{II.MGOe.5.21}) the following results: 
\begin{align}
0& \leq \frac{\left\langle \sum_{j=1}^{n}p_{j}f\left( A_{j}\right) g\left(
A_{j}\right) x,x\right\rangle }{\left\langle \sum_{j=1}^{n}p_{j}f\left(
A_{j}\right) x,x\right\rangle \cdot \left\langle \sum_{j=1}^{n}p_{j}g\left(
A_{j}\right) x,x\right\rangle }-1  \label{II.MGOe.5.23} \\
& \leq \frac{1}{4}\cdot \frac{\left( \Gamma -\gamma \right) \left( \Delta
-\delta \right) }{\sqrt{\Gamma \gamma \Delta \delta }}  \notag
\end{align}%
and 
\begin{align}
0& \leq \frac{\left\langle \sum_{j=1}^{n}p_{j}f\left( A_{j}\right) g\left(
A_{j}\right) x,x\right\rangle }{\left[ \left\langle
\sum_{j=1}^{n}p_{j}f\left( A_{j}\right) x,x\right\rangle \cdot \left\langle
\sum_{j=1}^{n}p_{j}g\left( A_{j}\right) x,x\right\rangle \right] ^{1/2}}
\label{II.MGOe.5.24} \\
& -\left[ \left\langle \sum_{j=1}^{n}p_{j}f\left( A_{j}\right)
x,x\right\rangle \cdot \left\langle \sum_{j=1}^{n}p_{j}g\left( A_{j}\right)
x,x\right\rangle \right] ^{1/2}  \notag \\
& \leq \left( \sqrt{\Gamma }-\sqrt{\gamma }\right) \left( \sqrt{\Delta }-%
\sqrt{\delta }\right)  \notag
\end{align}%
for each $x\in H,$ with $\left\Vert x\right\Vert =1,$ respectively.

If $f,g$ are asynchronous then 
\begin{align}
0& \leq 1-\frac{\left\langle \sum_{j=1}^{n}p_{j}f\left( A_{j}\right) g\left(
A_{j}\right) x,x\right\rangle }{\left\langle \sum_{j=1}^{n}p_{j}f\left(
A_{j}\right) x,x\right\rangle \cdot \left\langle \sum_{j=1}^{n}p_{j}g\left(
A_{j}\right) x,x\right\rangle }  \label{II.MGOe.5.25} \\
& \leq \frac{1}{4}\cdot \frac{\left( \Gamma -\gamma \right) \left( \Delta
-\delta \right) }{\sqrt{\Gamma \gamma \Delta \delta }}  \notag
\end{align}%
and 
\begin{align}
0& \leq \left[ \left\langle \sum_{j=1}^{n}p_{j}f\left( A_{j}\right)
x,x\right\rangle \cdot \left\langle \sum_{j=1}^{n}p_{j}g\left( A_{j}\right)
x,x\right\rangle \right] ^{1/2}  \label{II.MGOe.5.26} \\
& -\frac{\left\langle \sum_{j=1}^{n}p_{j}f\left( A_{j}\right) g\left(
A_{j}\right) x,x\right\rangle }{\left[ \left\langle
\sum_{j=1}^{n}p_{j}f\left( A_{j}\right) x,x\right\rangle \cdot \left\langle
\sum_{j=1}^{n}p_{j}g\left( A_{j}\right) x,x\right\rangle \right] ^{1/2}} 
\notag \\
& \leq \left( \sqrt{\Gamma }-\sqrt{\gamma }\right) \left( \sqrt{\Delta }-%
\sqrt{\delta }\right)  \notag
\end{align}%
for each $x\in H,$ with $\left\Vert x\right\Vert =1,$ respectively.
\end{remark}

The above inequalities (\ref{II.MGOe.5.23}) - (\ref{II.MGOe.5.26}) can be
used to state various particular inequalities as in the previous examples,
however the details are left to the interested reader.

\section{More Inequalities for the \v{C}eby\v{s}ev Functional}

\subsection{A Refinement and Some Related Results}

The following result can be stated:

\begin{theorem}[Dragomir, 2008, \protect\cite{II.SSDG2}]
\label{II.mcft.2.1}Let $A$ be a selfadjoint operator with $Sp\left( A\right)
\subseteq \left[ m,M\right] $ for some real numbers $m<M.$ If $f,g:\left[ m,M%
\right] \longrightarrow \mathbb{R}$ are continuous with $\delta :=\min_{t\in %
\left[ m,M\right] }g\left( t\right) $ and $\Delta :=\max_{t\in \left[ m,M%
\right] }g\left( t\right) ,$ then 
\begin{align}
\left\vert C\left( f,g;A;x\right) \right\vert & \leq \frac{1}{2}\left(
\Delta -\delta \right) \left\langle \left\vert f\left( A\right)
-\left\langle f\left( A\right) x,x\right\rangle \cdot 1_{H}\right\vert
x,x\right\rangle  \label{II.mcfe.2.1.1} \\
& \leq \frac{1}{2}\left( \Delta -\delta \right) C^{1/2}\left( f,f;A;x\right)
,  \notag
\end{align}%
for any $x\in H$ with $\left\Vert x\right\Vert =1.$
\end{theorem}

\begin{proof}
Since $\delta :=\min_{t\in \left[ m,M\right] }g\left( t\right) $ and $\Delta
:=\max_{t\in \left[ m,M\right] }g\left( t\right) ,$ we have 
\begin{equation}
\left\vert g\left( t\right) -\frac{\Delta +\delta }{2}\right\vert \leq \frac{%
1}{2}\left( \Delta -\delta \right) ,  \label{II.mcfe.2.7}
\end{equation}%
for any $t\in \left[ m,M\right] $ and for any $x\in H$ with $\left\Vert
x\right\Vert =1.$

If we multiply the inequality (\ref{II.mcfe.2.7}) with $\left\vert f\left(
t\right) -\left\langle f\left( A\right) x,x\right\rangle \right\vert $ we
get 
\begin{align}
& \left\vert f\left( t\right) g\left( t\right) -\left\langle f\left(
A\right) x,x\right\rangle g\left( t\right) -\frac{\Delta +\delta }{2}f\left(
t\right) +\frac{\Delta +\delta }{2}\left\langle f\left( A\right)
x,x\right\rangle \right\vert  \label{II.mcfe.2.8} \\
& \leq \frac{1}{2}\left( \Delta -\delta \right) \left\vert f\left( t\right)
-\left\langle f\left( A\right) x,x\right\rangle \right\vert ,  \notag
\end{align}%
for any $t\in \left[ m,M\right] $ and for any $x\in H$ with $\left\Vert
x\right\Vert =1.$

Now, if we apply the property (\ref{P}) for the inequality (\ref{II.mcfe.2.8}%
) and a selfadjoint operator $B$ with $Sp\left( B\right) \subset \left[ m,M%
\right] ,$ then we get the following inequality of interest in itself: 
\begin{align}
& \left\vert \left\langle f\left( B\right) g\left( B\right) y,y\right\rangle
-\left\langle f\left( A\right) x,x\right\rangle \left\langle g\left(
B\right) y,y\right\rangle \right.  \label{II.mcfe.2.9} \\
& \left. -\frac{\Delta +\delta }{2}\left\langle f\left( B\right)
y,y\right\rangle +\frac{\Delta +\delta }{2}\left\langle f\left( A\right)
x,x\right\rangle \right\vert  \notag \\
& \leq \frac{1}{2}\left( \Delta -\delta \right) \left\langle \left\vert
f\left( B\right) -\left\langle f\left( A\right) x,x\right\rangle \cdot
1_{H}\right\vert y,y\right\rangle ,  \notag
\end{align}%
for any $x,y\in H$ with $\left\Vert x\right\Vert =\left\Vert y\right\Vert
=1. $

If we choose in (\ref{II.mcfe.2.9}) $y=x$ and $B=A,$ then we deduce the
first inequality in (\ref{II.mcfe.2.1.1}).

Now, by the Schwarz inequality in $H$ we have 
\begin{align*}
\left\langle \left\vert f\left( A\right) -\left\langle f\left( A\right)
x,x\right\rangle \cdot 1_{H}\right\vert x,x\right\rangle & \leq \left\Vert
\left\vert f\left( A\right) -\left\langle f\left( A\right) x,x\right\rangle
\cdot 1_{H}\right\vert x\right\Vert \\
& =\left\Vert f\left( A\right) x-\left\langle f\left( A\right)
x,x\right\rangle \cdot x\right\Vert \\
& =\left[ \left\Vert f\left( A\right) x\right\Vert ^{2}-\left\langle f\left(
A\right) x,x\right\rangle ^{2}\right] ^{1/2} \\
& =C^{1/2}\left( f,f;A;x\right) ,
\end{align*}%
for any $x\in H$ with $\left\Vert x\right\Vert =1,$ and the second part of (%
\ref{II.mcfe.2.1.1}) is also proved.
\end{proof}

Let $U$ be a selfadjoint operator on the Hilbert space $\left(
H,\left\langle .,.\right\rangle \right) $ with the spectrum $Sp\left(
U\right) $ included in the interval $\left[ m,M\right] $ for some real
numbers $m<M$ and let $\left\{ E_{\lambda }\right\} _{\lambda \in \mathbb{R}%
} $ be its \textit{spectral family}. Then for any continuous function $f:%
\left[ m,M\right] \rightarrow \mathbb{R}$, it is well known that we have the
following representation in terms of the Riemann-Stieltjes integral: 
\begin{equation}
\left\langle f\left( U\right) x,x\right\rangle =\int_{m-0}^{M}f\left(
\lambda \right) d\left( \left\langle E_{\lambda }x,x\right\rangle \right) ,
\label{II.mcfe.2.10}
\end{equation}%
for any $x\in H$ with $\left\Vert x\right\Vert =1.$ The function $%
g_{x}\left( \lambda \right) :=\left\langle E_{\lambda }x,x\right\rangle $ is 
\textit{monotonic nondecreasing} on the interval $\left[ m,M\right] $ and 
\begin{equation}
g_{x}\left( m-0\right) =0\text{ \quad and\quad\ }g_{x}\left( M\right) =1
\label{II.mcfe.2.10.a}
\end{equation}%
for any $x\in H$ with $\left\Vert x\right\Vert =1.$

The following result is of interest:

\begin{theorem}[Dragomir, 2008, \protect\cite{II.SSDG2}]
\label{II.mcft.2.2}Let $A$ and $B$ be selfadjoint operators with $Sp\left(
A\right) ,Sp\left( B\right) \subseteq \left[ m,M\right] $ for some real
numbers $m<M.$ If $f:\left[ m,M\right] \longrightarrow \mathbb{R}$ is of $%
r-L-$H\"{o}lder type, i.e., for a given $r\in (0,1]$ and $L>0$ we have 
\begin{equation*}
\left\vert f\left( s\right) -f\left( t\right) \right\vert \leq L\left\vert
s-t\right\vert ^{r}\text{ \quad for any }s,t\in \left[ m,M\right] ,
\end{equation*}%
then we have the Ostrowski type inequality for selfadjoint operators: 
\begin{equation}
\left\vert f\left( s\right) -\left\langle f\left( A\right) x,x\right\rangle
\right\vert \leq L\left[ \frac{1}{2}\left( M-m\right) +\left\vert s-\frac{m+M%
}{2}\right\vert \right] ^{r},  \label{II.mcfe.2.11.0}
\end{equation}%
for any $s\in \left[ m,M\right] $ and any $x\in H$ with $\left\Vert
x\right\Vert =1$.

Moreover, we have 
\begin{align}
\left\vert \left\langle f\left( B\right) y,y\right\rangle -\left\langle
f\left( A\right) x,x\right\rangle \right\vert & \leq \left\langle \left\vert
f\left( B\right) -\left\langle f\left( A\right) x,x\right\rangle \cdot
1_{H}\right\vert y,y\right\rangle  \label{II.mcfe.2.11} \\
& \leq L\left[ \frac{1}{2}\left( M-m\right) +\left\langle \left\vert B-\frac{%
m+M}{2}\cdot 1_{H}\right\vert y,y\right\rangle \right] ^{r},  \notag
\end{align}%
for any $x,y\in H$ with $\left\Vert x\right\Vert =\left\Vert y\right\Vert
=1. $
\end{theorem}

\begin{proof}
We use the following Ostrowski type inequality for the Riemann-Stieltjes
integral obtained by the author in \cite{II.mcfSSD1}: 
\begin{align}
& \left\vert f\left( s\right) \left[ u\left( b\right) -u\left( a\right) %
\right] -\int_{a}^{b}f\left( t\right) du\left( t\right) \right\vert
\label{II.mcfe.2.12} \\
& \leq L\left[ \frac{1}{2}\left( b-a\right) +\left\vert s-\frac{a+b}{2}%
\right\vert \right] ^{r}\dbigvee\limits_{a}^{b}\left( u\right)  \notag
\end{align}%
for any $s\in \left[ a,b\right] ,$ provided that $f$ is of $r-L-$H\"{o}lder
type on $\left[ a,b\right] ,$ $u$ is of bounded variation on $\left[ a,b%
\right] $ and $\dbigvee_{a}^{b}\left( u\right) $ denotes the total variation
of $u$ on $\left[ a,b\right] .$

Now, applying this inequality for $u\left( \lambda \right) =g_{x}\left(
\lambda \right) :=\left\langle E_{\lambda }x,x\right\rangle $ where $x\in H$
with $\left\Vert x\right\Vert =1$ we get 
\begin{align}
& \left\vert f\left( s\right) -\int_{m-0}^{M}f\left( \lambda \right) d\left(
\left\langle E_{\lambda }x,x\right\rangle \right) \right\vert
\label{II.mcfe.2.13} \\
& \leq L\left[ \frac{1}{2}\left( M-m\right) +\left\vert s-\frac{m+M}{2}%
\right\vert \right] ^{r}\dbigvee\limits_{m-0}^{M}\left( g_{x}\right)  \notag
\end{align}%
which, by (\ref{II.mcfe.2.10}) and (\ref{II.mcfe.2.10.a}) is equivalent with
(\ref{II.mcfe.2.11.0}).

By applying the property (\ref{P}) for the inequality (\ref{II.mcfe.2.11.0})
and the operator $B$ we have 
\begin{align*}
\left\langle \left\vert f\left( B\right) -\left\langle f\left( A\right)
x,x\right\rangle \cdot 1_{H}\right\vert y,y\right\rangle & \leq
L\left\langle \left[ \frac{1}{2}\left( M-m\right) +\left\vert B-\frac{m+M}{2}%
\cdot 1_{H}\right\vert \right] ^{r}y,y\right\rangle \\
& \leq L\left\langle \left[ \frac{1}{2}\left( M-m\right) +\left\vert B-\frac{%
m+M}{2}\right\vert \cdot 1_{H}\right] y,y\right\rangle ^{r} \\
& =L\left[ \frac{1}{2}\left( M-m\right) +\left\langle \left\vert B-\frac{m+M%
}{2}\cdot 1_{H}\right\vert y,y\right\rangle \right] ^{r}
\end{align*}%
for any $x,y\in H$ with $\left\Vert x\right\Vert =\left\Vert y\right\Vert
=1, $ which proves the second inequality in (\ref{II.mcfe.2.11}).

Further, by the Jensen inequality for convex functions of selfadjoint
operators (see for instance \cite[p. 5]{II.FMPS}) applied for the modulus,
we can state that 
\begin{equation}
\left\vert \left\langle h\left( A\right) x,x\right\rangle \right\vert \leq
\left\langle \left\vert h\left( A\right) \right\vert x,x\right\rangle 
\tag{M}  \label{II.mcfM}
\end{equation}%
for any $x\in H$ with $\left\Vert x\right\Vert =1,$ where $h$ is a
continuous function on $\left[ m,M\right] .$

Now, if we apply the inequality (\ref{II.mcfM}), then we have 
\begin{equation*}
\left\vert \left\langle \left[ f\left( B\right) -\left\langle f\left(
A\right) x,x\right\rangle \cdot 1_{H}\right] y,y\right\rangle \right\vert
\leq \left\langle \left\vert f\left( B\right) -\left\langle f\left( A\right)
x,x\right\rangle \cdot 1_{H}\right\vert y,y\right\rangle
\end{equation*}%
which shows the first part of (\ref{II.mcfe.2.11}), and the proof is
complete.
\end{proof}

\begin{remark}
\label{II.mcfr.2.1}With the above assumptions for $f,A$ and $B$ we have the
following particular inequalities of interest: 
\begin{equation}
\left\vert f\left( \frac{m+M}{2}\right) -\left\langle f\left( A\right)
x,x\right\rangle \right\vert \leq \frac{1}{2^{r}}L\left( M-m\right) ^{r}
\label{II.mcfe.2.14}
\end{equation}%
and 
\begin{equation}
\left\vert f\left( \left\langle Ax,x\right\rangle \right) -\left\langle
f\left( A\right) x,x\right\rangle \right\vert \leq L\left[ \frac{1}{2}\left(
M-m\right) +\left\vert \left\langle Ax,x\right\rangle -\frac{m+M}{2}%
\right\vert \right] ^{r},  \label{II.mcfe.2.15}
\end{equation}%
for any $x\in H$ with $\left\Vert x\right\Vert =1$.

We also have the inequalities: 
\begin{align}
& \left\vert \left\langle f\left( A\right) y,y\right\rangle -\left\langle
f\left( A\right) x,x\right\rangle \right\vert  \label{II.mcfe.2.16} \\
& \leq \left\langle \left\vert f\left( A\right) -\left\langle f\left(
A\right) x,x\right\rangle \cdot 1_{H}\right\vert y,y\right\rangle  \notag \\
& \leq L\left[ \frac{1}{2}\left( M-m\right) +\left\langle \left\vert A-\frac{%
m+M}{2}\cdot 1_{H}\right\vert y,y\right\rangle \right] ^{r},  \notag
\end{align}%
for any $x,y\in H$ with $\left\Vert x\right\Vert =\left\Vert y\right\Vert
=1, $%
\begin{align}
\left\vert \left\langle \left[ f\left( B\right) -f\left( A\right) \right]
x,x\right\rangle \right\vert & \leq \left\langle \left\vert f\left( B\right)
-\left\langle f\left( A\right) x,x\right\rangle \cdot 1_{H}\right\vert
x,x\right\rangle  \label{II.mcfe.2.17} \\
& \leq L\left[ \frac{1}{2}\left( M-m\right) +\left\langle \left\vert B-\frac{%
m+M}{2}\cdot 1_{H}\right\vert x,x\right\rangle \right] ^{r}  \notag
\end{align}%
and, more particularly, 
\begin{align}
& \left\langle \left\vert f\left( A\right) -\left\langle f\left( A\right)
x,x\right\rangle \cdot 1_{H}\right\vert x,x\right\rangle
\label{II.mcfe.2.18} \\
& \leq L\left[ \frac{1}{2}\left( M-m\right) +\left\langle \left\vert A-\frac{%
m+M}{2}\cdot 1_{H}\right\vert x,x\right\rangle \right] ^{r},  \notag
\end{align}%
for any $x\in H$ with $\left\Vert x\right\Vert =1.$

We also have the norm inequality 
\begin{equation}
\left\Vert f\left( B\right) -f\left( A\right) \right\Vert \leq L\left[ \frac{%
1}{2}\left( M-m\right) +\left\Vert B-\frac{m+M}{2}\cdot 1_{H}\right\Vert %
\right] ^{r}.  \label{II.mcfe.2.19}
\end{equation}
\end{remark}

The following corollary of the above Theorem \ref{II.mcft.2.2} can be useful
for applications:

\begin{corollary}[Dragomir, 2008, \protect\cite{II.SSDG2}]
\label{II.mcfc.2.1}Let $A$ and $B$ be selfadjoint operators with $Sp\left(
A\right) ,Sp\left( B\right) \subseteq \left[ m,M\right] $ for some real
numbers $m<M.$ If $f:\left[ m,M\right] \longrightarrow \mathbb{R}$ is
absolutely continuous then we have the Ostrowski type inequality for
selfadjoint operators: 
\begin{multline}
\left\vert f\left( s\right) -\left\langle f\left( A\right) x,x\right\rangle
\right\vert  \label{II.mcfe.2.20} \\
\leq \left\{ 
\begin{array}{ll}
\left[ \frac{1}{2}\left( M-m\right) +\left\vert s-\frac{m+M}{2}\right\vert %
\right] \left\Vert f^{\prime }\right\Vert _{\infty ,\left[ m,M\right] } & 
\text{if \ }f^{\prime }\in L_{\infty }\left[ m,M\right] ; \\ 
&  \\ 
\left[ \frac{1}{2}\left( M-m\right) +\left\vert s-\frac{m+M}{2}\right\vert %
\right] ^{1/q}\left\Vert f^{\prime }\right\Vert _{p,\left[ m,M\right] } & 
\begin{array}{l}
\text{if \ }f^{\prime }\in L_{p}\left[ m,M\right] , \\ 
p,q>1,\frac{1}{p}+\frac{1}{q}=1,%
\end{array}%
\end{array}%
\right.
\end{multline}%
for any $s\in \left[ m,M\right] $ and any $x\in H$ with $\left\Vert
x\right\Vert =1$, where $\left\Vert \cdot \right\Vert _{p,\left[ m,M\right]
} $ are the Lebesgue norms, i.e., 
\begin{equation*}
\left\Vert h\right\Vert _{\infty ,\left[ m,M\right] }:=ess\sup_{t\in \left[
m,M\right] }\left\Vert h\left( t\right) \right\Vert
\end{equation*}%
and 
\begin{equation*}
\left\Vert h\right\Vert _{p,\left[ m,M\right] }:=\left(
\int_{m}^{M}\left\vert h\left( t\right) \right\vert ^{p}\right) ^{1/p},\quad
p\geq 1.
\end{equation*}

Moreover, we have 
\begin{align}
& \left\vert \left\langle f\left( B\right) y,y\right\rangle -\left\langle
f\left( A\right) x,x\right\rangle \right\vert  \label{II.mcfe.2.21} \\
& \leq \left\langle \left\vert f\left( B\right) -\left\langle f\left(
A\right) x,x\right\rangle \cdot 1_{H}\right\vert y,y\right\rangle  \notag \\
& \leq \hspace{-0.07in}\left\{ \hspace{-0.07in}%
\begin{array}{ll}
\left[ \frac{M-m}{2}+\left\langle \left\vert B-\frac{m+M}{2}\cdot
1_{H}\right\vert y,y\right\rangle \right] \left\Vert f^{\prime }\right\Vert
_{\infty ,\left[ m,M\right] } & \hspace{-0.07in}\text{if \ }f^{\prime }\in
L_{\infty }\left[ m,M\right] ; \\ 
& \hspace{-0.07in} \\ 
\left[ \frac{M-m}{2}+\left\langle \left\vert B-\frac{m+M}{2}\cdot
1_{H}\right\vert y,y\right\rangle \right] ^{\frac{1}{q}}\left\Vert f^{\prime
}\right\Vert _{p,\left[ m,M\right] } & \hspace{-0.07in}\hspace{-0.07in}%
\begin{array}{l}
\text{if \ }f^{\prime }\in L_{p}\left[ m,M\right] , \\ 
p,q>1,\frac{1}{p}+\frac{1}{q}=1,%
\end{array}%
\end{array}%
\right.  \notag
\end{align}%
for any $x,y\in H$ with $\left\Vert x\right\Vert =\left\Vert y\right\Vert
=1. $
\end{corollary}

Now, on utilising Theorem \ref{II.mcft.2.1} we can provide the following
upper bound for the \v{C}eby\v{s}ev functional that may be\ more useful in
applications:

\begin{corollary}[Dragomir, 2008, \protect\cite{II.SSDG2}]
\label{II.mcfc.2.2}Let $A$ be a selfadjoint operator with $Sp\left( A\right)
\subseteq \left[ m,M\right] $ for some real numbers $m<M.$ If $g:\left[ m,M%
\right] \longrightarrow \mathbb{R}$ is continuous with $\delta :=\min_{t\in %
\left[ m,M\right] }g\left( t\right) $ and $\Delta :=\max_{t\in \left[ m,M%
\right] }g\left( t\right) ,$ then for any $f:\left[ m,M\right]
\longrightarrow \mathbb{R}$ of $r-L-$H\"{o}lder type we have the inequality: 
\begin{equation}
\left\vert C\left( f,g;A;x\right) \right\vert \leq \frac{1}{2}\left( \Delta
-\delta \right) L\left[ \frac{1}{2}\left( M-m\right) +\left\langle
\left\vert A-\frac{m+M}{2}\cdot 1_{H}\right\vert x,x\right\rangle \right]
^{r},  \label{II.mcfe.2.22}
\end{equation}%
for any $x\in H$ with $\left\Vert x\right\Vert =1.$
\end{corollary}

\begin{remark}
\label{II.mcfr.2.2}With the assumptions from Corollary \ref{II.mcfc.2.2} for 
$g$ and $A$ and if $f$ is absolutely continuos on $\left[ m,M\right] ,$ then
we have the inequalities: 
\begin{multline}
\left\vert C\left( f,g;A;x\right) \right\vert \leq \frac{1}{2}\left( \Delta
-\delta \right)  \label{II.mcfe.2.23} \\
\times \left\{ \hspace{-0.07in}%
\begin{array}{ll}
\left[ \frac{1}{2}\left( M-m\right) +\left\langle \left\vert A-\frac{m+M}{2}%
\cdot 1_{H}\right\vert x,x\right\rangle \right] \left\Vert f^{\prime
}\right\Vert _{\infty ,\left[ m,M\right] } & \hspace{-0.07in}\text{if \ }%
f^{\prime }\in L_{\infty }\left[ m,M\right] ; \\ 
& \hspace{-0.07in} \\ 
\left[ \frac{1}{2}\left( M-m\right) +\left\langle \left\vert A-\frac{m+M}{2}%
\cdot 1_{H}\right\vert x,x\right\rangle \right] ^{1/q}\left\Vert f^{\prime
}\right\Vert _{p,\left[ m,M\right] } & \hspace{-0.07in}%
\begin{array}{l}
\text{if \ }f^{\prime }\in L_{\infty }\left[ m,M\right] , \\ 
p,q>1,\frac{1}{p}+\frac{1}{q}=1%
\end{array}%
\end{array}%
\right.
\end{multline}%
for any $x\in H$ with $\left\Vert x\right\Vert =1.$
\end{remark}

\subsection{Some Inequalities for Sequences of Operators}

Consider the sequence of selfadjoint operators $\mathbf{A}=\left(
A_{1},\dots ,A_{n}\right) $ with $Sp\left( A_{j}\right) \subseteq \left[ m,M%
\right] $ for $j\in \left\{ 1,\dots ,n\right\} $ and for some scalars $m<M.$
If $\mathbf{x}=\left( x_{1},\dots ,x_{n}\right) \in H^{n}$ are such that $%
\sum_{j=1}^{n}\left\Vert x_{j}\right\Vert ^{2}=1,$ then we can consider the
following \v{C}eby\v{s}ev type functional 
\begin{equation*}
C\left( f,g;\mathbf{A,x}\right) :=\sum_{j=1}^{n}\left\langle f\left(
A_{j}\right) g\left( A_{j}\right) x_{j},x_{j}\right\rangle
-\sum_{j=1}^{n}\left\langle f\left( A_{j}\right) x_{j},x_{j}\right\rangle
\cdot \sum_{j=1}^{n}\left\langle g\left( A_{j}\right)
x_{j},x_{j}\right\rangle .
\end{equation*}%
As a particular case of the above functional and for a probability sequence $%
\mathbf{p=}\left( p_{1},\dots ,p_{n}\right) ,$ i.e., $p_{j}\geq 0$ for $j\in
\left\{ 1,\dots ,n\right\} $ and $\sum_{j=1}^{n}p_{j}=1,$ we can also
consider the functional 
\begin{align*}
C\left( f,g;\mathbf{A,p,}x\right) & :=\left\langle
\sum_{j=1}^{n}p_{j}f\left( A_{j}\right) g\left( A_{j}\right) x,x\right\rangle
\\
& -\left\langle \sum_{j=1}^{n}p_{j}f\left( A_{j}\right) x,x\right\rangle
\cdot \left\langle \sum_{j=1}^{n}p_{j}g\left( A_{j}\right) x,x\right\rangle
\end{align*}%
where $x\in H,$ $\left\Vert x\right\Vert =1.$

We know, from \cite{II.SSD} that for the sequence of selfadjoint operators $%
\mathbf{A}=\left( A_{1},\dots ,A_{n}\right) $ with $Sp\left( A_{j}\right)
\subseteq \left[ m,M\right] $ for $j\in \left\{ 1,\dots ,n\right\} $ and for
the synchronous (asynchronous) functions $f,g:\left[ m,M\right]
\longrightarrow \mathbb{R}$ we have the inequality 
\begin{equation}
C\left( f,g;\mathbf{A,x}\right) \geq \left( \leq \right) 0\text{ }
\label{II.mcfe.4.1}
\end{equation}%
for any $\mathbf{x}=\left( x_{1},\dots ,x_{n}\right) \in H^{n}$ with $%
\sum_{j=1}^{n}\left\Vert x_{j}\right\Vert ^{2}=1.$ Also, for any probability
distribution $\mathbf{p=}\left( p_{1},\dots ,p_{n}\right) $ and any $x\in H,$
$\left\Vert x\right\Vert =1$ we have 
\begin{equation}
C\left( f,g;\mathbf{A,p,}x\right) \geq \left( \leq \right) 0.
\label{II.mcfe.4.2}
\end{equation}

On the other hand, the following Gr\"{u}ss' type inequality is valid as well 
\cite{II.SSDG}: 
\begin{equation}
\left\vert C\left( f,g;\mathbf{A,x}\right) \right\vert \leq \frac{1}{2}\cdot
\left( \Gamma -\gamma \right) \left[ C\left( g,g;\mathbf{A,x}\right) \right]
^{1/2}\left( \leq \frac{1}{4}\left( \Gamma -\gamma \right) \left( \Delta
-\delta \right) \right)  \label{II.mcfe.4.3}
\end{equation}%
for any $\mathbf{x}=\left( x_{1},\dots ,x_{n}\right) \in H^{n}$ with $%
\sum_{j=1}^{n}\left\Vert x_{j}\right\Vert ^{2}=1,$ where $f$ and $g$ are
continuous on $\left[ m,M\right] $ and $\gamma :=\min_{t\in \left[ m,M\right]
}f\left( t\right) $, $\Gamma :=\max_{t\in \left[ m,M\right] }f\left(
t\right) $, $\delta :=\min_{t\in \left[ m,M\right] }g\left( t\right) $ and $%
\Delta :=\max_{t\in \left[ m,M\right] }g\left( t\right) .$

Similarly, for any probability distribution $\mathbf{p=}\left( p_{1},\dots
,p_{n}\right) $ and any $x\in H,$ $\left\Vert x\right\Vert =1$ we also have
the inequality: 
\begin{equation}
\left\vert C\left( f,g;\mathbf{A,p,}x\right) \right\vert \leq \frac{1}{2}%
\cdot \left( \Gamma -\gamma \right) \left[ C\left( g,g;\mathbf{A,p,}x\right) %
\right] ^{1/2}\left( \leq \frac{1}{4}\left( \Gamma -\gamma \right) \left(
\Delta -\delta \right) \right) .  \label{II.mcfe.4.4}
\end{equation}

We can state now the following new result:

\begin{theorem}[Dragomir, 2008, \protect\cite{II.SSDG2}]
\label{II.mcft.4.1}Consider the sequence of selfadjoint operators $\mathbf{A}%
=\left( A_{1},\dots ,A_{n}\right) $ with $Sp\left( A_{j}\right) \subseteq %
\left[ m,M\right] $ for $j\in \left\{ 1,\dots ,n\right\} $ and for some
scalars $m<M.$ If $f,g:\left[ m,M\right] \longrightarrow \mathbb{R}$ are
continuous with $\delta :=\min_{t\in \left[ m,M\right] }g\left( t\right) $
and $\Delta :=\max_{t\in \left[ m,M\right] }g\left( t\right) ,$ then 
\begin{align}
\left\vert C\left( f,g;\mathbf{A};\mathbf{x}\right) \right\vert & \leq \frac{%
1}{2}\left( \Delta -\delta \right) \sum_{j=1}^{n}\left\langle \left\vert
f\left( A_{j}\right) -\sum_{k=1}^{n}\left\langle f\left( A_{k}\right)
x_{k},x_{k}\right\rangle \cdot 1_{H}\right\vert x_{j},x_{j}\right\rangle
\label{II.mcfe.4.5} \\
& \leq \frac{1}{2}\left( \Delta -\delta \right) C^{1/2}\left( f,f;\mathbf{A};%
\mathbf{x}\right) ,  \notag
\end{align}%
for any $\mathbf{x}=\left( x_{1},\dots ,x_{n}\right) \in H^{n}$ such that $%
\sum_{j=1}^{n}\left\Vert x_{j}\right\Vert ^{2}=1.$
\end{theorem}

\begin{proof}
Follows from Theorem \ref{II.mcft.2.1} and the details are omitted.
\end{proof}

The following particular results is of interest for applications:

\begin{corollary}[Dragomir, 2008, \protect\cite{II.SSDG2}]
\label{II.mcfc.4.1}Consider the sequence of selfadjoint operators $\mathbf{A}%
=\left( A_{1},\dots ,A_{n}\right) $ with $Sp\left( A_{j}\right) \subseteq %
\left[ m,M\right] $ for $j\in \left\{ 1,\dots ,n\right\} $ and for some
scalars $m<M.$ If $f,g:\left[ m,M\right] \longrightarrow \mathbb{R}$ are
continuous with $\delta :=\min_{t\in \left[ m,M\right] }g\left( t\right) $
and $\Delta :=\max_{t\in \left[ m,M\right] }g\left( t\right) ,$ then for any 
$p_{j}\geq 0,j\in \left\{ 1,\dots ,n\right\} $ with $\sum_{j=1}^{n}p_{j}=1$
and $x\in H$ with $\left\Vert x\right\Vert =1$ we have 
\begin{align}
& \left\vert C\left( f,g;\mathbf{A,p,}x\right) \right\vert
\label{II.mcfe.4.6} \\
& \leq \frac{1}{2}\left( \Delta -\delta \right) \left\langle
\sum_{j=1}^{n}p_{j}\left\vert f\left( A_{j}\right) -\left\langle
\sum_{k=1}^{n}p_{k}f\left( A_{k}\right) x,x\right\rangle \cdot
1_{H}\right\vert x,x\right\rangle  \notag \\
& \leq \frac{1}{2}\left( \Delta -\delta \right) C^{1/2}\left( f,f;\mathbf{%
A,p,}x\right) .  \notag
\end{align}
\end{corollary}

\begin{proof}
In we choose in Theorem \ref{II.mcft.4.1} $x_{j}=\sqrt{p_{j}}\cdot x,$ $j\in
\left\{ 1,\dots ,n\right\} ,$ where $p_{j}\geq 0,j\in \left\{ 1,\dots
,n\right\} ,$ $\sum_{j=1}^{n}p_{j}=1$ and $x\in H,$ with $\left\Vert
x\right\Vert =1$ then a simple calculation shows that the inequality (\ref%
{II.mcfe.4.5}) becomes (\ref{II.mcfe.4.6}). The details are omitted.
\end{proof}

In a similar manner we can prove the following result as well:

\begin{theorem}[Dragomir, 2008, \protect\cite{II.SSDG2}]
\label{II.mcft.4.2}Consider the sequences of selfadjoint operators $\mathbf{A%
}=\left( A_{1},\dots ,A_{n}\right) ,$ $\mathbf{B}=\left( B_{1},\dots
,B_{n}\right) $ with $Sp\left( A_{j}\right) ,Sp\left( B_{j}\right) \subseteq %
\left[ m,M\right] $ for $j\in \left\{ 1,\dots ,n\right\} $ and for some
scalars $m<M.$ If $f:\left[ m,M\right] \longrightarrow \mathbb{R}$ is of $%
r-L-$H\"{o}lder type, then we have the Ostrowski type inequality for
sequences of selfadjoint operators: 
\begin{equation}
\left\vert f\left( s\right) -\sum_{j=1}^{n}\left\langle f\left( A_{j}\right)
x_{j},x_{j}\right\rangle \right\vert \leq L\left[ \frac{1}{2}\left(
M-m\right) +\left\vert s-\frac{m+M}{2}\right\vert \right] ^{r},
\label{II.mcfe.4.7}
\end{equation}%
for any $s\in \left[ m,M\right] $ and any $\mathbf{x}=\left( x_{1},\dots
,x_{n}\right) \in H^{n}$ such that $\sum_{j=1}^{n}\left\Vert
x_{j}\right\Vert ^{2}=1$.

Moreover, we have 
\begin{align}
& \left\vert \sum_{j=1}^{n}\left\langle f\left( B_{j}\right)
y_{j},y_{j}\right\rangle -\sum_{k=1}^{n}\left\langle f\left( A_{k}\right)
x_{k},x_{k}\right\rangle \right\vert  \label{II.mcfe.4.8} \\
& \leq \sum_{j=1}^{n}\left\langle \left\vert f\left( B_{j}\right)
-\sum_{k=1}^{n}\left\langle f\left( A_{k}\right) x_{k},x_{k}\right\rangle
\cdot 1_{H}\right\vert y_{j},y_{j}\right\rangle  \notag \\
& \leq L\left[ \frac{1}{2}\left( M-m\right) +\sum_{j=1}^{n}\left\langle
\left\vert B_{j}-\frac{m+M}{2}\cdot 1_{H}\right\vert
y_{j},y_{j}\right\rangle \right] ^{r},  \notag
\end{align}%
for any $\mathbf{x}=\left( x_{1},\dots ,x_{n}\right) ,\mathbf{y}=\left(
y_{1},\dots ,y_{n}\right) \in H^{n}$ such that $\sum_{j=1}^{n}\left\Vert
x_{j}\right\Vert ^{2}=\sum_{j=1}^{n}\left\Vert y_{j}\right\Vert ^{2}=1.$
\end{theorem}

\begin{corollary}[Dragomir, 2008, \protect\cite{II.SSDG2}]
\label{II.mcfc.4.2}Consider the sequences of selfadjoint operators $\mathbf{A%
}=\left( A_{1},\dots ,A_{n}\right) ,$ $\mathbf{B}=\left( B_{1},\dots
,B_{n}\right) $ with $Sp\left( A_{j}\right) ,Sp\left( B_{j}\right) \subseteq %
\left[ m,M\right] $ for $j\in \left\{ 1,\dots ,n\right\} $ and for some
scalars $m<M.$ If $f:\left[ m,M\right] \longrightarrow \mathbb{R}$ is of $%
r-L-$H\"{o}lder type, then for any $p_{j}\geq 0,j\in \left\{ 1,\dots
,n\right\} $ with $\sum_{j=1}^{n}p_{j}=1$ and $x\in H$ with $\left\Vert
x\right\Vert =1$ we have the weighted Ostrowski type inequality for
sequences of selfadjoint operators: 
\begin{equation}
\left\vert f\left( s\right) -\left\langle \sum_{j=1}^{n}p_{j}f\left(
A_{j}\right) x,x\right\rangle \right\vert \leq L\left[ \frac{1}{2}\left(
M-m\right) +\left\vert s-\frac{m+M}{2}\right\vert \right] ^{r},
\label{II.mcfe.4.9}
\end{equation}%
for any $s\in \left[ m,M\right] $.

Moreover, we have 
\begin{align}
& \left\vert \left\langle \sum_{j=1}^{n}q_{j}f\left( B_{j}\right)
y,y\right\rangle -\left\langle \sum_{k=1}^{n}p_{k}f\left( A_{k}\right)
x,x\right\rangle \right\vert  \label{II.mcfe.4.10} \\
& \leq \left\langle \sum_{j=1}^{n}q_{j}\left\vert f\left( B_{j}\right)
-\left\langle \sum_{k=1}^{n}p_{k}f\left( A_{k}\right) x,x\right\rangle \cdot
1_{H}\right\vert y,y\right\rangle  \notag \\
& \leq L\left[ \frac{1}{2}\left( M-m\right) +\left\langle
\sum_{j=1}^{n}q_{j}\left\vert B_{j}-\frac{m+M}{2}\cdot 1_{H}\right\vert
y,y\right\rangle \right] ^{r},  \notag
\end{align}%
for any $q_{k}\geq 0,k\in \left\{ 1,\dots ,n\right\} $ with $%
\sum_{k=1}^{n}q_{k}=1$ and $x,y\in H$ with $\left\Vert x\right\Vert
=\left\Vert y\right\Vert =1.$
\end{corollary}

\subsection{Some Reverses of Jensen's Inequality}

It is clear that all the above inequalities can be applied for various
particular instances of functions $f$ and $g.$ However, in the following we
only consider the inequalities 
\begin{equation}
\left\vert f\left( \left\langle Ax,x\right\rangle \right) -\left\langle
f\left( A\right) x,x\right\rangle \right\vert \leq L\left[ \frac{1}{2}\left(
M-m\right) +\left\vert \left\langle Ax,x\right\rangle -\frac{m+M}{2}%
\right\vert \right] ^{r}  \label{II.mcfe.3.1}
\end{equation}%
for any $x\in H$ with$\left\Vert x\right\Vert =1,$ where the function $f:%
\left[ m,M\right] \rightarrow \mathbb{R}$ is of $r-L-$H\"{o}lder type, and 
\begin{multline}
\left\vert f\left( \left\langle Ax,x\right\rangle \right) -\left\langle
f\left( A\right) x,x\right\rangle \right\vert  \label{II.mcfe.3.2} \\
\leq \left\{ 
\begin{array}{ll}
\left[ \frac{1}{2}\left( M-m\right) +\left\vert \left\langle
Ax,x\right\rangle -\frac{m+M}{2}\right\vert \right] \left\Vert f^{\prime
}\right\Vert _{\infty ,\left[ m,M\right] }, & \text{if }f^{\prime }\in
L_{\infty }\left[ m,M\right] \\ 
&  \\ 
\left[ \frac{1}{2}\left( M-m\right) +\left\vert \left\langle
Ax,x\right\rangle -\frac{m+M}{2}\right\vert \right] ^{q}\left\Vert f^{\prime
}\right\Vert _{p,\left[ m,M\right] }, & 
\begin{array}{l}
\text{if }f^{\prime }\in L_{p}\left[ m,M\right] ; \\ 
p>1,\frac{1}{p}+\frac{1}{q}=1%
\end{array}%
\end{array}%
\right.
\end{multline}%
for any $x\in H$ with $\left\Vert x\right\Vert =1,$ where the function $f:%
\left[ m,M\right] \rightarrow \mathbb{R}$ is absolutely continuous on $\left[
m,M\right] ,$ which are related to the \textit{Jensen's inequality} for
convex functions.

\textbf{1.} Now, if we consider the concave function $f:\left[ m,M\right]
\subset \lbrack 0,\infty )\rightarrow \mathbb{R}$, $f\left( t\right) =t^{r}$
with $r\in \left( 0,1\right) $ and take into account that it is of $r-L-$H%
\"{o}lder type with the constant $L=1,$ then from (\ref{II.mcfe.3.1}) we
derive the following reverse for the \textit{H\"{o}lder-McCarthy inequality 
\cite{II.Mc}} 
\begin{equation}
0\leq \left\langle A^{r}x,x\right\rangle -\left\langle Ax,x\right\rangle
^{r}\leq \left[ \frac{1}{2}\left( M-m\right) +\left\vert \left\langle
Ax,x\right\rangle -\frac{m+M}{2}\right\vert \right] ^{r}  \label{II.mcfe.3.3}
\end{equation}%
for any $x\in H$ with $\left\Vert x\right\Vert =1.$

\textbf{2.} Now, if we consider the functions $f:\left[ m,M\right] \subset
\left( 0,\infty \right) \rightarrow \mathbb{R}$ with $f\left( t\right)
=t^{s} $ and $s\in \left( -\infty ,0\right) \cup \left( 0,\infty \right) ,$
then they are absolutely continuous and 
\begin{equation*}
\left\Vert f^{\prime }\right\Vert _{\infty ,\left[ m,M\right] }=\left\{ 
\begin{array}{ll}
sM^{s-1} & \text{for }s\in \lbrack 1,\infty ), \\ 
&  \\ 
\left\vert s\right\vert m^{s-1} & \text{for }s\in \left( -\infty ,0\right)
\cup \left( 0,1\right) .%
\end{array}%
\right.
\end{equation*}%
If $p\geq 1,$ then 
\begin{align*}
\left\Vert f^{\prime }\right\Vert _{p,\left[ m,M\right] }& =\left\vert
s\right\vert \left( \int_{m}^{M}t^{p\left( s-1\right) }dt\right) ^{1/p} \\
& =\left\vert s\right\vert \times \left\{ 
\begin{array}{ll}
\left( \frac{M^{p\left( s-1\right) +1}-m^{p\left( s-1\right) +1}}{p\left(
s-1\right) +1}\right) ^{1/p} & \text{if }s\neq 1-\frac{1}{p} \\ 
&  \\ 
\left[ \ln \left( \frac{M}{m}\right) \right] ^{1/p} & \text{if }s=1-\frac{1}{%
p}.%
\end{array}%
\right.
\end{align*}

On making use of the first inequality from (\ref{II.mcfe.3.2}) we deduce for
a given $s\in \left( -\infty ,0\right) \cup \left( 0,\infty \right) $ that 
\begin{align}
\left\vert \left\langle Ax,x\right\rangle ^{s}-\left\langle
A^{s}x,x\right\rangle \right\vert & \leq \left[ \frac{1}{2}\left( M-m\right)
+\left\vert \left\langle Ax,x\right\rangle -\frac{m+M}{2}\right\vert \right]
\label{II.mcfe.3.4} \\
& \times \left\{ 
\begin{array}{ll}
sM^{s-1} & \text{for }s\in \lbrack 1,\infty ), \\ 
&  \\ 
\left\vert s\right\vert m^{s-1} & \text{for }s\in \left( -\infty ,0\right)
\cup \left( 0,1\right) .%
\end{array}%
\right.  \notag
\end{align}%
for any $x\in H$ with $\left\Vert x\right\Vert =1.$

The second part of (\ref{II.mcfe.3.2}) will produce the following reverse of
the \textit{H\"{o}lder-McCarthy inequality }as well: 
\begin{align}
\left\vert \left\langle Ax,x\right\rangle ^{s}-\left\langle
A^{s}x,x\right\rangle \right\vert & \leq \left\vert s\right\vert \left[ 
\frac{1}{2}\left( M-m\right) +\left\vert \left\langle Ax,x\right\rangle -%
\frac{m+M}{2}\right\vert \right] ^{q}  \label{II.mcfe.3.5} \\
& \times \left\{ 
\begin{array}{ll}
\left( \frac{M^{p\left( s-1\right) +1}-m^{p\left( s-1\right) +1}}{p\left(
s-1\right) +1}\right) ^{1/p} & \text{if }s\neq 1-\frac{1}{p} \\ 
&  \\ 
\left[ \ln \left( \frac{M}{m}\right) \right] ^{1/p} & \text{if }s=1-\frac{1}{%
p}%
\end{array}%
\right.  \notag
\end{align}%
for any $x\in H$ with $\left\Vert x\right\Vert =1,$ where $s\in \left(
-\infty ,0\right) \cup \left( 0,\infty \right) ,$ $p>1$ and $\frac{1}{p}+%
\frac{1}{q}=1.$

\textbf{3.} Now, if we consider the function $f\left( t\right) =\ln t$
defined on the interval $\left[ m,M\right] \subset \left( 0,\infty \right) ,$
then $f$ is also absolutely continuous and 
\begin{equation*}
\left\Vert f^{\prime }\right\Vert _{p,\left[ m,M\right] }=\left\{ 
\begin{array}{ll}
m^{-1} & \text{for }p=\infty , \\ 
&  \\ 
\left( \frac{M^{p-1}-m^{p-1}}{\left( p-1\right) M^{p-1}m^{p-1}}\right) ^{1/p}
& \text{for }p>1, \\ 
&  \\ 
\ln \left( \frac{M}{m}\right) & \text{for }p=1.%
\end{array}%
\right.
\end{equation*}%
Making use of the first inequality in (\ref{II.mcfe.3.2}) we deduce 
\begin{equation}
0\leq \ln \left( \left\langle Ax,x\right\rangle \right) -\left\langle \ln
\left( A\right) x,x\right\rangle \leq \left[ \frac{1}{2}\left( M-m\right)
+\left\vert \left\langle Ax,x\right\rangle -\frac{m+M}{2}\right\vert \right]
m^{-1}  \label{II.mcfe.3.6}
\end{equation}%
and 
\begin{align}
0& \leq \ln \left( \left\langle Ax,x\right\rangle \right) -\left\langle \ln
\left( A\right) x,x\right\rangle  \label{II.mcfe.3.7} \\
& \leq \left[ \frac{1}{2}\left( M-m\right) +\left\vert \left\langle
Ax,x\right\rangle -\frac{m+M}{2}\right\vert \right] ^{q}\left( \frac{%
M^{p-1}-m^{p-1}}{\left( p-1\right) M^{p-1}m^{p-1}}\right) ^{1/p}  \notag
\end{align}%
for any $x\in H$ with $\left\Vert x\right\Vert =1,$ where $p>1$ and $\frac{1%
}{p}+\frac{1}{q}=1.$

Similar results can be stated for sequences of operators, however the
details are left to the interested reader.

\subsection{Some Particular Gr\"{u}ss' Type Inequalities}

In this last section we provide some particular cases that can be obtained
via the Gr\"{u}ss' type inequalities established before. For this purpose we
select only two examples as follows.

Let $A$ be a selfadjoint operator with $Sp\left( A\right) \subseteq \left[
m,M\right] $ for some real numbers $m<M.$ If $g:\left[ m,M\right]
\longrightarrow \mathbb{R}$ is continuous with $\delta :=\min_{t\in \left[
m,M\right] }g\left( t\right) $ and $\Delta :=\max_{t\in \left[ m,M\right]
}g\left( t\right) ,$ then for any $f:\left[ m,M\right] \longrightarrow 
\mathbb{R}$ of $r-L-$H\"{o}lder type we have the inequality: 
\begin{align}
& \left\vert \left\langle f\left( A\right) g\left( A\right) x,x\right\rangle
-\left\langle f\left( A\right) x,x\right\rangle \cdot \left\langle g\left(
A\right) x,x\right\rangle \right\vert  \label{II.mcfe.5.1} \\
& \leq \frac{1}{2}\left( \Delta -\delta \right) L\left[ \frac{1}{2}\left(
M-m\right) +\left\langle \left\vert A-\frac{m+M}{2}\cdot 1_{H}\right\vert
x,x\right\rangle \right] ^{r},  \notag
\end{align}%
for any $x\in H$ with $\left\Vert x\right\Vert =1.$

Moreover, if $f$ is absolutely continuos on $\left[ m,M\right] ,$ then we
have the inequalities: 
\begin{multline}
\left\vert \left\langle f\left( A\right) g\left( A\right) x,x\right\rangle
-\left\langle f\left( A\right) x,x\right\rangle \cdot \left\langle g\left(
A\right) x,x\right\rangle \right\vert \leq \frac{1}{2}\left( \Delta -\delta
\right)  \label{II.mcfe.5.2} \\
\times \left\{ 
\begin{array}{ll}
\left[ \frac{1}{2}\left( M-m\right) +\left\langle \left\vert A-\frac{m+M}{2}%
\cdot 1_{H}\right\vert x,x\right\rangle \right] \left\Vert f^{\prime
}\right\Vert _{\infty ,\left[ m,M\right] } & \text{ if \ }f^{\prime }\in
L_{\infty }\left[ m,M\right] ; \\ 
&  \\ 
\left[ \frac{1}{2}\left( M-m\right) +\left\langle \left\vert A-\frac{m+M}{2}%
\cdot 1_{H}\right\vert x,x\right\rangle \right] ^{1/q}\left\Vert f^{\prime
}\right\Vert _{p,\left[ m,M\right] } & 
\begin{array}{l}
\text{if \ }f^{\prime }\in L_{p}\left[ m,M\right] , \\ 
p,q>1,\frac{1}{p}+\frac{1}{q}=1%
\end{array}%
\end{array}%
\right.
\end{multline}%
for any $x\in H$ with $\left\Vert x\right\Vert =1.$

\textbf{1.} If we consider the concave function $f:\left[ m,M\right] \subset
\lbrack 0,\infty )\rightarrow \mathbb{R}$, $f\left( t\right) =t^{r}$ with $%
r\in \left( 0,1\right) $ and take into account that it is of $r-L-$H\"{o}%
lder type with the constant $L=1,$ then from (\ref{II.mcfe.5.1}) we derive
the following result: 
\begin{align}
& \left\vert \left\langle A^{r}g\left( A\right) x,x\right\rangle
-\left\langle A^{r}x,x\right\rangle \cdot \left\langle g\left( A\right)
x,x\right\rangle \right\vert  \label{II.mcfe.5.3} \\
& \leq \frac{1}{2}\left( \Delta -\delta \right) \left[ \frac{1}{2}\left(
M-m\right) +\left\langle \left\vert A-\frac{m+M}{2}\cdot 1_{H}\right\vert
x,x\right\rangle \right] ^{r},  \notag
\end{align}%
for any $x\in H$ with $\left\Vert x\right\Vert =1,$ where $g:\left[ m,M%
\right] \longrightarrow \mathbb{R}$ is continuous with $\delta :=\min_{t\in %
\left[ m,M\right] }g\left( t\right) $ and $\Delta :=\max_{t\in \left[ m,M%
\right] }g\left( t\right) .$

Now, consider the function $g:\left[ m,M\right] \subset \left( 0,\infty
\right) \rightarrow \mathbb{R}$, $g\left( t\right) =t^{p}$ with $p\in \left(
-\infty ,0\right) \cup \left( 0,\infty \right) $. Obviously, 
\begin{equation*}
\Delta -\delta =\left\{ 
\begin{array}{ll}
M^{p}-m^{p} & \text{ if }p>0, \\ 
&  \\ 
\frac{M^{-p}-m^{-p}}{M^{-p}m^{-p}} & \text{if }p<0,%
\end{array}%
\right.
\end{equation*}%
and by (\ref{II.mcfe.5.3}) we get for any $x\in H$ with $\left\Vert
x\right\Vert =1$ that 
\begin{align}
0& \leq \left\langle A^{r+p}x,x\right\rangle -\left\langle
A^{r}x,x\right\rangle \cdot \left\langle A^{p}x,x\right\rangle
\label{II.mcfe.5.4} \\
& \leq \frac{1}{2}\left( M^{p}-m^{p}\right) \left[ \frac{1}{2}\left(
M-m\right) +\left\langle \left\vert A-\frac{m+M}{2}\cdot 1_{H}\right\vert
x,x\right\rangle \right] ^{r},  \notag
\end{align}%
when $p>0$ and 
\begin{align}
0& \leq \left\langle A^{r}x,x\right\rangle \cdot \left\langle
A^{p}x,x\right\rangle -\left\langle A^{r+p}x,x\right\rangle
\label{II.mcfe.5.5} \\
& \leq \frac{1}{2}\cdot \frac{M^{-p}-m^{-p}}{M^{-p}m^{-p}}\left[ \frac{1}{2}%
\left( M-m\right) +\left\langle \left\vert A-\frac{m+M}{2}\cdot
1_{H}\right\vert x,x\right\rangle \right] ^{r},  \notag
\end{align}%
when $p<0.$

If $g:\left[ m,M\right] \subset \left( 0,\infty \right) \rightarrow \mathbb{R%
}$, $g\left( t\right) =\ln t,$ then by (\ref{II.mcfe.5.3}) we also get the
inequality for logarithm: 
\begin{align}
0& \leq \left\langle A^{r}\ln Ax,x\right\rangle -\left\langle
A^{r}x,x\right\rangle \cdot \left\langle \ln Ax,x\right\rangle
\label{II.mcfe.5.6} \\
& \leq \ln \sqrt{\frac{M}{m}}\cdot \left[ \frac{1}{2}\left( M-m\right)
+\left\langle \left\vert A-\frac{m+M}{2}\cdot 1_{H}\right\vert
x,x\right\rangle \right] ^{r},  \notag
\end{align}%
for any $x\in H$ with $\left\Vert x\right\Vert =1.$

\textbf{2. }Now consider the functions $f,g:\left[ m,M\right] \subset \left(
0,\infty \right) \rightarrow \mathbb{R}$, with $f\left( t\right) =t^{s}$ and 
$g\left( t\right) =t^{w}$ with $s,w\in \left( -\infty ,0\right) \cup \left(
0,\infty \right) .$ We have 
\begin{equation*}
\left\Vert f^{\prime }\right\Vert _{\infty ,\left[ m,M\right] }=\left\{ 
\begin{array}{ll}
sM^{s-1} & \text{for }s\in \lbrack 1,\infty ), \\ 
&  \\ 
\left\vert s\right\vert m^{s-1} & \text{for }s\in \left( -\infty ,0\right)
\cup \left( 0,1\right) .%
\end{array}%
\right.
\end{equation*}%
and, for $p\geq 1,$ 
\begin{equation*}
\left\Vert f^{\prime }\right\Vert _{p,\left[ m,M\right] }=\left\vert
s\right\vert \times \left\{ 
\begin{array}{ll}
\left( \frac{M^{p\left( s-1\right) +1}-m^{p\left( s-1\right) +1}}{p\left(
s-1\right) +1}\right) ^{1/p} & \text{if }s\neq 1-\frac{1}{p} \\ 
&  \\ 
\left[ \ln \left( \frac{M}{m}\right) \right] ^{1/p} & \text{if }s=1-\frac{1}{%
p}.%
\end{array}%
\right.
\end{equation*}%
If $w>0,$ then by the first inequality in (\ref{II.mcfe.5.2}) we have 
\begin{align}
& \left\vert \left\langle A^{s+w}x,x\right\rangle -\left\langle
A^{s}x,x\right\rangle \cdot \left\langle A^{w}x,x\right\rangle \right\vert
\label{II.mcfe.5.7} \\
& \leq \frac{1}{2}\left( M^{w}-m^{w}\right) \left[ \frac{1}{2}\left(
M-m\right) +\left\langle \left\vert A-\frac{m+M}{2}\cdot 1_{H}\right\vert
x,x\right\rangle \right]  \notag \\
& \times \left\{ 
\begin{array}{ll}
sM^{s-1} & \text{for }s\in \lbrack 1,\infty ), \\ 
&  \\ 
\left\vert s\right\vert m^{s-1} & \text{for }s\in \left( -\infty ,0\right)
\cup \left( 0,1\right) .%
\end{array}%
\right.  \notag
\end{align}%
for any $x\in H$ with $\left\Vert x\right\Vert =1.$

If $w<0,$ then by the same inequality we also have 
\begin{align}
& \left\vert \left\langle A^{s+w}x,x\right\rangle -\left\langle
A^{s}x,x\right\rangle \cdot \left\langle A^{w}x,x\right\rangle \right\vert
\label{II.mcfe.5.8} \\
& \leq \frac{1}{2}\cdot \frac{M^{-w}-m^{-w}}{M^{-w}m^{-w}}\left[ \frac{1}{2}%
\left( M-m\right) +\left\langle \left\vert A-\frac{m+M}{2}\cdot
1_{H}\right\vert x,x\right\rangle \right]  \notag \\
& \times \left\{ 
\begin{array}{ll}
sM^{s-1} & \text{for }s\in \lbrack 1,\infty ), \\ 
&  \\ 
\left\vert s\right\vert m^{s-1} & \text{for }s\in \left( -\infty ,0\right)
\cup \left( 0,1\right) ,%
\end{array}%
\right.  \notag
\end{align}%
for any $x\in H$ with $\left\Vert x\right\Vert =1.$

Finally, if we assume that $p>1$ and $w>0,$ then by the second inequality in
(\ref{II.mcfe.5.2}) we have 
\begin{align}
& \left\vert \left\langle A^{s+w}x,x\right\rangle -\left\langle
A^{s}x,x\right\rangle \cdot \left\langle A^{w}x,x\right\rangle \right\vert
\label{II.mcfe.5.9} \\
& \leq \frac{1}{2}\left\vert s\right\vert \left( M^{w}-m^{w}\right) \left[ 
\frac{1}{2}\left( M-m\right) +\left\langle \left\vert A-\frac{m+M}{2}\cdot
1_{H}\right\vert x,x\right\rangle \right] ^{1/q}  \notag \\
& \times \left\{ 
\begin{array}{ll}
\left( \frac{M^{p\left( s-1\right) +1}-m^{p\left( s-1\right) +1}}{p\left(
s-1\right) +1}\right) ^{1/p} & \text{if }s\neq 1-\frac{1}{p} \\ 
&  \\ 
\left[ \ln \left( \frac{M}{m}\right) \right] ^{1/p} & \text{if }s=1-\frac{1}{%
p},%
\end{array}%
\right.  \notag
\end{align}%
while for $w<0,$ we also have 
\begin{align}
& \left\vert \left\langle A^{s+w}x,x\right\rangle -\left\langle
A^{s}x,x\right\rangle \cdot \left\langle A^{w}x,x\right\rangle \right\vert
\label{II.mcfe.5.10} \\
& \leq \frac{1}{2}\left\vert s\right\vert \cdot \frac{M^{-w}-m^{-w}}{%
M^{-w}m^{-w}}\left[ \frac{1}{2}\left( M-m\right) +\left\langle \left\vert A-%
\frac{m+M}{2}\cdot 1_{H}\right\vert x,x\right\rangle \right] ^{1/q}  \notag
\\
& \times \left\{ 
\begin{array}{ll}
\left( \frac{M^{p\left( s-1\right) +1}-m^{p\left( s-1\right) +1}}{p\left(
s-1\right) +1}\right) ^{1/p} & \text{if }s\neq 1-\frac{1}{p} \\ 
&  \\ 
\left[ \ln \left( \frac{M}{m}\right) \right] ^{1/p} & \text{if }s=1-\frac{1}{%
p},%
\end{array}%
\right.  \notag
\end{align}%
where $q>1$ with $\frac{1}{p}+\frac{1}{q}=1$ and $x\in H$ with $\left\Vert
x\right\Vert =1.$

\section{Bounds for the \v Ceby\v sev Functional of Lipschitzian Functions}

\subsection{The Case of Lipschitzian Functions}

The following result can be stated:

\begin{theorem}[Dragomir, 2008, \protect\cite{II.SSDG3}]
\label{II.op2funct.2.1}Let $A$ be a selfadjoint operator with $Sp\left(
A\right) \subseteq \left[ m,M\right] $ for some real numbers $m<M.$ If $f:%
\left[ m,M\right] \longrightarrow \mathbb{R}$ is Lipschitzian with the
constant $L>0$ and $g:\left[ m,M\right] \longrightarrow \mathbb{R}$ is
continuous with $\delta :=\min_{t\in \left[ m,M\right] }g\left( t\right) $
and $\Delta :=\max_{t\in \left[ m,M\right] }g\left( t\right) ,$ then 
\begin{equation}
\left\vert C\left( f,g;A;x\right) \right\vert \leq \frac{1}{2}\left( \Delta
-\delta \right) L\left\langle \ell _{A,x}\left( A\right) x,x\right\rangle
\leq \left( \Delta -\delta \right) LC\left( e,e;A;x\right)
\label{II.op2funce.2.1.1}
\end{equation}%
for any $x\in H$ with $\left\Vert x\right\Vert =1,$ where 
\begin{equation*}
\ell _{A,x}\left( t\right) :=\left\langle \left\vert t\cdot
1_{H}-A\right\vert x,x\right\rangle
\end{equation*}%
is a continuous function on $\left[ m,M\right] ,$ $e\left( t\right) =t$ and 
\begin{equation}
C\left( e,e;A;x\right) =\left\Vert Ax\right\Vert ^{2}-\left\langle
Ax,x\right\rangle ^{2}\left( \geq 0\right) .  \label{II.op2funce.2.2}
\end{equation}
\end{theorem}

\begin{proof}
First of all, by the Jensen inequality for convex functions of selfadjoint
operators (see for instance \cite[p. 5]{II.FMPS}) applied for the modulus,
we can state that 
\begin{equation}
\left\vert \left\langle h\left( A\right) x,x\right\rangle \right\vert \leq
\left\langle \left\vert h\left( A\right) \right\vert x,x\right\rangle 
\tag{M}  \label{II.op2funcM}
\end{equation}%
for any $x\in H$ with $\left\Vert x\right\Vert =1,$ where $h$ is a
continuous function on $\left[ m,M\right] .$

Since $f$ is Lipschitzian with the constant $L>0,$ then for any $t,s\in %
\left[ m,M\right] $ we have 
\begin{equation}
\left\vert f\left( t\right) -f\left( s\right) \right\vert \leq L\left\vert
t-s\right\vert .  \label{II.op2funce.2.3}
\end{equation}%
Now, if we fix $t\in \left[ m,M\right] $ and apply the property (\ref{P})
for the inequality (\ref{II.op2funce.2.3}) and the operator $A$ we get 
\begin{equation}
\left\langle \left\vert f\left( t\right) \cdot 1_{H}-f\left( A\right)
\right\vert x,x\right\rangle \leq L\left\langle \left\vert t\cdot
1_{H}-A\right\vert x,x\right\rangle ,  \label{II.op2funce.2.5}
\end{equation}%
for any $x\in H$ with $\left\Vert x\right\Vert =1.$

Utilising the property (\ref{II.op2funcM}) we get 
\begin{equation*}
\left\vert f\left( t\right) -\left\langle f\left( A\right) x,x\right\rangle
\right\vert =\left\vert \left\langle f\left( t\right) \cdot 1_{H}-f\left(
A\right) x,x\right\rangle \right\vert \leq \left\langle \left\vert f\left(
t\right) \cdot 1_{H}-f\left( A\right) \right\vert x,x\right\rangle
\end{equation*}%
which together with (\ref{II.op2funce.2.5}) gives 
\begin{equation}
\left\vert f\left( t\right) -\left\langle f\left( A\right) x,x\right\rangle
\right\vert \leq L\ell _{A,x}\left( t\right)  \label{II.op2funce.2.6}
\end{equation}%
for any $t\in \left[ m,M\right] $ and for any $x\in H$ with $\left\Vert
x\right\Vert =1.$

Since $\delta :=\min_{t\in \left[ m,M\right] }g\left( t\right) $ and $\Delta
:=\max_{t\in \left[ m,M\right] }g\left( t\right) ,$ we also have 
\begin{equation}
\left\vert g\left( t\right) -\frac{\Delta +\delta }{2}\right\vert \leq \frac{%
1}{2}\left( \Delta -\delta \right)  \label{II.op2funce.2.7}
\end{equation}%
for any $t\in \left[ m,M\right] $ and for any $x\in H$ with $\left\Vert
x\right\Vert =1.$

If we multiply the inequality (\ref{II.op2funce.2.6}) with (\ref%
{II.op2funce.2.7}) we get 
\begin{align}
& \left\vert f\left( t\right) g\left( t\right) -\left\langle f\left(
A\right) x,x\right\rangle g\left( t\right) -\frac{\Delta +\delta }{2}f\left(
t\right) +\frac{\Delta +\delta }{2}\left\langle f\left( A\right)
x,x\right\rangle \right\vert  \label{II.op2funce.2.8} \\
& \leq \frac{1}{2}\left( \Delta -\delta \right) L\ell _{A,x}\left( t\right) =%
\frac{1}{2}\left( \Delta -\delta \right) L\left\langle \left\vert t\cdot
1_{H}-A\right\vert x,x\right\rangle  \notag \\
& \leq \frac{1}{2}\left( \Delta -\delta \right) L\left\langle \left\vert
t\cdot 1_{H}-A\right\vert ^{2}x,x\right\rangle ^{1/2}  \notag \\
& =\frac{1}{2}\left( \Delta -\delta \right) L\left( \left\langle
A^{2}x,x\right\rangle -2\left\langle Ax,x\right\rangle t+t^{2}\right) ^{1/2},
\notag
\end{align}%
for any $t\in \left[ m,M\right] $ and for any $x\in H$ with $\left\Vert
x\right\Vert =1.$

Now, if we apply the property (\ref{P}) for the inequality (\ref%
{II.op2funce.2.8}) and a selfadjoint operator $B$ with $Sp\left( B\right)
\subset \left[ m,M\right] ,$ then we get the following inequality of
interest in itself: 
\begin{align}
& \left\langle f\left( B\right) g\left( B\right) y,y\right\rangle
-\left\langle f\left( A\right) x,x\right\rangle \left\langle g\left(
B\right) y,y\right\rangle  \label{II.op2funce.2.9} \\
& \left. -\frac{\Delta +\delta }{2}\left\langle f\left( B\right)
y,y\right\rangle +\frac{\Delta +\delta }{2}\left\langle f\left( A\right)
x,x\right\rangle \right\vert  \notag \\
& \leq \frac{1}{2}\left( \Delta -\delta \right) L\left\langle \ell
_{A,x}\left( B\right) y,y\right\rangle  \notag \\
& \leq \frac{1}{2}\left( \Delta -\delta \right) L\left\langle \left(
\left\langle A^{2}x,x\right\rangle 1_{H}-2\left\langle Ax,x\right\rangle
B+B^{2}\right) ^{1/2}y,y\right\rangle  \notag \\
& \leq \frac{1}{2}\left( \Delta -\delta \right) L\left( \left\langle
A^{2}x,x\right\rangle -2\left\langle Ax,x\right\rangle \left\langle
By,y\right\rangle +\left\langle B^{2}y,y\right\rangle \right) ,  \notag
\end{align}%
for any $x,y\in H$ with $\left\Vert x\right\Vert =\left\Vert y\right\Vert
=1. $

Finally, if we choose in (\ref{II.op2funce.2.9}) $y=x$ and $B=A,$ then we
deduce the desired result (\ref{II.op2funce.2.1.1}).
\end{proof}

In the case of two Lipschitzian functions, the following result may be
stated as well:

\begin{theorem}[Dragomir, 2008, \protect\cite{II.SSDG3}]
\label{II.op2funct.2.2}Let $A$ be a selfadjoint operator with $Sp\left(
A\right) \subseteq \left[ m,M\right] $ for some real numbers $m<M.$ If $f,g:%
\left[ m,M\right] \longrightarrow \mathbb{R}$ are Lipschitzian with the
constants $L,K>0,$ then 
\begin{equation}
\left\vert C\left( f,g;A;x\right) \right\vert \leq LKC\left( e,e;A;x\right) ,
\label{II.op2funce.2.10}
\end{equation}%
for any $x\in H$ with $\left\Vert x\right\Vert =1.$
\end{theorem}

\begin{proof}
Since $f,g:\left[ m,M\right] \longrightarrow \mathbb{R}$ are Lipschitzian,
then 
\begin{equation*}
\left| f\left( t\right) -f\left( s\right) \right| \leq L\left| t-s\right| 
\text{ and }\left| g\left( t\right) -g\left( s\right) \right| \leq K\left|
t-s\right|
\end{equation*}
for any $t,s\in \left[ m,M\right] ,$ which gives the inequality 
\begin{equation*}
\left| f\left( t\right) g\left( t\right) -f\left( t\right) g\left( s\right)
-f\left( s\right) g\left( t\right) +f\left( s\right) g\left( s\right)
\right| \leq KL\left( t^{2}-2ts+s^{2}\right)
\end{equation*}
for any $t,s\in \left[ m,M\right] .$

Now, fix $t\in \left[ m,M\right] $ and if we apply the properties (\ref{P})
and (\ref{II.op2funcM}) for the operator $A$ we get successively 
\begin{align}
& \left\vert f\left( t\right) g\left( t\right) -\left\langle g\left(
A\right) x,x\right\rangle f\left( t\right) -\left\langle f\left( A\right)
x,x\right\rangle g\left( t\right) +\left\langle f\left( A\right) g\left(
A\right) x,x\right\rangle \right\vert  \label{II.op2funce.2.11} \\
& =\left\vert \left\langle \left[ f\left( t\right) g\left( t\right) \cdot
1_{H}-f\left( t\right) g\left( A\right) -f\left( A\right) g\left( t\right)
+f\left( A\right) g\left( A\right) \right] x,x\right\rangle \right\vert 
\notag \\
& \leq \left\langle \left\vert f\left( t\right) g\left( t\right) \cdot
1_{H}-f\left( t\right) g\left( A\right) -f\left( A\right) g\left( t\right)
+f\left( A\right) g\left( A\right) \right\vert x,x\right\rangle  \notag \\
& \leq KL\left\langle \left( t^{2}\cdot 1_{H}-2tA+A^{2}\right)
x,x\right\rangle =KL\left( t^{2}-2t\left\langle Ax,x\right\rangle
+\left\langle A^{2}x,x\right\rangle \right)  \notag
\end{align}%
for any $x\in H$ with $\left\Vert x\right\Vert =1.$

Further, fix $x\in H$ with $\left\Vert x\right\Vert =1.$ On applying the
same properties for the inequality (\ref{II.op2funce.2.11}) and another
selfadjoint operator $B$ with $Sp\left( B\right) \subset \left[ m,M\right] ,$
we have 
\begin{align}
& \left\vert \left\langle f\left( B\right) g\left( B\right) y,y\right\rangle
-\left\langle g\left( A\right) x,x\right\rangle \left\langle f\left(
B\right) y,y\right\rangle \right.  \label{II.op2funce.2.12.1} \\
& \left. -\left\langle f\left( A\right) x,x\right\rangle \left\langle
g\left( B\right) y,y\right\rangle +\left\langle f\left( A\right) g\left(
A\right) x,x\right\rangle \right\vert  \notag
\end{align}%
\begin{align*}
& =\left\vert \left\langle \left[ f\left( B\right) g\left( B\right)
-\left\langle g\left( A\right) x,x\right\rangle f\left( B\right)
-\left\langle f\left( A\right) x,x\right\rangle g\left( B\right)
+\left\langle f\left( A\right) g\left( A\right) x,x\right\rangle 1_{H}\right]
y,y\right\rangle \right\vert \\
& \leq \left\langle \left\vert f\left( B\right) g\left( B\right)
-\left\langle g\left( A\right) x,x\right\rangle f\left( B\right)
-\left\langle f\left( A\right) x,x\right\rangle g\left( B\right)
+\left\langle f\left( A\right) g\left( A\right) x,x\right\rangle
1_{H}\right\vert y,y\right\rangle \\
& \leq KL\left\langle \left( B^{2}-2\left\langle Ax,x\right\rangle
B+\left\langle A^{2}x,x\right\rangle 1_{H}\right) y,y\right\rangle \\
& =KL\left( \left\langle B^{2}y,y\right\rangle -2\left\langle
Ax,x\right\rangle \left\langle By,y\right\rangle +\left\langle
A^{2}x,x\right\rangle \right)
\end{align*}%
for any $x,y\in H$ with $\left\Vert x\right\Vert =\left\Vert y\right\Vert
=1, $ which is an inequality of interest in its own right.

Finally, on making $B=A$ and $y=x$ in (\ref{II.op2funce.2.12.1}) we deduce
the desired result (\ref{II.op2funce.2.10}).
\end{proof}

\subsection{Some Inequalities for Sequences of Operators}

Consider the sequence of selfadjoint operators $\mathbf{A}=\left(
A_{1},\dots ,A_{n}\right) $ with $Sp\left( A_{j}\right) \subseteq \left[ m,M%
\right] $ for $j\in \left\{ 1,\dots ,n\right\} $ and for some scalars $m<M.$
If $\mathbf{x}=\left( x_{1},\dots ,x_{n}\right) \in H^{n}$ are such that $%
\sum_{j=1}^{n}\left\Vert x_{j}\right\Vert ^{2}=1,$ then we can consider the
following \v{C}eby\v{s}ev type functional 
\begin{equation*}
C\left( f,g;\mathbf{A,x}\right) :=\sum_{j=1}^{n}\left\langle f\left(
A_{j}\right) g\left( A_{j}\right) x_{j},x_{j}\right\rangle
-\sum_{j=1}^{n}\left\langle f\left( A_{j}\right) x_{j},x_{j}\right\rangle
\cdot \sum_{j=1}^{n}\left\langle g\left( A_{j}\right)
x_{j},x_{j}\right\rangle .
\end{equation*}%
As a particular case of the above functional and for a probability sequence $%
\mathbf{p=}\left( p_{1},\dots ,p_{n}\right) ,$ i.e., $p_{j}\geq 0$ for $j\in
\left\{ 1,\dots ,n\right\} $ and $\sum_{j=1}^{n}p_{j}=1,$ we can also
consider the functional 
\begin{align*}
C\left( f,g;\mathbf{A,p,}x\right) & :=\left\langle
\sum_{j=1}^{n}p_{j}f\left( A_{j}\right) g\left( A_{j}\right) x,x\right\rangle
\\
& -\left\langle \sum_{j=1}^{n}p_{j}f\left( A_{j}\right) x,x\right\rangle
\cdot \left\langle \sum_{j=1}^{n}p_{j}g\left( A_{j}\right) x,x\right\rangle
\end{align*}%
where $x\in H,$ $\left\Vert x\right\Vert =1.$

We know, from \cite{II.SSD} that for the sequence of selfadjoint operators $%
\mathbf{A}=\left( A_{1},\dots ,A_{n}\right) $ with $Sp\left( A_{j}\right)
\subseteq \left[ m,M\right] $ for $j\in \left\{ 1,\dots ,n\right\} $ and for
the synchronous (asynchronous) functions $f,g:\left[ m,M\right]
\longrightarrow \mathbb{R}$ we have the inequality 
\begin{equation}
C\left( f,g;\mathbf{A,x}\right) \geq \left( \leq \right) 0\text{ }
\label{II.op2funce.3.1}
\end{equation}%
for any $\mathbf{x}=\left( x_{1},\dots ,x_{n}\right) \in H^{n}$ with $%
\sum_{j=1}^{n}\left\Vert x_{j}\right\Vert ^{2}=1.$ Also, for any probability
distribution $\mathbf{p=}\left( p_{1},\dots ,p_{n}\right) $ and any $x\in H,$
$\left\Vert x\right\Vert =1$ we have 
\begin{equation}
C\left( f,g;\mathbf{A,p,}x\right) \geq \left( \leq \right) 0.
\label{II.op2funce.3.2}
\end{equation}

On the other hand, the following Gr\"{u}ss' type inequality is valid as well 
\cite{II.SSD}: 
\begin{equation}
\left\vert C\left( f,g;\mathbf{A,x}\right) \right\vert \leq \frac{1}{2}\cdot
\left( \Gamma -\gamma \right) \left[ C\left( g,g;\mathbf{A,x}\right) \right]
^{1/2}\left( \leq \frac{1}{4}\left( \Gamma -\gamma \right) \left( \Delta
-\delta \right) \right)  \label{II.op2funce.3.3}
\end{equation}%
for any $\mathbf{x}=\left( x_{1},\dots ,x_{n}\right) \in H^{n}$ with $%
\sum_{j=1}^{n}\left\Vert x_{j}\right\Vert ^{2}=1,$ where $f$ and $g$ are
continuous on $\left[ m,M\right] $ and $\gamma :=\min_{t\in \left[ m,M\right]
}f\left( t\right) $, $\Gamma :=\max_{t\in \left[ m,M\right] }f\left(
t\right) $, $\delta :=\min_{t\in \left[ m,M\right] }g\left( t\right) $ and $%
\Delta :=\max_{t\in \left[ m,M\right] }g\left( t\right) .$

Similarly, for any probability distribution $\mathbf{p=}\left( p_{1},\dots
,p_{n}\right) $ and any $x\in H,$ $\left\Vert x\right\Vert =1$ we also have
the inequality: 
\begin{equation}
\left\vert C\left( f,g;\mathbf{A,p,}x\right) \right\vert \leq \frac{1}{2}%
\cdot \left( \Gamma -\gamma \right) \left[ C\left( g,g;\mathbf{A,p,}x\right) %
\right] ^{1/2}\left( \leq \frac{1}{4}\left( \Gamma -\gamma \right) \left(
\Delta -\delta \right) \right) .  \label{II.op2funce.3.4}
\end{equation}

We can state now the following new result:

\begin{theorem}[Dragomir, 2008, \protect\cite{II.SSDG3}]
\label{II.op2funct.3.1}Let $\mathbf{A}=\left( A_{1},\dots ,A_{n}\right) $ be
a sequence of selfadjoint operators with $Sp\left( A_{j}\right) \subseteq %
\left[ m,M\right] $ for $j\in \left\{ 1,\dots ,n\right\} $ and for some
scalars $m<M.$ If $f:\left[ m,M\right] \longrightarrow \mathbb{R}$ is
Lipschitzian with the constant $L>0$ and $g:\left[ m,M\right]
\longrightarrow \mathbb{R}$ is continuous with $\delta :=\min_{t\in \left[
m,M\right] }g\left( t\right) $ and $\Delta :=\max_{t\in \left[ m,M\right]
}g\left( t\right) ,$ then 
\begin{align}
\left\vert C\left( f,g;\mathbf{A,x}\right) \right\vert & \leq \frac{1}{2}%
\left( \Delta -\delta \right) L\sum_{k=1}^{n}\left\langle \ell _{\mathbf{A},%
\mathbf{x}}\left( A_{k}\right) x_{k},x_{k}\right\rangle
\label{II.op2funce.3.5} \\
& \leq \left( \Delta -\delta \right) LC\left( e,e;\mathbf{A};\mathbf{x}%
\right)  \notag
\end{align}%
for any $\mathbf{x}=\left( x_{1},\dots ,x_{n}\right) \in H^{n}$ with $%
\sum_{j=1}^{n}\left\Vert x_{j}\right\Vert ^{2}=1,$ where 
\begin{equation*}
\ell _{\mathbf{A},\mathbf{x}}\left( t\right) :=\sum_{j=1}^{n}\left\langle
\left\vert t\cdot 1_{H}-A_{j}\right\vert x_{j},x_{j}\right\rangle
\end{equation*}%
is a continuous function on $\left[ m,M\right] ,$ $e\left( t\right) =t$ and 
\begin{equation*}
C\left( e,e;\mathbf{A};\mathbf{x}\right) =\sum_{j=1}^{n}\left\Vert
Ax_{j}\right\Vert ^{2}-\left( \sum_{j=1}^{n}\left\langle
A_{j}x_{j},x_{j}\right\rangle \right) ^{2}\left( \geq 0\right) .
\end{equation*}
\end{theorem}

\begin{proof}
Follows from Theorem \ref{II.op2funct.2.1}. The details are omitted.
\end{proof}

As a particular case we have:

\begin{corollary}[Dragomir, 2008, \protect\cite{II.SSDG3}]
\label{op2funcc.3.1}Let $\mathbf{A}=\left( A_{1},\dots ,A_{n}\right) $ be a
sequence of selfadjoint operators with $Sp\left( A_{j}\right) \subseteq %
\left[ m,M\right] $ for $j\in \left\{ 1,\dots ,n\right\} $ and for some
scalars $m<M.$ If $f:\left[ m,M\right] \longrightarrow \mathbb{R}$ is
Lipschitzian with the constant $L>0$ and $g:\left[ m,M\right]
\longrightarrow \mathbb{R}$ is continuous with $\delta :=\min_{t\in \left[
m,M\right] }g\left( t\right) $ and $\Delta :=\max_{t\in \left[ m,M\right]
}g\left( t\right) ,$ then for any $p_{j}\geq 0,j\in \left\{ 1,\dots
,n\right\} $ with $\sum_{j=1}^{n}p_{j}=1$ and $x\in H$ with $\left\Vert
x\right\Vert =1$ we have 
\begin{align}
\left\vert C\left( f,g;\mathbf{A,p,}x\right) \right\vert & \leq \frac{1}{2}%
\left( \Delta -\delta \right) L\left\langle \sum_{k=1}^{n}p_{k}\ell _{%
\mathbf{A},\mathbf{p,}x}\left( A_{k}\right) x,x\right\rangle
\label{II.op2funce.3.6} \\
& \leq \left( \Delta -\delta \right) LC\left( e,e;\mathbf{A,p,}x\right) 
\notag
\end{align}%
where 
\begin{equation*}
\ell _{\mathbf{A},\mathbf{p,}x}\left( t\right) :=\left\langle
\sum_{j=1}^{n}p_{j}\left\vert t\cdot 1_{H}-A_{j}\right\vert x,x\right\rangle
\end{equation*}%
is a continuous function on $\left[ m,M\right] $ and 
\begin{equation*}
C\left( e,e;\mathbf{A,p,}x\right) =\sum_{j=1}^{n}p_{j}\left\Vert
Ax_{j}\right\Vert ^{2}-\left\langle \sum_{j=1}^{n}p_{j}A_{j}x,x\right\rangle
^{2}\left( \geq 0\right) .
\end{equation*}
\end{corollary}

\begin{proof}
In we choose in Theorem \ref{II.op2funct.3.1} $x_{j}=\sqrt{p_{j}}\cdot x,$ $%
j\in \left\{ 1,\dots ,n\right\} ,$ where $p_{j}\geq 0,j\in \left\{ 1,\dots
,n\right\} ,$ $\sum_{j=1}^{n}p_{j}=1$ and $x\in H,$ with $\left\Vert
x\right\Vert =1$ then a simple calculation shows that the inequality (\ref%
{II.op2funce.3.5}) becomes (\ref{II.op2funce.3.6}). The details are omitted.
\end{proof}

In a similar way we obtain the following results as well:

\begin{theorem}[Dragomir, 2008, \protect\cite{II.SSDG3}]
\label{II.op2funct.3.2}Let $\mathbf{A}=\left( A_{1},\dots ,A_{n}\right) $ be
a sequence of selfadjoint operators with $Sp\left( A_{j}\right) \subseteq %
\left[ m,M\right] $ for $j\in \left\{ 1,\dots ,n\right\} $ and for some
scalars $m<M.$ If $f,g:\left[ m,M\right] \longrightarrow \mathbb{R}$ are
Lipschitzian with the constants $L,K>0,$ then 
\begin{equation}
\left\vert C\left( f,g;\mathbf{A,x}\right) \right\vert \leq LKC\left( e,e;%
\mathbf{A,x}\right) ,  \label{II.op2funce.3.7}
\end{equation}%
for any $\mathbf{x}=\left( x_{1},\dots ,x_{n}\right) \in H^{n}$ with $%
\sum_{j=1}^{n}\left\Vert x_{j}\right\Vert ^{2}=1.$
\end{theorem}

\begin{corollary}
\label{II.op2funcc.3.2}Let $\mathbf{A}=\left( A_{1},\dots ,A_{n}\right) $ be
a sequence of selfadjoint operators with $Sp\left( A_{j}\right) \subseteq %
\left[ m,M\right] $ for $j\in \left\{ 1,\dots ,n\right\} $ and for some
scalars $m<M.$ If $f,g:\left[ m,M\right] \longrightarrow \mathbb{R}$ are
Lipschitzian with the constants $L,K>0,$ then for any $p_{j}\geq 0,j\in
\left\{ 1,\dots ,n\right\} $ with $\sum_{j=1}^{n}p_{j}=1$ we have 
\begin{equation}
\left\vert C\left( f,g;\mathbf{A,p,}x\right) \right\vert \leq LKC\left( e,e;%
\mathbf{A,p,}x\right) ,  \label{II.op2funce.3.8}
\end{equation}%
for any $x\in H$ with $\left\Vert x\right\Vert =1.$
\end{corollary}

\subsection{The Case of $\left( \protect\varphi ,\Phi \right) -$Lipschitzian
Functions}

The following lemma may be stated.

\begin{lemma}
\label{II.op2funcl.1}Let $u:\left[ a,b\right] \rightarrow \mathbb{R}$ and $%
\varphi ,\Phi \in \mathbb{R}$ with $\Phi >\varphi .$ The following
statements are equivalent:

\begin{enumerate}
\item[(i)] The function $u-\frac{\varphi +\Phi }{2}\cdot e,$ where $e\left(
t\right) =t,$ $t\in \left[ a,b\right] ,$ is $\frac{1}{2}\left( \Phi -\varphi
\right) -$Lipschitzian;

\item[(ii)] We have the inequality: 
\begin{equation}
\varphi \leq \frac{u\left( t\right) -u\left( s\right) }{t-s}\leq \Phi \quad 
\text{for each}\quad t,s\in \left[ a,b\right] \quad \text{with }t\neq s;
\label{II.op2funce.4.1}
\end{equation}

\item[(iii)] We have the inequality: 
\begin{equation}
\varphi \left( t-s\right) \leq u\left( t\right) -u\left( s\right) \leq \Phi
\left( t-s\right) \quad \text{for each}\quad t,s\in \left[ a,b\right] \quad 
\text{with }t>s.  \label{II.op2funce.4.2}
\end{equation}
\end{enumerate}
\end{lemma}

Following \cite{II.L}, we can introduce the concept:

\begin{definition}
\label{II.op2funcd.1}The function $u:\left[ a,b\right] \rightarrow \mathbb{R}
$ which satisfies one of the equivalent conditions (i) -- (iii) is said to
be $\left( \varphi ,\Phi \right) -$Lipschitzian on $\left[ a,b\right] .$
\end{definition}

Notice that in \cite{II.L}, the definition was introduced on utilising the
statement (iii) and only the equivalence (i) $\Leftrightarrow $ (iii) was
considered.

Utilising \textit{Lagrange's mean value theorem}, we can state the following
result that provides practical examples of $\left( \varphi ,\Phi \right) -$%
Lipschitzian functions.

\begin{proposition}
\label{II.op2funcp.1}Let $u:\left[ a,b\right] \rightarrow \mathbb{R}$ be
continuous on $\left[ a,b\right] $ and differentiable on $\left( a,b\right)
. $ If 
\begin{equation}
-\infty <\gamma :=\inf_{t\in \left( a,b\right) }u^{\prime }\left( t\right)
,\qquad \sup_{t\in \left( a,b\right) }u^{\prime }\left( t\right) =:\Gamma
<\infty  \label{II.op2funce.4.3}
\end{equation}%
then $u$ is $\left( \gamma ,\Gamma \right) -$Lipschitzian on $\left[ a,b%
\right] .$
\end{proposition}

The following result can be stated:

\begin{theorem}[Dragomir, 2008, \protect\cite{II.SSDG3}]
\label{II.op2funct.4.1}Let $A$ be a selfadjoint operator with $Sp\left(
A\right) \subseteq \left[ m,M\right] $ for some real numbers $m<M.$ If $f:%
\left[ m,M\right] \longrightarrow \mathbb{R}$ is $\left( \varphi ,\Phi
\right) -$Lipschitzian on $\left[ a,b\right] $ and $g:\left[ m,M\right]
\longrightarrow \mathbb{R}$ is continuous with $\delta :=\min_{t\in \left[
m,M\right] }g\left( t\right) $ and $\Delta :=\max_{t\in \left[ m,M\right]
}g\left( t\right) ,$ then 
\begin{align}
\left\vert C\left( f,g;A;x\right) -\frac{\varphi +\Phi }{2}C\left(
e,g;A;x\right) \right\vert & \leq \frac{1}{4}\left( \Delta -\delta \right)
\left( \Phi -\varphi \right) \left\langle \ell _{A,x}\left( A\right)
x,x\right\rangle  \label{II.op2funce.4.4} \\
& \leq \frac{1}{2}\left( \Delta -\delta \right) \left( \Phi -\varphi \right)
C\left( e,e;A;x\right)  \notag
\end{align}%
for any $x\in H$ with $\left\Vert x\right\Vert =1.$
\end{theorem}

The proof follows by Theorem \ref{II.op2funct.2.1} applied for the $\frac{1}{%
2}\left( \Phi -\varphi \right) -$Lipschitzian function $f-\frac{\varphi
+\Phi }{2}\cdot e$ (see Lemma \ref{II.op2funcl.1}) and the details are
omitted.

\begin{theorem}[Dragomir, 2008, \protect\cite{II.SSDG3}]
\label{II.op2funct.4.2}Let $A$ be a selfadjoint operator with $Sp\left(
A\right) \subseteq \left[ m,M\right] $ for some real numbers $m<M$ and $f,g:%
\left[ m,M\right] \longrightarrow \mathbb{R}$. If $f$ is $\left( \varphi
,\Phi \right) -$Lipschitzian and $g$ is $\left( \psi ,\Psi \right) -$%
Lipschitzian on $\left[ a,b\right] ,$ then 
\begin{align}
& \left\vert C\left( f,g;A;x\right) -\frac{\Phi +\varphi }{2}C\left(
e,g;A;x\right) \right.  \label{II.op2funce.4.5} \\
& \left. -\frac{\Psi +\psi }{2}C\left( f,e;A;x\right) +\frac{\Phi +\varphi }{%
2}\cdot \frac{\Psi +\psi }{2}C\left( e,e;A;x\right) \right\vert  \notag \\
& \leq \frac{1}{4}\left( \Phi -\varphi \right) \left( \Psi -\psi \right)
C\left( e,e;A;x\right) ,  \notag
\end{align}%
for any $x\in H$ with $\left\Vert x\right\Vert =1.$
\end{theorem}

The proof follows by Theorem \ref{II.op2funct.2.2} applied for the $\frac{1}{%
2}\left( \Phi -\varphi \right) -$Lipschitzian function $f-\frac{\varphi
+\Phi }{2}\cdot e$ and the $\frac{1}{2}\left( \Psi -\psi \right) -$%
Lipschitzian function $g-\frac{\Psi +\psi }{2}\cdot e.$ The details are
omitted.

Similar results can be derived for sequences of operators, however they will
not be presented here.

\subsection{Some Applications}

It is clear that all the inequalities obtained in the previous sections can
be applied to obtain particular inequalities of interest for different
selections of the functions $f$ and $g$ involved. However we will present
here only some particular results that can be derived from the inequality 
\begin{equation}
\left\vert C\left( f,g;A;x\right) \right\vert \leq LKC\left( e,e;A;x\right) ,
\label{II.op2funce.5.1}
\end{equation}%
that holds for the Lipschitzian functions $f$ and $g,$ the first with the
constant $L>0$ and the second with the constant $K>0.$

\textbf{1.} Now, if we consider the functions $f,g:\left[ m,M\right] \subset
\left( 0,\infty \right) \rightarrow \mathbb{R}$ with $f\left( t\right)
=t^{p},g\left( t\right) =t^{q}$ and $p,q\in \left( -\infty ,0\right) \cup
\left( 0,\infty \right) $ then they are Lipschitzian with the constants $%
L=\left\Vert f^{\prime }\right\Vert _{\infty }$ and $K=\left\Vert g^{\prime
}\right\Vert _{\infty }.$ Since $f^{\prime }\left( t\right)
=pt^{p-1},g\left( t\right) =qt^{q-1},$ hence 
\begin{equation*}
\left\Vert f^{\prime }\right\Vert _{\infty }=\left\{ 
\begin{array}{ll}
pM^{p-1} & \text{for }p\in \lbrack 1,\infty ), \\ 
&  \\ 
\left\vert p\right\vert m^{p-1} & \text{for }p\in \left( -\infty ,0\right)
\cup \left( 0,1\right)%
\end{array}%
\right. \text{ }
\end{equation*}%
and 
\begin{equation*}
\left\Vert g^{\prime }\right\Vert _{\infty }=\left\{ 
\begin{array}{ll}
qM^{q-1} & \text{for }q\in \lbrack 1,\infty ), \\ 
&  \\ 
\left\vert q\right\vert m^{q-1} & \text{for }q\in \left( -\infty ,0\right)
\cup \left( 0,1\right)%
\end{array}%
\right. .
\end{equation*}

Therefore we can state the following inequalities for the powers of a
positive definite operator $A$ with $Sp\left( A\right) \subset \left[ m,M%
\right] \subset \left( 0,\infty \right) .$

If $p,q\geq 1,$ then 
\begin{equation}
(0\leq )\left\langle A^{p+q}x,x\right\rangle -\left\langle
A^{p}x,x\right\rangle \cdot \left\langle A^{q}x,x\right\rangle \leq
pqM^{p+q-2}\left( \left\Vert Ax\right\Vert ^{2}-\left\langle
Ax,x\right\rangle ^{2}\right)  \label{II.op2funce.5.2}
\end{equation}%
for each $x\in H$ with $\left\Vert x\right\Vert =1.$

If $p\geq 1$ and $q\in \left( -\infty ,0\right) \cup \left( 0,1\right) ,$
then 
\begin{equation}
\left\vert \left\langle A^{p+q}x,x\right\rangle -\left\langle
A^{p}x,x\right\rangle \cdot \left\langle A^{q}x,x\right\rangle \right\vert
\leq p\left\vert q\right\vert M^{p-1}m^{q-1}\left( \left\Vert Ax\right\Vert
^{2}-\left\langle Ax,x\right\rangle ^{2}\right)  \label{II.op2funce.5.3}
\end{equation}%
for each $x\in H$ with $\left\Vert x\right\Vert =1.$

If $p\in \left( -\infty ,0\right) \cup \left( 0,1\right) $ and $q\geq 1,$
then 
\begin{equation}
\left\vert \left\langle A^{p+q}x,x\right\rangle -\left\langle
A^{p}x,x\right\rangle \cdot \left\langle A^{q}x,x\right\rangle \right\vert
\leq \left\vert p\right\vert qM^{q-1}m^{p-1}\left( \left\Vert Ax\right\Vert
^{2}-\left\langle Ax,x\right\rangle ^{2}\right)  \label{II.op2funce.5.4}
\end{equation}%
for each $x\in H$ with $\left\Vert x\right\Vert =1.$

If $p,q\in \left( -\infty ,0\right) \cup \left( 0,1\right) ,$ then 
\begin{equation}
\left\vert \left\langle A^{p+q}x,x\right\rangle -\left\langle
A^{p}x,x\right\rangle \cdot \left\langle A^{q}x,x\right\rangle \right\vert
\leq \left\vert pq\right\vert m^{p+q-2}\left( \left\Vert Ax\right\Vert
^{2}-\left\langle Ax,x\right\rangle ^{2}\right)  \label{II.op2funce.5.5}
\end{equation}%
for each $x\in H$ with $\left\Vert x\right\Vert =1.$

Moreover, if we take $p=1$ and $q=-1$ in (\ref{II.op2funce.5.3}), then we
get the following result 
\begin{equation}
(0\leq )\left\langle Ax,x\right\rangle \cdot \left\langle
A^{-1}x,x\right\rangle -1\leq m^{-2}\left( \left\Vert Ax\right\Vert
^{2}-\left\langle Ax,x\right\rangle ^{2}\right)  \label{II.op2funce.5.6}
\end{equation}%
for each $x\in H$ with $\left\Vert x\right\Vert =1.$

\textbf{2.} Consider now the functions $f,g:\left[ m,M\right] \subset \left(
0,\infty \right) \rightarrow \mathbb{R}$ with $f\left( t\right) =t^{p},p\in
\left( -\infty ,0\right) \cup \left( 0,\infty \right) $ and $g\left(
t\right) =\ln t.$ Then $g$ is also Lipschitzian with the constant $%
K=\left\Vert g^{\prime }\right\Vert _{\infty }=m^{-1}.$ Applying the
inequality (\ref{II.op2funce.5.1}) we then have for any $x\in H$ with $%
\left\Vert x\right\Vert =1$ that 
\begin{equation}
(0\leq )\left\langle A^{p}\ln Ax,x\right\rangle -\left\langle
A^{p}x,x\right\rangle \cdot \left\langle \ln Ax,x\right\rangle \leq
pM^{p-1}m^{-1}\left( \left\Vert Ax\right\Vert ^{2}-\left\langle
Ax,x\right\rangle ^{2}\right)  \label{II.op2funce.5.7}
\end{equation}%
if $p\geq 1,$%
\begin{equation}
(0\leq )\left\langle A^{p}\ln Ax,x\right\rangle -\left\langle
A^{p}x,x\right\rangle \cdot \left\langle \ln Ax,x\right\rangle \leq
pm^{p-2}\left( \left\Vert Ax\right\Vert ^{2}-\left\langle Ax,x\right\rangle
^{2}\right)  \label{II.op2funce.5.8}
\end{equation}%
if $p\in \left( 0,1\right) $ and 
\begin{equation}
(0\leq )\left\langle A^{p}x,x\right\rangle \cdot \left\langle \ln
Ax,x\right\rangle -\left\langle A^{p}\ln Ax,x\right\rangle \leq \left(
-p\right) m^{p-2}\left( \left\Vert Ax\right\Vert ^{2}-\left\langle
Ax,x\right\rangle ^{2}\right)  \label{II.op2funce.5.9}
\end{equation}%
if $p\in \left( -\infty ,0\right) .$

\textbf{3. }Now consider the functions $f,g:\left[ m,M\right] \subset 
\mathbb{R\rightarrow R}$ given by $f\left( t\right) =\exp \left( \alpha
t\right) $ and $g\left( t\right) =\exp \left( \beta t\right) $ with $\alpha
,\beta $ nonzero real numbers. It is obvious that%
\begin{equation*}
\left\Vert f^{\prime }\right\Vert _{\infty }=\left\vert \alpha \right\vert
\times \left\{ 
\begin{array}{ll}
\exp \left( \alpha M\right) & \text{for }\alpha >0, \\ 
&  \\ 
\exp \left( \alpha m\right) & \text{for }\alpha <0%
\end{array}%
\right. \text{ }
\end{equation*}%
and 
\begin{equation*}
\left\Vert g^{\prime }\right\Vert _{\infty }=\left\vert \beta \right\vert
\times \left\{ 
\begin{array}{ll}
\exp \left( \beta M\right) & \text{for }\beta >0, \\ 
&  \\ 
\exp \left( \beta m\right) & \text{for }\beta <0%
\end{array}%
\right. .
\end{equation*}

Finally, on applying the inequality (\ref{II.op2funce.5.1}) we get%
\begin{align*}
(0& \leq )\left\langle \exp \left[ \left( \alpha +\beta \right) A\right]
x,x\right\rangle -\left\langle \exp \left( \alpha A\right) x,x\right\rangle
\cdot \left\langle \exp \left( \beta A\right) x,x\right\rangle \\
& \leq \left\vert \alpha \beta \right\vert \left( \left\Vert Ax\right\Vert
^{2}-\left\langle Ax,x\right\rangle ^{2}\right) \times \left\{ 
\begin{array}{ll}
\exp \left[ \left( \alpha +\beta \right) M\right] & \text{for }\alpha ,\beta
>0, \\ 
&  \\ 
\exp \left[ \left( \alpha +\beta \right) m\right] & \text{for }\alpha ,\beta
<0%
\end{array}%
\right. \text{ }
\end{align*}%
and%
\begin{align*}
(0& \leq )\left\langle \exp \left( \alpha A\right) x,x\right\rangle \cdot
\left\langle \exp \left( \beta A\right) x,x\right\rangle -\left\langle \exp 
\left[ \left( \alpha +\beta \right) A\right] x,x\right\rangle \\
& \leq \left\vert \alpha \beta \right\vert \left( \left\Vert Ax\right\Vert
^{2}-\left\langle Ax,x\right\rangle ^{2}\right) \times \left\{ 
\begin{array}{ll}
\exp \left( \alpha M+\beta m\right) & \text{for }\alpha >0,\beta <0 \\ 
&  \\ 
\exp \left( \alpha m+\beta M\right) & \text{for }\alpha <0,\beta >0%
\end{array}%
\right. \text{ }
\end{align*}%
for each $x\in H$ with $\left\Vert x\right\Vert =1.$

\section{Quasi Gr\"{u}ss' Type Inequalities}

\subsection{Introduction}

In \cite{II.2b}, \ in order to generalize the above result in abstract
structures the author has proved the following Gr\"{u}ss' type inequality in
real or complex inner product spaces.

\begin{theorem}[Dragomir, 1999, \protect\cite{II.2b}]
\label{II.a.t}Let $\left( H,\left\langle .,.\right\rangle \right) $ be an
inner product space over $\mathbb{K}\left( \mathbb{K=R}\text{,}\mathbb{C}%
\right) $ and $e\in H,\left\Vert e\right\Vert =1.$ If $\varphi ,\gamma ,\Phi
,\Gamma $ are real or complex numbers and $x,y$ are vectors in $H$ such that
the conditions 
\begin{equation}
\func{Re}\left\langle \Phi e-x,x-\varphi e\right\rangle \geq 0\text{ and }%
\func{Re}\left\langle \Gamma e-y,y-\gamma e\right\rangle \geq 0
\label{II.a.i.1}
\end{equation}%
hold, then we have the inequality 
\begin{equation}
\left\vert \left\langle x,y\right\rangle -\left\langle x,e\right\rangle
\left\langle e,y\right\rangle \right\vert \leq \frac{1}{4}\left\vert \Phi
-\varphi \right\vert \cdot \left\vert \Gamma -\gamma \right\vert .
\label{II.a.i.2}
\end{equation}%
The constant $\frac{1}{4}$ is best possible in the sense that it can not be
replaced by a smaller constant.
\end{theorem}

For other results of this type, see the recent monograph \cite{II.SSDM} and
the references therein.

Let $U$ be a selfadjoint operator on the complex Hilbert space $\left(
H,\left\langle .,.\right\rangle \right) $ with the spectrum $Sp\left(
U\right) $ included in the interval $\left[ m,M\right] $ for some real
numbers $m<M$ and let $\left\{ E_{\lambda }\right\} _{\lambda }$ be its 
\textit{spectral family}. Then for any continuous function $f:\left[ m,M%
\right] \rightarrow \mathbb{C}$, it is well known that we have the following 
\textit{spectral representation theorem in terms of the Riemann-Stieltjes
integral}: 
\begin{equation}
f\left( U\right) =\int_{m-0}^{M}f\left( \lambda \right) dE_{\lambda },
\label{II.a.e.1.1}
\end{equation}%
which in terms of vectors can be written as 
\begin{equation}
\left\langle f\left( U\right) x,y\right\rangle =\int_{m-0}^{M}f\left(
\lambda \right) d\left\langle E_{\lambda }x,y\right\rangle ,
\label{II.a.e.1.2}
\end{equation}%
for any $x,y\in H.$ The function $g_{x,y}\left( \lambda \right)
:=\left\langle E_{\lambda }x,y\right\rangle $ is of \textit{bounded variation%
} on the interval $\left[ m,M\right] $ and 
\begin{equation*}
g_{x,y}\left( m-0\right) =0\text{ and }g_{x,y}\left( M\right) =\left\langle
x,y\right\rangle
\end{equation*}%
for any $x,y\in H.$ It is also well known that $g_{x}\left( \lambda \right)
:=\left\langle E_{\lambda }x,x\right\rangle $ is \textit{monotonic
nondecreasing} and \textit{right continuous} on $\left[ m,M\right] $.

\subsection{Vector Inequalities}

In this section we provide various bounds for the magnitude of the difference%
\begin{equation*}
\left\langle f\left( A\right) x,y\right\rangle -\left\langle
x,y\right\rangle \left\langle f\left( A\right) x,x\right\rangle
\end{equation*}%
under different assumptions on the continuous function, the selfadjoint
operator $A:H\rightarrow H$ and the vectors $x,y\in H$ with $\left\Vert
x\right\Vert =1.$

\begin{theorem}[Dragomir, 2010, \protect\cite{II.a.SSD2}]
\label{II.a.t.1.1}Let $A$ be a selfadjoint operator in the Hilbert space $H$
with the spectrum $Sp\left( A\right) \subseteq \left[ m,M\right] $ for some
real numbers $m<M$ and let $\left\{ E_{\lambda }\right\} _{\lambda }$ be its 
\textit{spectral family.} Assume that $x,y\in H,\left\Vert x\right\Vert =1$
are such that there exists $\gamma ,\Gamma \in \mathbb{C}$ with either%
\begin{equation}
\func{Re}\left\langle \Gamma x-y,y-\gamma x\right\rangle \geq 0
\label{II.a.C}
\end{equation}%
or, equivalently%
\begin{equation*}
\left\Vert y-\frac{\gamma +\Gamma }{2}x\right\Vert \leq \frac{1}{2}%
\left\vert \Gamma -\gamma \right\vert .
\end{equation*}

1. If $f:\left[ m,M\right] \rightarrow \mathbb{C}$ is a continuous function
of bounded variation on $\left[ m,M\right] $, then we have the inequality%
\begin{align}
& \left\vert \left\langle f\left( A\right) x,y\right\rangle -\left\langle
x,y\right\rangle \left\langle f\left( A\right) x,x\right\rangle \right\vert
\label{II.a.e.3.1} \\
& \leq \max_{\lambda \in \left[ m,M\right] }\left\vert \left\langle
E_{\lambda }x,y\right\rangle -\left\langle E_{\lambda }x,x\right\rangle
\left\langle x,y\right\rangle \right\vert \dbigvee\limits_{m}^{M}\left(
f\right)  \notag \\
& \leq \max_{\lambda \in \left[ m,M\right] }\left( \left\langle E_{\lambda
}x,x\right\rangle \left\langle \left( 1_{H}-E_{\lambda }\right)
x,x\right\rangle \right) ^{1/2}\left( \left\Vert y\right\Vert
^{2}-\left\vert \left\langle y,x\right\rangle \right\vert ^{2}\right)
^{1/2}\dbigvee\limits_{m}^{M}\left( f\right)  \notag \\
& \leq \frac{1}{2}\left( \left\Vert y\right\Vert ^{2}-\left\vert
\left\langle y,x\right\rangle \right\vert ^{2}\right)
^{1/2}\dbigvee\limits_{m}^{M}\left( f\right) \leq \frac{1}{4}\left\vert
\Gamma -\gamma \right\vert \dbigvee\limits_{m}^{M}\left( f\right) .  \notag
\end{align}

2. If $f:\left[ m,M\right] \rightarrow \mathbb{C}$ is a Lipschitzian
function with the constant $L>0$ on $\left[ m,M\right] $, then we have the
inequality%
\begin{align}
& \left\vert \left\langle f\left( A\right) x,y\right\rangle -\left\langle
x,y\right\rangle \left\langle f\left( A\right) x,x\right\rangle \right\vert
\label{II.a.e.3.2} \\
& \leq L\int_{m-0}^{M}\left\vert \left\langle E_{\lambda }x,y\right\rangle
-\left\langle E_{\lambda }x,x\right\rangle \left\langle x,y\right\rangle
\right\vert d\lambda  \notag \\
& \leq L\left( \left\Vert y\right\Vert ^{2}-\left\vert \left\langle
y,x\right\rangle \right\vert ^{2}\right) ^{1/2}\int_{m-0}^{M}\left(
\left\langle E_{\lambda }x,x\right\rangle \left\langle \left(
1_{H}-E_{\lambda }\right) x,x\right\rangle \right) ^{1/2}d\lambda  \notag \\
& \leq L\left( \left\Vert y\right\Vert ^{2}-\left\vert \left\langle
y,x\right\rangle \right\vert ^{2}\right) ^{1/2}\left\langle \left(
M1_{H}-A\right) x,x\right\rangle ^{1/2}\left\langle \left( A-m1_{H}\right)
x,x\right\rangle ^{1/2}  \notag \\
& \leq \frac{1}{2}\left( M-m\right) L\left( \left\Vert y\right\Vert
^{2}-\left\vert \left\langle y,x\right\rangle \right\vert ^{2}\right)
^{1/2}\leq \frac{1}{4}\left\vert \Gamma -\gamma \right\vert \left(
M-m\right) L.  \notag
\end{align}

3. If $f:\left[ m,M\right] \rightarrow \mathbb{R}$ is a continuous monotonic
nondecreasing function on $\left[ m,M\right] $, then we have the inequality%
\begin{align}
& \left\vert \left\langle f\left( A\right) x,y\right\rangle -\left\langle
x,y\right\rangle \left\langle f\left( A\right) x,x\right\rangle \right\vert
\label{II.a.e.3.3} \\
& \leq \int_{m-0}^{M}\left\vert \left\langle E_{\lambda }x,y\right\rangle
-\left\langle E_{\lambda }x,x\right\rangle \left\langle x,y\right\rangle
\right\vert df\left( \lambda \right)  \notag \\
& \leq \left( \left\Vert y\right\Vert ^{2}-\left\vert \left\langle
y,x\right\rangle \right\vert ^{2}\right) ^{1/2}\int_{m-0}^{M}\left(
\left\langle E_{\lambda }x,x\right\rangle \left\langle \left(
1_{H}-E_{\lambda }\right) x,x\right\rangle \right) ^{1/2}df\left( \lambda
\right)  \notag \\
& \leq \left( \left\Vert y\right\Vert ^{2}-\left\vert \left\langle
y,x\right\rangle \right\vert ^{2}\right) ^{1/2}  \notag \\
& \times \left\langle \left( f\left( M\right) 1_{H}-f\left( A\right) \right)
x,x\right\rangle ^{1/2}\left\langle \left( f\left( A\right) -f\left(
m\right) 1_{H}\right) x,x\right\rangle ^{1/2}  \notag \\
& \leq \frac{1}{2}\left[ f\left( M\right) -f\left( m\right) \right] \left(
\left\Vert y\right\Vert ^{2}-\left\vert \left\langle y,x\right\rangle
\right\vert ^{2}\right) ^{1/2}\leq \frac{1}{4}\left\vert \Gamma -\gamma
\right\vert \left[ f\left( M\right) -f\left( m\right) \right] .  \notag
\end{align}
\end{theorem}

\begin{proof}
First of all, we notice that by the Schwarz inequality in $H$ we have for
any $u,v,e\in H$ with $\left\Vert e\right\Vert =1$ that%
\begin{equation}
\left\vert \left\langle u,v\right\rangle -\left\langle u,e\right\rangle
\left\langle e,v\right\rangle \right\vert \leq \left( \left\Vert
u\right\Vert ^{2}-\left\vert \left\langle u,e\right\rangle \right\vert
^{2}\right) ^{1/2}\left( \left\Vert v\right\Vert ^{2}-\left\vert
\left\langle v,e\right\rangle \right\vert ^{2}\right) ^{1/2}.
\label{II.a.e.3.4}
\end{equation}

Now on utilizing (\ref{II.a.e.3.4}), we can state that%
\begin{align}
& \left\vert \left\langle E_{\lambda }x,y\right\rangle -\left\langle
E_{\lambda }x,x\right\rangle \left\langle x,y\right\rangle \right\vert
\label{II.a.e.3.5} \\
& \leq \left( \left\Vert E_{\lambda }x\right\Vert ^{2}-\left\vert
\left\langle E_{\lambda }x,x\right\rangle \right\vert ^{2}\right)
^{1/2}\left( \left\Vert y\right\Vert ^{2}-\left\vert \left\langle
y,x\right\rangle \right\vert ^{2}\right) ^{1/2}  \notag
\end{align}%
for any $\lambda \in \left[ m,M\right] .$

Since $E_{\lambda }$ are projections and $E_{\lambda }\geq 0$ then 
\begin{align}
\left\Vert E_{\lambda }x\right\Vert ^{2}-\left\vert \left\langle E_{\lambda
}x,x\right\rangle \right\vert ^{2}& =\left\langle E_{\lambda
}x,x\right\rangle -\left\langle E_{\lambda }x,x\right\rangle ^{2}
\label{II.a.e.3.6} \\
& =\left\langle E_{\lambda }x,x\right\rangle \left\langle \left(
1_{H}-E_{\lambda }\right) x,x\right\rangle \leq \frac{1}{4}  \notag
\end{align}%
for any $\lambda \in \left[ m,M\right] $ and $x\in H$ with $\left\Vert
x\right\Vert =1.$

Also, by making use of the Gr\"{u}ss' type inequality in inner product
spaces obtained by the author in \cite{II.2b} we have%
\begin{equation}
\left( \left\Vert y\right\Vert ^{2}-\left\vert \left\langle y,x\right\rangle
\right\vert ^{2}\right) ^{1/2}\leq \frac{1}{2}\left\vert \Gamma -\gamma
\right\vert .  \label{II.a.e.3.7}
\end{equation}%
Combining the relations (\ref{II.a.e.3.5})-(\ref{II.a.e.3.7}) we deduce the
following inequality that is of interest in itself%
\begin{align}
& \left\vert \left\langle E_{\lambda }x,y\right\rangle -\left\langle
E_{\lambda }x,x\right\rangle \left\langle x,y\right\rangle \right\vert
\label{II.a.e.3.8} \\
& \leq \left( \left\langle E_{\lambda }x,x\right\rangle \left\langle \left(
1_{H}-E_{\lambda }\right) x,x\right\rangle \right) ^{1/2}\left( \left\Vert
y\right\Vert ^{2}-\left\vert \left\langle y,x\right\rangle \right\vert
^{2}\right) ^{1/2}  \notag \\
& \leq \frac{1}{2}\left( \left\Vert y\right\Vert ^{2}-\left\vert
\left\langle y,x\right\rangle \right\vert ^{2}\right) ^{1/2}\leq \frac{1}{4}%
\left\vert \Gamma -\gamma \right\vert  \notag
\end{align}%
for any $\lambda \in \left[ m,M\right] .$

It is well known that if $p:\left[ a,b\right] \rightarrow \mathbb{C}$ is a
continuous function, $v:\left[ a,b\right] \rightarrow \mathbb{C}$ is of
bounded variation then the Riemann-Stieltjes integral $\int_{a}^{b}p\left(
t\right) dv\left( t\right) $ exists and the following inequality holds%
\begin{equation}
\left\vert \int_{a}^{b}p\left( t\right) dv\left( t\right) \right\vert \leq
\max_{t\in \left[ a,b\right] }\left\vert p\left( t\right) \right\vert
\dbigvee\limits_{a}^{b}\left( v\right) ,  \label{II.a.e.3.9}
\end{equation}%
where $\dbigvee\limits_{a}^{b}\left( v\right) $ denotes the total variation
of $v$ on $\left[ a,b\right] .$

Utilising this property of the Riemann-Stieltjes integral and the inequality
(\ref{II.a.e.3.8}) we have%
\begin{align}
& \left\vert \int_{m-0}^{M}\left[ \left\langle E_{\lambda }x,y\right\rangle
-\left\langle E_{\lambda }x,x\right\rangle \left\langle x,y\right\rangle %
\right] df\left( \lambda \right) \right\vert  \label{II.a.e.3.10} \\
& \leq \max_{\lambda \in \left[ m,M\right] }\left\vert \left\langle
E_{\lambda }x,y\right\rangle -\left\langle E_{\lambda }x,x\right\rangle
\left\langle x,y\right\rangle \right\vert \dbigvee\limits_{m}^{M}\left(
f\right)  \notag \\
& \leq \max_{\lambda \in \left[ m,M\right] }\left( \left\langle E_{\lambda
}x,x\right\rangle \left\langle \left( 1_{H}-E_{\lambda }\right)
x,x\right\rangle \right) ^{1/2}\left( \left\Vert y\right\Vert
^{2}-\left\vert \left\langle y,x\right\rangle \right\vert ^{2}\right)
^{1/2}\dbigvee\limits_{m}^{M}\left( f\right)  \notag \\
& \leq \frac{1}{2}\left( \left\Vert y\right\Vert ^{2}-\left\vert
\left\langle y,x\right\rangle \right\vert ^{2}\right)
^{1/2}\dbigvee\limits_{m}^{M}\left( f\right) \leq \frac{1}{4}\left\vert
\Gamma -\gamma \right\vert \dbigvee\limits_{m}^{M}\left( f\right)  \notag
\end{align}%
for $x$ and $y$ as in the assumptions of the theorem.

Now, integrating by parts in the Riemann-Stieltjes integral and making use
of the spectral representation theorem we have%
\begin{align}
& \int_{m-0}^{M}\left[ \left\langle E_{\lambda }x,y\right\rangle
-\left\langle E_{\lambda }x,x\right\rangle \left\langle x,y\right\rangle %
\right] df\left( \lambda \right)  \label{II.a.e.3.11} \\
& =\left. \left[ \left\langle E_{\lambda }x,y\right\rangle -\left\langle
E_{\lambda }x,x\right\rangle \left\langle x,y\right\rangle \right] f\left(
\lambda \right) \right\vert _{m-0}^{M}  \notag \\
& -\int_{m-0}^{M}f\left( \lambda \right) d\left[ \left\langle E_{\lambda
}x,y\right\rangle -\left\langle E_{\lambda }x,x\right\rangle \left\langle
x,y\right\rangle \right]  \notag \\
& =\left\langle x,y\right\rangle \int_{m-0}^{M}f\left( \lambda \right)
d\left\langle E_{\lambda }x,x\right\rangle -\int_{m-0}^{M}f\left( \lambda
\right) d\left\langle E_{\lambda }x,y\right\rangle  \notag \\
& =\left\langle x,y\right\rangle \left\langle f\left( A\right)
x,x\right\rangle -\left\langle f\left( A\right) x,y\right\rangle  \notag
\end{align}%
which together with (\ref{II.a.e.3.10}) produces the desired result (\ref%
{II.a.e.3.1}).

Now, recall that if $p:\left[ a,b\right] \rightarrow \mathbb{C}$ is a
Riemann integrable function and $v:\left[ a,b\right] \rightarrow \mathbb{C}$
is Lipschitzian with the constant $L>0$, i.e.,%
\begin{equation*}
\left\vert f\left( s\right) -f\left( t\right) \right\vert \leq L\left\vert
s-t\right\vert \text{ for any }t,s\in \left[ a,b\right] ,
\end{equation*}%
then the Riemann-Stieltjes integral $\int_{a}^{b}p\left( t\right) dv\left(
t\right) $ exists and the following inequality holds%
\begin{equation*}
\left\vert \int_{a}^{b}p\left( t\right) dv\left( t\right) \right\vert \leq
L\int_{a}^{b}\left\vert p\left( t\right) \right\vert dt.
\end{equation*}

Now, on applying this property of the Riemann-Stieltjes integral we have
from (\ref{II.a.e.3.8}) that 
\begin{align}
& \left\vert \int_{m-0}^{M}\left[ \left\langle E_{\lambda }x,y\right\rangle
-\left\langle E_{\lambda }x,x\right\rangle \left\langle x,y\right\rangle %
\right] df\left( \lambda \right) \right\vert  \label{II.a.e.3.12} \\
& \leq L\int_{m-0}^{M}\left\vert \left\langle E_{\lambda }x,y\right\rangle
-\left\langle E_{\lambda }x,x\right\rangle \left\langle x,y\right\rangle
\right\vert d\lambda  \notag \\
& \leq L\left( \left\Vert y\right\Vert ^{2}-\left\vert \left\langle
y,x\right\rangle \right\vert ^{2}\right) ^{1/2}\int_{m-0}^{M}\left(
\left\langle E_{\lambda }x,x\right\rangle \left\langle \left(
1_{H}-E_{\lambda }\right) x,x\right\rangle \right) ^{1/2}d\lambda .  \notag
\end{align}

If we use the Cauchy-Bunyakovsky-Schwarz integral inequality and the
spectral representation theorem we have successively 
\begin{align}
& \int_{m-0}^{M}\left( \left\langle E_{\lambda }x,x\right\rangle
\left\langle \left( 1_{H}-E_{\lambda }\right) x,x\right\rangle \right)
^{1/2}d\lambda  \label{II.a.e.3.13} \\
& \leq \left[ \int_{m-0}^{M}\left\langle E_{\lambda }x,x\right\rangle
d\lambda \right] ^{1/2}\left[ \int_{m-0}^{M}\left\langle \left(
1_{H}-E_{\lambda }\right) x,x\right\rangle d\lambda \right] ^{1/2}  \notag \\
& =\left[ \left. \left\langle E_{\lambda }x,x\right\rangle \lambda
\right\vert _{m-0}^{M}-\int_{m-0}^{M}\lambda d\left\langle E_{\lambda
}x,x\right\rangle \right] ^{1/2}  \notag \\
& \times \left[ \left. \left\langle \left( 1_{H}-E_{\lambda }\right)
x,x\right\rangle \lambda \right\vert _{m-0}^{M}-\int_{m-0}^{M}\lambda
d\left\langle \left( 1_{H}-E_{\lambda }\right) x,x\right\rangle \right] 
\notag \\
& =\left\langle \left( M1_{H}-A\right) x,x\right\rangle ^{1/2}\left\langle
\left( A-m1_{H}\right) x,x\right\rangle ^{1/2}.  \notag
\end{align}%
On utilizing (\ref{II.a.e.3.13}), (\ref{II.a.e.3.12}) and (\ref{II.a.e.3.11}%
) we deduce the first three\ inequalities in (\ref{II.a.e.3.2}).

The fourth inequality follows from the fact that 
\begin{align*}
& \left\langle \left( M1_{H}-A\right) x,x\right\rangle \left\langle \left(
A-m1_{H}\right) x,x\right\rangle \\
& \leq \frac{1}{4}\left[ \left\langle \left( M1_{H}-A\right)
x,x\right\rangle +\left\langle \left( A-m1_{H}\right) x,x\right\rangle %
\right] ^{2}=\frac{1}{4}\left( M-m\right) ^{2}.
\end{align*}%
The last part follows from (\ref{II.a.e.3.7}).

Further, from the theory of Riemann-Stieltjes integral it is also well known
that if $p:\left[ a,b\right] \rightarrow \mathbb{C}$ is of bounded variation
and $v:\left[ a,b\right] \rightarrow \mathbb{R}$ is continuous and monotonic
nondecreasing, then the Riemann-Stieltjes integrals $\int_{a}^{b}p\left(
t\right) dv\left( t\right) $ and $\int_{a}^{b}\left\vert p\left( t\right)
\right\vert dv\left( t\right) $ exist and%
\begin{equation}
\left\vert \int_{a}^{b}p\left( t\right) dv\left( t\right) \right\vert \leq
\int_{a}^{b}\left\vert p\left( t\right) \right\vert dv\left( t\right) .
\label{II.a.e.3.14}
\end{equation}%
Utilising this property and the inequality (\ref{II.a.e.3.8}) we have
successively 
\begin{align}
& \left\vert \int_{m-0}^{M}\left[ \left\langle E_{\lambda }x,y\right\rangle
-\left\langle E_{\lambda }x,x\right\rangle \left\langle x,y\right\rangle %
\right] df\left( \lambda \right) \right\vert  \label{II.a.e.3.15} \\
& \leq \int_{m-0}^{M}\left\vert \left\langle E_{\lambda }x,y\right\rangle
-\left\langle E_{\lambda }x,x\right\rangle \left\langle x,y\right\rangle
\right\vert df\left( \lambda \right)  \notag \\
& \leq \left( \left\Vert y\right\Vert ^{2}-\left\vert \left\langle
y,x\right\rangle \right\vert ^{2}\right) ^{1/2}\int_{m-0}^{M}\left(
\left\langle E_{\lambda }x,x\right\rangle \left\langle \left(
1_{H}-E_{\lambda }\right) x,x\right\rangle \right) ^{1/2}df\left( \lambda
\right) .  \notag
\end{align}%
Applying the Cauchy-Bunyakovsky-Schwarz integral inequality for the
Riemann-Stieltjes integral with monotonic integrators and the spectral
representation theorem we have 
\begin{align}
& \int_{m-0}^{M}\left( \left\langle E_{\lambda }x,x\right\rangle
\left\langle \left( 1_{H}-E_{\lambda }\right) x,x\right\rangle \right)
^{1/2}df\left( \lambda \right)  \label{II.a.e.3.16} \\
& \leq \left[ \int_{m-0}^{M}\left\langle E_{\lambda }x,x\right\rangle
df\left( \lambda \right) \right] ^{1/2}\left[ \int_{m-0}^{M}\left\langle
\left( 1_{H}-E_{\lambda }\right) x,x\right\rangle df\left( \lambda \right) %
\right] ^{1/2}  \notag \\
& =\left[ \left. \left\langle E_{\lambda }x,x\right\rangle f\left( \lambda
\right) \right\vert _{m-0}^{M}-\int_{m-0}^{M}f\left( \lambda \right)
d\left\langle E_{\lambda }x,x\right\rangle \right] ^{1/2}  \notag \\
& \times \left[ \left. \left\langle \left( 1_{H}-E_{\lambda }\right)
x,x\right\rangle f\left( \lambda \right) \right\vert
_{m-0}^{M}-\int_{m-0}^{M}f\left( \lambda \right) d\left\langle \left(
1_{H}-E_{\lambda }\right) x,x\right\rangle \right] ^{1/2}  \notag \\
& =\left\langle \left( f\left( M\right) 1_{H}-f\left( A\right) \right)
x,x\right\rangle ^{1/2}\left\langle \left( f\left( A\right) -f\left(
m\right) 1_{H}\right) x,x\right\rangle ^{1/2}  \notag \\
& \leq \frac{1}{2}\left[ f\left( M\right) -f\left( m\right) \right]  \notag
\end{align}%
and the proof is complete.
\end{proof}

\begin{remark}
\label{II.a.r.3.1}If we drop the conditions on $x,y,$ we can obtain from the
inequalities (\ref{II.a.e.3.1})-(\ref{II.a.e.3.2}) the following results
that can be easily applied for particular functions:

1. If $f:\left[ m,M\right] \rightarrow \mathbb{C}$ is a continuous function
of bounded variation on $\left[ m,M\right] $, then we have the inequality%
\begin{align}
& \left\vert \left\langle f\left( A\right) x,y\right\rangle \left\Vert
x\right\Vert ^{2}-\left\langle x,y\right\rangle \left\langle f\left(
A\right) x,x\right\rangle \right\vert  \label{II.a.e.3.17} \\
& \leq \frac{1}{2}\left\Vert x\right\Vert ^{2}\left( \left\Vert y\right\Vert
^{2}\left\Vert x\right\Vert ^{2}-\left\vert \left\langle y,x\right\rangle
\right\vert ^{2}\right) ^{1/2}\dbigvee\limits_{m}^{M}\left( f\right)  \notag
\end{align}%
for any $x,y\in H,x\neq 0.$

2. If $f:\left[ m,M\right] \rightarrow \mathbb{C}$ is a Lipschitzian
function with the constant $L>0$ on $\left[ m,M\right] $, then we have the
inequality%
\begin{align}
& \left\vert \left\langle f\left( A\right) x,y\right\rangle \left\Vert
x\right\Vert ^{2}-\left\langle x,y\right\rangle \left\langle f\left(
A\right) x,x\right\rangle \right\vert \leq L\left( \left\Vert y\right\Vert
^{2}\left\Vert x\right\Vert ^{2}-\left\vert \left\langle y,x\right\rangle
\right\vert ^{2}\right) ^{1/2}  \label{II.a.e.3.18} \\
& \times \left[ \left\langle \left( M1_{H}-A\right) x,x\right\rangle
\left\langle \left( A-m1_{H}\right) x,x\right\rangle \right] ^{1/2}  \notag
\\
& \leq \frac{1}{2}\left( M-m\right) L\left\Vert x\right\Vert ^{2}\left(
\left\Vert y\right\Vert ^{2}\left\Vert x\right\Vert ^{2}-\left\vert
\left\langle y,x\right\rangle \right\vert ^{2}\right) ^{1/2}  \notag
\end{align}%
for any $x,y\in H,x\neq 0.$

3. If $f:\left[ m,M\right] \rightarrow \mathbb{R}$ is a continuous monotonic
nondecreasing function on $\left[ m,M\right] $, then we have the inequality%
\begin{align}
& \left\vert \left\langle f\left( A\right) x,y\right\rangle \left\Vert
x\right\Vert ^{2}-\left\langle x,y\right\rangle \left\langle f\left(
A\right) x,x\right\rangle \right\vert \leq \left( \left\Vert y\right\Vert
^{2}\left\Vert x\right\Vert ^{2}-\left\vert \left\langle y,x\right\rangle
\right\vert ^{2}\right) ^{1/2}  \label{II.a.e.3.19} \\
& \times \left[ \left\langle \left( f\left( M\right) 1_{H}-f\left( A\right)
\right) x,x\right\rangle \left\langle \left( f\left( A\right) -f\left(
m\right) 1_{H}\right) x,x\right\rangle \right] ^{1/2}  \notag \\
& \leq \frac{1}{2}\left[ f\left( M\right) -f\left( m\right) \right]
\left\Vert x\right\Vert ^{2}\left( \left\Vert y\right\Vert ^{2}\left\Vert
x\right\Vert ^{2}-\left\vert \left\langle y,x\right\rangle \right\vert
^{2}\right) ^{1/2}  \notag
\end{align}%
for any $x,y\in H,x\neq 0.$
\end{remark}

The following lemma may be stated.

\begin{lemma}
\label{II.a.l.2.1}Let $u:\left[ a,b\right] \rightarrow \mathbb{R}$ and $%
\varphi ,\Phi \in \mathbb{R}$ with $\Phi >\varphi .$ The following
statements are equivalent:

\begin{enumerate}
\item[(i)] The function $u-\frac{\varphi +\Phi }{2}\cdot e,$ where $e\left(
t\right) =t,$ $t\in \left[ a,b\right] ,$ is $\frac{1}{2}\left( \Phi -\varphi
\right) -$Lipschitzian;

\item[(ii)] We have the inequality: 
\begin{equation}
\varphi \leq \frac{u\left( t\right) -u\left( s\right) }{t-s}\leq \Phi \quad 
\text{for each}\quad t,s\in \left[ a,b\right] \quad \text{with }t\neq s;
\label{II.a.e.2.10}
\end{equation}

\item[(iii)] We have the inequality: 
\begin{equation}
\varphi \left( t-s\right) \leq u\left( t\right) -u\left( s\right) \leq \Phi
\left( t-s\right) \quad \text{for each}\quad t,s\in \left[ a,b\right] \quad 
\text{with }t>s.  \label{II.a.e.2.11}
\end{equation}
\end{enumerate}
\end{lemma}

Following \cite{II.L}, we can introduce the concept:

\begin{definition}
\label{II.a.d.2.1}The function $u:\left[ a,b\right] \rightarrow \mathbb{R}$
which satisfies one of the equivalent conditions (i) -- (iii) is said to be $%
\left( \varphi ,\Phi \right) -$Lipschitzian on $\left[ a,b\right] .$
\end{definition}

Notice that in \cite{II.L}, the definition was introduced on utilizing the
statement (iii) and only the equivalence (i) $\Leftrightarrow $ (iii) was
considered.

Utilising \textit{Lagrange's mean value theorem}, we can state the following
result that provides practical examples of $\left( \varphi ,\Phi \right) -$%
Lipschitzian functions.

\begin{proposition}
\label{II.a.p.2.1}Let $u:\left[ a,b\right] \rightarrow \mathbb{R}$ be
continuous on $\left[ a,b\right] $ and differentiable on $\left( a,b\right)
. $ If 
\begin{equation}
-\infty <\gamma :=\inf_{t\in \left( a,b\right) }u^{\prime }\left( t\right)
,\qquad \sup_{t\in \left( a,b\right) }u^{\prime }\left( t\right) =:\Gamma
<\infty  \label{II.a.e.2.12}
\end{equation}%
then $u$ is $\left( \gamma ,\Gamma \right) -$Lipschitzian on $\left[ a,b%
\right] .$
\end{proposition}

We are able now to provide the following corollary:

\begin{corollary}[Dragomir, 2010, \protect\cite{II.a.SSD2}]
\label{II.a.c.2.1}With the assumptions of Theorem \ref{II.a.t.1.1} and if $f:%
\left[ m,M\right] \rightarrow \mathbb{R}$ is a $\left( \varphi ,\Phi \right) 
$-Lipschitzian function then we have%
\begin{align}
& \left\vert \left\langle f\left( A\right) x,y\right\rangle -\left\langle
x,y\right\rangle \left\langle f\left( A\right) x,x\right\rangle \right\vert
\label{II.a.e.2.13} \\
& \leq \frac{1}{2}\left( \Phi -\varphi \right) \int_{m-0}^{M}\left\vert
\left\langle E_{\lambda }x,y\right\rangle -\left\langle E_{\lambda
}x,x\right\rangle \left\langle x,y\right\rangle \right\vert d\lambda  \notag
\\
& \leq \frac{1}{2}\left( \Phi -\varphi \right) \left( \left\Vert
y\right\Vert ^{2}-\left\vert \left\langle y,x\right\rangle \right\vert
^{2}\right) ^{1/2}\int_{m-0}^{M}\left( \left\langle E_{\lambda
}x,x\right\rangle \left\langle \left( 1_{H}-E_{\lambda }\right)
x,x\right\rangle \right) ^{1/2}d\lambda  \notag \\
& \leq \frac{1}{2}\left( \Phi -\varphi \right) \left( \left\Vert
y\right\Vert ^{2}-\left\vert \left\langle y,x\right\rangle \right\vert
^{2}\right) ^{1/2}  \notag \\
& \times \left\langle \left( M1_{H}-A\right) x,x\right\rangle
^{1/2}\left\langle \left( A-m1_{H}\right) x,x\right\rangle ^{1/2}  \notag \\
& \leq \frac{1}{4}\left( M-m\right) \left( \Phi -\varphi \right) \left(
\left\Vert y\right\Vert ^{2}-\left\vert \left\langle y,x\right\rangle
\right\vert ^{2}\right) ^{1/2}  \notag \\
& \leq \frac{1}{8}\left\vert \Gamma -\gamma \right\vert \left( M-m\right)
\left( \Phi -\varphi \right) .  \notag
\end{align}
\end{corollary}

The proof follows from the second part of Theorem \ref{II.a.t.1.1} applied
for the $\frac{1}{2}\left( \Phi -\varphi \right) $-Lipschitzian function $f-%
\frac{\Phi +\varphi }{2}\cdot e$ by performing the required calculations in
the first term of the inequality. The details are omitted.

\subsection{Applications for Gr\"{u}ss' Type Inequalities}

The following result provides some Gr\"{u}ss' type inequalities for two
function of two selfadjoint operators.

\begin{proposition}[Dragomir, 2010, \protect\cite{II.a.SSD2}]
\label{II.a.p.4.1}Let $A,B$ be two selfadjoint operators in the Hilbert
space $H$ with the spectra $Sp\left( A\right) ,Sp\left( B\right) \subseteq %
\left[ m,M\right] $ for some real numbers $m<M$ and let $\left\{ E_{\lambda
}\right\} _{\lambda }$ be the \textit{spectral family of }$A$\textit{.}
Assume that $g:\left[ m,M\right] \rightarrow \mathbb{R}$ is a continuous
function and denote $n:=\min_{t\in \left[ m,M\right] }g\left( t\right) $ and 
$N:=\max_{t\in \left[ m,M\right] }g\left( t\right) .$

1. If $f:\left[ m,M\right] \rightarrow \mathbb{C}$ is a continuous function
of bounded variation on $\left[ m,M\right] $, then we have the inequality%
\begin{align}
& \left\vert \left\langle f\left( A\right) x,g\left( B\right) x\right\rangle
-\left\langle f\left( A\right) x,x\right\rangle \left\langle g\left(
B\right) x,x\right\rangle \right\vert  \label{II.a.e.4.1} \\
& \leq \max_{\lambda \in \left[ m,M\right] }\left\vert \left\langle
E_{\lambda }x,g\left( B\right) x\right\rangle -\left\langle E_{\lambda
}x,x\right\rangle \left\langle x,g\left( B\right) x\right\rangle \right\vert
\dbigvee\limits_{m}^{M}\left( f\right)  \notag \\
& \leq \max_{\lambda \in \left[ m,M\right] }\left( \left\langle E_{\lambda
}x,x\right\rangle \left\langle \left( 1_{H}-E_{\lambda }\right)
x,x\right\rangle \right) ^{1/2}  \notag \\
& \times \left( \left\Vert g\left( B\right) x\right\Vert ^{2}-\left\vert
\left\langle g\left( B\right) x,x\right\rangle \right\vert ^{2}\right)
^{1/2}\dbigvee\limits_{m}^{M}\left( f\right)  \notag \\
& \leq \frac{1}{2}\left( \left\Vert g\left( B\right) x\right\Vert
^{2}-\left\vert \left\langle g\left( B\right) x,x\right\rangle \right\vert
^{2}\right) ^{1/2}\dbigvee\limits_{m}^{M}\left( f\right) \leq \frac{1}{4}%
\left( N-n\right) \dbigvee\limits_{m}^{M}\left( f\right)  \notag
\end{align}%
for any $x\in H,\left\Vert x\right\Vert =1.$

2. If $f:\left[ m,M\right] \rightarrow \mathbb{C}$ is a Lipschitzian
function with the constant $L>0$ on $\left[ m,M\right] $, then we have the
inequality%
\begin{align}
& \left\vert \left\langle f\left( A\right) x,g\left( B\right) x\right\rangle
-\left\langle f\left( A\right) x,x\right\rangle \left\langle g\left(
B\right) x,x\right\rangle \right\vert  \label{II.a.e.4.2} \\
& \leq L\int_{m-0}^{M}\left\vert \left\langle E_{\lambda }x,g\left( B\right)
x\right\rangle -\left\langle E_{\lambda }x,x\right\rangle \left\langle
x,g\left( B\right) x\right\rangle \right\vert d\lambda  \notag \\
& \leq L\left( \left\Vert g\left( B\right) x\right\Vert ^{2}-\left\vert
\left\langle g\left( B\right) x,x\right\rangle \right\vert ^{2}\right) ^{1/2}
\notag \\
& \times \int_{m-0}^{M}\left( \left\langle E_{\lambda }x,x\right\rangle
\left\langle \left( 1_{H}-E_{\lambda }\right) x,x\right\rangle \right)
^{1/2}d\lambda  \notag \\
& \leq L\left( \left\Vert g\left( B\right) x\right\Vert ^{2}-\left\vert
\left\langle g\left( B\right) x,x\right\rangle \right\vert ^{2}\right) ^{1/2}
\notag \\
& \times \left\langle \left( M1_{H}-A\right) x,x\right\rangle
^{1/2}\left\langle \left( A-m1_{H}\right) x,x\right\rangle ^{1/2}  \notag \\
& \leq \frac{1}{2}\left( M-m\right) L\left( \left\Vert g\left( B\right)
x\right\Vert ^{2}-\left\vert \left\langle g\left( B\right) x,x\right\rangle
\right\vert ^{2}\right) ^{1/2}  \notag \\
& \leq \frac{1}{4}\left( N-n\right) \left( M-m\right) L  \notag
\end{align}%
for any $x\in H,\left\Vert x\right\Vert =1.$

3. If $f:\left[ m,M\right] \rightarrow \mathbb{R}$ is a continuous monotonic
nondecreasing function on $\left[ m,M\right] $, then we have the inequality%
\begin{align}
& \left\vert \left\langle f\left( A\right) x,g\left( B\right) x\right\rangle
-\left\langle f\left( A\right) x,x\right\rangle \left\langle g\left(
B\right) x,x\right\rangle \right\vert  \label{II.a.e.4.3} \\
& \leq \int_{m-0}^{M}\left\vert \left\langle E_{\lambda }x,g\left( B\right)
x\right\rangle -\left\langle E_{\lambda }x,x\right\rangle \left\langle
x,g\left( B\right) x\right\rangle \right\vert df\left( \lambda \right) 
\notag \\
& \leq \left( \left\Vert g\left( B\right) x\right\Vert ^{2}-\left\vert
\left\langle g\left( B\right) x,x\right\rangle \right\vert ^{2}\right) ^{1/2}
\notag \\
& \times \int_{m-0}^{M}\left( \left\langle E_{\lambda }x,x\right\rangle
\left\langle \left( 1_{H}-E_{\lambda }\right) x,x\right\rangle \right)
^{1/2}df\left( \lambda \right)  \notag \\
& \leq \left( \left\Vert g\left( B\right) x\right\Vert ^{2}-\left\vert
\left\langle g\left( B\right) x,x\right\rangle \right\vert ^{2}\right) ^{1/2}
\notag \\
& \times \left\langle \left( f\left( M\right) 1_{H}-f\left( A\right) \right)
x,x\right\rangle ^{1/2}\left\langle \left( f\left( A\right) -f\left(
m\right) 1_{H}\right) x,x\right\rangle ^{1/2}  \notag \\
& \leq \frac{1}{2}\left[ f\left( M\right) -f\left( m\right) \right] \left(
\left\Vert g\left( B\right) x\right\Vert ^{2}-\left\vert \left\langle
g\left( B\right) x,x\right\rangle \right\vert ^{2}\right) ^{1/2}  \notag \\
& \leq \frac{1}{4}\left( N-n\right) \left[ f\left( M\right) -f\left(
m\right) \right]  \notag
\end{align}%
for any $x\in H,\left\Vert x\right\Vert =1.$
\end{proposition}

\begin{proof}
We notice that, since $n:=\min_{t\in \left[ m,M\right] }g\left( t\right) $
and $N:=\max_{t\in \left[ m,M\right] }g\left( t\right) ,$ then $n\leq
\left\langle g\left( B\right) x,x\right\rangle \leq N$ which implies that $%
\left\langle g\left( B\right) x-nx,Mx-g\left( B\right) x\right\rangle \geq 0$
for any $x\in H,\left\Vert x\right\Vert =1.$ On applying Theorem \ref%
{II.a.t.1.1} for $y=Bx,\Gamma =N$ and $\gamma =n$ we deduce the desired
result.
\end{proof}

\begin{remark}
\label{II.a.r.3.2}We observe that if the function $f$ takes real values and
is a $\left( \varphi ,\Phi \right) $-Lipschitzian function on $\left[ m,M%
\right] $, then the inequality (\ref{II.a.e.4.2}) can be improved as follows%
\begin{align}
& \left\vert \left\langle f\left( A\right) x,g\left( B\right) x\right\rangle
-\left\langle f\left( A\right) x,x\right\rangle \left\langle g\left(
B\right) x,x\right\rangle \right\vert  \label{II.a.e.4.4} \\
& \leq \frac{1}{2}\left( \Phi -\varphi \right) \int_{m-0}^{M}\left\vert
\left\langle E_{\lambda }x,g\left( B\right) x\right\rangle -\left\langle
E_{\lambda }x,x\right\rangle \left\langle x,g\left( B\right) x\right\rangle
\right\vert d\lambda  \notag \\
& \leq \frac{1}{2}\left( \Phi -\varphi \right) \left( \left\Vert g\left(
B\right) x\right\Vert ^{2}-\left\vert \left\langle g\left( B\right)
x,x\right\rangle \right\vert ^{2}\right) ^{1/2}  \notag \\
& \times \int_{m-0}^{M}\left( \left\langle E_{\lambda }x,x\right\rangle
\left\langle \left( 1_{H}-E_{\lambda }\right) x,x\right\rangle \right)
^{1/2}d\lambda  \notag \\
& \leq \frac{1}{2}\left( \Phi -\varphi \right) \left( \left\Vert g\left(
B\right) x\right\Vert ^{2}-\left\vert \left\langle g\left( B\right)
x,x\right\rangle \right\vert ^{2}\right) ^{1/2}  \notag \\
& \times \left\langle \left( M1_{H}-A\right) x,x\right\rangle
^{1/2}\left\langle \left( A-m1_{H}\right) x,x\right\rangle ^{1/2}  \notag \\
& \leq \frac{1}{4}\left( M-m\right) \left( \Phi -\varphi \right) \left(
\left\Vert g\left( B\right) x\right\Vert ^{2}-\left\vert \left\langle
g\left( B\right) x,x\right\rangle \right\vert ^{2}\right) ^{1/2}  \notag \\
& \leq \frac{1}{8}\left( N-n\right) \left( M-m\right) \left( \Phi -\varphi
\right)  \notag
\end{align}%
for any $x\in H,\left\Vert x\right\Vert =1.$
\end{remark}

\subsection{Applications}

By choosing different examples of elementary functions into the above
inequalities, one can obtain various Gr\"{u}ss' type inequalities of
interest.

For instance, if we choose $f,g:\left( 0,\infty \right) \rightarrow \left(
0,\infty \right) $ with $f\left( t\right) =t^{p},g\left( t\right) =t^{q}$
and $p,q>0,$ then for any selfadjoint operators $A,B$ with $Sp\left(
A\right) ,Sp\left( B\right) \subseteq \left[ m,M\right] \subset \left(
0,\infty \right) $ we get from (\ref{II.a.e.4.3}) the inequalities

\begin{align}
& \left\vert \left\langle A^{p}x,B^{q}x\right\rangle -\left\langle
A^{p}x,x\right\rangle \left\langle B^{q}x,x\right\rangle \right\vert
\label{II.a.e.5.1} \\
& \leq p\left( \left\Vert B^{q}x\right\Vert ^{2}-\left\vert \left\langle
B^{q}x,x\right\rangle \right\vert ^{2}\right) ^{1/2}\int_{m-0}^{M}\left(
\left\langle E_{\lambda }x,x\right\rangle \left\langle \left(
1_{H}-E_{\lambda }\right) x,x\right\rangle \right) ^{1/2}\lambda
^{p-1}d\lambda  \notag \\
& \leq \left( \left\Vert B^{q}x\right\Vert ^{2}-\left\vert \left\langle
B^{q}x,x\right\rangle \right\vert ^{2}\right) ^{1/2}\left\langle \left(
M^{p}1_{H}-A^{p}\right) x,x\right\rangle ^{1/2}\left\langle \left(
A^{p}-m^{p}1_{H}\right) x,x\right\rangle ^{1/2}  \notag \\
& \leq \frac{1}{2}\left( M^{p}-m^{p}\right) \left( \left\Vert
B^{q}x\right\Vert ^{2}-\left\vert \left\langle B^{q}x,x\right\rangle
\right\vert ^{2}\right) ^{1/2}\leq \frac{1}{4}\left( M^{q}-m^{q}\right)
\left( M^{p}-m^{p}\right)  \notag
\end{align}%
for any $x\in H$ with $\left\Vert x\right\Vert =1,$ where $\left\{
E_{\lambda }\right\} _{\lambda }$ is the spectral family\textit{\ of }$A$%
\textit{.}

The same choice of functions considered in the inequality (\ref{II.a.e.4.4})
produce the result 
\begin{align}
& \left\vert \left\langle A^{p}x,B^{q}x\right\rangle -\left\langle
A^{p}x,x\right\rangle \left\langle B^{q}x,x\right\rangle \right\vert
\label{II.a.e.5.2} \\
& \leq \frac{1}{2}\Delta _{p}\left( \left\Vert B^{q}x\right\Vert
^{2}-\left\vert \left\langle B^{q}x,x\right\rangle \right\vert ^{2}\right)
^{1/2}  \notag \\
& \times \int_{m-0}^{M}\left( \left\langle E_{\lambda }x,x\right\rangle
\left\langle \left( 1_{H}-E_{\lambda }\right) x,x\right\rangle \right)
^{1/2}d\lambda  \notag \\
& \leq \frac{1}{2}\Delta _{p}\left( \left\Vert B^{q}x\right\Vert
^{2}-\left\vert \left\langle B^{q}x,x\right\rangle \right\vert ^{2}\right)
^{1/2}  \notag \\
& \times \left\langle \left( M^{p}1_{H}-A^{p}\right) x,x\right\rangle
^{1/2}\left\langle \left( A^{p}-m^{p}1_{H}\right) x,x\right\rangle ^{1/2} 
\notag \\
& \leq \frac{1}{4}\left( M-m\right) \Delta _{p}\left( \left\Vert
B^{q}x\right\Vert ^{2}-\left\vert \left\langle B^{q}x,x\right\rangle
\right\vert ^{2}\right) ^{1/2}  \notag \\
& \leq \frac{1}{8}\left( M^{q}-m^{q}\right) \left( M-m\right) \Delta _{p} 
\notag
\end{align}%
where%
\begin{equation}
\Delta _{p}:=p\times \left\{ 
\begin{array}{c}
M^{p-1}-m^{p-1}\text{ if }p\geq 1 \\ 
\\ 
\frac{M^{1-p}-m^{1-p}}{M^{1-p}m^{1-p}}\text{ if }0<p<1.%
\end{array}%
\right.  \label{II.a.e.5.3}
\end{equation}%
for any $x\in H$ with $\left\Vert x\right\Vert =1.$

Now, if we choose $f\left( t\right) =\ln t,t>0$ and keep the same $g$ then
we have the inequalities%
\begin{align}
& \left\vert \left\langle \ln Ax,B^{q}x\right\rangle -\left\langle \ln
Ax,x\right\rangle \left\langle B^{q}x,x\right\rangle \right\vert
\label{II.a.e.5.4} \\
& \leq \left( \left\Vert B^{q}x\right\Vert ^{2}-\left\vert \left\langle
B^{q}x,x\right\rangle \right\vert ^{2}\right) ^{1/2}  \notag \\
& \times \int_{m-0}^{M}\left( \left\langle E_{\lambda }x,x\right\rangle
\left\langle \left( 1_{H}-E_{\lambda }\right) x,x\right\rangle \right)
^{1/2}\lambda ^{-1}d\lambda  \notag \\
& \leq \left( \left\Vert B^{q}x\right\Vert ^{2}-\left\vert \left\langle
B^{q}x,x\right\rangle \right\vert ^{2}\right) ^{1/2}  \notag \\
& \times \left\langle \left( \ln M1_{H}-\ln A\right) x,x\right\rangle
^{1/2}\left\langle \left( \ln A-\ln m1_{H}\right) x,x\right\rangle ^{1/2} 
\notag \\
& \leq \left( \left\Vert B^{q}x\right\Vert ^{2}-\left\vert \left\langle
B^{q}x,x\right\rangle \right\vert ^{2}\right) ^{1/2}\ln \sqrt{\frac{M}{m}} 
\notag \\
& \leq \frac{1}{2}\left( M^{q}-m^{q}\right) \ln \sqrt{\frac{M}{m}}  \notag
\end{align}%
and%
\begin{align}
& \left\vert \left\langle \ln Ax,B^{q}x\right\rangle -\left\langle \ln
Ax,x\right\rangle \left\langle B^{q}x,x\right\rangle \right\vert
\label{II.a.e.5.5} \\
& \leq \frac{1}{2}\left( \frac{M-m}{mM}\right) \left( \left\Vert
B^{q}x\right\Vert ^{2}-\left\vert \left\langle B^{q}x,x\right\rangle
\right\vert ^{2}\right) ^{1/2}  \notag \\
& \times \int_{m-0}^{M}\left( \left\langle E_{\lambda }x,x\right\rangle
\left\langle \left( 1_{H}-E_{\lambda }\right) x,x\right\rangle \right)
^{1/2}d\lambda  \notag \\
& \leq \frac{1}{2}\left( \frac{M-m}{mM}\right) \left( \left\Vert
B^{q}x\right\Vert ^{2}-\left\vert \left\langle B^{q}x,x\right\rangle
\right\vert ^{2}\right) ^{1/2}  \notag \\
& \times \left\langle \left( M1_{H}-A\right) x,x\right\rangle
^{1/2}\left\langle \left( A-m1_{H}\right) x,x\right\rangle ^{1/2}  \notag \\
& \leq \frac{1}{4}\frac{\left( M-m\right) ^{2}}{mM}\left( \left\Vert
B^{q}x\right\Vert ^{2}-\left\vert \left\langle B^{q}x,x\right\rangle
\right\vert ^{2}\right) ^{1/2}  \notag \\
& \leq \frac{1}{8}\left( M^{q}-m^{q}\right) \frac{\left( M-m\right) ^{2}}{mM}
\notag
\end{align}%
for any $x\in H$ with $\left\Vert x\right\Vert =1.$

\section{Two Operators Gr\"{u}ss' Type Inequalities}

\subsection{Some Representation Results}

We start with the following representation result that will play a key role
in obtaining various bounds for different choices of functions including
continuous functions of bounded variation, Lipschitzian functions or
monotonic and continuous functions.

\begin{theorem}[Dragomir, 2010, \protect\cite{II.b.SSD3}]
\label{II.b.t.2.1}Let $A,B$ be two selfadjoint operators in the Hilbert
space $H$ with the spectra $Sp\left( A\right) ,Sp\left( B\right) \subseteq %
\left[ m,M\right] $ for some real numbers $m<M$ and let $\left\{ E_{\lambda
}\right\} _{\lambda }$ be the \textit{spectral family of }$A$\textit{\ and }$%
\left\{ F_{\mu }\right\} _{\mu }$ the \textit{spectral family of }$B.$ If $%
f,g:\left[ m,M\right] \rightarrow \mathbb{C}$ are continuous, then we have
the representation%
\begin{align}
& \left\langle f\left( A\right) x,g\left( B\right) x\right\rangle
-\left\langle f\left( A\right) x,x\right\rangle \left\langle x,g\left(
B\right) x\right\rangle  \label{II.b.e.2.1} \\
& =\int_{m-0}^{M}\left( \int_{m-0}^{M}\left[ \left\langle E_{\lambda
}x,x\right\rangle \left\langle x,F_{\mu }x\right\rangle -\left\langle
E_{\lambda }x,F_{\mu }x\right\rangle \right] d\left( g\left( \mu \right)
\right) \right) d\left( f\left( \lambda \right) \right)  \notag
\end{align}%
for any $x\in H$ with $\left\Vert x\right\Vert =1.$
\end{theorem}

\begin{proof}
Integrating by parts in the Riemann-Stieltjes integral and making use of the
spectral representation theorem we have%
\begin{align}
& \int_{m-0}^{M}\left[ \left\langle E_{\lambda }x,y\right\rangle
-\left\langle E_{\lambda }x,x\right\rangle \left\langle x,y\right\rangle %
\right] df\left( \lambda \right)  \label{II.b.e.2.2} \\
& =\left. \left[ \left\langle E_{\lambda }x,y\right\rangle -\left\langle
E_{\lambda }x,x\right\rangle \left\langle x,y\right\rangle \right] f\left(
\lambda \right) \right\vert _{m-0}^{M}  \notag \\
& -\int_{m-0}^{M}f\left( \lambda \right) d\left[ \left\langle E_{\lambda
}x,y\right\rangle -\left\langle E_{\lambda }x,x\right\rangle \left\langle
x,y\right\rangle \right]  \notag \\
& =\left\langle x,y\right\rangle \int_{m-0}^{M}f\left( \lambda \right)
d\left\langle E_{\lambda }x,x\right\rangle -\int_{m-0}^{M}f\left( \lambda
\right) d\left\langle E_{\lambda }x,y\right\rangle  \notag \\
& =\left\langle x,y\right\rangle \left\langle f\left( A\right)
x,x\right\rangle -\left\langle f\left( A\right) x,y\right\rangle  \notag
\end{align}%
for any $x,y\in H$ with $\left\Vert x\right\Vert =1.$

Now, if we chose $y=g\left( B\right) x$ in (\ref{II.b.e.2.2}) then we get
that%
\begin{align}
& \int_{m-0}^{M}\left[ \left\langle E_{\lambda }x,g\left( B\right)
x\right\rangle -\left\langle E_{\lambda }x,x\right\rangle \left\langle
x,g\left( B\right) x\right\rangle \right] df\left( \lambda \right)
\label{II.b.e.2.3} \\
& =\left\langle x,g\left( B\right) x\right\rangle \left\langle f\left(
A\right) x,x\right\rangle -\left\langle f\left( A\right) x,g\left( B\right)
x\right\rangle  \notag
\end{align}%
for any $x\in H$ with $\left\Vert x\right\Vert =1.$

Utilising the spectral representation theorem for $B$ we also have for each
fixed $\lambda \in \left[ m,M\right] $ that%
\begin{align}
& \left\langle E_{\lambda }x,g\left( B\right) x\right\rangle -\left\langle
E_{\lambda }x,x\right\rangle \left\langle x,g\left( B\right) x\right\rangle
\label{II.b.e.2.4} \\
& =\left\langle E_{\lambda }x,\int_{m-0}^{M}g\left( \mu \right) dF_{\mu
}x\right\rangle -\left\langle E_{\lambda }x,x\right\rangle \left\langle
x,\int_{m-0}^{M}g\left( \mu \right) dF_{\mu }x\right\rangle  \notag \\
& =\int_{m-0}^{M}g\left( \mu \right) d\left( \left\langle E_{\lambda
}x,F_{\mu }x\right\rangle \right) -\left\langle E_{\lambda }x,x\right\rangle
\int_{m-0}^{M}g\left( \mu \right) d\left( \left\langle x,F_{\mu
}x\right\rangle \right)  \notag
\end{align}%
for any $x\in H$ with $\left\Vert x\right\Vert =1.$

Integrating by parts in the Riemann-Stieltjes integral we have%
\begin{eqnarray*}
\int_{m-0}^{M}g\left( \mu \right) d\left( \left\langle E_{\lambda }x,F_{\mu
}x\right\rangle \right) &=&\left. g\left( \mu \right) \left\langle
E_{\lambda }x,F_{\mu }x\right\rangle \right] _{m-0}^{M}-\int_{m-0}^{M}\left%
\langle E_{\lambda }x,F_{\mu }x\right\rangle dg\left( \mu \right) \\
&=&g\left( M\right) \left\langle E_{\lambda }x,x\right\rangle
-\int_{m-0}^{M}\left\langle E_{\lambda }x,F_{\mu }x\right\rangle d\left(
g\left( \mu \right) \right)
\end{eqnarray*}%
and 
\begin{eqnarray*}
\int_{m-0}^{M}g\left( \mu \right) d\left( \left\langle x,F_{\mu
}x\right\rangle \right) &=&\left. g\left( \mu \right) \left\langle x,F_{\mu
}x\right\rangle \right] _{m-0}^{M}-\int_{m-0}^{M}\left\langle x,F_{\mu
}x\right\rangle d\left( g\left( \mu \right) \right) \\
&=&g\left( M\right) -\int_{m-0}^{M}\left\langle x,F_{\mu }x\right\rangle
d\left( g\left( \mu \right) \right) ,
\end{eqnarray*}%
therefore%
\begin{align}
& \int_{m-0}^{M}g\left( \mu \right) d\left( \left\langle E_{\lambda
}x,F_{\mu }x\right\rangle \right) -\left\langle E_{\lambda }x,x\right\rangle
\int_{m-0}^{M}g\left( \mu \right) d\left( \left\langle x,F_{\mu
}x\right\rangle \right)  \label{II.b.e.2.5} \\
& =g\left( M\right) \left\langle E_{\lambda }x,x\right\rangle
-\int_{m-0}^{M}\left\langle E_{\lambda }x,F_{\mu }x\right\rangle d\left(
g\left( \mu \right) \right)  \notag \\
& -\left\langle E_{\lambda }x,x\right\rangle \left( g\left( M\right)
-\int_{m-0}^{M}\left\langle x,F_{\mu }x\right\rangle d\left( g\left( \mu
\right) \right) \right)  \notag \\
& =\left\langle E_{\lambda }x,x\right\rangle \int_{m-0}^{M}\left\langle
x,F_{\mu }x\right\rangle d\left( g\left( \mu \right) \right)
-\int_{m-0}^{M}\left\langle E_{\lambda }x,F_{\mu }x\right\rangle d\left(
g\left( \mu \right) \right)  \notag \\
& =\int_{m-0}^{M}\left[ \left\langle E_{\lambda }x,x\right\rangle
\left\langle x,F_{\mu }x\right\rangle -\left\langle E_{\lambda }x,F_{\mu
}x\right\rangle \right] d\left( g\left( \mu \right) \right)  \notag
\end{align}%
for any $x\in H$ with $\left\Vert x\right\Vert =1$ and $\lambda \in \left[
m,M\right] .$

Utilising (\ref{II.b.e.2.3})-(\ref{II.b.e.2.5}) we deduce the desired result
(\ref{II.b.e.2.1}).
\end{proof}

\begin{remark}
\label{II.b.r.2.1}In particular, if we take $B=A,$ then we get from (\ref%
{II.b.e.2.1}) the equality%
\begin{align}
& \left\langle f\left( A\right) x,g\left( A\right) x\right\rangle
-\left\langle f\left( A\right) x,x\right\rangle \left\langle x,g\left(
A\right) x\right\rangle  \label{II.b.e.2.5.a} \\
& =\int_{m-0}^{M}\left( \int_{m-0}^{M}\left[ \left\langle E_{\lambda
}x,x\right\rangle \left\langle x,E_{\mu }x\right\rangle -\left\langle
E_{\lambda }x,E_{\mu }x\right\rangle \right] d\left( g\left( \mu \right)
\right) \right) d\left( f\left( \lambda \right) \right)  \notag
\end{align}%
for any $x\in H$ with $\left\Vert x\right\Vert =1,$ which for $g=f$ produces
the representation result for the variance of the selfadjoint operator $%
f\left( A\right) ,$%
\begin{align}
& \left\Vert f\left( A\right) x\right\Vert ^{2}-\left\langle f\left(
A\right) x,x\right\rangle ^{2}  \label{II.b.e.2.5.b} \\
& =\int_{m-0}^{M}\left( \int_{m-0}^{M}\left[ \left\langle E_{\lambda
}x,x\right\rangle \left\langle x,E_{\mu }x\right\rangle -\left\langle
E_{\lambda }x,E_{\mu }x\right\rangle \right] d\left( f\left( \mu \right)
\right) \right) d\left( f\left( \lambda \right) \right)  \notag
\end{align}%
for any $x\in H$ with $\left\Vert x\right\Vert =1.$
\end{remark}

\subsection{Bounds for $f$ of Bounded Variation}

The first vectorial Gr\"{u}ss' type inequality when one of the functions is
of bounded variation is as follows:

\begin{theorem}[Dragomir, 2010, \protect\cite{II.b.SSD3}]
\label{II.b.t.2.2} Let $A,B$ be two selfadjoint operators in the Hilbert
space $H$ with the spectra $Sp\left( A\right) ,Sp\left( B\right) \subseteq %
\left[ m,M\right] $ for some real numbers $m<M$ and let $\left\{ E_{\lambda
}\right\} _{\lambda }$ be the \textit{spectral family of }$A$\textit{\ and }$%
\left\{ F_{\mu }\right\} _{\mu }$ the \textit{spectral family of }$B.$ Also,
assume that $f:\left[ m,M\right] \rightarrow \mathbb{C}$ is continuous and
of bounded variation on $\left[ m,M\right] .$

1. If $g:\left[ m,M\right] \rightarrow \mathbb{C}$ is continuous and of
bounded variation on $\left[ m,M\right] ,$ then we have the inequality 
\begin{align}
& \left\vert \left\langle f\left( A\right) x,g\left( B\right) x\right\rangle
-\left\langle f\left( A\right) x,x\right\rangle \left\langle x,g\left(
B\right) x\right\rangle \right\vert  \label{II.b.e.2.6} \\
& \leq \max_{\left( \lambda ,\mu \right) \in \left[ m,M\right]
^{2}}\left\vert \left\langle E_{\lambda }x,x\right\rangle \left\langle
x,F_{\mu }x\right\rangle -\left\langle E_{\lambda }x,F_{\mu }x\right\rangle
\right\vert \dbigvee\limits_{m}^{M}\left( g\right)
\dbigvee\limits_{m}^{M}\left( f\right)  \notag \\
& \leq \max_{\lambda \in \left[ m,M\right] }\left[ \left\langle E_{\lambda
}x,x\right\rangle \left\langle \left( 1_{H}-E_{\lambda }\right)
x,x\right\rangle \right] ^{1/2}  \notag \\
& \times \max_{\mu \in \left[ m,M\right] }\left[ \left\langle F_{\mu
}x,x\right\rangle \left\langle \left( 1_{H}-F_{\mu }\right) x,x\right\rangle %
\right] ^{1/2}\dbigvee\limits_{m}^{M}\left( g\right)
\dbigvee\limits_{m}^{M}\left( f\right) \leq \frac{1}{4}\dbigvee%
\limits_{m}^{M}\left( g\right) \dbigvee\limits_{m}^{M}\left( f\right)  \notag
\end{align}%
for any $x\in H$ with $\left\Vert x\right\Vert =1.$

2. If $g:\left[ m,M\right] \rightarrow \mathbb{C}$ is Lipschitzian with the
constant $K>0$ on $\left[ m,M\right] ,$ then we have the inequality 
\begin{align}
& \left\vert \left\langle f\left( A\right) x,g\left( B\right) x\right\rangle
-\left\langle f\left( A\right) x,x\right\rangle \left\langle x,g\left(
B\right) x\right\rangle \right\vert  \label{II.b.e.2.7} \\
& \leq K\max_{\lambda \in \left[ m,M\right] }\left[ \int_{m-0}^{M}\left\vert
\left\langle E_{\lambda }x,x\right\rangle \left\langle x,F_{\mu
}x\right\rangle -\left\langle E_{\lambda }x,F_{\mu }x\right\rangle
\right\vert d\mu \right] \dbigvee\limits_{m}^{M}\left( f\right)  \notag \\
& \leq K\dbigvee\limits_{m}^{M}\left( f\right) \max_{\lambda \in \left[ m,M%
\right] }\left[ \left\langle E_{\lambda }x,x\right\rangle \left\langle
\left( 1_{H}-E_{\lambda }\right) x,x\right\rangle \right] ^{1/2}  \notag \\
& \times \int_{m-0}^{M}\left[ \left\langle F_{\mu }x,x\right\rangle
\left\langle \left( 1_{H}-F_{\mu }\right) x,x\right\rangle \right] ^{1/2}d\mu
\notag \\
& \leq \frac{1}{2}K\dbigvee\limits_{m}^{M}\left( f\right) \left\langle
\left( M1_{H}-B\right) x,x\right\rangle ^{1/2}\left\langle \left(
B-m1_{H}\right) x,x\right\rangle ^{1/2}  \notag \\
& \leq \frac{1}{4}K\left( M-m\right) \dbigvee\limits_{m}^{M}\left( f\right) 
\notag
\end{align}%
for any $x\in H$ with $\left\Vert x\right\Vert =1.$

3. If $g:\left[ m,M\right] \rightarrow \mathbb{R}$ is continuous and
monotonic nondecreasing on $\left[ m,M\right] ,$ then we have the inequality 
\begin{align}
& \left\vert \left\langle f\left( A\right) x,g\left( B\right) x\right\rangle
-\left\langle f\left( A\right) x,x\right\rangle \left\langle x,g\left(
B\right) x\right\rangle \right\vert  \label{II.b.e.2.8} \\
& \leq \max_{\lambda \in \left[ m,M\right] }\left[ \int_{m-0}^{M}\left\vert
\left\langle E_{\lambda }x,x\right\rangle \left\langle x,F_{\mu
}x\right\rangle -\left\langle E_{\lambda }x,F_{\mu }x\right\rangle
\right\vert dg\left( \mu \right) \right] \dbigvee\limits_{m}^{M}\left(
f\right)  \notag \\
& \leq \dbigvee\limits_{m}^{M}\left( f\right) \max_{\lambda \in \left[ m,M%
\right] }\left[ \left\langle E_{\lambda }x,x\right\rangle \left\langle
\left( 1_{H}-E_{\lambda }\right) x,x\right\rangle \right] ^{1/2}  \notag \\
& \times \int_{m-0}^{M}\left[ \left\langle F_{\mu }x,x\right\rangle
\left\langle \left( 1_{H}-F_{\mu }\right) x,x\right\rangle \right]
^{1/2}dg\left( \mu \right)  \notag \\
& \leq \frac{1}{2}\dbigvee\limits_{m}^{M}\left( f\right) \left\langle \left(
g\left( M\right) 1_{H}-g\left( B\right) \right) x,x\right\rangle
^{1/2}\left\langle \left( g\left( B\right) -g\left( m\right) 1_{H}\right)
x,x\right\rangle ^{1/2}  \notag \\
& \leq \frac{1}{4}\left[ g\left( M\right) -g\left( m\right) \right]
\dbigvee\limits_{m}^{M}\left( f\right)  \notag
\end{align}%
for any $x\in H$ with $\left\Vert x\right\Vert =1.$
\end{theorem}

\begin{proof}
1. It is well known that if $p:\left[ a,b\right] \rightarrow \mathbb{C}$ is
a continuous function, $v:\left[ a,b\right] \rightarrow \mathbb{C}$ is of
bounded variation then the Riemann-Stieltjes integral $\int_{a}^{b}p\left(
t\right) dv\left( t\right) $ exists and the following inequality holds%
\begin{equation}
\left\vert \int_{a}^{b}p\left( t\right) dv\left( t\right) \right\vert \leq
\max_{t\in \left[ a,b\right] }\left\vert p\left( t\right) \right\vert
\dbigvee\limits_{a}^{b}\left( v\right) ,  \label{II.b.e.2.9}
\end{equation}%
where $\dbigvee\limits_{a}^{b}\left( v\right) $ denotes the total variation
of $v$ on $\left[ a,b\right] .$

Now, on utilizing the property (\ref{II.b.e.2.9}) and the identity (\ref%
{II.b.e.2.1}) we have%
\begin{align}
& \left\vert \left\langle f\left( A\right) x,g\left( B\right) x\right\rangle
-\left\langle f\left( A\right) x,x\right\rangle \left\langle x,g\left(
B\right) x\right\rangle \right\vert  \label{II.b.e.2.10} \\
& \leq \max_{\lambda \in \left[ m,M\right] }\left\vert \int_{m-0}^{M}\left[
\left\langle E_{\lambda }x,x\right\rangle \left\langle x,F_{\mu
}x\right\rangle -\left\langle E_{\lambda }x,F_{\mu }x\right\rangle \right]
d\left( g\left( \mu \right) \right) \right\vert
\dbigvee\limits_{m}^{M}\left( f\right)  \notag
\end{align}%
for any $x\in \left[ m,M\right] .$

The same inequality (\ref{II.b.e.2.9}) produces the bound%
\begin{align}
& \max_{\lambda \in \left[ m,M\right] }\left\vert \int_{m-0}^{M}\left[
\left\langle E_{\lambda }x,x\right\rangle \left\langle x,F_{\mu
}x\right\rangle -\left\langle E_{\lambda }x,F_{\mu }x\right\rangle \right]
d\left( g\left( \mu \right) \right) \right\vert  \label{II.b.e.2.11} \\
& \leq \max_{\lambda \in \left[ m,M\right] }\left[ \max_{\mu \in \left[ m,M%
\right] }\left\vert \left\langle E_{\lambda }x,x\right\rangle \left\langle
x,F_{\mu }x\right\rangle -\left\langle E_{\lambda }x,F_{\mu }x\right\rangle
\right\vert \right] \dbigvee\limits_{m}^{M}\left( f\right)  \notag \\
& =\max_{\left( \lambda ,\mu \right) \in \left[ m,M\right] ^{2}}\left\vert
\left\langle E_{\lambda }x,x\right\rangle \left\langle x,F_{\mu
}x\right\rangle -\left\langle E_{\lambda }x,F_{\mu }x\right\rangle
\right\vert \dbigvee\limits_{m}^{M}\left( f\right)  \notag
\end{align}%
for any $x\in \left[ m,M\right] .$

By making use of (\ref{II.b.e.2.10}) and (\ref{II.b.e.2.11}) we deduce the
first part of (\ref{II.b.e.2.6}).

Further, we notice that by the Schwarz inequality in $H$ we have for any $%
u,v,e\in H$ with $\left\Vert e\right\Vert =1$ that%
\begin{equation}
\left\vert \left\langle u,v\right\rangle -\left\langle u,e\right\rangle
\left\langle e,v\right\rangle \right\vert \leq \left( \left\Vert
u\right\Vert ^{2}-\left\vert \left\langle u,e\right\rangle \right\vert
^{2}\right) ^{1/2}\left( \left\Vert v\right\Vert ^{2}-\left\vert
\left\langle v,e\right\rangle \right\vert ^{2}\right) ^{1/2}.
\label{II.b.e.2.12}
\end{equation}%
Indeed, if we write Schwarz's inequality for the vectors $u-\left\langle
u,e\right\rangle e$ and $v-\left\langle v,e\right\rangle e$ we have%
\begin{equation*}
\left\vert \left\langle u-\left\langle u,e\right\rangle e,v-\left\langle
v,e\right\rangle e\right\rangle \right\vert \leq \left\Vert u-\left\langle
u,e\right\rangle e\right\Vert \left\Vert v-\left\langle v,e\right\rangle
e\right\Vert
\end{equation*}%
which, by performing the calculations, is equivalent with (\ref{II.b.e.2.12}%
).

Now, on utilizing (\ref{II.b.e.2.12}), we can state that%
\begin{align}
& \left\vert \left\langle E_{\lambda }x,x\right\rangle \left\langle x,F_{\mu
}x\right\rangle -\left\langle E_{\lambda }x,F_{\mu }x\right\rangle
\right\vert  \label{II.b.e.2.13} \\
& \leq \left( \left\Vert E_{\lambda }x\right\Vert ^{2}-\left\vert
\left\langle E_{\lambda }x,x\right\rangle \right\vert ^{2}\right)
^{1/2}\left( \left\Vert F_{\mu }x\right\Vert ^{2}-\left\vert \left\langle
F_{\mu }x,x\right\rangle \right\vert ^{2}\right) ^{1/2}  \notag
\end{align}%
for any $\lambda ,\mu \in \left[ m,M\right] .$

Since $E_{\lambda }$ and $F_{\mu }$ are projections and $E_{\lambda },F_{\mu
}\geq 0$ then 
\begin{align}
\left\Vert E_{\lambda }x\right\Vert ^{2}-\left\vert \left\langle E_{\lambda
}x,x\right\rangle \right\vert ^{2}& =\left\langle E_{\lambda
}x,x\right\rangle -\left\langle E_{\lambda }x,x\right\rangle ^{2}
\label{II.b.e.2.14} \\
& =\left\langle E_{\lambda }x,x\right\rangle \left\langle \left(
1_{H}-E_{\lambda }\right) x,x\right\rangle \leq \frac{1}{4}  \notag
\end{align}%
and 
\begin{equation}
\left\Vert F_{\mu }x\right\Vert ^{2}-\left\vert \left\langle F_{\mu
}x,x\right\rangle \right\vert ^{2}=\left\langle F_{\mu }x,x\right\rangle
\left\langle \left( 1_{H}-F_{\mu }\right) x,x\right\rangle \leq \frac{1}{4}
\label{II.b.e.2.15}
\end{equation}

for any $\lambda ,\mu \in \left[ m,M\right] $ and $x\in H$ with $\left\Vert
x\right\Vert =1.$

Now, if we use (\ref{II.b.e.2.13})-(\ref{II.b.e.2.15}) then we get the
second part of (\ref{II.b.e.2.6}).

2. Further, recall that if $p:\left[ a,b\right] \rightarrow \mathbb{C}$ is a
Riemann integrable function and $v:\left[ a,b\right] \rightarrow \mathbb{C}$
is Lipschitzian with the constant $L>0$, i.e.,%
\begin{equation*}
\left\vert f\left( s\right) -f\left( t\right) \right\vert \leq L\left\vert
s-t\right\vert \text{ for any }t,s\in \left[ a,b\right] ,
\end{equation*}%
then the Riemann-Stieltjes integral $\int_{a}^{b}p\left( t\right) dv\left(
t\right) $ exists and the following inequality holds%
\begin{equation}
\left\vert \int_{a}^{b}p\left( t\right) dv\left( t\right) \right\vert \leq
L\int_{a}^{b}\left\vert p\left( t\right) \right\vert dt.  \label{II.b.e.2.16}
\end{equation}

If we use the inequality (\ref{II.b.e.2.16}), then we have in the case when $%
g$ is Lipschitzian with the constant $K>0$ that%
\begin{align}
& \max_{\lambda \in \left[ m,M\right] }\left\vert \int_{m-0}^{M}\left[
\left\langle E_{\lambda }x,x\right\rangle \left\langle x,F_{\mu
}x\right\rangle -\left\langle E_{\lambda }x,F_{\mu }x\right\rangle \right]
d\left( g\left( \mu \right) \right) \right\vert  \label{II.b.e.2.17} \\
& \leq K\max_{\lambda \in \left[ m,M\right] }\left[ \int_{m-0}^{M}\left\vert
\left\langle E_{\lambda }x,x\right\rangle \left\langle x,F_{\mu
}x\right\rangle -\left\langle E_{\lambda }x,F_{\mu }x\right\rangle
\right\vert d\mu \right]  \notag
\end{align}%
for any $x\in H$ with $\left\Vert x\right\Vert =1$ and the first part of (%
\ref{II.b.e.2.7}) is proved.

Further, by employing (\ref{II.b.e.2.13})-(\ref{II.b.e.2.15}) we also get
that%
\begin{align}
& \max_{\lambda \in \left[ m,M\right] }\int_{m-0}^{M}\left\vert \left\langle
E_{\lambda }x,x\right\rangle \left\langle x,F_{\mu }x\right\rangle
-\left\langle E_{\lambda }x,F_{\mu }x\right\rangle \right\vert d\mu
\label{II.b.e.2.18} \\
& \leq \max_{\lambda \in \left[ m,M\right] }\left[ \left\langle E_{\lambda
}x,x\right\rangle \left\langle \left( 1_{H}-E_{\lambda }\right)
x,x\right\rangle \right] ^{1/2}  \notag \\
& \times \int_{m-0}^{M}\left[ \left\langle F_{\mu }x,x\right\rangle
\left\langle \left( 1_{H}-F_{\mu }\right) x,x\right\rangle \right] ^{1/2}d\mu
\notag
\end{align}%
for any $x\in H$ with $\left\Vert x\right\Vert =1.$

If we use the Cauchy-Bunyakovsky-Schwarz integral inequality and the
spectral representation theorem, then we have successively 
\begin{align}
& \int_{m-0}^{M}\left( \left\langle F_{\mu }x,x\right\rangle \left\langle
\left( 1_{H}-F_{\mu }\right) x,x\right\rangle \right) ^{1/2}d\mu
\label{II.b.e.2.19} \\
& \leq \left[ \int_{m-0}^{M}\left\langle F_{\mu }x,x\right\rangle d\mu %
\right] ^{1/2}\left[ \int_{m-0}^{M}\left\langle \left( 1_{H}-F_{\mu }\right)
x,x\right\rangle d\mu \right] ^{1/2}  \notag \\
& =\left[ \left. \left\langle F_{\mu }x,x\right\rangle \mu \right\vert
_{m-0}^{M}-\int_{m-0}^{M}\mu d\left\langle F_{\mu }x,x\right\rangle \right]
^{1/2}  \notag \\
& \times \left[ \left. \left\langle \left( 1_{H}-F_{\mu }\right)
x,x\right\rangle \mu \right\vert _{m-0}^{M}-\int_{m-0}^{M}\mu d\left\langle
\left( 1_{H}-F_{\mu }\right) x,x\right\rangle \right]  \notag \\
& =\left\langle \left( M1_{H}-B\right) x,x\right\rangle ^{1/2}\left\langle
\left( B-m1_{H}\right) x,x\right\rangle ^{1/2},  \notag
\end{align}%
for any $x\in H$ with $\left\Vert x\right\Vert =1.$

On employing now (\ref{II.b.e.2.17})-(\ref{II.b.e.2.19}) we deduce the
second part of (\ref{II.b.e.2.7}).

The last part of (\ref{II.b.e.2.7}) follows by the elementary inequality 
\begin{equation}
\alpha \beta \leq \frac{1}{4}\left( \alpha +\beta \right) ^{2},\alpha \beta
\geq 0  \label{II.b.e.2.19.1}
\end{equation}%
for the choice $\alpha =\left\langle \left( M1_{H}-B\right) x,x\right\rangle 
$ and $\beta =\left\langle \left( B-m1_{H}\right) x,x\right\rangle $ and the
details are omitted.

3. Further, from the theory of Riemann-Stieltjes integral it is also well
known that if $p:\left[ a,b\right] \rightarrow \mathbb{C}$ is of bounded
variation and $v:\left[ a,b\right] \rightarrow \mathbb{R}$ is continuous and
monotonic nondecreasing, then the Riemann-Stieltjes integrals $%
\int_{a}^{b}p\left( t\right) dv\left( t\right) $ and $\int_{a}^{b}\left\vert
p\left( t\right) \right\vert dv\left( t\right) $ exist and%
\begin{equation}
\left\vert \int_{a}^{b}p\left( t\right) dv\left( t\right) \right\vert \leq
\int_{a}^{b}\left\vert p\left( t\right) \right\vert dv\left( t\right) .
\label{II.b.e.2.20}
\end{equation}

Now, if we assume that $g$ is monotonic nondecreasing on $\left[ m,M\right]
, $ then by (\ref{II.b.e.2.20}) we have that%
\begin{align}
& \max_{\lambda \in \left[ m,M\right] }\left\vert \int_{m-0}^{M}\left[
\left\langle E_{\lambda }x,x\right\rangle \left\langle x,F_{\mu
}x\right\rangle -\left\langle E_{\lambda }x,F_{\mu }x\right\rangle \right]
d\left( g\left( \mu \right) \right) \right\vert  \label{II.b.e.2.21} \\
& \leq \max_{\lambda \in \left[ m,M\right] }\left[ \int_{m-0}^{M}\left\vert
\left\langle E_{\lambda }x,x\right\rangle \left\langle x,F_{\mu
}x\right\rangle -\left\langle E_{\lambda }x,F_{\mu }x\right\rangle
\right\vert dg\left( \mu \right) \right]  \notag
\end{align}%
for any $x\in H$ with $\left\Vert x\right\Vert =1.$

Further, by employing (\ref{II.b.e.2.13})-(\ref{II.b.e.2.15}) we also get
that%
\begin{align}
& \max_{\lambda \in \left[ m,M\right] }\int_{m-0}^{M}\left\vert \left\langle
E_{\lambda }x,x\right\rangle \left\langle x,F_{\mu }x\right\rangle
-\left\langle E_{\lambda }x,F_{\mu }x\right\rangle \right\vert dg\left( \mu
\right)  \label{II.b.e.2.22} \\
& \leq \max_{\lambda \in \left[ m,M\right] }\left[ \left\langle E_{\lambda
}x,x\right\rangle \left\langle \left( 1_{H}-E_{\lambda }\right)
x,x\right\rangle \right] ^{1/2}  \notag \\
& \times \int_{m-0}^{M}\left[ \left\langle F_{\mu }x,x\right\rangle
\left\langle \left( 1_{H}-F_{\mu }\right) x,x\right\rangle \right]
^{1/2}dg\left( \mu \right)  \notag
\end{align}%
for any $x\in H$ with $\left\Vert x\right\Vert =1.$ These prove the first
part of (\ref{II.b.e.2.8}).

If we use the Cauchy-Bunyakovsky-Schwarz integral inequality for the
Riemann-Stieltjes integral with monotonic nondecreasing integrators and the
spectral representation theorem, then we have successively 
\begin{align}
& \int_{m-0}^{M}\left( \left\langle F_{\mu }x,x\right\rangle \left\langle
\left( 1_{H}-F_{\mu }\right) x,x\right\rangle \right) ^{1/2}dg\left( \mu
\right)  \label{II.b.e.2.23} \\
& \leq \left[ \int_{m-0}^{M}\left\langle F_{\mu }x,x\right\rangle dg\left(
\mu \right) \right] ^{1/2}\left[ \int_{m-0}^{M}\left\langle \left(
1_{H}-F_{\mu }\right) x,x\right\rangle dg\left( \mu \right) \right] ^{1/2} 
\notag \\
& =\left[ \left. \left\langle F_{\mu }x,x\right\rangle g\left( \mu \right)
\right\vert _{m-0}^{M}-\int_{m-0}^{M}g\left( \mu \right) d\left\langle
F_{\mu }x,x\right\rangle \right] ^{1/2}  \notag \\
& \times \left[ \left. \left\langle \left( 1_{H}-F_{\mu }\right)
x,x\right\rangle g\left( \mu \right) \right\vert
_{m-0}^{M}-\int_{m-0}^{M}g\left( \mu \right) d\left\langle \left(
1_{H}-F_{\mu }\right) x,x\right\rangle \right] ^{1/2}  \notag \\
& =\left\langle \left( g\left( M\right) 1_{H}-g\left( B\right) \right)
x,x\right\rangle ^{1/2}\left\langle \left( g\left( B\right) -g\left(
m\right) 1_{H}\right) x,x\right\rangle ^{1/2},  \notag
\end{align}%
for any $x\in H$ with $\left\Vert x\right\Vert =1.$

Utilising (\ref{II.b.e.2.23}) we then deduce the last part of (\ref%
{II.b.e.2.8}). The details are omitted.
\end{proof}

Now, in order to provide other results that are similar to the Gr\"{u}ss'
type inequalities stated in the introduction, we can state the following
corollary:

\begin{corollary}[Dragomir, 2010, \protect\cite{II.b.SSD3}]
\label{II.b.c.2.1}Let $A$ be a selfadjoint operators in the Hilbert space $H$
with the spectrum $Sp\left( A\right) \subseteq \left[ m,M\right] $ for some
real numbers $m<M$ and let $\left\{ E_{\lambda }\right\} _{\lambda }$ be the 
\textit{spectral family of }$A.$ Also, assume that $f:\left[ m,M\right]
\rightarrow \mathbb{C}$ is continuous and of bounded variation on $\left[ m,M%
\right] .$

1. If $g:\left[ m,M\right] \rightarrow \mathbb{C}$ is continuous and of
bounded variation on $\left[ m,M\right] ,$ then we have the inequality 
\begin{align}
& \left\vert \left\langle f\left( A\right) x,g\left( A\right) x\right\rangle
-\left\langle f\left( A\right) x,x\right\rangle \left\langle x,g\left(
A\right) x\right\rangle \right\vert  \label{II.b.e.2.24} \\
& \leq \max_{\left( \lambda ,\mu \right) \in \left[ m,M\right]
^{2}}\left\vert \left\langle E_{\lambda }x,x\right\rangle \left\langle
x,E_{\mu }x\right\rangle -\left\langle E_{\lambda }x,E_{\mu }x\right\rangle
\right\vert \dbigvee\limits_{m}^{M}\left( g\right)
\dbigvee\limits_{m}^{M}\left( f\right)  \notag \\
& \leq \max_{\lambda \in \left[ m,M\right] }\left[ \left\langle E_{\lambda
}x,x\right\rangle \left\langle \left( 1_{H}-E_{\lambda }\right)
x,x\right\rangle \right] \dbigvee\limits_{m}^{M}\left( g\right)
\dbigvee\limits_{m}^{M}\left( f\right) \leq \frac{1}{4}\dbigvee%
\limits_{m}^{M}\left( g\right) \dbigvee\limits_{m}^{M}\left( f\right)  \notag
\end{align}%
for any $x\in H$ with $\left\Vert x\right\Vert =1.$

2. If $g:\left[ m,M\right] \rightarrow \mathbb{C}$ is Lipschitzian with the
constant $K>0$ on $\left[ m,M\right] ,$ then we have the inequality 
\begin{align}
& \left\vert \left\langle f\left( A\right) x,g\left( A\right) x\right\rangle
-\left\langle f\left( A\right) x,x\right\rangle \left\langle x,g\left(
A\right) x\right\rangle \right\vert  \label{II.b.e.2.25} \\
& \leq K\max_{\lambda \in \left[ m,M\right] }\left[ \int_{m-0}^{M}\left\vert
\left\langle E_{\lambda }x,x\right\rangle \left\langle x,E_{\mu
}x\right\rangle -\left\langle E_{\lambda }x,E_{\mu }x\right\rangle
\right\vert d\mu \right] \dbigvee\limits_{m}^{M}\left( f\right)  \notag \\
& \leq K\dbigvee\limits_{m}^{M}\left( f\right) \max_{\lambda \in \left[ m,M%
\right] }\left[ \left\langle E_{\lambda }x,x\right\rangle \left\langle
\left( 1_{H}-E_{\lambda }\right) x,x\right\rangle \right] ^{1/2}  \notag \\
& \times \int_{m-0}^{M}\left[ \left\langle E_{\mu }x,x\right\rangle
\left\langle \left( 1_{H}-E_{\mu }\right) x,x\right\rangle \right] ^{1/2}d\mu
\notag \\
& \leq \frac{1}{2}K\dbigvee\limits_{m}^{M}\left( f\right) \left\langle
\left( M1_{H}-A\right) x,x\right\rangle ^{1/2}\left\langle \left(
A-m1_{H}\right) x,x\right\rangle ^{1/2}  \notag \\
& \leq \frac{1}{4}K\left( M-m\right) \dbigvee\limits_{m}^{M}\left( f\right) 
\notag
\end{align}%
for any $x\in H$ with $\left\Vert x\right\Vert =1.$

3. If $g:\left[ m,M\right] \rightarrow \mathbb{R}$ is continuous and
monotonic nondecreasing on $\left[ m,M\right] ,$ then we have the inequality 
\begin{align}
& \left\vert \left\langle f\left( A\right) x,g\left( A\right) x\right\rangle
-\left\langle f\left( A\right) x,x\right\rangle \left\langle x,g\left(
A\right) x\right\rangle \right\vert  \label{II.b.e.2.26} \\
& \leq \max_{\lambda \in \left[ m,M\right] }\left[ \int_{m-0}^{M}\left\vert
\left\langle E_{\lambda }x,x\right\rangle \left\langle x,E_{\mu
}x\right\rangle -\left\langle E_{\lambda }x,E_{\mu }x\right\rangle
\right\vert dg\left( \mu \right) \right] \dbigvee\limits_{m}^{M}\left(
f\right)  \notag \\
& \leq \dbigvee\limits_{m}^{M}\left( f\right) \max_{\lambda \in \left[ m,M%
\right] }\left[ \left\langle E_{\lambda }x,x\right\rangle \left\langle
\left( 1_{H}-E_{\lambda }\right) x,x\right\rangle \right] ^{1/2}  \notag \\
& \times \int_{m-0}^{M}\left[ \left\langle E_{\mu }x,x\right\rangle
\left\langle \left( 1_{H}-E_{\mu }\right) x,x\right\rangle \right]
^{1/2}dg\left( \mu \right)  \notag \\
& \leq \frac{1}{2}\dbigvee\limits_{m}^{M}\left( f\right) \left\langle \left(
g\left( M\right) 1_{H}-g\left( A\right) \right) x,x\right\rangle
^{1/2}\left\langle \left( g\left( A\right) -g\left( m\right) 1_{H}\right)
x,x\right\rangle ^{1/2}  \notag \\
& \leq \frac{1}{4}\left[ g\left( M\right) -g\left( m\right) \right]
\dbigvee\limits_{m}^{M}\left( f\right)  \notag
\end{align}%
for any $x\in H$ with $\left\Vert x\right\Vert =1.$
\end{corollary}

\begin{remark}
\label{II.b.r.2.2}The following inequality for the variance of $f\left(
A\right) $ under the assumptions that $A$ is a selfadjoint operators in the
Hilbert space $H$ with the spectrum $Sp\left( A\right) \subseteq \left[ m,M%
\right] $ for some real numbers $m<M$, $\left\{ E_{\lambda }\right\}
_{\lambda }$ is the \textit{spectral family of }$A$ and $f:\left[ m,M\right]
\rightarrow \mathbb{C}$ is continuous and of bounded variation on $\left[ m,M%
\right] $ can be stated%
\begin{align}
0& \leq \left\Vert f\left( A\right) x\right\Vert ^{2}-\left\langle f\left(
A\right) x,x\right\rangle ^{2}  \label{II.b.e.2.27} \\
& \leq \max_{\left( \lambda ,\mu \right) \in \left[ m,M\right]
^{2}}\left\vert \left\langle E_{\lambda }x,x\right\rangle \left\langle
x,E_{\mu }x\right\rangle -\left\langle E_{\lambda }x,E_{\mu }x\right\rangle
\right\vert \left[ \dbigvee\limits_{m}^{M}\left( f\right) \right] ^{2} 
\notag \\
& \leq \max_{\lambda \in \left[ m,M\right] }\left[ \left\langle E_{\lambda
}x,x\right\rangle \left\langle \left( 1_{H}-E_{\lambda }\right)
x,x\right\rangle \right] \left[ \dbigvee\limits_{m}^{M}\left( f\right) %
\right] ^{2}\leq \frac{1}{4}\left[ \dbigvee\limits_{m}^{M}\left( f\right) %
\right] ^{2}  \notag
\end{align}%
for any $x\in H$ with $\left\Vert x\right\Vert =1.$
\end{remark}

\subsection{Bounds for $f$ Lipschitzian}

The case when the first function is Lipschitzian is as follows:

\begin{theorem}[Dragomir, 2010, \protect\cite{II.b.SSD3}]
\label{II.b.t.3.1}Let $A,B$ be two selfadjoint operators in the Hilbert
space $H$ with the spectra $Sp\left( A\right) ,Sp\left( B\right) \subseteq %
\left[ m,M\right] $ for some real numbers $m<M$ and let $\left\{ E_{\lambda
}\right\} _{\lambda }$ be the \textit{spectral family of }$A$\textit{\ and }$%
\left\{ F_{\mu }\right\} _{\mu }$ the \textit{spectral family of }$B.$ Also,
assume that $f:\left[ m,M\right] \rightarrow \mathbb{C}$ is Lipschitzian
with the constant $L>0$ on $\left[ m,M\right] .$

1. If $g:\left[ m,M\right] \rightarrow \mathbb{C}$ is Lipschitzian with the
constant $K>0$ on $\left[ m,M\right] ,$ then we have the inequality 
\begin{align}
& \left\vert \left\langle f\left( A\right) x,g\left( B\right) x\right\rangle
-\left\langle f\left( A\right) x,x\right\rangle \left\langle x,g\left(
B\right) x\right\rangle \right\vert  \label{II.b.e.3.1} \\
& \leq LK\int_{m-0}^{M}\int_{m-0}^{M}\left\vert \left\langle E_{\lambda
}x,x\right\rangle \left\langle x,F_{\mu }x\right\rangle -\left\langle
E_{\lambda }x,F_{\mu }x\right\rangle \right\vert d\mu d\lambda  \notag \\
& \leq LK\int_{m-0}^{M}\left[ \left\langle E_{\lambda }x,x\right\rangle
\left\langle \left( 1_{H}-E_{\lambda }\right) x,x\right\rangle \right]
^{1/2}d\lambda  \notag \\
& \times \int_{m-0}^{M}\left[ \left\langle F_{\mu }x,x\right\rangle
\left\langle \left( 1_{H}-F_{\mu }\right) x,x\right\rangle \right] ^{1/2}d\mu
\notag \\
& \leq LK\left[ \left\langle \left( M1_{H}-A\right) x,x\right\rangle
\left\langle \left( A-m1_{H}\right) x,x\right\rangle \right] ^{1/2}  \notag
\\
& \times \left[ \left\langle \left( M1_{H}-B\right) x,x\right\rangle
\left\langle \left( B-m1_{H}\right) x,x\right\rangle \right] ^{1/2}\leq 
\frac{1}{4}LK\left( M-m\right) ^{2}  \notag
\end{align}%
for any $x\in H$ with $\left\Vert x\right\Vert =1.$

2. If $g:\left[ m,M\right] \rightarrow \mathbb{R}$ is continuous and
monotonic nondecreasing on $\left[ m,M\right] ,$ then we have the inequality 
\begin{align}
& \left\vert \left\langle f\left( A\right) x,g\left( B\right) x\right\rangle
-\left\langle f\left( A\right) x,x\right\rangle \left\langle x,g\left(
B\right) x\right\rangle \right\vert  \label{II.b.e.3.2} \\
& \leq L\int_{m-0}^{M}\int_{m-0}^{M}\left\vert \left\langle E_{\lambda
}x,x\right\rangle \left\langle x,F_{\mu }x\right\rangle -\left\langle
E_{\lambda }x,F_{\mu }x\right\rangle \right\vert dg\left( \mu \right)
d\lambda  \notag \\
& \leq L\int_{m-0}^{M}\left[ \left\langle E_{\lambda }x,x\right\rangle
\left\langle \left( 1_{H}-E_{\lambda }\right) x,x\right\rangle \right]
^{1/2}d\lambda  \notag \\
& \times \int_{m-0}^{M}\left[ \left\langle F_{\mu }x,x\right\rangle
\left\langle \left( 1_{H}-F_{\mu }\right) x,x\right\rangle \right]
^{1/2}dg\left( \mu \right)  \notag \\
& \leq L\left[ \left\langle \left( M1_{H}-A\right) x,x\right\rangle
\left\langle \left( A-m1_{H}\right) x,x\right\rangle \right] ^{1/2}  \notag
\\
& \times \left[ \left\langle \left( g\left( M\right) 1_{H}-g\left( B\right)
\right) x,x\right\rangle \left\langle \left( g\left( B\right) -g\left(
m\right) 1_{H}\right) x,x\right\rangle \right] ^{1/2}  \notag \\
& \leq \frac{1}{4}L\left( M-m\right) \left[ g\left( M\right) -g\left(
m\right) \right]  \notag
\end{align}%
for any $x\in H$ with $\left\Vert x\right\Vert =1.$
\end{theorem}

\begin{proof}
1. We observe that, on utilizing the property (\ref{II.b.e.2.16}) and the
identity (\ref{II.b.e.2.1}) we have%
\begin{align}
& \left\vert \left\langle f\left( A\right) x,g\left( B\right) x\right\rangle
-\left\langle f\left( A\right) x,x\right\rangle \left\langle x,g\left(
B\right) x\right\rangle \right\vert  \label{II.b.e.3.3} \\
& \leq L\int_{m-0}^{M}\left\vert \int_{m-0}^{M}\left[ \left\langle
E_{\lambda }x,x\right\rangle \left\langle x,F_{\mu }x\right\rangle
-\left\langle E_{\lambda }x,F_{\mu }x\right\rangle \right] d\left( g\left(
\mu \right) \right) \right\vert d\lambda  \notag
\end{align}%
for any $x\in H,\left\Vert x\right\Vert =1.$

By the same property (\ref{II.b.e.2.16}) we also have%
\begin{align}
& \left\vert \int_{m-0}^{M}\left[ \left\langle E_{\lambda }x,x\right\rangle
\left\langle x,F_{\mu }x\right\rangle -\left\langle E_{\lambda }x,F_{\mu
}x\right\rangle \right] d\left( g\left( \mu \right) \right) \right\vert
\label{II.b.e.3.4} \\
& \leq K\int_{m-0}^{M}\left\vert \left\langle E_{\lambda }x,x\right\rangle
\left\langle x,F_{\mu }x\right\rangle -\left\langle E_{\lambda }x,F_{\mu
}x\right\rangle \right\vert d\mu  \notag
\end{align}%
for any $x\in H,\left\Vert x\right\Vert =1$ and $\lambda \in \left[ m,M%
\right] .$

Therefore, by (\ref{II.b.e.3.3}) and (\ref{II.b.e.3.4}) we get%
\begin{align}
& \left\vert \left\langle f\left( A\right) x,g\left( B\right) x\right\rangle
-\left\langle f\left( A\right) x,x\right\rangle \left\langle x,g\left(
B\right) x\right\rangle \right\vert  \label{II.b.e.3.5} \\
& \leq LK\int_{m-0}^{M}\int_{m-0}^{M}\left\vert \left\langle E_{\lambda
}x,x\right\rangle \left\langle x,F_{\mu }x\right\rangle -\left\langle
E_{\lambda }x,F_{\mu }x\right\rangle \right\vert d\mu d\lambda  \notag
\end{align}%
for any $x\in H,\left\Vert x\right\Vert =1,$ which proves the first
inequality in (\ref{II.b.e.3.1}).

From (\ref{II.b.e.2.13})-(\ref{II.b.e.2.15}) we have%
\begin{align}
& \left\vert \left\langle E_{\lambda }x,x\right\rangle \left\langle x,F_{\mu
}x\right\rangle -\left\langle E_{\lambda }x,F_{\mu }x\right\rangle
\right\vert  \label{II.b.e.3.6} \\
& \leq \left[ \left\langle E_{\lambda }x,x\right\rangle \left\langle \left(
1_{H}-E_{\lambda }\right) x,x\right\rangle \right] ^{1/2}\left[ \left\langle
F_{\mu }x,x\right\rangle \left\langle \left( 1_{H}-F_{\mu }\right)
x,x\right\rangle \right] ^{1/2}  \notag
\end{align}%
for any $x\in H,\left\Vert x\right\Vert =1$ and $\lambda ,\mu \in \left[ m,M%
\right] .$

Integrating on $\left[ m,M\right] ^{2}$ the inequality (\ref{II.b.e.3.6})
and utilizing the Cauchy-Bunyakowsky-Schwarz integral inequality for the
Riemann integral we have%
\begin{align}
& \int_{m-0}^{M}\int_{m-0}^{M}\left\vert \left\langle E_{\lambda
}x,x\right\rangle \left\langle x,F_{\mu }x\right\rangle -\left\langle
E_{\lambda }x,F_{\mu }x\right\rangle \right\vert d\mu d\lambda
\label{II.b.e.3.7} \\
& \leq \int_{m-0}^{M}\left[ \left\langle E_{\lambda }x,x\right\rangle
\left\langle \left( 1_{H}-E_{\lambda }\right) x,x\right\rangle \right]
^{1/2}d\lambda  \notag \\
& \times \int_{m-0}^{M}\left[ \left\langle F_{\mu }x,x\right\rangle
\left\langle \left( 1_{H}-F_{\mu }\right) x,x\right\rangle \right] ^{1/2}d\mu
\notag \\
& \leq \left[ \int_{m-0}^{M}\left\langle E_{\lambda }x,x\right\rangle
d\lambda \right] ^{1/2}\left[ \int_{m-0}^{M}\left\langle \left(
1_{H}-E_{\lambda }\right) x,x\right\rangle d\lambda \right] ^{1/2}  \notag \\
& \times \left[ \int_{m-0}^{M}\left\langle F_{\mu }x,x\right\rangle d\mu %
\right] ^{1/2}\left[ \int_{m-0}^{M}\left\langle \left( 1_{H}-F_{\mu }\right)
x,x\right\rangle d\mu \right] ^{1/2}.  \notag
\end{align}%
Integrating by parts and utilizing the spectral representation theorem we
have%
\begin{eqnarray*}
\int_{m-0}^{M}\left\langle E_{\lambda }x,x\right\rangle d\lambda &=&\left.
\left\langle E_{\lambda }x,x\right\rangle \lambda \right\vert
_{m-0}^{M}-\int_{m-0}^{M}\lambda d\left\langle E_{\lambda }x,x\right\rangle
\\
&=&M-\left\langle Ax,x\right\rangle =\left\langle \left( M1_{H}-A\right)
x,x\right\rangle ,
\end{eqnarray*}%
\begin{equation*}
\int_{m-0}^{M}\left\langle \left( 1_{H}-E_{\lambda }\right) x,x\right\rangle
d\lambda =\left\langle \left( A-m1_{H}\right) x,x\right\rangle
\end{equation*}%
and the similar equalities for $B,$ providing the second part of (\ref%
{II.b.e.3.1}).

The last part follows from (\ref{II.b.e.2.19.1}) and we omit the details.

2. Utilising the inequality (\ref{II.b.e.2.20}) we have%
\begin{align}
& \left\vert \int_{m-0}^{M}\left[ \left\langle E_{\lambda }x,x\right\rangle
\left\langle x,F_{\mu }x\right\rangle -\left\langle E_{\lambda }x,F_{\mu
}x\right\rangle \right] d\left( g\left( \mu \right) \right) \right\vert
\label{II.b.e.3.8} \\
& \leq \int_{m-0}^{M}\left\vert \left\langle E_{\lambda }x,x\right\rangle
\left\langle x,F_{\mu }x\right\rangle -\left\langle E_{\lambda }x,F_{\mu
}x\right\rangle \right\vert dg\left( \mu \right)  \notag
\end{align}%
which, together with (\ref{II.b.e.3.3}), produces the inequality%
\begin{align}
& \left\vert \left\langle f\left( A\right) x,g\left( B\right) x\right\rangle
-\left\langle f\left( A\right) x,x\right\rangle \left\langle x,g\left(
B\right) x\right\rangle \right\vert  \label{II.b.e.3.9} \\
& \leq L\int_{m-0}^{M}\int_{m-0}^{M}\left\vert \left\langle E_{\lambda
}x,x\right\rangle \left\langle x,F_{\mu }x\right\rangle -\left\langle
E_{\lambda }x,F_{\mu }x\right\rangle \right\vert dg\left( \mu \right)
d\lambda  \notag
\end{align}%
for any $x\in H,\left\Vert x\right\Vert =1.$

Now, by utilizing (\ref{II.b.e.3.6}) and a similar argument to the one
outlined above, we deduce the desired result (\ref{II.b.e.3.2}) and the
details are omitted.
\end{proof}

The case of one operator is incorporated in

\begin{corollary}[Dragomir, 2010, \protect\cite{II.b.SSD3}]
\label{II.b.c.3.1}Let $A$ be a selfadjoint operators in the Hilbert space $H$
with the spectrum $Sp\left( A\right) \subseteq \left[ m,M\right] $ for some
real numbers $m<M$ and let $\left\{ E_{\lambda }\right\} _{\lambda }$ be the 
\textit{spectral family of }$A.$ Also, assume that $f:\left[ m,M\right]
\rightarrow \mathbb{C}$ is Lipschitzian with the constant $L>0$ on $\left[
m,M\right] .$

1. If $g:\left[ m,M\right] \rightarrow \mathbb{C}$ is Lipschitzian with the
constant $K>0$ on $\left[ m,M\right] ,$ then we have the inequality 
\begin{align}
& \left\vert \left\langle f\left( A\right) x,g\left( A\right) x\right\rangle
-\left\langle f\left( A\right) x,x\right\rangle \left\langle x,g\left(
A\right) x\right\rangle \right\vert  \label{II.b.e.3.10} \\
& \leq LK\int_{m-0}^{M}\int_{m-0}^{M}\left\vert \left\langle E_{\lambda
}x,x\right\rangle \left\langle x,E_{\mu }x\right\rangle -\left\langle
E_{\lambda }x,E_{\mu }x\right\rangle \right\vert d\mu d\lambda  \notag \\
& \leq LK\left( \int_{m-0}^{M}\left[ \left\langle E_{\lambda
}x,x\right\rangle \left\langle \left( 1_{H}-E_{\lambda }\right)
x,x\right\rangle \right] ^{1/2}d\lambda \right) ^{2}  \notag \\
& \leq LK\left[ \left\langle \left( M1_{H}-A\right) x,x\right\rangle
\left\langle \left( A-m1_{H}\right) x,x\right\rangle \right] \leq \frac{1}{4}%
LK\left( M-m\right) ^{2}  \notag
\end{align}%
for any $x\in H$ with $\left\Vert x\right\Vert =1.$

2. If $g:\left[ m,M\right] \rightarrow \mathbb{R}$ is continuous and
monotonic nondecreasing on $\left[ m,M\right] ,$ then we have the inequality 
\begin{align}
& \left\vert \left\langle f\left( A\right) x,g\left( A\right) x\right\rangle
-\left\langle f\left( A\right) x,x\right\rangle \left\langle x,g\left(
A\right) x\right\rangle \right\vert  \label{II.b.e.3.11} \\
& \leq L\int_{m-0}^{M}\int_{m-0}^{M}\left\vert \left\langle E_{\lambda
}x,x\right\rangle \left\langle x,F_{\mu }x\right\rangle -\left\langle
E_{\lambda }x,F_{\mu }x\right\rangle \right\vert dg\left( \mu \right)
d\lambda  \notag \\
& \leq L\int_{m-0}^{M}\left[ \left\langle E_{\lambda }x,x\right\rangle
\left\langle \left( 1_{H}-E_{\lambda }\right) x,x\right\rangle \right]
^{1/2}d\lambda  \notag \\
& \times \int_{m-0}^{M}\left[ \left\langle E_{\mu }x,x\right\rangle
\left\langle \left( 1_{H}-E_{\mu }\right) x,x\right\rangle \right]
^{1/2}dg\left( \mu \right)  \notag \\
& \leq L\left[ \left\langle \left( M1_{H}-A\right) x,x\right\rangle
\left\langle \left( A-m1_{H}\right) x,x\right\rangle \right] ^{1/2}  \notag
\\
& \times \left[ \left\langle \left( g\left( M\right) 1_{H}-g\left( A\right)
\right) x,x\right\rangle \left\langle \left( g\left( A\right) -g\left(
m\right) 1_{H}\right) x,x\right\rangle \right] ^{1/2}  \notag \\
& \leq \frac{1}{4}L\left( M-m\right) \left[ g\left( M\right) -g\left(
m\right) \right]  \notag
\end{align}%
for any $x\in H$ with $\left\Vert x\right\Vert =1.$
\end{corollary}

\begin{remark}
\label{II.b.r.3.1}The following inequality for the variance of $f\left(
A\right) $ under the assumptions that $A$ is a selfadjoint operators in the
Hilbert space $H$ with the spectrum $Sp\left( A\right) \subseteq \left[ m,M%
\right] $ for some real numbers $m<M$, $\left\{ E_{\lambda }\right\}
_{\lambda }$ is the \textit{spectral family of }$A$ and $f:\left[ m,M\right]
\rightarrow \mathbb{C}$ is Lipschitzian with the constant $L>0$ on $\left[
m,M\right] $ can be stated%
\begin{align}
& 0\leq \left\Vert f\left( A\right) x\right\Vert ^{2}-\left\langle f\left(
A\right) x,x\right\rangle ^{2}  \label{II.b.e.3.12} \\
& \leq L^{2}\int_{m-0}^{M}\int_{m-0}^{M}\left\vert \left\langle E_{\lambda
}x,x\right\rangle \left\langle x,E_{\mu }x\right\rangle -\left\langle
E_{\lambda }x,E_{\mu }x\right\rangle \right\vert d\mu d\lambda  \notag \\
& \leq L^{2}\left( \int_{m-0}^{M}\left[ \left\langle E_{\lambda
}x,x\right\rangle \left\langle \left( 1_{H}-E_{\lambda }\right)
x,x\right\rangle \right] ^{1/2}d\lambda \right) ^{2}  \notag \\
& \leq L^{2}\left[ \left\langle \left( M1_{H}-A\right) x,x\right\rangle
\left\langle \left( A-m1_{H}\right) x,x\right\rangle \right]  \notag \\
& \leq \frac{1}{4}L^{2}\left( M-m\right) ^{2}  \notag
\end{align}%
for any $x\in H$ with $\left\Vert x\right\Vert =1.$
\end{remark}

\subsection{Bounds for $f$ Monotonic Nondecreasing}

Finally, for the case of two monotonic functions we have the following
result as well:

\begin{theorem}[Dragomir, 2010, \protect\cite{II.b.SSD3}]
\label{II.b.t.4.1}Let $A,B$ be two selfadjoint operators in the Hilbert
space $H$ with the spectra $Sp\left( A\right) ,Sp\left( B\right) \subseteq %
\left[ m,M\right] $ for some real numbers $m<M$ and let $\left\{ E_{\lambda
}\right\} _{\lambda }$ be the \textit{spectral family of }$A$\textit{\ and }$%
\left\{ F_{\mu }\right\} _{\mu }$ the \textit{spectral family of }$B.$ If $%
f,g:\left[ m,M\right] \rightarrow \mathbb{C}$ are continuous and monotonic
nondecreasing on $\left[ m,M\right] ,$ then%
\begin{align}
& \left\vert \left\langle f\left( A\right) x,g\left( B\right) x\right\rangle
-\left\langle f\left( A\right) x,x\right\rangle \left\langle x,g\left(
B\right) x\right\rangle \right\vert  \label{II.b.e.4.1} \\
& \leq \int_{m-0}^{M}\int_{m-0}^{M}\left\vert \left\langle E_{\lambda
}x,x\right\rangle \left\langle x,F_{\mu }x\right\rangle -\left\langle
E_{\lambda }x,F_{\mu }x\right\rangle \right\vert dg\left( \mu \right)
df\left( \lambda \right)  \notag \\
& \leq \int_{m-0}^{M}\left[ \left\langle E_{\lambda }x,x\right\rangle
\left\langle \left( 1_{H}-E_{\lambda }\right) x,x\right\rangle \right]
^{1/2}df\left( \lambda \right)  \notag \\
& \times \int_{m-0}^{M}\left[ \left\langle F_{\mu }x,x\right\rangle
\left\langle \left( 1_{H}-F_{\mu }\right) x,x\right\rangle \right]
^{1/2}dg\left( \mu \right)  \notag \\
& \leq \left[ \left\langle \left( f\left( M\right) 1_{H}-f\left( A\right)
\right) x,x\right\rangle \left\langle \left( f\left( A\right) -f\left(
m\right) 1_{H}\right) x,x\right\rangle \right] ^{1/2}  \notag \\
& \times \left[ \left\langle \left( g\left( M\right) 1_{H}-g\left( B\right)
\right) x,x\right\rangle \left\langle \left( g\left( B\right) -g\left(
m\right) 1_{H}\right) x,x\right\rangle \right] ^{1/2}  \notag \\
& \leq \frac{1}{4}\left[ f\left( M\right) -f\left( m\right) \right] \left[
g\left( M\right) -g\left( m\right) \right]  \notag
\end{align}%
for any $x\in H,\left\Vert x\right\Vert =1.$
\end{theorem}

The details of the proof are omitted.

In particular we have:

\begin{corollary}[Dragomir, 2010, \protect\cite{II.b.SSD3}]
\label{II.b.c.4.1}Let $A$ be a selfadjoint operators in the Hilbert space $H$
with the spectrum $Sp\left( A\right) \subseteq \left[ m,M\right] $ for some
real numbers $m<M$ and let $\left\{ E_{\lambda }\right\} _{\lambda }$ be the 
\textit{spectral family of }$A.$ If $f,g:\left[ m,M\right] \rightarrow 
\mathbb{C}$ are continuous and monotonic nondecreasing on $\left[ m,M\right]
,$ then%
\begin{align}
& \left\vert \left\langle f\left( A\right) x,g\left( A\right) x\right\rangle
-\left\langle f\left( A\right) x,x\right\rangle \left\langle x,g\left(
A\right) x\right\rangle \right\vert  \label{II.b.e.4.2} \\
& \leq \int_{m-0}^{M}\int_{m-0}^{M}\left\vert \left\langle E_{\lambda
}x,x\right\rangle \left\langle x,E_{\mu }x\right\rangle -\left\langle
E_{\lambda }x,E_{\mu }x\right\rangle \right\vert dg\left( \mu \right)
df\left( \lambda \right)  \notag \\
& \leq \int_{m-0}^{M}\left[ \left\langle E_{\lambda }x,x\right\rangle
\left\langle \left( 1_{H}-E_{\lambda }\right) x,x\right\rangle \right]
^{1/2}df\left( \lambda \right)  \notag \\
& \times \int_{m-0}^{M}\left[ \left\langle E_{\mu }x,x\right\rangle
\left\langle \left( 1_{H}-E_{\mu }\right) x,x\right\rangle \right]
^{1/2}dg\left( \mu \right)  \notag \\
& \leq \left[ \left\langle \left( f\left( M\right) 1_{H}-f\left( A\right)
\right) x,x\right\rangle \left\langle \left( f\left( A\right) -f\left(
m\right) 1_{H}\right) x,x\right\rangle \right] ^{1/2}  \notag \\
& \times \left[ \left\langle \left( g\left( M\right) 1_{H}-g\left( A\right)
\right) x,x\right\rangle \left\langle \left( g\left( A\right) -g\left(
m\right) 1_{H}\right) x,x\right\rangle \right] ^{1/2}  \notag \\
& \leq \frac{1}{4}\left[ f\left( M\right) -f\left( m\right) \right] \left[
g\left( M\right) -g\left( m\right) \right]  \notag
\end{align}%
for any $x\in H,\left\Vert x\right\Vert =1.$

In particular, the following inequality for the variance of $f\left(
A\right) $ in the case of monotonic nondecreasing functions $f$ holds: 
\begin{align}
0& \leq \left\Vert f\left( A\right) x\right\Vert ^{2}-\left\langle f\left(
A\right) x,x\right\rangle ^{2}  \label{II.b.e.4.3} \\
& \leq \int_{m-0}^{M}\int_{m-0}^{M}\left\vert \left\langle E_{\lambda
}x,x\right\rangle \left\langle x,E_{\mu }x\right\rangle -\left\langle
E_{\lambda }x,E_{\mu }x\right\rangle \right\vert df\left( \mu \right)
df\left( \lambda \right)  \notag \\
& \leq \left( \int_{m-0}^{M}\left[ \left\langle E_{\lambda }x,x\right\rangle
\left\langle \left( 1_{H}-E_{\lambda }\right) x,x\right\rangle \right]
^{1/2}df\left( \lambda \right) \right) ^{2}  \notag \\
& \leq \left[ \left\langle \left( f\left( M\right) 1_{H}-f\left( A\right)
\right) x,x\right\rangle \left\langle \left( f\left( A\right) -f\left(
m\right) 1_{H}\right) x,x\right\rangle \right]  \notag \\
& \leq \frac{1}{4}\left[ f\left( M\right) -f\left( m\right) \right] ^{2} 
\notag
\end{align}%
for any $x\in H,\left\Vert x\right\Vert =1.$
\end{corollary}

\subsection{Applications}

By choosing different examples of elementary functions into the above
inequalities, one can obtain various Gr\"{u}ss' type inequalities of
interest.

For instance, if we choose $f,g:\left( 0,\infty \right) \rightarrow \left(
0,\infty \right) $ with $f\left( t\right) =t^{p},g\left( t\right) =t^{q}$
and $p,q>0,$ then for any selfadjoint operators $A,B$ with $Sp\left(
A\right) ,Sp\left( B\right) \subseteq \left[ m,M\right] \subset \left(
0,\infty \right) $ we get from (\ref{II.b.e.4.1}) the inequalities:

\begin{align}
& \left\vert \left\langle A^{p}x,B^{q}x\right\rangle -\left\langle
A^{p}x,x\right\rangle \left\langle B^{q}x,x\right\rangle \right\vert
\label{II.b.e.5.1} \\
& \leq pq\int_{m-0}^{M}\int_{m-0}^{M}\left\vert \left\langle E_{\lambda
}x,x\right\rangle \left\langle x,F_{\mu }x\right\rangle -\left\langle
E_{\lambda }x,F_{\mu }x\right\rangle \right\vert \mu ^{q-1}\lambda
^{p-1}d\mu d\lambda  \notag \\
& \leq pq\int_{m-0}^{M}\left[ \left\langle E_{\lambda }x,x\right\rangle
\left\langle \left( 1_{H}-E_{\lambda }\right) x,x\right\rangle \right]
^{1/2}\lambda ^{p-1}d\lambda  \notag \\
& \times \int_{m-0}^{M}\left[ \left\langle F_{\mu }x,x\right\rangle
\left\langle \left( 1_{H}-F_{\mu }\right) x,x\right\rangle \right] ^{1/2}\mu
^{q-1}d\mu  \notag \\
& \leq \left[ \left\langle \left( M^{p}1_{H}-A^{p}\right) x,x\right\rangle
\left\langle \left( A^{p}-m^{p}1_{H}\right) x,x\right\rangle \right] ^{1/2} 
\notag \\
& \times \left[ \left\langle \left( M^{q}1_{H}-B^{q}\right) x,x\right\rangle
\left\langle \left( B^{q}-m^{q}1_{H}\right) x,x\right\rangle \right] ^{1/2} 
\notag \\
& \leq \frac{1}{4}\left( M^{p}-m^{p}\right) \left( M^{q}-m^{q}\right)  \notag
\end{align}%
for any $x\in H$ with $\left\Vert x\right\Vert =1,$ where $\left\{
E_{\lambda }\right\} _{\lambda }$ is the spectral family of $A$\ and $%
\left\{ F_{\mu }\right\} _{\mu }$ is the spectral family of $B.$

When $B=A$ then by the \v{C}eby\v{s}ev's inequality for functions of same
monotonicity the inequality (\ref{II.b.e.5.1}) becomes%
\begin{align}
& 0\leq \left\langle A^{p}x,A^{q}x\right\rangle -\left\langle
A^{p}x,x\right\rangle \left\langle A^{q}x,x\right\rangle  \label{II.b.e.5.2}
\\
& \leq pq\int_{m-0}^{M}\int_{m-0}^{M}\left\vert \left\langle E_{\lambda
}x,x\right\rangle \left\langle x,E_{\mu }x\right\rangle -\left\langle
E_{\lambda }x,E_{\mu }x\right\rangle \right\vert \mu ^{q-1}\lambda
^{p-1}d\mu d\lambda  \notag \\
& \leq pq\int_{m-0}^{M}\left[ \left\langle E_{\lambda }x,x\right\rangle
\left\langle \left( 1_{H}-E_{\lambda }\right) x,x\right\rangle \right]
^{1/2}\lambda ^{p-1}d\lambda  \notag \\
& \times \int_{m-0}^{M}\left[ \left\langle E_{\mu }x,x\right\rangle
\left\langle \left( 1_{H}-E_{\mu }\right) x,x\right\rangle \right] ^{1/2}\mu
^{q-1}d\mu  \notag \\
& \leq \left[ \left\langle \left( M^{p}1_{H}-A^{p}\right) x,x\right\rangle
\left\langle \left( A^{p}-m^{p}1_{H}\right) x,x\right\rangle \right] ^{1/2} 
\notag \\
& \times \left[ \left\langle \left( M^{q}1_{H}-B^{q}\right) x,x\right\rangle
\left\langle \left( B^{q}-m^{q}1_{H}\right) x,x\right\rangle \right] ^{1/2} 
\notag \\
& \leq \frac{1}{4}\left( M^{p}-m^{p}\right) \left( M^{q}-m^{q}\right)  \notag
\end{align}%
for any $x\in H$ with $\left\Vert x\right\Vert =1$ and $p,q>0.$

Now. define the coefficients

\begin{equation}
\Delta _{p}:=p\times \left\{ 
\begin{array}{c}
M^{p-1}-m^{p-1}\text{ if }p\geq 1 \\ 
\\ 
\frac{M^{1-p}-m^{1-p}}{M^{1-p}m^{1-p}}\text{ if }0<p<1.%
\end{array}%
\right.  \label{II.b.e.5.3}
\end{equation}

On utilizing the inequality (\ref{II.b.e.3.1}) for the same power functions
considered above, we can state the inequality%
\begin{align}
& \left\vert \left\langle A^{p}x,B^{q}x\right\rangle -\left\langle
A^{p}x,x\right\rangle \left\langle B^{q}x,x\right\rangle \right\vert
\label{II.b.e.5.4} \\
& \leq \Delta _{p}\Delta _{q}\int_{m-0}^{M}\int_{m-0}^{M}\left\vert
\left\langle E_{\lambda }x,x\right\rangle \left\langle x,F_{\mu
}x\right\rangle -\left\langle E_{\lambda }x,F_{\mu }x\right\rangle
\right\vert d\mu d\lambda  \notag \\
& \leq \Delta _{p}\Delta _{q}\int_{m-0}^{M}\left[ \left\langle E_{\lambda
}x,x\right\rangle \left\langle \left( 1_{H}-E_{\lambda }\right)
x,x\right\rangle \right] ^{1/2}d\lambda  \notag \\
& \times \int_{m-0}^{M}\left[ \left\langle F_{\mu }x,x\right\rangle
\left\langle \left( 1_{H}-F_{\mu }\right) x,x\right\rangle \right] ^{1/2}d\mu
\notag \\
& \leq \Delta _{p}\Delta _{q}\left[ \left\langle \left( M1_{H}-A\right)
x,x\right\rangle \left\langle \left( A-m1_{H}\right) x,x\right\rangle \right]
^{1/2}  \notag \\
& \times \left[ \left\langle \left( M1_{H}-B\right) x,x\right\rangle
\left\langle \left( B-m1_{H}\right) x,x\right\rangle \right] ^{1/2}\leq 
\frac{1}{4}\Delta _{p}\Delta _{q}\left( M-m\right) ^{2}  \notag
\end{align}%
for any $x\in H$ with $\left\Vert x\right\Vert =1$ and $p,q>0.$

In particular, for $B=A$ we have%
\begin{align}
& 0\leq \left\langle A^{p}x,A^{q}x\right\rangle -\left\langle
A^{p}x,x\right\rangle \left\langle A^{q}x,x\right\rangle  \label{II.b.e.5.5}
\\
& \leq \Delta _{p}\Delta _{q}\int_{m-0}^{M}\int_{m-0}^{M}\left\vert
\left\langle E_{\lambda }x,x\right\rangle \left\langle x,E_{\mu
}x\right\rangle -\left\langle E_{\lambda }x,E_{\mu }x\right\rangle
\right\vert d\mu d\lambda  \notag \\
& \leq \Delta _{p}\Delta _{q}\left( \int_{m-0}^{M}\left[ \left\langle
E_{\lambda }x,x\right\rangle \left\langle \left( 1_{H}-E_{\lambda }\right)
x,x\right\rangle \right] ^{1/2}d\lambda \right) ^{2}  \notag \\
& \leq \Delta _{p}\Delta _{q}\left[ \left\langle \left( M1_{H}-A\right)
x,x\right\rangle \left\langle \left( A-m1_{H}\right) x,x\right\rangle \right]
\leq \frac{1}{4}\Delta _{p}\Delta _{q}\left( M-m\right) ^{2}  \notag
\end{align}%
for any $x\in H$ with $\left\Vert x\right\Vert =1$ and $p,q>0.$

Similar results can be stated if $p<0$ or $q<0.$ However the details are
left to the interest reader.

\bigskip

\chapter{Inequalities of Ostrowski Type}

\section{Introduction}

Ostrowski's type inequalities provide sharp error estimates in approximating
the value of a function by its integral mean. They can be utilized to obtain
a priory error bounds for different quadrature rules in approximating the
Riemann integral by different Riemann sums. They also shows, in general,
that the mid-point rule provides the best approximation in the class of all
Riemann sums sampled in the interior points of a given partition.

As revealed by a simple search in the data base \textit{MathSciNet} of the 
\textit{American Mathematical Society} with the key words "Ostrowski" and
"inequality" in the title, an exponential evolution of research papers
devoted to this result has been registered in the last decade. There are now
at least 280 papers that can be found by performing the above search.
Numerous extensions, generalisations in both the integral and discrete case
have been discovered. More general versions for $n$-time differentiable
functions, the corresponding versions on time scales, for vector valued
functions or multiple integrals have been established as well. Numerous
applications in Numerical Analysis, Probability Theory and other fields have
been also given.

In the present chapter we present some recent results obtained by the author
in extending Ostrowski inequality in various directions for continuous
functions of selfadjoint operators in complex Hilbert spaces. As far as we
know, the obtained results are new with no previous similar results ever
obtained in the literature.

Applications for mid-point inequalities and some elementary functions of
operators such as the power function, the logarithmic and exponential
functions are provided as well.

\section{Scalar Ostrowski's Type Inequalities}

In the scalar case, comparison between functions and integral means are
incorporated in Ostrowski type inequalities as mentioned below.

The first result in this direction is known in the literature as Ostrowski's
inequality \cite{III.a.O}.

\begin{theorem}
\label{III.b.ta}Let $f:\left[ a,b\right] \rightarrow \mathbb{R}$ be a
differentiable function on $\left( a,b\right) $ with the property that $%
\left\vert f^{\prime }\left( t\right) \right\vert \leq M$ for all $t\in
\left( a,b\right) $. Then 
\begin{equation}
\left\vert f\left( x\right) -\frac{1}{b-a}\int_{a}^{b}f\left( t\right)
dt\right\vert \leq \left[ \frac{1}{4}+\left( \frac{x-\frac{a+b}{2}}{b-a}%
\right) ^{2}\right] \left( b-a\right) M  \label{III.b.1.1}
\end{equation}%
for all $x\in \left[ a,b\right] $. The constant $\frac{1}{4}$ is the best
possible in the sense that it cannot be replaced by a smaller quantity..
\end{theorem}

The following Ostrowski type result for absolutely continuous functions
holds (see \cite{III.b.DRW1} -- \cite{III.b.DRW3}).

\begin{theorem}
\label{III.b.tb}Let $f:\left[ a,b\right] \rightarrow \mathbb{R}$ be
absolutely continuous on $\left[ a,b\right] $. Then, for all $x\in \left[ a,b%
\right] $, we have: 
\begin{multline}
\left\vert f\left( x\right) -\frac{1}{b-a}\int_{a}^{b}f\left( t\right)
dt\right\vert  \label{III.b.1.3} \\
\leq \left\{ 
\begin{array}{lll}
\left[ \frac{1}{4}+\left( \frac{x-\frac{a+b}{2}}{b-a}\right) ^{2}\right]
\left( b-a\right) \left\Vert f^{\prime }\right\Vert _{\infty } & \text{if} & 
f^{\prime }\in L_{\infty }\left[ a,b\right] ; \\ 
&  &  \\ 
\frac{1}{\left( p+1\right) ^{\frac{1}{p}}}\left[ \left( \frac{x-a}{b-a}%
\right) ^{p+1}+\left( \frac{b-x}{b-a}\right) ^{p+1}\right] ^{\frac{1}{p}%
}\left( b-a\right) ^{\frac{1}{p}}\left\Vert f^{\prime }\right\Vert _{q} & 
\text{if} & f^{\prime }\in L_{q}\left[ a,b\right] , \\ 
&  & \frac{1}{p}+\frac{1}{q}=1,\;p>1; \\ 
\left[ \frac{1}{2}+\left\vert \frac{x-\frac{a+b}{2}}{b-a}\right\vert \right]
\left\Vert f^{\prime }\right\Vert _{1}; &  & 
\end{array}%
\right.
\end{multline}%
where $\left\Vert \cdot \right\Vert _{r}$ \ ($r\in \left[ 1,\infty \right] $%
) are the usual Lebesgue norms on $L_{r}\left[ a,b\right] $, i.e., 
\begin{equation*}
\left\Vert g\right\Vert _{\infty }:=ess\sup\limits_{t\in \left[ a,b\right]
}\left\vert g\left( t\right) \right\vert
\end{equation*}%
and 
\begin{equation*}
\left\Vert g\right\Vert _{r}:=\left( \int_{a}^{b}\left\vert g\left( t\right)
\right\vert ^{r}dt\right) ^{\frac{1}{r}},\;r\in \lbrack 1,\infty ).
\end{equation*}%
The constants $\frac{1}{4}$, $\frac{1}{\left( p+1\right) ^{\frac{1}{p}}}$
and $\frac{1}{2}$ respectively are sharp in the sense presented in Theorem %
\ref{III.b.ta}.
\end{theorem}

The above inequalities can also be obtained from the Fink result in \cite%
{III.a.F} on choosing $n=1$ and performing some appropriate computations.

If one drops the condition of absolute continuity and assumes that $f$ is H%
\"{o}lder continuous, then one may state the result (see for instance \cite%
{III.b.DCRW} and the references therein for earlier contributions):

\begin{theorem}
\label{III.b.tc}Let $f:\left[ a,b\right] \rightarrow \mathbb{R}$ be of $r-H-$%
H\"{o}lder type, i.e., 
\begin{equation}
\left\vert f\left( x\right) -f\left( y\right) \right\vert \leq H\left\vert
x-y\right\vert ^{r},\;\text{for all \ }x,y\in \left[ a,b\right] ,
\label{III.b.1.4}
\end{equation}%
where $r\in (0,1]$ and $H>0$ are fixed. Then, for all $x\in \left[ a,b\right]
,$ we have the inequality: 
\begin{align}
& \left\vert f\left( x\right) -\frac{1}{b-a}\int_{a}^{b}f\left( t\right)
dt\right\vert  \label{III.b.1.5} \\
& \leq \frac{H}{r+1}\left[ \left( \frac{b-x}{b-a}\right) ^{r+1}+\left( \frac{%
x-a}{b-a}\right) ^{r+1}\right] \left( b-a\right) ^{r}.  \notag
\end{align}%
The constant $\frac{1}{r+1}$ is also sharp in the above sense.
\end{theorem}

Note that if $r=1$, i.e., $f$ is Lipschitz continuous, then we get the
following version of Ostrowski's inequality for Lipschitzian functions (with 
$L$ instead of $H$) (see for instance \cite{III.b.DR1}) 
\begin{equation}
\left\vert f\left( x\right) -\frac{1}{b-a}\int_{a}^{b}f\left( t\right)
dt\right\vert \leq \left[ \frac{1}{4}+\left( \frac{x-\frac{a+b}{2}}{b-a}%
\right) ^{2}\right] \left( b-a\right) L.  \label{III.b.1.6}
\end{equation}%
Here the constant $\frac{1}{4}$ is also best.

Moreover, if one drops the condition of the continuity of the function, and
assumes that it is of bounded variation, then the following result may be
stated (see \cite{III.b.DR2}).

\begin{theorem}
\label{III.b.td}Assume that $f:\left[ a,b\right] \rightarrow \mathbb{R}$ is
of bounded variation and denote by $\bigvee\limits_{a}^{b}\left( f\right) $
its total variation. Then 
\begin{equation}
\left\vert f\left( x\right) -\frac{1}{b-a}\int_{a}^{b}f\left( t\right)
dt\right\vert \leq \left[ \frac{1}{2}+\left\vert \frac{x-\frac{a+b}{2}}{b-a}%
\right\vert \right] \bigvee\limits_{a}^{b}\left( f\right)  \label{III.b.1.7}
\end{equation}%
for all $x\in \left[ a,b\right] $. The constant $\frac{1}{2}$ is the best
possible.
\end{theorem}

If we assume more about $f$, i.e., $f$ is monotonically increasing, then the
inequality (\ref{III.b.1.7}) may be improved in the following manner \cite%
{III.d.DR3} (see also the monograph \cite{III.b.DR}).

\begin{theorem}
\label{III.b.te}Let $f:\left[ a,b\right] \rightarrow \mathbb{R}$ be
monotonic nondecreasing. Then for all $x\in \left[ a,b\right] $, we have the
inequality: 
\begin{align}
& \left\vert f\left( x\right) -\frac{1}{b-a}\int_{a}^{b}f\left( t\right)
dt\right\vert  \label{III.b.1.8} \\
& \leq \frac{1}{b-a}\left\{ \left[ 2x-\left( a+b\right) \right] f\left(
x\right) +\int_{a}^{b}sgn\left( t-x\right) f\left( t\right) dt\right\} 
\notag \\
& \leq \frac{1}{b-a}\left\{ \left( x-a\right) \left[ f\left( x\right)
-f\left( a\right) \right] +\left( b-x\right) \left[ f\left( b\right)
-f\left( x\right) \right] \right\}  \notag \\
& \leq \left[ \frac{1}{2}+\left\vert \frac{x-\frac{a+b}{2}}{b-a}\right\vert %
\right] \left[ f\left( b\right) -f\left( a\right) \right] .  \notag
\end{align}%
All the inequalities in (\ref{III.b.1.8}) are sharp and the constant $\frac{1%
}{2}$ is the best possible.
\end{theorem}

For other scalar Ostrowski's type inequalities, see \cite{III.b.A}-\cite%
{III.b.CD} and \cite{III.b.7ab}.

\section{Ostrowski's type Inequalities for H\"{o}lder Continuous Functions}

\subsection{Introduction}

Let $U$ be a selfadjoint operator on the Hilbert space $\left(
H,\left\langle .,.\right\rangle \right) $ with the spectrum $Sp\left(
U\right) $ included in the interval $\left[ m,M\right] $ for some real
numbers $m<M$ and let $\left\{ E_{\lambda }\right\} _{\lambda \in \mathbb{R}%
} $ be its \textit{spectral family}. Then for any continuous function $f:%
\left[ m,M\right] \rightarrow \mathbb{C}$, it is well known that we have the
following \textit{spectral representation} theorem in terms of the\textit{\
Riemann-Stieltjes integral}: 
\begin{equation}
\left\langle f\left( U\right) x,x\right\rangle =\int_{m-0}^{M}f\left(
\lambda \right) d\left( \left\langle E_{\lambda }x,x\right\rangle \right) ,
\label{III.e.1.1}
\end{equation}%
for any $x\in H$ with $\left\Vert x\right\Vert =1.$ The function $%
g_{x}\left( \lambda \right) :=\left\langle E_{\lambda }x,x\right\rangle $ is 
\textit{monotonic nondecreasing} on the interval $\left[ m,M\right] $ and 
\begin{equation}
g_{x}\left( m-0\right) =0\text{ and }g_{x}\left( M\right) =1
\label{III.e.1.2}
\end{equation}%
for any $x\in H$ with $\left\Vert x\right\Vert =1.$

Utilising the representation (\ref{III.e.1.1}) and the following Ostrowski's
type inequality for the Riemann-Stieltjes integral obtained by the author in 
\cite{III.SSD1}: 
\begin{align}
& \left\vert f\left( s\right) \left[ u\left( b\right) -u\left( a\right) %
\right] -\int_{a}^{b}f\left( t\right) du\left( t\right) \right\vert
\label{III.e.1.3} \\
& \leq L\left[ \frac{1}{2}\left( b-a\right) +\left\vert s-\frac{a+b}{2}%
\right\vert \right] ^{r}\dbigvee\limits_{a}^{b}\left( u\right)  \notag
\end{align}%
for any $s\in \left[ a,b\right] ,$ provided that $f$ is of $r-L-$H\"{o}lder
type on $\left[ a,b\right] $ (see (\ref{III.e.1.4}) below), $u$ is of 
\textit{bounded variation} on $\left[ a,b\right] $ and $\dbigvee_{a}^{b}%
\left( u\right) $ denotes the \textit{total variation} of $u$ on $\left[ a,b%
\right] ,$ we obtained the following inequality of Ostrowski type for
selfadjoint operators:

\begin{theorem}[Dragomir, 2008, \protect\cite{III.SSD2}]
\label{III.t.1.1}Let $A$ and $B$ be selfadjoint operators with $Sp\left(
A\right) ,$ $Sp\left( B\right) \subseteq \left[ m,M\right] $ for some real
numbers $m<M.$ If $f:\left[ m,M\right] \longrightarrow \mathbb{R}$ is of $%
r-L-$H\"{o}lder type, i.e., for a given $r\in (0,1]$ and $L>0$ we have 
\begin{equation}
\left\vert f\left( s\right) -f\left( t\right) \right\vert \leq L\left\vert
s-t\right\vert ^{r}\text{ for any }s,t\in \left[ m,M\right] ,
\label{III.e.1.4}
\end{equation}%
then we have the inequality: 
\begin{equation}
\left\vert f\left( s\right) -\left\langle f\left( A\right) x,x\right\rangle
\right\vert \leq L\left[ \frac{1}{2}\left( M-m\right) +\left\vert s-\frac{m+M%
}{2}\right\vert \right] ^{r},  \label{III.e.1.5}
\end{equation}%
for any $s\in \left[ m,M\right] $ and any $x\in H$ with $\left\Vert
x\right\Vert =1$.

Moreover, we have 
\begin{align}
& \left\vert \left\langle f\left( B\right) y,y\right\rangle -\left\langle
f\left( A\right) x,x\right\rangle \right\vert  \label{III.e.1.6} \\
& \leq \left\langle \left\vert f\left( B\right) -\left\langle f\left(
A\right) x,x\right\rangle \cdot 1_{H}\right\vert y,y\right\rangle  \notag \\
& \leq L\left[ \frac{1}{2}\left( M-m\right) +\left\langle \left\vert B-\frac{%
m+M}{2}\cdot 1_{H}\right\vert y,y\right\rangle \right] ^{r},  \notag
\end{align}%
for any $x,y\in H$ with $\left\Vert x\right\Vert =\left\Vert y\right\Vert
=1. $
\end{theorem}

With the above assumptions for $f,A$ and $B$ we have the following
particular inequalities of interest: 
\begin{equation}
\left\vert f\left( \frac{m+M}{2}\right) -\left\langle f\left( A\right)
x,x\right\rangle \right\vert \leq \frac{1}{2^{r}}L\left( M-m\right) ^{r}
\label{III.e.1.7}
\end{equation}%
and 
\begin{equation}
\left\vert f\left( \left\langle Ax,x\right\rangle \right) -\left\langle
f\left( A\right) x,x\right\rangle \right\vert \leq L\left[ \frac{1}{2}\left(
M-m\right) +\left\vert \left\langle Ax,x\right\rangle -\frac{m+M}{2}%
\right\vert \right] ^{r},  \label{III.e.1.8}
\end{equation}%
for any $x\in H$ with $\left\Vert x\right\Vert =1$.

We also have the inequalities: 
\begin{align}
& \left\vert \left\langle f\left( A\right) y,y\right\rangle -\left\langle
f\left( A\right) x,x\right\rangle \right\vert  \label{III.e.1.9} \\
& \leq \left\langle \left\vert f\left( A\right) -\left\langle f\left(
A\right) x,x\right\rangle \cdot 1_{H}\right\vert y,y\right\rangle  \notag \\
& \leq L\left[ \frac{1}{2}\left( M-m\right) +\left\langle \left\vert A-\frac{%
m+M}{2}\cdot 1_{H}\right\vert y,y\right\rangle \right] ^{r},  \notag
\end{align}%
for any $x,y\in H$ with $\left\Vert x\right\Vert =\left\Vert y\right\Vert
=1, $

\begin{align}
& \left\vert \left\langle \left[ f\left( B\right) -f\left( A\right) \right]
x,x\right\rangle \right\vert  \label{III.e.1.10} \\
& \leq \left\langle \left\vert f\left( B\right) -\left\langle f\left(
A\right) x,x\right\rangle \cdot 1_{H}\right\vert x,x\right\rangle  \notag \\
& \leq L\left[ \frac{1}{2}\left( M-m\right) +\left\langle \left\vert B-\frac{%
m+M}{2}\cdot 1_{H}\right\vert x,x\right\rangle \right] ^{r}  \notag
\end{align}%
and, more particularly, 
\begin{align}
& \left\langle \left\vert f\left( A\right) -\left\langle f\left( A\right)
x,x\right\rangle \cdot 1_{H}\right\vert x,x\right\rangle  \label{III.e.1.11}
\\
& \leq L\left[ \frac{1}{2}\left( M-m\right) +\left\langle \left\vert A-\frac{%
m+M}{2}\cdot 1_{H}\right\vert x,x\right\rangle \right] ^{r},  \notag
\end{align}%
for any $x\in H$ with $\left\Vert x\right\Vert =1.$

We also have the norm inequality 
\begin{equation}
\left\Vert f\left( B\right) -f\left( A\right) \right\Vert \leq L\left[ \frac{%
1}{2}\left( M-m\right) +\left\Vert B-\frac{m+M}{2}\cdot 1_{H}\right\Vert %
\right] ^{r}.  \label{III.e.1.12}
\end{equation}

For various generalizations, extensions and related Ostrowski type
inequalities for functions of one or several variables see the monograph 
\cite{III.DR} and the references therein.

\subsection{More Inequalities of Ostrowski's Type}

The following result holds:

\begin{theorem}[Dragomir, 2010, \protect\cite{III.SSD3}]
\label{III.t.2.1}Let $A$ be a selfadjoint operator with $Sp\left( A\right)
\subseteq \left[ m,M\right] $ for some real numbers $m<M.$ If $f:\left[ m,M%
\right] \longrightarrow \mathbb{R}$ is of $r-L-$H\"{o}lder type with $r\in
(0,1]$, then we have the inequality: 
\begin{align}
\left\vert f\left( s\right) -\left\langle f\left( A\right) x,x\right\rangle
\right\vert & \leq L\left\langle \left\vert s\cdot 1_{H}-A\right\vert
x,x\right\rangle ^{r}  \label{III.e.2.1} \\
& \leq L\left[ \left( s-\left\langle Ax,x\right\rangle \right)
^{2}+D^{2}\left( A;x\right) \right] ^{r/2},  \notag
\end{align}%
for any $s\in \left[ m,M\right] $ and any $x\in H$ with $\left\Vert
x\right\Vert =1$, where $D\left( A;x\right) $ is the variance of the
selfadjoint operator $A$ in $x$ and is defined by 
\begin{equation*}
D\left( A;x\right) :=\left( \left\Vert Ax\right\Vert ^{2}-\left\langle
Ax,x\right\rangle ^{2}\right) ^{1/2},
\end{equation*}%
where $x\in H$ with $\left\Vert x\right\Vert =1.$
\end{theorem}

\begin{proof}
First of all, by the Jensen inequality for convex functions of selfadjoint
operators (see for instance \cite[p. 5]{III.FMPS}) applied for the modulus,
we can state that 
\begin{equation}
\left\vert \left\langle h\left( A\right) x,x\right\rangle \right\vert \leq
\left\langle \left\vert h\left( A\right) \right\vert x,x\right\rangle 
\tag{M}  \label{III.M}
\end{equation}%
for any $x\in H$ with $\left\Vert x\right\Vert =1,$ where $h$ is a
continuous function on $\left[ m,M\right] .$

Utilising the property (\ref{III.M}) we then get 
\begin{equation}
\left\vert f\left( s\right) -\left\langle f\left( A\right) x,x\right\rangle
\right\vert =\left\vert \left\langle f\left( s\right) \cdot 1_{H}-f\left(
A\right) x,x\right\rangle \right\vert \leq \left\langle \left\vert f\left(
s\right) \cdot 1_{H}-f\left( A\right) \right\vert x,x\right\rangle
\label{III.e.2.2}
\end{equation}%
for any $x\in H$ with $\left\Vert x\right\Vert =1$ and any $s\in \left[ m,M%
\right] .$

Since $f$ is of $r-L-$H\"{o}lder type$,$ then for any $t,s\in \left[ m,M%
\right] $ we have 
\begin{equation}
\left\vert f\left( s\right) -f\left( t\right) \right\vert \leq L\left\vert
s-t\right\vert ^{r}.  \label{III.e.2.3}
\end{equation}%
If we fix $s\in \left[ m,M\right] $ and apply the property (\ref{P}) for the
inequality (\ref{III.e.2.3}) and the operator $A$ we get 
\begin{equation}
\left\langle \left\vert f\left( s\right) \cdot 1_{H}-f\left( A\right)
\right\vert x,x\right\rangle \leq L\left\langle \left\vert s\cdot
1_{H}-A\right\vert ^{r}x,x\right\rangle \leq L\left\langle \left\vert s\cdot
1_{H}-A\right\vert x,x\right\rangle ^{r}  \label{III.e.2.4}
\end{equation}%
for any $x\in H$ with $\left\Vert x\right\Vert =1$ and any $s\in \left[ m,M%
\right] ,$ where, for the last inequality we have used the fact that if $P$
is a positive operator and $r\in \left( 0,1\right) $ then, by the H\"{o}%
lder-McCarthy inequality \cite{III.Mc}, 
\begin{equation}
\left\langle P^{r}x,x\right\rangle \leq \left\langle Px,x\right\rangle ^{r} 
\tag{HM}  \label{III.HM}
\end{equation}%
for any $x\in H$ with $\left\Vert x\right\Vert =1.$ This proves the fist
inequality in (\ref{III.e.2.1}).

Now, observe that for any bounded linear operator $T$ we have 
\begin{equation*}
\left\langle \left\vert T\right\vert x,x\right\rangle =\left\langle \left(
T^{\ast }T\right) ^{1/2}x,x\right\rangle \leq \left\langle \left( T^{\ast
}T\right) x,x\right\rangle ^{1/2}=\left\Vert Tx\right\Vert
\end{equation*}%
for any $x\in H$ with $\left\Vert x\right\Vert =1$ which implies that 
\begin{align}
\left\langle \left\vert s\cdot 1_{H}-A\right\vert x,x\right\rangle ^{r}&
\leq \left\Vert sx-Ax\right\Vert ^{r}  \label{III.e.2.5} \\
& =\left( s^{2}-2s\left\langle Ax,x\right\rangle +\left\Vert Ax\right\Vert
^{2}\right) ^{r/2}  \notag \\
& =\left[ \left( s-\left\langle Ax,x\right\rangle \right) ^{2}+\left\Vert
Ax\right\Vert ^{2}-\left\langle Ax,x\right\rangle ^{2}\right] ^{r/2}  \notag
\end{align}%
for any $x\in H$ with $\left\Vert x\right\Vert =1$ and any $s\in \left[ m,M%
\right] .$

Finally, on making use of (\ref{III.e.2.2}), (\ref{III.e.2.4}) and (\ref%
{III.e.2.5}) we deduce the desired result (\ref{III.e.2.1}).
\end{proof}

\begin{remark}
\label{III.r.2.0}If we choose in (\ref{III.e.2.1}) $s=\frac{m+M}{2},$ then
we get the sequence of inequalities 
\begin{align}
& \left\vert f\left( \frac{m+M}{2}\right) -\left\langle f\left( A\right)
x,x\right\rangle \right\vert  \label{III.e.2.5.a} \\
& \leq L\left\langle \left\vert \frac{m+M}{2}\cdot 1_{H}-A\right\vert
x,x\right\rangle ^{r}  \notag \\
& \leq L\left[ \left( \frac{m+M}{2}-\left\langle Ax,x\right\rangle \right)
^{2}+D^{2}\left( A;x\right) \right] ^{r/2}  \notag \\
& \leq L\left[ \frac{1}{4}\left( M-m\right) ^{2}+D^{2}\left( A;x\right) %
\right] ^{r/2}\leq \frac{1}{2^{r}}L\left( M-m\right) ^{r}  \notag
\end{align}%
for any $x\in H$ with $\left\Vert x\right\Vert =1,$ since, obviously, 
\begin{equation*}
\left( \frac{m+M}{2}-\left\langle Ax,x\right\rangle \right) ^{2}\leq \frac{1%
}{4}\left( M-m\right) ^{2}
\end{equation*}%
and 
\begin{equation*}
D^{2}\left( A;x\right) \leq \frac{1}{4}\left( M-m\right) ^{2}
\end{equation*}%
for any $x\in H$ with $\left\Vert x\right\Vert =1.$

We notice that the inequality (\ref{III.e.2.5.a}) provides a refinement for
the result (\ref{III.e.1.7}) above.
\end{remark}

The best inequality we can get from (\ref{III.e.2.1}) is incorporated in the
following:

\begin{corollary}[Dragomir, 2010, \protect\cite{III.SSD3}]
\label{III.c.2.1}Let $A$ be a selfadjoint operator with $Sp\left( A\right)
\subseteq \left[ m,M\right] $ for some real numbers $m<M.$ If $f:\left[ m,M%
\right] \longrightarrow \mathbb{R}$ is of $r-L-$H\"{o}lder type with $r\in
(0,1]$, then we have the inequality 
\begin{equation}
\left\vert f\left( \left\langle Ax,x\right\rangle \right) -\left\langle
f\left( A\right) x,x\right\rangle \right\vert \leq L\left\langle \left\vert
\left\langle Ax,x\right\rangle \cdot 1_{H}-A\right\vert x,x\right\rangle
^{r}\leq LD^{r}\left( A;x\right) ,  \label{III.e.2.6}
\end{equation}%
for any $x\in H$ with $\left\Vert x\right\Vert =1$.
\end{corollary}

The inequality (\ref{III.e.2.1}) may be used to obtain other inequalities
for two selfadjoint operators as follows:

\begin{corollary}[Dragomir, 2010, \protect\cite{III.SSD3}]
\label{III.c.2.2}Let $A$ and $B$ be selfadjoint operators with $Sp\left(
A\right) ,Sp\left( B\right) \subseteq \left[ m,M\right] $ for some real
numbers $m<M.$ If $f:\left[ m,M\right] \longrightarrow \mathbb{R}$ is of $%
r-L-$H\"{o}lder type with $r\in (0,1]$, then we have the inequality 
\begin{align}
& \left\vert \left\langle f\left( B\right) y,y\right\rangle -\left\langle
f\left( A\right) x,x\right\rangle \right\vert  \label{III.e.2.7} \\
& \leq L\left[ \left( \left\langle By,y\right\rangle -\left\langle
Ax,x\right\rangle \right) ^{2}+D^{2}\left( A;x\right) +D^{2}\left(
B;y\right) \right] ^{r/2}  \notag
\end{align}%
for any $x,y\in H$ with $\left\Vert x\right\Vert =\left\Vert y\right\Vert
=1. $
\end{corollary}

\begin{proof}
If we apply the property (\ref{P}) to the inequality (\ref{III.e.2.1}) and
for the operator $B,$ then we get 
\begin{align}
& \left\langle \left\vert f\left( B\right) -\left\langle f\left( A\right)
x,x\right\rangle \cdot 1_{H}\right\vert y,y\right\rangle  \label{III.e.2.8}
\\
& \leq L\left\langle \left[ \left( B-\left\langle Ax,x\right\rangle \cdot
1_{H}\right) ^{2}+D^{2}\left( A;x\right) \cdot 1_{H}\right]
^{r/2}y,y\right\rangle  \notag
\end{align}%
for any $x,y\in H$ with $\left\Vert x\right\Vert =\left\Vert y\right\Vert
=1. $

Utilising the inequality (\ref{III.M}) we also have that 
\begin{equation}
\left\vert f\left( \left\langle By,y\right\rangle \right) -\left\langle
f\left( A\right) x,x\right\rangle \right\vert \leq \left\langle \left\vert
f\left( B\right) -\left\langle f\left( A\right) x,x\right\rangle \cdot
1_{H}\right\vert y,y\right\rangle  \label{III.e.2.9}
\end{equation}%
for any $x,y\in H$ with $\left\Vert x\right\Vert =\left\Vert y\right\Vert
=1. $

Now, by the H\"{o}lder-McCarthy inequality (\ref{III.HM}) we also have 
\begin{align}
& \left\langle \left[ \left( B-\left\langle Ax,x\right\rangle \cdot
1_{H}\right) ^{2}+D^{2}\left( A;x\right) \cdot 1_{H}\right]
^{r/2}y,y\right\rangle  \label{III.e.2.10} \\
& \leq \left\langle \left[ \left( B-\left\langle Ax,x\right\rangle \cdot
1_{H}\right) ^{2}+D^{2}\left( A;x\right) \cdot 1_{H}\right] y,y\right\rangle
^{r/2}  \notag \\
& =\left( \left( \left\langle By,y\right\rangle -\left\langle
Ax,x\right\rangle \right) ^{2}+D^{2}\left( A;x\right) +D^{2}\left(
B;y\right) \right) ^{r/2}  \notag
\end{align}%
for any $x,y\in H$ with $\left\Vert x\right\Vert =\left\Vert y\right\Vert
=1. $

On making use of (\ref{III.e.2.8})-(\ref{III.e.2.10}) we deduce the desired
result (\ref{III.e.2.7}).
\end{proof}

\begin{remark}
\label{III.r.2.1} Since 
\begin{equation}
D^{2}\left( A;x\right) \leq \frac{1}{4}\left( M-m\right) ^{2},
\label{III.e.2.11}
\end{equation}%
then we obtain from (\ref{III.e.2.7}) the following vector inequalities 
\begin{align}
& \left\vert \left\langle f\left( A\right) y,y\right\rangle -\left\langle
f\left( A\right) x,x\right\rangle \right\vert  \label{III.e.2.12} \\
& \leq L\left[ \left( \left\langle Ay,y\right\rangle -\left\langle
Ax,x\right\rangle \right) ^{2}+D^{2}\left( A;x\right) +D^{2}\left(
A;y\right) \right] ^{r/2}  \notag \\
& \leq L\left[ \left( \left\langle Ay,y\right\rangle -\left\langle
Ax,x\right\rangle \right) ^{2}+\frac{1}{2}\left( M-m\right) ^{2}\right]
^{r/2},  \notag
\end{align}%
and 
\begin{align}
& \left\vert \left\langle \left[ f\left( B\right) -f\left( A\right) \right]
x,x\right\rangle \right\vert  \label{III.e.2.13} \\
& \leq L\left[ \left\langle \left( B-A\right) x,x\right\rangle
^{2}+D^{2}\left( A;x\right) +D^{2}\left( B;x\right) \right] ^{r/2}  \notag \\
& \leq L\left[ \left\langle \left( B-A\right) x,x\right\rangle ^{2}+\frac{1}{%
2}\left( M-m\right) ^{2}\right] ^{r/2}.  \notag
\end{align}%
In particular, we have the norm inequality 
\begin{equation}
\left\Vert f\left( B\right) -f\left( A\right) \right\Vert \leq L\left[
\left\Vert B-A\right\Vert ^{2}+\frac{1}{2}\left( M-m\right) ^{2}\right]
^{r/2}.  \label{III.e.2.14}
\end{equation}
\end{remark}

The following result provides convenient examples for applications:

\begin{corollary}[Dragomir, 2010, \protect\cite{III.SSD3}]
\label{III.c.2.3}Let $A$ be a selfadjoint operator with $Sp\left( A\right)
\subseteq \left[ m,M\right] $ for some real numbers $m<M.$ If $f:\left[ m,M%
\right] \longrightarrow \mathbb{R}$ is absolutely continuous on $\left[ m,M%
\right] $, then we have the inequality: 
\begin{align}
& \left\vert f\left( s\right) -\left\langle f\left( A\right)
x,x\right\rangle \right\vert  \label{III.e.2.15} \\
& \leq \left\{ 
\begin{array}{cc}
\left\langle \left\vert s\cdot 1_{H}-A\right\vert x,x\right\rangle
\left\Vert f^{\prime }\right\Vert _{\left[ m,M\right] ,\infty } & \text{if }%
f^{\prime }\in L_{\infty }\left[ m,M\right] , \\ 
&  \\ 
\left\langle \left\vert s\cdot 1_{H}-A\right\vert x,x\right\rangle
^{1/q}\left\Vert f^{\prime }\right\Vert _{\left[ m,M\right] ,p} & 
\begin{array}{c}
\text{if }f^{\prime }\in L_{p}\left[ m,M\right] , \\ 
p>1,\frac{1}{p}+\frac{1}{q}=1,%
\end{array}%
\end{array}%
\right.  \notag \\
& \leq \left\{ 
\begin{array}{cc}
\left[ \left( s-\left\langle Ax,x\right\rangle \right) ^{2}+D^{2}\left(
A;x\right) \right] ^{1/2}\left\Vert f^{\prime }\right\Vert _{\left[ m,M%
\right] ,\infty } & \text{if }f^{\prime }\in L_{\infty }\left[ m,M\right] ,
\\ 
&  \\ 
\left[ \left( s-\left\langle Ax,x\right\rangle \right) ^{2}+D^{2}\left(
A;x\right) \right] ^{\frac{1}{2q}}\left\Vert f^{\prime }\right\Vert _{\left[
m,M\right] ,p} & 
\begin{array}{c}
\text{if }f^{\prime }\in L_{p}\left[ m,M\right] , \\ 
p>1,\frac{1}{p}+\frac{1}{q}=1,%
\end{array}%
\end{array}%
\right.  \notag
\end{align}%
for any $s\in \left[ m,M\right] $ and any $x\in H$ with $\left\Vert
x\right\Vert =1,$ where $\left\Vert f^{\prime }\right\Vert _{\left[ m,M%
\right] ,\ell }$ are the Lebesgue norms, i.e., 
\begin{equation*}
\left\Vert f^{\prime }\right\Vert _{\left[ m,M\right] ,\ell }:=\left\{ 
\begin{array}{cc}
ess\sup_{t\in \left[ m,M\right] }\left\vert f^{\prime }\left( t\right)
\right\vert & \text{if }\ell =\infty \\ 
&  \\ 
\left( \int_{m}^{M}\left\vert f^{\prime }\left( t\right) \right\vert
^{p}dt\right) ^{1/p} & \text{if }\ell =p\geq 1.%
\end{array}%
\right.
\end{equation*}
\end{corollary}

\begin{proof}
Follows from Theorem \ref{III.t.2.1} and on tacking into account that if $f:%
\left[ m,M\right] \longrightarrow \mathbb{R}$ is absolutely continuous on $%
\left[ m,M\right] ,$ then for any $s,t\in \left[ m,M\right] $ we have 
\begin{align*}
& \left\vert f\left( s\right) -f\left( t\right) \right\vert \\
& =\left\vert \int_{t}^{s}f^{\prime }\left( u\right) du\right\vert \\
& \leq \left\{ 
\begin{array}{cc}
\left\vert s-t\right\vert ess\sup_{t\in \left[ m,M\right] }\left\vert
f^{\prime }\left( t\right) \right\vert & \text{if }f^{\prime }\in L_{\infty }%
\left[ m,M\right] \\ 
&  \\ 
\left\vert s-t\right\vert ^{1/q}\left( \int_{m}^{M}\left\vert f^{\prime
}\left( t\right) \right\vert ^{p}dt\right) ^{1/p} & \text{if }f^{\prime }\in
L_{p}\left[ m,M\right] ,p>1,\frac{1}{p}+\frac{1}{q}=1.%
\end{array}%
\right.
\end{align*}
\end{proof}

\begin{remark}
\label{III.r.2.2}It is clear that all the inequalities from Corollaries \ref%
{III.c.2.1}, \ref{III.c.2.2} and Remark \ref{III.r.2.1} may be stated for
absolutely continuous functions. However, we mention here only one, namely 
\begin{align}
& \left\vert f\left( \left\langle Ax,x\right\rangle \right) -\left\langle
f\left( A\right) x,x\right\rangle \right\vert  \label{III.e.2.16} \\
& \leq \left\{ 
\begin{array}{cc}
\left\langle \left\vert \left\langle Ax,x\right\rangle \cdot
1_{H}-A\right\vert x,x\right\rangle \left\Vert f^{\prime }\right\Vert _{ 
\left[ m,M\right] ,\infty } & \text{if }f^{\prime }\in L_{\infty }\left[ m,M%
\right] \\ 
&  \\ 
\left\langle \left\vert \left\langle Ax,x\right\rangle \cdot
1_{H}-A\right\vert x,x\right\rangle ^{1/q}\left\Vert f^{\prime }\right\Vert
_{\left[ m,M\right] ,p} & 
\begin{array}{c}
\text{if }f^{\prime }\in L_{p}\left[ m,M\right] , \\ 
p>1,\frac{1}{p}+\frac{1}{q}=1,%
\end{array}%
\end{array}%
\right.  \notag \\
& \leq \left\{ 
\begin{array}{cc}
D\left( A;x\right) \left\Vert f^{\prime }\right\Vert _{\left[ m,M\right]
,\infty } & \text{if }f^{\prime }\in L_{\infty }\left[ m,M\right] \\ 
&  \\ 
D^{1/q}\left( A;x\right) \left\Vert f^{\prime }\right\Vert _{\left[ m,M%
\right] ,p} & 
\begin{array}{c}
\text{if }f^{\prime }\in L_{p}\left[ m,M\right] , \\ 
p>1,\frac{1}{p}+\frac{1}{q}=1.%
\end{array}%
\end{array}%
\right.  \notag
\end{align}
\end{remark}

\subsection{The Case of $\left( \protect\varphi ,\Phi \right) -$Lipschitzian
Functions}

The following lemma may be stated.

\begin{lemma}
\label{III.l.4.1}Let $u:\left[ a,b\right] \rightarrow \mathbb{R}$ and $%
\varphi ,\Phi \in \mathbb{R}$ be such that $\Phi >\varphi .$ The following
statements are equivalent:

\begin{enumerate}
\item[(i)] The function $u-\frac{\varphi +\Phi }{2}\cdot e,$ where $e\left(
t\right) =t,$ $t\in \left[ a,b\right] ,$ is $\frac{1}{2}\left( \Phi -\varphi
\right) -$Lipschitzian;

\item[(ii)] We have the inequality: 
\begin{equation}
\varphi \leq \frac{u\left( t\right) -u\left( s\right) }{t-s}\leq \Phi \quad 
\text{for each}\quad t,s\in \left[ a,b\right] \quad \text{with }t\neq s;
\label{III.e.4.1}
\end{equation}

\item[(iii)] We have the inequality: 
\begin{equation}
\varphi \left( t-s\right) \leq u\left( t\right) -u\left( s\right) \leq \Phi
\left( t-s\right) \quad \text{for each}\quad t,s\in \left[ a,b\right] \quad 
\text{with }t>s.  \label{III.e.4.2}
\end{equation}
\end{enumerate}
\end{lemma}

We can introduce the following class of functions, see also \cite{III.L}:

\begin{definition}
\label{III.d.4.1}The function $u:\left[ a,b\right] \rightarrow \mathbb{R}$
which satisfies one of the equivalent conditions (i) -- (iii) is said to be $%
\left( \varphi ,\Phi \right) -$Lipschitzian on $\left[ a,b\right] .$
\end{definition}

Utilising \textit{Lagrange's mean value theorem}, we can state the following
result that provides practical examples of $\left( \varphi ,\Phi \right) -$%
Lipschitzian functions.

\begin{proposition}
\label{III.p.4.1}Let $u:\left[ a,b\right] \rightarrow \mathbb{R}$ be
continuous on $\left[ a,b\right] $ and differentiable on $\left( a,b\right)
. $ If 
\begin{equation}
-\infty <\gamma :=\inf_{t\in \left( a,b\right) }u^{\prime }\left( t\right)
,\qquad \sup_{t\in \left( a,b\right) }u^{\prime }\left( t\right) =:\Gamma
<\infty  \label{III.e.4.3}
\end{equation}%
then $u$ is $\left( \gamma ,\Gamma \right) -$Lipschitzian on $\left[ a,b%
\right] .$
\end{proposition}

The following result can be stated:

\begin{proposition}[Dragomir, 2010, \protect\cite{III.SSD3}]
\label{III.p.4.2}Let $A$ be a selfadjoint operator with $Sp\left( A\right)
\subseteq \left[ m,M\right] $ for some real numbers $m<M.$ If $f:\left[ m,M%
\right] \longrightarrow \mathbb{R}$ is $\left( \gamma ,\Gamma \right) -$%
Lipschitzian on $\left[ m,M\right] ,$ then we have the inequality 
\begin{align}
\left\vert f\left( \left\langle Ax,x\right\rangle \right) -\left\langle
f\left( A\right) x,x\right\rangle \right\vert & \leq \frac{1}{2}\left(
\Gamma -\gamma \right) \left\langle \left\vert \left\langle
Ax,x\right\rangle \cdot 1_{H}-A\right\vert x,x\right\rangle
\label{III.e.4.4} \\
& \leq \frac{1}{2}\left( \Gamma -\gamma \right) D\left( A;x\right) ,  \notag
\end{align}%
for any $x\in H$ with $\left\Vert x\right\Vert =1$.
\end{proposition}

\begin{proof}
Follows by Corollary \ref{III.c.2.1} on taking into account that in this
case we have $r=1$ and $L=\frac{1}{2}\left( \Gamma -\gamma \right) .$
\end{proof}

We can use the result (\ref{III.e.4.4}) for the particular case of convex
functions to provide an interesting reverse inequality for the Jensen's type
operator inequality due to Mond and Pe\v{c}ari\'{c} \cite{III.MP1} (see also 
\cite[p. 5]{III.FMPS}):

\begin{theorem}[Mond-Pe\v{c}ari\'{c}, 1993, \protect\cite{III.MP1}]
\label{III.t.4.1} Let $A$ be a selfadjoint operator on the Hilbert space $H$
and assume that $Sp\left( A\right) \subseteq \left[ m,M\right] $ for some
scalars $m,M$ with $m<M.$ If $f$ is a convex function on $\left[ m,M\right]
, $ then 
\begin{equation}
f\left( \left\langle Ax,x\right\rangle \right) \leq \left\langle f\left(
A\right) x,x\right\rangle  \tag{MP}  \label{III.MPI}
\end{equation}%
for each $x\in H$ with $\left\Vert x\right\Vert =1.$
\end{theorem}

\begin{corollary}[Dragomir, 2010, \protect\cite{III.SSD3}]
\label{III.c.4.1}With the assumptions of Theorem \ref{III.t.4.1} we have the
inequality%
\begin{align}
& \left( 0\leq \right) \left\langle f\left( A\right) x,x\right\rangle
-f\left( \left\langle Ax,x\right\rangle \right)  \label{III.e.4.5} \\
& \leq \frac{1}{2}\left( f_{-}^{\prime }\left( M\right) -f_{+}^{\prime
}\left( m\right) \right) \left\langle \left\vert \left\langle
Ax,x\right\rangle \cdot 1_{H}-A\right\vert x,x\right\rangle  \notag \\
& \leq \frac{1}{2}\left( f_{-}^{\prime }\left( M\right) -f_{+}^{\prime
}\left( m\right) \right) D\left( A;x\right) \leq \frac{1}{4}\left(
f_{-}^{\prime }\left( M\right) -f_{+}^{\prime }\left( m\right) \right)
\left( M-m\right)  \notag
\end{align}%
for each $x\in H$ with $\left\Vert x\right\Vert =1.$
\end{corollary}

\begin{proof}
Follows by Proposition \ref{III.p.4.2} on taking into account that 
\begin{equation*}
f_{+}^{\prime }\left( m\right) \left( t-s\right) \leq f\left( t\right)
-f\left( s\right) \leq f_{-}^{\prime }\left( M\right) \left( t-s\right) \quad
\end{equation*}%
for each $s,t$ with the property that $M>$ $t>s>m.$
\end{proof}

The following result may be stated as well:

\begin{proposition}[Dragomir, 2010, \protect\cite{III.SSD3}]
\label{III.p.4.5}Let $A$ be a selfadjoint operator with $Sp\left( A\right)
\subseteq \left[ m,M\right] $ for some real numbers $m<M.$ If $f:\left[ m,M%
\right] \longrightarrow \mathbb{R}$ is $\left( \gamma ,\Gamma \right) -$%
Lipschitzian on $\left[ m,M\right] ,$ then we have the inequality 
\begin{align}
& \left\vert f\left( \left\langle Ax,x\right\rangle \right) -\left\langle
f\left( A\right) x,x\right\rangle \right\vert  \label{III.e.4.6} \\
& \leq \frac{1}{2}\left( \Gamma -\gamma \right) \left[ \frac{1}{2}\left(
M-m\right) +\left\vert \left\langle Ax,x\right\rangle -\frac{m+M}{2}%
\right\vert \right]  \notag
\end{align}%
for any $x\in H$ with $\left\Vert x\right\Vert =1$.
\end{proposition}

The following particular case for convex functions holds:

\begin{corollary}[Dragomir, 2010, \protect\cite{III.SSD3}]
\label{III.c.4.2}With the assumptions of Theorem \ref{III.t.4.1} we have the
inequality%
\begin{align}
& \left( 0\leq \right) \left\langle f\left( A\right) x,x\right\rangle
-f\left( \left\langle Ax,x\right\rangle \right)  \label{III.e.4.7} \\
& \leq \frac{1}{2}\left( f_{-}^{\prime }\left( M\right) -f_{+}^{\prime
}\left( m\right) \right) \left[ \frac{1}{2}\left( M-m\right) +\left\vert
\left\langle Ax,x\right\rangle -\frac{m+M}{2}\right\vert \right]  \notag
\end{align}%
for each $x\in H$ with $\left\Vert x\right\Vert =1.$
\end{corollary}

\subsection{Related Results}

In the previous sections we have compared amongst other the following
quantities 
\begin{equation*}
f\left( \frac{m+M}{2}\right) \text{ and }f\left( \left\langle
Ax,x\right\rangle \right) \text{ }
\end{equation*}%
with $\left\langle f\left( A\right) x,x\right\rangle $ for a selfadjoint
operator $A$ on the Hilbert space $H$ with $Sp\left( A\right) \subseteq %
\left[ m,M\right] $ for some real numbers $m<M,$ $f:\left[ m,M\right]
\longrightarrow \mathbb{R}$ a function of $r-L-$H\"{o}lder type with $r\in
(0,1]$ and $x\in H$ with $\left\Vert x\right\Vert =1.$

Since, obviously,%
\begin{equation*}
m\leq \frac{1}{M-m}\int_{m}^{M}f\left( t\right) dt\leq M,
\end{equation*}%
then is also natural to compare $\frac{1}{M-m}\int_{m}^{M}f\left( t\right)
dt $ with $\left\langle f\left( A\right) x,x\right\rangle $ under the same
assumptions for $f,A$ and $x.$

The following result holds:

\begin{theorem}[Dragomir, 2010, \protect\cite{III.SSD3}]
\label{III.t.5.1}Let $A$ be a selfadjoint operator with $Sp\left( A\right)
\subseteq \left[ m,M\right] $ for some real numbers $m<M.$ If $f:\left[ m,M%
\right] \longrightarrow \mathbb{R}$ is of $r-L-$H\"{o}lder type with $r\in
(0,1]$, then we have the inequality: 
\begin{align}
& \left\vert \frac{1}{M-m}\int_{m}^{M}f\left( s\right) dt-\left\langle
f\left( A\right) x,x\right\rangle \right\vert  \label{III.e.5.1} \\
& \leq \frac{1}{r+1}L\left( M-m\right) ^{r}  \notag \\
& \times \left[ \left\langle \left( \frac{M\cdot 1_{H}-A}{M-m}\right)
^{r+1}x,x\right\rangle +\left\langle \left( \frac{A-m\cdot 1_{H}}{M-m}%
\right) ^{r+1}x,x\right\rangle \right]  \notag \\
& \leq \frac{1}{r+1}L\left( M-m\right) ^{r},  \notag
\end{align}%
for any $x\in H$ with $\left\Vert x\right\Vert =1.$

In particular, if $f:\left[ m,M\right] \longrightarrow \mathbb{R}$ is
Lipschitzian with a constant $K,$ then 
\begin{align}
& \left\vert \frac{1}{M-m}\int_{m}^{M}f\left( s\right) dt-\left\langle
f\left( A\right) x,x\right\rangle \right\vert  \label{III.e.5.1.a} \\
& \leq K\left( M-m\right) \left[ \frac{1}{4}+\frac{1}{\left( M-m\right) ^{2}}%
\left( D^{2}\left( A;x\right) +\left( \left\langle Ax,x\right\rangle -\frac{%
m+M}{2}\right) ^{2}\right) \right]  \notag \\
& \leq \frac{1}{2}K\left( M-m\right)  \notag
\end{align}%
for any $x\in H$ with $\left\Vert x\right\Vert =1.$
\end{theorem}

\begin{proof}
We use the following Ostrowski's type result (see for instance \cite[p. 3]%
{III.DR}) written for the function $f$ that is of $r-L-$H\"{o}lder type on
the interval $\left[ m,M\right] :$%
\begin{align}
& \left\vert \frac{1}{M-m}\int_{m}^{M}f\left( s\right) dt-f\left( t\right)
\right\vert  \label{III.e.5.2} \\
& \leq \frac{L}{r+1}\left( M-m\right) ^{r}\left[ \left( \frac{M-t}{M-m}%
\right) ^{r+1}+\left( \frac{t-m}{M-m}\right) ^{r+1}\right]  \notag
\end{align}%
for any $t\in \left[ m,M\right] .$

If we apply the properties (\ref{P}) and (\ref{III.M}) then we have
successively%
\begin{align}
& \left\vert \frac{1}{M-m}\int_{m}^{M}f\left( s\right) dt-\left\langle
f\left( A\right) x,x\right\rangle \right\vert  \label{III.e.5.3} \\
& \leq \left\langle \left\vert \frac{1}{M-m}\int_{m}^{M}f\left( s\right)
dt-f\left( A\right) \right\vert x,x\right\rangle  \notag \\
& \leq \frac{L}{r+1}\left( M-m\right) ^{r}  \notag \\
& \times \left[ \left\langle \left( \frac{M\cdot 1_{H}-A}{M-m}\right)
^{r+1}x,x\right\rangle +\left\langle \left( \frac{A-m\cdot 1_{H}}{M-m}%
\right) ^{r+1}x,x\right\rangle \right]  \notag
\end{align}%
which proves the first inequality in (\ref{III.e.5.1}).

Utilising the Lah-Ribari\'{c} inequality version for selfadjoint operators $%
A $ with $Sp\left( A\right) \subseteq \left[ m,M\right] $ for some real
numbers $m<M$ and convex functions $g:\left[ m,M\right] \rightarrow \mathbb{R%
},$ namely (see for instance \cite[p. 57]{III.FMPS}): 
\begin{equation*}
\left\langle g\left( A\right) x,x\right\rangle \leq \frac{M-\left\langle
Ax,x\right\rangle }{M-m}g\left( m\right) +\frac{\left\langle
Ax,x\right\rangle -m}{M-m}g\left( M\right)
\end{equation*}%
for any $x\in H$ with $\left\Vert x\right\Vert =1,$ then we get for the
convex function $g\left( t\right) :=\left( \frac{M-t}{M-m}\right) ^{r+1},$%
\begin{equation*}
\left\langle \left( \frac{M\cdot 1_{H}-A}{M-m}\right) ^{r+1}x,x\right\rangle
\leq \frac{M-\left\langle Ax,x\right\rangle }{M-m}
\end{equation*}%
and for the convex function $g\left( t\right) :=\left( \frac{t-m}{M-m}%
\right) ^{r+1},$%
\begin{equation*}
\left\langle \left( \frac{A-m\cdot 1_{H}}{M-m}\right) ^{r+1}x,x\right\rangle
\leq \frac{\left\langle Ax,x\right\rangle -m}{M-m}
\end{equation*}%
for any $x\in H$ with $\left\Vert x\right\Vert =1.$

Now, on making use of the\ last two inequalities, we deduce the second part
of (\ref{III.e.5.1}).

Since%
\begin{eqnarray*}
&&\frac{1}{2}\left\langle \left( \frac{M\cdot 1_{H}-A}{M-m}\right)
^{2}x,x\right\rangle +\left\langle \left( \frac{A-m\cdot 1_{H}}{M-m}\right)
^{2}x,x\right\rangle \\
&=&\frac{1}{4}+\frac{1}{\left( M-m\right) ^{2}}\left( D^{2}\left( A;x\right)
+\left( \left\langle Ax,x\right\rangle -\frac{m+M}{2}\right) ^{2}\right)
\end{eqnarray*}%
for any $x\in H$ with $\left\Vert x\right\Vert =1,$ then on choosing $r=1$
in (\ref{III.e.5.1}) we deduce the desired result (\ref{III.e.5.1.a}).
\end{proof}

\begin{remark}
\label{III.r.2.3}We should notice from the proof of the above theorem, we
also have the following inequalities in the operator order of $B\left(
H\right) $%
\begin{align}
& \left\vert f\left( A\right) -\left( \frac{1}{M-m}\int_{m}^{M}f\left(
s\right) dt\right) \cdot 1_{H}\right\vert  \label{III.5.3.a} \\
& \leq \frac{L}{r+1}\left( M-m\right) ^{r}\left[ \left( \frac{M\cdot 1_{H}-A%
}{M-m}\right) ^{r+1}+\left( \frac{A-m\cdot 1_{H}}{M-m}\right) ^{r+1}\right] 
\notag \\
& \leq \frac{1}{r+1}L\left( M-m\right) ^{r}\cdot 1_{H}.  \notag
\end{align}
\end{remark}

The following particular case is of interest:

\begin{corollary}[Dragomir, 2010, \protect\cite{III.SSD3}]
\label{III.c.5.1}Let $A$ be a selfadjoint operator with $Sp\left( A\right)
\subseteq \left[ m,M\right] $ for some real numbers $m<M.$ If $f:\left[ m,M%
\right] \longrightarrow \mathbb{R}$ is $\left( \gamma ,\Gamma \right) -$%
Lipschitzian on $\left[ m,M\right] ,$ then we have the inequality 
\begin{align}
& \left\vert \left\langle f\left( A\right) x,x\right\rangle -\frac{\Gamma
+\gamma }{2}-\frac{1}{M-m}\int_{m}^{M}f\left( s\right) dt+\frac{\Gamma
+\gamma }{2}\cdot \frac{m+M}{2}\right\vert  \label{III.e.5.4} \\
& \leq \frac{1}{2}\left( \Gamma -\gamma \right) \left( M-m\right)  \notag \\
& \times \left[ \frac{1}{4}+\frac{1}{\left( M-m\right) ^{2}}\left(
D^{2}\left( A;x\right) +\left( \left\langle Ax,x\right\rangle -\frac{m+M}{2}%
\right) ^{2}\right) \right]  \notag \\
& \leq \frac{1}{4}\left( \Gamma -\gamma \right) \left( M-m\right) .  \notag
\end{align}
\end{corollary}

\begin{proof}
Follows by (\ref{III.e.5.1.a}) applied for the $\frac{1}{2}\left( \Gamma
-\gamma \right) $-Lipshitzian function $f-\frac{\Gamma +\gamma }{2}\cdot e.$
\end{proof}

\subsection{Applications for Some Particular Functions}

\textbf{1.} We have the following important inequality in Operator Theory
that is well known as the H\"{o}lder-McCarthy inequality:

\begin{theorem}[H\"{o}lder-McCarthy, 1967, \protect\cite{III.Mc}]
\label{III.t.3.1} Let $A$ be a selfadjoint positive operator on a Hilbert
space $H$. Then

(i) \ \ $\left\langle A^{r}x,x\right\rangle \geq \left\langle
Ax,x\right\rangle ^{r}$ for all $r>1$ and $x\in H$ with $\left\| x\right\|
=1;$

(ii) \ $\left\langle A^{r}x,x\right\rangle \leq \left\langle
Ax,x\right\rangle ^{r}$ for all $0<r<1$ and $x\in H$ with $\left\| x\right\|
=1;$

(iii) If $A$ is invertible, then $\left\langle A^{-r}x,x\right\rangle \geq
\left\langle Ax,x\right\rangle ^{-r}$ for all $r>0$ and $x\in H$ with $%
\left\| x\right\| =1.$
\end{theorem}

We can provide the following reverse inequalities:

\begin{proposition}
\label{III.p.3.1}Let $A$ be a selfadjoint positive operator on a Hilbert
space $H$ and $0<r<1.$ Then 
\begin{equation}
\left( 0\leq \right) \left\langle Ax,x\right\rangle ^{r}-\left\langle
A^{r}x,x\right\rangle \leq \left\langle \left\vert \left\langle
Ax,x\right\rangle \cdot 1_{H}-A\right\vert x,x\right\rangle ^{r}\leq
D^{r}\left( A;x\right)  \label{III.e.3.1}
\end{equation}%
for all $x\in H$ with $\left\Vert x\right\Vert =1.$
\end{proposition}

\begin{proof}
Follows from Corollary \ref{III.c.2.1} by taking into account that the
function $f\left( t\right) =t^{r}$ is of $r-L-$H\"{o}lder type with $L=1$ on
any compact interval of $\left( 0,\infty \right) .$
\end{proof}

On making use of Corollary \ref{III.c.4.1} we can state the following result
as well:

\begin{proposition}
\label{III.p.3.2}Let $A$ be a selfadjoint positive operator on a Hilbert
space $H$. Assume that $Sp\left( A\right) \subseteq \left[ m,M\right]
\subseteq \lbrack 0,\infty ).$

(i) We have%
\begin{align}
0& \leq \left\langle A^{r}x,x\right\rangle -\left\langle Ax,x\right\rangle
^{r}  \label{III.e.3.2} \\
& \leq \frac{1}{2}r\left( M^{r-1}-m^{r-1}\right) \left\langle \left\vert
\left\langle Ax,x\right\rangle \cdot 1_{H}-A\right\vert x,x\right\rangle 
\notag \\
& \leq \frac{1}{2}r\left( M^{r-1}-m^{r-1}\right) D\left( A;x\right) \leq 
\frac{1}{4}r\left( M^{r-1}-m^{r-1}\right) \left( M-m\right)  \notag
\end{align}%
for all $r>1$ and $x\in H$ with $\left\Vert x\right\Vert =1;$

(ii) We also have%
\begin{align}
0& \leq \left\langle Ax,x\right\rangle ^{r}-\left\langle
A^{r}x,x\right\rangle  \label{III.e.3.3} \\
& \leq \frac{1}{2}r\left( \frac{M^{1-r}-m^{1-r}}{m^{1-r}M^{1-r}}\right)
\left\langle \left\vert \left\langle Ax,x\right\rangle \cdot
1_{H}-A\right\vert x,x\right\rangle  \notag \\
& \leq \frac{1}{2}r\left( \frac{M^{1-r}-m^{1-r}}{m^{1-r}M^{1-r}}\right)
D\left( A;x\right) \leq \frac{1}{4}r\left( \frac{M^{1-r}-m^{1-r}}{%
m^{1-r}M^{1-r}}\right) \left( M-m\right)  \notag
\end{align}%
for all $0<r<1$ and $x\in H$ with $\left\Vert x\right\Vert =1;$

(iii) If $A$ is invertible, then 
\begin{align}
0& \leq \left\langle A^{-r}x,x\right\rangle -\left\langle Ax,x\right\rangle
^{-r}  \label{III.e.3.4} \\
& \leq \frac{1}{2}r\left( \frac{M^{r+1}-m^{r+1}}{M^{r+1}m^{r+1}}\right)
\left\langle \left\vert \left\langle Ax,x\right\rangle \cdot
1_{H}-A\right\vert x,x\right\rangle  \notag \\
& \leq \frac{1}{2}r\left( \frac{M^{r+1}-m^{r+1}}{M^{r+1}m^{r+1}}\right)
D\left( A;x\right) \leq \frac{1}{4}r\left( \frac{M^{r+1}-m^{r+1}}{%
M^{r+1}m^{r+1}}\right) \left( M-m\right)  \notag
\end{align}%
for all $r>0$ and $x\in H$ with $\left\Vert x\right\Vert =1.$
\end{proposition}

\textbf{2.} Consider the convex function $f:\left( 0,\infty \right)
\rightarrow \mathbb{R}$, $f\left( x\right) =-\ln x.$ On utilizing the
inequality (\ref{III.e.4.5}), we can state the following result:

\begin{proposition}
\label{III.p.4.3}For any positive definite operator $A$ on the Hilbert space 
$H$ with $Sp\left( A\right) \subseteq \left[ m,M\right] \subseteq \lbrack
0,\infty )$ we have the inequality%
\begin{align}
& \left( 0\leq \right) \ln \left( \left\langle Ax,x\right\rangle \right)
-\left\langle \ln \left( A\right) x,x\right\rangle  \label{III.e.3.5} \\
& \leq \frac{1}{2}\cdot \frac{M-m}{mM}\left\langle \left\vert \left\langle
Ax,x\right\rangle \cdot 1_{H}-A\right\vert x,x\right\rangle  \notag \\
& \leq \frac{1}{2}\cdot \frac{M-m}{mM}D\left( A;x\right) \leq \frac{1}{4}%
\cdot \frac{\left( M-m\right) ^{2}}{mM}  \notag
\end{align}%
for any $x\in H$ with $\left\Vert x\right\Vert =1.$
\end{proposition}

Finally, the following result for logarithms also holds:

\begin{proposition}
\label{III.p.4.4}Under the assumptions of Proposition \ref{III.p.4.3} we
have the inequality%
\begin{align}
& \left( 0\leq \right) \left\langle A\ln \left( A\right) x,x\right\rangle
-\left\langle Ax,x\right\rangle \ln \left( \left\langle Ax,x\right\rangle
\right)  \label{III.e.3.6} \\
& \leq \ln \sqrt{\frac{M}{m}}\left\langle \left\vert \left\langle
Ax,x\right\rangle \cdot 1_{H}-A\right\vert x,x\right\rangle  \notag \\
& \leq \ln \sqrt{\frac{M}{m}}\cdot D\left( A;x\right) \leq \frac{1}{2}\left(
M-m\right) \ln \sqrt{\frac{M}{m}}  \notag
\end{align}%
for any $x\in H$ with $\left\Vert x\right\Vert =1.$
\end{proposition}

\begin{remark}
\label{III.r.3.1}On utilizing the results from the previous sections for
other convex functions of interest such as $f\left( x\right) =\ln \left[
\left( 1-x\right) /x\right] ,$ $x\in \left( 0,1/2\right) $ or $f\left(
x\right) =\ln \left( 1+\exp x\right) ,x\in \left( -\infty ,\infty \right) $
we can get other interesting operator inequalities. However, the details are
left to the interested reader.
\end{remark}

\section{Other Ostrowski Inequalities for Continuous Functions}

\subsection{Inequalities for Absolutely Continuous Functions of Selfadjoint
Operators}

We start with the following scalar inequality that is of interest in itself
since it provides a generalization of the Ostrowski inequality when upper
and lower bounds for the derivative are provided:

\begin{lemma}[Dragomir, 2010, \protect\cite{III.a.SSD3}]
\label{III.a.l.2.1}Let $f:\left[ a,b\right] \rightarrow \mathbb{R}$ be an
absolutely continuous function whose derivative is bounded above and below
on $\left[ a,b\right] ,$ i.e., there exists the real constants $\gamma $ and 
$\Gamma ,\gamma <\Gamma $ with the property that $\gamma \leq f^{\prime
}\left( s\right) \leq \Gamma $ for almost every $s\in \left[ a,b\right] .$
Then we have the double inequality 
\begin{align}
& -\frac{1}{2}\cdot \frac{\Gamma -\gamma }{b-a}\left[ \left( s-\frac{b\Gamma
-a\gamma }{\Gamma -\gamma }\right) ^{2}-\Gamma \gamma \left( \frac{b-a}{%
\Gamma -\gamma }\right) ^{2}\right]  \label{III.a.e.2.1} \\
& \leq f\left( s\right) -\frac{1}{b-a}\int_{a}^{b}f\left( t\right) dt  \notag
\\
& \leq \frac{1}{2}\cdot \frac{\Gamma -\gamma }{b-a}\left[ \left( s-\frac{%
a\Gamma -b\gamma }{\Gamma -\gamma }\right) ^{2}-\Gamma \gamma \left( \frac{%
b-a}{\Gamma -\gamma }\right) ^{2}\right]  \notag
\end{align}%
for any $s\in \left[ a,b\right] .$ The inequalities are sharp.
\end{lemma}

\begin{proof}
We start with \textit{Montgomery's identity}%
\begin{align}
& f\left( s\right) -\frac{1}{b-a}\int_{a}^{b}f\left( t\right) dt
\label{III.a.MO} \\
& =\frac{1}{b-a}\int_{a}^{s}\left( t-a\right) f^{\prime }\left( t\right) dt+%
\frac{1}{b-a}\int_{s}^{b}\left( t-b\right) f^{\prime }\left( t\right) dt 
\notag
\end{align}%
that holds for any $s\in \left[ a,b\right] .$

Since $\gamma \leq f^{\prime }\left( t\right) \leq \Gamma $ for almost every 
$t\in \left[ a,b\right] ,$ then%
\begin{equation*}
\frac{\gamma }{b-a}\int_{a}^{s}\left( t-a\right) dt\leq \frac{1}{b-a}%
\int_{a}^{s}\left( t-a\right) f^{\prime }\left( t\right) dt\leq \frac{\Gamma 
}{b-a}\int_{a}^{s}\left( t-a\right) dt
\end{equation*}%
and%
\begin{equation*}
\frac{\Gamma }{b-a}\int_{s}^{b}\left( b-t\right) dt\leq \frac{1}{b-a}%
\int_{s}^{b}\left( b-t\right) f^{\prime }\left( t\right) dt\leq \frac{\Gamma 
}{b-a}\int_{s}^{b}\left( b-t\right) dt
\end{equation*}%
for any $s\in \left[ a,b\right] .$

Now, due to the fact that%
\begin{equation*}
\int_{a}^{s}\left( t-a\right) dt=\frac{1}{2}\left( s-a\right) ^{2}\text{ and 
}\int_{s}^{b}\left( b-t\right) dt=\frac{1}{2}\left( b-s\right) ^{2}
\end{equation*}%
then by (\ref{III.a.MO}) we deduce the following inequality that is of
interest in itself:%
\begin{align}
& -\frac{1}{2\left( b-a\right) }\left[ \Gamma \left( b-s\right) ^{2}-\gamma
\left( s-a\right) ^{2}\right]  \label{III.a.e.2.2} \\
& \leq f\left( s\right) -\frac{1}{b-a}\int_{a}^{b}f\left( t\right) dt  \notag
\\
& \leq \frac{1}{2\left( b-a\right) }\left[ \Gamma \left( s-a\right)
^{2}-\gamma \left( b-s\right) ^{2}\right]  \notag
\end{align}%
for any $s\in \left[ a,b\right] .$

Further on, if we denote by%
\begin{equation*}
A:=\gamma \left( s-a\right) ^{2}-\Gamma \left( b-s\right) ^{2}\text{ and }%
B:=\Gamma \left( s-a\right) ^{2}-\gamma \left( b-s\right) ^{2}
\end{equation*}%
then, after some elementary calculations, we derive that%
\begin{equation*}
A=-\left( \Gamma -\gamma \right) \left( s-\frac{b\Gamma -a\gamma }{\Gamma
-\gamma }\right) ^{2}+\frac{\Gamma \gamma }{\Gamma -\gamma }\left(
b-a\right) ^{2}
\end{equation*}%
and%
\begin{equation*}
B=\left( \Gamma -\gamma \right) \left( s-\frac{a\Gamma -b\gamma }{\Gamma
-\gamma }\right) ^{2}-\frac{\Gamma \gamma }{\Gamma -\gamma }\left(
b-a\right) ^{2}
\end{equation*}%
which, together with (\ref{III.a.e.2.2}), produces the desired result (\ref%
{III.a.e.2.1}).

The sharpness of the inequalities follow from the sharpness of some
particular cases outlined below. The details are omitted.
\end{proof}

\begin{corollary}
\label{III.a.c.2.1}With the assumptions of Lemma \ref{III.a.l.2.1} we have
the inequalities%
\begin{equation}
\frac{1}{2}\gamma \left( b-a\right) \leq \frac{1}{b-a}\int_{a}^{b}f\left(
t\right) dt-f\left( a\right) \leq \frac{1}{2}\Gamma \left( b-a\right)
\label{III.a.e.2.3}
\end{equation}%
and%
\begin{equation}
\frac{1}{2}\gamma \left( b-a\right) \leq f\left( b\right) -\frac{1}{b-a}%
\int_{a}^{b}f\left( t\right) dt\leq \frac{1}{2}\Gamma \left( b-a\right)
\label{III.a.e.2.4}
\end{equation}%
and%
\begin{equation}
\left\vert f\left( \frac{a+b}{2}\right) -\frac{1}{b-a}\int_{a}^{b}f\left(
t\right) dt\right\vert \leq \frac{1}{8}\left( \Gamma -\gamma \right) \left(
b-a\right)  \label{III.a.e.2.5}
\end{equation}%
respectively. The constant $\frac{1}{8}$ is best possible in (\ref%
{III.a.e.2.5}).
\end{corollary}

The proof is obvious from (\ref{III.a.e.2.1}) on choosing $s=a,s=b$ and $s=%
\frac{a+b}{2},$ respectively.

\begin{corollary}[Dragomir, 2010, \protect\cite{III.a.SSD3}]
\label{III.a.c.2.2}With the assumptions of Lemma \ref{III.a.l.2.1} and if,
in addition $\gamma =-\alpha $ and $\Gamma =\beta $ with $\alpha ,\beta >0$
then%
\begin{equation}
\frac{1}{b-a}\int_{a}^{b}f\left( t\right) dt-f\left( \frac{b\beta +a\alpha }{%
\beta +\alpha }\right) \leq \frac{1}{2}\cdot \alpha \beta \left( \frac{b-a}{%
\beta +\alpha }\right)  \label{III.a.e.2.6}
\end{equation}%
and%
\begin{equation}
f\left( \frac{a\beta +b\alpha }{\beta +\alpha }\right) -\frac{1}{b-a}%
\int_{a}^{b}f\left( t\right) dt\leq \frac{1}{2}\cdot \alpha \beta \left( 
\frac{b-a}{\beta +\alpha }\right) .  \label{III.a.e.2.7}
\end{equation}
\end{corollary}

The proof follows from (\ref{III.a.e.2.1}) on choosing $s=\frac{b\beta
+a\alpha }{\beta +\alpha }\in \left[ a,b\right] $ and $s=\frac{a\beta
+b\alpha }{\beta +\alpha }\in \left[ a,b\right] ,$ respectively.

\begin{remark}
\label{III.a.r.2.1}If $f:\left[ a,b\right] \rightarrow \mathbb{R}$ is
absolutely continuous and 
\begin{equation*}
\left\Vert f^{\prime }\right\Vert _{\infty }:=ess\sup_{t\in \left[ a,b\right]
}\left\vert f^{\prime }\left( t\right) \right\vert <\infty ,
\end{equation*}
then by choosing $\gamma =-\left\Vert f^{\prime }\right\Vert _{\infty }$ and 
$\Gamma =\left\Vert f^{\prime }\right\Vert _{\infty }$ in (\ref{III.a.e.2.1}%
) we deduce the classical Ostrowski's inequality for absolutely continuous
functions. The constant $\frac{1}{4}$ in Ostrowski's inequality is best
possible.
\end{remark}

We are able now to state the following result providing upper and lower
bounds for absolutely convex functions of selfadjoint operators in Hilbert
spaces whose derivatives are bounded below and above:

\begin{theorem}[Dragomir, 2010, \protect\cite{III.a.SSD3}]
\label{III.a.t.2.1}Let $A$ be a selfadjoint operator in the Hilbert space $H$
with the spectrum $Sp\left( A\right) \subseteq \left[ m,M\right] $ for some
real numbers $m<M.$ If $f:\left[ m,M\right] \rightarrow \mathbb{R}$ is an
absolutely continuous function such that there exists the real constants $%
\gamma $ and $\Gamma ,\gamma <\Gamma $ with the property that $\gamma \leq
f^{\prime }\left( s\right) \leq \Gamma $ for almost every $s\in \left[ m,M%
\right] ,$ then we have the following double inequality in the operator
order of $B\left( H\right) :$ 
\begin{align}
& -\frac{1}{2}\cdot \frac{\Gamma -\gamma }{M-m}\left[ \left( A-\frac{M\Gamma
-m\gamma }{\Gamma -\gamma }\cdot 1_{H}\right) ^{2}-\Gamma \gamma \left( 
\frac{M-m}{\Gamma -\gamma }\right) ^{2}\cdot 1_{H}\right]
\label{III.a.e.2.9} \\
& \leq f\left( A\right) -\left( \frac{1}{M-m}\int_{m}^{M}f\left( t\right)
dt\right) \cdot 1_{H}  \notag \\
& \leq \frac{1}{2}\cdot \frac{\Gamma -\gamma }{M-m}\left[ \left( A-\frac{%
m\Gamma -M\gamma }{\Gamma -\gamma }\cdot 1_{H}\right) ^{2}-\Gamma \gamma
\left( \frac{M-m}{\Gamma -\gamma }\right) ^{2}\cdot 1_{H}\right] .  \notag
\end{align}
\end{theorem}

The proof follows by the property (\ref{P}) applied for the inequality (\ref%
{III.a.e.2.1}) in Lemma \ref{III.a.l.2.1}.

\begin{theorem}[Dragomir, 2010, \protect\cite{III.a.SSD3}]
\label{III.a.t.2.1.a}With the assumptions in Theorem \ref{III.a.t.2.1} we
have in the operator order the following inequalities%
\begin{align}
& \left\vert f\left( A\right) -\left( \frac{1}{M-m}\int_{m}^{M}f\left(
t\right) dt\right) \cdot 1_{H}\right\vert  \label{III.a.e.2.10} \\
& \leq \left\{ 
\begin{array}{l}
\left[ \frac{1}{4}1_{H}+\left( \frac{A-\frac{m+M}{2}1_{H}}{M-m}\right) ^{2}%
\right] \left( M-m\right) \left\Vert f^{\prime }\right\Vert _{\infty }\quad 
\text{if \ }f^{\prime }\in L_{\infty }\left[ m,M\right] ; \\ 
\\ 
\frac{1}{\left( p+1\right) ^{\frac{1}{p}}}\left[ \left( \frac{A-m1_{H}}{M-m}%
\right) ^{p+1}+\left( \frac{M1_{H}-A}{M-m}\right) ^{p+1}\right] \left(
M-m\right) ^{\frac{1}{q}}\left\Vert f^{\prime }\right\Vert _{q} \\ 
\hfill \text{if \ }f^{\prime }\in L_{p}\left[ m,M\right] ,\ \frac{1}{p}+%
\frac{1}{q}=1,\ p>1; \\ 
\\ 
\left[ \frac{1}{2}1_{H}+\left\vert \frac{A-\frac{m+M}{2}1_{H}}{M-m}%
\right\vert \right] \left\Vert f^{\prime }\right\Vert _{1}.%
\end{array}%
\right.  \notag
\end{align}
\end{theorem}

The proof is obvious by the scalar inequalities from Theorem \ref{III.b.tb}
and the property (\ref{P}).

The third inequality in (\ref{III.a.e.2.10}) can be naturally generalized
for functions of bounded variation as follows:

\begin{theorem}[Dragomir, 2010, \protect\cite{III.a.SSD3}]
\label{III.a.t.2.1.b}Let $A$ be a selfadjoint operator in the Hilbert space $%
H$ with the spectrum $Sp\left( A\right) \subseteq \left[ m,M\right] $ for
some real numbers $m<M.$ If $f:\left[ m,M\right] \rightarrow \mathbb{R}$ is
a continuous function of bounded variation on $\left[ m,M\right] ,$ then we
have the inequality%
\begin{align}
& \left\vert f\left( A\right) -\left( \frac{1}{M-m}\int_{m}^{M}f\left(
t\right) dt\right) \cdot 1_{H}\right\vert  \label{III.a.e.2.10.1} \\
& \leq \left[ \frac{1}{2}1_{H}+\left\vert \frac{A-\frac{m+M}{2}1_{H}}{M-m}%
\right\vert \right] \dbigvee\limits_{m}^{M}\left( f\right)  \notag
\end{align}%
where $\dbigvee\limits_{m}^{M}\left( f\right) $ denotes the total variation
of $f$ on $\left[ m,M\right] .$ The constant $\frac{1}{2}$ is best possible
in (\ref{III.a.e.2.10.1}).
\end{theorem}

\begin{proof}
Follows from the scalar inequality obtained by the author in \cite{III.b.DR2}%
, namely%
\begin{equation}
\left\vert f\left( s\right) -\frac{1}{b-a}\int_{a}^{b}f\left( t\right)
dt\right\vert \leq \left[ \frac{1}{2}+\left\vert \frac{s-\frac{a+b}{2}}{b-a}%
\right\vert \right] \dbigvee\limits_{a}^{b}\left( f\right)
\label{III.a.e.2.10.2}
\end{equation}%
for any $s\in \left[ a,b\right] ,$ where $f$ is a function of bounded
variation on $\left[ a,b\right] .$ The constant $\frac{1}{2}$ is best
possible in (\ref{III.a.e.2.10.2}).
\end{proof}

\subsection{Inequalities for Convex Functions of Selfadjoint Operators}

The case of convex functions is important for applications.

We need the following lemma.

\begin{lemma}[Dragomir, 2010, \protect\cite{III.a.SSD3}]
\label{III.a.l.2.2}Let $f:\left[ a,b\right] \rightarrow \mathbb{R}$ be a
differentiable convex function such that the derivative $f^{\prime }$ is
continuous on $\left( a,b\right) $ and with the lateral derivative finite
and $f_{-}^{\prime }\left( b\right) \neq f_{+}^{\prime }\left( a\right) $.
Then we have the following double inequality%
\begin{align}
& -\frac{1}{2}\cdot \frac{f_{-}^{\prime }\left( b\right) -f_{+}^{\prime
}\left( a\right) }{b-a}  \label{III.a.e.2.11} \\
& \times \left[ \left( s-\frac{bf_{-}^{\prime }\left( b\right)
-af_{+}^{\prime }\left( a\right) }{f_{-}^{\prime }\left( b\right)
-f_{+}^{\prime }\left( a\right) }\right) ^{2}-f_{-}^{\prime }\left( b\right)
f_{+}^{\prime }\left( a\right) \left( \frac{b-a}{f_{-}^{\prime }\left(
b\right) -f_{+}^{\prime }\left( a\right) }\right) ^{2}\right]  \notag \\
& \leq f\left( s\right) -\frac{1}{b-a}\int_{a}^{b}f\left( t\right) dt\leq
f^{\prime }\left( s\right) \left( s-\frac{a+b}{2}\right)  \notag
\end{align}%
for any $s\in \left[ a,b\right] .$
\end{lemma}

\begin{proof}
Since $f$ is convex, then by the fact that $f^{\prime }$ is monotonic
nondecreasing, we have 
\begin{equation*}
\frac{f_{+}^{\prime }\left( a\right) }{b-a}\int_{a}^{s}\left( t-a\right)
dt\leq \frac{1}{b-a}\int_{a}^{s}\left( t-a\right) f^{\prime }\left( t\right)
dt\leq \frac{f^{\prime }\left( s\right) }{b-a}\int_{a}^{s}\left( t-a\right)
dt
\end{equation*}%
and%
\begin{equation*}
\frac{f^{\prime }\left( s\right) }{b-a}\int_{s}^{b}\left( b-t\right) dt\leq 
\frac{1}{b-a}\int_{s}^{b}\left( b-t\right) f^{\prime }\left( t\right) dt\leq 
\frac{f_{-}^{\prime }\left( b\right) }{b-a}\int_{s}^{b}\left( b-t\right) dt
\end{equation*}%
for any $s\in \left[ a,b\right] ,$ where $f_{+}^{\prime }\left( a\right) $
and $f_{-}^{\prime }\left( b\right) $ are the lateral derivatives in $a$ and 
$b$ respectively.

Utilising the Montgomery identity (\ref{III.a.MO}) we then have%
\begin{align*}
& \frac{f_{+}^{\prime }\left( a\right) }{b-a}\int_{a}^{s}\left( t-a\right)
dt-\frac{f_{-}^{\prime }\left( b\right) }{b-a}\int_{s}^{b}\left( b-t\right)
dt \\
& \leq f\left( s\right) -\frac{1}{b-a}\int_{a}^{b}f\left( t\right) dt \\
& \leq \frac{f^{\prime }\left( s\right) }{b-a}\int_{a}^{s}\left( t-a\right)
dt-\frac{f^{\prime }\left( s\right) }{b-a}\int_{s}^{b}\left( b-t\right) dt
\end{align*}%
which is equivalent with the following inequality that is of interest in
itself%
\begin{align}
& \frac{1}{2\left( b-a\right) }\left[ f_{+}^{\prime }\left( a\right) \left(
s-a\right) ^{2}-f_{-}^{\prime }\left( b\right) \left( b-s\right) ^{2}\right]
\label{III.a.e.2.12} \\
& \leq f\left( s\right) -\frac{1}{b-a}\int_{a}^{b}f\left( t\right) dt\leq
f^{\prime }\left( s\right) \left( s-\frac{a+b}{2}\right)  \notag
\end{align}%
for any $s\in \left[ a,b\right] .$

A simple calculation reveals now that the left side of (\ref{III.a.e.2.12})
coincides with the same side of the desired inequality (\ref{III.a.e.2.11}).
\end{proof}

We are able now to sate our result for convex functions of selfadjoint
operators:

\begin{theorem}[Dragomir, 2010, \protect\cite{III.a.SSD3}]
\label{III.a.t.2.2}Let $A$ be a selfadjoint operator in the Hilbert space $H$
with the spectrum $Sp\left( A\right) \subseteq \left[ m,M\right] $ for some
real numbers $m<M.$ If $f:\left[ m,M\right] \rightarrow \mathbb{R}$ is a
differentiable convex function such that the derivative $f^{\prime }$ is
continuous on $\left( m,M\right) $ and with the lateral derivative finite
and $f_{-}^{\prime }\left( M\right) \neq f_{+}^{\prime }\left( m\right) ,$
then we have the double inequality in the operator order of $B\left(
H\right) $%
\begin{multline}
-\frac{1}{2}\cdot \frac{f_{-}^{\prime }\left( M\right) -f_{+}^{\prime
}\left( m\right) }{M-m}  \label{III.a.e.2.13} \\
\times \left[ \left( A-\frac{Mf_{-}^{\prime }\left( M\right) -mf_{+}^{\prime
}\left( m\right) }{f_{-}^{\prime }\left( M\right) -f_{+}^{\prime }\left(
m\right) }\cdot 1_{H}\right) ^{2}-f_{-}^{\prime }\left( M\right)
f_{+}^{\prime }\left( m\right) \left( \frac{M-m}{f_{-}^{\prime }\left(
M\right) -f_{+}^{\prime }\left( m\right) }\right) ^{2}\cdot 1_{H}\right] \\
\leq f\left( A\right) -\left( \frac{1}{M-m}\int_{m}^{M}f\left( t\right)
dt\right) \cdot 1_{H}\leq \left( A-\frac{m+M}{2}\cdot 1_{H}\right) f^{\prime
}\left( A\right) .
\end{multline}
\end{theorem}

The proof follows from the scalar case in Lemma \ref{III.a.l.2.2}.

\begin{remark}
\label{III.a.r.2.2}We observe that one can drop the assumption of
differentiability of the convex function and will still have the first
inequality in (\ref{III.a.e.2.13}). This follows from the fact that the
class of differentiable convex functions is dense in the class of all convex
functions defined on a given interval.
\end{remark}

A different lower bound for the quantity%
\begin{equation*}
f\left( A\right) -\left( \frac{1}{M-m}\int_{m}^{M}f\left( t\right) dt\right)
\cdot 1_{H}
\end{equation*}%
expressed only in terms of the operator $A$ and not its second power as
above, also holds:

\begin{theorem}[Dragomir, 2010, \protect\cite{III.a.SSD3}]
\label{III.a.t.2.3}Let $A$ be a selfadjoint operator in the Hilbert space $H$
with the spectrum $Sp\left( A\right) \subseteq \left[ m,M\right] $ for some
real numbers $m<M.$ If $f:\left[ m,M\right] \rightarrow \mathbb{R}$ is a
convex function on $\left[ m,M\right] ,$ then we have the following
inequality in the operator order of $B\left( H\right) $%
\begin{align}
& f\left( A\right) -\left( \frac{1}{M-m}\int_{m}^{M}f\left( t\right)
dt\right) \cdot 1_{H}  \label{III.a.e.2.14} \\
& \geq \left( \frac{1}{M-m}\int_{m}^{M}f\left( t\right) dt\right) \cdot 1_{H}
\notag \\
& -\frac{f\left( M\right) \left( M\cdot 1_{H}-A\right) +f\left( m\right)
\left( A-m\cdot 1_{H}\right) }{M-m}.  \notag
\end{align}
\end{theorem}

\begin{proof}
It suffices to prove for the case of differentiable convex functions defined
on $\left( m,M\right) .$

So, by the gradient inequality we have that%
\begin{equation*}
f\left( t\right) -f\left( s\right) \geq \left( t-s\right) f^{\prime }\left(
s\right)
\end{equation*}%
for any $t,s\in \left( m,M\right) .$

Now, if we integrate this inequality over $s\in \left[ m,M\right] $ we get%
\begin{align}
& \left( M-m\right) f\left( t\right) -\int_{m}^{M}f\left( s\right) ds
\label{III.a.e.2.15} \\
& \geq \int_{m}^{M}\left( t-s\right) f^{\prime }\left( s\right) ds  \notag \\
& =\int_{m}^{M}f\left( s\right) ds-\left( M-t\right) f\left( M\right)
-\left( t-m\right) f\left( m\right)  \notag
\end{align}%
for each $s\in \left[ m,M\right] .$

Finally, if we apply to the inequality (\ref{III.a.e.2.15}) the property (%
\ref{P}), we deduce the desired result (\ref{III.a.e.2.14}).
\end{proof}

\begin{corollary}[Dragomir, 2010, \protect\cite{III.a.SSD3}]
\label{III.a.c.2.3}With the assumptions of Theorem \ref{III.a.t.2.3} we have
the following double inequality in the operator order%
\begin{align}
& \frac{f\left( m\right) +f\left( M\right) }{2}\cdot 1_{H}
\label{III.a.e.2.16} \\
& \geq \frac{1}{2}\left[ f\left( A\right) +\frac{f\left( M\right) \left(
M\cdot 1_{H}-A\right) +f\left( m\right) \left( A-m\cdot 1_{H}\right) }{M-m}%
\right]  \notag \\
& \geq \left( \frac{1}{M-m}\int_{m}^{M}f\left( t\right) dt\right) \cdot
1_{H}.  \notag
\end{align}
\end{corollary}

\begin{proof}
The second inequality is equivalent with (\ref{III.a.e.2.14}).

For the first inequality, we observe, by the convexity of $f$ we have that%
\begin{equation*}
\frac{f\left( M\right) \left( t-m\right) +f\left( m\right) \left( M-t\right) 
}{M-m}\geq f\left( t\right)
\end{equation*}%
for any $t\in \left[ m,M\right] ,$ which produces the operator inequality%
\begin{equation}
\frac{f\left( M\right) \left( A-m\cdot 1_{H}\right) +f\left( m\right) \left(
M\cdot 1_{H}-A\right) }{M-m}\geq f\left( A\right) .  \label{III.a.e.2.17}
\end{equation}%
Now, if in both sides of (\ref{III.a.e.2.17}) we add the same quantity 
\begin{equation*}
\frac{f\left( M\right) \left( M\cdot 1_{H}-A\right) +f\left( m\right) \left(
A-m\cdot 1_{H}\right) }{M-m}
\end{equation*}%
and perform the calculations, then we obtain the first part of (\ref%
{III.a.e.2.16}) and the proof is complete.
\end{proof}

\subsection{Some Vector Inequalities}

The following result holds:

\begin{theorem}[Dragomir, 2010, \protect\cite{III.a.SSD3}]
\label{III.a.t.3.1}Let $A$ be a selfadjoint operator in the Hilbert space $H$
with the spectrum $Sp\left( A\right) \subseteq \left[ m,M\right] $ for some
real numbers $m<M$ and let $\left\{ E_{\lambda }\right\} _{\lambda }$ be its 
\textit{spectral family.} If $f:\left[ m,M\right] \rightarrow \mathbb{R}$ is
an absolutely continuous function on $\left[ m,M\right] $, then we have the
inequalities%
\begin{align}
& \left\vert f\left( s\right) \left\langle x,y\right\rangle -\left\langle
f\left( A\right) x,y\right\rangle \right\vert   \label{III.a.e.3.2} \\
& \leq \dbigvee\limits_{m-0}^{M}\left( \left\langle E_{\left( \cdot \right)
}x,y\right\rangle \right)   \notag \\
& \times \left\{ 
\begin{array}{cc}
\left[ \frac{1}{2}\left( M-m\right) +\left\vert s-\frac{m+M}{2}\right\vert %
\right] \left\Vert f^{\prime }\right\Vert _{\infty } & \text{if }f^{\prime
}\in L_{\infty }\left[ m,M\right]  \\ 
&  \\ 
\left[ \frac{1}{2}\left( M-m\right) +\left\vert s-\frac{m+M}{2}\right\vert %
\right] ^{1/q}\left\Vert f^{\prime }\right\Vert _{p} & 
\begin{array}{c}
\text{if }f^{\prime }\in L_{p}\left[ m,M\right] ,p>1, \\ 
\frac{1}{p}+\frac{1}{q}=1,%
\end{array}%
\end{array}%
\right.   \notag \\
& \leq \left\Vert x\right\Vert \left\Vert y\right\Vert   \notag \\
& \times \left\{ 
\begin{array}{cc}
\left[ \frac{1}{2}\left( M-m\right) +\left\vert s-\frac{m+M}{2}\right\vert %
\right] \left\Vert f^{\prime }\right\Vert _{\infty } & \text{if }f^{\prime
}\in L_{\infty }\left[ m,M\right]  \\ 
&  \\ 
\left[ \frac{1}{2}\left( M-m\right) +\left\vert s-\frac{m+M}{2}\right\vert %
\right] ^{1/q}\left\Vert f^{\prime }\right\Vert _{p} & 
\begin{array}{c}
\text{if }f^{\prime }\in L_{p}\left[ m,M\right] ,p>1, \\ 
\frac{1}{p}+\frac{1}{q}=1,%
\end{array}%
\end{array}%
\right.   \notag
\end{align}%
for any $x,y\in H$ and $s\in \left[ m,M\right] .$
\end{theorem}

\begin{proof}
Since $f$ is absolutely continuous, then we have%
\begin{align}
& \left\vert f\left( s\right) -f\left( t\right) \right\vert
\label{III.a.e.3.3} \\
& =\left\vert \int_{s}^{t}f^{\prime }\left( u\right) du\right\vert \leq
\left\vert \int_{s}^{t}\left\vert f^{\prime }\left( u\right) \right\vert
du\right\vert  \notag \\
& \leq \left\{ 
\begin{array}{cc}
\left\vert t-s\right\vert \left\Vert f^{\prime }\right\Vert _{\infty } & 
\text{if }f^{\prime }\in L_{\infty }\left[ m,M\right] \\ 
&  \\ 
\left\vert t-s\right\vert ^{1/q}\left\Vert f^{\prime }\right\Vert _{p} & 
\text{if }f^{\prime }\in L_{p}\left[ m,M\right] ,p>1,\frac{1}{p}+\frac{1}{q}%
=1,%
\end{array}%
\right.  \notag
\end{align}%
for any $s,t\in \left[ m,M\right] .$

It is well known that if $p:\left[ a,b\right] \rightarrow \mathbb{C}$ is a
continuous functions and $v:\left[ a,b\right] \rightarrow \mathbb{C}$ is of
bounded variation, then the Riemann-Stieltjes integral $\int_{a}^{b}p\left(
t\right) dv\left( t\right) $ exists and the following inequality holds%
\begin{equation*}
\left\vert \int_{a}^{b}p\left( t\right) dv\left( t\right) \right\vert \leq
\max_{t\in \left[ a,b\right] }\left\vert p\left( t\right) \right\vert
\dbigvee\limits_{a}^{b}\left( v\right) ,
\end{equation*}%
where $\dbigvee\limits_{a}^{b}\left( v\right) $ denotes the total variation
of $v$ on $\left[ a,b\right] .$

Now, by the above property of the Riemann-Stieltjes integral we have from
the representation (\ref{III.a.e.3.10}) that%
\begin{align}
& \left\vert f\left( s\right) \left\langle x,y\right\rangle -\left\langle
f\left( A\right) x,y\right\rangle \right\vert   \label{III.a.e.3.4} \\
& =\left\vert \int_{m-0}^{M}\left[ f\left( s\right) -f\left( t\right) \right]
d\left( \left\langle E_{t}x,y\right\rangle \right) \right\vert   \notag \\
& \leq \max_{t\in \left[ m,M\right] }\left\vert f\left( s\right) -f\left(
t\right) \right\vert \dbigvee\limits_{m-0}^{M}\left( \left\langle E_{\left(
\cdot \right) }x,y\right\rangle \right)   \notag \\
& \leq \dbigvee\limits_{m-0}^{M}\left( \left\langle E_{\left( \cdot \right)
}x,y\right\rangle \right)   \notag \\
& \times \left\{ 
\begin{array}{cc}
\max_{t\in \left[ m,M\right] }\left\vert t-s\right\vert \left\Vert f^{\prime
}\right\Vert _{\infty } & \text{if }f^{\prime }\in L_{\infty }\left[ m,M%
\right]  \\ 
&  \\ 
\max_{t\in \left[ m,M\right] }\left\vert t-s\right\vert ^{1/q}\left\Vert
f^{\prime }\right\Vert _{p} & 
\begin{array}{c}
\text{if }f^{\prime }\in L_{p}\left[ m,M\right] ,p>1, \\ 
\frac{1}{p}+\frac{1}{q}=1,%
\end{array}%
\end{array}%
\right.   \notag \\
& :=F  \notag
\end{align}%
where $\dbigvee\limits_{m-0}^{M}\left( \left\langle E_{\left( \cdot \right)
}x,y\right\rangle \right) $ denotes the total variation of $\left\langle
E_{\left( \cdot \right) }x,y\right\rangle $ and $x,y\in H.$

Since, obviously, we have $\max_{t\in \left[ m,M\right] }\left\vert
t-s\right\vert =\frac{1}{2}\left( M-m\right) +\left\vert s-\frac{m+M}{2}%
\right\vert ,$ then 
\begin{align}
F& =\dbigvee\limits_{m-0}^{M}\left( \left\langle E_{\left( \cdot \right)
}x,y\right\rangle \right)   \label{III.a.e.3.4.a} \\
& \times \left\{ 
\begin{array}{cc}
\left[ \frac{1}{2}\left( M-m\right) +\left\vert s-\frac{m+M}{2}\right\vert %
\right] \left\Vert f^{\prime }\right\Vert _{\infty } & \text{if }f^{\prime
}\in L_{\infty }\left[ m,M\right]  \\ 
&  \\ 
\left[ \frac{1}{2}\left( M-m\right) +\left\vert s-\frac{m+M}{2}\right\vert %
\right] ^{1/q}\left\Vert f^{\prime }\right\Vert _{p} & 
\begin{array}{c}
\text{if }f^{\prime }\in L_{p}\left[ m,M\right] ,p>1, \\ 
\frac{1}{p}+\frac{1}{q}=1,%
\end{array}%
\end{array}%
\right.   \notag
\end{align}%
for any $x,y\in H.$

The last part follows by the Total Variation Schwarz's inequality and the
details are omitted.
\end{proof}

\begin{corollary}[Dragomir, 2010, \protect\cite{III.a.SSD3}]
\label{III.a.c.3.1}With the assumptions of Theorem \ref{III.a.t.3.1} we have
the following inequalities%
\begin{align}
& \left\vert f\left( \frac{\left\langle Ax,x\right\rangle }{\left\Vert
x\right\Vert ^{2}}\right) \left\langle x,y\right\rangle -\left\langle
f\left( A\right) x,y\right\rangle \right\vert  \label{III.a.e.3.7} \\
& \leq \left\Vert x\right\Vert \left\Vert y\right\Vert  \notag \\
& \times \left\{ 
\begin{array}{cc}
\left[ \frac{1}{2}\left( M-m\right) +\left\vert \frac{\left\langle
Ax,x\right\rangle }{\left\Vert x\right\Vert ^{2}}-\frac{m+M}{2}\right\vert %
\right] \left\Vert f^{\prime }\right\Vert _{\infty } & \text{if }f^{\prime
}\in L_{\infty }\left[ m,M\right] \\ 
&  \\ 
\left[ \frac{1}{2}\left( M-m\right) +\left\vert \frac{\left\langle
Ax,x\right\rangle }{\left\Vert x\right\Vert ^{2}}-\frac{m+M}{2}\right\vert %
\right] ^{1/q}\left\Vert f^{\prime }\right\Vert _{p} & 
\begin{array}{c}
\text{if }f^{\prime }\in L_{p}\left[ m,M\right] ,p>1, \\ 
\frac{1}{p}+\frac{1}{q}=1,%
\end{array}%
\end{array}%
\right.  \notag
\end{align}%
and%
\begin{align}
& \left\vert f\left( \frac{m+M}{2}\right) \left\langle x,y\right\rangle
-\left\langle f\left( A\right) x,y\right\rangle \right\vert
\label{III.a.e.3.8} \\
& \leq \left\Vert x\right\Vert \left\Vert y\right\Vert  \notag \\
& \times \left\{ 
\begin{array}{cc}
\frac{1}{2}\left( M-m\right) \left\Vert f^{\prime }\right\Vert _{\infty } & 
\text{if }f^{\prime }\in L_{\infty }\left[ m,M\right] \\ 
&  \\ 
\frac{1}{2^{1/q}}\left( M-m\right) ^{1/q}\left\Vert f^{\prime }\right\Vert
_{p} & 
\begin{array}{c}
\text{if }f^{\prime }\in L_{p}\left[ m,M\right] ,p>1, \\ 
\frac{1}{p}+\frac{1}{q}=1,%
\end{array}%
\end{array}%
\right.  \notag
\end{align}%
for any $x,y\in H.$
\end{corollary}

\begin{remark}
\label{III.a.r.3.1}In particular, we obtain from (\ref{III.a.e.2.7}) the
following inequalities%
\begin{align}
& \left\vert f\left( \left\langle Ax,x\right\rangle \right) -\left\langle
f\left( A\right) x,x\right\rangle \right\vert  \label{III.a.e.3.9} \\
& \leq \left\{ 
\begin{array}{cc}
\left[ \frac{1}{2}\left( M-m\right) +\left\vert \left\langle
Ax,x\right\rangle -\frac{m+M}{2}\right\vert \right] \left\Vert f^{\prime
}\right\Vert _{\infty } & \text{if }f^{\prime }\in L_{\infty }\left[ m,M%
\right] \\ 
&  \\ 
\left[ \frac{1}{2}\left( M-m\right) +\left\vert \left\langle
Ax,x\right\rangle -\frac{m+M}{2}\right\vert \right] ^{1/q}\left\Vert
f^{\prime }\right\Vert _{p} & 
\begin{array}{c}
\text{if }f^{\prime }\in L_{p}\left[ m,M\right] , \\ 
p>1,\frac{1}{p}+\frac{1}{q}=1,%
\end{array}%
\end{array}%
\right.  \notag
\end{align}%
and%
\begin{align}
& \left\vert f\left( \frac{m+M}{2}\right) -\left\langle f\left( A\right)
x,x\right\rangle \right\vert  \label{III.a.e.3.10} \\
& \leq \left\{ 
\begin{array}{cc}
\frac{1}{2}\left( M-m\right) \left\Vert f^{\prime }\right\Vert _{\infty } & 
\text{if }f^{\prime }\in L_{\infty }\left[ m,M\right] \\ 
&  \\ 
\frac{1}{2^{1/q}}\left( M-m\right) ^{1/q}\left\Vert f^{\prime }\right\Vert
_{p} & 
\begin{array}{c}
\text{if }f^{\prime }\in L_{p}\left[ m,M\right] ,p>1, \\ 
\frac{1}{p}+\frac{1}{q}=1,%
\end{array}%
\end{array}%
\right.  \notag
\end{align}%
for any $x\in H$ with $\left\Vert x\right\Vert =1.$
\end{remark}

\begin{theorem}[Dragomir, 2010, \protect\cite{III.a.SSD3}]
\label{III.a.t.3.2}Let $A$ be a selfadjoint operator in the Hilbert space $H$
with the spectrum $Sp\left( A\right) \subseteq \left[ m,M\right] $ for some
real numbers $m<M$ and let $\left\{ E_{\lambda }\right\} _{\lambda }$ be its 
\textit{spectral family.} If $f:\left[ m,M\right] \rightarrow \mathbb{R}$ is 
$r-H$-H\"{o}lder continuous on $\left[ m,M\right] $, then we have the
inequality%
\begin{align}
& \left\vert f\left( s\right) \left\langle x,y\right\rangle -\left\langle
f\left( A\right) x,y\right\rangle \right\vert   \label{III.a.e.3.11} \\
& \leq H\dbigvee\limits_{m-0}^{M}\left( \left\langle E_{\left( \cdot \right)
}x,y\right\rangle \right) \left[ \frac{1}{2}\left( M-m\right) +\left\vert s-%
\frac{m+M}{2}\right\vert \right] ^{r}  \notag \\
& \leq H\left\Vert x\right\Vert \left\Vert y\right\Vert \left[ \frac{1}{2}%
\left( M-m\right) +\left\vert s-\frac{m+M}{2}\right\vert \right] ^{r}  \notag
\end{align}%
for any $x,y\in H$ and $s\in \left[ m,M\right] .$

In particular, we have the inequalities%
\begin{align}
& \left\vert f\left( \frac{\left\langle Ax,x\right\rangle }{\left\Vert
x\right\Vert ^{2}}\right) \left\langle x,y\right\rangle -\left\langle
f\left( A\right) x,y\right\rangle \right\vert  \label{III.a.e.3.11.a} \\
& \leq H\left\Vert x\right\Vert \left\Vert y\right\Vert \left[ \frac{1}{2}%
\left( M-m\right) +\left\vert \frac{\left\langle Ax,x\right\rangle }{%
\left\Vert x\right\Vert ^{2}}-\frac{m+M}{2}\right\vert \right] ^{r}  \notag
\end{align}%
and%
\begin{equation}
\left\vert f\left( \frac{m+M}{2}\right) \left\langle x,y\right\rangle
-\left\langle f\left( A\right) x,y\right\rangle \right\vert \leq \frac{1}{%
2^{r}}H\left\Vert x\right\Vert \left\Vert y\right\Vert \left( M-m\right) ^{r}
\label{III.a.e.3.11.b}
\end{equation}%
for any $x,y\in H.$
\end{theorem}

\begin{proof}
Utilising the inequality (\ref{III.a.e.3.4}) and the fact that $f$ is $r-H$-H%
\"{o}lder continuous we have successively 
\begin{align}
& \left\vert f\left( s\right) \left\langle x,y\right\rangle -\left\langle
f\left( A\right) x,y\right\rangle \right\vert   \label{III.a.e.3.12} \\
& =\left\vert \int_{m-0}^{M}\left[ f\left( s\right) -f\left( t\right) \right]
d\left( \left\langle E_{t}x,y\right\rangle \right) \right\vert   \notag \\
& \leq \max_{t\in \left[ m,M\right] }\left\vert f\left( s\right) -f\left(
t\right) \right\vert \dbigvee\limits_{m-0}^{M}\left( \left\langle E_{\left(
\cdot \right) }x,y\right\rangle \right)   \notag \\
& \leq H\max_{t\in \left[ m,M\right] }\left\vert s-t\right\vert
^{r}\dbigvee\limits_{m-0}^{M}\left( \left\langle E_{\left( \cdot \right)
}x,y\right\rangle \right)   \notag \\
& =H\left[ \frac{1}{2}\left( M-m\right) +\left\vert s-\frac{m+M}{2}%
\right\vert \right] ^{r}\dbigvee\limits_{m-0}^{M}\left( \left\langle
E_{\left( \cdot \right) }x,y\right\rangle \right)   \notag
\end{align}%
for any $x,y\in H$ and $s\in \left[ m,M\right] .$

The argument follows now as in the proof of Theorem \ref{III.a.t.3.1} and
the details are omitted.
\end{proof}

\subsection{Logarithmic Inequalities}

Consider \textit{the identric mean}\emph{\ } 
\begin{equation*}
I=I\left( a,b\right) :=\left\{ 
\begin{array}{lll}
a & \text{if } & a=b, \\ 
&  &  \\ 
\frac{1}{e}\left( \frac{b^{b}}{a^{a}}\right) ^{\frac{1}{b-a}} & \text{if} & 
a\neq b,%
\end{array}%
\right. \;a,b>0;
\end{equation*}%
and observe that%
\begin{equation*}
\frac{1}{b-a}\int_{a}^{b}\ln tdt=\ln \left[ I\left( a,b\right) \right] .
\end{equation*}

If we apply Theorem \ref{III.a.t.2.2} for the convex function $f\left(
t\right) =-\ln t,t>0,$ then we can state:

\begin{proposition}
\label{III.a.p.4.1}Let $A$ be a positive selfadjoint operator in the Hilbert
space $H$ with the spectrum $Sp\left( A\right) \subseteq \left[ m,M\right] $
for some positive numbers $0<$ $m<M.$ Then we have the double inequality in
the operator order of $B\left( H\right) $%
\begin{equation}
-\frac{1}{2mM}\left( A^{2}-mM\right) \leq \ln I\left( m,M\right) \cdot
1_{H}-\ln A\leq \frac{m+M}{2}\cdot A^{-1}-1_{H}.  \label{III.a.e.4.1}
\end{equation}
\end{proposition}

If we denote by $G\left( a,b\right) :=\sqrt{ab}$ the geometric mean of the
positive numbers $a,b,$ then we can state the following result as well:

\begin{proposition}
\label{III.a.p.4.2}With the assumptions of Proposition \ref{III.a.p.4.1}, we
have the inequalities in the operator order of $B\left( H\right) $%
\begin{align}
& \ln G\left( m,M\right) \cdot 1_{H}  \label{III.a.e.4.2} \\
& \leq \frac{1}{2}\left[ \ln A+\frac{\ln M\cdot \left( M\cdot 1_{H}-A\right)
+\ln m\cdot \left( A-m\cdot 1_{H}\right) }{M-m}\right]  \notag \\
& \leq \ln I\left( m,M\right) \cdot 1_{H}.  \notag
\end{align}
\end{proposition}

The inequality follows by Corollary \ref{III.a.c.2.3} applied for the convex
function $f\left( t\right) =-\ln t,t>0.$

Finally, the following vector inequality may be stated

\begin{proposition}
\label{III.a.p.4.3}With the assumptions of Proposition \ref{III.a.p.4.1},
for any $x,y\in H$ we have the inequalities%
\begin{align}
& \left\vert \left\langle x,y\right\rangle \ln s-\left\langle \ln
Ax,y\right\rangle \right\vert  \label{III.a.e.4.3} \\
& \leq \left\Vert x\right\Vert \left\Vert y\right\Vert \left\{ 
\begin{array}{ll}
\left[ \frac{1}{2}\left( M-m\right) +\left\vert s-\frac{m+M}{2}\right\vert %
\right] \frac{1}{m}, &  \\ 
&  \\ 
\left[ \frac{1}{2}\left( M-m\right) +\left\vert s-\frac{m+M}{2}\right\vert %
\right] ^{1/q}\frac{M^{p-1}-m^{p-1}}{\left( p-1\right) M^{p-1}m^{p-1}}, & 
\end{array}%
\right.  \notag
\end{align}%
for any $s\in \left[ m,M\right] ,$ where $p>1,\frac{1}{p}+\frac{1}{q}=1.$
\end{proposition}

\section{More Ostrowski's Type Inequalities}

\subsection{Some Vector Inequalities for Functions of Bounded Variation}

The following result holds:

\begin{theorem}[Dragomir, 2010, \protect\cite{III.d.SSD4}]
\label{III.d.t.2.1}Let $A$ be a selfadjoint operator in the Hilbert space $H$
with the spectrum $Sp\left( A\right) \subseteq \left[ m,M\right] $ for some
real numbers $m<M$ and let $\left\{ E_{\lambda }\right\} _{\lambda }$ be its 
\textit{spectral family.} If $f:\left[ m,M\right] \rightarrow \mathbb{R}$ is
a continuous function of bounded variation on $\left[ m,M\right] $, then we
have the inequality%
\begin{align}
& \left\vert f\left( s\right) \left\langle x,y\right\rangle -\left\langle
f\left( A\right) x,y\right\rangle \right\vert   \label{III.d.e.2.2} \\
& \leq \left\langle E_{s}x,x\right\rangle ^{1/2}\left\langle
E_{s}y,y\right\rangle ^{1/2}\dbigvee\limits_{m}^{s}\left( f\right)   \notag
\\
& +\left\langle \left( 1_{H}-E_{s}\right) x,x\right\rangle
^{1/2}\left\langle \left( 1_{H}-E_{s}\right) y,y\right\rangle
^{1/2}\dbigvee\limits_{s}^{M}\left( f\right)   \notag \\
& \leq \left\Vert x\right\Vert \left\Vert y\right\Vert \left( \frac{1}{2}%
\dbigvee\limits_{m}^{M}\left( f\right) +\frac{1}{2}\left\vert
\dbigvee\limits_{m}^{s}\left( f\right) -\dbigvee\limits_{s}^{M}\left(
f\right) \right\vert \right) \left( \leq \left\Vert x\right\Vert \left\Vert
y\right\Vert \dbigvee\limits_{m}^{M}\left( f\right) \right)   \notag
\end{align}%
for any $x,y\in H$ and for any $s\in \left[ m,M\right] .$
\end{theorem}

\begin{proof}
We use the following identity for the Riemann-Stieltjes integral established
by the author in 2000 in \cite{III.d.DR5} (see also \cite[p. 452]{III.DR}):%
\begin{align}
& \left[ u\left( b\right) -u\left( a\right) \right] f\left( s\right)
-\int_{a}^{b}f\left( t\right) du\left( t\right)  \label{III.d.e.2.3} \\
& =\int_{a}^{s}\left[ u\left( t\right) -u\left( a\right) \right] df\left(
t\right) +\int_{s}^{b}\left[ u\left( t\right) -u\left( b\right) \right]
df\left( t\right) ,  \notag
\end{align}%
for any $s\in \left[ a,b\right] ,$ provided the Riemann-Stieltjes integral $%
\int_{a}^{b}f\left( t\right) du\left( t\right) $ exists.

A simple proof can be done by utilizing the integration by parts formula and
starting from the right hand side of (\ref{III.d.e.2.3}).

If we choose in (\ref{III.d.e.2.3}) $a=m,b=M$ and $u\left( t\right)
=\left\langle E_{t}x,y\right\rangle ,$ then we have the following identity
of interest in itself%
\begin{equation}
f\left( s\right) \left\langle x,y\right\rangle -\left\langle f\left(
A\right) x,y\right\rangle =\int_{m-0}^{s}\left\langle E_{t}x,y\right\rangle
df\left( t\right) +\int_{s}^{M}\left\langle \left( E_{t}-1_{H}\right)
x,y\right\rangle df\left( t\right)  \label{III.d.e.2.4}
\end{equation}%
for any $x,y\in H$ and for any $s\in \left[ m,M\right] .$

It is well known that if $p:\left[ a,b\right] \rightarrow \mathbb{C}$ is a
continuous function and $v:\left[ a,b\right] \rightarrow \mathbb{C}$ is of
bounded variation, then the Riemann-Stieltjes integral $\int_{a}^{b}p\left(
t\right) dv\left( t\right) $ exists and the following inequality holds%
\begin{equation*}
\left\vert \int_{a}^{b}p\left( t\right) dv\left( t\right) \right\vert \leq
\max_{t\in \left[ a,b\right] }\left\vert p\left( t\right) \right\vert
\dbigvee\limits_{a}^{b}\left( v\right)
\end{equation*}%
where $\dbigvee\limits_{a}^{b}\left( v\right) $ denotes the total variation
of $v$ on $\left[ a,b\right] .$

Utilising this property we have from (\ref{III.d.e.2.4}) that%
\begin{align}
& \left\vert f\left( s\right) \left\langle x,y\right\rangle -\left\langle
f\left( A\right) x,y\right\rangle \right\vert  \label{III.d.e.2.4.1} \\
& \leq \left\vert \int_{m-0}^{s}\left\langle E_{t}x,y\right\rangle df\left(
t\right) \right\vert +\left\vert \int_{s}^{M}\left\langle \left(
E_{t}-1_{H}\right) x,y\right\rangle df\left( t\right) \right\vert  \notag \\
& \leq \max_{t\in \left[ m,s\right] }\left\vert \left\langle
E_{t}x,y\right\rangle \right\vert \dbigvee\limits_{m}^{s}\left( f\right)
+\max_{t\in \left[ s,M\right] }\left\vert \left\langle \left(
E_{t}-1_{H}\right) x,y\right\rangle \right\vert
\dbigvee\limits_{s}^{M}\left( f\right) :=T  \notag
\end{align}%
for any $x,y\in H$ and for any $s\in \left[ m,M\right] .$

If $P$ is a nonnegative operator on $H,$ i.e., $\left\langle
Px,x\right\rangle \geq 0$ for any $x\in H,$ then the following inequality is
a generalization of the Schwarz inequality in $H$%
\begin{equation}
\left\vert \left\langle Px,y\right\rangle \right\vert ^{2}\leq \left\langle
Px,x\right\rangle \left\langle Py,y\right\rangle  \label{III.d.e.2.4.a}
\end{equation}%
for any $x,y\in H.$

On applying the inequality (\ref{III.d.e.2.4.a}) we have 
\begin{equation*}
\left\vert \left\langle E_{t}x,y\right\rangle \right\vert \leq \left\langle
E_{t}x,x\right\rangle ^{1/2}\left\langle E_{t}y,y\right\rangle ^{1/2}
\end{equation*}%
and%
\begin{equation*}
\left\vert \left\langle \left( 1_{H}-E_{t}\right) x,y\right\rangle
\right\vert \leq \left\langle \left( 1_{H}-E_{t}\right) x,x\right\rangle
^{1/2}\left\langle \left( 1_{H}-E_{t}\right) y,y\right\rangle ^{1/2}
\end{equation*}%
for any $x,y\in H$ and $t\in \left[ m,M\right] .$

Therefore%
\begin{align}
T& \leq \max_{t\in \left[ m,s\right] }\left[ \left\langle
E_{t}x,x\right\rangle ^{1/2}\left\langle E_{t}y,y\right\rangle ^{1/2}\right]
\dbigvee\limits_{m}^{s}\left( f\right)  \label{III.d.e.2.4.b} \\
& +\max_{t\in \left[ s,M\right] }\left[ \left\langle \left(
1_{H}-E_{t}\right) x,x\right\rangle ^{1/2}\left\langle \left(
1_{H}-E_{t}\right) y,y\right\rangle ^{1/2}\right] \dbigvee\limits_{s}^{M}%
\left( f\right)  \notag \\
& \leq \max_{t\in \left[ m,s\right] }\left\langle E_{t}x,x\right\rangle
^{1/2}\max_{t\in \left[ m,s\right] }\left\langle E_{t}y,y\right\rangle
^{1/2}\dbigvee\limits_{m}^{s}\left( f\right)  \notag \\
& +\max_{t\in \left[ s,M\right] }\left\langle \left( 1_{H}-E_{t}\right)
x,x\right\rangle ^{1/2}\max_{t\in \left[ s,M\right] }\left\langle \left(
1_{H}-E_{t}\right) y,y\right\rangle ^{1/2}\dbigvee\limits_{s}^{M}\left(
f\right)  \notag \\
& =\left\langle E_{s}x,x\right\rangle ^{1/2}\left\langle
E_{s}y,y\right\rangle ^{1/2}\dbigvee\limits_{m}^{s}\left( f\right)  \notag \\
& +\left\langle \left( 1_{H}-E_{s}\right) x,x\right\rangle
^{1/2}\left\langle \left( 1_{H}-E_{s}\right) y,y\right\rangle
^{1/2}\dbigvee\limits_{s}^{M}\left( f\right)  \notag \\
& :=V  \notag
\end{align}%
for any $x,y\in H$ and for any $s\in \left[ m,M\right] ,$ proving the first
inequality in (\ref{III.d.e.2.2}).

Now, observe that%
\begin{align*}
V& \leq \max \left\{ \dbigvee\limits_{m}^{s}\left( f\right)
,\dbigvee\limits_{s}^{M}\left( f\right) \right\} \\
& \times \left[ \left\langle E_{s}x,x\right\rangle ^{1/2}\left\langle
E_{s}y,y\right\rangle ^{1/2}+\left\langle \left( 1_{H}-E_{s}\right)
x,x\right\rangle ^{1/2}\left\langle \left( 1_{H}-E_{s}\right)
y,y\right\rangle ^{1/2}\right] .
\end{align*}%
Since%
\begin{equation*}
\max \left\{ \dbigvee\limits_{m}^{s}\left( f\right)
,\dbigvee\limits_{s}^{M}\left( f\right) \right\} =\frac{1}{2}%
\dbigvee\limits_{m}^{M}\left( f\right) +\frac{1}{2}\left\vert
\dbigvee\limits_{m}^{s}\left( f\right) -\dbigvee\limits_{s}^{M}\left(
f\right) \right\vert
\end{equation*}%
and by the Cauchy-Buniakovski-Schwarz inequality for positive real numbers $%
a_{1},b_{1},a_{2},b_{2}$%
\begin{equation}
a_{1}b_{1}+a_{2}b_{2}\leq \left( a_{1}^{2}+a_{2}^{2}\right) ^{1/2}\left(
b_{1}^{2}+b_{2}^{2}\right) ^{1/2}  \label{III.d.e.2.4.c}
\end{equation}%
we have%
\begin{align*}
& \left\langle E_{s}x,x\right\rangle ^{1/2}\left\langle
E_{s}y,y\right\rangle ^{1/2}+\left\langle \left( 1_{H}-E_{s}\right)
x,x\right\rangle ^{1/2}\left\langle \left( 1_{H}-E_{s}\right)
y,y\right\rangle ^{1/2} \\
& \leq \left[ \left\langle E_{s}x,x\right\rangle +\left\langle \left(
1_{H}-E_{s}\right) x,x\right\rangle \right] ^{1/2}\left[ \left\langle
E_{s}y,y\right\rangle +\left\langle \left( 1_{H}-E_{s}\right)
y,y\right\rangle \right] ^{1/2} \\
& =\left\Vert x\right\Vert \left\Vert y\right\Vert
\end{align*}%
for any $x,y\in H$ and $s\in \left[ m,M\right] ,$ then the last part of (\ref%
{III.d.e.2.2}) is proven as well.
\end{proof}

\begin{remark}
\label{III.d.r.2.1}For the continuous function with bounded variation $f:%
\left[ m,M\right] \rightarrow \mathbb{R}$ if $p\in \left[ m,M\right] $ is a
point with the property that%
\begin{equation*}
\dbigvee\limits_{m}^{p}\left( f\right) =\dbigvee\limits_{p}^{M}\left(
f\right)
\end{equation*}%
then from (\ref{III.d.e.2.2}) we get the interesting inequality%
\begin{equation}
\left\vert f\left( p\right) \left\langle x,y\right\rangle -\left\langle
f\left( A\right) x,y\right\rangle \right\vert \leq \frac{1}{2}\left\Vert
x\right\Vert \left\Vert y\right\Vert \dbigvee\limits_{m}^{M}\left( f\right)
\label{III.d.e.2.4.d}
\end{equation}%
for any $x,y\in H.$

If the continuous function $f:\left[ m,M\right] \rightarrow \mathbb{R}$ is
monotonic nondecreasing and therefore of bounded variation, we get from (\ref%
{III.d.e.2.2}) the following inequality as well%
\begin{align}
& \left\vert f\left( s\right) \left\langle x,y\right\rangle -\left\langle
f\left( A\right) x,y\right\rangle \right\vert  \label{III.d.e.2.4.e} \\
& \leq \left\langle E_{s}x,x\right\rangle ^{1/2}\left\langle
E_{s}y,y\right\rangle ^{1/2}\left( f\left( s\right) -f\left( m\right) \right)
\notag \\
& +\left\langle \left( 1_{H}-E_{s}\right) x,x\right\rangle
^{1/2}\left\langle \left( 1_{H}-E_{s}\right) y,y\right\rangle ^{1/2}\left(
f\left( M\right) -f\left( s\right) \right)  \notag \\
& \leq \left\Vert x\right\Vert \left\Vert y\right\Vert \left( \frac{1}{2}%
\left( f\left( M\right) -f\left( m\right) \right) +\left\vert f\left(
s\right) -\frac{f\left( m\right) +f\left( M\right) }{2}\right\vert \right) 
\notag \\
& \left( \leq \left\Vert x\right\Vert \left\Vert y\right\Vert f\left(
M\right) -f\left( m\right) \right)  \notag
\end{align}%
for any $x,y\in H$ and $s\in \left[ m,M\right] .$

Moreover, if the continuous function $f:\left[ m,M\right] \rightarrow 
\mathbb{R}$ is nondecreasing on $\left[ m,M\right] ,$ then the equation 
\begin{equation*}
f\left( s\right) =\frac{f\left( m\right) +f\left( M\right) }{2}
\end{equation*}%
has got at least a solution in $\left[ m,M\right] .$ In his case we get from
(\ref{III.d.e.2.4.e}) the following trapezoidal type inequality%
\begin{equation}
\left\vert \frac{f\left( m\right) +f\left( M\right) }{2}\left\langle
x,y\right\rangle -\left\langle f\left( A\right) x,y\right\rangle \right\vert
\leq \frac{1}{2}\left\Vert x\right\Vert \left\Vert y\right\Vert \left(
f\left( M\right) -f\left( m\right) \right)  \label{III.d.e.2.4.f}
\end{equation}%
for any $x,y\in H.$
\end{remark}

\subsection{Some Vector Inequalities for Lipshitzian Functions}

The following result that incorporates the case of Lipschitzian functions
also holds

\begin{theorem}[Dragomir, 2010, \protect\cite{III.d.SSD4}]
\label{III.d.t.2.2}Let $A$ be a selfadjoint operator in the Hilbert space $H$
with the spectrum $Sp\left( A\right) \subseteq \left[ m,M\right] $ for some
real numbers $m<M$ and let $\left\{ E_{\lambda }\right\} _{\lambda }$ be its 
\textit{spectral family.} If $f:\left[ m,M\right] \rightarrow \mathbb{R}$ is
Lipschitzian with the constant $L>0$ on $\left[ m,M\right] $, i.e., 
\begin{equation*}
\left\vert f\left( s\right) -f\left( t\right) \right\vert \leq L\left\vert
s-t\right\vert \text{ for any }s,t\in \left[ m,M\right] ,
\end{equation*}%
then we have the inequality%
\begin{align}
& \left\vert f\left( s\right) \left\langle x,y\right\rangle -\left\langle
f\left( A\right) x,y\right\rangle \right\vert  \label{III.d.e.2.5} \\
& \leq L\left[ \left( \int_{m-0}^{s}\left\langle E_{t}x,x\right\rangle
dt\right) ^{1/2}\left( \int_{m-0}^{s}\left\langle E_{t}y,y\right\rangle
dt\right) ^{1/2}\right.  \notag \\
& \left. +\left( \int_{s}^{M}\left\langle \left( 1_{H}-E_{t}\right)
x,x\right\rangle dt\right) ^{1/2}\left( \int_{s}^{M}\left\langle \left(
1_{H}-E_{t}\right) y,y\right\rangle dt\right) ^{1/2}\right]  \notag \\
& \leq L\left\langle \left\vert A-s1_{H}\right\vert x,x\right\rangle
^{1/2}\left\langle \left\vert A-s1_{H}\right\vert y,y\right\rangle ^{1/2} 
\notag \\
& \leq L\left[ D^{2}\left( A;x\right) +\left( s\left\Vert x\right\Vert
^{2}-\left\langle Ax,x\right\rangle \right) ^{2}\right] ^{1/4}  \notag \\
& \times \left[ D^{2}\left( A;y\right) +\left( s\left\Vert y\right\Vert
^{2}-\left\langle Ay,y\right\rangle \right) ^{2}\right] ^{1/4}  \notag
\end{align}%
for any $x,y\in H$ and $s\in \left[ m,M\right] ,$ where $D\left( A;x\right) $
is the variance of the selfadjoint operator $A$ in $x$ and is defined by 
\begin{equation*}
D\left( A;x\right) :=\left( \left\Vert Ax\right\Vert ^{2}\left\Vert
x\right\Vert ^{2}-\left\langle Ax,x\right\rangle ^{2}\right) ^{1/2}.
\end{equation*}
\end{theorem}

\begin{proof}
It is well known that if $p:\left[ a,b\right] \rightarrow \mathbb{C}$ is a
Riemann integrable function and $v:\left[ a,b\right] \rightarrow \mathbb{C}$
is Lipschitzian with the constant $L>0$, i.e.,%
\begin{equation*}
\left\vert f\left( s\right) -f\left( t\right) \right\vert \leq L\left\vert
s-t\right\vert \text{ for any }t,s\in \left[ a,b\right] ,
\end{equation*}%
then the Riemann-Stieltjes integral $\int_{a}^{b}p\left( t\right) dv\left(
t\right) $ exists and the following inequality holds%
\begin{equation*}
\left\vert \int_{a}^{b}p\left( t\right) dv\left( t\right) \right\vert \leq
L\int_{a}^{b}\left\vert p\left( t\right) \right\vert dt.
\end{equation*}

Now, on applying this property of the Riemann-Stieltjes integral, we have
from the representation (\ref{III.d.e.2.4}) that%
\begin{align*}
& \left\vert f\left( s\right) \left\langle x,y\right\rangle -\left\langle
f\left( A\right) x,y\right\rangle \right\vert \\
& \leq \left\vert \int_{m-0}^{s}\left\langle E_{t}x,y\right\rangle df\left(
t\right) \right\vert +\left\vert \int_{s}^{M}\left\langle \left(
E_{t}-1_{H}\right) x,y\right\rangle df\left( t\right) \right\vert \\
& \leq L\left[ \int_{m-0}^{s}\left\vert \left\langle E_{t}x,y\right\rangle
\right\vert dt+\int_{s}^{M}\left\vert \left\langle \left( E_{t}-1_{H}\right)
x,y\right\rangle \right\vert dt\right] :=LW
\end{align*}%
for any $x,y\in H$ and $s\in \left[ m,M\right] .$

By utilizing the generalized Schwarz inequality for nonnegative operators (%
\ref{III.d.e.2.4.a}) and the Cauchy-Buniakovski-Schwarz inequality for the
Riemann integral we have%
\begin{align}
W& \leq \int_{m-0}^{s}\left\langle E_{t}x,x\right\rangle ^{1/2}\left\langle
E_{t}y,y\right\rangle ^{1/2}dt  \label{III.d.e.2.7} \\
& +\int_{s}^{M}\left\langle \left( 1_{H}-E_{t}\right) x,x\right\rangle
^{1/2}\left\langle \left( 1_{H}-E_{t}\right) y,y\right\rangle ^{1/2}dt 
\notag \\
& \leq \left( \int_{m-0}^{s}\left\langle E_{t}x,x\right\rangle dt\right)
^{1/2}\left( \int_{m-0}^{s}\left\langle E_{t}y,y\right\rangle dt\right)
^{1/2}  \notag \\
& +\left( \int_{s}^{M}\left\langle \left( 1_{H}-E_{t}\right)
x,x\right\rangle dt\right) ^{1/2}\left( \int_{s}^{M}\left\langle \left(
1_{H}-E_{t}\right) y,y\right\rangle dt\right) ^{1/2}  \notag \\
& :=Z  \notag
\end{align}%
for any $x,y\in H$ and $s\in \left[ m,M\right] .$

On the other hand, by making use of the elementary inequality (\ref%
{III.d.e.2.4.c}) we also have%
\begin{align}
Z& \leq \left( \int_{m-0}^{s}\left\langle E_{t}x,x\right\rangle
dt+\int_{s}^{M}\left\langle \left( 1_{H}-E_{t}\right) x,x\right\rangle
dt\right) ^{1/2}  \label{III.d.e.2.8} \\
& \times \left( \int_{m-0}^{s}\left\langle E_{t}y,y\right\rangle
dt+\int_{s}^{M}\left\langle \left( 1_{H}-E_{t}\right) y,y\right\rangle
dt\right) ^{1/2}  \notag
\end{align}%
for any $x,y\in H$ and $s\in \left[ m,M\right] .$

Now, observe that, by the use of the representation (\ref{III.d.e.2.4}) for
the continuous function $f:\left[ m,M\right] \rightarrow \mathbb{R}$, $%
f\left( t\right) =\left\vert t-s\right\vert $ where $s$ is fixed in $\left[
m,M\right] $ we have the following identity that is of interest in itself%
\begin{equation}
\left\langle \left\vert A-s\cdot 1_{H}\right\vert x,y\right\rangle
=\int_{m-0}^{s}\left\langle E_{t}x,y\right\rangle
dt+\int_{s}^{M}\left\langle \left( 1_{H}-E_{t}\right) x,y\right\rangle dt
\label{III.d.e.2.9}
\end{equation}%
for any $x,y\in H.$

On utilizing (\ref{III.d.e.2.9}) for $x$ and then for $y$ we deduce the
second part of (\ref{III.d.e.2.5}).

Finally, by the well known inequality for the modulus of a bounded linear
operator%
\begin{equation*}
\left\langle \left\vert T\right\vert x,x\right\rangle \leq \left\Vert
Tx\right\Vert \left\Vert x\right\Vert ,x\in H
\end{equation*}%
we have%
\begin{align*}
\left\langle \left\vert A-s\cdot 1_{H}\right\vert x,x\right\rangle ^{1/2}&
\leq \left\Vert Ax-sx\right\Vert ^{1/2}\left\Vert x\right\Vert ^{1/2} \\
& =\left( \left\Vert Ax\right\Vert ^{2}-2s\left\langle Ax,x\right\rangle
+s^{2}\left\Vert x\right\Vert ^{2}\right) ^{1/4}\left\Vert x\right\Vert
^{1/2} \\
& =\left[ \left\Vert Ax\right\Vert ^{2}\left\Vert x\right\Vert
^{2}-\left\langle Ax,x\right\rangle ^{2}+\left( s\left\Vert x\right\Vert
^{2}-\left\langle Ax,x\right\rangle \right) ^{2}\right] ^{1/4} \\
& =\left[ D^{2}\left( A;x\right) +\left( s\left\Vert x\right\Vert
^{2}-\left\langle Ax,x\right\rangle \right) ^{2}\right] ^{1/4}
\end{align*}%
and a similar relation for $y$. The proof is thus complete.
\end{proof}

\begin{remark}
\label{III.d.r.2.2}Since $A$ is a selfadjoint operator in the Hilbert space $%
H$ with the spectrum $Sp\left( A\right) \subseteq \left[ m,M\right] ,$ then 
\begin{equation*}
\left\vert A-\frac{m+M}{2}\cdot 1_{H}\right\vert \leq \frac{M-m}{2}1_{H}
\end{equation*}%
giving from (\ref{III.d.e.2.5}) that%
\begin{align}
& \left\vert f\left( \frac{m+M}{2}\right) \left\langle x,y\right\rangle
-\left\langle f\left( A\right) x,y\right\rangle \right\vert
\label{III.d.e.2.9.a} \\
& \leq L\left[ \left( \int_{m-0}^{\frac{m+M}{2}}\left\langle
E_{t}x,x\right\rangle dt\right) ^{1/2}\left( \int_{m-0}^{\frac{m+M}{2}%
}\left\langle E_{t}y,y\right\rangle dt\right) ^{1/2}\right.  \notag \\
& \left. +\left( \int_{\frac{m+M}{2}}^{M}\left\langle \left(
1_{H}-E_{t}\right) x,x\right\rangle dt\right) ^{1/2}\left( \int_{\frac{m+M}{2%
}}^{M}\left\langle \left( 1_{H}-E_{t}\right) y,y\right\rangle dt\right)
^{1/2}\right]  \notag \\
& \leq L\left\langle \left\vert A-\frac{m+M}{2}\cdot 1_{H}\right\vert
x,x\right\rangle ^{1/2}\left\langle \left\vert A-\frac{m+M}{2}\cdot
1_{H}\right\vert y,y\right\rangle ^{1/2}  \notag \\
& \leq \frac{1}{2}L\left( M-m\right) \left\Vert x\right\Vert \left\Vert
y\right\Vert  \notag
\end{align}%
for any $x,y\in H.$
\end{remark}

The particular case of equal vectors is of interest:

\begin{corollary}[Dragomir, 2010, \protect\cite{III.d.SSD4}]
\label{III.d.c.2.0}Let $A$ be a selfadjoint operator in the Hilbert space $H$
with the spectrum $Sp\left( A\right) \subseteq \left[ m,M\right] $ for some
real numbers $m<M$\textit{.} If $f:\left[ m,M\right] \rightarrow \mathbb{R}$
is Lipschitzian with the constant $L>0$ on $\left[ m,M\right] $, then we
have the inequality%
\begin{align}
\left\vert f\left( s\right) \left\Vert x\right\Vert ^{2}-\left\langle
f\left( A\right) x,x\right\rangle \right\vert & \leq L\left\langle
\left\vert A-s\cdot 1_{H}\right\vert x,x\right\rangle  \label{III.d.e.2.9.b}
\\
& \leq L\left[ D^{2}\left( A;x\right) +\left( s\left\Vert x\right\Vert
^{2}-\left\langle Ax,x\right\rangle \right) ^{2}\right] ^{1/2}  \notag
\end{align}%
for any $x\in H$ and $s\in \left[ m,M\right] .$
\end{corollary}

\begin{remark}
\label{III.d.r.2.3}An important particular case that can be obtained from (%
\ref{III.d.e.2.9.b}) is the one when $s=\frac{\left\langle Ax,x\right\rangle 
}{\left\Vert x\right\Vert ^{2}},x\neq 0,$ giving the inequality%
\begin{align}
\left\vert f\left( \frac{\left\langle Ax,x\right\rangle }{\left\Vert
x\right\Vert ^{2}}\right) \left\Vert x\right\Vert ^{2}-\left\langle f\left(
A\right) x,x\right\rangle \right\vert & \leq L\left\langle \left\vert A-%
\frac{\left\langle Ax,x\right\rangle }{\left\Vert x\right\Vert ^{2}}\cdot
1_{H}\right\vert x,x\right\rangle  \label{III.d.e.2.9.c} \\
& \leq LD\left( A;x\right) \leq \frac{1}{2}L\left( M-m\right) \left\Vert
x\right\Vert ^{2}  \notag
\end{align}%
for any $x\in H,x\neq 0.$
\end{remark}

We are able now to provide the following corollary:

\begin{corollary}[Dragomir, 2010, \protect\cite{III.d.SSD4}]
\label{III.d.c.2.1}Let $A$ be a selfadjoint operator in the Hilbert space $H$
with the spectrum $Sp\left( A\right) \subseteq \left[ m,M\right] $ for some
real numbers $m<M$ and let $\left\{ E_{\lambda }\right\} _{\lambda }$ be its 
\textit{spectral family.} If $f:\left[ m,M\right] \rightarrow \mathbb{R}$ is
a $\left( \varphi ,\Phi \right) -$Lipschitzian functions on $\left[ m,M%
\right] $ with $\Phi >\varphi ,$ then we have the inequality%
\begin{align}
& \left\vert \left\langle f\left( A\right) x,y\right\rangle -\frac{\Phi
+\varphi }{2}\left\langle Ax,y\right\rangle +\frac{\Phi +\varphi }{2}%
s\left\langle x,y\right\rangle -f\left( s\right) \left\langle
x,y\right\rangle \right\vert  \label{III.d.e.2.13} \\
& \leq \frac{1}{2}\left( \Phi -\varphi \right) \left[ \left(
\int_{m-0}^{s}\left\langle E_{t}x,x\right\rangle dt\right) ^{1/2}\left(
\int_{m-0}^{s}\left\langle E_{t}y,y\right\rangle dt\right) ^{1/2}\right. 
\notag \\
& \left. +\left( \int_{s}^{M}\left\langle \left( 1_{H}-E_{t}\right)
x,x\right\rangle dt\right) ^{1/2}\left( \int_{s}^{M}\left\langle \left(
1_{H}-E_{t}\right) y,y\right\rangle dt\right) ^{1/2}\right]  \notag \\
& \leq \frac{1}{2}\left( \Phi -\varphi \right) \left\langle \left\vert
A-s1_{H}\right\vert x,x\right\rangle ^{1/2}\left\langle \left\vert
A-s1_{H}\right\vert y,y\right\rangle ^{1/2}  \notag \\
& \leq \frac{1}{2}\left( \Phi -\varphi \right) \left[ D^{2}\left( A;x\right)
+\left( s\left\Vert x\right\Vert ^{2}-\left\langle Ax,x\right\rangle \right)
^{2}\right] ^{1/4}  \notag \\
& \times \left[ D^{2}\left( A;y\right) +\left( s\left\Vert y\right\Vert
^{2}-\left\langle Ay,y\right\rangle \right) ^{2}\right] ^{1/4}  \notag
\end{align}%
for any $x,y\in H.$
\end{corollary}

\begin{remark}
\label{III.d.r.2.4}Various particular cases can be stated by utilizing the
inequality (\ref{III.d.e.2.13}), however the details are left to the
interested reader.
\end{remark}

\section{Some Vector Inequalities for Monotonic Functions}

The case of monotonic functions is of interest as well. The corresponding
result is incorporated in the following

\begin{theorem}[Dragomir, 2010, \protect\cite{III.d.SSD4}]
\label{III.d.t.3.1}Let $A$ be a selfadjoint operator in the Hilbert space $H$
with the spectrum $Sp\left( A\right) \subseteq \left[ m,M\right] $ for some
real numbers $m<M$ and let $\left\{ E_{\lambda }\right\} _{\lambda }$ be its 
\textit{spectral family.} If $f:\left[ m,M\right] \rightarrow \mathbb{R}$ is
a continuous monotonic nondecreasing function on $\left[ m,M\right] $, then
we have the inequality%
\begin{align}
& \left\vert f\left( s\right) \left\langle x,y\right\rangle -\left\langle
f\left( A\right) x,y\right\rangle \right\vert  \label{III.d.e.3.1} \\
& \leq \left( \int_{m-0}^{s}\left\langle E_{t}x,x\right\rangle df\left(
t\right) \right) ^{1/2}\left( \int_{m-0}^{s}\left\langle
E_{t}y,y\right\rangle df\left( t\right) \right) ^{1/2}  \notag \\
& +\left( \int_{s}^{M}\left\langle \left( 1_{H}-E_{t}\right)
x,x\right\rangle df\left( t\right) \right) ^{1/2}\left(
\int_{s}^{M}\left\langle \left( 1_{H}-E_{t}\right) y,y\right\rangle df\left(
t\right) \right) ^{1/2}  \notag \\
& \leq \left\langle \left\vert f\left( A\right) -f\left( s\right)
1_{H}\right\vert x,x\right\rangle ^{1/2}\left\langle \left\vert f\left(
A\right) -f\left( s\right) 1_{H}\right\vert y,y\right\rangle ^{1/2}  \notag
\\
& \leq \left[ D^{2}\left( f\left( A\right) ;x\right) +\left( f\left(
s\right) \left\Vert x\right\Vert ^{2}-\left\langle f\left( A\right)
x,x\right\rangle \right) ^{2}\right] ^{1/4}  \notag \\
& \times \left[ D^{2}\left( f\left( A\right) ;y\right) +\left( f\left(
s\right) \left\Vert y\right\Vert ^{2}-\left\langle f\left( A\right)
y,y\right\rangle \right) ^{2}\right] ^{1/4}  \notag
\end{align}%
for any $x,y\in H$ and $s\in \left[ m,M\right] ,$ where, as above $D\left(
f\left( A\right) ;x\right) $ is the variance of the selfadjoint operator $%
f\left( A\right) $ in $x.$
\end{theorem}

\begin{proof}
From the theory of Riemann-Stieltjes integral is well known that if $p:\left[
a,b\right] \rightarrow \mathbb{C}$ is of bounded variation and $v:\left[ a,b%
\right] \rightarrow \mathbb{R}$ is continuous and monotonic nondecreasing,
then the Riemann-Stieltjes integrals $\int_{a}^{b}p\left( t\right) dv\left(
t\right) $ and $\int_{a}^{b}\left\vert p\left( t\right) \right\vert dv\left(
t\right) $ exist and%
\begin{equation*}
\left\vert \int_{a}^{b}p\left( t\right) dv\left( t\right) \right\vert \leq
\int_{a}^{b}\left\vert p\left( t\right) \right\vert dv\left( t\right) .
\end{equation*}%
On utilizing this property and the representation (\ref{III.d.e.2.4}) we
have successively%
\begin{align}
& \left\vert f\left( s\right) \left\langle x,y\right\rangle -\left\langle
f\left( A\right) x,y\right\rangle \right\vert  \label{III.d.e.3.2} \\
& \leq \left\vert \int_{m-0}^{s}\left\langle E_{t}x,y\right\rangle df\left(
t\right) \right\vert +\left\vert \int_{s}^{M}\left\langle \left(
E_{t}-1_{H}\right) x,y\right\rangle df\left( t\right) \right\vert  \notag \\
& \leq \int_{m-0}^{s}\left\vert \left\langle E_{t}x,y\right\rangle
\right\vert df\left( t\right) +\int_{s}^{M}\left\vert \left\langle \left(
E_{t}-1_{H}\right) x,y\right\rangle \right\vert df\left( t\right)  \notag \\
& \leq \int_{m-0}^{s}\left\langle E_{t}x,x\right\rangle ^{1/2}\left\langle
E_{t}y,y\right\rangle ^{1/2}df\left( t\right)  \notag \\
& +\int_{s}^{M}\left\langle \left( 1_{H}-E_{t}\right) x,x\right\rangle
^{1/2}\left\langle \left( 1_{H}-E_{t}\right) y,y\right\rangle ^{1/2}df\left(
t\right)  \notag \\
& :=Y,  \notag
\end{align}%
for any $x,y\in H$ and $s\in \left[ m,M\right] .$

We use now the following version of the Cauchy-Buniakovski-Schwarz
inequality for the Riemann-Stieltjes integral with monotonic nondecreasing
integrators 
\begin{equation*}
\left( \int_{a}^{b}p\left( t\right) q\left( t\right) dv\left( t\right)
\right) ^{2}\leq \int_{a}^{b}p^{2}\left( t\right) dv\left( t\right)
\int_{a}^{b}q^{2}\left( t\right) dv\left( t\right)
\end{equation*}%
to get that%
\begin{equation*}
\int_{m-0}^{s}\left\langle E_{t}x,x\right\rangle ^{1/2}\left\langle
E_{t}y,y\right\rangle ^{1/2}df\left( t\right) \leq \left(
\int_{m-0}^{s}\left\langle E_{t}x,x\right\rangle df\left( t\right) \right)
^{1/2}\left( \int_{m-0}^{s}\left\langle E_{t}y,y\right\rangle df\left(
t\right) \right) ^{1/2}
\end{equation*}%
and%
\begin{align*}
& \int_{s}^{M}\left\langle \left( 1_{H}-E_{t}\right) x,x\right\rangle
^{1/2}\left\langle \left( 1_{H}-E_{t}\right) y,y\right\rangle ^{1/2}df\left(
t\right) \\
& \leq \left( \int_{s}^{M}\left\langle \left( 1_{H}-E_{t}\right)
x,x\right\rangle df\left( t\right) \right) ^{1/2}\left(
\int_{s}^{M}\left\langle \left( 1_{H}-E_{t}\right) y,y\right\rangle df\left(
t\right) \right) ^{1/2}
\end{align*}%
for any $x,y\in H$ and $s\in \left[ m,M\right] .$

Therefore%
\begin{align*}
Y& \leq \left( \int_{m-0}^{s}\left\langle E_{t}x,x\right\rangle df\left(
t\right) \right) ^{1/2}\left( \int_{m-0}^{s}\left\langle
E_{t}y,y\right\rangle df\left( t\right) \right) ^{1/2} \\
& +\left( \int_{s}^{M}\left\langle \left( 1_{H}-E_{t}\right)
x,x\right\rangle df\left( t\right) \right) ^{1/2}\left(
\int_{s}^{M}\left\langle \left( 1_{H}-E_{t}\right) y,y\right\rangle df\left(
t\right) \right) ^{1/2} \\
& \leq \left( \int_{m-0}^{s}\left\langle E_{t}x,x\right\rangle df\left(
t\right) +\int_{s}^{M}\left\langle \left( 1_{H}-E_{t}\right)
x,x\right\rangle df\left( t\right) \right) ^{1/2} \\
& \times \left( \int_{m-0}^{s}\left\langle E_{t}y,y\right\rangle df\left(
t\right) +\int_{s}^{M}\left\langle \left( 1_{H}-E_{t}\right)
y,y\right\rangle df\left( t\right) \right) ^{1/2}
\end{align*}%
for any $x,y\in H$ and $s\in \left[ m,M\right] ,$ where, to get the last
inequality we have used the elementary inequality (\ref{III.d.e.2.4.c}).

Now, since $f$ is monotonic nondecreasing, on applying the representation (%
\ref{III.d.e.2.4}) for the function $\left\vert f\left( \cdot \right)
-f\left( s\right) \right\vert $ with $s$ fixed in $\left[ m,M\right] $ we
deduce the following identity that is of interest in itself as well:%
\begin{equation}
\left\langle \left\vert f\left( A\right) -f\left( s\right) \right\vert
x,y\right\rangle =\int_{m-0}^{s}\left\langle E_{t}x,y\right\rangle df\left(
t\right) +\int_{s}^{M}\left\langle \left( 1_{H}-E_{t}\right)
x,y\right\rangle df\left( t\right)  \label{III.d.e.3.3}
\end{equation}%
for any $x,y\in H.$

The second part of (\ref{III.d.e.3.1}) follows then by writing (\ref%
{III.d.e.3.3}) for $x$ then by $y$ and utilizing the relevant inequalities
from above.

The last part is similar to the corresponding one from the proof of Theorem %
\ref{III.d.t.2.2} and the details are omitted.
\end{proof}

The following corollary is of interest:

\begin{corollary}[Dragomir, 2010, \protect\cite{III.d.SSD4}]
\label{III.d.c.3.1}With the assumption of Theorem \ref{III.d.t.3.1} we have
the inequalities%
\begin{align}
& \left\vert \frac{f\left( m\right) +f\left( M\right) }{2}\left\langle
x,y\right\rangle -\left\langle f\left( A\right) x,y\right\rangle \right\vert
\label{III.d.e.3.4} \\
& \leq \left\langle \left\vert f\left( A\right) -\frac{f\left( m\right)
+f\left( M\right) }{2}\cdot 1_{H}\right\vert x,x\right\rangle ^{1/2}  \notag
\\
& \times \left\langle \left\vert f\left( A\right) -\frac{f\left( m\right)
+f\left( M\right) }{2}\cdot 1_{H}\right\vert y,y\right\rangle ^{1/2}  \notag
\\
& \leq \frac{1}{2}\left( f\left( M\right) -f\left( m\right) \right)
\left\Vert x\right\Vert \left\Vert y\right\Vert ,  \notag
\end{align}%
for any $x,y\in H.$
\end{corollary}

\begin{proof}
Since $f$ is monotonic nondecreasing, then $f\left( u\right) \in \left[
f\left( m\right) ,f\left( M\right) \right] $ for any $u\in \left[ m,M\right]
.$ By the continuity of $f$ it follows that there exists at list one $s\in %
\left[ m,M\right] $ such that 
\begin{equation*}
f\left( s\right) =\frac{f\left( m\right) +f\left( M\right) }{2}.
\end{equation*}%
Now, on utilizing the inequality (\ref{III.d.e.3.1}) for this $s$ we deduce
the first inequality in (\ref{III.d.e.3.4}). The second part follows as
above and the details are omitted.
\end{proof}

\subsection{Power Inequalities}

We consider the power function $f\left( t\right) :=t^{p}$ where $p\in 
\mathbb{R}\diagdown \left\{ 0\right\} $ and $t>0.$ The following mid-point
inequalities hold:

\begin{proposition}
\label{III.d.p.4.1}Let $A$ be a selfadjoint operator in the Hilbert space $H$
with the spectrum $Sp\left( A\right) \subseteq \left[ m,M\right] $ for some
real numbers with $0\leq m<M$.

If $p>0,$ then for any $x,y\in H$%
\begin{align}
& \left\vert \left( \frac{m+M}{2}\right) ^{p}\left\langle x,y\right\rangle
-\left\langle A^{p}x,y\right\rangle \right\vert  \label{III.d.e.4.1} \\
& \leq B_{p}\left\langle \left\vert A-\frac{m+M}{2}\cdot 1_{H}\right\vert
x,x\right\rangle ^{1/2}\left\langle \left\vert A-\frac{m+M}{2}\cdot
1_{H}\right\vert y,y\right\rangle ^{1/2}  \notag \\
& \leq \frac{1}{2}B_{p}\left( M-m\right) \left\Vert x\right\Vert \left\Vert
y\right\Vert  \notag
\end{align}%
where%
\begin{equation*}
B_{p}=p\times \left\{ 
\begin{array}{cc}
M^{p-1} & \text{if }p\geq 1 \\ 
&  \\ 
m^{p-1} & \text{if }0<p<1,m>0.%
\end{array}%
\right.
\end{equation*}%
and%
\begin{align}
& \left\vert \left( \frac{m+M}{2}\right) ^{-p}\left\langle x,y\right\rangle
-\left\langle A^{-p}x,y\right\rangle \right\vert  \label{III.d.e.4.2} \\
& \leq C_{p}\left\langle \left\vert A-\frac{m+M}{2}\cdot 1_{H}\right\vert
x,x\right\rangle ^{1/2}\left\langle \left\vert A-\frac{m+M}{2}\cdot
1_{H}\right\vert y,y\right\rangle ^{1/2}  \notag \\
& \leq \frac{1}{2}C_{p}\left( M-m\right) \left\Vert x\right\Vert \left\Vert
y\right\Vert  \notag
\end{align}%
where%
\begin{equation*}
C_{p}=pm^{-p-1}\text{ and }m>0.
\end{equation*}
\end{proposition}

The proof follows from (\ref{III.d.e.2.9.a}).

We can also state the following trapezoidal type inequalities:

\begin{proposition}
\label{III.d.p.4.2}With the assumption of Proposition \ref{III.d.p.4.1} and
if $p>0$ we have the inequalities%
\begin{align}
& \left\vert \frac{m^{p}+M^{p}}{2}\left\langle x,y\right\rangle
-\left\langle A^{p}x,y\right\rangle \right\vert  \label{III.d.e.4.3} \\
& \leq \left\langle \left\vert A^{p}-\frac{m^{p}+M^{p}}{2}\cdot
1_{H}\right\vert x,x\right\rangle ^{1/2}\left\langle \left\vert A^{p}-\frac{%
m^{p}+M^{p}}{2}\cdot 1_{H}\right\vert y,y\right\rangle ^{1/2}  \notag \\
& \leq \frac{1}{2}\left( M^{p}-m^{p}\right) \left\Vert x\right\Vert
\left\Vert y\right\Vert ,  \notag
\end{align}%
and, for $m>0,$%
\begin{align}
& \left\vert \frac{m^{p}+M^{p}}{2m^{p}M^{p}}\left\langle x,y\right\rangle
-\left\langle A^{-p}x,y\right\rangle \right\vert  \label{III.d.e.4.4} \\
& \leq \left\langle \left\vert A^{-p}-\frac{m^{p}+M^{p}}{2m^{p}M^{p}}\cdot
1_{H}\right\vert x,x\right\rangle ^{1/2}\left\langle \left\vert A^{-p}-\frac{%
m^{p}+M^{p}}{2m^{p}M^{p}}\cdot 1_{H}\right\vert y,y\right\rangle ^{1/2} 
\notag \\
& \leq \frac{1}{2}\left( \frac{M^{p}-m^{p}}{M^{p}m^{p}}\right) \left\Vert
x\right\Vert \left\Vert y\right\Vert ,  \notag
\end{align}%
for any $x,y\in H.$
\end{proposition}

The proof follows from Corollary \ref{III.d.c.3.1}.

\subsection{Logarithmic Inequalities}

Consider the function $f\left( t\right) =\ln t,t>0.$ Denote by $A\left(
a,b\right) :=\frac{a+b}{2}$ the arithmetic mean of $a,b>0$ and $G\left(
a,b\right) :=\sqrt{ab}$ the geometric mean of these numbers. We have the
following result:

\begin{proposition}
\label{III.d.p.5.1}Let $A$ be a selfadjoint operator in the Hilbert space $H$
with the spectrum $Sp\left( A\right) \subseteq \left[ m,M\right] $ for some
real numbers with $0<m<M$. For any $x,y\in H$ we have 
\begin{align}
& \left\vert \ln A\left( m,M\right) \cdot \left\langle x,y\right\rangle
-\left\langle \ln Ax,y\right\rangle \right\vert  \label{III.d.e.5.1} \\
& \leq \frac{1}{m}\left\langle \left\vert A-\frac{m+M}{2}\cdot
1_{H}\right\vert x,x\right\rangle ^{1/2}\left\langle \left\vert A-\frac{m+M}{%
2}\cdot 1_{H}\right\vert y,y\right\rangle ^{1/2}  \notag \\
& \leq \frac{1}{2}\left( \frac{M}{m}-1\right) \left\Vert x\right\Vert
\left\Vert y\right\Vert  \notag
\end{align}%
and%
\begin{align}
& \left\vert \ln G\left( m,M\right) \cdot \left\langle x,y\right\rangle
-\left\langle \ln Ax,y\right\rangle \right\vert  \label{III.d.e.5.2} \\
& \leq \left\langle \left\vert \ln A-\ln G\left( m,M\right) \cdot
1_{H}\right\vert x,x\right\rangle ^{1/2}\left\langle \left\vert \ln A-\ln
G\left( m,M\right) \cdot 1_{H}\right\vert y,y\right\rangle ^{1/2}  \notag \\
& \leq \ln \sqrt{\frac{M}{m}}\cdot \left\Vert x\right\Vert \left\Vert
y\right\Vert .  \notag
\end{align}
\end{proposition}

The proof follows by (\ref{III.d.e.2.9.a}) and (\ref{III.d.e.3.4}).

\section{Ostrowski's Type Vector Inequalities}

\subsection{Some Vector Inequalities}

The following result holds:

\begin{theorem}[Dragomir, 2010, \protect\cite{III.b.SSD5}]
\label{III.b.t.3.1}Let $A$ be a selfadjoint operator in the Hilbert space $H$
with the spectrum $Sp\left( A\right) \subseteq \left[ m,M\right] $ for some
real numbers $m<M$ and let $\left\{ E_{\lambda }\right\} _{\lambda }$ be its 
\textit{spectral family.} If $f:\left[ m,M\right] \rightarrow \mathbb{C}$ is
a continuous function of bounded variation on $\left[ m,M\right] $, then we
have the inequality%
\begin{align}
& \left\vert \left\langle x,y\right\rangle \frac{1}{M-m}\int_{m}^{M}f\left(
s\right) ds-\left\langle f\left( A\right) x,y\right\rangle \right\vert
\label{III.b.e.3.1} \\
& \leq \frac{1}{M-m}\dbigvee\limits_{m}^{M}\left( f\right) \max_{t\in \left[
m,M\right] }\left[ \left( M-t\right) \left\langle E_{t}x,x\right\rangle
^{1/2}\left\langle E_{t}y,y\right\rangle ^{1/2}\right.  \notag \\
& \left. +\left( t-m\right) \left\langle \left( 1_{H}-E_{t}\right)
x,x\right\rangle ^{1/2}\left\langle \left( 1_{H}-E_{t}\right)
y,y\right\rangle ^{1/2}\right]  \notag \\
& \leq \left\Vert x\right\Vert \left\Vert y\right\Vert
\dbigvee\limits_{m}^{M}\left( f\right)  \notag
\end{align}%
for any $x,y\in H.$
\end{theorem}

\begin{proof}
Assume that $f:\left[ m,M\right] \rightarrow \mathbb{C}$ is a continuous
function on $\left[ m,M\right] .$ Then under the assumptions of the theorem
for $A$ and $\left\{ E_{\lambda }\right\} _{\lambda },$ we have the
following representation%
\begin{align}
& \left\langle x,y\right\rangle \frac{1}{M-m}\int_{m}^{M}f\left( s\right)
ds-\left\langle f\left( A\right) x,y\right\rangle  \label{III.b.e.3.2} \\
& =\frac{1}{M-m}\int_{m-0}^{M}\left\langle \left[ \left( M-t\right)
E_{t}+\left( t-m\right) \left( E_{t}-1_{H}\right) \right] x,y\right\rangle
df\left( t\right)  \notag
\end{align}%
for any $x,y\in H.$

Indeed, integrating by parts in the Riemann-Stieltjes integral and using the
spectral representation theorem we have%
\begin{align*}
& \frac{1}{M-m}\int_{m-0}^{M}\left\langle \left[ \left( M-t\right)
E_{t}+\left( t-m\right) \left( E_{t}-1_{H}\right) \right] x,y\right\rangle
df\left( t\right) \\
& =\int_{m-0}^{M}\left( \left\langle E_{t}x,y\right\rangle -\frac{t-m}{M-m}%
\left\langle x,y\right\rangle \right) df\left( t\right) \\
& =\left. \left( \left\langle E_{t}x,y\right\rangle -\frac{t-m}{M-m}%
\left\langle x,y\right\rangle \right) f\left( t\right) \right\vert _{m-0}^{M}
\\
& -\int_{m-0}^{M}f\left( t\right) d\left( \left\langle E_{t}x,y\right\rangle
-\frac{t-m}{M-m}\left\langle x,y\right\rangle \right) \\
& =-\int_{m-0}^{M}f\left( t\right) d\left\langle E_{t}x,y\right\rangle
+\left\langle x,y\right\rangle \frac{1}{M-m}\int_{m}^{M}f\left( t\right) dt
\\
& =\left\langle x,y\right\rangle \frac{1}{M-m}\int_{m}^{M}f\left( t\right)
dt-\left\langle f\left( A\right) x,y\right\rangle
\end{align*}%
for any $x,y\in H$ and the equality (\ref{III.b.e.3.2}) is proved.

It is well known that if $p:\left[ a,b\right] \rightarrow \mathbb{C}$ is a
continuous function and $v:\left[ a,b\right] \rightarrow \mathbb{C}$ is of
bounded variation, then the Riemann-Stieltjes integral $\int_{a}^{b}p\left(
t\right) dv\left( t\right) $ exists and the following inequality holds%
\begin{equation*}
\left\vert \int_{a}^{b}p\left( t\right) dv\left( t\right) \right\vert \leq
\max_{t\in \left[ a,b\right] }\left\vert p\left( t\right) \right\vert
\dbigvee\limits_{a}^{b}\left( v\right)
\end{equation*}%
where $\dbigvee\limits_{a}^{b}\left( v\right) $ denotes the total variation
of $v$ on $\left[ a,b\right] .$

Utilising this property we have from (\ref{III.b.e.3.2}) that%
\begin{align}
& \left\vert \left\langle x,y\right\rangle \frac{1}{M-m}\int_{m}^{M}f\left(
s\right) ds-\left\langle f\left( A\right) x,y\right\rangle \right\vert
\label{III.b.e.3.3} \\
& \leq \frac{1}{M-m}\max_{t\in \left[ m,M\right] }\left\vert \left\langle %
\left[ \left( M-t\right) E_{t}+\left( t-m\right) \left( E_{t}-1_{H}\right) %
\right] x,y\right\rangle \right\vert \dbigvee\limits_{m}^{M}\left( f\right) 
\notag
\end{align}%
for any $x,y\in H.$

Now observe that%
\begin{align}
& \left\vert \left\langle \left[ \left( M-t\right) E_{t}+\left( t-m\right)
\left( E_{t}-1_{H}\right) \right] x,y\right\rangle \right\vert
\label{III.b.e.3.3.a} \\
& =\left\vert \left( M-t\right) \left\langle E_{t}x,y\right\rangle +\left(
t-m\right) \left\langle \left( E_{t}-1_{H}\right) x,y\right\rangle
\right\vert  \notag \\
& \leq \left( M-t\right) \left\vert \left\langle E_{t}x,y\right\rangle
\right\vert +\left( t-m\right) \left\vert \left\langle \left(
E_{t}-1_{H}\right) x,y\right\rangle \right\vert  \notag
\end{align}%
for any $x,y\in H$ and $t\in \left[ m,M\right] .$

If $P$ is a nonnegative operator on $H,$ i.e., $\left\langle
Px,x\right\rangle \geq 0$ for any $x\in H,$ then the following inequality is
a generalization of the Schwarz inequality in $H$%
\begin{equation}
\left\vert \left\langle Px,y\right\rangle \right\vert ^{2}\leq \left\langle
Px,x\right\rangle \left\langle Py,y\right\rangle  \label{III.b.e.3.4}
\end{equation}%
for any $x,y\in H.$

On applying the inequality (\ref{III.b.e.3.4}) we have%
\begin{align}
& \left( M-t\right) \left\vert \left\langle E_{t}x,y\right\rangle
\right\vert +\left( t-m\right) \left\vert \left\langle \left(
E_{t}-1_{H}\right) x,y\right\rangle \right\vert  \label{III.b.e.3.5} \\
& \leq \left( M-t\right) \left\langle E_{t}x,x\right\rangle
^{1/2}\left\langle E_{t}y,y\right\rangle ^{1/2}  \notag \\
& +\left( t-m\right) \left\langle \left( 1_{H}-E_{t}\right) x,x\right\rangle
^{1/2}\left\langle \left( 1_{H}-E_{t}\right) y,y\right\rangle ^{1/2}  \notag
\\
& \leq \max \left\{ M-t,t-m\right\}  \notag \\
& \times \left[ \left\langle E_{t}x,x\right\rangle ^{1/2}\left\langle
E_{t}y,y\right\rangle ^{1/2}+\left\langle \left( 1_{H}-E_{t}\right)
x,x\right\rangle ^{1/2}\left\langle \left( 1_{H}-E_{t}\right)
y,y\right\rangle ^{1/2}\right]  \notag \\
& \leq \max \left\{ M-t,t-m\right\}  \notag \\
& \times \left[ \left\langle E_{s}x,x\right\rangle +\left\langle \left(
1_{H}-E_{s}\right) x,x\right\rangle \right] ^{1/2}\left[ \left\langle
E_{s}y,y\right\rangle +\left\langle \left( 1_{H}-E_{s}\right)
y,y\right\rangle \right] ^{1/2}  \notag \\
& =\max \left\{ M-t,t-m\right\} \left\Vert x\right\Vert \left\Vert
y\right\Vert ,  \notag
\end{align}%
where for the last inequality we used the elementary fact%
\begin{equation}
a_{1}b_{1}+a_{2}b_{2}\leq \left( a_{1}^{2}+a_{2}^{2}\right) ^{1/2}\left(
b_{1}^{2}+b_{2}^{2}\right) ^{1/2}  \label{III.b.e.3.6}
\end{equation}%
that holds for $a_{1},b_{1},a_{2},b_{2}$ positive real numbers.

Utilising the inequalities (\ref{III.b.e.3.3}), (\ref{III.b.e.3.3.a}) and (%
\ref{III.b.e.3.5}) we deduce the desired result (\ref{III.b.e.3.1}).
\end{proof}

The case of Lipschitzian functions is embodied in the following result:

\begin{theorem}[Dragomir, 2010, \protect\cite{III.b.SSD5}]
\label{III.b.t.3.2}Let $A$ be a selfadjoint operator in the Hilbert space $H$
with the spectrum $Sp\left( A\right) \subseteq \left[ m,M\right] $ for some
real numbers $m<M$ and let $\left\{ E_{\lambda }\right\} _{\lambda }$ be its 
\textit{spectral family.} If $f:\left[ m,M\right] \rightarrow \mathbb{C}$ is
a Lipschitzian function with the constant $L>0$ on $\left[ m,M\right] $,
then we have the inequality%
\begin{align}
& \left\vert \left\langle x,y\right\rangle \frac{1}{M-m}\int_{m}^{M}f\left(
s\right) ds-\left\langle f\left( A\right) x,y\right\rangle \right\vert
\label{III.b.e.3.7} \\
& \leq \frac{L}{M-m}\int_{m}^{M}\left[ \left( M-t\right) \left\langle
E_{t}x,x\right\rangle ^{1/2}\left\langle E_{t}y,y\right\rangle ^{1/2}\right.
\notag \\
& \left. +\left( t-m\right) \left\langle \left( 1_{H}-E_{t}\right)
x,x\right\rangle ^{1/2}\left\langle \left( 1_{H}-E_{t}\right)
y,y\right\rangle ^{1/2}\right] dt  \notag \\
& \leq \frac{3}{4}L\left( M-m\right) \left\Vert x\right\Vert \left\Vert
y\right\Vert  \notag
\end{align}%
for any $x,y\in H.$
\end{theorem}

\begin{proof}
It is well known that if $p:\left[ a,b\right] \rightarrow \mathbb{C}$ is a
Riemann integrable function and $v:\left[ a,b\right] \rightarrow \mathbb{C}$
is Lipschitzian with the constant $L>0$, i.e.,%
\begin{equation*}
\left\vert f\left( s\right) -f\left( t\right) \right\vert \leq L\left\vert
s-t\right\vert \text{ for any }t,s\in \left[ a,b\right] ,
\end{equation*}%
then the Riemann-Stieltjes integral $\int_{a}^{b}p\left( t\right) dv\left(
t\right) $ exists and the following inequality holds%
\begin{equation*}
\left\vert \int_{a}^{b}p\left( t\right) dv\left( t\right) \right\vert \leq
L\int_{a}^{b}\left\vert p\left( t\right) \right\vert dt.
\end{equation*}

Now, on applying this property of the Riemann-Stieltjes integral, we have
from the representation (\ref{III.b.e.3.2}) that%
\begin{align}
& \left\vert \left\langle x,y\right\rangle \frac{1}{M-m}\int_{m}^{M}f\left(
s\right) ds-\left\langle f\left( A\right) x,y\right\rangle \right\vert
\label{III.b.e.3.8} \\
& \leq \frac{L}{M-m}\int_{m-0}^{M}\left\vert \left\langle \left[ \left(
M-t\right) E_{t}+\left( t-m\right) \left( E_{t}-1_{H}\right) \right]
x,y\right\rangle \right\vert dt.  \notag
\end{align}%
Since, from the proof of Theorem \ref{III.b.t.3.1}, we have%
\begin{align}
& \left\vert \left\langle \left[ \left( M-t\right) E_{t}+\left( t-m\right)
\left( E_{t}-1_{H}\right) \right] x,y\right\rangle \right\vert
\label{III.b.e.3.9} \\
& \leq \left( M-t\right) \left\langle E_{t}x,x\right\rangle
^{1/2}\left\langle E_{t}y,y\right\rangle ^{1/2}  \notag \\
& +\left( t-m\right) \left\langle \left( 1_{H}-E_{t}\right) x,x\right\rangle
^{1/2}\left\langle \left( 1_{H}-E_{t}\right) y,y\right\rangle ^{1/2}  \notag
\\
& \leq \max \left\{ M-t,t-m\right\} \left\Vert x\right\Vert \left\Vert
y\right\Vert  \notag \\
& =\left[ \frac{1}{2}\left( M-m\right) +\left\vert t-\frac{m+M}{2}%
\right\vert \right] \left\Vert x\right\Vert \left\Vert y\right\Vert  \notag
\end{align}%
for any $x,y\in H$ and $t\in \left[ m,M\right] ,$ then integrating (\ref%
{III.b.e.3.9}) and taking into account that 
\begin{equation*}
\int_{m}^{M}\left\vert t-\frac{m+M}{2}\right\vert dt=\frac{1}{4}\left(
M-m\right) ^{2}
\end{equation*}%
we deduce the desired result (\ref{III.b.e.3.7}).
\end{proof}

Finally for the section, we provide here the case of monotonic nondecreasing
functions as well:

\begin{theorem}[Dragomir, 2010, \protect\cite{III.b.SSD5}]
\label{III.b.t.3.3}Let $A$ be a selfadjoint operator in the Hilbert space $H$
with the spectrum $Sp\left( A\right) \subseteq \left[ m,M\right] $ for some
real numbers $m<M$ and let $\left\{ E_{\lambda }\right\} _{\lambda }$ be its 
\textit{spectral family.} If $f:\left[ m,M\right] \rightarrow \mathbb{R}$ is
a continuous monotonic nondecreasing function on $\left[ m,M\right] $, then
we have the inequality%
\begin{align}
& \left\vert \left\langle x,y\right\rangle \frac{1}{M-m}\int_{m}^{M}f\left(
s\right) ds-\left\langle f\left( A\right) x,y\right\rangle \right\vert
\label{III.b.e.3.10} \\
& \leq \frac{1}{M-m}\int_{m}^{M}\left[ \left( M-t\right) \left\langle
E_{t}x,x\right\rangle ^{1/2}\left\langle E_{t}y,y\right\rangle ^{1/2}\right.
\notag \\
& \left. +\left( t-m\right) \left\langle \left( 1_{H}-E_{t}\right)
x,x\right\rangle ^{1/2}\left\langle \left( 1_{H}-E_{t}\right)
y,y\right\rangle ^{1/2}\right] df\left( t\right)  \notag \\
& \leq \left[ f\left( M\right) -f\left( m\right) -\frac{1}{M-m}%
\int_{m}^{M}sgn\left( t-\frac{m+M}{2}\right) f\left( t\right) dt\right]
\left\Vert x\right\Vert \left\Vert y\right\Vert  \notag \\
& \leq \left[ f\left( M\right) -f\left( m\right) \right] \left\Vert
x\right\Vert \left\Vert y\right\Vert  \notag
\end{align}%
for any $x,y\in H.$
\end{theorem}

\begin{proof}
From the theory of Riemann-Stieltjes integral is well known that if $p:\left[
a,b\right] \rightarrow \mathbb{C}$ is of bounded variation and $v:\left[ a,b%
\right] \rightarrow \mathbb{R}$ is continuous and monotonic nondecreasing,
then the Riemann-Stieltjes integrals $\int_{a}^{b}p\left( t\right) dv\left(
t\right) $ and $\int_{a}^{b}\left\vert p\left( t\right) \right\vert dv\left(
t\right) $ exist and%
\begin{equation*}
\left\vert \int_{a}^{b}p\left( t\right) dv\left( t\right) \right\vert \leq
\int_{a}^{b}\left\vert p\left( t\right) \right\vert dv\left( t\right) .
\end{equation*}%
Now, on applying this property of the Riemann-Stieltjes integral, we have
from the representation (\ref{III.b.e.3.2}) that%
\begin{align}
& \left\vert \left\langle x,y\right\rangle \frac{1}{M-m}\int_{m}^{M}f\left(
s\right) ds-\left\langle f\left( A\right) x,y\right\rangle \right\vert
\label{III.b.e.3.11} \\
& \leq \frac{1}{M-m}\int_{m-0}^{M}\left\vert \left\langle \left[ \left(
M-t\right) E_{t}+\left( t-m\right) \left( E_{t}-1_{H}\right) \right]
x,y\right\rangle \right\vert df\left( t\right) .  \notag
\end{align}%
Further on, by utilizing the inequality (\ref{III.b.e.3.9}) we also have
that 
\begin{align}
& \int_{m-0}^{M}\left\vert \left\langle \left[ \left( M-t\right)
E_{t}+\left( t-m\right) \left( E_{t}-1_{H}\right) \right] x,y\right\rangle
\right\vert df\left( t\right)  \label{III.b.e.3.12} \\
& \leq \int_{m}^{M}\left[ \left( M-t\right) \left\langle
E_{t}x,x\right\rangle ^{1/2}\left\langle E_{t}y,y\right\rangle ^{1/2}\right.
\notag \\
& \left. +\left( t-m\right) \left\langle \left( 1_{H}-E_{t}\right)
x,x\right\rangle ^{1/2}\left\langle \left( 1_{H}-E_{t}\right)
y,y\right\rangle ^{1/2}\right] df\left( t\right)  \notag \\
& \leq \left[ \frac{1}{2}\left( M-m\right) \left[ f\left( M\right) -f\left(
m\right) \right] +\int_{m}^{M}\left\vert t-\frac{m+M}{2}\right\vert df\left(
t\right) \right] \left\Vert x\right\Vert \left\Vert y\right\Vert .  \notag
\end{align}

Now, integrating by parts in the Riemann-Stieltjes integral we have%
\begin{align*}
& \int_{m}^{M}\left\vert t-\frac{m+M}{2}\right\vert df\left( t\right) \\
& =\int_{m}^{\frac{M+m}{2}}\left( \frac{m+M}{2}-t\right) df\left( t\right)
+\int_{\frac{m+M}{2}}^{M}\left( t-\frac{m+M}{2}\right) df\left( t\right) \\
& =\left. \left( \frac{m+M}{2}-t\right) f\left( t\right) \right\vert _{m}^{%
\frac{M+m}{2}}+\int_{m}^{\frac{M+m}{2}}f\left( t\right) dt \\
& +\left. \left( t-\frac{m+M}{2}\right) f\left( t\right) \right\vert _{\frac{%
m+M}{2}}^{M}-\int_{\frac{m+M}{2}}^{M}f\left( t\right) dt \\
& =\frac{1}{2}\left( M-m\right) \left[ f\left( M\right) -f\left( m\right) %
\right] -\int_{m}^{M}sgn\left( t-\frac{m+M}{2}\right) f\left( t\right) dt,
\end{align*}%
which together with (\ref{III.b.e.3.12}) produces the second inequality in (%
\ref{III.b.e.3.10}).

Since the functions $sgn\left( \cdot -\frac{m+M}{2}\right) $ and $f\left(
\cdot \right) $ have the same monotonicity, then by the \v{C}eby\v{s}ev
inequality we have%
\begin{align*}
& \int_{m}^{M}sgn\left( t-\frac{m+M}{2}\right) f\left( t\right) dt \\
& \geq \frac{1}{M-m}\int_{m}^{M}sgn\left( t-\frac{m+M}{2}\right)
dt\int_{m}^{M}f\left( t\right) dt=0
\end{align*}%
and the last part of (\ref{III.b.e.3.10}) is proved.
\end{proof}

\subsection{Applications for Particular Functions}

It is obvious that the above results can be applied for various particular
functions. However, we will restrict here only to the power and logarithmic
functions.

\textbf{1.} Consider now the power function $f:\left( 0,\infty \right)
\rightarrow \mathbb{R}$, $f\left( t\right) =t^{p}$ with $p>0.$ This function
is monotonic increasing on $\left( 0,\infty \right) $ and applying Theorem %
\ref{III.b.t.3.3} we can state the following proposition:

\begin{proposition}
\label{III.b.p.4.1}Let $A$ be a selfadjoint operator in the Hilbert space $H$
with the spectrum $Sp\left( A\right) \subseteq \left[ m,M\right] $ for some
real numbers $0<m<M$ and let $\left\{ E_{\lambda }\right\} _{\lambda }$ be
its \textit{spectral family. Then }for any $x,y\in H$ we have the
inequalities%
\begin{align}
& \left\vert \left\langle A^{p}x,y\right\rangle -\frac{M^{p+1}-m^{p+1}}{%
\left( p+1\right) \left( M-m\right) }\left\langle x,y\right\rangle
\right\vert  \label{III.b.e.4.1} \\
& \leq \frac{p}{M-m}\int_{m}^{M}\left[ \left( M-t\right) \left\langle
E_{t}x,x\right\rangle ^{1/2}\left\langle E_{t}y,y\right\rangle ^{1/2}\right.
\notag \\
& \left. +\left( t-m\right) \left\langle \left( 1_{H}-E_{t}\right)
x,x\right\rangle ^{1/2}\left\langle \left( 1_{H}-E_{t}\right)
y,y\right\rangle ^{1/2}\right] t^{p-1}dt  \notag \\
& \leq \left[ M^{p}-m^{p}-\frac{M^{p+1}+m^{p+1}-2^{p}\left( M+m\right) ^{p+1}%
}{\left( p+1\right) \left( M-m\right) }\right] \left\Vert x\right\Vert
\left\Vert y\right\Vert .  \notag
\end{align}
\end{proposition}

On applying now Theorem \ref{III.b.t.3.2} to the same power function, then
we can state the following result as well:

\begin{proposition}
\label{III.b.p.4.2}With the same assumptions from Proposition \ref%
{III.b.p.4.1} we have%
\begin{align}
& \left\vert \left\langle A^{p}x,y\right\rangle -\frac{M^{p+1}-m^{p+1}}{%
\left( p+1\right) \left( M-m\right) }\left\langle x,y\right\rangle
\right\vert  \label{III.b.e.4.2} \\
& \leq \frac{B_{p}}{M-m}\int_{m}^{M}\left[ \left( M-t\right) \left\langle
E_{t}x,x\right\rangle ^{1/2}\left\langle E_{t}y,y\right\rangle ^{1/2}\right.
\notag \\
& \left. +\left( t-m\right) \left\langle \left( 1_{H}-E_{t}\right)
x,x\right\rangle ^{1/2}\left\langle \left( 1_{H}-E_{t}\right)
y,y\right\rangle ^{1/2}\right] dt  \notag \\
& \leq \frac{3}{4}B_{p}\left( M-m\right) \left\Vert x\right\Vert \left\Vert
y\right\Vert  \notag
\end{align}%
for any $x,y\in H$, where%
\begin{equation*}
B_{p}=p\times \left\{ 
\begin{array}{cc}
M^{p-1} & \text{if }p\geq 1 \\ 
&  \\ 
m^{p-1} & \text{if }0<p<1,m>0.%
\end{array}%
\right.
\end{equation*}
\end{proposition}

The case of negative powers except $p=-1$ goes likewise and we omit the
details.

Now, if we apply Theorem \ref{III.b.t.3.3} and \ref{III.b.t.3.2} for the
increasing function $f\left( t\right) =-\frac{1}{t}$ with $t>0,$ then we can
state the following proposition:

\begin{proposition}
\label{III.b.p.4.3}Let $A$ be a selfadjoint operator in the Hilbert space $H$
with the spectrum $Sp\left( A\right) \subseteq \left[ m,M\right] $ for some
real numbers $0<m<M$ and let $\left\{ E_{\lambda }\right\} _{\lambda }$ be
its \textit{spectral family. Then }for any $x,y\in H$ we have the
inequalities%
\begin{align}
& \left\vert \left\langle A^{-1}x,y\right\rangle -\frac{\ln M-\ln m}{M-m}%
\left\langle x,y\right\rangle \right\vert  \label{III.b.e.4.3} \\
& \leq \frac{1}{M-m}\int_{m}^{M}\left[ \left( M-t\right) \left\langle
E_{t}x,x\right\rangle ^{1/2}\left\langle E_{t}y,y\right\rangle ^{1/2}\right.
\notag \\
& \left. +\left( t-m\right) \left\langle \left( 1_{H}-E_{t}\right)
x,x\right\rangle ^{1/2}\left\langle \left( 1_{H}-E_{t}\right)
y,y\right\rangle ^{1/2}\right] t^{2}dt  \notag \\
& \leq \left[ \frac{M-m}{mM}-\frac{\ln \left[ \left( \frac{m+M}{2}\right)
^{2}\right] -\ln \left( mM\right) }{M-m}\right] \left\Vert x\right\Vert
\left\Vert y\right\Vert  \notag
\end{align}%
and%
\begin{align}
& \left\vert \left\langle A^{-1}x,y\right\rangle -\frac{\ln M-\ln m}{M-m}%
\left\langle x,y\right\rangle \right\vert  \label{III.b.e.4.4} \\
& \leq \frac{1}{m^{2}\left( M-m\right) }\int_{m}^{M}\left[ \left( M-t\right)
\left\langle E_{t}x,x\right\rangle ^{1/2}\left\langle E_{t}y,y\right\rangle
^{1/2}\right.  \notag \\
& \left. +\left( t-m\right) \left\langle \left( 1_{H}-E_{t}\right)
x,x\right\rangle ^{1/2}\left\langle \left( 1_{H}-E_{t}\right)
y,y\right\rangle ^{1/2}\right] dt  \notag \\
& \leq \frac{3}{4}\frac{M-m}{m^{2}}\left\Vert x\right\Vert \left\Vert
y\right\Vert .  \notag
\end{align}
\end{proposition}

\textbf{2. }Now, if we apply Theorems \ref{III.b.t.3.3} and \ref{III.b.t.3.2}
to the function $f:\left( 0,\infty \right) \rightarrow \mathbb{R}$, $f\left(
t\right) =\ln t$, then we can state

\begin{proposition}
\label{III.b.p.4.4}Let $A$ be a selfadjoint operator in the Hilbert space $H$
with the spectrum $Sp\left( A\right) \subseteq \left[ m,M\right] $ for some
real numbers $0<m<M$ and let $\left\{ E_{\lambda }\right\} _{\lambda }$ be
its \textit{spectral family. Then }for any $x,y\in H$ we have the
inequalities%
\begin{align}
& \left\vert \left\langle \ln Ax,y\right\rangle -\left\langle
x,y\right\rangle \ln I\left( m,M\right) \right\vert  \label{III.b.e.4.5} \\
& \leq \frac{1}{M-m}\int_{m}^{M}\left[ \left( M-t\right) \left\langle
E_{t}x,x\right\rangle ^{1/2}\left\langle E_{t}y,y\right\rangle ^{1/2}\right.
\notag \\
& \left. +\left( t-m\right) \left\langle \left( 1_{H}-E_{t}\right)
x,x\right\rangle ^{1/2}\left\langle \left( 1_{H}-E_{t}\right)
y,y\right\rangle ^{1/2}\right] tdt  \notag \\
& \leq \left[ \ln \left( \frac{M}{m}\right) -\ln \left( \sqrt{\frac{I\left( 
\frac{m+M}{2},M\right) }{I\left( m,\frac{m+M}{2}\right) }}\right) \right]
\left\Vert x\right\Vert \left\Vert y\right\Vert  \notag
\end{align}%
and%
\begin{align}
& \left\vert \left\langle \ln Ax,y\right\rangle -\left\langle
x,y\right\rangle \ln I\left( m,M\right) \right\vert  \label{III.b.e.4.6} \\
& \leq \frac{1}{m\left( M-m\right) }\int_{m}^{M}\left[ \left( M-t\right)
\left\langle E_{t}x,x\right\rangle ^{1/2}\left\langle E_{t}y,y\right\rangle
^{1/2}\right.  \notag \\
& \left. +\left( t-m\right) \left\langle \left( 1_{H}-E_{t}\right)
x,x\right\rangle ^{1/2}\left\langle \left( 1_{H}-E_{t}\right)
y,y\right\rangle ^{1/2}\right] dt  \notag \\
& \leq \frac{3}{4}\left( \frac{M}{m}-1\right) \left\Vert x\right\Vert
\left\Vert y\right\Vert ,  \notag
\end{align}%
where $I\left( m,M\right) $ is the identric mean of $m$ and $M$ and is
defined by%
\begin{equation*}
I\left( m,M\right) =\frac{1}{e}\left( \frac{M^{M}}{m^{m}}\right) ^{1/(M-m)}.
\end{equation*}
\end{proposition}

\section{Bounds for the Difference Between Functions and Integral Means}

\subsection{Vector Inequalities Via Ostrowski's Type Bounds}

The following result holds:

\begin{theorem}[Dragomir, 2010, \protect\cite{III.c.SSD5}]
\label{III.c.t.2.1}Let $A$ be a selfadjoint operator in the Hilbert space $H$
with the spectrum $Sp\left( A\right) \subseteq \left[ m,M\right] $ for some
real numbers $m<M$ and let $\left\{ E_{\lambda }\right\} _{\lambda }$ be its 
\textit{spectral family.} If $f:\left[ m,M\right] \rightarrow \mathbb{R}$ is
a continuous function on $\left[ m,M\right] $, then we have the inequality%
\begin{align}
& \left\vert \left\langle f\left( A\right) x,y\right\rangle -\left\langle
x,y\right\rangle \frac{1}{M-m}\int_{m}^{M}f\left( s\right) ds\right\vert 
\label{III.c.e.2.1} \\
& \leq \max_{t\in \left[ m,M\right] }\left\vert f\left( t\right) -\frac{1}{%
M-m}\int_{m}^{M}f\left( s\right) ds\right\vert
\dbigvee\limits_{m-0}^{M}\left( \left\langle E_{\left( \cdot \right)
}x,y\right\rangle \right)   \notag \\
& \leq \max_{t\in \left[ m,M\right] }\left\vert f\left( t\right) -\frac{1}{%
M-m}\int_{m}^{M}f\left( s\right) ds\right\vert \left\Vert x\right\Vert
\left\Vert y\right\Vert   \notag
\end{align}%
for any $x,y\in H.$
\end{theorem}

\begin{proof}
Utilising the spectral representation theorem we have the following equality
of interest 
\begin{align}
& \left\langle f\left( A\right) x,y\right\rangle -\left\langle
x,y\right\rangle \frac{1}{M-m}\int_{m}^{M}f\left( s\right) ds
\label{III.c.e.2.2} \\
& =\int_{m-0}^{M}\left[ f\left( t\right) -\frac{1}{M-m}\int_{m}^{M}f\left(
s\right) ds\right] d\left( \left\langle E_{t}x,y\right\rangle \right)  \notag
\end{align}%
for any $x,y\in H.$

It is well known that if $p:\left[ a,b\right] \rightarrow \mathbb{C}$ is a
continuous function and $v:\left[ a,b\right] \rightarrow \mathbb{C}$ is of
bounded variation, then the Riemann-Stieltjes integral $\int_{a}^{b}p\left(
t\right) dv\left( t\right) $ exists and the following inequality holds%
\begin{equation}
\left\vert \int_{a}^{b}p\left( t\right) dv\left( t\right) \right\vert \leq
\max_{t\in \left[ a,b\right] }\left\vert p\left( t\right) \right\vert
\dbigvee\limits_{a}^{b}\left( v\right) ,  \label{III.c.e.2.2.1}
\end{equation}%
where $\dbigvee\limits_{a}^{b}\left( v\right) $ denotes the total variation
of $v$ on $\left[ a,b\right] .$

Utilising these two facts we get the first part of (\ref{III.c.e.2.1}).

The last part follows by the Total Variation Schwarz's inequality and we
omit the details.
\end{proof}

For particular classes of continuous functions $f:\left[ m,M\right]
\rightarrow \mathbb{C}$ we are able to provide simpler bounds as
incorporated in the following corollary:

\begin{corollary}[Dragomir, 2010, \protect\cite{III.c.SSD5}]
\label{III.c.c.2.1}Let $A$ be a selfadjoint operator in the Hilbert space $H$
with the spectrum $Sp\left( A\right) \subseteq \left[ m,M\right] $ for some
real numbers $m<M,$ $\left\{ E_{\lambda }\right\} _{\lambda }$ be its 
\textit{spectral family and} $f:\left[ m,M\right] \rightarrow \mathbb{C}$ a
continuous function on $\left[ m,M\right] .$

1. If $f$ is of bounded variation on $\left[ m,M\right] ,$ then%
\begin{align}
& \left\vert \left\langle f\left( A\right) x,y\right\rangle -\left\langle
x,y\right\rangle \frac{1}{M-m}\int_{m}^{M}f\left( s\right) ds\right\vert 
\label{III.c.e.2.6} \\
& \leq \dbigvee\limits_{m}^{M}\left( f\right)
\dbigvee\limits_{m-0}^{M}\left( \left\langle E_{\left( \cdot \right)
}x,y\right\rangle \right) \leq \left\Vert x\right\Vert \left\Vert
y\right\Vert \dbigvee\limits_{m}^{M}\left( f\right)   \notag
\end{align}%
for any $x,y\in H.$

2. If $f:\left[ m,M\right] \longrightarrow \mathbb{C}$ is of $r-H-$H\"{o}%
lder type, i.e., for a given $r\in (0,1]$ and $H>0$ we have 
\begin{equation}
\left\vert f\left( s\right) -f\left( t\right) \right\vert \leq H\left\vert
s-t\right\vert ^{r}\text{ for any }s,t\in \left[ m,M\right] ,
\label{III.c.e.2.7}
\end{equation}%
then we have the inequality:%
\begin{align}
& \left\vert \left\langle f\left( A\right) x,y\right\rangle -\left\langle
x,y\right\rangle \frac{1}{M-m}\int_{m}^{M}f\left( s\right) ds\right\vert 
\label{III.c.e.2.8} \\
& \leq \frac{1}{r+1}H\left( M-m\right) ^{r}\dbigvee\limits_{m-0}^{M}\left(
\left\langle E_{\left( \cdot \right) }x,y\right\rangle \right) \leq \frac{1}{%
r+1}H\left( M-m\right) ^{r}\left\Vert x\right\Vert \left\Vert y\right\Vert  
\notag
\end{align}%
for any $x,y\in H.$

In particular, if $f:\left[ m,M\right] \longrightarrow \mathbb{C}$ is
Lipschitzian with the constant $L>0,$ then%
\begin{align}
& \left\vert \left\langle f\left( A\right) x,y\right\rangle -\left\langle
x,y\right\rangle \frac{1}{M-m}\int_{m}^{M}f\left( s\right) ds\right\vert 
\label{III.c.e.2.9} \\
& \leq \frac{1}{2}L\left( M-m\right) \dbigvee\limits_{m-0}^{M}\left(
\left\langle E_{\left( \cdot \right) }x,y\right\rangle \right) \leq \frac{1}{%
2}L\left( M-m\right) \left\Vert x\right\Vert \left\Vert y\right\Vert   \notag
\end{align}%
for any $x,y\in H.$

3. If $f:\left[ m,M\right] \longrightarrow \mathbb{C}$ is absolutely
continuous, then 
\begin{align}
& \left\vert \left\langle f\left( A\right) x,y\right\rangle -\left\langle
x,y\right\rangle \frac{1}{M-m}\int_{m}^{M}f\left( s\right) ds\right\vert 
\label{III.c.e.2.10} \\
& \leq \dbigvee\limits_{m-0}^{M}\left( \left\langle E_{\left( \cdot \right)
}x,y\right\rangle \right)   \notag \\
& \times \left\{ 
\begin{array}{ll}
\frac{1}{2}\left( M-m\right) \left\Vert f^{\prime }\right\Vert _{\infty } & 
\text{if }f^{\prime }\in L_{\infty }\left[ m,M\right]  \\ 
&  \\ 
\frac{1}{\left( q+1\right) ^{1/q}}\left( M-m\right) ^{1/q}\left\Vert
f^{\prime }\right\Vert _{p} & 
\begin{array}{c}
\text{if }f^{\prime }\in L_{p}\left[ m,M\right]  \\ 
p>1,1/p+1/q=1;%
\end{array}
\\ 
&  \\ 
\left\Vert f^{\prime }\right\Vert _{1} & 
\end{array}%
\right.   \notag \\
& \leq \left\Vert x\right\Vert \left\Vert y\right\Vert   \notag \\
& \times \left\{ 
\begin{array}{ll}
\frac{1}{2}\left( M-m\right) \left\Vert f^{\prime }\right\Vert _{\infty } & 
\text{if }f^{\prime }\in L_{\infty }\left[ m,M\right]  \\ 
&  \\ 
\frac{1}{\left( q+1\right) ^{1/q}}\left( M-m\right) ^{1/q}\left\Vert
f^{\prime }\right\Vert _{p} & 
\begin{array}{c}
\text{if }f^{\prime }\in L_{p}\left[ m,M\right]  \\ 
p>1,1/p+1/q=1;%
\end{array}
\\ 
&  \\ 
\left\Vert f^{\prime }\right\Vert _{1} & 
\end{array}%
\right.   \notag
\end{align}%
for any $x,y\in H,$ where $\left\Vert f^{\prime }\right\Vert _{p}$ are the
Lebesgue norms, i.e., we recall that%
\begin{equation*}
\left\Vert f^{\prime }\right\Vert _{p}:=\left\{ 
\begin{array}{cc}
ess\sup_{s\in \left[ m,M\right] }\left\vert f^{\prime }\left( s\right)
\right\vert  & \text{if }p=\infty ; \\ 
&  \\ 
\left( \int_{m}^{M}\left\vert f\left( s\right) \right\vert ^{p}ds\right)
^{1/p} & \text{if }p\geq 1.%
\end{array}%
\right. 
\end{equation*}
\end{corollary}

\begin{proof}
We use the Ostrowski type inequalities in order to provide upper bounds for
the quantity%
\begin{equation*}
\max_{t\in \left[ m,M\right] }\left\vert f\left( t\right) -\frac{1}{M-m}%
\int_{m}^{M}f\left( s\right) ds\right\vert
\end{equation*}%
where $f:\left[ m,M\right] \longrightarrow \mathbb{C}$ is a continuous
function.

The following result may be stated (see \cite{III.b.DR2}) for functions of
bounded variation:

\begin{lemma}
\label{III.c.l.2.1}Assume that $f:\left[ m,M\right] \rightarrow \mathbb{C}$
is of bounded variation and denote by $\bigvee\limits_{m}^{M}\left( f\right) 
$ its total variation. Then 
\begin{equation}
\left\vert f\left( t\right) -\frac{1}{M-m}\int_{m}^{M}f\left( s\right)
ds\right\vert \leq \left[ \frac{1}{2}+\left\vert \frac{t-\frac{m+M}{2}}{M-m}%
\right\vert \right] \bigvee\limits_{m}^{M}\left( f\right)
\label{III.c.e.2.11}
\end{equation}%
for all $t\in \left[ m,M\right] $. The constant $\frac{1}{2}$ is the best
possible.
\end{lemma}

Now, taking the maximum over $x\in \left[ m,M\right] $ in (\ref{III.c.e.2.11}%
) we deduce (\ref{III.c.e.2.6}).

If $f$ is H\"{o}lder continuous, then one may state the result:

\begin{lemma}
\label{III.c.tc}Let $f:\left[ m,M\right] \rightarrow \mathbb{C}$ be of $r-H-$%
H\"{o}lder type, where $r\in (0,1]$ and $H>0$ are fixed, then, for all $x\in %
\left[ m,M\right] ,$ we have the inequality: 
\begin{align}
& \left\vert f\left( t\right) -\frac{1}{M-m}\int_{m}^{M}f\left( s\right)
ds\right\vert  \label{III.c.e.2.12} \\
& \leq \frac{H}{r+1}\left[ \left( \frac{M-t}{M-m}\right) ^{r+1}+\left( \frac{%
t-m}{M-m}\right) ^{r+1}\right] \left( M-m\right) ^{r}.  \notag
\end{align}%
The constant $\frac{1}{r+1}$ is also sharp in the above sense.
\end{lemma}

Note that if $r=1$, i.e., $f$ is Lipschitz continuous, then we get the
following version of Ostrowski's inequality for Lipschitzian functions (with 
$L$ instead of $H$) (see for instance \cite{III.c.DR1}) 
\begin{equation}
\left\vert f\left( t\right) -\frac{1}{M-m}\int_{m}^{M}f\left( s\right)
ds\right\vert \leq \left[ \frac{1}{4}+\left( \frac{t-\frac{m+M}{2}}{M-m}%
\right) ^{2}\right] \left( M-m\right) L,  \label{III.c.e.2.13}
\end{equation}%
for any $x\in \left[ m,M\right] .$ Here the constant $\frac{1}{4}$ is also
best.

Taking the maximum over $x\in \left[ m,M\right] $ in (\ref{III.c.e.2.12}) we
deduce (\ref{III.c.e.2.8}) and the second part of the corollary is proved.

The following Ostrowski type result for absolutely continuous functions
holds.

\begin{lemma}
\label{III.c.tb}Let $f:\left[ a,b\right] \rightarrow \mathbb{R}$ be
absolutely continuous on $\left[ a,b\right] $. Then, for all $t\in \left[ a,b%
\right] $, we have: 
\begin{multline}
\left\vert f\left( t\right) -\frac{1}{M-m}\int_{m}^{M}f\left( s\right)
ds\right\vert  \label{III.c.e.2.14} \\
\leq \left\{ 
\begin{array}{lll}
\left[ \frac{1}{4}+\left( \frac{t-\frac{m+M}{2}}{M-m}\right) ^{2}\right]
\left( M-m\right) \left\Vert f^{\prime }\right\Vert _{\infty } & \text{if} & 
f^{\prime }\in L_{\infty }\left[ m,M\right] ; \\ 
&  &  \\ 
\frac{1}{\left( q+1\right) ^{\frac{1}{q}}}\left[ \left( \frac{t-m}{M-m}%
\right) ^{q+1}+\left( \frac{M-t}{M-m}\right) ^{q+1}\right] ^{\frac{1}{q}%
}\left( M-m\right) ^{\frac{1}{q}}\left\Vert f^{\prime }\right\Vert _{p} & 
\text{if} & f^{\prime }\in L_{p}\left[ m,M\right] , \\ 
&  & \frac{1}{p}+\frac{1}{q}=1,\;p>1; \\ 
\left[ \frac{1}{2}+\left\vert \frac{t-\frac{m+M}{2}}{M-m}\right\vert \right]
\left\Vert f^{\prime }\right\Vert _{1}. &  & 
\end{array}%
\right.
\end{multline}%
The constants $\frac{1}{4}$, $\frac{1}{\left( p+1\right) ^{\frac{1}{p}}}$
and $\frac{1}{2}$ respectively are sharp in the sense presented above.
\end{lemma}

The above inequalities can also be obtained from the Fink result in \cite%
{III.a.F} on choosing $n=1$ and performing some appropriate computations.

Taking the maximum in these inequalities we deduce (\ref{III.c.e.2.10}).
\end{proof}

For other scalar Ostrowski's type inequalities, see \cite{III.c.A} and \cite%
{III.c.7ab}.

\subsection{Other Vector Inequalities}

In \cite{III.c.DF1}, the authors have considered the following functional%
\begin{equation}
D\left( f;u\right) :=\int_{a}^{b}f\left( s\right) du\left( s\right) -\left[
u\left( b\right) -u\left( a\right) \right] \cdot \frac{1}{b-a}%
\int_{a}^{b}f\left( t\right) dt,  \label{III.c.e.3.1}
\end{equation}%
provided that the Stieltjes integral $\int_{a}^{b}f\left( s\right) du\left(
s\right) $ exists.

This functional plays an important role in approximating the Stieltjes
integral $\int_{a}^{b}f\left( s\right) du\left( s\right) $ in terms of the
Riemann integral $\int_{a}^{b}f\left( t\right) dt$ and the divided
difference of the integrator $u.$

In \cite{III.c.DF1}, the following result in estimating the above functional 
$D\left( f;u\right) $ has been obtained:%
\begin{equation}
\left\vert D\left( f;u\right) \right\vert \leq \frac{1}{2}L\left( M-m\right)
\left( b-a\right) ,  \label{III.c.e.3.2}
\end{equation}%
provided $u$ is $L-$\textit{Lipschitzian} and $f$ is\textit{\ Riemann
integrable} and with the property that there exists the constants $m,M\in 
\mathbb{R}$ such that%
\begin{equation}
m\leq f\left( t\right) \leq M\quad \text{for any\quad }t\in \left[ a,b\right]
.  \label{III.c.e.3.3}
\end{equation}%
The constant $\frac{1}{2}$ is best possible in (\ref{III.c.e.3.2}) in the
sense that it cannot be replaced by a smaller quantity.

If one assumes that $u$ is of \textit{bounded variation} and $f$ is $K-$%
\textit{Lipschitzian}, then $D\left( f,u\right) $ satisfies the inequality 
\cite{III.c.DF2}%
\begin{equation}
\left\vert D\left( f;u\right) \right\vert \leq \frac{1}{2}K\left( b-a\right)
\bigvee_{a}^{b}\left( u\right) .  \label{III.c.e.3.4}
\end{equation}%
Here the constant $\frac{1}{2}$ is also best possible.

Now, for the function $u:\left[ a,b\right] \rightarrow \mathbb{C},$ consider
the following auxiliary mappings $\Phi ,\Gamma $ and $\Delta $ \cite{III.c.D}%
:%
\begin{align*}
\Phi \left( t\right) & :=\frac{\left( t-a\right) u\left( b\right) +\left(
b-t\right) u\left( a\right) }{b-a}-u\left( t\right) ,\qquad t\in \left[ a,b%
\right] , \\
\Gamma \left( t\right) & :=\left( t-a\right) \left[ u\left( b\right)
-u\left( t\right) \right] -\left( b-t\right) \left[ u\left( t\right)
-u\left( a\right) \right] ,\qquad t\in \left[ a,b\right] , \\
\Delta \left( t\right) & :=\left[ u;b,t\right] -\left[ u;t,a\right] ,\qquad
t\in \left( a,b\right) ,
\end{align*}%
where $\left[ u;\alpha ,\beta \right] $ is the \textit{divided difference}
of $u$ in $\alpha ,\beta ,$ i.e.,%
\begin{equation*}
\left[ u;\alpha ,\beta \right] :=\frac{u\left( \alpha \right) -u\left( \beta
\right) }{\alpha -\beta }.
\end{equation*}

The following representation of $D\left( f,u\right) $ may be stated, see 
\cite{III.c.D} and \cite{III.c.D1}. Due to its importance in proving our new
results we present here a short proof as well.

\begin{lemma}
\label{III.c.l.3.1}Let $f,u:\left[ a,b\right] \rightarrow \mathbb{C}$ be
such that the Stieltjes integral $\int_{a}^{b}f\left( t\right) du\left(
t\right) $ and the Riemann integral $\int_{a}^{b}f\left( t\right) dt$ exist.
Then%
\begin{align}
D\left( f,u\right) & =\int_{a}^{b}\Phi \left( t\right) df\left( t\right) =%
\frac{1}{b-a}\int_{a}^{b}\Gamma \left( t\right) df\left( t\right)
\label{III.c.e.3.5} \\
& =\frac{1}{b-a}\int_{a}^{b}\left( t-a\right) \left( b-t\right) \Delta
\left( t\right) df\left( t\right) .  \notag
\end{align}
\end{lemma}

\begin{proof}
Since $\int_{a}^{b}f\left( t\right) du\left( t\right) $ exists, hence $%
\int_{a}^{b}\Phi \left( t\right) df\left( t\right) $ also exists, and the
integration by parts formula for Riemann-Stieltjes integrals gives that%
\begin{align*}
\int_{a}^{b}\Phi \left( t\right) df\left( t\right) & =\int_{a}^{b}\left[ 
\frac{\left( t-a\right) u\left( b\right) +\left( b-t\right) u\left( a\right) 
}{b-a}-u\left( t\right) \right] df\left( t\right) \\
& =\left. \left[ \frac{\left( t-a\right) u\left( b\right) +\left( b-t\right)
u\left( a\right) }{b-a}-u\left( t\right) \right] f\left( t\right)
\right\vert _{a}^{b} \\
& -\int_{a}^{b}f\left( t\right) d\left[ \frac{\left( t-a\right) u\left(
b\right) +\left( b-t\right) u\left( a\right) }{b-a}-u\left( t\right) \right]
\\
& =-\int_{a}^{b}f\left( t\right) \left[ \frac{u\left( b\right) -u\left(
a\right) }{b-a}dt-du\left( t\right) \right] =D\left( f,u\right) ,
\end{align*}%
proving the required identity.
\end{proof}

For recent inequalities related to $D\left( f;u\right) $ for various pairs
of functions $\left( f,u\right) ,$ see \cite{III.c.D2}.

The following representation for a continuous function of selfadjoint
operator may be stated:

\begin{lemma}[Dragomir, 2010, \protect\cite{III.c.SSD5}]
\label{III.c.l.3.2}Let $A$ be a selfadjoint operator in the Hilbert space $H$
with the spectrum $Sp\left( A\right) \subseteq \left[ m,M\right] $ for some
real numbers $m<M,$ $\left\{ E_{\lambda }\right\} _{\lambda }$ be its 
\textit{spectral family and} $f:\left[ m,M\right] \rightarrow \mathbb{C}$ a
continuous function on $\left[ m,M\right] .$ If $x,y\in H,$ then we have the
representation%
\begin{align}
\left\langle f\left( A\right) x,y\right\rangle & =\left\langle
x,y\right\rangle \frac{1}{M-m}\int_{m}^{M}f\left( s\right) ds
\label{III.c.e.3.6} \\
& +\frac{1}{M-m}\int_{m-0}^{M}\left\langle \left[ \left( t-m\right) \left(
1_{H}-E_{t}\right) -\left( M-t\right) E_{t}\right] x,y\right\rangle df\left(
t\right) .  \notag
\end{align}
\end{lemma}

\begin{proof}
Utilising Lemma \ref{III.c.l.3.1} we have 
\begin{align}
\int_{m}^{M}f\left( t\right) du\left( t\right) & =\left[ u\left( M\right)
-u\left( m\right) \right] \cdot \frac{1}{M-m}\int_{m}^{M}f\left( s\right) ds
\label{III.c.e.3.7} \\
& +\int_{m}^{M}\left[ \frac{\left( t-m\right) u\left( M\right) +\left(
M-t\right) u\left( m\right) }{M-m}-u\left( t\right) \right] df\left(
t\right) ,  \notag
\end{align}%
for any continuous function $f:\left[ m,M\right] \rightarrow \mathbb{C}$ and
any function of bounded variation $u:\left[ m,M\right] \rightarrow \mathbb{C}
$.

Now, if we write the equality (\ref{III.c.e.3.7}) for $u\left( t\right)
=\left\langle E_{t}x,y\right\rangle $ with $x,y\in H,$ then we get%
\begin{align}
\int_{m-0}^{M}f\left( t\right) d\left\langle E_{t}x,y\right\rangle &
=\left\langle x,y\right\rangle \cdot \frac{1}{M-m}\int_{m}^{M}f\left(
s\right) ds  \label{III.c.e.3.8} \\
& +\int_{m-0}^{M}\left[ \frac{\left( t-m\right) \left\langle
x,y\right\rangle }{M-m}-\left\langle E_{t}x,y\right\rangle \right] df\left(
t\right) ,  \notag
\end{align}%
which, by the spectral representation theorem, produces the desired result (%
\ref{III.c.e.3.6}).
\end{proof}

The following result may be stated:

\begin{theorem}[Dragomir, 2010, \protect\cite{III.c.SSD5}]
\label{III.c.t.3.1}Let $A$ be a selfadjoint operator in the Hilbert space $H$
with the spectrum $Sp\left( A\right) \subseteq \left[ m,M\right] $ for some
real numbers $m<M$ $\left\{ E_{\lambda }\right\} _{\lambda }$ be its \textit{%
spectral family and} $f:\left[ m,M\right] \rightarrow \mathbb{C}$ a
continuous function on $\left[ m,M\right] .$

1. If $f$ is of bounded variation, then 
\begin{align}
& \left\vert \left\langle f\left( A\right) x,y\right\rangle -\left\langle
x,y\right\rangle \frac{1}{M-m}\int_{m}^{M}f\left( s\right) ds\right\vert
\label{III.c.e.3.9} \\
& \leq \left\Vert y\right\Vert \dbigvee\limits_{m}^{M}\left( f\right)  \notag
\\
& \times \max_{t\in \left[ m,M\right] }\left[ \left( \frac{t-m}{M-m}\right)
^{2}\left\Vert \left( 1_{H}-E_{t}\right) x\right\Vert ^{2}+\left( \frac{M-t}{%
M-m}\right) ^{2}\left\Vert E_{t}x\right\Vert ^{2}\right] ^{1/2}  \notag \\
& \leq \left\Vert x\right\Vert \left\Vert y\right\Vert
\dbigvee\limits_{m}^{M}\left( f\right)  \notag
\end{align}%
for any $x,y\in H.$

2. If $f$ is Lipschitzian with the constant $L>0$, then 
\begin{align}
& \left\vert \left\langle f\left( A\right) x,y\right\rangle -\left\langle
x,y\right\rangle \frac{1}{M-m}\int_{m}^{M}f\left( s\right) ds\right\vert
\label{III.c.e.3.10} \\
& \leq \frac{L\left\Vert y\right\Vert }{M-m}\int_{m-0}^{M}\left[ \left(
t-m\right) ^{2}\left\Vert \left( 1_{H}-E_{t}\right) x\right\Vert ^{2}+\left(
M-t\right) ^{2}\left\Vert E_{t}x\right\Vert ^{2}\right] ^{1/2}dt  \notag \\
& \leq \frac{1}{2}\left[ 1+\frac{\sqrt{2}}{2}\ln \left( \sqrt{2}+1\right) %
\right] \left( M-m\right) L\left\Vert y\right\Vert \left\Vert x\right\Vert 
\notag
\end{align}%
for any $x,y\in H.$

3. If $f:\left[ m,M\right] \rightarrow \mathbb{R}$ is monotonic
nondecreasing, then 
\begin{align}
& \left\vert \left\langle f\left( A\right) x,y\right\rangle -\left\langle
x,y\right\rangle \frac{1}{M-m}\int_{m}^{M}f\left( s\right) ds\right\vert
\label{III.c.e.3.11} \\
& \leq \frac{\left\Vert y\right\Vert }{M-m}\int_{m-0}^{M}\left[ \left(
t-m\right) ^{2}\left\Vert \left( 1_{H}-E_{t}\right) x\right\Vert ^{2}+\left(
M-t\right) ^{2}\left\Vert E_{t}x\right\Vert ^{2}\right] ^{1/2}df\left(
t\right)  \notag \\
& \leq \left\Vert y\right\Vert \left\Vert x\right\Vert \int_{m}^{M}\left[
\left( \frac{t-m}{M-m}\right) ^{2}+\left( \frac{M-t}{M-m}\right) ^{2}\right]
^{1/2}df\left( t\right)  \notag \\
& \leq \left\Vert y\right\Vert \left\Vert x\right\Vert \left[ f\left(
M\right) -f\left( m\right) \right] ^{1/2}  \notag \\
& \times \left[ f\left( M\right) -f\left( m\right) -\frac{4}{M-m}%
\int_{m}^{M}\left( t-\frac{m+M}{2}\right) f\left( t\right) dt\right] ^{1/2} 
\notag \\
& \leq \left\Vert y\right\Vert \left\Vert x\right\Vert \left[ f\left(
M\right) -f\left( m\right) \right]  \notag
\end{align}%
for any $x,y\in H.$
\end{theorem}

\begin{proof}
If we assume that $f$ is of bounded variation, then on applying the property
(\ref{III.c.e.2.2.1}) to the representation (\ref{III.c.e.3.6}) we get%
\begin{align}
& \left\vert \left\langle f\left( A\right) x,y\right\rangle -\left\langle
x,y\right\rangle \frac{1}{M-m}\int_{m}^{M}f\left( s\right) ds\right\vert
\label{III.c.e.3.12} \\
& \leq \frac{1}{M-m}\max_{t\in \left[ m,M\right] }\left\vert \left\langle %
\left[ \left( t-m\right) \left( 1_{H}-E_{t}\right) -\left( M-t\right) E_{t}%
\right] x,y\right\rangle \right\vert \dbigvee\limits_{m}^{M}\left( f\right) .
\notag
\end{align}

Now, on utilizing the Schwarz inequality and the fact that $E_{t}$ is a
projector for any $t\in \left[ m,M\right] ,$ then we have%
\begin{align}
& \left\vert \left\langle \left[ \left( t-m\right) \left( 1_{H}-E_{t}\right)
-\left( M-t\right) E_{t}\right] x,y\right\rangle \right\vert
\label{III.c.e.3.13} \\
& \leq \left\Vert \left[ \left( t-m\right) \left( 1_{H}-E_{t}\right) -\left(
M-t\right) E_{t}\right] x\right\Vert \left\Vert y\right\Vert  \notag \\
& =\left[ \left( t-m\right) ^{2}\left\Vert \left( 1_{H}-E_{t}\right)
x\right\Vert ^{2}+\left( M-t\right) ^{2}\left\Vert E_{t}x\right\Vert ^{2}%
\right] ^{1/2}\left\Vert y\right\Vert  \notag \\
& \leq \left[ \left( t-m\right) ^{2}+\left( M-t\right) ^{2}\right]
^{1/2}\left\Vert x\right\Vert \left\Vert y\right\Vert  \notag
\end{align}%
for any $x,y\in H$ and for any $t\in \left[ m,M\right] .$

Taking the maximum in (\ref{III.c.e.3.13}) we deduce the desired inequality (%
\ref{III.c.e.3.9}).

It is well known that if $p:\left[ a,b\right] \rightarrow \mathbb{C}$ is a
Riemann integrable function and $v:\left[ a,b\right] \rightarrow \mathbb{C}$
is Lipschitzian with the constant $L>0$, i.e.,%
\begin{equation*}
\left\vert f\left( s\right) -f\left( t\right) \right\vert \leq L\left\vert
s-t\right\vert \text{ for any }t,s\in \left[ a,b\right] ,
\end{equation*}%
then the Riemann-Stieltjes integral $\int_{a}^{b}p\left( t\right) dv\left(
t\right) $ exists and the following inequality holds%
\begin{equation*}
\left\vert \int_{a}^{b}p\left( t\right) dv\left( t\right) \right\vert \leq
L\int_{a}^{b}\left\vert p\left( t\right) \right\vert dt.
\end{equation*}

Now, on applying this property of the Riemann-Stieltjes integral to the
representation (\ref{III.c.e.3.6}), we get%
\begin{align}
& \left\vert \left\langle f\left( A\right) x,y\right\rangle -\left\langle
x,y\right\rangle \frac{1}{M-m}\int_{m}^{M}f\left( s\right) ds\right\vert
\label{III.c.e.3.14} \\
& \leq \frac{L}{M-m}\int_{m-0}^{M}\left\vert \left\langle \left[ \left(
t-m\right) \left( 1_{H}-E_{t}\right) -\left( M-t\right) E_{t}\right]
x,y\right\rangle \right\vert dt  \notag \\
& \leq \frac{L\left\Vert y\right\Vert }{M-m}\int_{m-0}^{M}\left[ \left(
t-m\right) ^{2}\left\Vert \left( 1_{H}-E_{t}\right) x\right\Vert ^{2}+\left(
M-t\right) ^{2}\left\Vert E_{t}x\right\Vert ^{2}\right] ^{1/2}dt  \notag \\
& \leq L\left\Vert y\right\Vert \left\Vert x\right\Vert \int_{m}^{M}\left[
\left( \frac{t-m}{M-m}\right) ^{2}+\left( \frac{M-t}{M-m}\right) ^{2}\right]
^{1/2}dt,  \notag
\end{align}%
for any $x,y\in H.$

Now, if we change the variable in the integral by choosing $u=\frac{t-m}{M-m}
$ then we get%
\begin{align*}
& \int_{m}^{M}\left[ \left( \frac{t-m}{M-m}\right) ^{2}+\left( \frac{M-t}{M-m%
}\right) ^{2}\right] ^{1/2}dt \\
& =\left( M-m\right) \int_{0}^{1}\left[ u^{2}+\left( 1-u\right) ^{2}\right]
^{1/2}du \\
& =\frac{1}{2}\left( M-m\right) \left[ 1+\frac{\sqrt{2}}{2}\ln \left( \sqrt{2%
}+1\right) \right] ,
\end{align*}%
which together with (\ref{III.c.e.3.14}) produces the desired result (\ref%
{III.c.e.3.10}).

From the theory of Riemann-Stieltjes integral is well known that if $p:\left[
a,b\right] \rightarrow \mathbb{C}$ is of bounded variation and $v:\left[ a,b%
\right] \rightarrow \mathbb{R}$ is continuous and monotonic nondecreasing,
then the Riemann-Stieltjes integrals $\int_{a}^{b}p\left( t\right) dv\left(
t\right) $ and $\int_{a}^{b}\left\vert p\left( t\right) \right\vert dv\left(
t\right) $ exist and%
\begin{equation*}
\left\vert \int_{a}^{b}p\left( t\right) dv\left( t\right) \right\vert \leq
\int_{a}^{b}\left\vert p\left( t\right) \right\vert dv\left( t\right) .
\end{equation*}%
Now, on applying this property of the Riemann-Stieltjes integral, we have
from the representation (\ref{III.c.e.3.6})%
\begin{align}
& \left\vert \left\langle f\left( A\right) x,y\right\rangle -\left\langle
x,y\right\rangle \frac{1}{M-m}\int_{m}^{M}f\left( s\right) ds\right\vert
\label{III.c.e.3.15} \\
& \leq \frac{1}{M-m}\int_{m-0}^{M}\left\vert \left\langle \left[ \left(
t-m\right) \left( 1_{H}-E_{t}\right) -\left( M-t\right) E_{t}\right]
x,y\right\rangle \right\vert df\left( t\right)  \notag \\
& \leq \frac{\left\Vert y\right\Vert }{M-m}\int_{m-0}^{M}\left[ \left(
t-m\right) ^{2}\left\Vert \left( 1_{H}-E_{t}\right) x\right\Vert ^{2}+\left(
M-t\right) ^{2}\left\Vert E_{t}x\right\Vert ^{2}\right] ^{1/2}df\left(
t\right)  \notag \\
& \leq \left\Vert y\right\Vert \left\Vert x\right\Vert \int_{m}^{M}\left[
\left( \frac{t-m}{M-m}\right) ^{2}+\left( \frac{M-t}{M-m}\right) ^{2}\right]
^{1/2}df\left( t\right) ,  \notag
\end{align}%
for any $x,y\in H$ and the proof of the first and second inequality in (\ref%
{III.c.e.3.11}) is completed.

For the last part we use the following Cauchy-Buniakowski-Schwarz integral
inequality for the Riemann-Stieltjes integral with monotonic nondecreasing
integrator $v$%
\begin{equation*}
\left\vert \int_{a}^{b}p\left( t\right) q\left( t\right) dv\left( t\right)
\right\vert \leq \left[ \int_{a}^{b}\left\vert p\left( t\right) \right\vert
^{2}dv\left( t\right) \right] ^{1/2}\left[ \int_{a}^{b}\left\vert q\left(
t\right) \right\vert ^{2}dv\left( t\right) \right] ^{1/2}
\end{equation*}%
where $p,q:\left[ a,b\right] \rightarrow \mathbb{C}$ are continuous on $%
\left[ a,b\right] .$

By applying this inequality we conclude that%
\begin{align}
& \int_{m}^{M}\left[ \left( \frac{t-m}{M-m}\right) ^{2}+\left( \frac{M-t}{M-m%
}\right) ^{2}\right] ^{1/2}df\left( t\right)  \label{III.c.e.3.16} \\
& \leq \left[ \int_{m}^{M}df\left( t\right) \right] ^{1/2}\left[ \int_{m}^{M}%
\left[ \left( \frac{t-m}{M-m}\right) ^{2}+\left( \frac{M-t}{M-m}\right) ^{2}%
\right] df\left( t\right) \right] ^{1/2}.  \notag
\end{align}%
Further, integrating by parts in the Riemann-Stieltjes integral we also have
that%
\begin{align}
& \int_{m}^{M}\left[ \left( \frac{t-m}{M-m}\right) ^{2}+\left( \frac{M-t}{M-m%
}\right) ^{2}\right] df\left( t\right)  \label{III.c.e.3.17} \\
& =f\left( M\right) -f\left( m\right) -\frac{4}{M-m}\int_{m}^{M}\left( t-%
\frac{m+M}{2}\right) f\left( t\right) dt  \notag \\
& \leq f\left( M\right) -f\left( m\right)  \notag
\end{align}%
where for the last part we used the fact that by the \v{C}eby\v{s}ev
integral inequality for monotonic functions with the same monotonicity we
have that 
\begin{align*}
& \int_{m}^{M}\left( t-\frac{m+M}{2}\right) f\left( t\right) dt \\
& \geq \frac{1}{M-m}\int_{m}^{M}\left( t-\frac{m+M}{2}\right)
dt\int_{m}^{M}f\left( t\right) dt=0.
\end{align*}
\end{proof}

\subsection{Some Applications for Particular Functions}

\textbf{1.} Consider the function $f:\left( 0,\infty \right) \rightarrow 
\mathbb{R}$ given by $f\left( t\right) =t^{r}$ with $r\in (0,1].$ This
function is $r$-H\"{o}lder continuous with the constant $H>0.$ Then, by
applying Corollary \ref{III.c.c.2.1} we can state the following result

\begin{proposition}
\label{III.c.p.4.1}Let $A$ be a selfadjoint operator in the Hilbert space $H$
with the spectrum $Sp\left( A\right) \subseteq \left[ m,M\right] $ for some
real numbers $\,0<m<M$ and $\left\{ E_{\lambda }\right\} _{\lambda }$ be its 
\textit{spectral family. Then for all }$r$\textit{\ }with $r\in (0,1]$%
\textit{\ we have the inequality }%
\begin{align}
& \left\vert \left\langle A^{r}x,y\right\rangle -\left\langle
x,y\right\rangle \frac{M^{r+1}-m^{r+1}}{\left( r+1\right) \left( M-m\right) }%
\right\vert   \label{III.c.e.4.1} \\
& \leq \frac{1}{r+1}\left( M-m\right) ^{r}\dbigvee\limits_{m-0}^{M}\left(
\left\langle E_{\left( \cdot \right) }x,y\right\rangle \right) \leq \frac{1}{%
r+1}\left( M-m\right) ^{r}\left\Vert x\right\Vert \left\Vert y\right\Vert  
\notag
\end{align}%
for any $x,y\in H.$
\end{proposition}

The case of $p>1$ is incorporated in the following proposition:

\begin{proposition}
\label{III.c.p.4.2}With the same assumptions from Proposition \ref%
{III.c.p.4.1} and if $p>1,$ then we have%
\begin{align}
& \left\vert \left\langle A^{p}x,y\right\rangle -\frac{M^{p+1}-m^{p+1}}{%
\left( p+1\right) \left( M-m\right) }\left\langle x,y\right\rangle
\right\vert   \label{III.c.e.4.2} \\
& \leq \frac{1}{2}pM^{p-1}\left( M-m\right) \dbigvee\limits_{m-0}^{M}\left(
\left\langle E_{\left( \cdot \right) }x,y\right\rangle \right) \leq \frac{1}{%
2}pM^{p-1}\left( M-m\right) \left\Vert x\right\Vert \left\Vert y\right\Vert 
\notag
\end{align}%
for any $x,y\in H$.
\end{proposition}

The case of negative powers except $p=-1$ goes likewise and we omit the
details.

Now, if we apply Corollary \ref{III.c.c.2.1} for the function $f\left(
t\right) =-\frac{1}{t}$ with $t>0,$ then we can state the following
proposition:

\begin{proposition}
\label{III.c.p.4.3}Let $A$ be a selfadjoint operator in the Hilbert space $H$
with the spectrum $Sp\left( A\right) \subseteq \left[ m,M\right] $ for some
real numbers $0<m<M$ and let $\left\{ E_{\lambda }\right\} _{\lambda }$ be
its \textit{spectral family. Then }for any $x,y\in H$ we have the
inequalities%
\begin{align}
& \left\vert \left\langle A^{-1}x,y\right\rangle -\frac{\ln M-\ln m}{M-m}%
\left\langle x,y\right\rangle \right\vert   \label{III.c.e.4.4} \\
& \leq \frac{1}{2}\frac{M-m}{m^{2}}\dbigvee\limits_{m-0}^{M}\left(
\left\langle E_{\left( \cdot \right) }x,y\right\rangle \right) \leq \frac{1}{%
2}\frac{M-m}{m^{2}}\left\Vert x\right\Vert \left\Vert y\right\Vert .  \notag
\end{align}
\end{proposition}

\textbf{2. }Now, if we apply Corollary \ref{III.c.c.2.1} to the function $%
f:\left( 0,\infty \right) \rightarrow \mathbb{R}$, $f\left( t\right) =\ln t$%
, then we can state

\begin{proposition}
\label{III.c.p.4.4}Let $A$ be a selfadjoint operator in the Hilbert space $H$
with the spectrum $Sp\left( A\right) \subseteq \left[ m,M\right] $ for some
real numbers $0<m<M$ and let $\left\{ E_{\lambda }\right\} _{\lambda }$ be
its \textit{spectral family. Then }for any $x,y\in H$ we have the
inequalities%
\begin{align}
& \left\vert \left\langle \ln Ax,y\right\rangle -\left\langle
x,y\right\rangle \ln I\left( m,M\right) \right\vert   \label{III.c.e.4.6} \\
& \leq \frac{1}{2}\left( \frac{M}{m}-1\right)
\dbigvee\limits_{m-0}^{M}\left( \left\langle E_{\left( \cdot \right)
}x,y\right\rangle \right) \leq \frac{1}{2}\left( \frac{M}{m}-1\right)
\left\Vert x\right\Vert \left\Vert y\right\Vert ,  \notag
\end{align}%
where $I\left( m,M\right) $ is the identric mean of $m$ and $M$ and is
defined by%
\begin{equation*}
I\left( m,M\right) =\frac{1}{e}\left( \frac{M^{M}}{m^{m}}\right) ^{1/(M-m)}.
\end{equation*}
\end{proposition}

\section{Ostrowski's Type Inequalities for $n$-Time Differentiable Functions}

\subsection{Some Identities}

In \cite{III.e.CDJ}, the authors have pointed out the following integral
identity:

\begin{lemma}[ Cerone-Dragomir-Roumeliotis, 1999, \protect\cite{III.e.CDJ}]
\label{III.e.l2.1}Let $f:\left[ a,b\right] \rightarrow \mathbb{R}$ be a
mapping such that the $\left( n-1\right) $-derivative $f^{\left( n-1\right)
} $ (where $n\geq 1)$ is absolutely continuous on $\left[ a,b\right] $. Then
for all $x\in \left[ a,b\right] $, we have the identity: 
\begin{align}
\int_{a}^{b}f\left( t\right) dt& =\sum_{k=0}^{n-1}\left[ \frac{\left(
b-x\right) ^{k+1}+\left( -1\right) ^{k}\left( x-a\right) ^{k+1}}{\left(
k+1\right) !}\right] f^{\left( k\right) }\left( x\right)  \label{III.e.e.2.1}
\\
& +\left( -1\right) ^{n}\int_{a}^{b}K_{n}\left( x,t\right) f^{\left(
n\right) }\left( t\right) dt  \notag
\end{align}%
where the kernel $K_{n}:\left[ a,b\right] ^{2}\rightarrow \mathbb{R}$ is
given by 
\begin{equation}
K_{n}\left( x,t\right) :=\left\{ 
\begin{array}{ll}
\frac{\left( t-a\right) ^{n}}{n!}, & a\leq t\leq x\leq b \\ 
&  \\ 
\frac{\left( t-b\right) ^{n}}{n!}, & a\leq x<t\leq b.%
\end{array}%
\right.  \label{III.e.e.2.2}
\end{equation}
\end{lemma}

The identity (\ref{III.e.e.2.2}) can be written in the following equivalent
form as: 
\begin{align}
f\left( z\right) & =\frac{1}{b-a}\int_{a}^{b}f\left( t\right) dt
\label{III.e.e.2.3} \\
& -\frac{1}{b-a}\sum_{k=1}^{n-1}\frac{1}{\left( k+1\right) !}\left[ \left(
b-z\right) ^{k+1}+\left( -1\right) ^{k}\left( z-a\right) ^{k+1}\right]
f^{\left( k\right) }\left( z\right)  \notag \\
& +\frac{\left( -1\right) ^{n-1}}{\left( b-a\right) n!}\left[
\int_{a}^{z}\left( t-a\right) ^{n}f^{\left( n\right) }\left( t\right)
dt+\int_{z}^{b}\left( t-b\right) ^{n}f^{\left( n\right) }\left( t\right) dt%
\right]  \notag
\end{align}%
for all $z\in \left[ a,b\right] $.

Note that for $n=1$, the sum $\sum_{k=1}^{n-1}$ is empty and we obtain the
well known\textit{\ Montgomery's identity} (see for example \cite{III.e.CD}) 
\begin{align}
f\left( z\right) & =\frac{1}{b-a}\int_{a}^{b}f\left( t\right) dt
\label{e.2.4} \\
& +\frac{1}{b-a}\left[ \int_{a}^{z}\left( t-a\right) f^{\left( 1\right)
}\left( t\right) dt+\int_{z}^{b}\left( t-b\right) f^{\left( 1\right) }\left(
t\right) dt\right] ,\;  \notag
\end{align}%
for any $z\in \left[ a,b\right] .$

In a slightly more general setting, by the use of the identity (\ref%
{III.e.e.2.3}), we can state the following result as well:

\begin{lemma}[Dragomir, 2010, \protect\cite{III.e.SSD7}]
\label{III.l.2.2}Let $f:\left[ a,b\right] \rightarrow \mathbb{R}$ be a
mapping such that the $n$-derivative $f^{\left( n\right) }$ (where $n\geq 1)$
is of bounded variation on $\left[ a,b\right] $. Then for all $\lambda \in %
\left[ a,b\right] $, we have the identity:%
\begin{align}
f\left( \lambda \right) & =\frac{1}{b-a}\int_{a}^{b}f\left( t\right) dt
\label{III.e.e.2.5} \\
& -\frac{1}{b-a}\sum_{k=1}^{n}\frac{1}{\left( k+1\right) !}\left[ \left(
b-\lambda \right) ^{k+1}+\left( -1\right) ^{k}\left( \lambda -a\right) ^{k+1}%
\right] f^{\left( k\right) }\left( \lambda \right)  \notag \\
& +\frac{\left( -1\right) ^{n}}{\left( b-a\right) \left( n+1\right) !} 
\notag \\
& \times \left[ \int_{a}^{\lambda }\left( t-a\right) ^{n+1}d\left( f^{\left(
n\right) }\left( t\right) \right) +\int_{\lambda }^{b}\left( t-b\right)
^{n+1}d\left( f^{\left( n\right) }\left( t\right) \right) \right] .  \notag
\end{align}
\end{lemma}

Now we can state the following representation result for functions of
selfadjoint operators:

\begin{theorem}[Dragomir, 2010, \protect\cite{III.e.SSD7}]
\label{III.e.t.2.1}Let $A$ be a selfadjoint operator in the Hilbert space $H$
with the spectrum $Sp\left( A\right) \subseteq \left[ m,M\right] $ for some
real numbers $m<M$, $\left\{ E_{\lambda }\right\} _{\lambda }$ be its 
\textit{spectral family,} $I$ be a closed subinterval on $\mathbb{R}$ with $%
\left[ m,M\right] \subset \mathring{I}$ (the interior of $I)$ and let $n$ be
an integer with $n\geq 1.$ If $f:I\rightarrow \mathbb{C}$ is such that the $%
n $-th derivative $f^{\left( n\right) }$ is of bounded variation on the
interval $\left[ m,M\right] $, then we have the representation%
\begin{align}
f\left( A\right) & =\left( \frac{1}{M-m}\int_{m}^{M}f\left( t\right)
dt\right) 1_{H}-\frac{1}{M-m}  \label{III.e.e.2.6} \\
& \times \sum_{k=1}^{n}\frac{1}{\left( k+1\right) !}\left[ \left(
M1_{H}-A\right) ^{k+1}+\left( -1\right) ^{k}\left( A-m1_{H}\right) ^{k+1}%
\right] f^{\left( k\right) }\left( A\right)  \notag \\
& +T_{n}\left( A,m,M\right)  \notag
\end{align}%
where the remainder is given by%
\begin{align}
T_{n}\left( A,m,M\right) & :=\frac{\left( -1\right) ^{n}}{\left( M-m\right)
\left( n+1\right) !}  \label{III.e.e.2.7} \\
& \times \left[ \int_{m-0}^{M}\left( \int_{m}^{\lambda }\left( t-m\right)
^{n+1}d\left( f^{\left( n\right) }\left( t\right) \right) \right)
dE_{\lambda }\right.  \notag \\
& \left. +\int_{m-0}^{M}\left( \int_{\lambda }^{M}\left( t-M\right)
^{n+1}d\left( f^{\left( n\right) }\left( t\right) \right) \right)
dE_{\lambda }\right] .  \notag
\end{align}

In particular, if the $n$-th derivative $f^{\left( n\right) }$ is absolutely
continuous on $\left[ m,M\right] $, then the remainder can be represented as 
\begin{align}
& T_{n}\left( A,m,M\right)  \label{III.e.e.2.8} \\
& =\frac{\left( -1\right) ^{n}}{\left( M-m\right) \left( n+1\right) !} 
\notag \\
& \times \int_{m-0}^{M}\left[ \left( \lambda -m\right) ^{n+1}\left(
1_{H}-E_{\lambda }\right) +\left( \lambda -M\right) ^{n+1}E_{\lambda }\right]
f^{\left( n+1\right) }\left( \lambda \right) d\lambda .  \notag
\end{align}
\end{theorem}

\begin{proof}
By Lemma \ref{III.l.2.2} we have%
\begin{align}
f\left( \lambda \right) & =\frac{1}{M-m}\int_{m}^{M}f\left( t\right) dt-%
\frac{1}{M-m}  \label{III.e.e.2.9} \\
& \times \sum_{k=1}^{n}\frac{1}{\left( k+1\right) !}\left[ \left( M-\lambda
\right) ^{k+1}+\left( -1\right) ^{k}\left( \lambda -m\right) ^{k+1}\right]
f^{\left( k\right) }\left( \lambda \right)  \notag \\
& +\frac{\left( -1\right) ^{n}}{\left( M-m\right) \left( n+1\right) !} 
\notag \\
& \times \left[ \int_{m}^{\lambda }\left( t-m\right) ^{n+1}d\left( f^{\left(
n\right) }\left( t\right) \right) +\int_{\lambda }^{M}\left( t-M\right)
^{n+1}d\left( f^{\left( n\right) }\left( t\right) \right) \right]  \notag
\end{align}%
for any $\lambda \in \left[ m,M\right] .$

Integrating the identity (\ref{III.e.e.2.9}) in the Riemann-Stieltjes sense
with the integrator $E_{\lambda }$ we get 
\begin{align}
& \int_{m}^{M}f\left( \lambda \right) dE_{\lambda }  \label{III.e.e.2.10} \\
& =\frac{1}{M-m}\int_{m}^{M}f\left( t\right) dt\int_{m}^{M}dE_{\lambda }-%
\frac{1}{M-m}  \notag \\
& \times \sum_{k=1}^{n}\frac{1}{\left( k+1\right) !}\int_{m}^{M}\left[
\left( M-\lambda \right) ^{k+1}+\left( -1\right) ^{k}\left( \lambda
-m\right) ^{k+1}\right] f^{\left( k\right) }\left( \lambda \right)
dE_{\lambda }  \notag \\
& +T_{n}\left( A,m,M\right) .  \notag
\end{align}

Since, by the spectral representation theorem we have%
\begin{equation*}
\int_{m-0}^{M}f\left( \lambda \right) dE_{\lambda }=f\left( A\right)
,\int_{m-0}^{M}dE_{\lambda }=1_{H}
\end{equation*}%
and 
\begin{align*}
& \int_{m-0}^{M}\left[ \left( M-\lambda \right) ^{k+1}+\left( -1\right)
^{k}\left( \lambda -m\right) ^{k+1}\right] f^{\left( k\right) }\left(
\lambda \right) dE_{\lambda } \\
& =\left[ \left( M1_{H}-A\right) ^{k+1}+\left( -1\right) ^{k}\left(
A-m1_{H}\right) ^{k+1}\right] f^{\left( k\right) }\left( A\right) ,
\end{align*}%
then by (\ref{III.e.e.2.10}) we deduce the representation (\ref{III.e.e.2.6}%
).

Now, if the $n$-th derivative $f^{\left( n\right) }$ is absolutely
continuous on $\left[ m,M\right] ,$ then 
\begin{equation*}
\int_{m}^{\lambda }\left( t-m\right) ^{n+1}d\left( f^{\left( n\right)
}\left( t\right) \right) =\int_{m}^{\lambda }\left( t-m\right)
^{n+1}f^{\left( n+1\right) }\left( t\right) dt
\end{equation*}%
and 
\begin{equation*}
\int_{\lambda }^{M}\left( t-M\right) ^{n+1}d\left( f^{\left( n\right)
}\left( t\right) \right) =\int_{\lambda }^{M}\left( t-M\right)
^{n+1}f^{\left( n+1\right) }\left( t\right) dt
\end{equation*}%
where the integrals in the right hand side are taken in the Lebesgue sense.

Utilising the integration by parts formula for the Riemann-Stieltjes
integral and the differentiation rule for the Stieltjes integral we have
successively 
\begin{eqnarray*}
&&\int_{m-0}^{M}\left( \int_{m}^{\lambda }\left( t-m\right) ^{n+1}f^{\left(
n+1\right) }\left( t\right) dt\right) dE_{\lambda } \\
&=&\left. \left( \int_{m}^{\lambda }\left( t-m\right) ^{n+1}f^{\left(
n+1\right) }\left( t\right) dt\right) E_{\lambda }\right\vert
_{m-0}^{M}-\int_{m-0}^{M}\left( \lambda -m\right) ^{n+1}f^{\left( n+1\right)
}\left( \lambda \right) E_{\lambda }d\lambda \\
&=&\left( \int_{m}^{M}\left( t-m\right) ^{n+1}f^{\left( n+1\right) }\left(
t\right) dt\right) 1_{H}-\int_{m-0}^{M}\left( \lambda -m\right)
^{n+1}f^{\left( n+1\right) }\left( \lambda \right) E_{\lambda }d\lambda \\
&=&\int_{m-0}^{M}\left( \lambda -m\right) ^{n+1}f^{\left( n+1\right) }\left(
\lambda \right) \left( 1_{H}-E_{\lambda }\right) d\lambda
\end{eqnarray*}%
and%
\begin{eqnarray*}
&&\int_{m-0}^{M}\left( \int_{\lambda }^{M}\left( t-M\right) ^{n+1}f^{\left(
n+1\right) }\left( t\right) dt\right) dE_{\lambda } \\
&=&\left. \left( \int_{\lambda }^{M}\left( t-M\right) ^{n+1}f^{\left(
n+1\right) }\left( t\right) dt\right) E_{\lambda }\right\vert
_{m-0}^{M}+\int_{m-0}^{M}\left( \lambda -M\right) ^{n+1}f^{\left( n+1\right)
}\left( \lambda \right) E_{\lambda }d\lambda \\
&=&\int_{m-0}^{M}\left( \lambda -M\right) ^{n+1}f^{\left( n+1\right) }\left(
\lambda \right) E_{\lambda }d\lambda
\end{eqnarray*}%
and the representation (\ref{III.e.e.2.8}) is thus obtained.
\end{proof}

\begin{remark}
\label{III.e.r.2.1}Let $A$ be a positive selfadjoint operator in the Hilbert
space $H$ with the spectrum $Sp\left( A\right) \subseteq \left[ m,M\right] $
for some positive real numbers $0<m<M$ and $\left\{ E_{\lambda }\right\}
_{\lambda }$ be its \textit{spectral family. Then, for }$n\geq 1,$\textit{\
we have the equality}%
\begin{align}
\ln A& =\left[ \ln I\left( m,M\right) \right] 1_{H}+\frac{1}{M-m}
\label{III.e.e.2.11} \\
& \times \sum_{k=1}^{n}\frac{1}{k\left( k+1\right) }\left[ \left(
A-m1_{H}\right) ^{k+1}+\left( -1\right) ^{k}\left( M1_{H}-A\right) ^{k+1}%
\right] A^{-k}  \notag \\
& +\frac{1}{\left( M-m\right) \left( n+1\right) }  \notag \\
& \times \left[ \int_{m-0}^{M}\left[ \left( \lambda -m\right) ^{n+1}\left(
1_{H}-E_{\lambda }\right) +\left( \lambda -M\right) ^{n+1}E_{\lambda }\right]
\lambda ^{-n-1}d\lambda \right] ,  \notag
\end{align}%
where $I\left( m,M\right) $ is the identric mean and is defined by%
\begin{equation*}
I\left( m,M\right) =\left\{ 
\begin{array}{c}
\frac{1}{e}\left( \frac{M^{M}}{m^{m}}\right) ^{1/\left( M-m\right) }\text{
if }M\neq m; \\ 
\\ 
M\text{ if }M=m.%
\end{array}%
\right.
\end{equation*}
\end{remark}

\begin{remark}
\label{III.e.r.2.1.a}If we introduce the exponential mean by%
\begin{equation*}
E\left( m,M\right) =\left\{ 
\begin{array}{c}
\frac{\exp M-\exp m}{M-m}\text{ if }M\neq m; \\ 
\\ 
M\text{ if }M=m%
\end{array}%
\right.
\end{equation*}%
and applying the identity (\ref{III.e.e.2.6}) for the exponential function,
we have 
\begin{align}
& \left[ 1_{H}+\frac{1}{M-m}\sum_{k=1}^{n}\frac{1}{\left( k+1\right) !}\left[
\left( M1_{H}-A\right) ^{k+1}+\left( -1\right) ^{k}\left( A-m1_{H}\right)
^{k+1}\right] \right]  \label{III.e.e.2.12} \\
& \times \exp A-E\left( m,M\right) 1_{H}  \notag \\
& =\frac{\left( -1\right) ^{n}}{\left( M-m\right) \left( n+1\right) !}%
\int_{m-0}^{M}\left[ \left( \lambda -m\right) ^{n+1}\left( 1_{H}-E_{\lambda
}\right) +\left( \lambda -M\right) ^{n+1}E_{\lambda }\right] e^{\lambda
}d\lambda  \notag
\end{align}%
where $A$ is a selfadjoint operator in the Hilbert space $H$ with the
spectrum $Sp\left( A\right) \subseteq \left[ m,M\right] $ for some real
numbers $m<M$ and $\left\{ E_{\lambda }\right\} _{\lambda }$ is its \textit{%
spectral family.}
\end{remark}

\subsection{Error Bounds for $f^{\left( n\right) }$ of Bounded Variation}

From the identity (\ref{III.e.e.2.6}), we define for any $x,y\in H$ 
\begin{align}
& T_{n}\left( A,m,M;x,y\right)  \label{III.e.e.3.1} \\
& :=\left\langle f\left( A\right) x,y\right\rangle +\frac{1}{M-m}%
\sum_{k=1}^{n}\frac{1}{\left( k+1\right) !}  \notag \\
& \times \left[ \left\langle \left( M1_{H}-A\right) ^{k+1}f^{\left( k\right)
}\left( A\right) x,y\right\rangle +\left( -1\right) ^{k}\left\langle \left(
A-m1_{H}\right) ^{k+1}f^{\left( k\right) }\left( A\right) x,y\right\rangle %
\right]  \notag \\
& -\left( \frac{1}{M-m}\int_{m}^{M}f\left( t\right) dt\right) \left\langle
x,y\right\rangle .  \notag
\end{align}

We have the following result concerning bounds for the absolute value of $%
T_{n}\left( A,m,M;x,y\right) $ when the $n$-th derivative $f^{\left(
n\right) }$ is of bounded variation:

\begin{theorem}[Dragomir, 2010, \protect\cite{III.e.SSD7}]
\label{III.e.t.3.1}Let $A$ be a selfadjoint operator in the Hilbert space $H$
with the spectrum $Sp\left( A\right) \subseteq \left[ m,M\right] $ for some
real numbers $m<M$, $\left\{ E_{\lambda }\right\} _{\lambda }$ be its 
\textit{spectral family,} $I$ be a closed subinterval on $\mathbb{R}$ with $%
\left[ m,M\right] \subset \mathring{I}$ and let $n$ be an integer with $%
n\geq 1.$

1. If $f:I\rightarrow \mathbb{C}$ is such that the $n$-th derivative $%
f^{\left( n\right) }$ is of bounded variation on the interval $\left[ m,M%
\right] $, then we have the inequalities%
\begin{align}
& \left\vert T_{n}\left( A,m,M;x,y\right) \right\vert   \label{III.e.e.3.2}
\\
& \leq \frac{1}{\left( M-m\right) \left( n+1\right) !}\dbigvee%
\limits_{m-0}^{M}\left( \left\langle E_{\left( \cdot \right)
}x,y\right\rangle \right)   \notag \\
& \times \max_{\lambda \in \left[ m,M\right] }\left[ \left( \lambda
-m\right) ^{n+1}\dbigvee\limits_{m}^{\lambda }\left( f^{\left( n\right)
}\right) +\left( M-\lambda \right) ^{n+1}\dbigvee\limits_{\lambda
}^{M}\left( f^{\left( n\right) }\right) \right]   \notag \\
& \leq \frac{\left( M-m\right) ^{n}}{\left( n+1\right) !}\dbigvee%
\limits_{m-0}^{M}\left( \left\langle E_{\left( \cdot \right)
}x,y\right\rangle \right) \dbigvee\limits_{m}^{M}\left( f^{\left( n\right)
}\right) \leq \frac{\left( M-m\right) ^{n}}{\left( n+1\right) !}%
\dbigvee\limits_{m}^{M}\left( f^{\left( n\right) }\right) \left\Vert
x\right\Vert \left\Vert y\right\Vert   \notag
\end{align}%
for any $x,y\in H.$

2. If $f:I\rightarrow \mathbb{C}$ is such that the $n$-th derivative $%
f^{\left( n\right) }$ is Lipschitzian with the constant $L_{n}>0$ on the
interval $\left[ m,M\right] $, then we have the inequalities%
\begin{align}
\left\vert T_{n}\left( A,m,M;x,y\right) \right\vert & \leq \frac{L_{n}\left(
M-m\right) ^{n+1}}{\left( n+2\right) !}\dbigvee\limits_{m-0}^{M}\left(
\left\langle E_{\left( \cdot \right) }x,y\right\rangle \right) 
\label{III.e.e.3.3} \\
& \leq \frac{L_{n}\left( M-m\right) ^{n+1}}{\left( n+2\right) !}\left\Vert
x\right\Vert \left\Vert y\right\Vert   \notag
\end{align}%
for any $x,y\in H.$

3. If $f:I\rightarrow \mathbb{R}$ is such that the $n$-th derivative $%
f^{\left( n\right) }$ is monotonic nondecreasing on the interval $\left[ m,M%
\right] $, then we have the inequalities%
\begin{align}
& \left\vert T_{n}\left( A,m,M;x,y\right) \right\vert   \label{III.e.e.3.3.a}
\\
& \leq \frac{1}{\left( M-m\right) \left( n+1\right) !}\dbigvee%
\limits_{m-0}^{M}\left( \left\langle E_{\left( \cdot \right)
}x,y\right\rangle \right)   \notag \\
& \times \max_{\lambda \in \left[ m,M\right] }\left[ f^{\left( n\right)
}\left( \lambda \right) \left( \left( \lambda -m\right) ^{n+1}-\left(
M-\lambda \right) ^{n+1}\right) \right.   \notag \\
& \left. +\left( n+1\right) \left[ \int_{\lambda }^{M}\left( M-t\right)
^{n}f^{\left( n\right) }\left( t\right) dt-\int_{m}^{\lambda }\left(
t-m\right) ^{n}f^{\left( n\right) }\left( t\right) dt\right] \right]   \notag
\\
& \leq \frac{1}{\left( M-m\right) \left( n+1\right) !}\max_{\lambda \in %
\left[ m,M\right] }\left[ \left( \lambda -m\right) ^{n+1}\left[ f^{\left(
n\right) }\left( \lambda \right) -f^{\left( n\right) }\left( m\right) \right]
\right.   \notag \\
& \left. +\left( M-\lambda \right) ^{n+1}\left[ f^{\left( n\right) }\left(
M\right) -f^{\left( n\right) }\left( \lambda \right) \right] \right]
\dbigvee\limits_{m-0}^{M}\left( \left\langle E_{\left( \cdot \right)
}x,y\right\rangle \right)   \notag \\
& \leq \frac{\left( M-m\right) ^{n}}{\left( n+1\right) !}\dbigvee%
\limits_{m-0}^{M}\left( \left\langle E_{\left( \cdot \right)
}x,y\right\rangle \right) \left[ f^{\left( n\right) }\left( M\right)
-f^{\left( n\right) }\left( m\right) \right]   \notag \\
& \leq \frac{\left( M-m\right) ^{n}}{\left( n+1\right) !}\left[ f^{\left(
n\right) }\left( M\right) -f^{\left( n\right) }\left( m\right) \right]
\left\Vert x\right\Vert \left\Vert y\right\Vert   \notag
\end{align}%
for any $x,y\in H.$
\end{theorem}

\begin{proof}
1. By the identity (\ref{III.e.e.2.7}) we have for any $x,y\in H$ that%
\begin{align}
T_{n}\left( A,m,M;x,y\right) & :=\frac{\left( -1\right) ^{n}}{\left(
M-m\right) \left( n+1\right) !}  \label{III.e.e.3.4} \\
& \times \left[ \int_{m-0}^{M}\left( \int_{m}^{\lambda }\left( t-m\right)
^{n+1}d\left( f^{\left( n\right) }\left( t\right) \right) \right)
d\left\langle E_{\lambda }x,y\right\rangle \right.  \notag \\
& \left. +\int_{m-0}^{M}\left( \int_{\lambda }^{M}\left( t-M\right)
^{n+1}d\left( f^{\left( n\right) }\left( t\right) \right) \right)
d\left\langle E_{\lambda }x,y\right\rangle \right] .  \notag
\end{align}

It is well known that if $p:\left[ a,b\right] \rightarrow \mathbb{C}$ is a
continuous function, $v:\left[ a,b\right] \rightarrow \mathbb{C}$ is of
bounded variation then the Riemann-Stieltjes integral $\int_{a}^{b}p\left(
t\right) dv\left( t\right) $ exists and the following inequality holds%
\begin{equation}
\left\vert \int_{a}^{b}p\left( t\right) dv\left( t\right) \right\vert \leq
\max_{t\in \left[ a,b\right] }\left\vert p\left( t\right) \right\vert
\dbigvee\limits_{a}^{b}\left( v\right) ,  \label{III.e.e.3.5}
\end{equation}%
where $\dbigvee\limits_{a}^{b}\left( v\right) $ denotes the total variation
of $v$ on $\left[ a,b\right] .$

Taking the modulus in (\ref{III.e.e.3.4}) and utilizing the property (\ref%
{III.e.e.3.5}), we have successively that%
\begin{multline}
\left\vert T_{n}\left( A,m,M;x,y\right) \right\vert =\frac{1}{\left(
M-m\right) \left( n+1\right) !}  \label{III.e.e.3.6} \\
\times \left\vert \int_{m-0}^{M}\left[ \left( \int_{m}^{\lambda }\left(
t-m\right) ^{n+1}d\left( f^{\left( n\right) }\left( t\right) \right) +\left(
\int_{\lambda }^{M}\left( t-M\right) ^{n+1}d\left( f^{\left( n\right)
}\left( t\right) \right) \right) \right) \right] d\left\langle E_{\lambda
}x,y\right\rangle \right\vert  \\
\leq \frac{1}{\left( M-m\right) \left( n+1\right) !}\dbigvee%
\limits_{m-0}^{M}\left( \left\langle E_{\left( \cdot \right)
}x,y\right\rangle \right)  \\
\times \max_{\lambda \in \left[ m,M\right] }\left\vert \int_{m}^{\lambda
}\left( t-m\right) ^{n+1}d\left( f^{\left( n\right) }\left( t\right) \right)
+\int_{\lambda }^{M}\left( t-M\right) ^{n+1}d\left( f^{\left( n\right)
}\left( t\right) \right) \right\vert 
\end{multline}%
for any $x,y\in H.$

By the same property (\ref{III.e.e.3.5}) we have for $\lambda \in \left(
m,M\right) $ that%
\begin{align*}
\left\vert \int_{m}^{\lambda }\left( t-m\right) ^{n+1}d\left( f^{\left(
n\right) }\left( t\right) \right) \right\vert & \leq \max_{t\in \left[
m,\lambda \right] }\left( t-m\right) ^{n+1}\dbigvee\limits_{m}^{\lambda
}\left( f^{\left( n\right) }\right) \\
& =\left( \lambda -m\right) ^{n+1}\dbigvee\limits_{m}^{\lambda }\left(
f^{\left( n\right) }\right)
\end{align*}%
and%
\begin{align*}
\left\vert \int_{\lambda }^{M}\left( t-M\right) ^{n+1}d\left( f^{\left(
n\right) }\left( t\right) \right) \right\vert & \leq \max_{t\in \left[
\lambda ,M\right] }\left( M-t\right) ^{n+1}\dbigvee\limits_{\lambda
}^{M}\left( f^{\left( n\right) }\right) \\
& =\left( M-\lambda \right) ^{n+1}\dbigvee\limits_{\lambda }^{M}\left(
f^{\left( n\right) }\right)
\end{align*}%
which produce the inequality%
\begin{align}
& \left\vert \int_{m}^{\lambda }\left( t-m\right) ^{n+1}d\left( f^{\left(
n\right) }\left( t\right) \right) +\int_{\lambda }^{M}\left( t-M\right)
^{n+1}d\left( f^{\left( n\right) }\left( t\right) \right) \right\vert
\label{III.e.e.3.7} \\
& \leq \left( \lambda -m\right) ^{n+1}\dbigvee\limits_{m}^{\lambda }\left(
f^{\left( n\right) }\right) +\left( M-\lambda \right)
^{n+1}\dbigvee\limits_{\lambda }^{M}\left( f^{\left( n\right) }\right) . 
\notag
\end{align}%
Taking the maximum over $\lambda \in \left[ m,M\right] $ in (\ref%
{III.e.e.3.7}) and utilizing (\ref{III.e.e.3.6}) we deduce the first
inequality in (\ref{III.e.e.3.2}).

Now observe that%
\begin{align*}
& \left( \lambda -m\right) ^{n+1}\dbigvee\limits_{m}^{\lambda }\left(
f^{\left( n\right) }\right) +\left( M-\lambda \right)
^{n+1}\dbigvee\limits_{\lambda }^{M}\left( f^{\left( n\right) }\right) \\
& \leq \max \left\{ \left( \lambda -m\right) ^{n+1},\left( M-\lambda \right)
^{n+1}\right\} \left[ \dbigvee\limits_{m}^{\lambda }\left( f^{\left(
n\right) }\right) +\dbigvee\limits_{\lambda }^{M}\left( f^{\left( n\right)
}\right) \right] \\
& =\max \left\{ \left( \lambda -m\right) ^{n+1},\left( M-\lambda \right)
^{n+1}\right\} \dbigvee\limits_{m}^{M}\left( f^{\left( n\right) }\right) \\
& =\left[ \frac{1}{2}\left( M-m\right) +\left\vert \lambda -\frac{m+M}{2}%
\right\vert \right] ^{n+1}\dbigvee\limits_{m}^{M}\left( f^{\left( n\right)
}\right)
\end{align*}%
giving that%
\begin{align*}
& \max_{\lambda \in \left[ m,M\right] }\left[ \left( \lambda -m\right)
^{n+1}\dbigvee\limits_{m}^{\lambda }\left( f^{\left( n\right) }\right)
+\left( M-\lambda \right) ^{n+1}\dbigvee\limits_{\lambda }^{M}\left(
f^{\left( n\right) }\right) \right] \\
& \leq \left( M-m\right) ^{n+1}\dbigvee\limits_{m}^{M}\left( f^{\left(
n\right) }\right)
\end{align*}%
and the second inequality in (\ref{III.e.e.3.2}) is proved.

The last part of (\ref{III.e.e.3.2}) follows by the Total Variation
Schwarz's inequality and we omit the details.

2. Now, recall that if $p:\left[ a,b\right] \rightarrow \mathbb{C}$ is a
Riemann integrable function and $v:\left[ a,b\right] \rightarrow \mathbb{C}$
is Lipschitzian with the constant $L>0$, i.e.,%
\begin{equation*}
\left\vert f\left( s\right) -f\left( t\right) \right\vert \leq L\left\vert
s-t\right\vert \text{ for any }t,s\in \left[ a,b\right] ,
\end{equation*}%
then the Riemann-Stieltjes integral $\int_{a}^{b}p\left( t\right) dv\left(
t\right) $ exists and the following inequality holds%
\begin{equation}
\left\vert \int_{a}^{b}p\left( t\right) dv\left( t\right) \right\vert \leq
L\int_{a}^{b}\left\vert p\left( t\right) \right\vert dt.  \label{III.e.e.3.9}
\end{equation}

By the property (\ref{III.e.e.3.9}) we have for $\lambda \in \left(
m,M\right) $ that%
\begin{equation*}
\left\vert \int_{m}^{\lambda }\left( t-m\right) ^{n+1}d\left( f^{\left(
n\right) }\left( t\right) \right) \right\vert \leq L_{n}\int_{m}^{\lambda
}\left( t-m\right) ^{n+1}d\left( t\right) =\frac{L_{n}}{n+2}\left( \lambda
-m\right) ^{n+2}
\end{equation*}%
and%
\begin{equation*}
\left\vert \int_{\lambda }^{M}\left( t-M\right) ^{n+1}d\left( f^{\left(
n\right) }\left( t\right) \right) \right\vert \leq L_{n}\int_{\lambda
}^{M}\left( M-t\right) ^{n+1}dt=\frac{L_{n}}{n+2}\left( M-\lambda \right)
^{n+2}.
\end{equation*}%
By the inequality (\ref{III.e.e.3.6}) we then have%
\begin{align}
& \left\vert T_{n}\left( A,m,M;x,y\right) \right\vert   \label{III.e.e.3.10}
\\
& \leq \frac{1}{\left( M-m\right) \left( n+1\right) !}\dbigvee%
\limits_{m-0}^{M}\left( \left\langle E_{\left( \cdot \right)
}x,y\right\rangle \right)   \notag \\
& \times \max_{\lambda \in \left[ m,M\right] }\left[ \frac{L_{n}}{n+2}\left(
\lambda -m\right) ^{n+2}+\frac{L_{n}}{n+2}\left( M-\lambda \right) ^{n+2}%
\right]   \notag \\
& =\frac{L_{n}\left( M-m\right) ^{n+1}}{\left( n+2\right) !}%
\dbigvee\limits_{m-0}^{M}\left( \left\langle E_{\left( \cdot \right)
}x,y\right\rangle \right) \leq \frac{L_{n}\left( M-m\right) ^{n+1}}{\left(
n+2\right) !}\left\Vert x\right\Vert \left\Vert y\right\Vert   \notag
\end{align}%
for any $x,y\in H$ and the inequality (\ref{III.e.e.3.3}) is proved.

3. Further, from the theory of Riemann-Stieltjes integral it is also well
known that if $p:\left[ a,b\right] \rightarrow \mathbb{C}$ is continuous and 
$v:\left[ a,b\right] \rightarrow \mathbb{R}$ is monotonic nondecreasing,
then the Riemann-Stieltjes integrals $\int_{a}^{b}p\left( t\right) dv\left(
t\right) $ and $\int_{a}^{b}\left\vert p\left( t\right) \right\vert dv\left(
t\right) $ exist and%
\begin{equation}
\left\vert \int_{a}^{b}p\left( t\right) dv\left( t\right) \right\vert \leq
\int_{a}^{b}\left\vert p\left( t\right) \right\vert dv\left( t\right) \leq
\max_{t\in \left[ a,b\right] }\left\vert p\left( t\right) \right\vert \left[
v\left( b\right) -v\left( a\right) \right] .  \label{III.e.e.3.11}
\end{equation}

On making use of (\ref{III.e.e.3.11}) we have%
\begin{align}
\left\vert \int_{m}^{\lambda }\left( t-m\right) ^{n+1}d\left( f^{\left(
n\right) }\left( t\right) \right) \right\vert & \leq \int_{m}^{\lambda
}\left( t-m\right) ^{n+1}d\left( f^{\left( n\right) }\left( t\right) \right)
\label{III.e.e.3.12} \\
& \leq \left( \lambda -m\right) ^{n+1}\left[ f^{\left( n\right) }\left(
\lambda \right) -f^{\left( n\right) }\left( m\right) \right]  \notag
\end{align}%
and%
\begin{align}
\left\vert \int_{\lambda }^{M}\left( t-M\right) ^{n+1}d\left( f^{\left(
n\right) }\left( t\right) \right) \right\vert & \leq \int_{\lambda
}^{M}\left( M-t\right) ^{n+1}d\left( f^{\left( n\right) }\left( t\right)
\right)  \label{III.e.e.3.13} \\
& \leq \left( M-\lambda \right) ^{n+1}\left[ f^{\left( n\right) }\left(
M\right) -f^{\left( n\right) }\left( \lambda \right) \right]  \notag
\end{align}%
for any $\lambda \in \left( m,M\right) .$

Integrating by parts in the Riemann-Stieltjes integral, we also have%
\begin{align*}
& \int_{m}^{\lambda }\left( t-m\right) ^{n+1}d\left( f^{\left( n\right)
}\left( t\right) \right) \\
& =\left( \lambda -m\right) ^{n+1}f^{\left( n\right) }\left( \lambda \right)
-\left( n+1\right) \int_{m}^{\lambda }\left( t-m\right) ^{n}f^{\left(
n\right) }\left( t\right) dt
\end{align*}%
and%
\begin{align*}
& \int_{\lambda }^{M}\left( M-t\right) ^{n+1}d\left( f^{\left( n\right)
}\left( t\right) \right) \\
& =\left( n+1\right) \int_{\lambda }^{M}\left( M-t\right) ^{n}f^{\left(
n\right) }\left( t\right) dt-\left( M-\lambda \right) ^{n+1}f^{\left(
n\right) }\left( \lambda \right)
\end{align*}%
for any $\lambda \in \left( m,M\right) .$

Therefore, by adding (\ref{III.e.e.3.12}) with (\ref{III.e.e.3.13}) we get%
\begin{align*}
& \left\vert \int_{m}^{\lambda }\left( t-m\right) ^{n+1}d\left( f^{\left(
n\right) }\left( t\right) \right) \right\vert +\left\vert \int_{\lambda
}^{M}\left( t-M\right) ^{n+1}d\left( f^{\left( n\right) }\left( t\right)
\right) \right\vert \\
& \leq \left[ f^{\left( n\right) }\left( \lambda \right) \left( \left(
\lambda -m\right) ^{n+1}-\left( M-\lambda \right) ^{n+1}\right) \right] \\
& +\left( n+1\right) \left[ \int_{\lambda }^{M}\left( M-t\right)
^{n}f^{\left( n\right) }\left( t\right) dt-\int_{m}^{\lambda }\left(
t-m\right) ^{n}f^{\left( n\right) }\left( t\right) dt\right] \\
& \leq \left( \lambda -m\right) ^{n+1}\left[ f^{\left( n\right) }\left(
\lambda \right) -f^{\left( n\right) }\left( m\right) \right] +\left(
M-\lambda \right) ^{n+1}\left[ f^{\left( n\right) }\left( M\right)
-f^{\left( n\right) }\left( \lambda \right) \right]
\end{align*}%
for any $\lambda \in \left( m,M\right) .$

Now, on making use of the inequality (\ref{III.e.e.3.6}) we deduce (\ref%
{III.e.e.3.3.a}).
\end{proof}

\begin{remark}
\label{III.e.r.2.2}If we use the inequality (\ref{III.e.e.3.2}) for the
function $\ln $, then we get the inequality%
\begin{align}
& \left\vert L_{n}\left( A,m,M;x,y\right) \right\vert   \label{III.e.e.3.14}
\\
& \leq \frac{1}{\left( M-m\right) n\left( n+1\right) }\dbigvee%
\limits_{m-0}^{M}\left( \left\langle E_{\left( \cdot \right)
}x,y\right\rangle \right)   \notag \\
& \times \max_{\lambda \in \left[ m,M\right] }\left[ \left( \lambda
-m\right) ^{n+1}\frac{\lambda ^{n}-m^{n}}{\lambda ^{n}m^{n}}+\left(
M-\lambda \right) ^{n+1}\frac{M^{n}-\lambda ^{n}}{M^{n}\lambda ^{n}}\right] 
\notag \\
& \leq \frac{\left( M-m\right) ^{n}\left( M^{n}-m^{n}\right) }{n\left(
n+1\right) M^{n}m^{n}}\dbigvee\limits_{m-0}^{M}\left( \left\langle E_{\left(
\cdot \right) }x,y\right\rangle \right)   \notag \\
& \leq \frac{\left( M-m\right) ^{n}\left( M^{n}-m^{n}\right) }{n\left(
n+1\right) M^{n}m^{n}}\left\Vert x\right\Vert \left\Vert y\right\Vert  
\notag
\end{align}%
for any $x,y\in H,$ where%
\begin{align}
& L_{n}\left( A,m,M;x,y\right)   \label{III.e.e.3.15} \\
& :=\left\langle \ln Ax,y\right\rangle -\left[ \ln I\left( m,M\right) \right]
\left\langle x,y\right\rangle   \notag \\
& -\frac{1}{M-m}\sum_{k=1}^{n}\frac{1}{k\left( k+1\right) }  \notag \\
& \times \left[ \left\langle \left( A-m1_{H}\right)
^{k+1}A^{-k}x,y\right\rangle +\left( -1\right) ^{k}\left\langle \left(
M1_{H}-A\right) ^{k+1}A^{-k}x,y\right\rangle \right] .  \notag
\end{align}%
If we use the inequality (\ref{III.e.e.3.3}) for the function $\ln $ we get
the following bound as well 
\begin{align}
& \left\vert L_{n}\left( A,m,M;x,y\right) \right\vert   \label{III.e.e.3.16}
\\
& \leq \frac{1}{\left( n+1\right) \left( n+2\right) }\left( \frac{M}{m}%
-1\right) ^{n+1}\dbigvee\limits_{m-0}^{M}\left( \left\langle E_{\left( \cdot
\right) }x,y\right\rangle \right)   \notag \\
& \leq \frac{1}{\left( n+1\right) \left( n+2\right) }\left( \frac{M}{m}%
-1\right) ^{n+1}\left\Vert x\right\Vert \left\Vert y\right\Vert   \notag
\end{align}%
for any $x,y\in H.$
\end{remark}

\begin{remark}
\label{III.e.r.2.3}If we define%
\begin{multline}
E_{n}\left( A,m,M;x,y\right)   \label{III.e.e.3.17} \\
:=\left\langle \left[ 1_{H}+\frac{1}{M-m}\sum_{k=1}^{n}\frac{1}{\left(
k+1\right) !}\left[ \left( M1_{H}-A\right) ^{k+1}+\left( -1\right)
^{k}\left( A-m1_{H}\right) ^{k+1}\right] \right] \exp Ax,y\right\rangle  \\
-E\left( m,M\right) \left\langle x,y\right\rangle ,
\end{multline}%
then by the inequality (\ref{III.e.e.3.2}) we have 
\begin{align}
& \left\vert E_{n}\left( A,m,M;x,y\right) \right\vert   \label{III.e.e.3.18}
\\
& \leq \frac{1}{\left( M-m\right) \left( n+1\right) !}\dbigvee%
\limits_{m-0}^{M}\left( \left\langle E_{\left( \cdot \right)
}x,y\right\rangle \right)   \notag \\
& \times \max_{\lambda \in \left[ m,M\right] }\left[ \left( \lambda
-m\right) ^{n+1}\left( e^{\lambda }-e^{m}\right) +\left( M-\lambda \right)
^{n+1}\left( e^{M}-e^{\lambda }\right) \right]   \notag \\
& \leq \frac{\left( M-m\right) ^{n}}{\left( n+1\right) !}\dbigvee%
\limits_{m-0}^{M}\left( \left\langle E_{\left( \cdot \right)
}x,y\right\rangle \right) \left( e^{M}-e^{m}\right) \leq \frac{\left(
M-m\right) ^{n}}{\left( n+1\right) !}\left( e^{M}-e^{m}\right) \left\Vert
x\right\Vert \left\Vert y\right\Vert   \notag
\end{align}%
for any $x,y\in H.$

If we use the inequality (\ref{III.e.e.3.3}) for the function $\exp $ we get
the following bound as well 
\begin{align}
\left\vert E_{n}\left( A,m,M;x,y\right) \right\vert & \leq \frac{e^{M}\left(
M-m\right) ^{n+1}}{\left( n+2\right) !}\dbigvee\limits_{m-0}^{M}\left(
\left\langle E_{\left( \cdot \right) }x,y\right\rangle \right) 
\label{III.e.e.3.19} \\
& \leq \frac{e^{M}\left( M-m\right) ^{n+1}}{\left( n+2\right) !}\left\Vert
x\right\Vert \left\Vert y\right\Vert   \notag
\end{align}%
for any $x,y\in H.$
\end{remark}

\subsection{Error Bounds for $f^{\left( n\right) }$ Absolutely Continuous}

We consider the Lebesgue norms defined by%
\begin{equation*}
\left\Vert g\right\Vert _{\left[ a,b\right] ,\infty }:=ess\sup_{t\in \left[
a,b\right] }\left\vert g\left( t\right) \right\vert \text{ if }g\in
L_{\infty }\left[ a,b\right]
\end{equation*}%
and%
\begin{equation*}
\left\Vert g\right\Vert _{\left[ a,b\right] ,p}:=\left(
\int_{a}^{b}\left\vert g\left( t\right) \right\vert ^{p}dt\right) ^{1/p}%
\text{ if }g\in L_{p}\left[ a,b\right] ,p\geq 1.
\end{equation*}

\begin{theorem}[Dragomir, 2010, \protect\cite{III.e.SSD7}]
\label{III.e.t.4.1}Let $A$ be a selfadjoint operator in the Hilbert space $H$
with the spectrum $Sp\left( A\right) \subseteq \left[ m,M\right] $ for some
real numbers $m<M$, $\left\{ E_{\lambda }\right\} _{\lambda }$ be its 
\textit{spectral family,} $I$ be a closed subinterval on $\mathbb{R}$ with $%
\left[ m,M\right] \subset \mathring{I}$ and let $n$ be an integer with $%
n\geq 1.$ If the $n$-th derivative $f^{\left( n\right) }$ is absolutely
continuous on $\left[ m,M\right] $, then%
\begin{multline}
\left\vert T_{n}\left( A,m,M;x,y\right) \right\vert \leq \frac{1}{\left(
M-m\right) \left( n+1\right) !}  \label{III.e.e.4.1} \\
\times \int_{m-0}^{M}\left\vert \left( \lambda -m\right) ^{n+1}\left\langle
\left( 1_{H}-E_{\lambda }\right) x,y\right\rangle +\left( \lambda -M\right)
^{n+1}\left\langle E_{\lambda }x,y\right\rangle \right\vert \left\vert
f^{\left( n+1\right) }\left( \lambda \right) \right\vert d\lambda . \\
\leq \frac{1}{\left( M-m\right) \left( n+1\right) !} \\
\times \left\{ 
\begin{array}{l}
B_{n,1}\left( A,m,M;x,y\right) \left\Vert f^{\left( n\right) }\right\Vert _{ 
\left[ m,M\right] ,\infty }\text{ if }f^{\left( n\right) }\in L_{\infty }%
\left[ m,M\right] , \\ 
\\ 
B_{n,p}\left( A,m,M;x,y\right) \left\Vert f^{\left( n\right) }\right\Vert _{ 
\left[ m,M\right] ,q}\text{ if }f^{\left( n\right) }\in L_{q}\left[ m,M%
\right] ,p>1,\frac{1}{p}+\frac{1}{q}=1, \\ 
\\ 
B_{n,\infty }\left( A,m,M;x,y\right) \left\Vert f^{\left( n\right)
}\right\Vert _{\left[ m,M\right] ,1},%
\end{array}%
\right.
\end{multline}%
for any $x,y\in H,$ where%
\begin{align*}
& B_{n,p}\left( A,m,M;x,y\right) \\
& :=\left( \int_{m-0}^{M}\left\vert \left( \lambda -m\right)
^{n+1}\left\langle \left( 1_{H}-E_{\lambda }\right) x,y\right\rangle +\left(
\lambda -M\right) ^{n+1}\left\langle E_{\lambda }x,y\right\rangle
\right\vert ^{p}d\lambda \right) ^{1/p},p\geq 1
\end{align*}%
and%
\begin{align*}
& B_{n,\infty }\left( A,m,M;x,y\right) \\
& :=\sup_{t\in \left[ m,M\right] }\left\vert \left( \lambda -m\right)
^{n+1}\left\langle \left( 1_{H}-E_{\lambda }\right) x,y\right\rangle +\left(
\lambda -M\right) ^{n+1}\left\langle E_{\lambda }x,y\right\rangle
\right\vert .
\end{align*}
\end{theorem}

\begin{proof}
Follows from the representation%
\begin{align*}
& T_{n}\left( A,m,M;x,y\right) \\
& =\frac{\left( -1\right) ^{n}}{\left( M-m\right) \left( n+1\right) !} \\
& \times \int_{m-0}^{M}\left[ \left( \lambda -m\right) ^{n+1}\left\langle
\left( 1_{H}-E_{\lambda }\right) x,y\right\rangle +\left( \lambda -M\right)
^{n+1}\left\langle E_{\lambda }x,y\right\rangle \right] f^{\left( n+1\right)
}\left( \lambda \right) d\lambda
\end{align*}%
for any $x,y\in H,$ by taking the modulus and utilizing the H\"{o}lder
integral inequality.

The details are omitted.
\end{proof}

The bounds provided by $B_{n,p}\left( A,m,M;x,y\right) $ are not useful for
applications, therefore we will establish in the following some simpler,
however coarser bounds.

\begin{proposition}[Dragomir, 2010, \protect\cite{III.e.SSD7}]
\label{III.e.p.4.1}With the above notations, we have%
\begin{equation}
B_{n,\infty }\left( A,m,M;x,y\right) \leq \left( M-m\right) ^{n+1}\left\Vert
x\right\Vert \left\Vert y\right\Vert ,  \label{III.e.e.4.2}
\end{equation}%
\begin{equation}
B_{n,1}\left( A,m,M;x,y\right) \leq \frac{\left( 2^{n+2}-1\right) }{\left(
n+2\right) 2^{n+1}}\left( M-m\right) ^{n+2}\left\Vert x\right\Vert
\left\Vert y\right\Vert  \label{III.e.e.4.3}
\end{equation}%
and for $p>1$%
\begin{equation}
B_{n,p}\left( A,m,M;x,y\right) \leq \frac{\left( 2^{\left( n+1\right)
p+1}-1\right) ^{1/p}}{2^{n+1}\left[ \left( n+1\right) p+1\right] ^{1/p}}%
\left( M-m\right) ^{n+1+1/p}\left\Vert x\right\Vert \left\Vert y\right\Vert
\label{III.e.e.4.4}
\end{equation}%
for any $x,y\in H.$
\end{proposition}

\begin{proof}
Utilising the triangle inequality for the modulus we have%
\begin{align}
& \left\vert \left( \lambda -m\right) ^{n+1}\left\langle \left(
1_{H}-E_{\lambda }\right) x,y\right\rangle +\left( \lambda -M\right)
^{n+1}\left\langle E_{\lambda }x,y\right\rangle \right\vert
\label{III.e.e.4.5} \\
& \leq \left( \lambda -m\right) ^{n+1}\left\vert \left\langle \left(
1_{H}-E_{\lambda }\right) x,y\right\rangle \right\vert +\left( M-\lambda
\right) ^{n+1}\left\vert \left\langle E_{\lambda }x,y\right\rangle
\right\vert  \notag \\
& \leq \max \left\{ \left( \lambda -m\right) ^{n+1},\left( M-\lambda \right)
^{n+1}\right\} \left[ \left\vert \left\langle \left( 1_{H}-E_{\lambda
}\right) x,y\right\rangle \right\vert +\left\vert \left\langle E_{\lambda
}x,y\right\rangle \right\vert \right]  \notag
\end{align}%
for any $x,y\in H.$

Utilising the generalization of Schwarz's inequality for nonnegative
selfadjoint operators we have%
\begin{equation*}
\left\vert \left\langle \left( 1_{H}-E_{\lambda }\right) x,y\right\rangle
\right\vert \leq \left\langle \left( 1_{H}-E_{\lambda }\right)
x,x\right\rangle ^{1/2}\left\langle \left( 1_{H}-E_{\lambda }\right)
y,y\right\rangle ^{1/2}
\end{equation*}%
and%
\begin{equation*}
\left\vert \left\langle E_{\lambda }x,y\right\rangle \right\vert \leq
\left\langle E_{\lambda }x,x\right\rangle ^{1/2}\left\langle E_{\lambda
}y,y\right\rangle ^{1/2}
\end{equation*}%
for any $x,y\in H$ and $\lambda \in \left[ m,M\right] .$

Further, by making use of the elementary inequality 
\begin{equation*}
ac+bd\leq \left( a^{2}+b^{2}\right) ^{1/2}\left( c^{2}+d^{2}\right)
^{1/2},a,b,c,d\geq 0
\end{equation*}%
we have%
\begin{align}
& \left\vert \left\langle \left( 1_{H}-E_{\lambda }\right) x,y\right\rangle
\right\vert +\left\vert \left\langle E_{\lambda }x,y\right\rangle \right\vert
\label{III.e.e.4.6} \\
& \leq \left\langle \left( 1_{H}-E_{\lambda }\right) x,x\right\rangle
^{1/2}\left\langle \left( 1_{H}-E_{\lambda }\right) y,y\right\rangle
^{1/2}+\left\langle E_{\lambda }x,x\right\rangle ^{1/2}\left\langle
E_{\lambda }y,y\right\rangle ^{1/2}  \notag \\
& \leq \left( \left\langle \left( 1_{H}-E_{\lambda }\right) x,x\right\rangle
+\left\langle E_{\lambda }x,x\right\rangle \right) ^{1/2}\left( \left\langle
\left( 1_{H}-E_{\lambda }\right) y,y\right\rangle +\left\langle E_{\lambda
}y,y\right\rangle \right) ^{1/2}  \notag \\
& =\left\Vert x\right\Vert \left\Vert y\right\Vert  \notag
\end{align}%
for any $x,y\in H$ and $\lambda \in \left[ m,M\right] .$

Combining (\ref{III.e.e.4.5}) with (\ref{III.e.e.4.6}) we deduce that%
\begin{align}
& \left\vert \left( \lambda -m\right) ^{n+1}\left\langle \left(
1_{H}-E_{\lambda }\right) x,y\right\rangle +\left( \lambda -M\right)
^{n+1}\left\langle E_{\lambda }x,y\right\rangle \right\vert
\label{III.e.e.4.7} \\
& \leq \max \left\{ \left( \lambda -m\right) ^{n+1},\left( M-\lambda \right)
^{n+1}\right\} \left\Vert x\right\Vert \left\Vert y\right\Vert  \notag
\end{align}%
for any $x,y\in H$ and $\lambda \in \left[ m,M\right] .$

Taking the supremum over $\lambda \in \left[ m,M\right] $ in (\ref%
{III.e.e.4.7}) we deduce the inequality (\ref{III.e.e.4.2}).

Now, if we take the power $r\geq 1$ in (\ref{III.e.e.4.7}) and integrate,
then we get%
\begin{align}
& \int_{m-0}^{M}\left\vert \left( \lambda -m\right) ^{n+1}\left\langle
\left( 1_{H}-E_{\lambda }\right) x,y\right\rangle +\left( \lambda -M\right)
^{n+1}\left\langle E_{\lambda }x,y\right\rangle \right\vert ^{r}d\lambda
\label{III.e.e.4.8} \\
& \leq \left\Vert x\right\Vert ^{r}\left\Vert y\right\Vert
^{r}\int_{m}^{M}\max \left\{ \left( \lambda -m\right) ^{\left( n+1\right)
r},\left( M-\lambda \right) ^{\left( n+1\right) r}\right\} d\lambda  \notag
\\
& =\left\Vert x\right\Vert ^{r}\left\Vert y\right\Vert ^{r}\left[ \int_{m}^{%
\frac{M+m}{2}}\left( M-\lambda \right) ^{\left( n+1\right) r}d\lambda +\int_{%
\frac{M+m}{2}}^{M}\left( \lambda -m\right) ^{\left( n+1\right) r}d\lambda %
\right]  \notag \\
& =\frac{\left( 2^{\left( n+1\right) r+1}-1\right) }{\left[ \left(
n+1\right) r+1\right] 2^{\left( n+1\right) r}}\left( M-m\right) ^{\left(
n+1\right) r+1}\left\Vert x\right\Vert ^{r}\left\Vert y\right\Vert ^{r} 
\notag
\end{align}%
for any $x,y\in H.$

Utilizing (\ref{III.e.e.4.8}) for $r=1$ we deduce the bound (\ref%
{III.e.e.4.3}). Also, by making $r=p$ and then taking the power $1/p,$ we
deduce the last inequality (\ref{III.e.e.4.4}).
\end{proof}

The following result provides refinements of the inequalities in Proposition %
\ref{III.e.p.4.1}:

\begin{proposition}[Dragomir, 2010, \protect\cite{III.e.SSD7}]
\label{III.e.p.4.2}With the above notations, we have%
\begin{multline}
B_{n,\infty }\left( A,m,M;x,y\right)  \label{III.e.e.4.9} \\
\leq \left\Vert y\right\Vert \max_{\lambda \in \left[ m,M\right] }\left[
\left( \lambda -m\right) ^{2\left( n+1\right) }\left\langle \left(
1_{H}-E_{\lambda }\right) x,x\right\rangle +\left( M-\lambda \right)
^{2\left( n+1\right) }\left\langle E_{\lambda }x,x\right\rangle \right]
^{1/2} \\
\leq \left( M-m\right) ^{n+1}\left\Vert x\right\Vert \left\Vert y\right\Vert
,
\end{multline}%
\begin{multline}
B_{n,1}\left( A,m,M;x,y\right)  \label{III.e.e.4.10} \\
\leq \left\Vert y\right\Vert \int_{m-0}^{M}\left[ \left( \lambda -m\right)
^{2\left( n+1\right) }\left\langle \left( 1_{H}-E_{\lambda }\right)
x,x\right\rangle +\left( M-\lambda \right) ^{2\left( n+1\right)
}\left\langle E_{\lambda }x,x\right\rangle \right] ^{1/2}d\lambda \\
\leq \frac{\left( 2^{n+2}-1\right) }{\left( n+2\right) 2^{n+1}}\left(
M-m\right) ^{n+2}\left\Vert x\right\Vert \left\Vert y\right\Vert
\end{multline}%
and for $p>1$%
\begin{multline}
B_{n,p}\left( A,m,M;x,y\right)  \label{III.e.e.4.11} \\
\leq \left\Vert y\right\Vert \left( \int_{m-0}^{M}\left[ \left( \lambda
-m\right) ^{2\left( n+1\right) }\left\langle \left( 1_{H}-E_{\lambda
}\right) x,x\right\rangle +\left( M-\lambda \right) ^{2\left( n+1\right)
}\left\langle E_{\lambda }x,x\right\rangle \right] ^{p/2}d\lambda \right)
^{1/p} \\
\leq \frac{\left( 2^{\left( n+1\right) p+1}-1\right) ^{1/p}}{2^{n+1}\left[
\left( n+1\right) p+1\right] ^{1/p}}\left( M-m\right) ^{n+1+1/p}\left\Vert
x\right\Vert \left\Vert y\right\Vert
\end{multline}%
for any $x,y\in H.$
\end{proposition}

\begin{proof}
Utilising the Schwarz inequality in $H$, we have%
\begin{align}
& \left\vert \left\langle \left( \lambda -m\right) ^{n+1}\left(
1_{H}-E_{\lambda }\right) x+\left( \lambda -M\right) ^{n+1}E_{\lambda
}x,y\right\rangle \right\vert  \label{III.e.e.4.12} \\
& \leq \left\Vert y\right\Vert \left\Vert \left( \lambda -m\right)
^{n+1}\left( 1_{H}-E_{\lambda }\right) x+\left( \lambda -M\right)
^{n+1}E_{\lambda }x\right\Vert  \notag
\end{align}%
for any $x,y\in H.$

Since $E_{\lambda }$ are projectors for each $\lambda \in \left[ m,M\right]
, $ then we have%
\begin{align}
& \left\Vert \left( \lambda -m\right) ^{n+1}\left( 1_{H}-E_{\lambda }\right)
x+\left( \lambda -M\right) ^{n+1}E_{\lambda }x\right\Vert ^{2}
\label{III.e.e.4.13} \\
& =\left( \lambda -m\right) ^{2\left( n+1\right) }\left\Vert \left(
1_{H}-E_{\lambda }\right) x\right\Vert ^{2}  \notag \\
& +2\left( \lambda -m\right) ^{n+1}\left( \lambda -M\right) ^{n+1}\func{Re}%
\left\langle \left( 1_{H}-E_{\lambda }\right) x,E_{\lambda }x\right\rangle 
\notag \\
& +\left( M-\lambda \right) ^{2\left( n+1\right) }\left\Vert E_{\lambda
}x\right\Vert ^{2}  \notag \\
& =\left( \lambda -m\right) ^{2\left( n+1\right) }\left\Vert \left(
1_{H}-E_{\lambda }\right) x\right\Vert ^{2}+\left( M-\lambda \right)
^{2\left( n+1\right) }\left\Vert E_{\lambda }x\right\Vert ^{2}  \notag \\
& =\left( \lambda -m\right) ^{2\left( n+1\right) }\left\langle \left(
1_{H}-E_{\lambda }\right) x,x\right\rangle +\left( M-\lambda \right)
^{2\left( n+1\right) }\left\langle E_{\lambda }x,x\right\rangle  \notag \\
& \leq \left\Vert x\right\Vert ^{2}\max \left\{ \left( \lambda -m\right)
^{2\left( n+1\right) },\left( M-\lambda \right) ^{2\left( n+1\right)
}\right\}  \notag
\end{align}%
for any $x,y\in H$ and $\lambda \in \left[ m,M\right] .$

On making use of (\ref{III.e.e.4.12}) and (\ref{III.e.e.4.13}) we obtain the
following refinement of (\ref{III.e.e.4.7})%
\begin{align}
& \left\vert \left\langle \left( \lambda -m\right) ^{n+1}\left(
1_{H}-E_{\lambda }\right) x+\left( \lambda -M\right) ^{n+1}E_{\lambda
}x,y\right\rangle \right\vert  \label{III.e.e.4.14} \\
& \leq \left\Vert y\right\Vert \left[ \left( \lambda -m\right) ^{2\left(
n+1\right) }\left\langle \left( 1_{H}-E_{\lambda }\right) x,x\right\rangle
+\left( M-\lambda \right) ^{2\left( n+1\right) }\left\langle E_{\lambda
}x,x\right\rangle \right] ^{1/2}  \notag \\
& \leq \max \left\{ \left( \lambda -m\right) ^{n+1},\left( M-\lambda \right)
^{n+1}\right\} \left\Vert x\right\Vert \left\Vert y\right\Vert  \notag
\end{align}%
for any $x,y\in H$ and $\lambda \in \left[ m,M\right] .$

The proof now follows the lines of the proof from Proposition \ref%
{III.e.p.4.1} and we omit the details.
\end{proof}

\begin{remark}
\label{III.e.r.4.1}One can apply Theorem \ref{III.e.t.4.1} and Proposition %
\ref{III.e.p.4.1} for particular functions including the exponential and
logarithmic function. However the details are left to the interested reader.
\end{remark}

\bigskip

\chapter{Inequalities of Trapezoidal Type}

\section{Introduction}

From a complementary viewpoint to Ostrowski/mid-point inequalities,
trapezoidal type inequality provide a priory error bounds in approximating
the Riemann integral by a (generalized) trapezoidal formula.

Just like in the case of Ostrowski's inequality the development of these
kind of results have registered a sharp growth in the last decade with more
than 50 papers published, as one can easily asses this by performing a
search with the key word "trapezoid" and "inequality" in the title of the
papers reviewed by MathSciNet data base of the American Mathematical Society.

Numerous extensions, generalisations in both the integral and discrete case
have been discovered. More general versions for $n$-time differentiable
functions, the corresponding versions on time scales, for vector valued
functions or multiple integrals have been established as well. Numerous
applications in Numerical Analysis, Probability Theory and other fields have
been also given.

In the present chapter we present some recent results obtained by the author
in extending trapezoidal type inequality in various directions for
continuous functions of selfadjoint operators in complex Hilbert spaces. As
far as we know, the obtained results are new with no previous similar
results ever obtained in the literature.

Applications for some elementary functions of operators such as the power
function, the logarithmic and exponential functions are provided as well.

\section{Scalar Trapezoidal Type Inequalities}

In Classical Analysis a \textit{trapezoidal type inequality} is an
inequality that provides upper and/or lower bounds for the quantity%
\begin{equation*}
\frac{f\left( a\right) +f\left( b\right) }{2}\left( b-a\right)
-\int_{a}^{b}f\left( t\right) dt,
\end{equation*}%
that is the error in approximating the integral by a trapezoidal rule, for
various classes of integrable functions $f$ defined on the compact interval $%
\left[ a,b\right] .$

In order to introduce the reader to some of the well known results and
prepare the background for considering a similar problem for functions of
selfadjoint operators in Hilbert spaces, we mention the following
inequalities.

The case of functions of bounded variation was obtained in \cite{IV.2b} (see
also \cite[p. 68]{IV.1b}):

\begin{theorem}
\label{IV.ta.1}Let $f:\left[ a,b\right] \rightarrow \mathbb{C}$ be a
function of bounded variation. We have the inequality 
\begin{equation}
\left\vert \int_{a}^{b}f\left( t\right) dt-\frac{f\left( a\right) +f\left(
b\right) }{2}\left( b-a\right) \right\vert \leq \frac{1}{2}\left( b-a\right)
\bigvee_{a}^{b}\left( f\right) ,  \label{IV.a.1}
\end{equation}%
where $\bigvee_{a}^{b}\left( f\right) $ denotes the total variation of $f$
on the interval $\left[ a,b\right] $. The constant $\frac{1}{2}$ is the best
possible one.
\end{theorem}

This result may be improved if one assumes the monotonicity of $f$ as
follows (see \cite[p. 76]{IV.1b}):

\begin{theorem}
\label{IV.ta.2}Let $f:\left[ a,b\right] \rightarrow \mathbb{R}$ be a
monotonic nondecreasing function on $\left[ a,b\right] $. Then we have the
inequalities: 
\begin{align}
& \left\vert \int_{a}^{b}f\left( t\right) dt-\frac{f\left( a\right) +f\left(
b\right) }{2}\left( b-a\right) \right\vert  \label{IV.a.2} \\
& \leq \frac{1}{2}\left( b-a\right) \left[ f\left( b\right) -f\left(
a\right) \right] -\int_{a}^{b}sgn\left( t-\frac{a+b}{2}\right) f\left(
t\right) dt  \notag \\
& \leq \frac{1}{2}\left( b-a\right) \left[ f\left( b\right) -f\left(
a\right) \right] .  \notag
\end{align}%
The above inequalities are sharp.
\end{theorem}

If the mapping is Lipschitzian, then the following result holds as well \cite%
{IV.3b} (see also \cite[p. 82]{IV.1b}).

\begin{theorem}
\label{IV.ta.3}Let $f:\left[ a,b\right] \rightarrow \mathbb{C}$ be an $L-$%
Lipschitzian function on $\left[ a,b\right] ,$ i.e., $f$ satisfies the
condition: 
\begin{equation}
\left\vert f\left( s\right) -f\left( t\right) \right\vert \leq L\left\vert
s-t\right\vert \text{ \ for any \ }s,t\in \left[ a,b\right] \;\;\;\text{(}L>0%
\text{ is given).}  \tag{L}  \label{IV.L}
\end{equation}%
Then we have the inequality: 
\begin{equation}
\left\vert \int_{a}^{b}f\left( t\right) dt-\frac{f\left( a\right) +f\left(
b\right) }{2}\left( b-a\right) \right\vert \leq \frac{1}{4}\left( b-a\right)
^{2}L.  \label{IV.a.3}
\end{equation}%
The constant $\frac{1}{4}$ is best in (\ref{IV.a.3}).
\end{theorem}

If we would assume absolute continuity for the function $f$, then the
following estimates in terms of the Lebesgue norms of the derivative $%
f^{\prime }$ hold \cite[p. 93]{IV.1b}.

\begin{theorem}
\label{IV.ta.4}Let $f:\left[ a,b\right] \rightarrow \mathbb{C}$ be an
absolutely continuous function on $\left[ a,b\right] $. Then we have 
\begin{align}
& \left\vert \int_{a}^{b}f\left( t\right) dt-\frac{f\left( a\right) +f\left(
b\right) }{2}\left( b-a\right) \right\vert  \label{IV.a.4} \\
& \leq \left\{ 
\begin{array}{lll}
\dfrac{1}{4}\left( b-a\right) ^{2}\left\Vert f^{\prime }\right\Vert _{\infty
} & \text{if} & f^{\prime }\in L_{\infty }\left[ a,b\right] ; \\ 
&  &  \\ 
\dfrac{1}{2\left( q+1\right) ^{\frac{1}{q}}}\left( b-a\right)
^{1+1/q}\left\Vert f^{\prime }\right\Vert _{p} & \text{if} & f^{\prime }\in
L_{p}\left[ a,b\right] , \\ 
&  & p>1,\;\frac{1}{p}+\frac{1}{q}=1; \\ 
\dfrac{1}{2}\left( b-a\right) \left\Vert f^{\prime }\right\Vert _{1}, &  & 
\end{array}%
\right.  \notag
\end{align}%
where $\left\Vert \cdot \right\Vert _{p}$ $\left( p\in \left[ 1,\infty %
\right] \right) $ are the Lebesgue norms, i.e., 
\begin{equation*}
\left\Vert f^{\prime }\right\Vert _{\infty }=ess\sup\limits_{s\in \left[ a,b%
\right] }\left\vert f^{\prime }\left( s\right) \right\vert
\end{equation*}%
and 
\begin{equation*}
\left\Vert f^{\prime }\right\Vert _{p}:=\left( \int_{a}^{b}\left\vert
f^{\prime }\left( s\right) \right\vert ds\right) ^{\frac{1}{p}},\;\;p\geq 1.
\end{equation*}
\end{theorem}

The case of convex functions is as follows \cite{IV.4b}:

\begin{theorem}
\label{IV.t1}Let $f:\left[ a,b\right] \rightarrow \mathbb{R}$ be a convex
function on $\left[ a,b\right] .$ Then we have the inequalities 
\begin{align}
& \frac{1}{8}\left( b-a\right) ^{2}\left[ f_{+}^{\prime }\left( \frac{a+b}{2}%
\right) -f_{-}^{\prime }\left( \frac{a+b}{2}\right) \right]  \label{IV.2.1}
\\
& \leq \frac{f\left( a\right) +f\left( b\right) }{2}\left( b-a\right)
-\int_{a}^{b}f\left( t\right) dt  \notag \\
& \leq \frac{1}{8}\left( b-a\right) ^{2}\left[ f_{-}^{\prime }\left(
b\right) -f_{+}^{\prime }\left( a\right) \right] .  \notag
\end{align}%
The constant $\frac{1}{8}$ is sharp in both sides of (\ref{IV.2.1}).
\end{theorem}

For other scalar trapezoidal type inequalities, see \cite{IV.1b}.

\section{Trapezoidal Vector Inequalities}

\subsection{Some General Results}

With the notations introduced above, we consider in this paper the problem
of bounding the error%
\begin{equation*}
\frac{f\left( M\right) +f\left( m\right) }{2}\cdot \left\langle
x,y\right\rangle -\left\langle f\left( A\right) x,y\right\rangle
\end{equation*}%
in approximating $\left\langle f\left( A\right) x,y\right\rangle $ by the
trapezoidal type formula $\frac{f\left( M\right) +f\left( m\right) }{2}\cdot
\left\langle x,y\right\rangle ,$ where $x,y$ are vectors in the Hilbert
space $H,$ $f$ is a continuous functions of the selfadjoint operator $A$
with the spectrum in the compact interval of real numbers $\left[ m,M\right]
.$ Applications for some particular elementary functions are also provided.
The following result holds:

\begin{theorem}[Dragomir, 2010, \protect\cite{IV.SSD3}]
\label{IV.t.2.1}Let $A$ be a selfadjoint operator in the Hilbert space $H$
with the spectrum $Sp\left( A\right) \subseteq \left[ m,M\right] $ for some
real numbers $m<M$ and let $\left\{ E_{\lambda }\right\} _{\lambda }$ be its 
\textit{spectral family.} If $f:\left[ m,M\right] \rightarrow \mathbb{C}$ is
a continuous function of bounded variation on $\left[ m,M\right] $, then we
have the inequality%
\begin{align}
& \left\vert \frac{f\left( M\right) +f\left( m\right) }{2}\cdot \left\langle
x,y\right\rangle -\left\langle f\left( A\right) x,y\right\rangle \right\vert
\label{IV.e.2.2} \\
& \leq \frac{1}{2}\max_{\lambda \in \left[ m,M\right] }\left[ \left\langle
E_{\lambda }x,x\right\rangle ^{1/2}\left\langle E_{\lambda }y,y\right\rangle
^{1/2}\right.  \notag \\
& \left. +\left\langle \left( 1_{H}-E_{\lambda }\right) x,x\right\rangle
^{1/2}\left\langle \left( 1_{H}-E_{\lambda }\right) y,y\right\rangle ^{1/2} 
\right] \dbigvee\limits_{m}^{M}\left( f\right)  \notag \\
& \leq \frac{1}{2}\left\Vert x\right\Vert \left\Vert y\right\Vert
\dbigvee\limits_{m}^{M}\left( f\right)  \notag
\end{align}%
for any $x,y\in H.$
\end{theorem}

\begin{proof}
If $f,u:\left[ m,M\right] \rightarrow \mathbb{C}$ are such that the
Riemann-Stieltjes integral $\int_{a}^{b}f\left( t\right) du\left( t\right) $
exists, then a simple integration by parts reveals the identity%
\begin{align}
\int_{a}^{b}f\left( t\right) du\left( t\right) & =\frac{f\left( a\right)
+f\left( b\right) }{2}\left[ u\left( b\right) -u\left( a\right) \right]
\label{IV.e.2.3} \\
& -\int_{a}^{b}\left[ u\left( t\right) -\frac{u\left( a\right) +u\left(
b\right) }{2}\right] df\left( t\right) .  \notag
\end{align}

If we write the identity (\ref{IV.e.2.3}) for $u\left( \lambda \right)
=\left\langle E_{\lambda }x,y\right\rangle ,$ then we get%
\begin{align*}
\int_{m-0}^{M}f\left( \lambda \right) d\left( \left\langle E_{\lambda
}x,y\right\rangle \right) & =\frac{f\left( m\right) +f\left( M\right) }{2}%
\cdot \left\langle x,y\right\rangle \\
& -\int_{m-0}^{M}\left( \left\langle E_{\lambda }x,y\right\rangle -\frac{1}{2%
}\left\langle x,y\right\rangle \right) df\left( \lambda \right)
\end{align*}%
which gives the following identity of interest in itself%
\begin{align}
& \frac{f\left( m\right) +f\left( M\right) }{2}\cdot \left\langle
x,y\right\rangle -\left\langle f\left( A\right) x,y\right\rangle
\label{IV.e.2.4} \\
& =\frac{1}{2}\int_{m-0}^{M}\left[ \left\langle E_{\lambda }x,y\right\rangle
+\left\langle \left( E_{\lambda }-1_{H}\right) x,y\right\rangle \right]
df\left( \lambda \right) ,  \notag
\end{align}%
for any $x,y\in H.$

It is well known that if $p:\left[ a,b\right] \rightarrow \mathbb{C}$ is a
continuous function and $v:\left[ a,b\right] \rightarrow \mathbb{C}$ is of
bounded variation, then the Riemann-Stieltjes integral $\int_{a}^{b}p\left(
t\right) dv\left( t\right) $ exists and the following inequality holds%
\begin{equation}
\left\vert \int_{a}^{b}p\left( t\right) dv\left( t\right) \right\vert \leq
\max_{t\in \left[ a,b\right] }\left\vert p\left( t\right) \right\vert
\dbigvee\limits_{a}^{b}\left( v\right)  \label{IV.e.2.5}
\end{equation}%
where $\dbigvee\limits_{a}^{b}\left( v\right) $ denotes the total variation
of $v$ on $\left[ a,b\right] .$

Utilising the property (\ref{IV.e.2.5}), we have from (\ref{IV.e.2.4}) that%
\begin{align}
& \left\vert \frac{f\left( m\right) +f\left( M\right) }{2}\cdot \left\langle
x,y\right\rangle -\left\langle f\left( A\right) x,y\right\rangle \right\vert
\label{IV.e.2.6} \\
& \leq \frac{1}{2}\max_{\lambda \in \left[ m,M\right] }\left\vert
\left\langle E_{\lambda }x,y\right\rangle +\left\langle \left( E_{\lambda
}-1_{H}\right) x,y\right\rangle \right\vert \dbigvee\limits_{m}^{M}\left(
f\right)  \notag \\
& \leq \frac{1}{2}\left[ \max_{\lambda \in \left[ m,M\right] }\left[
\left\vert \left\langle E_{\lambda }x,y\right\rangle \right\vert +\left\vert
\left\langle \left( 1_{H}-E_{\lambda }\right) x,y\right\rangle \right\vert %
\right] \right] \dbigvee\limits_{m}^{M}\left( f\right) .  \notag
\end{align}%
If $P$ is a nonnegative operator on $H,$ i.e., $\left\langle
Px,x\right\rangle \geq 0$ for any $x\in H,$ then the following inequality is
a generalization of the Schwarz inequality in the Hilbert space $H$%
\begin{equation}
\left\vert \left\langle Px,y\right\rangle \right\vert ^{2}\leq \left\langle
Px,x\right\rangle \left\langle Py,y\right\rangle ,  \label{IV.e.2.7}
\end{equation}%
for any $x,y\in H.$

On applying the inequality (\ref{IV.e.2.7}) we have%
\begin{equation*}
\left\vert \left\langle E_{\lambda }x,y\right\rangle \right\vert \leq
\left\langle E_{\lambda }x,x\right\rangle ^{1/2}\left\langle E_{\lambda
}y,y\right\rangle ^{1/2}\text{ }
\end{equation*}%
and%
\begin{equation*}
\left\vert \left\langle \left( 1_{H}-E_{\lambda }\right) x,y\right\rangle
\right\vert \leq \left\langle \left( 1_{H}-E_{\lambda }\right)
x,x\right\rangle ^{1/2}\left\langle \left( 1_{H}-E_{\lambda }\right)
y,y\right\rangle ^{1/2},
\end{equation*}%
which, together with the elementary inequality for $a,b,c,d\geq 0$ 
\begin{equation*}
ab+cd\leq \left( a^{2}+c^{2}\right) ^{1/2}\left( b^{2}+d^{2}\right) ^{1/2}
\end{equation*}%
produce the inequalities%
\begin{align}
& \left\vert \left\langle E_{\lambda }x,y\right\rangle \right\vert
+\left\vert \left\langle \left( 1_{H}-E_{\lambda }\right) x,y\right\rangle
\right\vert  \label{IV.e.2.8} \\
& \leq \left\langle E_{\lambda }x,x\right\rangle ^{1/2}\left\langle
E_{\lambda }y,y\right\rangle ^{1/2}+\left\langle \left( 1_{H}-E_{\lambda
}\right) x,x\right\rangle ^{1/2}\left\langle \left( 1_{H}-E_{\lambda
}\right) y,y\right\rangle ^{1/2}  \notag \\
& \leq \left( \left\langle E_{\lambda }x,x\right\rangle +\left\langle \left(
1_{H}-E_{\lambda }\right) x,x\right\rangle \right) \left( \left\langle
E_{\lambda }y,y\right\rangle +\left\langle \left( 1_{H}-E_{\lambda }\right)
y,y\right\rangle \right)  \notag \\
& =\left\Vert x\right\Vert \left\Vert y\right\Vert  \notag
\end{align}%
for any $x,y\in H.$

On utilizing (\ref{IV.e.2.6}) and taking the maximum in (\ref{IV.e.2.8}) we
deduce the desired result (\ref{IV.e.2.2}).
\end{proof}

The case of Lipschitzian functions may be useful for applications:

\begin{theorem}[Dragomir, 2010, \protect\cite{IV.SSD3}]
\label{IV.t.2.2}Let $A$ be a selfadjoint operator in the Hilbert space $H$
with the spectrum $Sp\left( A\right) \subseteq \left[ m,M\right] $ for some
real numbers $m<M$ and let $\left\{ E_{\lambda }\right\} _{\lambda }$ be its 
\textit{spectral family.} If $f:\left[ m,M\right] \rightarrow \mathbb{C}$ is
Lipschitzian with the constant $L>0$ on $\left[ m,M\right] $, then we have
the inequality%
\begin{align}
& \left\vert \frac{f\left( M\right) +f\left( m\right) }{2}\cdot \left\langle
x,y\right\rangle -\left\langle f\left( A\right) x,y\right\rangle \right\vert
\label{IV.e.2.9} \\
& \leq \frac{1}{2}L\int_{m-0}^{M}\left[ \left\langle E_{\lambda
}x,x\right\rangle ^{1/2}\left\langle E_{\lambda }y,y\right\rangle
^{1/2}\right.  \notag \\
& \left. +\left\langle \left( 1_{H}-E_{\lambda }\right) x,x\right\rangle
^{1/2}\left\langle \left( 1_{H}-E_{\lambda }\right) y,y\right\rangle ^{1/2} 
\right] d\lambda  \notag \\
& \leq \frac{1}{2}\left( M-m\right) L\left\Vert x\right\Vert \left\Vert
y\right\Vert  \notag
\end{align}%
for any $x,y\in H.$
\end{theorem}

\begin{proof}
It is well known that if $p:\left[ a,b\right] \rightarrow \mathbb{C}$ is a
Riemann integrable function and $v:\left[ a,b\right] \rightarrow \mathbb{C}$
is Lipschitzian with the constant $L>0$, i.e.,%
\begin{equation*}
\left\vert f\left( s\right) -f\left( t\right) \right\vert \leq L\left\vert
s-t\right\vert \text{ for any }t,s\in \left[ a,b\right] ,
\end{equation*}%
then the Riemann-Stieltjes integral $\int_{a}^{b}p\left( t\right) dv\left(
t\right) $ exists and the following inequality holds%
\begin{equation*}
\left\vert \int_{a}^{b}p\left( t\right) dv\left( t\right) \right\vert \leq
L\int_{a}^{b}\left\vert p\left( t\right) \right\vert dt.
\end{equation*}

Now, on applying this property of the Riemann-Stieltjes integral, we have
from the representation (\ref{IV.e.2.4}) that%
\begin{align}
& \left\vert \frac{f\left( m\right) +f\left( M\right) }{2}\cdot \left\langle
x,y\right\rangle -\left\langle f\left( A\right) x,y\right\rangle \right\vert
\label{IV.e.2.10} \\
& \leq \frac{1}{2}L\int_{m-0}^{M}\left\vert \left\langle E_{\lambda
}x,y\right\rangle +\left\langle \left( E_{\lambda }-1_{H}\right)
x,y\right\rangle \right\vert d\lambda ,  \notag \\
& \leq \frac{1}{2}L\int_{m-0}^{M}\left[ \left\vert \left\langle E_{\lambda
}x,y\right\rangle \right\vert +\left\vert \left\langle \left(
1_{H}-E_{\lambda }\right) x,y\right\rangle \right\vert \right] d\lambda , 
\notag
\end{align}%
for any $x,y\in H.$

Further, integrating (\ref{IV.e.2.8}) on $\left[ m,M\right] $ we have%
\begin{align}
& \int_{m-0}^{M}\left[ \left\vert \left\langle E_{\lambda }x,y\right\rangle
\right\vert +\left\vert \left\langle \left( 1_{H}-E_{\lambda }\right)
x,y\right\rangle \right\vert \right] d\lambda  \label{IV.e.2.11} \\
& \leq \int_{m-0}^{M}\left[ \left\langle E_{\lambda }x,x\right\rangle
^{1/2}\left\langle E_{\lambda }y,y\right\rangle ^{1/2}\right.  \notag \\
& \left. +\left\langle \left( 1_{H}-E_{\lambda }\right) x,x\right\rangle
^{1/2}\left\langle \left( 1_{H}-E_{\lambda }\right) y,y\right\rangle ^{1/2} 
\right] d\lambda  \notag \\
& \leq \left( M-m\right) \left\Vert x\right\Vert \left\Vert y\right\Vert 
\notag
\end{align}%
which together with (\ref{IV.e.2.10}) produces the desired result (\ref%
{IV.e.2.9}).
\end{proof}

\subsection{Other Trapezoidal Vector Inequalities}

The following result provides a different perspective in bounding the error
in the trapezoidal approximation:

\begin{theorem}[Dragomir, 2010, \protect\cite{IV.SSD3}]
\label{IV.t.2.3}Let $A$ be a selfadjoint operator in the Hilbert space $H$
with the spectrum $Sp\left( A\right) \subseteq \left[ m,M\right] $ for some
real numbers $m<M$ and let $\left\{ E_{\lambda }\right\} _{\lambda }$ be its 
\textit{spectral family.} Assume that $f:\left[ m,M\right] \rightarrow 
\mathbb{C}$ is a continuous function on $\left[ m,M\right] $. Then we have
the inequalities%
\begin{align}
& \left\vert \frac{f\left( M\right) +f\left( m\right) }{2}\cdot \left\langle
x,y\right\rangle -\left\langle f\left( A\right) x,y\right\rangle \right\vert
\label{IV.e.2.12} \\
& \leq \left\{ 
\begin{array}{ll}
\max_{\lambda \in \left[ m,M\right] }\left\vert \left\langle E_{\lambda }x-%
\frac{1}{2}x,y\right\rangle \right\vert \dbigvee\limits_{m}^{M}\left(
f\right) & \text{if }f\text{ is of bounded variation} \\ 
&  \\ 
L\int_{m-0}^{M}\left\vert \left\langle E_{\lambda }x-\frac{1}{2}%
x,y\right\rangle \right\vert d\lambda & \text{if }f\text{ is }L\text{
Lipschitzian} \\ 
&  \\ 
\int_{m-0}^{M}\left\vert \left\langle E_{\lambda }x-\frac{1}{2}%
x,y\right\rangle \right\vert df\left( \lambda \right) & \text{if }f\text{ is
nondecreasing}%
\end{array}%
\right.  \notag \\
& \leq \frac{1}{2}\left\Vert x\right\Vert \left\Vert y\right\Vert \left\{ 
\begin{array}{ll}
\dbigvee\limits_{m}^{M}\left( f\right) & \text{if }f\text{ is of bounded
variation} \\ 
&  \\ 
L\left( M-m\right) & \text{if }f\text{ is }L\text{ Lipschitzian} \\ 
&  \\ 
\left( f\left( M\right) -f\left( m\right) \right) & \text{if }f\text{ is
nondecreasing}%
\end{array}%
\right.  \notag
\end{align}%
for any $x,y\in H.$
\end{theorem}

\begin{proof}
From (\ref{IV.e.2.6}) we have that%
\begin{align}
& \left\vert \frac{f\left( m\right) +f\left( M\right) }{2}\cdot \left\langle
x,y\right\rangle -\left\langle f\left( A\right) x,y\right\rangle \right\vert
\label{IV.e.2.13} \\
& \leq \frac{1}{2}\max_{\lambda \in \left[ m,M\right] }\left\vert
\left\langle E_{\lambda }x,y\right\rangle +\left\langle \left( E_{\lambda
}-1_{H}\right) x,y\right\rangle \right\vert \dbigvee\limits_{m}^{M}\left(
f\right)  \notag \\
& =\max_{\lambda \in \left[ m,M\right] }\left\vert \left\langle E_{\lambda
}x-\frac{1}{2}x,y\right\rangle \right\vert \dbigvee\limits_{m}^{M}\left(
f\right)  \notag
\end{align}%
for any $x,y\in H.$

Utilising the Schwarz inequality in $H$ and the fact that $E_{\lambda }$ are
projectors we have successively%
\begin{align}
\left\vert \left\langle E_{\lambda }x-\frac{1}{2}x,y\right\rangle
\right\vert & \leq \left\Vert E_{\lambda }x-\frac{1}{2}x\right\Vert
\left\Vert y\right\Vert  \label{IV.e.2.14} \\
& =\left[ \left\langle E_{\lambda }x,E_{\lambda }x\right\rangle
-\left\langle E_{\lambda }x,x\right\rangle +\frac{1}{4}\left\Vert
x\right\Vert ^{2}\right] ^{1/2}\left\Vert y\right\Vert  \notag \\
& =\frac{1}{2}\left\Vert x\right\Vert \left\Vert y\right\Vert  \notag
\end{align}%
for any $x,y\in H,$ which proves the first branch in (\ref{IV.e.2.12}).

The second inequality follows from (\ref{IV.e.2.10}).

From the theory of Riemann-Stieltjes integral is well known that if $p:\left[
a,b\right] \rightarrow \mathbb{C}$ is of bounded variation and $v:\left[ a,b%
\right] \rightarrow \mathbb{R}$ is continuous and monotonic nondecreasing,
then the Riemann-Stieltjes integrals $\int_{a}^{b}p\left( t\right) dv\left(
t\right) $ and $\int_{a}^{b}\left\vert p\left( t\right) \right\vert dv\left(
t\right) $ exist and%
\begin{equation}
\left\vert \int_{a}^{b}p\left( t\right) dv\left( t\right) \right\vert \leq
\int_{a}^{b}\left\vert p\left( t\right) \right\vert dv\left( t\right) .
\label{IV.e.2.14.a}
\end{equation}

From the representation (\ref{IV.e.2.4}) we then have%
\begin{align}
& \left\vert \frac{f\left( m\right) +f\left( M\right) }{2}\cdot \left\langle
x,y\right\rangle -\left\langle f\left( A\right) x,y\right\rangle \right\vert
\label{IV.e.2.15} \\
& \leq \frac{1}{2}\int_{m-0}^{M}\left\vert \left\langle E_{\lambda
}x,y\right\rangle +\left\langle \left( E_{\lambda }-1_{H}\right)
x,y\right\rangle \right\vert df\left( \lambda \right)  \notag \\
& =\int_{m-0}^{M}\left\vert \left\langle E_{\lambda }x-\frac{1}{2}%
x,y\right\rangle \right\vert df\left( \lambda \right)  \notag
\end{align}%
for any $x,y\in H,$ from which we obtain the last branch in (\ref{IV.e.2.12}%
).
\end{proof}

We recall that a function $f:\left[ a,b\right] \rightarrow \mathbb{C}$ is
called $r-H$-H\"{o}lder continuous with fixed $r\in (0,1]$ and $H>0$ if%
\begin{equation*}
\left\vert f\left( t\right) -f\left( s\right) \right\vert \leq H\left\vert
t-s\right\vert ^{r}\text{ for any }t,s\in \left[ a,b\right] .
\end{equation*}%
We have the following result concerning this class of functions.

\begin{theorem}[Dragomir, 2010, \protect\cite{IV.SSD3}]
\label{IV.t.2.4}Let $A$ be a selfadjoint operator in the Hilbert space $H$
with the spectrum $Sp\left( A\right) \subseteq \left[ m,M\right] $ for some
real numbers $m<M$ and let $\left\{ E_{\lambda }\right\} _{\lambda }$ be its 
\textit{spectral family.} If $f:\left[ m,M\right] \rightarrow \mathbb{C}$ is 
$r-H$-H\"{o}lder continuous on $\left[ m,M\right] $, then we have the
inequality%
\begin{align}
\left\vert \frac{f\left( m\right) +f\left( M\right) }{2}\cdot \left\langle
x,y\right\rangle -\left\langle f\left( A\right) x,y\right\rangle \right\vert
& \leq \frac{1}{2^{r}}H(M-m)^{r}\dbigvee\limits_{m-0}^{M}\left( \left\langle
E_{\left( \cdot \right) }x,y\right\rangle \right)   \label{IV.e.2.16} \\
& \leq \frac{1}{2^{r}}H(M-m)^{r}\left\Vert x\right\Vert \left\Vert
y\right\Vert   \notag
\end{align}%
for any $x,y\in H.$
\end{theorem}

\begin{proof}
We start with the equality%
\begin{align}
& \frac{f\left( M\right) +f\left( m\right) }{2}\cdot \left\langle
x,y\right\rangle -\left\langle f\left( A\right) x,y\right\rangle
\label{IV.e.2.17} \\
& =\int_{m-0}^{M}\left[ \frac{f\left( M\right) +f\left( m\right) }{2}%
-f\left( \lambda \right) \right] d\left( \left\langle E_{\lambda
}x,y\right\rangle \right)  \notag
\end{align}%
for any $x,y\in H,$ that follows from the spectral representation theorem.

Since the function $\left\langle E_{\left( \cdot \right) }x,y\right\rangle $
is of bounded variation for any vector $x,y\in H,$ by applying the
inequality (\ref{IV.e.2.5}) we conclude that%
\begin{align}
& \left\vert \frac{f\left( m\right) +f\left( M\right) }{2}\cdot \left\langle
x,y\right\rangle -\left\langle f\left( A\right) x,y\right\rangle \right\vert 
\label{IV.e.2.18} \\
& \leq \max_{\lambda \in \left[ m,M\right] }\left\vert \frac{f\left(
M\right) +f\left( m\right) }{2}-f\left( \lambda \right) \right\vert
\dbigvee\limits_{m-0}^{M}\left( \left\langle E_{\left( \cdot \right)
}x,y\right\rangle \right)   \notag
\end{align}%
for any $x,y\in H.$

As $f:\left[ m,M\right] \rightarrow \mathbb{C}$ is $r-H$-H\"{o}lder
continuous on $\left[ m,M\right] $, then we have%
\begin{align}
\left\vert \frac{f\left( M\right) +f\left( m\right) }{2}-f\left( \lambda
\right) \right\vert & \leq \frac{1}{2}\left\vert f\left( M\right) -f\left(
\lambda \right) \right\vert +\frac{1}{2}\left\vert f\left( \lambda \right)
-f\left( m\right) \right\vert  \label{IV.e.2.19} \\
& \leq \frac{1}{2}H\left[ \left( M-\lambda \right) ^{r}+\left( \lambda
-m\right) ^{r}\right]  \notag
\end{align}%
for any $\lambda \in \left[ m,M\right] .$

Since, obviously, the function $g_{r}\left( \lambda \right) :=\left(
M-\lambda \right) ^{r}+\left( \lambda -m\right) ^{r},r\in \left( 0,1\right) $
has the property that%
\begin{equation*}
\max_{\lambda \in \left[ m,M\right] }g_{r}\left( \lambda \right)
=g_{r}\left( \frac{m+M}{2}\right) =2^{1-r}\left( M-m\right) ^{r},
\end{equation*}%
then by (\ref{IV.e.2.18}) we deduce the first part of (\ref{IV.e.2.16}).

The last part follows by the Total Variation Schwarz's inequality and we
omit the details.
\end{proof}

\subsection{Applications for Some Particular Functions}

It is obvious that the results established above can be applied for various
particular functions of selfadjoint operators. We restrict ourselves here to
only two examples, namely the logarithm and the power functions.

\textbf{1.} If we consider the logarithmic function $f:\left( 0,\infty
\right) \rightarrow \mathbb{R}$, $f\left( t\right) =\ln t,$ then we can
state the following result:

\begin{proposition}
\label{IV.p.4.1}Let $A$ be a selfadjoint operator in the Hilbert space $H$
with the spectrum $Sp\left( A\right) \subseteq \left[ m,M\right] $ for some
real numbers with $0<m<M$ and let $\left\{ E_{\lambda }\right\} _{\lambda }$
be its \textit{spectral family. Then for any }$x,y\in H$ we have%
\begin{align}
& \left\vert \left\langle x,y\right\rangle \ln \sqrt{mM}-\left\langle \ln
Ax,y\right\rangle \right\vert  \label{IV.e.4.1} \\
& \leq \ln \left( \frac{M}{m}\right) \times \left\{ 
\begin{array}{ll}
\begin{array}{l}
\frac{1}{2}\max_{\lambda \in \left[ m,M\right] }\left[ \left\langle
E_{\lambda }x,x\right\rangle ^{1/2}\left\langle E_{\lambda }y,y\right\rangle
^{1/2}\right. \\ 
\left. +\left\langle \left( 1_{H}-E_{\lambda }\right) x,x\right\rangle
^{1/2}\left\langle \left( 1_{H}-E_{\lambda }\right) y,y\right\rangle ^{1/2} 
\right]%
\end{array}
&  \\ 
&  \\ 
\max_{\lambda \in \left[ m,M\right] }\left\vert \left\langle E_{\lambda }x-%
\frac{1}{2}x,y\right\rangle \right\vert & 
\end{array}%
\right.  \notag \\
& \leq \frac{1}{2}\left\Vert x\right\Vert \left\Vert y\right\Vert \ln \left( 
\frac{M}{m}\right)  \notag
\end{align}%
and%
\begin{align}
& \left\vert \left\langle x,y\right\rangle \ln \sqrt{mM}-\left\langle \ln
Ax,y\right\rangle \right\vert  \label{IV.e.4.2} \\
& \leq \frac{1}{m}\times \left\{ 
\begin{array}{ll}
\begin{array}{l}
\frac{1}{2}\int_{m-0}^{M}\left[ \left\langle E_{\lambda }x,x\right\rangle
^{1/2}\left\langle E_{\lambda }y,y\right\rangle ^{1/2}\right. \\ 
\left. +\left\langle \left( 1_{H}-E_{\lambda }\right) x,x\right\rangle
^{1/2}\left\langle \left( 1_{H}-E_{\lambda }\right) y,y\right\rangle ^{1/2} 
\right] d\lambda%
\end{array}
&  \\ 
&  \\ 
\int_{m-0}^{M}\left\vert \left\langle E_{\lambda }x-\frac{1}{2}%
x,y\right\rangle \right\vert d\lambda & 
\end{array}%
\right.  \notag \\
& \leq \frac{1}{2}\left\Vert x\right\Vert \left\Vert y\right\Vert \left( 
\frac{M}{m}-1\right)  \notag
\end{align}%
and%
\begin{eqnarray}
\left\vert \left\langle x,y\right\rangle \ln \sqrt{mM}-\left\langle \ln
Ax,y\right\rangle \right\vert &\leq &\int_{m-0}^{M}\left\vert \left\langle
E_{\lambda }x-\frac{1}{2}x,y\right\rangle \right\vert \lambda ^{-1}d\lambda
\label{IV.e.4.3} \\
&\leq &\frac{1}{2}\left\Vert x\right\Vert \left\Vert y\right\Vert \ln \left( 
\frac{M}{m}\right)  \notag
\end{eqnarray}%
respectively.
\end{proposition}

The proof is obvious from Theorems \ref{IV.t.2.1}, \ref{IV.t.2.2} and \ref%
{IV.t.2.3} applied for the logarithmic function. The details are omitted.

\textbf{2.} Consider now the power function $f:\left( 0,\infty \right)
\rightarrow \mathbb{R}$, $f\left( t\right) =t^{p}$ with $p\in \left( -\infty
,0\right) \cup \left( 0,\infty \right) .$ In the case when $p\in \left(
0,1\right) ,$ the function is $p-H$-H\"{o}lder continuous with $H=1$ on any
subinterval $\left[ m,M\right] $ of $[0,\infty ).$ By making use of Theorem %
\ref{IV.t.2.4} we can state the following result:

\begin{proposition}
\label{IV.p.4.2}Let $A$ be a selfadjoint operator in the Hilbert space $H$
with the spectrum $Sp\left( A\right) \subseteq \left[ m,M\right] $ for some
real numbers with $0\leq m<M$ and let $\left\{ E_{\lambda }\right\}
_{\lambda }$ be its \textit{spectral family. Then for }$p\in \left(
0,1\right) $ we have%
\begin{align}
\left\vert \frac{m^{p}+M^{p}}{2}\cdot \left\langle x,y\right\rangle
-\left\langle A^{p}x,y\right\rangle \right\vert & \leq \frac{1}{2^{p}}%
(M-m)^{p}\dbigvee\limits_{m-0}^{M}\left( \left\langle E_{\left( \cdot
\right) }x,y\right\rangle \right)   \label{IV.e.4.4} \\
& \leq \frac{1}{2^{p}}(M-m)^{p}\left\Vert x\right\Vert \left\Vert
y\right\Vert ,  \notag
\end{align}%
for any $x,y\in H.$
\end{proposition}

The case of powers $p\geq 1$ is embodied in the following:

\begin{proposition}
\label{IV.p.4.3}Let $A$ be a selfadjoint operator in the Hilbert space $H$
with the spectrum $Sp\left( A\right) \subseteq \left[ m,M\right] $ for some
real numbers with $0\leq m<M$ and let $\left\{ E_{\lambda }\right\}
_{\lambda }$ be its \textit{spectral family. Then for }$p\geq 1$ \textit{and
for any }$x,y\in H$ we have%
\begin{align}
& \left\vert \frac{m^{p}+M^{p}}{2}\cdot \left\langle x,y\right\rangle
-\left\langle A^{p}x,y\right\rangle \right\vert  \label{IV.e.4.5} \\
& \leq \left( M^{p}-m^{p}\right) \times \left\{ 
\begin{array}{ll}
\begin{array}{l}
\frac{1}{2}\max_{\lambda \in \left[ m,M\right] }\left[ \left\langle
E_{\lambda }x,x\right\rangle ^{1/2}\left\langle E_{\lambda }y,y\right\rangle
^{1/2}\right. \\ 
\left. +\left\langle \left( 1_{H}-E_{\lambda }\right) x,x\right\rangle
^{1/2}\left\langle \left( 1_{H}-E_{\lambda }\right) y,y\right\rangle ^{1/2} 
\right]%
\end{array}
&  \\ 
&  \\ 
\max_{\lambda \in \left[ m,M\right] }\left\vert \left\langle E_{\lambda }x-%
\frac{1}{2}x,y\right\rangle \right\vert & 
\end{array}%
\right.  \notag \\
& \leq \frac{1}{2}\left\Vert x\right\Vert \left\Vert y\right\Vert \left(
M^{p}-m^{p}\right)  \notag
\end{align}%
and%
\begin{align}
& \left\vert \frac{m^{p}+M^{p}}{2}\cdot \left\langle x,y\right\rangle
-\left\langle A^{p}x,y\right\rangle \right\vert  \label{IV.e.4.6} \\
& \leq pM^{p-1}\times \left\{ 
\begin{array}{ll}
\begin{array}{l}
\frac{1}{2}\int_{m-0}^{M}\left[ \left\langle E_{\lambda }x,x\right\rangle
^{1/2}\left\langle E_{\lambda }y,y\right\rangle ^{1/2}\right. \\ 
\left. +\left\langle \left( 1_{H}-E_{\lambda }\right) x,x\right\rangle
^{1/2}\left\langle \left( 1_{H}-E_{\lambda }\right) y,y\right\rangle ^{1/2} 
\right] d\lambda%
\end{array}
&  \\ 
&  \\ 
\int_{m-0}^{M}\left\vert \left\langle E_{\lambda }x-\frac{1}{2}%
x,y\right\rangle \right\vert d\lambda & 
\end{array}%
\right.  \notag \\
& \leq \frac{1}{2}p\left\Vert x\right\Vert \left\Vert y\right\Vert M^{p-1} 
\notag
\end{align}%
and%
\begin{align}
\left\vert \frac{m^{p}+M^{p}}{2}\cdot \left\langle x,y\right\rangle
-\left\langle A^{p}x,y\right\rangle \right\vert & \leq
p\int_{m-0}^{M}\left\vert \left\langle E_{\lambda }x-\frac{1}{2}%
x,y\right\rangle \right\vert \lambda ^{p-1}d\lambda  \label{IV.e.4.7} \\
& \leq \frac{1}{2}\left\Vert x\right\Vert \left\Vert y\right\Vert \left(
M^{p}-m^{p}\right)  \notag
\end{align}%
respectively.
\end{proposition}

The proof is obvious from Theorems \ref{IV.t.2.1}, \ref{IV.t.2.2} and \ref%
{IV.t.2.3} applied for the power function $f:\left( 0,\infty \right)
\rightarrow \mathbb{R}$, $f\left( t\right) =t^{p}$ with $p\geq 1.$ The
details are omitted.

The case of negative powers is similar. The details are left to the
interested reader.

\section{Generalised Trapezoidal Inequalities}

\subsection{Some Vector Inequalities}

In the present section we are interested in providing error bounds for
approximating $\left\langle f\left( A\right) x,y\right\rangle $ with the
quantity%
\begin{equation}
\frac{1}{M-m}\left[ f\left( m\right) \left( M\left\langle x,y\right\rangle
-\left\langle Ax,y\right\rangle \right) +f\left( M\right) \left(
\left\langle Ax,y\right\rangle -m\left\langle x,y\right\rangle \right) %
\right]  \label{IV.a.e.1.4}
\end{equation}%
where $x,y\in H,$ which is a generalized trapezoid formula. Applications for
some particular functions are provided as well. The following representation
is of interest in itself and will be useful in deriving our inequalities
later as well:

\begin{lemma}[Dragomir, 2010, \protect\cite{IV.a.SSD4}]
\label{IV.a.l.2.1}Let $A$ be a selfadjoint operator in the Hilbert space $H$
with the spectrum $Sp\left( A\right) \subseteq \left[ m,M\right] $ for some
real numbers $m<M$ and let $\left\{ E_{\lambda }\right\} _{\lambda }$ be its 
\textit{spectral family.} If $f:\left[ m,M\right] \rightarrow \mathbb{C}$ is
a continuous function on $\left[ m,M\right] $, then we have the
representation%
\begin{align}
& \left\langle \left[ \frac{f\left( m\right) \left( M1_{H}-A\right) +f\left(
M\right) \left( A-m1_{H}\right) }{M-m}\right] x,y\right\rangle -\left\langle
f\left( A\right) x,y\right\rangle  \label{IV.a.e.2.1} \\
& =\int_{m-0}^{M}\left\langle E_{t}x,y\right\rangle df\left( t\right) -\frac{%
f\left( M\right) -f\left( m\right) }{M-m}\int_{m-0}^{M}\left\langle
E_{t}x,y\right\rangle dt  \notag \\
& =\int_{m-0}^{M}\left[ \left\langle E_{t}x,y\right\rangle -\frac{1}{M-m}%
\int_{m-0}^{M}\left\langle E_{s}x,y\right\rangle ds\right] df\left( t\right)
\notag
\end{align}%
for any $x,y\in H.$
\end{lemma}

\begin{proof}
Integrating by parts and utilizing the spectral representation theorem we
have%
\begin{align*}
\int_{m-0}^{M}\left\langle E_{t}x,y\right\rangle df\left( t\right) &
=f\left( M\right) \left\langle x,y\right\rangle -\int_{m-0}^{M}f\left(
t\right) d\left\langle E_{t}x,y\right\rangle \\
& =f\left( M\right) \left\langle x,y\right\rangle -\left\langle f\left(
A\right) x,y\right\rangle
\end{align*}%
and%
\begin{equation*}
\int_{m-0}^{M}\left\langle E_{t}x,y\right\rangle dt=M\left\langle
x,y\right\rangle -\left\langle Ax,y\right\rangle
\end{equation*}%
for any $x,y\in H.$

Therefore%
\begin{align*}
& \int_{m-0}^{M}\left\langle E_{t}x,y\right\rangle df\left( t\right) -\frac{%
f\left( M\right) -f\left( m\right) }{M-m}\int_{m-0}^{M}\left\langle
E_{t}x,y\right\rangle dt \\
& =f\left( M\right) \left\langle x,y\right\rangle -\left\langle f\left(
A\right) x,y\right\rangle -\frac{f\left( M\right) -f\left( m\right) }{M-m}%
\left( M\left\langle x,y\right\rangle -\left\langle Ax,y\right\rangle \right)
\\
& =\frac{1}{M-m}\left[ f\left( m\right) \left( M\left\langle
x,y\right\rangle -\left\langle Ax,y\right\rangle \right) +f\left( M\right)
\left( \left\langle Ax,y\right\rangle -m\left\langle x,y\right\rangle
\right) \right] \\
& -\left\langle f\left( A\right) x,y\right\rangle
\end{align*}%
for any $x,y\in H,$ which proves the first equality in (\ref{IV.a.e.2.1}).

The second equality is obvious.
\end{proof}

The following result provides error bounds in approximating $\left\langle
f\left( A\right) x,y\right\rangle $ by the generalized trapezoidal rule (\ref%
{IV.a.e.1.4}):

\begin{theorem}[Dragomir, 2010, \protect\cite{IV.a.SSD4}]
\label{IV.a.t.2.1}Let $A$ be a selfadjoint operator in the Hilbert space $H$
with the spectrum $Sp\left( A\right) \subseteq \left[ m,M\right] $ for some
real numbers $m<M$ and let $\left\{ E_{\lambda }\right\} _{\lambda }$ be its 
\textit{spectral family.}

1. If $f:\left[ m,M\right] \rightarrow \mathbb{C}$ is of bounded variation
on $\left[ m,M\right] $, then%
\begin{align}
& \left\vert \left\langle \left[ \frac{f\left( m\right) \left(
M1_{H}-A\right) +f\left( M\right) \left( A-m1_{H}\right) }{M-m}\right]
x,y\right\rangle -\left\langle f\left( A\right) x,y\right\rangle \right\vert
\label{IV.a.e.2.2} \\
& \leq \sup_{t\in \left[ m,M\right] }\left[ \frac{t-m}{M-m}%
\dbigvee\limits_{m-0}^{t}\left( \left\langle E_{\left( \cdot \right)
}x,y\right\rangle \right) +\frac{M-t}{M-m}\dbigvee\limits_{t}^{M}\left(
\left\langle E_{\left( \cdot \right) }x,y\right\rangle \right) \right]
\dbigvee\limits_{m}^{M}\left( f\right)  \notag \\
& \leq \dbigvee\limits_{m-0}^{M}\left( \left\langle E_{\left( \cdot \right)
}x,y\right\rangle \right) \dbigvee\limits_{m}^{M}\left( f\right) \leq
\left\Vert x\right\Vert \left\Vert y\right\Vert
\dbigvee\limits_{m}^{M}\left( f\right)  \notag
\end{align}%
for any $x,y\in H.$

2. If $f:\left[ m,M\right] \rightarrow \mathbb{C}$ is Lipschitzian with the
constant $L>0$ on $\left[ m,M\right] $, then%
\begin{align}
& \left\vert \left\langle \left[ \frac{f\left( m\right) \left(
M1_{H}-A\right) +f\left( M\right) \left( A-m1_{H}\right) }{M-m}\right]
x,y\right\rangle -\left\langle f\left( A\right) x,y\right\rangle \right\vert
\label{IV.a.e.2.2.a} \\
& \leq L\int_{m}^{M}\left[ \frac{t-m}{M-m}\dbigvee\limits_{m-0}^{t}\left(
\left\langle E_{\left( \cdot \right) }x,y\right\rangle \right) +\frac{M-t}{%
M-m}\dbigvee\limits_{t}^{M}\left( \left\langle E_{\left( \cdot \right)
}x,y\right\rangle \right) \right] dt  \notag \\
& \leq L\left( M-m\right) \dbigvee\limits_{m-0}^{M}\left( \left\langle
E_{\left( \cdot \right) }x,y\right\rangle \right) \leq L\left( M-m\right)
\left\Vert x\right\Vert \left\Vert y\right\Vert  \notag
\end{align}%
for any $x,y\in H.$

3. If $f:\left[ m,M\right] \rightarrow \mathbb{R}$ is monotonic
nondecreasing on $\left[ m,M\right] $, then%
\begin{align}
& \left\vert \left\langle \left[ \frac{f\left( m\right) \left(
M1_{H}-A\right) +f\left( M\right) \left( A-m1_{H}\right) }{M-m}\right]
x,y\right\rangle -\left\langle f\left( A\right) x,y\right\rangle \right\vert
\label{IV.a.e.2.2.b} \\
& \leq \int_{m}^{M}\left[ \frac{t-m}{M-m}\dbigvee\limits_{m-0}^{t}\left(
\left\langle E_{\left( \cdot \right) }x,y\right\rangle \right) +\frac{M-t}{%
M-m}\dbigvee\limits_{t}^{M}\left( \left\langle E_{\left( \cdot \right)
}x,y\right\rangle \right) \right] df\left( t\right)  \notag \\
& \leq \dbigvee\limits_{m-0}^{M}\left( \left\langle E_{\left( \cdot \right)
}x,y\right\rangle \right) \left[ f\left( M\right) -f\left( m\right) \right]
\leq \left\Vert x\right\Vert \left\Vert y\right\Vert \left[ f\left( M\right)
-f\left( m\right) \right]  \notag
\end{align}%
for any $x,y\in H.$
\end{theorem}

\begin{proof}
It is well known that if $p:\left[ a,b\right] \rightarrow \mathbb{C}$ is a
bounded function, $v:\left[ a,b\right] \rightarrow \mathbb{C}$ is of bounded
variation and the Riemann-Stieltjes integral $\int_{a}^{b}p\left( t\right)
dv\left( t\right) $ exists, then the following inequality holds%
\begin{equation}
\left\vert \int_{a}^{b}p\left( t\right) dv\left( t\right) \right\vert \leq
\sup_{t\in \left[ a,b\right] }\left\vert p\left( t\right) \right\vert
\dbigvee\limits_{a}^{b}\left( v\right) ,  \label{IV.a.e.2.3}
\end{equation}%
where $\dbigvee\limits_{a}^{b}\left( v\right) $ denotes the total variation
of $v$ on $\left[ a,b\right] .$

Applying this property to the equality (\ref{IV.a.e.2.1}), we have%
\begin{align}
& \left\vert \left\langle \left[ \frac{f\left( m\right) \left(
M1_{H}-A\right) +f\left( M\right) \left( A-m1_{H}\right) }{M-m}\right]
x,y\right\rangle -\left\langle f\left( A\right) x,y\right\rangle \right\vert
\label{IV.a.e.2.4} \\
& \leq \sup_{t\in \left[ m,M\right] }\left\vert \left\langle
E_{t}x,y\right\rangle -\frac{1}{M-m}\int_{m-0}^{M}\left\langle
E_{s}x,y\right\rangle ds\right\vert \dbigvee\limits_{m}^{M}\left( f\right) 
\notag
\end{align}%
for any $x,y\in H.$

Now, a simple integration by parts in the Riemann-Stieltjes integral reveals
the following equality of interest%
\begin{align}
& \left\langle E_{t}x,y\right\rangle -\frac{1}{M-m}\int_{m-0}^{M}\left%
\langle E_{s}x,y\right\rangle ds  \label{IV.a.e.2.5} \\
& =\frac{1}{M-m}\left[ \int_{m-0}^{t}\left( s-m\right) d\left\langle
E_{s}x,y\right\rangle +\int_{t}^{M}\left( s-M\right) d\left\langle
E_{s}x,y\right\rangle \right]  \notag
\end{align}%
that holds for any $t\in \left[ m,M\right] $ and for any $x,y\in H.$

Since the function $v\left( s\right) :=\left\langle E_{s}x,y\right\rangle $
is of bounded variation on $\left[ m,M\right] $ for any $x,y\in H,$ then on
applying the inequality (\ref{IV.a.e.2.3}) once more, we get%
\begin{align}
& \left\vert \left\langle E_{t}x,y\right\rangle -\frac{1}{M-m}%
\int_{m-0}^{M}\left\langle E_{s}x,y\right\rangle ds\right\vert
\label{IV.a.e.2.6} \\
& \leq \frac{1}{M-m}\left[ \left\vert \int_{m-0}^{t}\left( s-m\right)
d\left\langle E_{s}x,y\right\rangle \right\vert +\left\vert
\int_{t}^{M}\left( s-M\right) d\left\langle E_{s}x,y\right\rangle
\right\vert \right]  \notag \\
& \leq \frac{t-m}{M-m}\dbigvee\limits_{m-0}^{t}\left( \left\langle E_{\left(
\cdot \right) }x,y\right\rangle \right) +\frac{M-t}{M-m}\dbigvee%
\limits_{t}^{M}\left( \left\langle E_{\left( \cdot \right) }x,y\right\rangle
\right)  \notag
\end{align}%
that holds for any $t\in \left[ m,M\right] $ and for any $x,y\in H.$

Now, taking the supremum in (\ref{IV.a.e.2.6}) and taking into account that 
\begin{equation*}
\dbigvee\limits_{m-0}^{t}\left( \left\langle E_{\left( \cdot \right)
}x,y\right\rangle \right) ,\dbigvee\limits_{t}^{M}\left( \left\langle
E_{\left( \cdot \right) }x,y\right\rangle \right) \leq
\dbigvee\limits_{m-0}^{M}\left( \left\langle E_{\left( \cdot \right)
}x,y\right\rangle \right)
\end{equation*}%
for any $t\in \left[ m,M\right] $ and for any $x,y\in H,$ we deduce the
first and the second inequality in (\ref{IV.a.e.2.2}).

The last part of (\ref{IV.a.e.2.2}) follows by the Total Variation Schwarz's
inequality and we omit the details.

Now, recall that if $p:\left[ a,b\right] \rightarrow \mathbb{C}$ is a
Riemann integrable function and $v:\left[ a,b\right] \rightarrow \mathbb{C}$
is Lipschitzian with the constant $L>0$, i.e.,%
\begin{equation*}
\left\vert f\left( s\right) -f\left( t\right) \right\vert \leq L\left\vert
s-t\right\vert \text{ for any }t,s\in \left[ a,b\right] ,
\end{equation*}%
then the Riemann-Stieltjes integral $\int_{a}^{b}p\left( t\right) dv\left(
t\right) $ exists and the following inequality holds%
\begin{equation*}
\left\vert \int_{a}^{b}p\left( t\right) dv\left( t\right) \right\vert \leq
L\int_{a}^{b}\left\vert p\left( t\right) \right\vert dt.
\end{equation*}

Now, on applying this property of the Riemann-Stieltjes integral, we have
from the representation (\ref{IV.a.e.2.1}) that%
\begin{align}
& \left\vert \left\langle \left[ \frac{f\left( m\right) \left(
M1_{H}-A\right) +f\left( M\right) \left( A-m1_{H}\right) }{M-m}\right]
x,y\right\rangle -\left\langle f\left( A\right) x,y\right\rangle \right\vert
\label{IV.a.e.2.8} \\
& \leq L\int_{m-0}^{M}\left\vert \left\langle E_{t}x,y\right\rangle -\frac{1%
}{M-m}\int_{m-0}^{M}\left\langle E_{s}x,y\right\rangle ds\right\vert dt 
\notag
\end{align}%
for any $x,y\in H.$

Further on, by utilizing (\ref{IV.a.e.2.5}) we can state that%
\begin{align*}
& \int_{m-0}^{M}\left\vert \left\langle E_{t}x,y\right\rangle -\frac{1}{M-m}%
\int_{m-0}^{M}\left\langle E_{s}x,y\right\rangle ds\right\vert dt \\
& \leq \frac{1}{M-m}\int_{m-0}^{M}\left[ \left\vert \int_{m-0}^{t}\left(
s-m\right) d\left\langle E_{s}x,y\right\rangle \right\vert +\left\vert
\int_{t}^{M}\left( s-M\right) d\left\langle E_{s}x,y\right\rangle
\right\vert \right] dt \\
& \leq \int_{m-0}^{M}\left[ \frac{t-m}{M-m}\dbigvee\limits_{m-0}^{t}\left(
\left\langle E_{\left( \cdot \right) }x,y\right\rangle \right) +\frac{M-t}{%
M-m}\dbigvee\limits_{t}^{M}\left( \left\langle E_{\left( \cdot \right)
}x,y\right\rangle \right) \right] dt \\
& \leq \left( M-m\right) \dbigvee\limits_{m-0}^{M}\left( \left\langle
E_{\left( \cdot \right) }x,y\right\rangle \right)
\end{align*}%
for any $x,y\in H,$ which proves the desired result (\ref{IV.a.e.2.2.a}).

From the theory of Riemann-Stieltjes integral it is also well known that if $%
p:\left[ a,b\right] \rightarrow \mathbb{C}$ is of bounded variation and $v:%
\left[ a,b\right] \rightarrow \mathbb{R}$ is continuous and monotonic
nondecreasing, then the Riemann-Stieltjes integrals $\int_{a}^{b}p\left(
t\right) dv\left( t\right) $ and $\int_{a}^{b}\left\vert p\left( t\right)
\right\vert dv\left( t\right) $ exist and%
\begin{equation*}
\left\vert \int_{a}^{b}p\left( t\right) dv\left( t\right) \right\vert \leq
\int_{a}^{b}\left\vert p\left( t\right) \right\vert dv\left( t\right) .
\end{equation*}

From the representation (\ref{IV.a.e.2.1}) we then have%
\begin{align}
& \left\vert \left\langle \left[ \frac{f\left( m\right) \left(
M1_{H}-A\right) +f\left( M\right) \left( A-m1_{H}\right) }{M-m}\right]
x,y\right\rangle -\left\langle f\left( A\right) x,y\right\rangle \right\vert
\label{IV.a.e.2.9} \\
& \leq \int_{m-0}^{M}\left\vert \left\langle E_{t}x,y\right\rangle -\frac{1}{%
M-m}\int_{m-0}^{M}\left\langle E_{s}x,y\right\rangle ds\right\vert df\left(
t\right)  \notag
\end{align}%
for any $x,y\in H.$

Further on, by utilizing (\ref{IV.a.e.2.5}) we can state that%
\begin{align*}
& \int_{m-0}^{M}\left\vert \left\langle E_{t}x,y\right\rangle -\frac{1}{M-m}%
\int_{m-0}^{M}\left\langle E_{s}x,y\right\rangle ds\right\vert df\left(
t\right) \\
& \leq \frac{1}{M-m}\int_{m-0}^{M}\left[ \left\vert \int_{m-0}^{t}\left(
s-m\right) d\left\langle E_{s}x,y\right\rangle \right\vert +\left\vert
\int_{t}^{M}\left( s-M\right) d\left\langle E_{s}x,y\right\rangle
\right\vert \right] df\left( t\right) \\
& \leq \int_{m-0}^{M}\left[ \frac{t-m}{M-m}\dbigvee\limits_{m-0}^{t}\left(
\left\langle E_{\left( \cdot \right) }x,y\right\rangle \right) +\frac{M-t}{%
M-m}\dbigvee\limits_{t}^{M}\left( \left\langle E_{\left( \cdot \right)
}x,y\right\rangle \right) \right] df\left( t\right) \\
& \leq \left( f\left( M\right) -f\left( m\right) \right)
\dbigvee\limits_{m-0}^{M}\left( \left\langle E_{\left( \cdot \right)
}x,y\right\rangle \right)
\end{align*}%
for any $x,y\in H,$ which proves the desired result (\ref{IV.a.e.2.2.b}).
\end{proof}

A different approach for Lipschitzian functions is incorporated in:

\begin{theorem}[Dragomir, 2010, \protect\cite{IV.a.SSD4}]
\label{IV.a.t.2.2}Let $A$ be a selfadjoint operator in the Hilbert space $H$
with the spectrum $Sp\left( A\right) \subseteq \left[ m,M\right] $ for some
real numbers $m<M$ and let $\left\{ E_{\lambda }\right\} _{\lambda }$ be its 
\textit{spectral family. }If $f:\left[ m,M\right] \rightarrow \mathbb{C}$ is
Lipschitzian with the constant $L>0$ on $\left[ m,M\right] $, then%
\begin{align}
& \left\vert \left\langle \left[ \frac{f\left( m\right) \left(
M1_{H}-A\right) +f\left( M\right) \left( A-m1_{H}\right) }{M-m}\right]
x,y\right\rangle -\left\langle f\left( A\right) x,y\right\rangle \right\vert
\label{IV.a.e.2.10} \\
& \leq L\left\Vert y\right\Vert \int_{m-0}^{M}\left\Vert E_{t}x-\frac{1}{M-m}%
\int_{m-0}^{M}E_{s}xds\right\Vert dt\leq \frac{1}{2}L\left( M-m\right)
\left\Vert x\right\Vert \left\Vert y\right\Vert  \notag
\end{align}%
for any $x,y\in H.$
\end{theorem}

\begin{proof}
We will use the inequality (\ref{IV.a.e.2.8}) for which a different upper
bound will be provided.

By the Schwarz inequality in $H$ we have that%
\begin{align}
& \int_{m-0}^{M}\left\vert \left\langle E_{t}x,y\right\rangle -\frac{1}{M-m}%
\int_{m-0}^{M}\left\langle E_{s}x,y\right\rangle ds\right\vert dt
\label{IV.a.e.2.11} \\
& =\int_{m-0}^{M}\left\vert \left\langle \left[ E_{t}x-\frac{1}{M-m}%
\int_{m-0}^{M}E_{s}xds\right] ,y\right\rangle \right\vert dt  \notag \\
& \leq \left\Vert y\right\Vert \int_{m-0}^{M}\left\Vert E_{t}x-\frac{1}{M-m}%
\int_{m-0}^{M}E_{s}xds\right\Vert dt  \notag
\end{align}%
for any $x,y\in H.$

On utilizing the Cauchy-Buniakovski-Schwarz integral inequality we may state
that%
\begin{align}
& \int_{m-0}^{M}\left\Vert E_{t}x-\frac{1}{M-m}\int_{m-0}^{M}E_{s}xds\right%
\Vert dt  \label{IV.a.e.2.12} \\
& \leq \left( M-m\right) ^{1/2}\left( \int_{m-0}^{M}\left\Vert E_{t}x-\frac{1%
}{M-m}\int_{m-0}^{M}E_{s}xds\right\Vert ^{2}dt\right) ^{1/2}  \notag
\end{align}%
for any $x\in H.$

Observe that the following equalities of interest hold and they can be
easily proved by direct calculations%
\begin{align}
& \frac{1}{M-m}\int_{m-0}^{M}\left\Vert E_{t}x-\frac{1}{M-m}%
\int_{m-0}^{M}E_{s}xds\right\Vert ^{2}dt  \label{IV.a.e.2.13} \\
& =\frac{1}{M-m}\int_{m-0}^{M}\left\Vert E_{t}x\right\Vert ^{2}dt-\left\Vert 
\frac{1}{M-m}\int_{m-0}^{M}E_{s}xds\right\Vert ^{2}  \notag
\end{align}%
and%
\begin{align}
& \frac{1}{M-m}\int_{m-0}^{M}\left\Vert E_{t}x\right\Vert ^{2}dt-\left\Vert 
\frac{1}{M-m}\int_{m-0}^{M}E_{s}xds\right\Vert ^{2}  \label{IV.a.e.2.14} \\
& =\frac{1}{M-m}\int_{m-0}^{M}\left\langle E_{t}x-\frac{1}{M-m}%
\int_{m-0}^{M}E_{s}xds,E_{t}x-\frac{1}{2}x\right\rangle dt  \notag
\end{align}%
for any $x\in H.$

By (\ref{IV.a.e.2.12}), (\ref{IV.a.e.2.13}) and (\ref{IV.a.e.2.14}) we get

\begin{align}
& \int_{m-0}^{M}\left\Vert E_{t}x-\frac{1}{M-m}\int_{m-0}^{M}E_{s}xds\right%
\Vert dt  \label{IV.a.e.2.15} \\
& \leq \left( M-m\right) ^{1/2}\left( \int_{m-0}^{M}\left\langle E_{t}x-%
\frac{1}{M-m}\int_{m-0}^{M}E_{s}xds,E_{t}x-\frac{1}{2}x\right\rangle
dt\right) ^{1/2}  \notag
\end{align}%
for any $x\in H.$

On making use of the Schwarz inequality in $H$ we also have 
\begin{align}
& \int_{m-0}^{M}\left\langle E_{t}x-\frac{1}{M-m}%
\int_{m-0}^{M}E_{s}xds,E_{t}x-\frac{1}{2}x\right\rangle dt
\label{IV.a.e.2.16} \\
& \leq \int_{m-0}^{M}\left\Vert E_{t}x-\frac{1}{M-m}\int_{m-0}^{M}E_{s}xds%
\right\Vert \left\Vert E_{t}x-\frac{1}{2}x\right\Vert dt  \notag \\
& =\frac{1}{2}\left\Vert x\right\Vert \int_{m-0}^{M}\left\Vert E_{t}x-\frac{1%
}{M-m}\int_{m-0}^{M}E_{s}xds\right\Vert dt,  \notag
\end{align}%
where we used the fact that $E_{t}$ are projectors, and in this case we have%
\begin{align*}
\left\Vert E_{t}x-\frac{1}{2}x\right\Vert ^{2}& =\left\Vert
E_{t}x\right\Vert ^{2}-\left\langle E_{t}x,x\right\rangle +\frac{1}{4}%
\left\Vert x\right\Vert ^{2} \\
& =\left\langle E_{t}^{2}x,x\right\rangle -\left\langle
E_{t}x,x\right\rangle +\frac{1}{4}\left\Vert x\right\Vert ^{2}=\frac{1}{4}%
\left\Vert x\right\Vert ^{2}
\end{align*}%
for any $t\in \left[ m,M\right] $ for any $x\in H.$

From (\ref{IV.a.e.2.15}) and (\ref{IV.a.e.2.16}) we get%
\begin{align}
& \int_{m-0}^{M}\left\Vert E_{t}x-\frac{1}{M-m}\int_{m-0}^{M}E_{s}xds\right%
\Vert dt  \label{IV.a.e.2.17} \\
& \leq \left( M-m\right) ^{1/2}\left( \frac{1}{2}\left\Vert x\right\Vert
\int_{m-0}^{M}\left\Vert E_{t}x-\frac{1}{M-m}\int_{m-0}^{M}E_{s}xds\right%
\Vert dt\right) ^{1/2}  \notag
\end{align}%
which is clearly equivalent with the following inequality of interest in
itself%
\begin{equation}
\int_{m-0}^{M}\left\Vert E_{t}x-\frac{1}{M-m}\int_{m-0}^{M}E_{s}xds\right%
\Vert dt\leq \frac{1}{2}\left\Vert x\right\Vert \left( M-m\right)
\label{IV.a.e.2.18}
\end{equation}%
for any $x\in H.$

This proves the last part of (\ref{IV.a.e.2.10}).
\end{proof}

\subsection{Applications for Particular Functions}

It is obvious that the above results can be applied for various particular
functions. However, we will restrict here only to the power and logarithmic
functions.

\textbf{1.} Consider now the power function $f:\left( 0,\infty \right)
\rightarrow \mathbb{R}$, $f\left( t\right) =t^{p}$ with $p\neq 0.$ On
applying Theorem \ref{IV.a.t.2.2} we can state the following proposition:

\begin{proposition}
\label{IV.a.p.4.1}Let $A$ be a selfadjoint operator in the Hilbert space $H$
with the spectrum $Sp\left( A\right) \subseteq \left[ m,M\right] $ for some
real numbers $0\leq m<M$ and let $\left\{ E_{\lambda }\right\} _{\lambda }$
be its \textit{spectral family. Then }for any $x,y\in H$ we have the
inequalities%
\begin{align}
& \left\vert \left\langle \left[ \frac{m^{p}\left( M1_{H}-A\right)
+M^{p}\left( A-m1_{H}\right) }{M-m}\right] x,y\right\rangle -\left\langle
A^{p}x,y\right\rangle \right\vert  \label{IV.a.e.3.1} \\
& \leq B_{p}\left\Vert y\right\Vert \int_{m-0}^{M}\left\Vert E_{t}x-\frac{1}{%
M-m}\int_{m-0}^{M}E_{s}xds\right\Vert dt\leq \frac{1}{2}B_{p}\left(
M-m\right) \left\Vert x\right\Vert \left\Vert y\right\Vert  \notag
\end{align}%
where%
\begin{equation*}
B_{p}=p\times \left\{ 
\begin{array}{cc}
M^{p-1} & \text{if }p\geq 1 \\ 
&  \\ 
m^{p-1} & \text{if }0<p<1,m>0%
\end{array}%
\right.
\end{equation*}%
and 
\begin{equation*}
B_{p}=\left( -p\right) m^{p-1}\text{ if }p<0,m>0.
\end{equation*}
\end{proposition}

\textbf{2.} The case of logarithmic function is as follows:

\begin{proposition}
\label{IV.a.p.4.2}Let $A$ be a selfadjoint operator in the Hilbert space $H$
with the spectrum $Sp\left( A\right) \subseteq \left[ m,M\right] $ for some
real numbers $0<m<M$ and let $\left\{ E_{\lambda }\right\} _{\lambda }$ be
its \textit{spectral family. Then }for any $x,y\in H$ we have the
inequalities%
\begin{align}
& \left\vert \left\langle \left[ \frac{\left( M1_{H}-A\right) \ln m+\left(
A-m1_{H}\right) \ln M}{M-m}\right] x,y\right\rangle -\left\langle \ln
Ax,y\right\rangle \right\vert  \label{IV.a.e.3.2} \\
& \leq \frac{1}{m}\left\Vert y\right\Vert \int_{m-0}^{M}\left\Vert E_{t}x-%
\frac{1}{M-m}\int_{m-0}^{M}E_{s}xds\right\Vert dt\leq \frac{1}{2}\left( 
\frac{M}{m}-1\right) \left\Vert x\right\Vert \left\Vert y\right\Vert . 
\notag
\end{align}
\end{proposition}

\section{More Generalised Trapezoidal Inequalities}

\subsection{Other Vector Inequalities}

The following result for general continuous functions holds:

\begin{theorem}[Dragomir, 2010, \protect\cite{IV.b.SSD5}]
\label{IV.b.t.2.1}Let $A$ be a selfadjoint operator in the Hilbert space $H$
with the spectrum $Sp\left( A\right) \subseteq \left[ m,M\right] $ for some
real numbers $m<M$ and let $\left\{ E_{\lambda }\right\} _{\lambda }$ be its 
\textit{spectral family. If }$f:\left[ m,M\right] \rightarrow \mathbb{R}$ is
continuous on $\left[ m,M\right] ,$ then we have the inequalities:%
\begin{align}
& \left\vert \left\langle \left[ \frac{f\left( m\right) \left(
M1_{H}-A\right) +f\left( M\right) \left( A-m1_{H}\right) }{M-m}\right]
x,y\right\rangle -\left\langle f\left( A\right) x,y\right\rangle \right\vert 
\label{IV.b.e.2.1} \\
& \leq \left[ \max_{t\in \left[ m,M\right] }f\left( t\right) -\min_{t\in %
\left[ m,M\right] }f\left( t\right) \right] \dbigvee\limits_{m-0}^{M}\left(
\left\langle E_{\left( \cdot \right) }x,y\right\rangle \right)   \notag \\
& \leq \left[ \max_{t\in \left[ m,M\right] }f\left( t\right) -\min_{t\in %
\left[ m,M\right] }f\left( t\right) \right] \left\Vert x\right\Vert
\left\Vert y\right\Vert   \notag
\end{align}%
for any $x,y\in H.$
\end{theorem}

\begin{proof}
We observe that, by the spectral representation theorem, we have the
equality 
\begin{align}
& \left\langle \left[ \frac{f\left( m\right) \left( M1_{H}-A\right) +f\left(
M\right) \left( A-m1_{H}\right) }{M-m}\right] x,y\right\rangle -\left\langle
f\left( A\right) x,y\right\rangle  \label{IV.b.e.2.2} \\
& =\int_{m-0}^{M}\Phi _{f}\left( t\right) d\left( \left\langle
E_{t}x,y\right\rangle \right)  \notag
\end{align}%
for any $x,y\in H,$ where $\Phi _{f}:\left[ m,M\right] \rightarrow \mathbb{R}
$ is given by 
\begin{equation*}
\Phi _{f}\left( t\right) =\frac{1}{M-m}\left[ \left( M-t\right) f\left(
m\right) +\left( t-m\right) f\left( M\right) \right] -f\left( t\right) .
\end{equation*}

It is well known that if $p:\left[ a,b\right] \rightarrow \mathbb{C}$ is a
continuous function and $v:\left[ a,b\right] \rightarrow \mathbb{C}$ is of
bounded variation, then the Riemann-Stieltjes integral $\int_{a}^{b}p\left(
t\right) dv\left( t\right) $ exists and the following inequality holds%
\begin{equation}
\left\vert \int_{a}^{b}p\left( t\right) dv\left( t\right) \right\vert \leq
\sup_{t\in \left[ a,b\right] }\left\vert p\left( t\right) \right\vert
\dbigvee\limits_{a}^{b}\left( v\right) ,  \label{IV.b.e.2.3}
\end{equation}%
where $\dbigvee\limits_{a}^{b}\left( v\right) $ denotes the total variation
of $v$ on $\left[ a,b\right] .$

Now, if we denote by $\gamma :=\min_{t\in \left[ m,M\right] }f\left(
t\right) $ and by $\Gamma :=\max_{t\in \left[ m,M\right] }f\left( t\right) $
then we have%
\begin{align*}
\gamma \left( M-t\right) & \leq \left( M-t\right) f\left( m\right) \leq
\Gamma \left( M-t\right) , \\
\gamma \left( t-m\right) & \leq \left( t-m\right) f\left( M\right) \leq
\Gamma \left( t-m\right)
\end{align*}%
and%
\begin{equation*}
-\left( M-m\right) \Gamma \leq -\left( M-m\right) f\left( t\right) \leq
-\gamma \left( M-m\right)
\end{equation*}%
for any $t\in \left[ m,M\right] .$ If we add these three inequalities, then
we get%
\begin{equation*}
-\left( M-m\right) \left( \Gamma -\gamma \right) \leq \left( M-m\right) \Phi
_{f}\left( t\right) \leq \left( M-m\right) \left( \Gamma -\gamma \right)
\end{equation*}%
for any $t\in \left[ m,M\right] ,$ which shows that%
\begin{equation}
\left\vert \Phi _{f}\left( t\right) \right\vert \leq \Gamma -\gamma \text{
for any }t\in \left[ m,M\right] .  \label{IV.b.e.2.3.a}
\end{equation}

On applying the inequality (\ref{IV.b.e.2.3}) for the representation (\ref%
{IV.b.e.2.2}) we have from (\ref{IV.b.e.2.3.a}) that%
\begin{equation*}
\left\vert \int_{m-0}^{M}\Phi _{f}\left( t\right) d\left( \left\langle
E_{t}x,y\right\rangle \right) \right\vert \leq \left( \Gamma -\gamma \right)
\dbigvee\limits_{m-0}^{M}\left( \left\langle E_{\left( \cdot \right)
}x,y\right\rangle \right) 
\end{equation*}%
for any $x,y\in H,$ which proves the first part of (\ref{IV.b.e.2.1}).

The last part of (\ref{IV.b.e.2.1}) follows by the Total Variation Schwarz's
inequality and we omit the details.
\end{proof}

When the generating function is of bounded variation, we have the following
result.

\begin{theorem}[Dragomir, 2010, \protect\cite{IV.b.SSD5}]
\label{IV.b.t.2.2}Let $A$ be a selfadjoint operator in the Hilbert space $H$
with the spectrum $Sp\left( A\right) \subseteq \left[ m,M\right] $ for some
real numbers $m<M$ and let $\left\{ E_{\lambda }\right\} _{\lambda }$ be its 
\textit{spectral family. If }$f:\left[ m,M\right] \rightarrow \mathbb{C}$ is
continuous and of bounded variation on $\left[ m,M\right] ,$ then we have
the inequalities:%
\begin{align}
& \left\vert \left\langle \left[ \frac{f\left( m\right) \left(
M1_{H}-A\right) +f\left( M\right) \left( A-m1_{H}\right) }{M-m}\right]
x,y\right\rangle -\left\langle f\left( A\right) x,y\right\rangle \right\vert 
\label{IV.b.e.2.5} \\
& \leq \max_{t\in \left[ m,M\right] }\left[ \frac{M-t}{M-m}%
\dbigvee\limits_{m}^{t}\left( f\right) +\frac{t-m}{M-m}\dbigvee%
\limits_{t}^{M}\left( f\right) \right] \dbigvee\limits_{m-0}^{M}\left(
\left\langle E_{\left( \cdot \right) }x,y\right\rangle \right)   \notag \\
& \leq \dbigvee\limits_{m-0}^{M}\left( \left\langle E_{\left( \cdot \right)
}x,y\right\rangle \right) \dbigvee\limits_{m}^{M}\left( f\right) \leq
\dbigvee\limits_{m}^{M}\left( f\right) \left\Vert x\right\Vert \left\Vert
y\right\Vert   \notag
\end{align}%
for any $x,y\in H.$
\end{theorem}

\begin{proof}
First of all, observe that%
\begin{align}
\left( M-m\right) \Phi _{f}\left( t\right) & =\left( t-M\right) \left[
f\left( t\right) -f\left( m\right) \right] +\left( t-m\right) \left[ f\left(
M\right) -f\left( t\right) \right]  \label{IV.b.e.2.6} \\
& =\left( t-M\right) \int_{m}^{t}df\left( s\right) +\left( t-m\right)
\int_{t}^{M}df\left( s\right)  \notag
\end{align}%
for any $t\in \left[ m,M\right] .$

Therefore%
\begin{align}
\left\vert \Phi _{f}\left( t\right) \right\vert & \leq \frac{M-t}{M-m}%
\left\vert \int_{m}^{t}df\left( s\right) \right\vert +\frac{t-m}{M-m}%
\left\vert \int_{t}^{M}df\left( s\right) \right\vert  \label{IV.b.e.2.6.a} \\
& \leq \frac{M-t}{M-m}\dbigvee\limits_{m}^{t}\left( f\right) +\frac{t-m}{M-m}%
\dbigvee\limits_{t}^{M}\left( f\right)  \notag \\
& \leq \max \left\{ \frac{M-t}{M-m},\frac{t-m}{M-m}\right\} \left[
\dbigvee\limits_{m}^{t}\left( f\right) +\dbigvee\limits_{t}^{M}\left(
f\right) \right]  \notag \\
& =\left[ \frac{1}{2}+\frac{\left\vert t-\frac{m+M}{2}\right\vert }{M-m}%
\right] \dbigvee\limits_{m}^{M}\left( f\right)  \notag
\end{align}%
for any $t\in \left[ m,M\right] ,$ which implies that%
\begin{align}
\max_{t\in \left[ m,M\right] }\left\vert \Phi _{f}\left( t\right)
\right\vert & \leq \max_{t\in \left[ m,M\right] }\left[ \frac{M-t}{M-m}%
\dbigvee\limits_{m}^{t}\left( f\right) +\frac{t-m}{M-m}\dbigvee%
\limits_{t}^{M}\left( f\right) \right]  \label{IV.b.e.2.7} \\
& \leq \max_{t\in \left[ m,M\right] }\left[ \frac{1}{2}+\frac{\left\vert t-%
\frac{m+M}{2}\right\vert }{M-m}\right] \dbigvee\limits_{m}^{M}\left(
f\right) =\dbigvee\limits_{m}^{M}\left( f\right) .  \notag
\end{align}

On applying the inequality (\ref{IV.b.e.2.3}) for the representation (\ref%
{IV.b.e.2.2}) we have from (\ref{IV.b.e.2.7}) that%
\begin{align*}
& \left\vert \int_{m-0}^{M}\Phi _{f}\left( t\right) d\left( \left\langle
E_{t}x,y\right\rangle \right) \right\vert  \\
& \leq \max_{t\in \left[ m,M\right] }\left[ \frac{M-t}{M-m}%
\dbigvee\limits_{m}^{t}\left( f\right) +\frac{t-m}{M-m}\dbigvee%
\limits_{t}^{M}\left( f\right) \right] \dbigvee\limits_{m-0}^{M}\left(
\left\langle E_{\left( \cdot \right) }x,y\right\rangle \right)  \\
& \leq \dbigvee\limits_{m}^{M}\left( f\right)
\dbigvee\limits_{m-0}^{M}\left( \left\langle E_{\left( \cdot \right)
}x,y\right\rangle \right) 
\end{align*}%
for any $x,y\in H,$ which produces the desired result (\ref{IV.b.e.2.5}).
\end{proof}

The case of Lipschitzian functions is as follows:

\begin{theorem}[Dragomir, 2010, \protect\cite{IV.b.SSD5}]
\label{IV.b.t.2.3}Let $A$ be a selfadjoint operator in the Hilbert space $H$
with the spectrum $Sp\left( A\right) \subseteq \left[ m,M\right] $ for some
real numbers $m<M$ and let $\left\{ E_{\lambda }\right\} _{\lambda }$ be its 
\textit{spectral family. If }$f:\left[ m,M\right] \rightarrow \mathbb{C}$ is
Lipschitzian with the constant $L>0$ on $\left[ m,M\right] ,$ then we have
the inequalities:%
\begin{align}
& \left\vert \left\langle \left[ \frac{f\left( m\right) \left(
M1_{H}-A\right) +f\left( M\right) \left( A-m1_{H}\right) }{M-m}\right]
x,y\right\rangle -\left\langle f\left( A\right) x,y\right\rangle \right\vert 
\label{IV.b.e.2.8} \\
& \leq \dbigvee\limits_{m-0}^{M}\left( \left\langle E_{\left( \cdot \right)
}x,y\right\rangle \right)   \notag \\
& \times \max_{t\in \left[ m,M\right] }\left[ \frac{M-t}{M-m}\left\vert
f\left( t\right) -f\left( m\right) \right\vert +\frac{t-m}{M-m}\left\vert
f\left( M\right) -f\left( t\right) \right\vert \right]   \notag \\
& \leq \frac{1}{2}\left( M-m\right) L\dbigvee\limits_{m-0}^{M}\left(
\left\langle E_{\left( \cdot \right) }x,y\right\rangle \right) \leq \frac{1}{%
2}\left( M-m\right) L\left\Vert x\right\Vert \left\Vert y\right\Vert   \notag
\end{align}%
for any $x,y\in H.$
\end{theorem}

\begin{proof}
We have from the first part of the equality (\ref{IV.b.e.2.6}) that%
\begin{align}
\left\vert \Phi _{f}\left( t\right) \right\vert & \leq \frac{M-t}{M-m}%
\left\vert f\left( t\right) -f\left( m\right) \right\vert +\frac{t-m}{M-m}%
\left\vert f\left( M\right) -f\left( t\right) \right\vert
\label{IV.b.e.2.8.a} \\
& \leq \frac{2L}{M-m}\left( M-t\right) \left( t-m\right) \leq \frac{1}{2}%
\left( M-m\right) L  \notag
\end{align}%
for any $t\in \left[ m,M\right] ,$ which, by a similar argument to the one
from the above Theorem \ref{IV.b.t.2.2}, produces the desired result (\ref%
{IV.b.e.2.8}). The details are omitted.
\end{proof}

The following corollary holds:

\begin{corollary}[Dragomir, 2010, \protect\cite{IV.b.SSD5}]
\label{IV.b.c.2.1}Let $A$ be a selfadjoint operator in the Hilbert space $H$
with the spectrum $Sp\left( A\right) \subseteq \left[ m,M\right] $ for some
real numbers $m<M$ and let $\left\{ E_{\lambda }\right\} _{\lambda }$ be its 
\textit{spectral family. If }$l,L\in \mathbb{R}$ are such that $L>l$ and $f:%
\left[ m,M\right] \rightarrow \mathbb{R}$ is $\left( l,L\right) -$%
Lipschitzian on $\left[ m,M\right] ,$ then we have the inequalities:%
\begin{align}
& \left\vert \left\langle \left[ \frac{f\left( m\right) \left(
M1_{H}-A\right) +f\left( M\right) \left( A-m1_{H}\right) }{M-m}\right]
x,y\right\rangle -\left\langle f\left( A\right) x,y\right\rangle \right\vert 
\label{IV.b.e.2.13} \\
& \leq \frac{1}{4}\left( M-m\right) \left( L-l\right)
\dbigvee\limits_{m-0}^{M}\left( \left\langle E_{\left( \cdot \right)
}x,y\right\rangle \right) \leq \frac{1}{4}\left( M-m\right) \left(
L-l\right) \left\Vert x\right\Vert \left\Vert y\right\Vert   \notag
\end{align}%
for any $x,y\in H.$
\end{corollary}

\begin{proof}
Follows by applying the inequality (\ref{IV.b.e.2.8}) to the $\frac{1}{2}%
\left( L-l\right) $-Lipschitzian function $f-\frac{1}{2}\left( l+L\right) e,$
where $e\left( t\right) =t,$ $t\in \left[ m,M\right] .$ The details are
omitted.
\end{proof}

When the generating function is continuous convex, we can state the
following result as well:

\begin{theorem}[Dragomir, 2010, \protect\cite{IV.b.SSD5}]
\label{IV.b.t.2.4}Let $A$ be a selfadjoint operator in the Hilbert space $H$
with the spectrum $Sp\left( A\right) \subseteq \left[ m,M\right] $ for some
real numbers $m<M$ and let $\left\{ E_{\lambda }\right\} _{\lambda }$ be its 
\textit{spectral family. If }$f:\left[ m,M\right] \rightarrow \mathbb{R}$ is
continuous convex on $\left[ m,M\right] $ with finite lateral derivatives $%
f_{-}^{\prime }\left( M\right) $ and $f_{+}^{\prime }\left( m\right) ,$ then
we have the inequalities:%
\begin{align}
& \left\vert \left\langle \left[ \frac{f\left( m\right) \left(
M1_{H}-A\right) +f\left( M\right) \left( A-m1_{H}\right) }{M-m}\right]
x,y\right\rangle -\left\langle f\left( A\right) x,y\right\rangle \right\vert 
\label{IV.b.e.2.14} \\
& \leq \frac{1}{4}\left( M-m\right) \left[ f_{-}^{\prime }\left( M\right)
-f_{+}^{\prime }\left( m\right) \right] \dbigvee\limits_{m-0}^{M}\left(
\left\langle E_{\left( \cdot \right) }x,y\right\rangle \right)   \notag \\
& \leq \frac{1}{4}\left( M-m\right) \left[ f_{-}^{\prime }\left( M\right)
-f_{+}^{\prime }\left( m\right) \right] \left\Vert x\right\Vert \left\Vert
y\right\Vert   \notag
\end{align}%
for any $x,y\in H.$
\end{theorem}

\begin{proof}
By the convexity of $f$ on $\left[ m,M\right] $ we have%
\begin{equation*}
f\left( t\right) -f\left( M\right) \geq f_{-}^{\prime }\left( M\right)
\left( t-M\right)
\end{equation*}%
for any $t\in \left[ m,M\right] .$ If we multiply this inequality with $%
t-m\geq 0$ we deduce%
\begin{equation}
\left( t-m\right) f\left( t\right) -\left( t-m\right) f\left( M\right) \geq
f_{-}^{\prime }\left( M\right) \left( t-M\right) \left( t-m\right)
\label{IV.b.e.2.15}
\end{equation}%
for any $t\in \left[ m,M\right] .$

Similarly, we get%
\begin{equation}
\left( M-t\right) f\left( t\right) -\left( M-t\right) f\left( m\right) \geq
f_{+}^{\prime }\left( m\right) \left( M-t\right) \left( t-m\right)
\label{IV.b.e.2.16}
\end{equation}%
for any $t\in \left[ m,M\right] .$

Summing the above inequalities and dividing by $M-m$ we deduce the inequality%
\begin{align}
\Phi _{f}\left( t\right) & \leq \frac{\left( M-t\right) \left( t-m\right) }{%
M-m}\left[ f_{-}^{\prime }\left( M\right) -f_{+}^{\prime }\left( m\right) %
\right]  \label{IV.b.e.2.17} \\
& \leq \frac{1}{4}\left( M-m\right) \left[ f_{-}^{\prime }\left( M\right)
-f_{+}^{\prime }\left( m\right) \right]  \notag
\end{align}%
for any $t\in \left[ m,M\right] .$

By the convexity of $f$ we also have that%
\begin{align}
\frac{1}{M-m}\left[ \left( M-t\right) f\left( m\right) +\left( t-m\right)
f\left( M\right) \right] & \geq f\left( \frac{\left( M-t\right) m+\left(
t-m\right) M}{M-m}\right)  \label{IV.b.e.2.17.a} \\
& =f\left( t\right)  \notag
\end{align}%
giving that 
\begin{equation}
\Phi _{f}\left( t\right) \geq 0\text{ for any }t\in \left[ m,M\right] .
\label{IV.b.e.2.18}
\end{equation}

Utilising (\ref{IV.b.e.2.3}) for the representation (\ref{IV.b.e.2.2}) we
deduce from (\ref{IV.b.e.2.17}) and (\ref{IV.b.e.2.18}) the desired result (%
\ref{IV.b.e.2.14}).
\end{proof}

\subsection{Inequalities in the Operator Order}

The following result providing some inequalities in the operator order may
be stated:

\begin{theorem}[Dragomir, 2010, \protect\cite{IV.b.SSD5}]
\label{IV.b.t.3.1}Let $A$ be a selfadjoint operator in the Hilbert space $H$
with the spectrum $Sp\left( A\right) \subseteq \left[ m,M\right] $ for some
real numbers $m<M$.

1. \ \textit{If }$f:\left[ m,M\right] \rightarrow \mathbb{R}$ is continuous
on $\left[ m,M\right] ,$ then 
\begin{align}
& \left\vert \frac{f\left( m\right) \left( M1_{H}-A\right) +f\left( M\right)
\left( A-m1_{H}\right) }{M-m}-f\left( A\right) \right\vert
\label{IV.b.e.3.1} \\
& \leq \left[ \max_{t\in \left[ m,M\right] }f\left( t\right) -\min_{t\in %
\left[ m,M\right] }f\left( t\right) \right] 1_{H}.  \notag
\end{align}

2. \ \textit{If }$f:\left[ m,M\right] \rightarrow \mathbb{C}$ is continuous
and of bounded variation on $\left[ m,M\right] ,$ then 
\begin{align}
& \left\vert \frac{f\left( m\right) \left( M1_{H}-A\right) +f\left( M\right)
\left( A-m1_{H}\right) }{M-m}-f\left( A\right) \right\vert
\label{IV.b.e.3.2} \\
& \leq \frac{M1_{H}-A}{M-m}\dbigvee\limits_{m}^{A}\left( f\right) +\frac{%
A-m1_{H}}{M-m}\dbigvee\limits_{A}^{M}\left( f\right) \leq \left[ \frac{1}{2}+%
\frac{\left\vert A-\frac{m+M}{2}1_{H}\right\vert }{M-m}\right]
\dbigvee\limits_{m}^{M}\left( f\right) ,  \notag
\end{align}%
where $\dbigvee\limits_{m}^{A}\left( f\right) $ denotes the operator
generated by the scalar function $\left[ m,M\right] \ni t\longmapsto
\dbigvee\limits_{m}^{t}\left( f\right) \in \mathbb{R}$. The same notation
applies for $\dbigvee\limits_{A}^{M}\left( f\right) .$

3. \ \textit{If }$f:\left[ m,M\right] \rightarrow \mathbb{C}$ is
Lipschitzian with the constant $L>0$ on $\left[ m,M\right] ,$ then 
\begin{align}
& \left\vert \frac{f\left( m\right) \left( M1_{H}-A\right) +f\left( M\right)
\left( A-m1_{H}\right) }{M-m}-f\left( A\right) \right\vert
\label{IV.b.e.3.3.} \\
& \leq \frac{M1_{H}-A}{M-m}\left\vert f\left( A\right) -f\left( m\right)
1_{H}\right\vert +\frac{A-m1_{H}}{M-m}\left\vert f\left( M\right)
1_{H}-f\left( A\right) \right\vert  \notag \\
& \leq \frac{1}{2}\left( M-m\right) L1_{H}.  \notag
\end{align}

4. \textit{\ If }$f:\left[ m,M\right] \rightarrow \mathbb{R}$ is continuous
convex on $\left[ m,M\right] $ with finite lateral derivatives $%
f_{-}^{\prime }\left( M\right) $ and $f_{+}^{\prime }\left( m\right) ,$ then
we have the inequalities:%
\begin{align}
0& \leq \frac{f\left( m\right) \left( M1_{H}-A\right) +f\left( M\right)
\left( A-m1_{H}\right) }{M-m}-f\left( A\right)  \label{IV.b.e.3.4} \\
& \leq \frac{\left( M1_{H}-A\right) \left( A-m1_{H}\right) }{M-m}\left[
f_{-}^{\prime }\left( M\right) -f_{+}^{\prime }\left( m\right) \right] 
\notag \\
& \leq \frac{1}{4}\left( M-m\right) \left[ f_{-}^{\prime }\left( M\right)
-f_{+}^{\prime }\left( m\right) \right] 1_{H}.  \notag
\end{align}
\end{theorem}

\begin{proof}
Follows by applying the property (\ref{P}) to the scalar inequalities (\ref%
{IV.b.e.2.3.a}), (\ref{IV.b.e.2.6.a}), (\ref{IV.b.e.2.8.a}), (\ref%
{IV.b.e.2.17}) and (\ref{IV.b.e.2.18}). The details are omitted.
\end{proof}

The following particular case is perhaps more useful for applications:

\begin{corollary}[Dragomir, 2010, \protect\cite{IV.b.SSD5}]
\label{IV.b.c.3.1}Let $A$ be a selfadjoint operator in the Hilbert space $H$
with the spectrum $Sp\left( A\right) \subseteq \left[ m,M\right] $ for some
real numbers $m<M$\textit{. If }$l,L\in \mathbb{R}$ with $L>l$ and $f:\left[
m,M\right] \rightarrow \mathbb{R}$ is $\left( l,L\right) -$Lipschitzian on $%
\left[ m,M\right] ,$ then we have the inequalities:%
\begin{equation}
\left\vert \frac{f\left( m\right) \left( M1_{H}-A\right) +f\left( M\right)
\left( A-m1_{H}\right) }{M-m}-f\left( A\right) \right\vert \leq \frac{1}{4}%
\left( M-m\right) \left( L-l\right) 1_{H}.  \label{IV.b.e.3.5}
\end{equation}
\end{corollary}

\subsection{More Inequalities for Differentiable Functions}

The following result holds:

\begin{theorem}[Dragomir, 2010, \protect\cite{IV.b.SSD5}]
\label{IV.b.t.4.1}Let $A$ be a selfadjoint operator in the Hilbert space $H$
with the spectrum $Sp\left( A\right) \subseteq \left[ m,M\right] $ for some
real numbers $m<M$. Assume that the function $f:I\rightarrow \mathbb{C}$
with $\left[ m,M\right] \subset \mathring{I}$ (the interior of $I)$ is
differentiable on $\mathring{I}.$

1. If the derivative $f^{\prime }$ is continuous and of bounded variation on 
$\left[ m,M\right] ,$ then we have the inequality 
\begin{align}
& \left\vert \left\langle \left[ \frac{f\left( m\right) \left(
M1_{H}-A\right) +f\left( M\right) \left( A-m1_{H}\right) }{M-m}\right]
x,y\right\rangle -\left\langle f\left( A\right) x,y\right\rangle \right\vert 
\label{IV.b.e.5.1} \\
& \leq \frac{1}{4}\left( M-m\right) \bigvee\limits_{m}^{M}\left( f^{\prime
}\right) \dbigvee\limits_{m-0}^{M}\left( \left\langle E_{\left( \cdot
\right) }x,y\right\rangle \right) \leq \frac{1}{4}\left( M-m\right)
\bigvee\limits_{m}^{M}\left( f^{\prime }\right) \left\Vert x\right\Vert
\left\Vert y\right\Vert   \notag
\end{align}%
for any $x,y\in H.$

2. If the derivative $f^{\prime }$ is Lipschitzian with the constant $K>0$
on $\left[ m,M\right] ,$ then we have the inequality 
\begin{align}
& \left\vert \left\langle \left[ \frac{f\left( m\right) \left(
M1_{H}-A\right) +f\left( M\right) \left( A-m1_{H}\right) }{M-m}\right]
x,y\right\rangle -\left\langle f\left( A\right) x,y\right\rangle \right\vert 
\label{IV.b.e.5.1.1} \\
& \leq \frac{1}{8}\left( M-m\right) ^{2}K\dbigvee\limits_{m-0}^{M}\left(
\left\langle E_{\left( \cdot \right) }x,y\right\rangle \right) \leq \frac{1}{%
8}\left( M-m\right) ^{2}K\left\Vert x\right\Vert \left\Vert y\right\Vert  
\notag
\end{align}%
for any $x,y\in H.$
\end{theorem}

\begin{proof}
First of all we notice that if $f:\left[ m,M\right] \rightarrow \mathbb{C}$
is absolutely continuous on $\left[ m,M\right] $ and such that the
derivative $f^{\prime }$ is Riemann integrable on $\left[ m,M\right] ,$ then
we have the following representation in terms of the Riemann-Stieltjes
integral:%
\begin{equation}
\Phi _{f}\left( t\right) =\frac{1}{M-m}\int_{m}^{M}K\left( t,s\right)
df^{\prime }\left( s\right) ,\quad t\in \left[ m,M\right] ,
\label{IV.b.e.5.2}
\end{equation}%
where the kernel $K:\left[ m,M\right] ^{2}\rightarrow \mathbb{R}$ is given by%
\begin{equation}
K\left( t,s\right) :=\left\{ 
\begin{array}{ll}
\left( M-t\right) \left( s-m\right) & \text{if \ }m\leq s\leq t \\[5pt] 
\left( t-m\right) \left( M-s\right) & \text{if \ }t<s\leq M.%
\end{array}%
\right.  \label{IV.b.e.5.3}
\end{equation}%
Indeed, since $f^{\prime }$ is Riemann integrable on $\left[ m,M\right] ,$
it follows that the Riemann-Stieltjes integrals $\int_{m}^{t}\left(
s-m\right) df^{\prime }\left( s\right) $ and $\int_{t}^{M}\left( M-s\right)
df^{\prime }\left( s\right) $ exist for each $t\in \left[ m,M\right] .$ Now,
integrating by parts in the Riemann-Stieltjes integral, we have:%
\begin{align*}
\int_{m}^{M}K\left( t,s\right) df^{\prime }\left( s\right) & =\left(
M-t\right) \int_{m}^{t}\left( s-m\right) df^{\prime }\left( s\right) +\left(
t-m\right) \int_{t}^{M}\left( M-s\right) df^{\prime }\left( s\right) \\
& =\left( M-t\right) \left[ \left( s-m\right) f^{\prime }\left( s\right) %
\big|_{m}^{t}-\int_{m}^{t}f^{\prime }\left( s\right) ds\right] \\
& +\left( t-m\right) \left[ \left( M-s\right) f^{\prime }\left( s\right) %
\big|_{t}^{M}-\int_{t}^{M}f^{\prime }\left( s\right) ds\right] \\
& =\left( M-t\right) \left[ \left( t-m\right) f^{\prime }\left( t\right)
-\left( f\left( t\right) -f\left( m\right) \right) \right] \\
& +\left( t-m\right) \left[ -\left( M-t\right) f^{\prime }\left( t\right)
+f\left( M\right) -f\left( t\right) \right] \\
& =\left( t-m\right) \left[ f\left( M\right) -f\left( t\right) \right]
-\left( M-t\right) \left[ f\left( t\right) -f\left( m\right) \right] \\
& =\left( M-m\right) \Phi _{f}\left( t\right)
\end{align*}%
for any $t\in \left[ m,M\right] ,$ which provides the desired representation
(\ref{IV.b.e.5.2}).

Now, utilizing the representation (\ref{IV.b.e.5.2}) and the property (\ref%
{IV.b.e.2.3}), we have%
\begin{align}
& \left\vert \Phi _{f}\left( t\right) \right\vert  \label{IV.b.e.5.4} \\
& =\frac{1}{M-m}\left\vert \left( M-t\right) \int_{m}^{t}\left( s-m\right)
df^{\prime }\left( s\right) +\left( t-m\right) \int_{t}^{M}\left( M-s\right)
df^{\prime }\left( s\right) \right\vert  \notag \\
& \leq \frac{1}{M-m}\left[ \left( M-t\right) \left\vert \int_{m}^{t}\left(
s-m\right) df^{\prime }\left( s\right) \right\vert +\left( t-m\right)
\left\vert \int_{t}^{M}\left( M-s\right) df^{\prime }\left( s\right)
\right\vert \right]  \notag \\
& \leq \frac{1}{M-m}  \notag \\
& \times \left[ \left( M-t\right) \bigvee\limits_{m}^{t}\left( f^{\prime
}\right) \sup_{s\in \left[ m,t\right] }\left( s-m\right) +\left( t-m\right)
\bigvee\limits_{t}^{M}\left( f^{\prime }\right) \sup_{s\in \left[ t,M\right]
}\left( M-s\right) \right]  \notag \\
& =\frac{\left( t-m\right) \left( M-t\right) }{M-m}\left[ \bigvee%
\limits_{m}^{t}\left( f^{\prime }\right) +\bigvee\limits_{t}^{M}\left(
f^{\prime }\right) \right]  \notag \\
& =\frac{\left( t-m\right) \left( M-t\right) }{M-m}\bigvee_{m}^{M}\left(
f^{\prime }\right) \leq \frac{1}{4}\left( M-m\right) \bigvee_{m}^{M}\left(
f^{\prime }\right)  \notag
\end{align}%
for any $t\in \left[ m,M\right] .$

On making use of the representation (\ref{IV.b.e.2.2}) we deduce the desired
result (\ref{IV.b.e.5.1}).

Further, we utilize the fact that for an $L-$Lipschitzian function, $p:\left[
\alpha ,\beta \right] \rightarrow \mathbb{C}$ and a Riemann integrable
function $v:\left[ \alpha ,\beta \right] \rightarrow \mathbb{C}$, the
Riemann-Stieltjes integral $\int_{\alpha }^{\beta }p\left( s\right) dv\left(
s\right) $ exists and 
\begin{equation*}
\left\vert \int_{\alpha }^{\beta }p\left( s\right) dv\left( s\right)
\right\vert \leq L\int_{\alpha }^{\beta }\left\vert p\left( s\right)
\right\vert ds.
\end{equation*}%
Then, by utilizing (\ref{IV.b.e.5.4}) we have%
\begin{align}
& \left\vert \Phi _{f}\left( t\right) \right\vert  \label{e.5.5} \\
& \leq \frac{1}{M-m}\left[ \left( M-t\right) \left\vert \int_{m}^{t}\left(
s-m\right) df^{\prime }\left( s\right) \right\vert +\left( t-m\right)
\left\vert \int_{t}^{M}\left( M-s\right) df^{\prime }\left( s\right)
\right\vert \right]  \notag \\
& \leq \frac{K}{M-m}\left[ \left( M-t\right) \int_{m}^{t}\left( s-m\right)
ds+\left( t-m\right) \int_{t}^{M}\left( M-s\right) ds\right]  \notag \\
& =\frac{K}{M-m}\left[ \frac{\left( M-t\right) \left( t-m\right) ^{2}}{2}+%
\frac{\left( t-m\right) \left( M-t\right) ^{2}}{2}\right]  \notag \\
& =\frac{1}{2}\left( M-m\right) \left( t-m\right) \left( M-t\right) K\leq 
\frac{1}{8}\left( M-m\right) ^{2}K  \notag
\end{align}%
for any $t\in \left[ m,M\right] .$

On making use of the representation (\ref{IV.b.e.2.2}) we deduce the desired
result (\ref{IV.b.e.5.1.1}).
\end{proof}

The following inequalities in the operator order are of interest as well:

\begin{theorem}[Dragomir, 2010, \protect\cite{IV.b.SSD5}]
\label{IV.b.t.4.2}Let $A$ be a selfadjoint operator in the Hilbert space $H$
with the spectrum $Sp\left( A\right) \subseteq \left[ m,M\right] $ for some
real numbers $m<M$. Assume that the function $f:I\rightarrow \mathbb{C}$
with $\left[ m,M\right] \subset \mathring{I}$ (the interior of $I)$ is
differentiable on $\mathring{I}.$

1. If the derivative $f^{\prime }$ is continuous and of bounded variation on 
$\left[ m,M\right] ,$ then we have the inequality 
\begin{align}
& \left\vert \frac{f\left( m\right) \left( M1_{H}-A\right) +f\left( M\right)
\left( A-m1_{H}\right) }{M-m}-f\left( A\right) \right\vert
\label{IV.b.e.5.6} \\
& \leq \frac{\left( A-m1_{H}\right) \left( M1_{H}-A\right) }{M-m}%
\bigvee_{m}^{M}\left( f^{\prime }\right) \leq \frac{1}{4}\left( M-m\right)
\bigvee\limits_{m}^{M}\left( f^{\prime }\right) 1_{H}.  \notag
\end{align}

2. \ If the derivative $f^{\prime }$ is Lipschitzian with the constant $K>0$
on $\left[ m,M\right] ,$ then we have the inequality 
\begin{align}
& \left\vert \frac{f\left( m\right) \left( M1_{H}-A\right) +f\left( M\right)
\left( A-m1_{H}\right) }{M-m}-f\left( A\right) \right\vert
\label{IV.b.e.5.7} \\
& \leq \frac{1}{2}\left( M-m\right) \left( A-m1_{H}\right) \left(
M1_{H}-A\right) K\leq \frac{1}{8}\left( M-m\right) ^{2}K1_{H}.  \notag
\end{align}
\end{theorem}

\subsection{Applications for Particular Functions}

It is obvious that the above results can be applied for various particular
functions. However, we will restrict here only to the power and logarithmic
functions.

\textbf{1.} Consider now the power function $f:\left( 0,\infty \right)
\rightarrow \mathbb{R}$, $f\left( t\right) =t^{p}$ with $p\neq 0.$ On
applying Theorem \ref{IV.b.t.2.4} we can state the following proposition:

\begin{proposition}
\label{IV.b.p.4.1}Let $A$ be a selfadjoint operator in the Hilbert space $H$
with the spectrum $Sp\left( A\right) \subseteq \left[ m,M\right] $ for some
real numbers $0<m<M$ \textit{. Then }for any $x,y\in H$ we have the
inequalities%
\begin{align}
& \left\vert \left\langle \left[ \frac{m^{p}\left( M1_{H}-A\right)
+M^{p}\left( A-m1_{H}\right) }{M-m}\right] x,y\right\rangle -\left\langle
A^{p}x,y\right\rangle \right\vert  \label{IV.b.e.4.1} \\
& \leq \frac{1}{2}\left( M-m\right) \Delta _{p}\left\Vert x\right\Vert
\left\Vert y\right\Vert  \notag
\end{align}%
where%
\begin{equation*}
\Delta _{p}=p\times \left\{ 
\begin{array}{cc}
M^{p-1}-m^{p-1} & \text{if }p\in \left( -\infty ,0\right) \cup \lbrack
1,\infty ) \\ 
&  \\ 
m^{p-1}-M^{p-1} & \text{if }0<p<1.%
\end{array}%
\right.
\end{equation*}

In particular, 
\begin{align}
& \left\vert \left\langle \left[ \frac{M\left( M1_{H}-A\right) +m\left(
A-m1_{H}\right) }{mM\left( M-m\right) }\right] x,y\right\rangle
-\left\langle A^{-1}x,y\right\rangle \right\vert  \label{IV.b.e.4.1.1} \\
& \leq \frac{1}{2}\frac{\left( M-m\right) ^{2}\left( M+m\right) }{m^{2}M^{2}}%
\left\Vert x\right\Vert \left\Vert y\right\Vert  \notag
\end{align}%
for any $x,y\in H.$
\end{proposition}

The following inequalities in the operator order also hold:

\begin{proposition}
\label{IV.b.p.2.1.a}Let $A$ be a selfadjoint operator in the Hilbert space $%
H $ with the spectrum $Sp\left( A\right) \subseteq \left[ m,M\right] $ for
some real numbers $0<m<M$ \textit{. }

If $p\in \left( -\infty ,0\right) \cup \lbrack 1,\infty ),$ \textit{then}%
\begin{align}
0& \leq \frac{m^{p}\left( M1_{H}-A\right) +M^{p}\left( A-m1_{H}\right) }{M-m}%
-A^{p}  \label{IV.b.e.4.1.2} \\
& \leq p\frac{\left( M1_{H}-A\right) \left( A-m1_{H}\right) }{M-m}\left(
M^{p-1}-m^{p-1}\right)  \notag \\
& \leq \frac{1}{4}p\left( M-m\right) \left( M^{p-1}-m^{p-1}\right) 1_{H}. 
\notag
\end{align}

If $p\in \left( 0,1\right) ,$ then%
\begin{align}
0& \leq A^{p}-\frac{m^{p}\left( M1_{H}-A\right) +M^{p}\left( A-m1_{H}\right) 
}{M-m}  \label{IV.b.e.4.1.3} \\
& \leq p\frac{\left( M1_{H}-A\right) \left( A-m1_{H}\right) }{M-m}\left(
m^{p-1}-M^{p-1}\right)  \notag \\
& \leq \frac{1}{4}p\left( M-m\right) \left( m^{p-1}-M^{p-1}\right) 1_{H}. 
\notag
\end{align}

In particular, we have the inequalities%
\begin{align}
0& \leq \frac{M\left( M1_{H}-A\right) +m\left( A-m1_{H}\right) }{mM\left(
M-m\right) }-A^{-1}  \label{IV.b.e.4.1.4} \\
& \leq \frac{\left( M1_{H}-A\right) \left( A-m1_{H}\right) }{M-m}\cdot \frac{%
M^{2}-m^{2}}{m^{2}M^{2}}  \notag \\
& \leq \frac{1}{2}\frac{\left( M-m\right) ^{2}\left( M+m\right) }{m^{2}M^{2}}%
1_{H}.  \notag
\end{align}
\end{proposition}

The proof follows from (\ref{IV.b.e.3.4}) and the details are omitted.

\textbf{2.} The case of logarithmic function is as follows:

\begin{proposition}
\label{IV.b.p.4.3}Let $A$ be a selfadjoint operator in the Hilbert space $H$
with the spectrum $Sp\left( A\right) \subseteq \left[ m,M\right] $ for some
real numbers $0<m<M$\textit{. Then }for any $x,y\in H$ we have the
inequalities%
\begin{align}
& \left\vert \left\langle \left[ \frac{\left( M1_{H}-A\right) \ln m+\left(
A-m1_{H}\right) \ln M}{M-m}\right] x,y\right\rangle -\left\langle \ln
Ax,y\right\rangle \right\vert  \label{IV.b.e.4.2} \\
& \leq \frac{1}{4}\frac{\left( M-m\right) ^{2}}{mM}\left\Vert x\right\Vert
\left\Vert y\right\Vert .  \notag
\end{align}%
We also have the following inequality in the operator order%
\begin{align}
0& \leq \ln A-\frac{\left( M1_{H}-A\right) \ln m+\left( A-m1_{H}\right) \ln M%
}{M-m}  \label{IV.b.e.4.3} \\
& \leq \frac{\left( M1_{H}-A\right) \left( A-m1_{H}\right) }{Mm}\leq \frac{1%
}{4}\frac{\left( M-m\right) ^{2}}{mM}1_{H}.  \notag
\end{align}
\end{proposition}

\begin{remark}
\label{IV.b.r.4.1}Similar results can be obtained if ones uses the
inequalities from Theorem \ref{IV.b.t.4.1} and \ref{IV.b.t.4.2}. However the
details are left to the interested reader.
\end{remark}

\section{Product Inequalities}

\subsection{Some Vector Inequalities}

In this section we investigate the quantity 
\begin{equation*}
\left\vert \left\langle \left[ f\left( M\right) 1_{H}-f\left( A\right) %
\right] \left[ f\left( A\right) -f\left( m\right) 1_{H}\right]
x,y\right\rangle \right\vert
\end{equation*}%
where $x,y$ are vectors in the Hilbert space $H$ and $A$ is a selfadjoint
operator with $Sp\left( A\right) \subseteq \left[ m,M\right] ,$ and provide
different bounds for some classes of continuous functions $f:\left[ m,M%
\right] \rightarrow \mathbb{C}$. Applications for some particular cases
including the power and logarithmic functions are provided as well.

The following representation in terms of the spectral family is of interest
in itself:

\begin{lemma}[Dragomir, 2010, \protect\cite{IV.c.SSD7}]
\label{IV.c.l.2.1}Let $A$ be a selfadjoint operator in the Hilbert space $H$
with the spectrum $Sp\left( A\right) \subseteq \left[ m,M\right] $ for some
real numbers $m<M$ and let $\left\{ E_{\lambda }\right\} _{\lambda }$ be its 
\textit{spectral family.} If $f:\left[ m,M\right] \rightarrow \mathbb{C}$ is
a continuous function on $\left[ m,M\right] $ with $f\left( M\right) \neq
f\left( m\right) $ , then we have the representation%
\begin{align}
& \frac{1}{\left[ f\left( M\right) -f\left( m\right) \right] ^{2}}\left[
f\left( M\right) 1_{H}-f\left( A\right) \right] \left[ f\left( A\right)
-f\left( m\right) 1_{H}\right]  \label{IV.c.e.2.1} \\
& =\frac{1}{f\left( M\right) -f\left( m\right) }  \notag \\
& \times \int_{m-0}^{M}\left( E_{t}-\frac{1}{f\left( M\right) -f\left(
m\right) }\int_{m-0}^{M}E_{s}df\left( s\right) \right) \left( E_{t}-\frac{1}{%
2}1_{H}\right) df\left( t\right) .  \notag
\end{align}
\end{lemma}

\begin{proof}
We observe that, 
\begin{align}
& \frac{1}{f\left( M\right) -f\left( m\right) }\int_{m-0}^{M}\left( E_{t}-%
\frac{1}{f\left( M\right) -f\left( m\right) }\int_{m-0}^{M}E_{s}df\left(
s\right) \right)  \label{IV.c.e.2.2} \\
& \times \left( E_{t}-\frac{1}{2}1_{H}\right) df\left( t\right)  \notag \\
& =\frac{1}{f\left( M\right) -f\left( m\right) }\int_{m-0}^{M}E_{t}^{2}df%
\left( t\right)  \notag \\
& -\frac{1}{f\left( M\right) -f\left( m\right) }\int_{m-0}^{M}E_{s}df\left(
s\right) \frac{1}{f\left( M\right) -f\left( m\right) }\int_{m-0}^{M}E_{t}df%
\left( t\right)  \notag \\
& -\frac{1}{2}\int_{m-0}^{M}E_{t}df\left( t\right) +\frac{1}{2}%
\int_{m-0}^{M}E_{s}df\left( s\right)  \notag \\
& =\frac{1}{f\left( M\right) -f\left( m\right) }\int_{m-0}^{M}E_{t}^{2}df%
\left( t\right) -\left[ \frac{1}{f\left( M\right) -f\left( m\right) }%
\int_{m-0}^{M}E_{t}df\left( t\right) \right] ^{2}  \notag
\end{align}%
which is an equality of interest in itself.

Since $E_{t}$ are projections, we have $E_{t}^{2}=E_{t}$ for any $t\in \left[
m,M\right] $ and then we can write that 
\begin{align}
& \frac{1}{f\left( M\right) -f\left( m\right) }\int_{m-0}^{M}E_{t}^{2}df%
\left( t\right) -\left[ \frac{1}{f\left( M\right) -f\left( m\right) }%
\int_{m-0}^{M}E_{t}df\left( t\right) \right] ^{2}  \label{IV.c.e.2.3} \\
& =\frac{1}{f\left( M\right) -f\left( m\right) }\int_{m-0}^{M}E_{t}df\left(
t\right) -\left[ \frac{1}{f\left( M\right) -f\left( m\right) }%
\int_{m-0}^{M}E_{t}df\left( t\right) \right] ^{2}  \notag \\
& =\frac{1}{f\left( M\right) -f\left( m\right) }\int_{m-0}^{M}E_{t}df\left(
t\right) \left[ 1_{H}-\frac{1}{f\left( M\right) -f\left( m\right) }%
\int_{m-0}^{M}E_{t}df\left( t\right) \right] .  \notag
\end{align}

Integrating by parts in the Riemann-Stieltjes integral and utilizing the
spectral representation theorem we have%
\begin{equation*}
\int_{m-0}^{M}E_{t}df\left( t\right) =f\left( M\right) 1_{H}-f\left( A\right)
\end{equation*}%
and%
\begin{equation*}
1_{H}-\frac{1}{f\left( M\right) -f\left( m\right) }\int_{m-0}^{M}E_{t}df%
\left( t\right) =\frac{f\left( A\right) -f\left( m\right) 1_{H}}{f\left(
M\right) -f\left( m\right) },
\end{equation*}%
which together with (\ref{IV.c.e.2.3}) and (\ref{IV.c.e.2.2}) produce the
desired result (\ref{IV.c.e.2.1}).
\end{proof}

The following vector version may be stated as well:

\begin{corollary}[Dragomir, 2010, \protect\cite{IV.c.SSD7}]
\label{IV.c.c.2.1}With the assumptions of Lemma \ref{IV.c.l.2.1} we have the
equality%
\begin{align}
& \left\langle \left[ f\left( M\right) 1_{H}-f\left( A\right) \right] \left[
f\left( A\right) -f\left( m\right) 1_{H}\right] x,y\right\rangle
\label{IV.c.e.2.4} \\
& =\left[ f\left( M\right) -f\left( m\right) \right]  \notag \\
& \times \int_{m-0}^{M}\left\langle \left( E_{t}-\frac{1}{f\left( M\right)
-f\left( m\right) }\int_{m-0}^{M}E_{s}df\left( s\right) \right) x,\left(
E_{t}-\frac{1}{2}1_{H}\right) y\right\rangle df\left( t\right) ,  \notag
\end{align}%
for any $x,y\in \left[ m,M\right] .$
\end{corollary}

The following result that provides some bounds for continuous functions of
bounded variation may be stated as well:

\begin{theorem}[Dragomir, 2010, \protect\cite{IV.c.SSD7}]
\label{IV.c.t.2.1}Let $A$ be a selfadjoint operator in the Hilbert space $H$
with the spectrum $Sp\left( A\right) \subseteq \left[ m,M\right] $ for some
real numbers $m<M$ and let $\left\{ E_{\lambda }\right\} _{\lambda }$ be its 
\textit{spectral family.} If $f:\left[ m,M\right] \rightarrow \mathbb{C}$ is
a continuous function of bounded variation on $\left[ m,M\right] $ with $%
f\left( M\right) \neq f\left( m\right) $, then we have the inequality%
\begin{align}
& \left\vert \left\langle \left[ f\left( M\right) 1_{H}-f\left( A\right) %
\right] \left[ f\left( A\right) -f\left( m\right) 1_{H}\right]
x,y\right\rangle \right\vert  \label{IV.c.e.2.5} \\
& \leq \frac{1}{2}\left\Vert y\right\Vert \left\vert f\left( M\right)
-f\left( m\right) \right\vert \dbigvee\limits_{m}^{M}\left( f\right)  \notag
\\
& \times \sup_{t\in \left[ m,M\right] }\left\Vert E_{t}x-\frac{1}{f\left(
M\right) -f\left( m\right) }\int_{m-0}^{M}E_{s}df\left( s\right) \right\Vert
\leq \frac{1}{2}\left\Vert x\right\Vert \left\Vert y\right\Vert \left[
\dbigvee\limits_{m}^{M}\left( f\right) \right] ^{2},  \notag
\end{align}%
for any $x,y\in H.$
\end{theorem}

\begin{proof}
It is well known that if $p:\left[ a,b\right] \rightarrow \mathbb{C}$ is a
bounded function, $v:\left[ a,b\right] \rightarrow \mathbb{C}$ is of bounded
variation and the Riemann-Stieltjes integral $\int_{a}^{b}p\left( t\right)
dv\left( t\right) $ exists, then the following inequality holds%
\begin{equation}
\left\vert \int_{a}^{b}p\left( t\right) dv\left( t\right) \right\vert \leq
\sup_{t\in \left[ a,b\right] }\left\vert p\left( t\right) \right\vert
\dbigvee\limits_{a}^{b}\left( v\right) ,  \label{IV.c.e.2.6}
\end{equation}%
where $\dbigvee\limits_{a}^{b}\left( v\right) $ denotes the total variation
of $v$ on $\left[ a,b\right] .$

Utilising this property and the representation (\ref{IV.c.e.2.4}) we have by
the Schwarz inequality in Hilbert space $H$ that%
\begin{align}
& \left\vert \left\langle \left[ f\left( M\right) 1_{H}-f\left( A\right) %
\right] \left[ f\left( A\right) -f\left( m\right) 1_{H}\right]
x,y\right\rangle \right\vert  \label{IV.c.e.2.7} \\
& \leq \left\vert f\left( M\right) -f\left( m\right) \right\vert
\dbigvee\limits_{m}^{M}\left( f\right)  \notag \\
& \times \sup_{t\in \left[ m,M\right] }\left\vert \left\langle \left( E_{t}-%
\frac{1}{f\left( M\right) -f\left( m\right) }\int_{m-0}^{M}E_{s}df\left(
s\right) \right) x,\left( E_{t}-\frac{1}{2}1_{H}\right) y\right\rangle
\right\vert  \notag \\
& \leq \left\vert f\left( M\right) -f\left( m\right) \right\vert
\dbigvee\limits_{m}^{M}\left( f\right)  \notag \\
& \times \sup_{t\in \left[ m,M\right] }\left[ \left\Vert E_{t}x-\frac{1}{%
f\left( M\right) -f\left( m\right) }\int_{m-0}^{M}E_{s}xdf\left( s\right)
\right\Vert \left\Vert E_{t}y-\frac{1}{2}y\right\Vert \right]  \notag
\end{align}%
for any $x,y\in \left[ m,M\right] .$

Since $E_{t}$ are projections, and in this case we have%
\begin{eqnarray*}
\left\Vert E_{t}y-\frac{1}{2}y\right\Vert ^{2} &=&\left\Vert
E_{t}y\right\Vert ^{2}-\left\langle E_{t}y,y\right\rangle +\frac{1}{4}%
\left\Vert y\right\Vert ^{2} \\
&=&\left\langle E_{t}^{2}y,y\right\rangle -\left\langle
E_{t}y,y\right\rangle +\frac{1}{4}\left\Vert y\right\Vert ^{2}=\frac{1}{4}%
\left\Vert y\right\Vert ^{2},
\end{eqnarray*}%
then from (\ref{IV.c.e.2.7}) we deduce the first part of (\ref{IV.c.e.2.5}).

Now, by the same property (\ref{IV.c.e.2.6}) for vector valued functions $p$
with values in Hilbert spaces, we also have that%
\begin{align}
& \left\Vert \left[ f\left( M\right) -f\left( m\right) \right]
E_{t}x-\int_{m-0}^{M}E_{s}xdf\left( s\right) \right\Vert  \label{IV.c.e.2.8}
\\
& =\left\Vert \int_{m-0}^{M}\left( E_{t}x-E_{s}x\right) df\left( s\right)
\right\Vert \leq \dbigvee\limits_{m}^{M}\left( f\right) \sup_{s\in \left[ m,M%
\right] }\left\Vert E_{t}x-E_{s}x\right\Vert  \notag
\end{align}%
for any $t\in \left[ m,M\right] $ and $x\in H.$

Since $0\leq E_{t}\leq 1_{H}$ in the operator order, then $-1_{H}\leq
E_{t}-E_{s}\leq 1$ which gives that $-\left\Vert x\right\Vert ^{2}\leq
\left\langle \left( E_{t}-E_{s}\right) x,x\right\rangle \leq \left\Vert
x\right\Vert ^{2}$, i.e., $\left\vert \left\langle \left( E_{t}-E_{s}\right)
x,x\right\rangle \right\vert \leq \left\Vert x\right\Vert ^{2}$ for any $%
x\in H,$ which implies that $\left\Vert E_{t}-E_{s}\right\Vert \leq 1$ for
any $t,s\in \left[ m,M\right] .$ Therefore $\sup_{s\in \left[ m,M\right]
}\left\Vert E_{t}x-E_{s}x\right\Vert \leq \left\Vert x\right\Vert $ which
together with (\ref{IV.c.e.2.8}) prove the last part of (\ref{IV.c.e.2.5}).
\end{proof}

The case of Lipschitzian functions is as follows:

\begin{theorem}[Dragomir, 2010, \protect\cite{IV.c.SSD7}]
\label{IV.c.t.2.2}Let $A$ be a selfadjoint operator in the Hilbert space $H$
with the spectrum $Sp\left( A\right) \subseteq \left[ m,M\right] $ for some
real numbers $m<M$ and let $\left\{ E_{\lambda }\right\} _{\lambda }$ be its 
\textit{spectral family.} If $f:\left[ m,M\right] \rightarrow \mathbb{C}$ is
a Lipschitzian function with the constant $L>0$ on $\left[ m,M\right] $ and
with $f\left( M\right) \neq f\left( m\right) $, then we have the inequality%
\begin{align}
& \left\vert \left\langle \left[ f\left( M\right) 1_{H}-f\left( A\right) %
\right] \left[ f\left( A\right) -f\left( m\right) 1_{H}\right]
x,y\right\rangle \right\vert  \label{IV.c.e.2.10} \\
& \leq \frac{1}{2}L\left\Vert y\right\Vert \left\vert f\left( M\right)
-f\left( m\right) \right\vert  \notag \\
& \times \int_{m-0}^{M}\left\Vert E_{t}x-\frac{1}{f\left( M\right) -f\left(
m\right) }\int_{m-0}^{M}E_{s}xdf\left( s\right) \right\Vert dt  \notag \\
& \leq \frac{1}{2}L^{2}\left\Vert y\right\Vert
\int_{m-0}^{M}\int_{m-0}^{M}\left\Vert E_{t}x-E_{s}x\right\Vert dsdt  \notag
\\
& \leq \frac{\sqrt{2}}{2}L^{2}\left\Vert y\right\Vert \left( M-m\right)
\left\langle Ax-mx,Mx-Ax\right\rangle ^{1/2}\leq \frac{\sqrt{2}}{4}%
L^{2}\left\Vert y\right\Vert \left\Vert x\right\Vert \left( M-m\right) ^{2} 
\notag
\end{align}%
for any $x,y\in H.$
\end{theorem}

\begin{proof}
Recall that if $p:\left[ a,b\right] \rightarrow \mathbb{C}$ is a Riemann
integrable function and $v:\left[ a,b\right] \rightarrow \mathbb{C}$ is
Lipschitzian with the constant $L>0$, i.e.,%
\begin{equation*}
\left\vert f\left( s\right) -f\left( t\right) \right\vert \leq L\left\vert
s-t\right\vert \text{ for any }t,s\in \left[ a,b\right] ,
\end{equation*}%
then the Riemann-Stieltjes integral $\int_{a}^{b}p\left( t\right) dv\left(
t\right) $ exists and the following inequality holds%
\begin{equation}
\left\vert \int_{a}^{b}p\left( t\right) dv\left( t\right) \right\vert \leq
L\int_{a}^{b}\left\vert p\left( t\right) \right\vert dt.
\label{IV.c.e.2.10.a}
\end{equation}

Now, on applying this property of the Riemann-Stieltjes integral, then we
have from the representation (\ref{IV.c.e.2.4}) that%
\begin{align}
& \left\vert \left\langle \left[ f\left( M\right) 1_{H}-f\left( A\right) %
\right] \left[ f\left( A\right) -f\left( m\right) 1_{H}\right]
x,y\right\rangle \right\vert  \label{IV.c.e.2.11} \\
& \leq \left\vert f\left( M\right) -f\left( m\right) \right\vert  \notag \\
& \times \int_{m-0}^{M}\left\vert \left\langle \left( E_{t}-\frac{1}{f\left(
M\right) -f\left( m\right) }\int_{m-0}^{M}E_{s}df\left( s\right) \right)
x,\left( E_{t}-\frac{1}{2}1_{H}\right) y\right\rangle \right\vert df\left(
t\right) ,  \notag \\
& \leq L\left\vert f\left( M\right) -f\left( m\right) \right\vert  \notag \\
& \times \int_{m-0}^{M}\left\Vert E_{t}x-\frac{1}{f\left( M\right) -f\left(
m\right) }\int_{m-0}^{M}E_{s}xdf\left( s\right) \right\Vert \left\Vert
E_{t}y-\frac{1}{2}y\right\Vert dt  \notag \\
& =\frac{1}{2}L\left\Vert y\right\Vert \left\vert f\left( M\right) -f\left(
m\right) \right\vert  \notag \\
& \times \int_{m-0}^{M}\left\Vert E_{t}x-\frac{1}{f\left( M\right) -f\left(
m\right) }\int_{m-0}^{M}E_{s}xdf\left( s\right) \right\Vert dt  \notag
\end{align}%
for any $x,y\in H$ and the first inequality in (\ref{IV.c.e.2.10}) is proved.

Further, observe that%
\begin{align}
& \left\vert f\left( M\right) -f\left( m\right) \right\vert
\int_{m-0}^{M}\left\Vert E_{t}x-\frac{1}{f\left( M\right) -f\left( m\right) }%
\int_{m-0}^{M}E_{s}xdf\left( s\right) \right\Vert dt  \label{IV.c.e.2.12} \\
& =\int_{m-0}^{M}\left\Vert \left[ f\left( M\right) -f\left( m\right) \right]
E_{t}x-\int_{m-0}^{M}E_{s}xdf\left( s\right) \right\Vert dt  \notag \\
& =\int_{m-0}^{M}\left\Vert \int_{m-0}^{M}\left( E_{t}x-E_{s}x\right)
df\left( s\right) \right\Vert dt  \notag
\end{align}%
for any $x\in H.$

If we use the vector valued version of the property (\ref{IV.c.e.2.10.a}),
then we have%
\begin{equation}
\int_{m-0}^{M}\left\Vert \int_{m-0}^{M}\left( E_{t}x-E_{s}x\right) df\left(
s\right) \right\Vert dt\leq L\int_{m-0}^{M}\int_{m-0}^{M}\left\Vert
E_{t}x-E_{s}x\right\Vert dsdt  \label{IV.c.e.2.13}
\end{equation}%
for any $x\in H$ and the second part of (\ref{IV.c.e.2.10}) is proved.

Further on, by applying the double integral version of the
Cauchy-Buniakowski-Schwarz inequality we have%
\begin{align}
& \int_{m-0}^{M}\int_{m-0}^{M}\left\Vert E_{t}x-E_{s}x\right\Vert dsdt
\label{IV.c.e.2.14} \\
& \leq \left( M-m\right) \left( \int_{m-0}^{M}\int_{m-0}^{M}\left\Vert
E_{t}x-E_{s}x\right\Vert ^{2}dsdt\right) ^{1/2}  \notag
\end{align}%
for any $x\in H.$

Now, by utilizing the fact that $E_{s\text{ }}$are projections for each $%
s\in \left[ m,M\right] $, then we have%
\begin{align}
& \int_{m-0}^{M}\int_{m-0}^{M}\left\Vert E_{t}x-E_{s}x\right\Vert ^{2}dsdt
\label{IV.c.e.2.15} \\
& =2\left[ \left( M-m\right) \int_{m-0}^{M}\left\Vert E_{t}x\right\Vert
^{2}dt-\left\Vert \int_{m-0}^{M}E_{t}xdt\right\Vert ^{2}\right]  \notag \\
& =2\left[ \left( M-m\right) \int_{m-0}^{M}\left\langle
E_{t}x,x\right\rangle dt-\left\Vert \int_{m-0}^{M}E_{t}xdt\right\Vert ^{2}%
\right]  \notag
\end{align}%
for any $x\in H.$

If we integrate by parts and use the spectral representation theorem, then
we get%
\begin{equation*}
\int_{m-0}^{M}\left\langle E_{t}x,x\right\rangle dt=\left\langle
Mx-Ax,x\right\rangle \text{ and }\int_{m-0}^{M}E_{t}xdt=Mx-Ax
\end{equation*}%
and by (\ref{IV.c.e.2.15}) we then obtain the following equality of interest 
\begin{equation}
\int_{m-0}^{M}\int_{m-0}^{M}\left\Vert E_{t}x-E_{s}x\right\Vert
^{2}dsdt=2\left\langle Ax-mx,Mx-Ax\right\rangle  \label{IV.c.e.2.16}
\end{equation}%
for any $x\in H.$

On making use of (\ref{IV.c.e.2.16}) and (\ref{IV.c.e.2.14}) we then deduce
the third part of (\ref{IV.c.e.2.10}).

Finally, by utilizing the elementary inequality in inner product spaces%
\begin{equation}
\func{Re}\left\langle a,b\right\rangle \leq \frac{1}{4}\left\Vert
a+b\right\Vert ^{2},\text{ }a,b\in H,  \label{IV.c.e.2.16.a}
\end{equation}%
we also have that%
\begin{equation*}
\left\langle Ax-mx,Mx-Ax\right\rangle \leq \frac{1}{4}\left( M-m\right)
^{2}\left\Vert x\right\Vert ^{2}
\end{equation*}%
for any $x\in H,$ which proves the last inequality in (\ref{IV.c.e.2.10}).
\end{proof}

The case of nondecreasing monotonic functions is as follows:

\begin{theorem}[Dragomir, 2010, \protect\cite{IV.c.SSD7}]
\label{IV.c.t.2.3}Let $A$ be a selfadjoint operator in the Hilbert space $H$
with the spectrum $Sp\left( A\right) \subseteq \left[ m,M\right] $ for some
real numbers $m<M$ and let $\left\{ E_{\lambda }\right\} _{\lambda }$ be its 
\textit{spectral family.} If $f:\left[ m,M\right] \rightarrow \mathbb{R}$ is
a monotonic nondecreasing function on $\left[ m,M\right] $, then we have the
inequality%
\begin{align}
& \left\vert \left\langle \left[ f\left( M\right) 1_{H}-f\left( A\right) %
\right] \left[ f\left( A\right) -f\left( m\right) 1_{H}\right]
x,y\right\rangle \right\vert  \label{IV.c.e.2.17} \\
& \leq \frac{1}{2}\left\Vert y\right\Vert \left[ f\left( M\right) -f\left(
m\right) \right]  \notag \\
& \times \int_{m-0}^{M}\left\Vert E_{t}x-\frac{1}{f\left( M\right) -f\left(
m\right) }\int_{m-0}^{M}E_{s}xdf\left( s\right) \right\Vert df\left( t\right)
\notag \\
& \leq \frac{1}{2}\left\Vert y\right\Vert \left[ f\left( M\right) -f\left(
m\right) \right]  \notag \\
& \times \left\langle \left[ f\left( M\right) 1_{H}-f\left( A\right) \right] %
\left[ f\left( A\right) -f\left( m\right) 1_{H}\right] x,x\right\rangle
^{1/2}  \notag \\
& \leq \frac{1}{4}\left\Vert y\right\Vert \left\Vert x\right\Vert \left[
f\left( M\right) -f\left( m\right) \right] ^{2}  \notag
\end{align}%
for any $x,y\in H.$
\end{theorem}

\begin{proof}
From the theory of Riemann-Stieltjes integral it is also well known that if $%
p:\left[ a,b\right] \rightarrow \mathbb{C}$ is of bounded variation and $v:%
\left[ a,b\right] \rightarrow \mathbb{R}$ is continuous and monotonic
nondecreasing, then the Riemann-Stieltjes integrals $\int_{a}^{b}p\left(
t\right) dv\left( t\right) $ and $\int_{a}^{b}\left\vert p\left( t\right)
\right\vert dv\left( t\right) $ exist and%
\begin{equation*}
\left\vert \int_{a}^{b}p\left( t\right) dv\left( t\right) \right\vert \leq
\int_{a}^{b}\left\vert p\left( t\right) \right\vert dv\left( t\right) .
\end{equation*}

Now, on applying this property of the Riemann-Stieltjes integral, we have
from the representation (\ref{IV.c.e.2.4}) that%
\begin{align}
& \left\vert \left\langle \left[ f\left( M\right) 1_{H}-f\left( A\right) %
\right] \left[ f\left( A\right) -f\left( m\right) 1_{H}\right]
x,y\right\rangle \right\vert  \label{IV.c.e.2.18} \\
& \leq \left[ f\left( M\right) -f\left( m\right) \right]  \notag \\
& \times \int_{m-0}^{M}\left\vert \left\langle \left( E_{t}-\frac{1}{f\left(
M\right) -f\left( m\right) }\int_{m-0}^{M}E_{s}df\left( s\right) \right)
x,\left( E_{t}-\frac{1}{2}1_{H}\right) y\right\rangle \right\vert df\left(
t\right) ,  \notag \\
& \leq \left[ f\left( M\right) -f\left( m\right) \right]  \notag \\
& \times \int_{m-0}^{M}\left\Vert E_{t}x-\frac{1}{f\left( M\right) -f\left(
m\right) }\int_{m-0}^{M}E_{s}xdf\left( s\right) \right\Vert \left\Vert
E_{t}y-\frac{1}{2}y\right\Vert df\left( t\right)  \notag \\
& =\frac{1}{2}\left\Vert y\right\Vert \left[ f\left( M\right) -f\left(
m\right) \right]  \notag \\
& \times \int_{m-0}^{M}\left\Vert E_{t}x-\frac{1}{f\left( M\right) -f\left(
m\right) }\int_{m-0}^{M}E_{s}xdf\left( s\right) \right\Vert df\left( t\right)
\notag
\end{align}%
for any $x,y\in H,$ which proves the first inequality in (\ref{IV.c.e.2.17}).

On utilizing the Cauchy-Buniakowski-Schwarz type inequality for the
Riemann-Stieltjes integral of monotonic nondecreasing integrators, we have%
\begin{align}
& \int_{m-0}^{M}\left\Vert E_{t}x-\frac{1}{f\left( M\right) -f\left(
m\right) }\int_{m-0}^{M}E_{s}xdf\left( s\right) \right\Vert df\left( t\right)
\label{IV.c.e.2.19} \\
& \leq \left[ \int_{m-0}^{M}df\left( t\right) \right] ^{1/2}  \notag \\
& \times \left[ \int_{m-0}^{M}\left\Vert E_{t}x-\frac{1}{f\left( M\right)
-f\left( m\right) }\int_{m-0}^{M}E_{s}xdf\left( s\right) \right\Vert
^{2}df\left( t\right) \right] ^{1/2}  \notag
\end{align}%
for any $x,y\in H.$

Observe that%
\begin{align}
& \int_{m-0}^{M}\left\Vert E_{t}x-\frac{1}{f\left( M\right) -f\left(
m\right) }\int_{m-0}^{M}E_{s}xdf\left( s\right) \right\Vert ^{2}df\left(
t\right)  \label{IV.c.e.2.20} \\
& =\int_{m-0}^{M}\left[ \left\Vert E_{t}x\right\Vert ^{2}-2\func{Re}%
\left\langle E_{t}x,\frac{1}{f\left( M\right) -f\left( m\right) }%
\int_{m-0}^{M}E_{s}xdf\left( s\right) \right\rangle \right.  \notag \\
& \left. +\left\Vert \frac{1}{f\left( M\right) -f\left( m\right) }%
\int_{m-0}^{M}E_{s}xdf\left( s\right) \right\Vert ^{2}\right] df\left(
t\right)  \notag \\
& =\left[ f\left( M\right) -f\left( m\right) \right] \left[ \frac{1}{f\left(
M\right) -f\left( m\right) }\int_{m-0}^{M}\left\Vert E_{t}x\right\Vert
^{2}df\left( t\right) \right.  \notag \\
& \left. -\left\Vert \frac{1}{f\left( M\right) -f\left( m\right) }%
\int_{m-0}^{M}E_{s}xdf\left( s\right) \right\Vert ^{2}\right]  \notag
\end{align}%
and, integrating by parts in the Riemann-Stieltjes integral, we have%
\begin{align}
\int_{m-0}^{M}\left\Vert E_{t}x\right\Vert ^{2}df\left( t\right) &
=\int_{m-0}^{M}\left\langle E_{t}x,E_{t}x\right\rangle df\left( t\right)
=\int_{m-0}^{M}\left\langle E_{t}x,x\right\rangle df\left( t\right)
\label{IV.c.e.2.21} \\
& =f\left( M\right) \left\Vert x\right\Vert ^{2}-\int_{m-0}^{M}f\left(
t\right) d\left\langle E_{t}x,x\right\rangle  \notag \\
& =f\left( M\right) \left\Vert x\right\Vert ^{2}-\left\langle f\left(
A\right) x,x\right\rangle =\left\langle \left[ f\left( M\right)
1_{H}-f\left( A\right) \right] x,x\right\rangle  \notag
\end{align}%
and%
\begin{equation}
\int_{m-0}^{M}E_{s}xdf\left( s\right) =f\left( M\right) x-f\left( A\right) x
\label{IV.c.e.2.22}
\end{equation}%
for any $x\in H.$

On making use of the equalities (\ref{IV.c.e.2.21}) and (\ref{IV.c.e.2.22})
we have%
\begin{align}
& \frac{1}{f\left( M\right) -f\left( m\right) }\int_{m-0}^{M}\left\Vert
E_{t}x\right\Vert ^{2}df\left( t\right) -\left\Vert \frac{1}{f\left(
M\right) -f\left( m\right) }\int_{m-0}^{M}E_{s}xdf\left( s\right)
\right\Vert ^{2}  \label{IV.c.e.2.23} \\
& =\frac{1}{\left[ f\left( M\right) -f\left( m\right) \right] ^{2}}  \notag
\\
& \times \left[ \left[ f\left( M\right) -f\left( m\right) \right]
\left\langle \left[ f\left( M\right) 1_{H}-f\left( A\right) \right]
x,x\right\rangle -\left\Vert f\left( M\right) x-f\left( A\right)
x\right\Vert ^{2}\right]  \notag \\
& =\frac{\left[ f\left( M\right) -f\left( m\right) \right] \left\langle %
\left[ f\left( M\right) 1_{H}-f\left( A\right) \right] x,x\right\rangle
-\left\langle f\left( M\right) x-f\left( A\right) x,f\left( M\right)
x-f\left( A\right) x\right\rangle }{\left[ f\left( M\right) -f\left(
m\right) \right] ^{2}}  \notag \\
& =\frac{\left[ f\left( M\right) -f\left( m\right) \right] \left\langle %
\left[ f\left( M\right) 1_{H}-f\left( A\right) \right] x,x\right\rangle
-\left\langle f\left( M\right) x-f\left( A\right) x,f\left( M\right)
x-f\left( A\right) x\right\rangle }{\left[ f\left( M\right) -f\left(
m\right) \right] ^{2}}  \notag \\
& =\frac{\left\langle f\left( M\right) x-f\left( A\right) x,f\left( A\right)
x-f\left( m\right) x\right\rangle }{\left[ f\left( M\right) -f\left(
m\right) \right] ^{2}}  \notag
\end{align}%
for any $x\in H.$

Therefore, we obtain the following equality of interest in itself as well%
\begin{align}
& \frac{1}{f\left( M\right) -f\left( m\right) }\int_{m-0}^{M}\left\Vert
E_{t}x-\frac{1}{f\left( M\right) -f\left( m\right) }\int_{m-0}^{M}E_{s}xdf%
\left( s\right) \right\Vert ^{2}df\left( t\right)  \label{IV.c.e.2.24} \\
& =\frac{\left\langle f\left( M\right) x-f\left( A\right) x,f\left( A\right)
x-f\left( m\right) x\right\rangle }{\left[ f\left( M\right) -f\left(
m\right) \right] ^{2}}  \notag \\
& =\frac{\left\langle \left[ f\left( M\right) 1_{H}-f\left( A\right) \right] %
\left[ f\left( A\right) -f\left( m\right) 1_{H}\right] x,x\right\rangle }{%
\left[ f\left( M\right) -f\left( m\right) \right] ^{2}}  \notag
\end{align}%
for any $x\in H$

On making use of the inequality (\ref{IV.c.e.2.19}) we deduce the second
inequality in (\ref{IV.c.e.2.17}).

The last part follows by (\ref{IV.c.e.2.16.a}) and the details are omitted.
\end{proof}

\subsection{Applications}

We consider the power function $f\left( t\right) :=t^{p}$ where $p\in 
\mathbb{R}\diagdown \left\{ 0\right\} $ and $t>0.$ The following power
inequalities hold:

\begin{proposition}
\label{IV.c.p.3.1}Let $A$ be a selfadjoint operator in the Hilbert space $H$
with the spectrum $Sp\left( A\right) \subseteq \left[ m,M\right] $ for some
real numbers with $0\leq m<M$.

If $p>0,$ then for any $x,y\in H$%
\begin{align}
& \left\vert \left\langle \left( M^{p}1_{H}-A^{p}\right) \left(
A^{p}-m^{p}1_{H}\right) x,y\right\rangle \right\vert  \label{IV.c.e.3.1} \\
& \leq \frac{\sqrt{2}}{2}B_{p}^{2}\left\Vert y\right\Vert \left( M-m\right)
\left\langle Ax-mx,Mx-Ax\right\rangle ^{1/2}  \notag \\
& \leq \frac{\sqrt{2}}{4}B_{p}^{2}\left\Vert y\right\Vert \left\Vert
x\right\Vert \left( M-m\right) ^{2}  \notag
\end{align}%
where%
\begin{equation*}
B_{p}=p\times \left\{ 
\begin{array}{cc}
M^{p-1} & \text{if }p\geq 1 \\ 
&  \\ 
m^{p-1} & \text{if }0<p<1,m>0%
\end{array}%
\right.
\end{equation*}%
and%
\begin{align}
& \left\vert \left\langle \left( A^{-p}-M^{-p}1_{H}\right) \left(
m^{-p}1_{H}-A^{-p}\right) x,y\right\rangle \right\vert  \label{IV.c.e.3.2} \\
& \leq \frac{\sqrt{2}}{2}C_{p}^{2}\left\Vert y\right\Vert \left( M-m\right)
\left\langle Ax-mx,Mx-Ax\right\rangle ^{1/2}  \notag \\
& \leq \frac{\sqrt{2}}{4}C_{p}^{2}\left\Vert y\right\Vert \left\Vert
x\right\Vert \left( M-m\right) ^{2},  \notag
\end{align}%
where%
\begin{equation*}
C_{p}=pm^{-p-1}\text{ and }m>0.
\end{equation*}
\end{proposition}

The proof follows from Theorem \ref{IV.c.t.2.2} applied for the power
function.

\begin{proposition}
\label{IV.c.p.3.2}Let $A$ be a selfadjoint operator in the Hilbert space $H$
with the spectrum $Sp\left( A\right) \subseteq \left[ m,M\right] $ for some
real numbers with $0\leq m<M$.

If $p>0,$ then for any $x,y\in H$%
\begin{align}
& \left\vert \left\langle \left( M^{p}1_{H}-A^{p}\right) \left(
A^{p}-m^{p}1_{H}\right) x,y\right\rangle \right\vert  \label{IV.c.e.3.3} \\
& \leq \frac{1}{2}\left\Vert y\right\Vert \left( M^{p}-m^{p}\right)
\left\langle \left( M^{p}1_{H}-A^{p}\right) \left( A^{p}-m^{p}1_{H}\right)
x,x\right\rangle ^{1/2}  \notag \\
& \leq \frac{1}{4}\left\Vert y\right\Vert \left\Vert x\right\Vert \left(
M^{p}-m^{p}\right) ^{2}  \notag
\end{align}%
and%
\begin{align}
& \left\vert \left\langle \left( A^{-p}-M^{-p}1_{H}\right) \left(
m^{-p}1_{H}-A^{-p}\right) x,y\right\rangle \right\vert  \label{IV.c.e.3.4} \\
& \leq \frac{1}{2}\left\Vert y\right\Vert \left( m^{-p}-M^{-p}\right)
\left\langle \left( A^{-p}-M^{-p}1_{H}\right) \left(
m^{-p}1_{H}-A^{-p}\right) x,x\right\rangle ^{1/2}  \notag \\
& \leq \frac{1}{4}\left\Vert y\right\Vert \left\Vert x\right\Vert \left(
m^{-p}-M^{-p}\right) ^{2}.  \notag
\end{align}
\end{proposition}

The proof follows from Theorem \ref{IV.c.t.2.3}.

Now, consider the logarithmic function $f\left( t\right) =\ln t,t>0.$ We have

\begin{proposition}
\label{IV.c.p.3.3}Let $A$ be a selfadjoint operator in the Hilbert space $H$
with the spectrum $Sp\left( A\right) \subseteq \left[ m,M\right] $ for some
real numbers with $0<m<M$. Then we have the inequalities%
\begin{align}
& \left\vert \left\langle \left[ \left( \ln M\right) 1_{H}-\ln A\right] %
\left[ \ln A-\left( \ln m\right) 1_{H}\right] x,y\right\rangle \right\vert
\label{IV.c.e.3.5} \\
& \leq \frac{\sqrt{2}}{2m^{2}}\left\Vert y\right\Vert \left( M-m\right)
\left\langle Ax-mx,Mx-Ax\right\rangle ^{1/2}  \notag \\
& \leq \frac{\sqrt{2}}{4}\left\Vert y\right\Vert \left\Vert x\right\Vert
\left( \frac{M}{m}-1\right) ^{2}  \notag
\end{align}%
and%
\begin{align}
& \left\vert \left\langle \left[ \left( \ln M\right) 1_{H}-\ln A\right] %
\left[ \ln A-\left( \ln m\right) 1_{H}\right] x,y\right\rangle \right\vert
\label{IV.c.e.3.6} \\
& \leq \frac{1}{2}\left\Vert y\right\Vert \left\langle \left[ \left( \ln
M\right) 1_{H}-\ln A\right] \left[ \ln A-\left( \ln m\right) 1_{H}\right]
x,x\right\rangle ^{1/2}\ln \left( \frac{M}{m}\right)  \notag \\
& \leq \frac{1}{4}\left\Vert y\right\Vert \left\Vert x\right\Vert \left[ \ln
\left( \frac{M}{m}\right) \right] ^{2}.  \notag
\end{align}
\end{proposition}

The proof follows from Theorem \ref{IV.c.t.2.2} and \ref{IV.c.t.2.3} applied
for the logarithmic function.

\bigskip

\chapter{Inequalities of Taylor Type}

\section{Introduction}

In approximating $n$-time differentiable functions around a point, perhaps
the classical Taylor's expansion is one of the simplest and most convenient
and elegant methods that has been employed in the development of Mathematics
for the last three centuries. There is probably no field of Science where
Mathematical Modelling is used not to contain in a form or another Taylor's
expansion for functions that are differentiable in a certain sense.

In the present chapter, that is intended to be developed to a later stage,
we present some error bounds in approximating $n$-time differentiable
functions of selfadjoint operators by the use of operator Taylor's type
expansions around a point or two points from its spectrum for which the
remainder is known in an integral form.

Some applications for elementary functions including the exponential and
logarithmic functions are provided as well.

\section{Taylor's Type Inequalities}

\subsection{Some Identities}

In this section, by utilizing the spectral representation theorem of
selfadjoint operators in Hilbert spaces, some error bounds in approximating $%
n$-time differentiable functions of selfadjoint operators in Hilbert Spaces
via a Taylor's type expansion are given. Applications for some elementary
functions of interest including the exponential and logarithmic functions
are also provided.

The following result provides a Taylor's type representation for a function
of selfadjoint operators in Hilbert spaces with integral remainder.

\begin{theorem}[Dragomir, 2010, \protect\cite{V.SSD7}]
\label{V.t.2.1}Let $A$ be a selfadjoint operator in the Hilbert space $H$
with the spectrum $Sp\left( A\right) \subseteq \left[ m,M\right] $ for some
real numbers $m<M$, $\left\{ E_{\lambda }\right\} _{\lambda }$ be its 
\textit{spectral family,} $I$ be a closed subinterval on $\mathbb{R}$ with $%
\left[ m,M\right] \subset \mathring{I}$ (the interior of $I)$ and let $n$ be
an integer with $n\geq 1.$ If $f:I\rightarrow \mathbb{C}$ is such that the $%
n $-th derivative $f^{\left( n\right) }$ is of bounded variation on the
interval $\left[ m,M\right] $, then for any $c\in \left[ m,M\right] $ we
have the equalities%
\begin{equation}
f\left( A\right) =\sum_{k=0}^{n}\frac{1}{k!}f^{\left( k\right) }\left(
c\right) \left( A-c1_{H}\right) ^{k}+R_{n}\left( f,c,m,M\right)
\label{V.e.2.1}
\end{equation}%
where%
\begin{equation}
R_{n}\left( f,c,m,M\right) =\frac{1}{n!}\int_{m-0}^{M}\left(
\int_{c}^{\lambda }\left( \lambda -t\right) ^{n}d\left( f^{\left( n\right)
}\left( t\right) \right) \right) dE_{\lambda }.  \label{V.e.2.2}
\end{equation}
\end{theorem}

\begin{proof}
We utilize the Taylor formula for a function $f:I\rightarrow \mathbb{C}$
whose $n$-th derivative $f^{\left( n\right) }$ is locally of bounded
variation on the interval $I$ to write the equality%
\begin{equation}
f\left( \lambda \right) =\sum_{k=0}^{n}\frac{1}{k!}f^{\left( k\right)
}\left( c\right) \left( \lambda -c\right) ^{k}+\frac{1}{n!}\int_{c}^{\lambda
}\left( \lambda -t\right) ^{n}d\left( f^{\left( n\right) }\left( t\right)
\right)  \label{V.e.2.3}
\end{equation}%
for any $\lambda ,c\in \left[ m,M\right] $, where the integral is taken in
the Riemann-Stieltjes sense.

If we integrate the equality on $\left[ m,M\right] $ in the
Riemann-Stieltjes sense with the integrator $E_{\lambda }$ we get 
\begin{align*}
\int_{m-0}^{M}f\left( \lambda \right) dE_{\lambda }& =\sum_{k=0}^{n}\frac{1}{%
k!}f^{\left( k\right) }\left( c\right) \int_{m-0}^{M}\left( \lambda
-c\right) ^{k}dE_{\lambda } \\
& +\frac{1}{n!}\int_{m-0}^{M}\left( \int_{c}^{\lambda }\left( \lambda
-t\right) ^{n}d\left( f^{\left( n\right) }\left( t\right) \right) \right)
dE_{\lambda }
\end{align*}%
which, by the spectral representation theorem, produces the equality (\ref%
{V.e.2.1}) with the representation of the remainder from (\ref{V.e.2.2}).
\end{proof}

The following particular instances are of interest for applications:

\begin{corollary}[Dragomir, 2010, \protect\cite{V.SSD7}]
\label{V.c.2.1}With the assumptions of the above Theorem \ref{V.t.2.1}, we
have the equalities%
\begin{equation}
f\left( A\right) =\sum_{k=0}^{n}\frac{1}{k!}f^{\left( k\right) }\left(
m\right) \left( A-m1_{H}\right) ^{k}+L_{n}\left( f,c,m,M\right)
\label{V.e.2.6}
\end{equation}%
where%
\begin{equation*}
L_{n}\left( f,c,m,M\right) =\frac{1}{n!}\int_{m-0}^{M}\left(
\int_{m}^{\lambda }\left( \lambda -t\right) ^{n}d\left( f^{\left( n\right)
}\left( t\right) \right) \right) dE_{\lambda }
\end{equation*}%
and%
\begin{equation}
f\left( A\right) =\sum_{k=0}^{n}\frac{1}{k!}f^{\left( k\right) }\left( \frac{%
m+M}{2}\right) \left( A-\frac{m+M}{2}1_{H}\right) ^{k}+M_{n}\left(
f,c,m,M\right)  \label{V.e.2.8}
\end{equation}%
where%
\begin{equation*}
M_{n}\left( f,c,m,M\right) =\frac{1}{n!}\int_{m-0}^{M}\left( \int_{\frac{m+M%
}{2}}^{\lambda }\left( \lambda -t\right) ^{n}d\left( f^{\left( n\right)
}\left( t\right) \right) \right) dE_{\lambda }
\end{equation*}%
and%
\begin{equation}
f\left( A\right) =\sum_{k=0}^{n}\frac{\left( -1\right) ^{k}}{k!}f^{\left(
k\right) }\left( M\right) \left( M1_{H}-A\right) ^{k}+U_{n}\left(
f,c,m,M\right)  \label{V.e.2.10}
\end{equation}%
where%
\begin{equation}
U_{n}\left( f,c,m,M\right) =\frac{\left( -1\right) ^{n+1}}{n!}%
\int_{m-0}^{M}\left( \int_{\lambda }^{M}\left( t-\lambda \right) ^{n}d\left(
f^{\left( n\right) }\left( t\right) \right) \right) dE_{\lambda },
\label{V.e.2.11}
\end{equation}%
respectively.
\end{corollary}

\begin{remark}
\label{V.r.2.1}We remark that, if the $n$-th derivative of the function $f$
considered above is absolutely continuous on the interval $\left[ m,M\right]
,$ then we have the representation (\ref{V.e.2.1}) with the remainder 
\begin{equation}
R_{n}\left( f,c,m,M\right) =\frac{1}{n!}\int_{m-0}^{M}\left(
\int_{c}^{\lambda }\left( \lambda -t\right) ^{n}f^{\left( n+1\right) }\left(
t\right) dt\right) dE_{\lambda }.  \label{V.e.2.12}
\end{equation}%
Here the integral $\int_{c}^{\lambda }\left( \lambda -t\right) ^{n}f^{\left(
n+1\right) }\left( t\right) dt$ is considered in the Lebesgue sense. Similar
representations hold true when $c$ is taken the particular values $m,M$ or $%
\frac{m+M}{2}.$
\end{remark}

Now, if we consider the exponential function, then for any selfadjoint
operator $A$ in the Hilbert space $H$ with the spectrum $Sp\left( A\right)
\subseteq \left[ m,M\right] $ and with the spectral family $\left\{
E_{\lambda }\right\} _{\lambda }$ we have the representation%
\begin{equation}
e^{A-c1_{H}}=\sum_{k=0}^{n}\frac{1}{k!}\left( A-c1_{H}\right) ^{k}+\frac{1}{%
n!}\int_{m-0}^{M}\left( \int_{c}^{\lambda }\left( \lambda -t\right)
^{n}e^{t-c}dt\right) dE_{\lambda },  \label{V.e.2.13}
\end{equation}%
where $c$ is any real number.

Further, if we consider the logarithmic function, then for any positive
definite operator $A$ with $Sp\left( A\right) \subseteq \left[ m,M\right]
\subset \left( 0,\infty \right) $ and with the spectral family $\left\{
E_{\lambda }\right\} _{\lambda }$ we have%
\begin{align}
\ln A& =\left( \ln c\right) 1_{H}+\sum_{k=1}^{n}\frac{\left( -1\right)
^{k-1}\left( A-c1_{H}\right) ^{k}}{kc^{k}}  \label{V.e.2.14} \\
& +\left( -1\right) ^{n}\int_{m-0}^{M}\left( \int_{c}^{\lambda }\frac{\left(
\lambda -t\right) ^{n}}{t^{n+1}}dt\right) dE_{\lambda }  \notag
\end{align}%
for any $c>0.$

\subsection{Some Error Bounds}

We start with the following result that provides an approximation for an $n$%
-time differentiable function of selfadjoint operators in Hilbert spaces:

\begin{theorem}[Dragomir, 2010, \protect\cite{V.SSD7}]
\label{V.t.3.1}Let $A$ be a selfadjoint operator in the Hilbert space $H$
with the spectrum $Sp\left( A\right) \subseteq \left[ m,M\right] $ for some
real numbers $m<M$, $\left\{ E_{\lambda }\right\} _{\lambda }$ be its 
\textit{spectral family,} $I$ be a closed subinterval on $\mathbb{R}$ with $%
\left[ m,M\right] \subset \mathring{I}$ (the interior of $I)$ and let $n$ be
an integer with $n\geq 1.$ If $f:I\rightarrow \mathbb{C}$ is such that the $n
$-th derivative $f^{\left( n\right) }$ is of bounded variation on the
interval $\left[ m,M\right] $, then for any $c\in \left[ m,M\right] $ we
have the inequality%
\begin{align}
& \left\vert \left\langle R_{n}\left( f,c,m,M\right) x,y\right\rangle
\right\vert   \label{V.e.3.1} \\
& =\left\vert \left\langle f\left( A\right) x,y\right\rangle -\sum_{k=0}^{n}%
\frac{1}{k!}f^{\left( k\right) }\left( c\right) \left\langle \left(
A-c1_{H}\right) ^{k}x,y\right\rangle \right\vert   \notag \\
& \leq \frac{1}{n!}\left[ \left( c-m\right)
^{n}\dbigvee\limits_{m}^{c}\left( f^{\left( n\right) }\right)
\dbigvee\limits_{m-0}^{c}\left( \left\langle E_{\left( \cdot \right)
}x,y\right\rangle \right) \right.   \notag \\
& \left. +\left( M-c\right) ^{n}\dbigvee\limits_{c}^{M}\left( f^{\left(
n\right) }\right) \dbigvee\limits_{c}^{M}\left( \left\langle E_{\left( \cdot
\right) }x,y\right\rangle \right) \right]   \notag \\
& \leq \frac{1}{n!}\max \left\{ \left( M-c\right)
^{n}\dbigvee\limits_{c}^{M}\left( f^{\left( n\right) }\right) ,\left(
c-m\right) ^{n}\dbigvee\limits_{c}^{M}\left( f^{\left( n\right) }\right)
\right\} \dbigvee\limits_{m-0}^{M}\left( \left\langle E_{\left( \cdot
\right) }x,y\right\rangle \right)   \notag \\
& \leq \frac{1}{n!}\left( \frac{1}{2}\left( M-m\right) +\left\vert c-\frac{%
m+M}{2}\right\vert \right) ^{n}\dbigvee\limits_{m}^{M}\left( f^{\left(
n\right) }\right) \dbigvee\limits_{m-0}^{M}\left( \left\langle E_{\left(
\cdot \right) }x,y\right\rangle \right) ,  \notag
\end{align}%
for any $x,y\in H.$
\end{theorem}

\begin{proof}
From the identities (\ref{V.e.2.1}) and (\ref{V.e.2.2}) we have%
\begin{align}
& \left\langle R_{n}\left( f,c,m,M\right) x,y\right\rangle  \label{V.e.3.2}
\\
& =\left\langle f\left( A\right) x,y\right\rangle -\sum_{k=0}^{n}\frac{1}{k!}%
f^{\left( k\right) }\left( c\right) \left\langle \left( A-c1_{H}\right)
^{k}x,y\right\rangle  \notag \\
& =\frac{1}{n!}\int_{m-0}^{M}\left( \int_{c}^{\lambda }\left( \lambda
-t\right) ^{n}d\left( f^{\left( n\right) }\left( t\right) \right) \right)
d\left\langle E_{\lambda }x,y\right\rangle  \notag \\
& =\frac{1}{n!}\int_{m-0}^{c}\left( \int_{c}^{\lambda }\left( \lambda
-t\right) ^{n}d\left( f^{\left( n\right) }\left( t\right) \right) \right)
d\left\langle E_{\lambda }x,y\right\rangle  \notag \\
& +\frac{1}{n!}\int_{c}^{M}\left( \int_{c}^{\lambda }\left( \lambda
-t\right) ^{n}d\left( f^{\left( n\right) }\left( t\right) \right) \right)
d\left\langle E_{\lambda }x,y\right\rangle  \notag
\end{align}%
for any $x,y\in H.$

It is well known that if $p:\left[ a,b\right] \rightarrow \mathbb{C}$ is a
continuous function, $v:\left[ a,b\right] \rightarrow \mathbb{C}$ is of
bounded variation then the Riemann-Stieltjes integral $\int_{a}^{b}p\left(
t\right) dv\left( t\right) $ exists and the following inequality holds%
\begin{equation}
\left\vert \int_{a}^{b}p\left( t\right) dv\left( t\right) \right\vert \leq
\max_{t\in \left[ a,b\right] }\left\vert p\left( t\right) \right\vert
\dbigvee\limits_{a}^{b}\left( v\right) ,  \label{V.e.3.3}
\end{equation}%
where $\dbigvee\limits_{a}^{b}\left( v\right) $ denotes the total variation
of $v$ on $\left[ a,b\right] .$

Taking the modulus in (\ref{V.e.3.2}) and utilizing the inequality (\ref%
{V.e.3.3}) we have%
\begin{align}
& \left\vert \left\langle R_{n}\left( f,c,m,M\right) x,y\right\rangle
\right\vert   \label{V.e.3.4} \\
& \leq \frac{1}{n!}\left\vert \int_{m-0}^{c}\left( \int_{c}^{\lambda }\left(
\lambda -t\right) ^{n}d\left( f^{\left( n\right) }\left( t\right) \right)
\right) d\left\langle E_{\lambda }x,y\right\rangle \right\vert   \notag \\
& +\frac{1}{n!}\left\vert \int_{c}^{M}\left( \int_{c}^{\lambda }\left(
\lambda -t\right) ^{n}d\left( f^{\left( n\right) }\left( t\right) \right)
\right) d\left\langle E_{\lambda }x,y\right\rangle \right\vert   \notag \\
& \leq \frac{1}{n!}\max_{\lambda \in \left[ m,c\right] }\left\vert
\int_{c}^{\lambda }\left( \lambda -t\right) ^{n}d\left( f^{\left( n\right)
}\left( t\right) \right) \right\vert \dbigvee\limits_{m-0}^{c}\left(
\left\langle E_{\left( \cdot \right) }x,y\right\rangle \right)   \notag \\
& +\frac{1}{n!}\max_{\lambda \in \left[ c,M\right] }\left\vert
\int_{c}^{\lambda }\left( \lambda -t\right) ^{n}d\left( f^{\left( n\right)
}\left( t\right) \right) \right\vert \dbigvee\limits_{c}^{M}\left(
\left\langle E_{\left( \cdot \right) }x,y\right\rangle \right)   \notag
\end{align}%
for any $x,y\in H.$

By the same property (\ref{V.e.3.3}) we have%
\begin{equation}
\max_{\lambda \in \left[ m,c\right] }\left\vert \int_{c}^{\lambda }\left(
\lambda -t\right) ^{n}d\left( f^{\left( n\right) }\left( t\right) \right)
\right\vert \leq \left( c-m\right) ^{n}\dbigvee\limits_{m}^{c}\left(
f^{\left( n\right) }\right)  \label{V.e.3.5}
\end{equation}%
and%
\begin{equation}
\max_{\lambda \in \left[ c,M\right] }\left\vert \int_{c}^{\lambda }\left(
\lambda -t\right) ^{n}d\left( f^{\left( n\right) }\left( t\right) \right)
\right\vert \leq \left( M-c\right) ^{n}\dbigvee\limits_{c}^{M}\left(
f^{\left( n\right) }\right) .  \label{V.e.3.6}
\end{equation}

Now, on making \ use of (\ref{V.e.3.4})-(\ref{V.e.3.6}) we deduce%
\begin{align*}
& \left\vert \left\langle R_{n}\left( f,c,m,M\right) x,y\right\rangle
\right\vert  \\
& \leq \frac{1}{n!}\left[ \left( c-m\right)
^{n}\dbigvee\limits_{m}^{c}\left( f^{\left( n\right) }\right)
\dbigvee\limits_{m-0}^{c}\left( \left\langle E_{\left( \cdot \right)
}x,y\right\rangle \right) \right.  \\
& \left. +\left( M-c\right) ^{n}\dbigvee\limits_{c}^{M}\left( f^{\left(
n\right) }\right) \dbigvee\limits_{c}^{M}\left( \left\langle E_{\left( \cdot
\right) }x,y\right\rangle \right) \right]  \\
& \leq \frac{1}{n!}\max \left\{ \left( c-m\right)
^{n}\dbigvee\limits_{m}^{c}\left( f^{\left( n\right) }\right) ,\left(
M-c\right) ^{n}\dbigvee\limits_{c}^{M}\left( f^{\left( n\right) }\right)
\right\}  \\
& \times \left[ \dbigvee\limits_{m-0}^{c}\left( \left\langle E_{\left( \cdot
\right) }x,y\right\rangle \right) +\dbigvee\limits_{c}^{M}\left(
\left\langle E_{\left( \cdot \right) }x,y\right\rangle \right) \right]  \\
& \leq \frac{1}{n!}\max \left\{ \left( c-m\right) ^{n},\left( M-c\right)
^{n}\right\} \dbigvee\limits_{m}^{M}\left( f^{\left( n\right) }\right)
\dbigvee\limits_{m-0}^{M}\left( \left\langle E_{\left( \cdot \right)
}x,y\right\rangle \right)  \\
& =\frac{1}{n!}\left( \frac{1}{2}\left( M-m\right) +\left\vert c-\frac{m+M}{2%
}\right\vert \right) ^{n}\dbigvee\limits_{m}^{M}\left( f^{\left( n\right)
}\right) \dbigvee\limits_{m-0}^{M}\left( \left\langle E_{\left( \cdot
\right) }x,y\right\rangle \right) 
\end{align*}%
for any $x,y\in H$ and the proof is complete.
\end{proof}

The following particular cases are of interest for applications

\begin{corollary}[Dragomir, 2010, \protect\cite{V.SSD7}]
\label{V.c.3.1}With the assumption of Theorem \ref{V.t.3.1} we have the
inequalities%
\begin{align}
& \left\vert \left\langle f\left( A\right) x,y\right\rangle -\sum_{k=0}^{n}%
\frac{1}{k!}f^{\left( k\right) }\left( m\right) \left\langle \left(
A-m1_{H}\right) ^{k}x,y\right\rangle \right\vert   \label{V.e.3.7} \\
& \leq \frac{1}{n!}\left( M-m\right) ^{n}\dbigvee\limits_{m}^{M}\left(
f^{\left( n\right) }\right) \dbigvee\limits_{m-0}^{M}\left( \left\langle
E_{\left( \cdot \right) }x,y\right\rangle \right)   \notag \\
& \leq \frac{1}{n!}\left( M-m\right) ^{n}\dbigvee\limits_{m}^{M}\left(
f^{\left( n\right) }\right) \left\Vert x\right\Vert \left\Vert y\right\Vert ,
\notag
\end{align}%
\begin{align}
& \left\vert \left\langle f\left( A\right) x,y\right\rangle -\sum_{k=0}^{n}%
\frac{\left( -1\right) ^{k}}{k!}f^{\left( k\right) }\left( M\right)
\left\langle \left( M1_{H}-A\right) ^{k}x,y\right\rangle \right\vert 
\label{V.e.3.8} \\
& \leq \frac{1}{n!}\left( M-m\right) ^{n}\dbigvee\limits_{m}^{M}\left(
f^{\left( n\right) }\right) \dbigvee\limits_{m-0}^{M}\left( \left\langle
E_{\left( \cdot \right) }x,y\right\rangle \right)   \notag \\
& \leq \frac{1}{n!}\left( M-m\right) ^{n}\dbigvee\limits_{m}^{M}\left(
f^{\left( n\right) }\right) \left\Vert x\right\Vert \left\Vert y\right\Vert 
\notag
\end{align}%
and%
\begin{align}
& \left\vert \left\langle f\left( A\right) x,y\right\rangle -\sum_{k=0}^{n}%
\frac{1}{k!}f^{\left( k\right) }\left( \frac{m+M}{2}\right) \left\langle
\left( A-\frac{m+M}{2}1_{H}\right) ^{k}x,y\right\rangle \right\vert 
\label{V.e.3.9} \\
& \leq \frac{1}{2^{n}n!}\left( M-m\right) ^{n}\max \left\{ \dbigvee\limits_{%
\frac{m+M}{2}}^{M}\left( f^{\left( n\right) }\right) ,\dbigvee\limits_{m}^{%
\frac{m+M}{2}}\left( f^{\left( n\right) }\right) \right\}
\dbigvee\limits_{m-0}^{M}\left( \left\langle E_{\left( \cdot \right)
}x,y\right\rangle \right)   \notag \\
& \leq \frac{1}{2^{n}n!}\left( M-m\right) ^{n}\max \left\{ \dbigvee\limits_{%
\frac{m+M}{2}}^{M}\left( f^{\left( n\right) }\right) ,\dbigvee\limits_{m}^{%
\frac{m+M}{2}}\left( f^{\left( n\right) }\right) \right\} \left\Vert
x\right\Vert \left\Vert y\right\Vert   \notag
\end{align}%
respectively, for any $x,y\in H.$
\end{corollary}

\begin{proof}
The first part in the inequalities follow from (\ref{V.e.3.1}) by choosing $%
c=m,c=M$ and $c=\frac{m+M}{2}$ respectively.

The last part follows by the Total Variation Schwarz's inequality and we
omit the details.
\end{proof}

The following result also holds:

\begin{theorem}[Dragomir, 2010, \protect\cite{V.SSD7}]
\label{V.t.3.2}Let $A$ be a selfadjoint operator in the Hilbert space $H$
with the spectrum $Sp\left( A\right) \subseteq \left[ m,M\right] $ for some
real numbers $m<M$, $\left\{ E_{\lambda }\right\} _{\lambda }$ be its 
\textit{spectral family,} $I$ be a closed subinterval on $\mathbb{R}$ with $%
\left[ m,M\right] \subset \mathring{I}$ (the interior of $I)$ and let $n$ be
an integer with $n\geq 1.$ If $f:I\rightarrow \mathbb{C}$ is such that the $n
$-th derivative $f^{\left( n\right) }$ is Lipschitzian with the constant $%
L_{n}>0$ on the interval $\left[ m,M\right] $, then for any $c\in \left[ m,M%
\right] $ we have the inequality%
\begin{align}
& \left\vert \left\langle R_{n}\left( f,c,m,M\right) x,y\right\rangle
\right\vert   \label{V.e.3.11} \\
& \leq \frac{1}{\left( n+1\right) !}L_{n}\left[ \left( c-m\right)
^{n+1}\dbigvee\limits_{m-0}^{c}\left( \left\langle E_{\left( \cdot \right)
}x,y\right\rangle \right) +\left( M-c\right)
^{n+1}\dbigvee\limits_{c}^{M}\left( \left\langle E_{\left( \cdot \right)
}x,y\right\rangle \right) \right]   \notag \\
& \leq \frac{1}{\left( n+1\right) !}L_{n}\left( \frac{1}{2}\left( M-m\right)
+\left\vert c-\frac{m+M}{2}\right\vert \right)
^{n+1}\dbigvee\limits_{m-0}^{M}\left( \left\langle E_{\left( \cdot \right)
}x,y\right\rangle \right)   \notag \\
& \leq \frac{1}{\left( n+1\right) !}L_{n}\left( \frac{1}{2}\left( M-m\right)
+\left\vert c-\frac{m+M}{2}\right\vert \right) ^{n+1}\left\Vert x\right\Vert
\left\Vert y\right\Vert   \notag
\end{align}%
for any $x,y\in H.$
\end{theorem}

\begin{proof}
First of all, recall that if $p:\left[ a,b\right] \rightarrow \mathbb{C}$ is
a Riemann integrable function and $v:\left[ a,b\right] \rightarrow \mathbb{C}
$ is Lipschitzian with the constant $L>0$, i.e.,%
\begin{equation*}
\left\vert f\left( s\right) -f\left( t\right) \right\vert \leq L\left\vert
s-t\right\vert \text{ for any }t,s\in \left[ a,b\right] ,
\end{equation*}%
then the Riemann-Stieltjes integral $\int_{a}^{b}p\left( t\right) dv\left(
t\right) $ exists and the following inequality holds%
\begin{equation*}
\left\vert \int_{a}^{b}p\left( t\right) dv\left( t\right) \right\vert \leq
L\int_{a}^{b}\left\vert p\left( t\right) \right\vert dt.
\end{equation*}

Now, on applying this property of the Riemann-Stieltjes integral we have%
\begin{align}
\max_{\lambda \in \left[ m,c\right] }\left\vert \int_{\lambda }^{c}\left(
t-\lambda \right) ^{n}d\left( f^{\left( n\right) }\left( t\right) \right)
\right\vert & \leq \max_{\lambda \in \left[ m,c\right] }\left[
L_{n}\int_{\lambda }^{c}\left( t-\lambda \right) ^{n}dt\right] 
\label{V.e.3.12} \\
& =\frac{L_{n}}{n+1}\left( c-m\right) ^{n+1}  \notag
\end{align}%
and%
\begin{align}
\max_{\lambda \in \left[ c,M\right] }\left\vert \int_{c}^{\lambda }\left(
\lambda -t\right) ^{n}d\left( f^{\left( n\right) }\left( t\right) \right)
\right\vert & \leq \max_{\lambda \in \left[ c,M\right] }\left[
L_{n}\int_{c}^{\lambda }\left( \lambda -t\right) ^{n}dt\right] 
\label{V.e.3.13} \\
& =\frac{L_{n}}{n+1}\left( M-c\right) ^{n+1}.  \notag
\end{align}%
Now, on utilizing the inequality (\ref{V.e.3.4}), then we have from (\ref%
{V.e.3.12}) and (\ref{V.e.3.13})\ that%
\begin{align}
& \left\vert \left\langle R_{n}\left( f,c,m,M\right) x,y\right\rangle
\right\vert   \label{V.e.3.14} \\
& \leq \frac{1}{\left( n+1\right) !}L_{n}\left( c-m\right)
^{n+1}\dbigvee\limits_{m-0}^{c}\left( \left\langle E_{\left( \cdot \right)
}x,y\right\rangle \right)   \notag \\
& +\frac{1}{\left( n+1\right) !}L_{n}\left( M-c\right)
^{n+1}\dbigvee\limits_{c}^{M}\left( \left\langle E_{\left( \cdot \right)
}x,y\right\rangle \right)   \notag \\
& \leq \frac{1}{\left( n+1\right) !}L_{n}\max \left\{ \left( c-m\right)
^{n+1},\left( M-c\right) ^{n+1}\right\} \dbigvee\limits_{m-0}^{M}\left(
\left\langle E_{\left( \cdot \right) }x,y\right\rangle \right)   \notag \\
& =\frac{1}{\left( n+1\right) !}L_{n}\left( \frac{1}{2}\left( M-m\right)
+\left\vert c-\frac{m+M}{2}\right\vert \right)
^{n+1}\dbigvee\limits_{m-0}^{M}\left( \left\langle E_{\left( \cdot \right)
}x,y\right\rangle \right) ,  \notag
\end{align}%
and the proof is complete.
\end{proof}

The following particular cases are of interest for applications:

\begin{corollary}[Dragomir, 2010, \protect\cite{V.SSD7}]
\label{V.c.3.2}With the assumption of Theorem \ref{V.t.3.2} we have the
inequalities%
\begin{align}
& \left\vert \left\langle f\left( A\right) x,y\right\rangle -\sum_{k=0}^{n}%
\frac{1}{k!}f^{\left( k\right) }\left( m\right) \left\langle \left(
A-m1_{H}\right) ^{k}x,y\right\rangle \right\vert   \label{V.e.3.15} \\
& \leq \frac{1}{\left( n+1\right) !}\left( M-m\right)
^{n+1}L_{n}\dbigvee\limits_{m-0}^{M}\left( \left\langle E_{\left( \cdot
\right) }x,y\right\rangle \right)   \notag \\
& \leq \frac{1}{\left( n+1\right) !}\left( M-m\right) ^{n+1}L_{n}\left\Vert
x\right\Vert \left\Vert y\right\Vert   \notag
\end{align}%
and%
\begin{align}
& \left\vert \left\langle f\left( A\right) x,y\right\rangle -\sum_{k=0}^{n}%
\frac{\left( -1\right) ^{k}}{k!}f^{\left( k\right) }\left( M\right)
\left\langle \left( M1_{H}-A\right) ^{k}x,y\right\rangle \right\vert 
\label{V.e.3.16} \\
& \leq \frac{1}{\left( n+1\right) !}\left( M-m\right)
^{n+1}L_{n}\dbigvee\limits_{m-0}^{M}\left( \left\langle E_{\left( \cdot
\right) }x,y\right\rangle \right)   \notag \\
& \leq \frac{1}{\left( n+1\right) !}\left( M-m\right) ^{n+1}L_{n}\left\Vert
x\right\Vert \left\Vert y\right\Vert   \notag
\end{align}%
and%
\begin{align}
& \left\vert \left\langle f\left( A\right) x,y\right\rangle -\sum_{k=0}^{n}%
\frac{1}{k!}f^{\left( k\right) }\left( \frac{m+M}{2}\right) \left\langle
\left( A-\frac{m+M}{2}1_{H}\right) ^{k}x,y\right\rangle \right\vert 
\label{V.e.3.17} \\
& \leq \frac{1}{2^{n+1}\left( n+1\right) !}\left( M-m\right)
^{n+1}L_{n}\dbigvee\limits_{m-0}^{M}\left( \left\langle E_{\left( \cdot
\right) }x,y\right\rangle \right)   \notag \\
& \leq \frac{1}{2^{n+1}\left( n+1\right) !}\left( M-m\right)
^{n+1}L_{n}\left\Vert x\right\Vert \left\Vert y\right\Vert   \notag
\end{align}%
respectively, for any $x,y\in H.$
\end{corollary}

The following corollary that provides a perturbed version of Taylor's
expansion holds:

\begin{corollary}[Dragomir, 2010, \protect\cite{V.SSD7}]
\label{V.c.3.3}Let $A$ be a selfadjoint operator in the Hilbert space $H$
with the spectrum $Sp\left( A\right) \subseteq \left[ m,M\right] $ for some
real numbers $m<M$, $\left\{ E_{\lambda }\right\} _{\lambda }$ be its 
\textit{spectral family,} $I$ be a closed subinterval on $\mathbb{R}$ with $%
\left[ m,M\right] \subset \mathring{I}$ (the interior of $I)$ and let $n$ be
an integer with $n\geq 1.$ If $g:I\rightarrow \mathbb{R}$ is such that the $n
$-th derivative $g^{\left( n\right) }$ is $\left( l_{n},L_{n}\right) -$%
Lipschitzian with the constant $L_{n}>l_{n}>0$ on the interval $\left[ m,M%
\right] $, then for any $c\in \left[ m,M\right] $ we have the inequality%
\begin{align}
& \left\vert \left\langle g\left( A\right) x,y\right\rangle -g\left(
c\right) \left\langle x,y\right\rangle -\sum_{k=1}^{n}\frac{1}{k!}g^{\left(
k\right) }\left( c\right) \left\langle \left( A-c1_{H}\right)
^{k}x,y\right\rangle -\frac{l_{n}+L_{n}}{2}\right.   \label{V.e.3.21} \\
& \times \left[ \frac{1}{\left( n+1\right) !}\left\langle
A^{n+1}x,y\right\rangle -\frac{c^{n+1}}{\left( n+1\right) !}\left\langle
x,y\right\rangle \right.   \notag \\
& \left. \left. -\sum_{k=1}^{n}\frac{c^{n-k+1}}{k!\left( n-k+1\right) !}%
\left\langle \left( A-c1_{H}\right) ^{k}x,y\right\rangle \right] \right\vert 
\notag \\
& \leq \frac{1}{2\left( n+1\right) !}\left( L_{n}-l_{n}\right)   \notag \\
& \times \left[ \left( c-m\right) ^{n+1}\dbigvee\limits_{m-0}^{c}\left(
\left\langle E_{\left( \cdot \right) }x,y\right\rangle \right) +\left(
M-c\right) ^{n+1}\dbigvee\limits_{c}^{M}\left( \left\langle E_{\left( \cdot
\right) }x,y\right\rangle \right) \right]   \notag \\
& \leq \frac{1}{2\left( n+1\right) !}\left( L_{n}-l_{n}\right) \left( \frac{1%
}{2}\left( M-m\right) +\left\vert c-\frac{m+M}{2}\right\vert \right)
^{n+1}\dbigvee\limits_{m-0}^{M}\left( \left\langle E_{\left( \cdot \right)
}x,y\right\rangle \right)   \notag \\
& \leq \frac{1}{2\left( n+1\right) !}\left( L_{n}-l_{n}\right) \left( \frac{1%
}{2}\left( M-m\right) +\left\vert c-\frac{m+M}{2}\right\vert \right)
^{n+1}\left\Vert x\right\Vert \left\Vert y\right\Vert   \notag
\end{align}%
for any $x,y\in H.$
\end{corollary}

\begin{proof}
Consider the function $f:I\rightarrow \mathbb{R}$ defined by 
\begin{equation*}
f\left( t\right) :=g\left( t\right) -\frac{1}{\left( n+1\right) !}\frac{%
L_{n}+l_{n}}{2}\cdot t^{n+1}.
\end{equation*}%
Observe that 
\begin{equation*}
f^{\left( k\right) }\left( t\right) :=g^{\left( k\right) }\left( t\right) -%
\frac{1}{\left( n-k+1\right) !}\frac{L_{n}+l_{n}}{2}\cdot t^{n-k+1}
\end{equation*}%
for any $k=0,...,n.$

Since $g^{\left( n\right) }$ is $\left( l_{n},L_{n}\right) -$Lipschitzian it
follows that 
\begin{equation*}
f^{\left( n\right) }\left( t\right) :=g^{\left( n\right) }\left( t\right) -%
\frac{L_{n}+l_{n}}{2}\cdot t
\end{equation*}%
is $\frac{1}{2}\left( L_{n}-l_{n}\right) $-Lipschitzian and applying Theorem %
\ref{V.t.3.2} for the function $f,$ we deduce after required calculations
the desired result (\ref{V.e.3.1}).
\end{proof}

\begin{remark}
\label{V.r.3.1}In particular, we can state from (\ref{V.e.3.21}) the
following inequalities%
\begin{align}
& \left\vert \left\langle g\left( A\right) x,y\right\rangle -g\left(
m\right) \left\langle x,y\right\rangle -\sum_{k=1}^{n}\frac{1}{k!}g^{\left(
k\right) }\left( m\right) \left\langle \left( A-m1_{H}\right)
^{k}x,y\right\rangle -\frac{l_{n}+L_{n}}{2}\right.   \label{V.e.3.22} \\
& \times \left[ \frac{1}{\left( n+1\right) !}\left\langle
A^{n+1}x,y\right\rangle -\frac{m^{n+1}}{\left( n+1\right) !}\left\langle
x,y\right\rangle \right.   \notag \\
& \left. \left. -\sum_{k=1}^{n}\frac{m^{n-k+1}}{k!\left( n-k+1\right) !}%
\left\langle \left( A-m1_{H}\right) ^{k}x,y\right\rangle \right] \right\vert 
\notag \\
& \leq \frac{1}{2\left( n+1\right) !}\left( L_{n}-l_{n}\right) \left(
M-m\right) ^{n+1}\dbigvee\limits_{m-0}^{M}\left( \left\langle E_{\left(
\cdot \right) }x,y\right\rangle \right)   \notag \\
& \leq \frac{1}{2\left( n+1\right) !}\left( L_{n}-l_{n}\right) \left(
M-m\right) ^{n+1}\left\Vert x\right\Vert \left\Vert y\right\Vert   \notag
\end{align}%
and%
\begin{align}
& \left\vert \left\langle g\left( A\right) x,y\right\rangle -g\left(
M\right) \left\langle x,y\right\rangle -\sum_{k=1}^{n}\frac{\left( -1\right)
^{k}}{k!}g^{\left( k\right) }\left( M\right) \left\langle \left(
M1_{H}-A\right) ^{k}x,y\right\rangle \right.   \label{V.e.3.23} \\
& -\frac{l_{n}+L_{n}}{2}\left[ \frac{1}{\left( n+1\right) !}\left\langle
A^{n+1}x,y\right\rangle -\frac{M^{n+1}}{\left( n+1\right) !}\left\langle
x,y\right\rangle \right.   \notag \\
& \left. \left. -\sum_{k=1}^{n}\left( -1\right) ^{k}\frac{M^{n-k+1}}{%
k!\left( n-k+1\right) !}\left\langle \left( M1_{H}-A\right)
^{k}x,y\right\rangle \right] \right\vert   \notag \\
& \leq \frac{1}{2\left( n+1\right) !}\left( L_{n}-l_{n}\right) \left(
M-m\right) ^{n+1}\dbigvee\limits_{m-0}^{M}\left( \left\langle E_{\left(
\cdot \right) }x,y\right\rangle \right)   \notag \\
& \leq \frac{1}{2\left( n+1\right) !}\left( L_{n}-l_{n}\right) \left(
M-m\right) ^{n+1}\left\Vert x\right\Vert \left\Vert y\right\Vert   \notag
\end{align}%
and%
\begin{align}
& \left\vert \left\langle g\left( A\right) x,y\right\rangle -g\left( \frac{%
m+M}{2}\right) \left\langle x,y\right\rangle \right.   \label{V.e.3.24} \\
& -\sum_{k=1}^{n}\frac{1}{k!}g^{\left( k\right) }\left( \frac{m+M}{2}\right)
\left\langle \left( A-\frac{m+M}{2}1_{H}\right) ^{k}x,y\right\rangle   \notag
\\
& -\frac{l_{n}+L_{n}}{2}\left[ \frac{1}{\left( n+1\right) !}\left\langle
A^{n+1}x,y\right\rangle -\frac{1}{\left( n+1\right) !}\left\langle
x,y\right\rangle \left( \frac{m+M}{2}\right) ^{n+1}\right.   \notag \\
& \left. \left. -\sum_{k=1}^{n}\frac{1}{\left( n-k+1\right) !k!}\left( \frac{%
m+M}{2}\right) ^{n-k+1}\left\langle \left( A-\frac{m+M}{2}1_{H}\right)
^{k}x,y\right\rangle \right] \right\vert   \notag \\
& \leq \frac{1}{2^{n+2}\left( n+1\right) !}\left( L_{n}-l_{n}\right) \left(
M-m\right) ^{n+1}\dbigvee\limits_{m-0}^{M}\left( \left\langle E_{\left(
\cdot \right) }x,y\right\rangle \right)   \notag \\
& \leq \frac{1}{2^{n+2}\left( n+1\right) !}\left( L_{n}-l_{n}\right) \left(
M-m\right) ^{n+1}\left\Vert x\right\Vert \left\Vert y\right\Vert   \notag
\end{align}%
respectively, for any $x,y\in H.$
\end{remark}

\subsection{Applications}

By utilizing Theorem \ref{V.t.3.1} and \ref{V.t.3.2} for the exponential
function, we can state the following result:

\begin{proposition}
\label{V.p.4.1}Let $A$ be a selfadjoint operator in the Hilbert space $H$
with the spectrum $Sp\left( A\right) \subseteq \left[ m,M\right] $ for some
real numbers $m<M$ and $\left\{ E_{\lambda }\right\} _{\lambda }$ be its 
\textit{spectral family, }then for any $c\in \left[ m,M\right] $ we have the
inequality%
\begin{align}
& \left\vert \left\langle e^{A}x,y\right\rangle -e^{c}\sum_{k=0}^{n}\frac{1}{%
k!}\left\langle \left( A-c1_{H}\right) ^{k}x,y\right\rangle \right\vert 
\label{V.e.4.1} \\
& \leq \frac{1}{n!}\left[ \left( c-m\right) ^{n}\left( e^{c}-e^{m}\right)
\dbigvee\limits_{m-0}^{c}\left( \left\langle E_{\left( \cdot \right)
}x,y\right\rangle \right) \right.   \notag \\
& \left. +\left( M-c\right) ^{n}\left( e^{M}-e^{c}\right)
\dbigvee\limits_{c}^{M}\left( \left\langle E_{\left( \cdot \right)
}x,y\right\rangle \right) \right]   \notag \\
& \leq \frac{1}{n!}\max \left\{ \left( M-c\right) ^{n}\left(
e^{M}-e^{c}\right) ,\left( c-m\right) ^{n}\left( e^{c}-e^{m}\right) \right\}
\dbigvee\limits_{m-0}^{M}\left( \left\langle E_{\left( \cdot \right)
}x,y\right\rangle \right)   \notag \\
& \leq \frac{1}{n!}\left( \frac{1}{2}\left( M-m\right) +\left\vert c-\frac{%
m+M}{2}\right\vert \right) ^{n}\left( e^{M}-e^{m}\right)
\dbigvee\limits_{m-0}^{M}\left( \left\langle E_{\left( \cdot \right)
}x,y\right\rangle \right)   \notag \\
& \leq \frac{1}{n!}\left( \frac{1}{2}\left( M-m\right) +\left\vert c-\frac{%
m+M}{2}\right\vert \right) ^{n}\left( e^{M}-e^{m}\right) \left\Vert
x\right\Vert \left\Vert y\right\Vert   \notag
\end{align}%
and%
\begin{align}
& \left\vert \left\langle e^{A}x,y\right\rangle -e^{c}\sum_{k=0}^{n}\frac{1}{%
k!}\left\langle \left( A-c1_{H}\right) ^{k}x,y\right\rangle \right\vert 
\label{V.e.4.2} \\
& \leq \frac{1}{\left( n+1\right) !}e^{M}\left[ \left( c-m\right)
^{n+1}\dbigvee\limits_{m-0}^{c}\left( \left\langle E_{\left( \cdot \right)
}x,y\right\rangle \right) +\left( M-c\right)
^{n+1}\dbigvee\limits_{c}^{M}\left( \left\langle E_{\left( \cdot \right)
}x,y\right\rangle \right) \right]   \notag \\
& \leq \frac{1}{\left( n+1\right) !}e^{M}\left( \frac{1}{2}\left( M-m\right)
+\left\vert c-\frac{m+M}{2}\right\vert \right)
^{n+1}\dbigvee\limits_{m-0}^{M}\left( \left\langle E_{\left( \cdot \right)
}x,y\right\rangle \right)   \notag \\
& \leq \frac{1}{\left( n+1\right) !}e^{M}\left( \frac{1}{2}\left( M-m\right)
+\left\vert c-\frac{m+M}{2}\right\vert \right) ^{n+1}\left\Vert x\right\Vert
\left\Vert y\right\Vert   \notag
\end{align}%
for any $x,y\in H.$
\end{proposition}

\begin{remark}
\label{V.r.4.1}We observe that the best inequalities we can get from (\ref%
{V.e.4.1}) and (\ref{V.e.4.2}) are%
\begin{align}
& \left\vert \left\langle e^{A}x,y\right\rangle -e^{\frac{m+M}{2}%
}\sum_{k=0}^{n}\frac{1}{k!}\left\langle \left( A-\frac{m+M}{2}1_{H}\right)
^{k}x,y\right\rangle \right\vert   \label{V.e.4.3} \\
& \leq \frac{1}{2^{n}n!}\left( M-m\right) ^{n}\left( e^{M}-e^{m}\right)
\dbigvee\limits_{m-0}^{M}\left( \left\langle E_{\left( \cdot \right)
}x,y\right\rangle \right)   \notag \\
& \leq \frac{1}{2^{n}n!}\left( M-m\right) ^{n}\left( e^{M}-e^{m}\right)
\left\Vert x\right\Vert \left\Vert y\right\Vert   \notag
\end{align}%
and%
\begin{align}
& \left\vert \left\langle e^{A}x,y\right\rangle -e^{\frac{m+M}{2}%
}\sum_{k=0}^{n}\frac{1}{k!}\left\langle \left( A-\frac{m+M}{2}1_{H}\right)
^{k}x,y\right\rangle \right\vert   \label{V.e.4.4} \\
& \leq \frac{1}{2^{n+1}\left( n+1\right) !}e^{M}\left( M-m\right)
^{n+1}\dbigvee\limits_{m-0}^{M}\left( \left\langle E_{\left( \cdot \right)
}x,y\right\rangle \right)   \notag \\
& \leq \frac{1}{2^{n+1}\left( n+1\right) !}e^{M}\left( M-m\right)
^{n+1}\left\Vert x\right\Vert \left\Vert y\right\Vert   \notag
\end{align}%
for any $x,y\in H.$
\end{remark}

The same Theorems \ref{V.t.3.1} and \ref{V.t.3.2} applied for the
logarithmic function produce:

\begin{proposition}
\label{V.p.4.2}Let $A$ be a positive definite operator in the Hilbert space $%
H$ with the spectrum $Sp\left( A\right) \subseteq \left[ m,M\right] \subset
\left( 0,\infty \right) $ and $\left\{ E_{\lambda }\right\} _{\lambda }$ be
its \textit{spectral family, }then for any $c\in \left[ m,M\right] $ we have
the inequalities%
\begin{align}
& \left\vert \left\langle \ln Ax,y\right\rangle -\left\langle
x,y\right\rangle \ln c-\sum_{k=1}^{n}\frac{\left( -1\right)
^{k-1}\left\langle \left( A-c1_{H}\right) ^{k}x,y\right\rangle }{kc^{k}}%
\right\vert   \label{V.e.4.5} \\
& \leq \frac{1}{n}\left[ \frac{\left( c-m\right) ^{n}\left(
c^{n}-m^{n}\right) }{c^{n}m^{n}}\dbigvee\limits_{m-0}^{c}\left( \left\langle
E_{\left( \cdot \right) }x,y\right\rangle \right) \right.   \notag \\
& \left. +\frac{\left( M-c\right) ^{n}\left( M^{n}-c^{n}\right) }{M^{m}c^{m}}%
\dbigvee\limits_{c}^{M}\left( \left\langle E_{\left( \cdot \right)
}x,y\right\rangle \right) \right]   \notag \\
& \leq \frac{1}{n}\max \left\{ \frac{\left( c-m\right) ^{n}\left(
c^{n}-m^{n}\right) }{c^{n}m^{n}},\frac{\left( M-c\right) ^{n}\left(
M^{n}-c^{n}\right) }{M^{m}c^{m}}\right\} \dbigvee\limits_{m-0}^{M}\left(
\left\langle E_{\left( \cdot \right) }x,y\right\rangle \right)   \notag \\
& \leq \frac{1}{n}\left( \frac{1}{2}\left( M-m\right) +\left\vert c-\frac{m+M%
}{2}\right\vert \right) ^{n}\frac{\left( M^{n}-m^{n}\right) }{M^{m}m^{m}}%
\dbigvee\limits_{m-0}^{M}\left( \left\langle E_{\left( \cdot \right)
}x,y\right\rangle \right)   \notag \\
& \leq \frac{1}{n}\left( \frac{1}{2}\left( M-m\right) +\left\vert c-\frac{m+M%
}{2}\right\vert \right) ^{n}\frac{\left( M^{n}-m^{n}\right) }{M^{m}m^{m}}%
\left\Vert x\right\Vert \left\Vert y\right\Vert   \notag
\end{align}%
and%
\begin{align}
& \left\vert \left\langle \ln Ax,y\right\rangle -\left\langle
x,y\right\rangle \ln c-\sum_{k=1}^{n}\frac{\left( -1\right)
^{k-1}\left\langle \left( A-c1_{H}\right) ^{k}x,y\right\rangle }{kc^{k}}%
\right\vert   \label{V.e.4.6} \\
& \leq \frac{1}{\left( n+1\right) m^{n+1}}\left[ \left( c-m\right)
^{n+1}\dbigvee\limits_{m-0}^{c}\left( \left\langle E_{\left( \cdot \right)
}x,y\right\rangle \right) +\left( M-c\right)
^{n+1}\dbigvee\limits_{c}^{M}\left( \left\langle E_{\left( \cdot \right)
}x,y\right\rangle \right) \right]   \notag \\
& \leq \frac{1}{\left( n+1\right) m^{n+1}}\left( \frac{1}{2}\left(
M-m\right) +\left\vert c-\frac{m+M}{2}\right\vert \right)
^{n+1}\dbigvee\limits_{m-0}^{M}\left( \left\langle E_{\left( \cdot \right)
}x,y\right\rangle \right)   \notag \\
& \leq \frac{1}{\left( n+1\right) m^{n+1}}\left( \frac{1}{2}\left(
M-m\right) +\left\vert c-\frac{m+M}{2}\right\vert \right) ^{n+1}\left\Vert
x\right\Vert \left\Vert y\right\Vert   \notag
\end{align}%
for any $x,y\in H.$
\end{proposition}

\begin{remark}
\label{V.r.4.2}The best inequalities we can get from (\ref{V.e.4.5}) and (%
\ref{V.e.4.6}) are%
\begin{align}
& \left\vert \left\langle \ln Ax,y\right\rangle -\left\langle
x,y\right\rangle \ln \left( \frac{m+M}{2}\right) -\sum_{k=1}^{n}\frac{\left(
-1\right) ^{k-1}\left\langle \left( A-\frac{m+M}{2}1_{H}\right)
^{k}x,y\right\rangle }{k\left( \frac{m+M}{2}\right) ^{k}}\right\vert 
\label{V.e.4.7} \\
& \leq \frac{1}{2^{n}n}\left( M-m\right) ^{n}\frac{\left( M^{n}-m^{n}\right) 
}{M^{m}m^{m}}\dbigvee\limits_{m-0}^{M}\left( \left\langle E_{\left( \cdot
\right) }x,y\right\rangle \right)   \notag \\
& \leq \frac{1}{2^{n}n}\left( M-m\right) ^{n}\frac{\left( M^{n}-m^{n}\right) 
}{M^{m}m^{m}}\left\Vert x\right\Vert \left\Vert y\right\Vert   \notag
\end{align}%
and%
\begin{align}
& \left\vert \left\langle \ln Ax,y\right\rangle -\left\langle
x,y\right\rangle \ln \left( \frac{m+M}{2}\right) -\sum_{k=1}^{n}\frac{\left(
-1\right) ^{k-1}\left\langle \left( A-\frac{m+M}{2}1_{H}\right)
^{k}x,y\right\rangle }{k\left( \frac{m+M}{2}\right) ^{k}}\right\vert 
\label{V.e.4.8} \\
& \leq \frac{1}{2^{n+1}\left( n+1\right) }\left( \frac{M}{m}-1\right)
^{n+1}\dbigvee\limits_{m-0}^{M}\left( \left\langle E_{\left( \cdot \right)
}x,y\right\rangle \right)   \notag \\
& \leq \frac{1}{2^{n+1}\left( n+1\right) }\left( \frac{M}{m}-1\right)
^{n+1}\left\Vert x\right\Vert \left\Vert y\right\Vert   \notag
\end{align}%
for any $x,y\in H.$
\end{remark}

\section{Perturbed Version}

\subsection{Some Identities}

The following result provides a perturbed Taylor's type representation for a
function of selfadjoint operators in Hilbert spaces.

\begin{theorem}[Dragomir, 2010, \protect\cite{V.a.SSD7}]
\label{V.a.t.2.1}Let $A$ be a selfadjoint operator in the Hilbert space $H$
with the spectrum $Sp\left( A\right) \subseteq \left[ m,M\right] $ for some
real numbers $m<M$, $\left\{ E_{\lambda }\right\} _{\lambda }$ be its 
\textit{spectral family,} $I$ be a closed subinterval on $\mathbb{R}$ with $%
\left[ m,M\right] \subset \mathring{I}$ (the interior of $I)$ and let $n$ be
an integer with $n\geq 1.$ If $f:I\rightarrow \mathbb{C}$ is such that the $%
n $-th derivative $f^{\left( n\right) }$ is of bounded variation on the
interval $\left[ m,M\right] $, then for any $c\in \left[ m,M\right] $ we
have the equalities%
\begin{align}
f\left( A\right) & =\sum_{k=0}^{n}\frac{1}{k!}f^{\left( k\right) }\left(
c\right) \left( A-c1_{H}\right) ^{k}  \label{V.a.e.2.1} \\
& +\left[ f\left( M\right) -\sum_{k=0}^{n}\frac{1}{k!}f^{\left( k\right)
}\left( c\right) \left( M-c\right) ^{k}\right] 1_{H}+V_{n}\left(
f,c,m,M\right)  \notag
\end{align}%
where%
\begin{equation}
V_{n}\left( f,c,m,M\right) :=\frac{\left( -1\right) ^{n}}{\left( n-1\right) !%
}\int_{m-0}^{M}\left( \int_{c}^{\lambda }\left( t-\lambda \right)
^{n-1}d\left( f^{\left( n\right) }\left( t\right) \right) \right) E_{\lambda
}d\lambda .  \label{V.a.e.2.2}
\end{equation}
\end{theorem}

\begin{proof}
We utilize the Taylor's formula for functions $f:I\rightarrow \mathbb{C}$
whose $n$-th derivative $f^{\left( n\right) }$ is locally of bounded
variation on the interval $I$ to write the equality%
\begin{equation}
f\left( \lambda \right) =\sum_{k=0}^{n}\frac{1}{k!}f^{\left( k\right)
}\left( c\right) \left( \lambda -c\right) ^{k}+\frac{1}{n!}\int_{c}^{\lambda
}\left( \lambda -t\right) ^{n}d\left( f^{\left( n\right) }\left( t\right)
\right)  \label{V.a.e.2.3}
\end{equation}%
for any $\lambda ,c\in \left[ m,M\right] $, where the integral is taken in
the Riemann-Stieltjes sense.

If we integrate the equality on $\left[ m,M\right] $ in the
Riemann-Stieltjes sense with the integrator $E_{\lambda }$ we get 
\begin{align*}
\int_{m-0}^{M}f\left( \lambda \right) dE_{\lambda }& =\sum_{k=0}^{n}\frac{1}{%
k!}f^{\left( k\right) }\left( c\right) \int_{m-0}^{M}\left( \lambda
-c\right) ^{k}dE_{\lambda } \\
& +\frac{1}{n!}\int_{m-0}^{M}\left( \int_{c}^{\lambda }\left( \lambda
-t\right) ^{n}d\left( f^{\left( n\right) }\left( t\right) \right) \right)
dE_{\lambda }
\end{align*}%
which, by the spectral representation theorem, produces the equality%
\begin{equation}
f\left( A\right) =\sum_{k=0}^{n}\frac{1}{k!}f^{\left( k\right) }\left(
c\right) \left( A-c1_{H}\right) ^{k}+\frac{1}{n!}\int_{m-0}^{M}\left(
\int_{c}^{\lambda }\left( \lambda -t\right) ^{n}d\left( f^{\left( n\right)
}\left( t\right) \right) \right) dE_{\lambda }  \label{V.a.e.2.3.a}
\end{equation}%
that is of interest in itself as well.

Now, integrating by parts in the Riemann-Stieltjes integral and using the
Leibnitz formula for integrals with parameters, we have%
\begin{align}
& \int_{m-0}^{M}\left( \int_{c}^{\lambda }\left( \lambda -t\right)
^{n}d\left( f^{\left( n\right) }\left( t\right) \right) \right) dE_{\lambda }
\label{V.a.e.2.4} \\
& =\left. E_{\lambda }\left( \int_{c}^{\lambda }\left( \lambda -t\right)
^{n}d\left( f^{\left( n\right) }\left( t\right) \right) \right) \right\vert
_{m-0}^{M}  \notag \\
& -\int_{m-0}^{M}E_{\lambda }d\left( \int_{c}^{\lambda }\left( \lambda
-t\right) ^{n}d\left( f^{\left( n\right) }\left( t\right) \right) \right) 
\notag \\
& =\left( \int_{c}^{M}\left( M-t\right) ^{n}d\left( f^{\left( n\right)
}\left( t\right) \right) \right) 1_{H}  \notag \\
& -n\int_{m-0}^{M}\left( \int_{c}^{\lambda }\left( \lambda -t\right)
^{n-1}d\left( f^{\left( n\right) }\left( t\right) \right) \right) E_{\lambda
}d\lambda  \notag
\end{align}%
and, since by the Taylor's formula (\ref{V.a.e.2.3}) we have%
\begin{equation}
\frac{1}{n!}\int_{c}^{M}\left( M-t\right) ^{n}d\left( f^{\left( n\right)
}\left( t\right) \right) =f\left( M\right) -\sum_{k=0}^{n}\frac{1}{k!}%
f^{\left( k\right) }\left( c\right) \left( M-c\right) ^{k},
\label{V.a.e.2.5}
\end{equation}%
then, by (\ref{V.a.e.2.3.a}) and (\ref{V.a.e.2.5}), we deduce the equality (%
\ref{V.a.e.2.1}) with the integral representation for the remainder provided
by (\ref{V.a.e.2.2}).
\end{proof}

The following particular instances are of interest for applications:

\begin{corollary}[Dragomir, 2010, \protect\cite{V.a.SSD7}]
\label{V.a.c.2.1}With the assumptions of the above Theorem \ref{V.a.t.2.1},
we have the equalities%
\begin{align}
f\left( A\right) & =\sum_{k=0}^{n}\frac{1}{k!}f^{\left( k\right) }\left(
m\right) \left( A-m1_{H}\right) ^{k}  \label{V.a.e.2.6} \\
& +\left[ f\left( M\right) -\sum_{k=0}^{n}\frac{1}{k!}f^{\left( k\right)
}\left( m\right) \left( M-m\right) ^{k}\right] 1_{H}+T_{n}\left(
f,c,m,M\right)  \notag
\end{align}%
where%
\begin{equation}
T_{n}\left( f,m,M\right) :=-\frac{1}{\left( n-1\right) !}\int_{m-0}^{M}%
\left( \int_{m}^{\lambda }\left( \lambda -t\right) ^{n-1}d\left( f^{\left(
n\right) }\left( t\right) \right) \right) E_{\lambda }d\lambda
\label{V.a.e.2.7}
\end{equation}%
and%
\begin{align}
f\left( A\right) & =\sum_{k=0}^{n}\frac{1}{k!}f^{\left( k\right) }\left( 
\frac{m+M}{2}\right) \left( A-\frac{m+M}{2}1_{H}\right) ^{k}
\label{V.a.e.2.8} \\
& +\left[ f\left( M\right) -\sum_{k=0}^{n}\frac{1}{k!}f^{\left( k\right)
}\left( \frac{m+M}{2}\right) \left( \frac{M-m}{2}\right) ^{k}\right] 1_{H} 
\notag \\
& +W_{n}\left( f,c,m,M\right)  \notag
\end{align}%
where%
\begin{equation}
W_{n}\left( f,m,M\right) :=\frac{\left( -1\right) ^{n}}{\left( n-1\right) !}%
\int_{m-0}^{M}\left( \int_{\frac{m+M}{2}}^{\lambda }\left( t-\lambda \right)
^{n-1}d\left( f^{\left( n\right) }\left( t\right) \right) \right) E_{\lambda
}d\lambda  \label{V.a.e.2.9}
\end{equation}%
and%
\begin{equation}
f\left( A\right) =\sum_{k=0}^{n}\frac{\left( -1\right) ^{k}}{k!}f^{\left(
k\right) }\left( M\right) \left( M1_{H}-A\right) ^{k}+Y_{n}\left(
f,c,m,M\right)  \label{V.a.e.2.10}
\end{equation}%
where%
\begin{equation}
Y_{n}\left( f,m,M\right) :=\frac{\left( -1\right) ^{n+1}}{\left( n-1\right) !%
}\int_{m-0}^{M}\left( \int_{\lambda }^{M}\left( t-\lambda \right)
^{n-1}d\left( f^{\left( n\right) }\left( t\right) \right) \right) E_{\lambda
}d\lambda ,  \label{V.a.e.2.11}
\end{equation}%
respectively.
\end{corollary}

\begin{remark}
\label{V.a.r.2.1}In order to give some examples we use the simplest
representation, namely (\ref{V.a.e.2.10}) for the exponential and the
logarithmic functions.

Let $A$ be a selfadjoint operator in the Hilbert space $H$ with the spectrum 
$Sp\left( A\right) \subseteq \left[ m,M\right] $ for some real numbers $m<M$
and let $\left\{ E_{\lambda }\right\} _{\lambda }$ be its \textit{spectral
family}. Then we have the representation 
\begin{align}
e^{A}& =e^{M}\sum_{k=0}^{n}\frac{\left( -1\right) ^{k}}{k!}\left(
M1_{H}-A\right) ^{k}  \label{V.a.e.2.12} \\
& +\frac{\left( -1\right) ^{n+1}}{\left( n-1\right) !}\int_{m-0}^{M}\left(
\int_{\lambda }^{M}\left( t-\lambda \right) ^{n-1}e^{t}dt\right) E_{\lambda
}d\lambda .  \notag
\end{align}%
In the case when $A$ is positive definite, i.e., $m>0,$ then we have the
representation%
\begin{align}
\ln A& =\left( \ln M\right) 1_{H}-\sum_{k=1}^{n}\frac{\left( M1_{H}-A\right)
^{k}}{kM^{k}}  \label{V.a.e.2.13} \\
& -n\int_{m-0}^{M}\left( \int_{\lambda }^{M}\frac{\left( t-\lambda \right)
^{n-1}}{t^{n+1}}dt\right) E_{\lambda }d\lambda .  \notag
\end{align}
\end{remark}

\subsection{Error Bounds for $f^{\left( n\right) }$ of Bounded Variation}

We start with the following result that provides an approximation for an $n$%
-time differentiable function of selfadjoint operators in Hilbert spaces:

\begin{theorem}[Dragomir, 2010, \protect\cite{V.a.SSD7}]
\label{V.a.t.3.1}Let $A$ be a selfadjoint operator in the Hilbert space $H$
with the spectrum $Sp\left( A\right) \subseteq \left[ m,M\right] $ for some
real numbers $m<M$, $\left\{ E_{\lambda }\right\} _{\lambda }$ be its 
\textit{spectral family,} $I$ be a closed subinterval on $\mathbb{R}$ with $%
\left[ m,M\right] \subset \mathring{I}$ (the interior of $I)$ and let $n$ be
an integer with $n\geq 1.$ If $f:I\rightarrow \mathbb{C}$ is such that the $%
n $-th derivative $f^{\left( n\right) }$ is of bounded variation on the
interval $\left[ m,M\right] $, then for any $c\in \left[ m,M\right] $ we
have the inequalities%
\begin{align}
& \left\vert \left\langle V_{n}\left( f,c,m,M\right) x,y\right\rangle
\right\vert  \label{V.a.e.3.1} \\
& \leq \frac{1}{\left( n-1\right) !}\int_{m-0}^{c}\left( c-\lambda \right)
^{n-1}\dbigvee\limits_{\lambda }^{c}\left( f^{\left( n\right) }\right)
\left\vert \left\langle E_{\lambda }x,y\right\rangle \right\vert d\lambda 
\notag \\
& +\frac{1}{\left( n-1\right) !}\int_{c}^{M}\left( \lambda -c\right)
^{n-1}\dbigvee\limits_{c}^{\lambda }\left( f^{\left( n\right) }\right)
\left\vert \left\langle E_{\lambda }x,y\right\rangle \right\vert d\lambda 
\notag \\
& \leq \frac{1}{\left( n-1\right) !}\dbigvee\limits_{m}^{c}\left( f^{\left(
n\right) }\right) \int_{m-0}^{c}\left( c-\lambda \right) ^{n-1}\left\vert
\left\langle E_{\lambda }x,y\right\rangle \right\vert d\lambda  \notag \\
& +\frac{1}{\left( n-1\right) !}\dbigvee\limits_{c}^{M}\left( f^{\left(
n\right) }\right) \int_{c}^{M}\left( \lambda -c\right) ^{n-1}\left\vert
\left\langle E_{\lambda }x,y\right\rangle \right\vert d\lambda  \notag \\
& \leq \frac{1}{\left( n-1\right) !}\max \left\{
\dbigvee\limits_{m}^{c}\left( f^{\left( n\right) }\right)
,\dbigvee\limits_{c}^{M}\left( f^{\left( n\right) }\right) \right\}
\int_{m-0}^{M}\left\vert \lambda -c\right\vert ^{n-1}\left\vert \left\langle
E_{\lambda }x,y\right\rangle \right\vert d\lambda  \notag \\
& \leq \frac{1}{n!}\max \left\{ \dbigvee\limits_{m}^{c}\left( f^{\left(
n\right) }\right) ,\dbigvee\limits_{c}^{M}\left( f^{\left( n\right) }\right)
\right\} B_{n}(c,m,M,x,y),  \notag
\end{align}%
for any $x,y\in H,$ where%
\begin{equation}
B_{n}(c,m,M,x,y):=\left\{ 
\begin{array}{l}
\left[ \left( M-c\right) ^{n}+\left( c-m\right) ^{n}\right] \left\Vert
x\right\Vert \left\Vert y\right\Vert ; \\ 
\\ 
C_{n}(c,m,M,x,y); \\ 
\\ 
n\left[ \frac{1}{2}\left( M-m\right) +\left\vert c-\frac{m+M}{2}\right\vert %
\right] ^{n-1} \\ 
\times \left[ \left\langle \left( M1_{H}-A\right) x,x\right\rangle
\left\langle \left( M1_{H}-A\right) y,y\right\rangle \right] ^{1/2}%
\end{array}%
\right.  \label{V.a.e.3.1.a}
\end{equation}%
and%
\begin{align}
& C_{n}(c,m,M,x,y)  \label{V.a.e.3.1.b} \\
& :=\left[ \left\langle \left[ \left( M-c\right) ^{n}1_{H}-sgn\left(
A-c1_{H}\right) \left\vert A-c1_{H}\right\vert ^{n}\right] x,x\right\rangle %
\right] ^{1/2}  \notag \\
& \times \left[ \left\langle \left[ \left( M-c\right) ^{n}1_{H}-sgn\left(
A-c1_{H}\right) \left\vert A-c1_{H}\right\vert ^{n}\right] y,y\right\rangle %
\right] ^{1/2}.  \notag
\end{align}%
Here the operator function $sgn\left( A-c1_{H}\right) \left\vert
A-c1_{H}\right\vert ^{n}$ is generated by the continuous function $sgn\left(
\cdot -c\right) \left\vert \cdot -c\right\vert ^{n}$ defined on the interval 
$\left[ m,M\right] .$
\end{theorem}

\begin{proof}
From the identities (\ref{V.a.e.2.1}) and (\ref{V.a.e.2.2}) we have%
\begin{align}
& \left\vert \left\langle V_{n}\left( f,c,m,M\right) x,y\right\rangle
\right\vert  \label{V.a.e.3.2} \\
& =\left\vert \frac{1}{\left( n-1\right) !}\int_{m-0}^{M}\left(
\int_{c}^{\lambda }\left( t-\lambda \right) ^{n-1}d\left( f^{\left( n\right)
}\left( t\right) \right) \right) \left\langle E_{\lambda }x,y\right\rangle
d\lambda \right\vert  \notag \\
& \leq \frac{1}{\left( n-1\right) !}\left\vert \int_{m-0}^{c}\left(
\int_{c}^{\lambda }\left( t-\lambda \right) ^{n-1}d\left( f^{\left( n\right)
}\left( t\right) \right) \right) \left\langle E_{\lambda }x,y\right\rangle
d\lambda \right\vert  \notag \\
& +\frac{1}{\left( n-1\right) !}\left\vert \int_{c}^{M}\left(
\int_{c}^{\lambda }\left( t-\lambda \right) ^{n-1}d\left( f^{\left( n\right)
}\left( t\right) \right) \right) \left\langle E_{\lambda }x,y\right\rangle
d\lambda \right\vert  \notag \\
& \leq \frac{1}{\left( n-1\right) !}\int_{m-0}^{c}\left\vert
\int_{c}^{\lambda }\left( t-\lambda \right) ^{n-1}d\left( f^{\left( n\right)
}\left( t\right) \right) \right\vert \left\vert \left\langle E_{\lambda
}x,y\right\rangle \right\vert d\lambda  \notag \\
& +\frac{1}{\left( n-1\right) !}\int_{c}^{M}\left\vert \int_{c}^{\lambda
}\left( t-\lambda \right) ^{n-1}d\left( f^{\left( n\right) }\left( t\right)
\right) \right\vert \left\vert \left\langle E_{\lambda }x,y\right\rangle
\right\vert d\lambda  \notag
\end{align}%
for any $x,y\in H.$

It is well known that if $p:\left[ a,b\right] \rightarrow \mathbb{C}$ is a
continuous function, $v:\left[ a,b\right] \rightarrow \mathbb{C}$ is of
bounded variation, then the Riemann-Stieltjes integral $\int_{a}^{b}p\left(
t\right) dv\left( t\right) $ exists and the following inequality holds%
\begin{equation}
\left\vert \int_{a}^{b}p\left( t\right) dv\left( t\right) \right\vert \leq
\sup_{t\in \left[ a,b\right] }\left\vert p\left( t\right) \right\vert
\dbigvee\limits_{a}^{b}\left( v\right) ,  \label{V.a.e.3.3}
\end{equation}%
where $\dbigvee\limits_{a}^{b}\left( v\right) $ denotes the total variation
of $v$ on $\left[ a,b\right] .$

By the same property (\ref{V.a.e.3.3}) we have%
\begin{equation}
\left\vert \int_{c}^{\lambda }\left( t-\lambda \right) ^{n-1}d\left(
f^{\left( n\right) }\left( t\right) \right) \right\vert \leq \left(
c-\lambda \right) ^{n-1}\dbigvee\limits_{\lambda }^{c}\left( f^{\left(
n\right) }\right)  \label{V.a.e.3.5}
\end{equation}%
for $\lambda \in \left[ m,c\right] $ and%
\begin{equation}
\left\vert \int_{c}^{\lambda }\left( t-\lambda \right) ^{n-1}d\left(
f^{\left( n\right) }\left( t\right) \right) \right\vert \leq \left( \lambda
-c\right) ^{n-1}\dbigvee\limits_{c}^{\lambda }\left( f^{\left( n\right)
}\right)  \label{V.a.e.3.6}
\end{equation}%
for $\lambda \in \left[ c,M\right] .$

Now, on making \ use of (\ref{V.a.e.3.2}) and (\ref{V.a.e.3.5})-(\ref%
{V.a.e.3.6}) we deduce%
\begin{align*}
& \left\vert \left\langle V_{n}\left( f,c,m,M\right) x,y\right\rangle
\right\vert \\
& \leq \frac{1}{\left( n-1\right) !}\int_{m-0}^{c}\left( c-\lambda \right)
^{n-1}\dbigvee\limits_{\lambda }^{c}\left( f^{\left( n\right) }\right)
\left\vert \left\langle E_{\lambda }x,y\right\rangle \right\vert d\lambda \\
& +\frac{1}{\left( n-1\right) !}\int_{c}^{M}\left( \lambda -c\right)
^{n-1}\dbigvee\limits_{c}^{\lambda }\left( f^{\left( n\right) }\right)
\left\vert \left\langle E_{\lambda }x,y\right\rangle \right\vert d\lambda
\end{align*}%
for any $x,y\in H$ which proves the first part of (\ref{V.a.e.3.1}).

The second and the third inequalities follow by the properties of the
integral.

For the last part we observe that 
\begin{align*}
\int_{m-0}^{M}\left\vert \lambda -c\right\vert ^{n-1}\left\vert \left\langle
E_{\lambda }x,y\right\rangle \right\vert d\lambda & \leq \max_{\lambda \in 
\left[ m,M\right] }\left\vert \left\langle E_{\lambda }x,y\right\rangle
\right\vert \int_{m}^{M}\left\vert \lambda -c\right\vert ^{n-1}d\lambda \\
& \leq \frac{1}{n}\left\Vert x\right\Vert \left\Vert y\right\Vert \left[
\left( M-c\right) ^{n}+\left( c-m\right) ^{n}\right]
\end{align*}%
for any $x,y\in H,$and the proof for the first branch of $B(c,m,M,x,y)$ is
complete.

Now, to prove the inequality for the second branch of $B(c,m,M,x,y)$ we use
the fact that if $P$ is a nonnegative operator on $H,$ i.e., $\left\langle
Px,x\right\rangle \geq 0$ for any $x\in H,$ then the following inequality
that provides a generalization of the Schwarz inequality in $H$ can be stated%
\begin{equation}
\left\vert \left\langle Px,y\right\rangle \right\vert ^{2}\leq \left\langle
Px,x\right\rangle \left\langle Py,y\right\rangle  \label{V.a.e.3.7}
\end{equation}%
for any $x,y\in H.$

If we use (\ref{V.a.e.3.7}) and the Cauchy-Buniakowski-Schwarz weighted
integral inequality we can write that%
\begin{align}
& \int_{m-0}^{M}\left\vert \lambda -c\right\vert ^{n-1}\left\vert
\left\langle E_{\lambda }x,y\right\rangle \right\vert d\lambda
\label{V.a.e.3.8} \\
& \leq \int_{m-0}^{M}\left\vert \lambda -c\right\vert ^{n-1}\left\langle
E_{\lambda }x,x\right\rangle ^{1/2}\left\langle E_{\lambda }y,y\right\rangle
^{1/2}d\lambda  \notag \\
& \leq \left( \int_{m-0}^{M}\left\vert \lambda -c\right\vert
^{n-1}\left\langle E_{\lambda }x,x\right\rangle d\lambda \right)
^{1/2}\left( \int_{m-0}^{M}\left\vert \lambda -c\right\vert
^{n-1}\left\langle E_{\lambda }y,y\right\rangle d\lambda \right) ^{1/2} 
\notag
\end{align}%
for any $x,y\in H.$

Integrating by parts in the Riemann-Stieltjes integral, we have%
\begin{align}
& \int_{m-0}^{M}\left\vert \lambda -c\right\vert ^{n-1}\left\langle
E_{\lambda }x,x\right\rangle d\lambda  \label{V.a.e.3.9} \\
& =\int_{m-0}^{c}\left( c-\lambda \right) ^{n-1}\left\langle E_{\lambda
}x,x\right\rangle d\lambda +\int_{c}^{M}\left( \lambda -c\right)
^{n-1}\left\langle E_{\lambda }x,x\right\rangle d\lambda  \notag \\
& =\frac{1}{n}\left[ -\int_{m-0}^{c}\left\langle E_{\lambda
}x,x\right\rangle d\left( c-\lambda \right) ^{n}+\int_{c}^{M}\left\langle
E_{\lambda }x,x\right\rangle d\left( \lambda -c\right) ^{n}\right]  \notag \\
& =\frac{1}{n}\left[ \left. -\left( c-\lambda \right) ^{n}\left\langle
E_{\lambda }x,x\right\rangle \right\vert _{m-0}^{c}+\int_{m-0}^{c}\left(
c-\lambda \right) ^{n}d\left\langle E_{\lambda }x,x\right\rangle \right] 
\notag \\
& +\frac{1}{n}\left[ \left. \left\langle E_{\lambda }x,x\right\rangle \left(
\lambda -c\right) ^{n}\right\vert _{c}^{M}-\int_{c}^{M}\left( \lambda
-c\right) ^{n}d\left\langle E_{\lambda }x,x\right\rangle \right]  \notag \\
& =\frac{1}{n}\int_{m-0}^{c}\left( c-\lambda \right) ^{n}d\left\langle
E_{\lambda }x,x\right\rangle  \notag \\
& +\frac{1}{n}\left[ \left\Vert x\right\Vert ^{2}\left( M-c\right)
^{n}-\int_{c}^{M}\left( \lambda -c\right) ^{n}d\left\langle E_{\lambda
}x,x\right\rangle \right]  \notag \\
& =\frac{1}{n}\left\Vert x\right\Vert ^{2}\left( M-c\right) ^{n}  \notag \\
& +\frac{1}{n}\left[ \int_{m-0}^{c}\left( c-\lambda \right)
^{n}d\left\langle E_{\lambda }x,x\right\rangle -\int_{c}^{M}\left( \lambda
-c\right) ^{n}d\left\langle E_{\lambda }x,x\right\rangle \right]  \notag \\
& =\frac{1}{n}\left[ \left\Vert x\right\Vert ^{2}\left( M-c\right)
^{n}-\int_{m-0}^{M}sgn\left( \lambda -c\right) \left\vert \lambda
-c\right\vert ^{n}d\left\langle E_{\lambda }x,x\right\rangle \right]  \notag
\\
& =\frac{1}{n}\left[ \left\langle \left[ \left( M-c\right)
^{n}1_{H}-sgn\left( A-c1_{H}\right) \left\vert A-c1_{H}\right\vert ^{n}%
\right] x,x\right\rangle \right]  \notag
\end{align}%
for any $x\in H,$ and a similar relation for $y,$ namely%
\begin{align}
& \int_{m-0}^{M}\left\vert \lambda -c\right\vert ^{n-1}\left\langle
E_{\lambda }y,y\right\rangle d\lambda  \label{V.a.e.3.10} \\
& =\frac{1}{n}\left[ \left\langle \left[ \left( M-c\right)
^{n}1_{H}-sgn\left( A-c1_{H}\right) \left\vert A-c1_{H}\right\vert ^{n}%
\right] y,y\right\rangle \right]  \notag
\end{align}%
for any $y\in H.$

The inequality (\ref{V.a.e.3.8}) and the equalities (\ref{V.a.e.3.9}) and (%
\ref{V.a.e.3.10}) produce the second bound in (\ref{V.a.e.3.1.a}).

Finally, observe also that%
\begin{align}
& \int_{m-0}^{M}\left\vert \lambda -c\right\vert ^{n-1}\left\langle
E_{\lambda }x,x\right\rangle d\lambda  \label{V.a.e.3.11} \\
& =\int_{m-0}^{c}\left( c-\lambda \right) ^{n-1}\left\langle E_{\lambda
}x,x\right\rangle d\lambda +\int_{c}^{M}\left( \lambda -c\right)
^{n-1}\left\langle E_{\lambda }x,x\right\rangle d\lambda  \notag \\
& \leq \left( c-m\right) ^{n-1}\int_{m-0}^{c}\left\langle E_{\lambda
}x,x\right\rangle d\lambda +\left( M-c\right) ^{n-1}\int_{c}^{M}\left\langle
E_{\lambda }x,x\right\rangle d\lambda  \notag \\
& \leq \max \left\{ \left( c-m\right) ^{n-1},\left( M-c\right)
^{n-1}\right\} \int_{m-0}^{M}\left\langle E_{\lambda }x,x\right\rangle
d\lambda  \notag \\
& =\left[ \frac{1}{2}\left( M-m\right) +\left\vert c-\frac{m+M}{2}%
\right\vert \right] ^{n-1}  \notag \\
& \times \left[ \left. \left\langle E_{\lambda }x,x\right\rangle \lambda
\right\vert _{m-0}^{M}-\int_{m-0}^{M}\lambda d\left\langle E_{\lambda
}x,x\right\rangle \right]  \notag \\
& =\left[ \frac{1}{2}\left( M-m\right) +\left\vert c-\frac{m+M}{2}%
\right\vert \right] ^{n-1}\left\langle \left( M1_{H}-A\right)
x,x\right\rangle  \notag
\end{align}%
for any $x\in H$ and similarly,%
\begin{align}
& \int_{m-0}^{M}\left\vert \lambda -c\right\vert ^{n-1}\left\langle
E_{\lambda }x,x\right\rangle d\lambda  \label{V.a.e.3.12} \\
& \leq \left[ \frac{1}{2}\left( M-m\right) +\left\vert c-\frac{m+M}{2}%
\right\vert \right] ^{n-1}\left\langle \left( M1_{H}-A\right)
y,y\right\rangle  \notag
\end{align}%
for any $y\in H.$

On making use of (\ref{V.a.e.3.8}), (\ref{V.a.e.3.11}) and (\ref{V.a.e.3.12}%
) we deduce the last bound provided in (\ref{V.a.e.3.1.a}).
\end{proof}

The following particular cases are of interest for applications

\begin{corollary}[Dragomir, 2010, \protect\cite{V.a.SSD7}]
\label{V.a.c.2.2}With the assumption of Theorem \ref{V.a.t.3.1} we have the
inequalities%
\begin{align}
& \left\vert \left\langle T_{n}\left( f,m,M\right) x,y\right\rangle
\right\vert  \label{V.a.e.3.13} \\
& \leq \frac{1}{\left( n-1\right) !}\int_{m-0}^{M}\left( \lambda -m\right)
^{n-1}\dbigvee\limits_{m}^{\lambda }\left( f^{\left( n\right) }\right)
\left\vert \left\langle E_{\lambda }x,y\right\rangle \right\vert d\lambda 
\notag \\
& \leq \frac{1}{\left( n-1\right) !}\dbigvee\limits_{m}^{M}\left( f^{\left(
n\right) }\right) \int_{m-0}^{M}\left( \lambda -m\right) ^{n-1}\left\vert
\left\langle E_{\lambda }x,y\right\rangle \right\vert d\lambda  \notag \\
& \leq \frac{1}{n!}\dbigvee\limits_{m}^{M}\left( f^{\left( n\right) }\right)
B_{n}(m,M,x,y),  \notag
\end{align}%
for any $x,y\in H,$ where%
\begin{align}
& B_{n}(m,M,x,y)  \label{V.a.e.3.14} \\
&  \notag \\
& :=\left\{ 
\begin{array}{l}
\left( M-m\right) ^{n}\left\Vert x\right\Vert \left\Vert y\right\Vert ; \\ 
\\ 
C_{n}(m,M,x,y); \\ 
\\ 
n\left( M-m\right) ^{n-1}\left[ \left\langle \left( M1_{H}-A\right)
x,x\right\rangle \left\langle \left( M1_{H}-A\right) y,y\right\rangle \right]
^{1/2}%
\end{array}%
\right.  \notag
\end{align}%
and%
\begin{align}
& C_{n}(m,M,x,y):=\left[ \left\langle \left[ \left( M-m\right)
^{n}1_{H}-\left( A-m1_{H}\right) ^{n}\right] x,x\right\rangle \right] ^{1/2}
\label{V.a.e.3.15} \\
& \times \left[ \left\langle \left[ \left( M-m\right) ^{n}1_{H}-\left(
A-m1_{H}\right) ^{n}\right] y,y\right\rangle \right] ^{1/2}.  \notag
\end{align}
\end{corollary}

The proof follows from Theorem \ref{V.a.t.3.1} by choosing $c=m$ and
performing the corresponding calculations.

\begin{corollary}[Dragomir, 2010, \protect\cite{V.a.SSD7}]
\label{V.a.c.2.3}With the assumption of Theorem \ref{V.a.t.3.1} we have the
inequalities%
\begin{align}
& \left\vert \left\langle Y_{n}\left( f,m,M\right) x,y\right\rangle
\right\vert  \label{V.a.e.3.16} \\
& \leq \frac{1}{\left( n-1\right) !}\int_{m-0}^{M}\left( M-\lambda \right)
^{n-1}\dbigvee\limits_{\lambda }^{M}\left( f^{\left( n\right) }\right)
\left\vert \left\langle E_{\lambda }x,y\right\rangle \right\vert d\lambda 
\notag \\
& \leq \frac{1}{\left( n-1\right) !}\dbigvee\limits_{m}^{M}\left( f^{\left(
n\right) }\right) \int_{m-0}^{M}\left( M-\lambda \right) ^{n-1}\left\vert
\left\langle E_{\lambda }x,y\right\rangle \right\vert d\lambda  \notag \\
& \leq \frac{1}{n!}\dbigvee\limits_{m}^{M}\left( f^{\left( n\right) }\right) 
\tilde{B}_{n}(m,M,x,y),  \notag
\end{align}%
for any $x,y\in H,$ where%
\begin{align}
& \tilde{B}_{n}(m,M,x,y)  \label{V.a.e.3.17} \\
&  \notag \\
& :=\left\{ 
\begin{array}{l}
\left( M-m\right) ^{n}\left\Vert x\right\Vert \left\Vert y\right\Vert ; \\ 
\\ 
\tilde{C}_{n}(m,M,x,y); \\ 
\\ 
n\left( M-m\right) ^{n-1}\left[ \left\langle \left( M1_{H}-A\right)
x,x\right\rangle \left\langle \left( M1_{H}-A\right) y,y\right\rangle \right]
^{1/2}%
\end{array}%
\right.  \notag
\end{align}%
and%
\begin{equation}
\tilde{C}_{n}(m,M,x,y):=\left[ \left\langle \left( M1_{H}-A\right)
^{n}x,x\right\rangle \right] ^{1/2}\left[ \left\langle \left(
M1_{H}-A\right) ^{n}y,y\right\rangle \right] ^{1/2}.  \label{V.a.e.3.18}
\end{equation}
\end{corollary}

The proof follows from Theorem \ref{V.a.t.3.1} by choosing $c=M$ and
performing the corresponding calculations.

The best bound we can get is incorporated in

\begin{corollary}[Dragomir, 2010, \protect\cite{V.a.SSD7}]
\label{V.a.c.2.4}With the assumption of Theorem \ref{V.a.t.3.1} we have the
inequalities%
\begin{align}
& \left\vert \left\langle W_{n}\left( f,m,M\right) x,y\right\rangle
\right\vert  \label{V.a.e.3.19} \\
& \leq \frac{1}{\left( n-1\right) !}\int_{m-0}^{\frac{m+M}{2}}\left( \frac{%
m+M}{2}-\lambda \right) ^{n-1}\dbigvee\limits_{\lambda }^{\frac{m+M}{2}%
}\left( f^{\left( n\right) }\right) \left\vert \left\langle E_{\lambda
}x,y\right\rangle \right\vert d\lambda  \notag \\
& +\frac{1}{\left( n-1\right) !}\int_{\frac{m+M}{2}}^{M}\left( \lambda -%
\frac{m+M}{2}\right) ^{n-1}\dbigvee\limits_{\frac{m+M}{2}}^{\lambda }\left(
f^{\left( n\right) }\right) \left\vert \left\langle E_{\lambda
}x,y\right\rangle \right\vert d\lambda  \notag \\
& \leq \frac{1}{\left( n-1\right) !}\dbigvee\limits_{m}^{\frac{m+M}{2}%
}\left( f^{\left( n\right) }\right) \int_{m-0}^{\frac{m+M}{2}}\left( \frac{%
m+M}{2}-\lambda \right) ^{n-1}\left\vert \left\langle E_{\lambda
}x,y\right\rangle \right\vert d\lambda  \notag \\
& +\frac{1}{\left( n-1\right) !}\dbigvee\limits_{\frac{m+M}{2}}^{M}\left(
f^{\left( n\right) }\right) \int_{\frac{m+M}{2}}^{M}\left( \lambda -\frac{m+M%
}{2}\right) ^{n-1}\left\vert \left\langle E_{\lambda }x,y\right\rangle
\right\vert d\lambda  \notag \\
& \leq \frac{1}{\left( n-1\right) !}\max \left\{ \dbigvee\limits_{m}^{\frac{%
m+M}{2}}\left( f^{\left( n\right) }\right) ,\dbigvee\limits_{\frac{m+M}{2}%
}^{M}\left( f^{\left( n\right) }\right) \right\}  \notag \\
& \times \int_{m-0}^{M}\left\vert \lambda -\frac{m+M}{2}\right\vert
^{n-1}\left\vert \left\langle E_{\lambda }x,y\right\rangle \right\vert
d\lambda  \notag \\
& \leq \frac{1}{n!}\max \left\{ \dbigvee\limits_{m}^{\frac{m+M}{2}}\left(
f^{\left( n\right) }\right) ,\dbigvee\limits_{\frac{m+M}{2}}^{M}\left(
f^{\left( n\right) }\right) \right\} \breve{B}_{n}(m,M,x,y),  \notag
\end{align}%
for any $x,y\in H,$ where%
\begin{eqnarray}
&&\breve{B}_{n}(m,M,x,y)  \label{V.a.e.3.20} \\
&:&=\left\{ 
\begin{array}{l}
\frac{\left( M-m\right) ^{n}}{2^{n-1}}\left\Vert x\right\Vert \left\Vert
y\right\Vert ; \\ 
\\ 
\breve{C}(m,M,x,y) \\ 
\\ 
\frac{n}{2^{n-1}}\left( M-m\right) ^{n-1}\left[ \left\langle \left(
M1_{H}-A\right) x,x\right\rangle \left\langle \left( M1_{H}-A\right)
y,y\right\rangle \right] ^{1/2}%
\end{array}%
\right.  \notag
\end{eqnarray}%
and%
\begin{multline}
\breve{C}_{n}(m,M,x,y)  \label{V.a.e.3.21} \\
:=\left[ \left\langle \left[ \frac{\left( M-m\right) ^{n}}{2^{n}}%
1_{H}-sgn\left( A-\frac{m+M}{2}1_{H}\right) \left\vert A-\frac{m+M}{2}%
1_{H}\right\vert ^{n}\right] x,x\right\rangle \right] ^{1/2} \\
\times \left[ \left\langle \left[ \frac{\left( M-m\right) ^{n}}{2^{n}}%
1_{H}-sgn\left( A-\frac{m+M}{2}1_{H}\right) \left\vert A-\frac{m+M}{2}%
1_{H}\right\vert ^{n}\right] y,y\right\rangle \right] ^{1/2}.
\end{multline}
\end{corollary}

\subsection{Error Bounds for $f^{\left( n\right) }$ Lipschitzian}

The case when the $n$-th derivative is Lipschitzian is incorporated in the
following result:

\begin{theorem}[Dragomir, 2010, \protect\cite{V.a.SSD7}]
\label{V.a.t.3.2}Let $A$ be a selfadjoint operator in the Hilbert space $H$
with the spectrum $Sp\left( A\right) \subseteq \left[ m,M\right] $ for some
real numbers $m<M$, $\left\{ E_{\lambda }\right\} _{\lambda }$ be its 
\textit{spectral family,} $I$ be a closed subinterval on $\mathbb{R}$ with $%
\left[ m,M\right] \subset \mathring{I}$ (the interior of $I)$ and let $n$ be
an integer with $n\geq 1.$ If $f:I\rightarrow \mathbb{C}$ is such that the $%
n $-th derivative $f^{\left( n\right) }$ is Lipschitzian with the constant $%
L_{n}>0$ on the interval $\left[ m,M\right] $, then for any $c\in \left[ m,M%
\right] $ we have the inequalities%
\begin{align}
& \left\vert \left\langle V_{n}\left( f,c,m,M\right) x,y\right\rangle
\right\vert  \label{V.a.e.3.22} \\
& \leq \frac{1}{n!}L_{n}\int_{m-0}^{M}\left\vert \lambda -c\right\vert
^{n}\left\vert \left\langle E_{\lambda }x,y\right\rangle \right\vert d\lambda
\notag \\
& \leq \frac{1}{\left( n+1\right) !}L_{n}  \notag \\
& \times \left\{ 
\begin{array}{l}
\left[ \left( M-c\right) ^{n+1}+\left( c-m\right) ^{n+1}\right] \left\Vert
x\right\Vert \left\Vert y\right\Vert ; \\ 
\\ 
\left[ \left\langle \left[ \left( M-c\right) ^{n+1}1_{H}-sgn\left(
A-c1_{H}\right) \left\vert A-c1_{H}\right\vert ^{n+1}\right]
x,x\right\rangle \right] ^{1/2} \\ 
\times \left[ \left\langle \left[ \left( M-c\right) ^{n+1}1_{H}-sgn\left(
A-c1_{H}\right) \left\vert A-c1_{H}\right\vert ^{n+1}\right]
y,y\right\rangle \right] ^{1/2}; \\ 
\\ 
\left( n+1\right) \left[ \frac{1}{2}\left( M-m\right) +\left\vert c-\frac{m+M%
}{2}\right\vert \right] ^{n} \\ 
\times \left[ \left\langle \left( M1_{H}-A\right) x,x\right\rangle
\left\langle \left( M1_{H}-A\right) y,y\right\rangle \right] ^{1/2};%
\end{array}%
\right.  \notag
\end{align}%
for any $x,y\in H.$
\end{theorem}

\begin{proof}
From the inequality (\ref{V.a.e.3.2}) in the proof of Theorem \ref{V.a.t.3.1}
we have%
\begin{align}
& \left\vert \left\langle V_{n}\left( f,c,m,M\right) x,y\right\rangle
\right\vert  \label{V.a.e.3.25} \\
& \leq \frac{1}{\left( n-1\right) !}\int_{m-0}^{c}\left\vert \int_{\lambda
}^{c}\left( t-\lambda \right) ^{n-1}d\left( f^{\left( n\right) }\left(
t\right) \right) \right\vert \left\vert \left\langle E_{\lambda
}x,y\right\rangle \right\vert d\lambda  \notag \\
& +\frac{1}{\left( n-1\right) !}\int_{c}^{M}\left\vert \int_{c}^{\lambda
}\left( t-\lambda \right) ^{n-1}d\left( f^{\left( n\right) }\left( t\right)
\right) \right\vert \left\vert \left\langle E_{\lambda }x,y\right\rangle
\right\vert d\lambda  \notag
\end{align}%
for any $x,y\in H.$

Further, we utilize the fact that for an $L-$Lipschitzian function, $p:\left[
\alpha ,\beta \right] \rightarrow \mathbb{C}$ and a Riemann integrable
function $v:\left[ \alpha ,\beta \right] \rightarrow \mathbb{C}$, the
Riemann-Stieltjes integral $\int_{\alpha }^{\beta }p\left( s\right) dv\left(
s\right) $ exists and 
\begin{equation*}
\left\vert \int_{\alpha }^{\beta }p\left( s\right) dv\left( s\right)
\right\vert \leq L\int_{\alpha }^{\beta }\left\vert p\left( s\right)
\right\vert ds.
\end{equation*}

On making use of this property we have for $\lambda \in \left[ m,c\right] $
that 
\begin{equation*}
\left\vert \int_{\lambda }^{c}\left( t-\lambda \right) ^{n-1}d\left(
f^{\left( n\right) }\left( t\right) \right) \right\vert \leq
L_{n}\int_{\lambda }^{c}\left( t-\lambda \right) ^{n-1}dt=\frac{1}{n}%
L_{n}\left( c-\lambda \right) ^{n}
\end{equation*}%
and for $\lambda \in \left[ c,M\right] $ that%
\begin{equation*}
\left\vert \int_{c}^{\lambda }\left( t-\lambda \right) ^{n-1}d\left(
f^{\left( n\right) }\left( t\right) \right) \right\vert \leq
L_{n}\int_{c}^{\lambda }\left( \lambda -t\right) ^{n-1}dt=\frac{1}{n}%
L_{n}\left( \lambda -c\right) ^{n}
\end{equation*}%
which, by (\ref{V.a.e.3.25}) produces the inequality%
\begin{align}
& \left\vert \left\langle V_{n}\left( f,c,m,M\right) x,y\right\rangle
\right\vert  \label{V.a.e.3.26} \\
& \leq \frac{1}{n!}L_{n}\int_{m-0}^{c}\left( c-\lambda \right)
^{n}\left\vert \left\langle E_{\lambda }x,y\right\rangle \right\vert
d\lambda +\frac{1}{n!}L_{n}\int_{c}^{M}\left( \lambda -c\right)
^{n}\left\vert \left\langle E_{\lambda }x,y\right\rangle \right\vert d\lambda
\notag \\
& =\frac{1}{n!}L_{n}\int_{m-0}^{M}\left\vert \lambda -c\right\vert
^{n}\left\vert \left\langle E_{\lambda }x,y\right\rangle \right\vert
d\lambda ,  \notag
\end{align}%
for any $x,y\in H,$ and the first part of (\ref{V.a.e.3.22}) is proved.

Finally, we observe that the bounds for the integral $\int_{m-0}^{M}\left%
\vert \lambda -c\right\vert ^{n}\left\vert \left\langle E_{\lambda
}x,y\right\rangle \right\vert d\lambda $ can be obtained in a similar manner
to the ones from the proof of Theorem \ref{V.a.t.3.1} and the details are
omitted.
\end{proof}

The following result contains error bounds for the particular expansions
considered in Corollary \ref{V.a.c.2.1}:

\begin{corollary}[Dragomir, 2010, \protect\cite{V.a.SSD7}]
\label{V.a.c.3.2}With the assumptions in Theorem \ref{V.a.t.3.2} we have the
inequalities%
\begin{align}
& \left\vert \left\langle T_{n}\left( f,m,M\right) x,y\right\rangle
\right\vert  \label{V.a.e.3.27} \\
& \leq \frac{1}{n!}L_{n}\int_{m-0}^{M}\left( \lambda -m\right)
^{n}\left\vert \left\langle E_{\lambda }x,y\right\rangle \right\vert d\lambda
\notag \\
& \leq \frac{1}{\left( n+1\right) !}L_{n}  \notag \\
& \times \left\{ 
\begin{array}{l}
\left( M-m\right) ^{n+1}\left\Vert x\right\Vert \left\Vert y\right\Vert ; \\ 
\\ 
\left[ \left\langle \left[ \left( M-m\right) ^{n+1}1_{H}-\left(
A-m1_{H}\right) ^{n+1}\right] x,x\right\rangle \right] ^{1/2} \\ 
\times \left[ \left\langle \left[ \left( M-m\right) ^{n+1}1_{H}-\left(
A-m1_{H}\right) ^{n+1}\right] y,y\right\rangle \right] ^{1/2}; \\ 
\\ 
\left( n+1\right) \left( M-m\right) ^{n}\left[ \left\langle \left(
M1_{H}-A\right) x,x\right\rangle \left\langle \left( M1_{H}-A\right)
y,y\right\rangle \right] ^{1/2};%
\end{array}%
\right.  \notag
\end{align}%
and%
\begin{align}
& \left\vert \left\langle Y_{n}\left( f,m,M\right) x,y\right\rangle
\right\vert  \label{V.a.e.3.28} \\
& \leq \frac{1}{n!}L_{n}\int_{m-0}^{M}\left( M-\lambda \right)
^{n}\left\vert \left\langle E_{\lambda }x,y\right\rangle \right\vert d\lambda
\notag \\
& \leq \frac{1}{\left( n+1\right) !}L_{n}  \notag \\
& \times \left\{ 
\begin{array}{l}
\left( M-m\right) ^{n+1}\left\Vert x\right\Vert \left\Vert y\right\Vert ; \\ 
\\ 
\left[ \left\langle \left[ \left( M1_{H}-A\right) ^{n+1}\right]
x,x\right\rangle \right] ^{1/2}\left[ \left\langle \left[ \left(
M1_{H}-A\right) ^{n+1}\right] y,y\right\rangle \right] ^{1/2}; \\ 
\\ 
\left( n+1\right) \left[ \left( M-m\right) \right] ^{n}\left[ \left\langle
\left( M1_{H}-A\right) x,x\right\rangle \left\langle \left( M1_{H}-A\right)
y,y\right\rangle \right] ^{1/2};%
\end{array}%
\right.  \notag
\end{align}%
and%
\begin{multline}
\left\vert \left\langle W_{n}\left( f,m,M\right) x,y\right\rangle \right\vert
\label{V.a.e.3.29} \\
\leq \frac{1}{n!}L_{n}\int_{m-0}^{M}\left\vert \lambda -\frac{m+M}{2}%
\right\vert ^{n}\left\vert \left\langle E_{\lambda }x,y\right\rangle
\right\vert d\lambda \leq \frac{1}{\left( n+1\right) !}L_{n} \\
\times \left\{ 
\begin{array}{l}
\frac{\left( M-m\right) ^{n+1}}{2^{n}}\left\Vert x\right\Vert \left\Vert
y\right\Vert ; \\ 
\\ 
\left[ \left\langle \left[ \frac{\left( M-m\right) ^{n+1}}{2^{n}}%
1_{H}-sgn\left( A-\frac{m+M}{2}1_{H}\right) \left\vert A-\frac{m+M}{2}%
1_{H}\right\vert ^{n+1}\right] x,x\right\rangle \right] ^{1/2} \\ 
\times \left[ \left\langle \left[ \frac{\left( M-m\right) ^{n+1}}{2^{n}}%
1_{H}-sgn\left( A-\frac{m+M}{2}1_{H}\right) \left\vert A-\frac{m+M}{2}%
1_{H}\right\vert ^{n+1}\right] y,y\right\rangle \right] ^{1/2}; \\ 
\\ 
\frac{n+1}{2^{n}}\left( M-m\right) ^{n}\left[ \left\langle \left(
M1_{H}-A\right) x,x\right\rangle \left\langle \left( M1_{H}-A\right)
y,y\right\rangle \right] ^{1/2};%
\end{array}%
\right.
\end{multline}%
for any $x,y\in H,$ respectively.
\end{corollary}

\subsection{Applications}

In order to obtain various vectorial operator inequalities one can use the
above results for particular elementary functions. We restrict ourself \ to
only two examples of functions, namely the exponential and the logarithmic
functions.

If we apply Corollary \ref{V.a.c.2.3} for the exponential function, we can
state the following result:

\begin{proposition}
\label{V.a.p.4.1}Let $A$ be a selfadjoint operator in the Hilbert space $H$
with the spectrum $Sp\left( A\right) \subseteq \left[ m,M\right] $ for some
real numbers $m<M$ and $\left\{ E_{\lambda }\right\} _{\lambda }$ be its 
\textit{spectral family. Then we have}%
\begin{align}
& \left\vert \left\langle e^{A}x,y\right\rangle -e^{M}\sum_{k=0}^{n}\frac{%
\left( -1\right) ^{k}}{k!}\left\langle \left( M1_{H}-A\right)
^{k}x,y\right\rangle \right\vert  \label{V.a.e.4.1} \\
& \leq \frac{1}{\left( n-1\right) !}\int_{m-0}^{M}\left( M-\lambda \right)
^{n-1}\left( e^{M}-e^{\lambda }\right) \left\vert \left\langle E_{\lambda
}x,y\right\rangle \right\vert d\lambda  \notag \\
& \leq \frac{1}{\left( n-1\right) !}\left( e^{M}-e^{m}\right)
\int_{m-0}^{M}\left( M-\lambda \right) ^{n-1}\left\vert \left\langle
E_{\lambda }x,y\right\rangle \right\vert d\lambda  \notag \\
& \leq \frac{1}{n!}\left( e^{M}-e^{m}\right)  \notag \\
& \times \left\{ 
\begin{array}{l}
\left( M-m\right) ^{n}\left\Vert x\right\Vert \left\Vert y\right\Vert ; \\ 
\\ 
\left[ \left\langle \left( M1_{H}-A\right) ^{n}x,x\right\rangle \right]
^{1/2}\left[ \left\langle \left( M1_{H}-A\right) ^{n}y,y\right\rangle \right]
^{1/2} \\ 
\\ 
n\left( M-m\right) ^{n-1}\left[ \left\langle \left( M1_{H}-A\right)
x,x\right\rangle \left\langle \left( M1_{H}-A\right) y,y\right\rangle \right]
^{1/2}%
\end{array}%
\right.  \notag
\end{align}%
for any $x,y\in H.$
\end{proposition}

If we use Corollary \ref{V.a.c.3.2} then we can provide other bounds as
follows:

\begin{proposition}
\label{V.a.p.4.2}With the assumptions of Proposition \ref{V.a.p.4.1} we have%
\begin{align}
& \left\vert \left\langle e^{A}x,y\right\rangle -e^{M}\sum_{k=0}^{n}\frac{%
\left( -1\right) ^{k}}{k!}\left\langle \left( M1_{H}-A\right)
^{k}x,y\right\rangle \right\vert  \label{V.a.e.4.2} \\
& \leq \frac{1}{n!}e^{M}\int_{m-0}^{M}\left( M-\lambda \right)
^{n}\left\vert \left\langle E_{\lambda }x,y\right\rangle \right\vert d\lambda
\notag \\
& \leq \frac{1}{\left( n+1\right) !}e^{M}  \notag \\
& \times \left\{ 
\begin{array}{l}
\left( M-m\right) ^{n+1}\left\Vert x\right\Vert \left\Vert y\right\Vert ; \\ 
\\ 
\left[ \left\langle \left[ \left( M1_{H}-A\right) ^{n+1}\right]
x,x\right\rangle \right] ^{1/2}\left[ \left\langle \left[ \left(
M1_{H}-A\right) ^{n+1}\right] y,y\right\rangle \right] ^{1/2}; \\ 
\\ 
\left( n+1\right) \left[ \left( M-m\right) \right] ^{n}\left[ \left\langle
\left( M1_{H}-A\right) x,x\right\rangle \left\langle \left( M1_{H}-A\right)
y,y\right\rangle \right] ^{1/2};%
\end{array}%
\right.  \notag
\end{align}
\end{proposition}

Finally, the Corollaries \ref{V.a.c.2.3} and \ref{V.a.c.3.2} produce the
following results for the logarithmic function:

\begin{proposition}
\label{V.a.p.4.3}Let $A$ be a positive definite operator in the Hilbert
space $H$ with the spectrum $Sp\left( A\right) \subseteq \left[ m,M\right]
\subset \left( 0,\infty \right) $ and $\left\{ E_{\lambda }\right\}
_{\lambda }$ be its \textit{spectral family, }then%
\begin{align}
& \left\vert \left\langle \ln Ax,y\right\rangle -\left\langle
x,y\right\rangle \ln M+\sum_{k=1}^{n}\frac{\left\langle \left(
M1_{H}-A\right) ^{k}x,y\right\rangle }{kM^{k}}\right\vert  \label{V.a.e.4.3}
\\
& \leq \int_{m-0}^{M}\left( M-\lambda \right) ^{n-1}\frac{M^{n}-\lambda ^{n}%
}{M^{n}\lambda ^{n}}\left\vert \left\langle E_{\lambda }x,y\right\rangle
\right\vert d\lambda  \notag \\
& \leq \frac{M^{n}-m^{n}}{M^{n}m^{n}}\int_{m-0}^{M}\left( M-\lambda \right)
^{n-1}\left\vert \left\langle E_{\lambda }x,y\right\rangle \right\vert
d\lambda  \notag \\
& \leq \frac{M^{n}-m^{n}}{nM^{n}m^{n}}  \notag \\
& \times \left\{ 
\begin{array}{l}
\left( M-m\right) ^{n}\left\Vert x\right\Vert \left\Vert y\right\Vert ; \\ 
\\ 
\left[ \left\langle \left( M1_{H}-A\right) ^{n}x,x\right\rangle \right]
^{1/2}\left[ \left\langle \left( M1_{H}-A\right) ^{n}y,y\right\rangle \right]
^{1/2} \\ 
\\ 
n\left( M-m\right) ^{n-1}\left[ \left\langle \left( M1_{H}-A\right)
x,x\right\rangle \left\langle \left( M1_{H}-A\right) y,y\right\rangle \right]
^{1/2}%
\end{array}%
\right.  \notag
\end{align}%
and%
\begin{align}
& \left\vert \left\langle \ln Ax,y\right\rangle -\left\langle
x,y\right\rangle \ln M+\sum_{k=1}^{n}\frac{\left\langle \left(
M1_{H}-A\right) ^{k}x,y\right\rangle }{kM^{k}}\right\vert  \label{V.a.e.4.4}
\\
& \leq \frac{1}{m^{n+1}}\int_{m-0}^{M}\left( M-\lambda \right)
^{n}\left\vert \left\langle E_{\lambda }x,y\right\rangle \right\vert d\lambda
\notag \\
& \leq \frac{1}{\left( n+1\right) m^{n+1}}  \notag \\
& \times \left\{ 
\begin{array}{l}
\left( M-m\right) ^{n+1}\left\Vert x\right\Vert \left\Vert y\right\Vert ; \\ 
\\ 
\left[ \left\langle \left( M1_{H}-A\right) ^{n+1}x,x\right\rangle \right]
^{1/2}\left[ \left\langle \left( M1_{H}-A\right) ^{n+1}y,y\right\rangle %
\right] ^{1/2}; \\ 
\\ 
\left( n+1\right) \left[ \left( M-m\right) \right] ^{n}\left[ \left\langle
\left( M1_{H}-A\right) x,x\right\rangle \left\langle \left( M1_{H}-A\right)
y,y\right\rangle \right] ^{1/2};%
\end{array}%
\right.  \notag
\end{align}
\end{proposition}

\section{Two Points Taylor's Type Inequalities}

\subsection{Representation Results}

We start with the following identity that has been obtained in \cite%
{V.b.SSD1a}. For the sake of completeness we give here a short proof as well.

\begin{lemma}[Dragomir, 2010, \protect\cite{V.b.SSD1a}]
\label{V.b.l.2.1}Let $I$ be a closed subinterval on $\mathbb{R}$, let $%
a,b\in I$ with $a<b$ and let $n$ be a nonnegative integer. If $%
f:I\rightarrow \mathbb{R}$ is such that the $n$-th derivative $f^{\left(
n\right) }$ is of bounded variation on the interval $\left[ a,b\right] ,$
then, for any $x\in \left[ a,b\right] $ we have the representation 
\begin{align}
f\left( x\right) & =\frac{1}{b-a}\left[ \left( b-x\right) f\left( a\right)
+\left( x-a\right) f\left( b\right) \right]  \label{V.b.e.2.1} \\
& +\frac{\left( b-x\right) \left( x-a\right) }{b-a}  \notag \\
& \times \sum_{k=1}^{n}\frac{1}{k!}\left\{ \left( x-a\right) ^{k-1}f^{\left(
k\right) }\left( a\right) +\left( -1\right) ^{k}\left( b-x\right)
^{k-1}f^{\left( k\right) }\left( b\right) \right\}  \notag \\
& +\frac{1}{b-a}\int_{a}^{b}S_{n}\left( x,t\right) d\left( f^{\left(
n\right) }\left( t\right) \right) ,  \notag
\end{align}%
where the kernel $S_{n}:\left[ a,b\right] ^{2}\rightarrow \mathbb{R}$ is
given by 
\begin{equation}
S_{n}\left( x,t\right) =\frac{1}{n!}\times \left\{ 
\begin{array}{ll}
\left( x-t\right) ^{n}\left( b-x\right) & \text{ if }a\leq t\leq x; \\[10pt] 
\left( -1\right) ^{n+1}\left( t-x\right) ^{n}\left( x-a\right) & \text{if }%
x<t\leq b%
\end{array}%
\right.  \label{V.b.e.2.2}
\end{equation}%
and the integral in the remainder is taken in the Riemann-Stieltjes sense.
\end{lemma}

\begin{proof}
We utilize the following Taylor's representation formula for functions $%
f:I\rightarrow \mathbb{R}$ such that the $n$-th derivatives $f^{\left(
n\right) }$ are of locally bounded variation on the interval $I,$%
\begin{equation}
f\left( x\right) =\sum_{k=0}^{n}\frac{1}{k!}\left( x-c\right) ^{k}f^{\left(
k\right) }\left( c\right) +\frac{1}{n!}\int_{c}^{x}\left( x-t\right)
^{n}d\left( f^{\left( n\right) }\left( t\right) \right) ,  \label{V.b.e.2.3}
\end{equation}%
where $x$ and $c$ are in $I$ and the integral in the remainder is taken in
the Riemann-Stieltjes sense.

Choosing $c=a$ and then $c=b$ in (\ref{V.b.e.2.3}) we can write that 
\begin{equation}
f\left( x\right) =\sum_{k=0}^{n}\frac{1}{k!}\left( x-a\right) ^{k}f^{\left(
k\right) }\left( a\right) +\frac{1}{n!}\int_{a}^{x}\left( x-t\right)
^{n}d\left( f^{\left( n\right) }\left( t\right) \right) ,  \label{V.b.e.2.4}
\end{equation}

and 
\begin{equation}
f\left( x\right) =\sum_{k=0}^{n}\frac{\left( -1\right) ^{k}}{k!}\left(
b-x\right) ^{k}f^{\left( k\right) }\left( b\right) +\frac{\left( -1\right)
^{n+1}}{n!}\int_{x}^{b}\left( t-x\right) ^{n}d\left( f^{\left( n\right)
}\left( t\right) \right) ,  \label{V.b.e.2.5}
\end{equation}%
for any $x\in \left[ a,b\right] .$

Now, by multiplying (\ref{V.b.e.2.4}) with $\left( b-x\right) $ and (\ref%
{V.b.e.2.5}) with $\left( x-a\right) $ we get 
\begin{align}
\left( b-x\right) f\left( x\right) & =\left( b-x\right) f\left( a\right)
+\left( b-x\right) \left( x-a\right) \sum_{k=1}^{n}\frac{1}{k!}\left(
x-a\right) ^{k-1}f^{\left( k\right) }\left( a\right)  \label{V.b.e.2.6} \\
& +\frac{1}{n!}\left( b-x\right) \int_{a}^{x}\left( x-t\right) ^{n}d\left(
f^{\left( n\right) }\left( t\right) \right)  \notag
\end{align}%
and 
\begin{align}
\left( x-a\right) f\left( x\right) & =\left( x-a\right) f\left( b\right)
+\left( b-x\right) \left( x-a\right) \sum_{k=1}^{n}\frac{\left( -1\right)
^{k}}{k!}\left( b-x\right) ^{k-1}f^{\left( k\right) }\left( b\right)
\label{V.b.e.2.7} \\
& +\frac{\left( -1\right) ^{n+1}}{n!}\left( x-a\right) \int_{x}^{b}\left(
t-x\right) ^{n}d\left( f^{\left( n\right) }\left( t\right) \right)  \notag
\end{align}%
respectively.

Finally, by adding the equalities (\ref{V.b.e.2.6}) and (\ref{V.b.e.2.7})
and dividing the sum with $\left( b-a\right) ,$ we obtain the desired
representation (\ref{V.b.e.2.2}).
\end{proof}

\begin{remark}
\label{V.b.r.2.1}The case $n=0$ provides the representation 
\begin{equation}
f\left( x\right) =\frac{1}{b-a}\left[ \left( b-x\right) f\left( a\right)
+\left( x-a\right) f\left( b\right) \right] +\frac{1}{b-a}%
\int_{a}^{b}S\left( x,t\right) d\left( f\left( t\right) \right)
\label{V.b.e.2.8}
\end{equation}%
for any $x\in \left[ a,b\right] ,$ where 
\begin{equation*}
S\left( x,t\right) =\left\{ 
\begin{array}{ll}
b-x & \text{ if }a\leq t\leq x, \\[6pt] 
a-x & \text{if }x<t\leq b,%
\end{array}%
\right.
\end{equation*}%
and $f$ is of bounded variation on $\left[ a,b\right] .$ This result was
obtained by a different approach in \cite{V.b.SSD1}.

The case $n=1$ provides the representation 
\begin{equation}
f\left( x\right) =\frac{1}{b-a}\left[ \left( b-x\right) f\left( a\right)
+\left( x-a\right) f\left( b\right) \right] +\frac{1}{b-a}%
\int_{a}^{b}Q\left( x,t\right) d\left( f^{\prime }\left( t\right) \right) ,
\label{V.b.e.2.10}
\end{equation}%
where 
\begin{equation*}
Q\left( x,t\right) =\left\{ 
\begin{array}{ll}
\left( a-t\right) \left( b-x\right) & \text{ if }a\leq t\leq x, \\[6pt] 
\left( t-b\right) \left( x-a\right) & \text{if }x\leq t\leq b.%
\end{array}%
\right.
\end{equation*}%
Notice that the representation (\ref{V.b.e.2.10}) was obtained by a
different approach in \cite{V.b.SSD1}.
\end{remark}

\begin{theorem}[Dragomir, 2010, \protect\cite{V.b.SSD8}]
\label{V.b.t.2.1}Let $A$ be a selfadjoint operator in the Hilbert space $H$
with the spectrum $Sp\left( A\right) \subseteq \left[ m,M\right] $ for some
real numbers $m<M$, $\left\{ E_{\lambda }\right\} _{\lambda }$ be its 
\textit{spectral family,} $I$ be a closed subinterval on $\mathbb{R}$ with $%
\left[ m,M\right] \subset \mathring{I}$ and let $n$ be an integer with $%
n\geq 1.$ If $f:I\rightarrow \mathbb{C}$ is such that the $n$-th derivative $%
f^{\left( n\right) }$ is of bounded variation on the interval $\left[ m,M%
\right] $, then we have the representation%
\begin{align}
f\left( A\right) & =\frac{1}{M-m}\left[ f\left( m\right) \left(
M1_{H}-A\right) +f\left( M\right) \left( A-m1_{H}\right) \right]
\label{V.b.e.2.11} \\
& +\frac{\left( M1_{H}-A\right) \left( A-m1_{H}\right) }{M-m}  \notag \\
& \times \sum_{k=1}^{n}\frac{1}{k!}\left\{ f^{\left( k\right) }\left(
m\right) \left( A-m1_{H}\right) ^{k-1}+\left( -1\right) ^{k}f^{\left(
k\right) }\left( M\right) \left( M1_{H}-A\right) ^{k-1}\right\}  \notag \\
& +T_{n}\left( f,m,M\right) ,  \notag
\end{align}%
where the remainder $T_{n}\left( f,m,M\right) $ is given by%
\begin{equation}
T_{n}\left( f,m,M\right) :=\frac{1}{\left( M-m\right) n!}\int_{m-0}^{M}K_{n}%
\left( m,M,f;\lambda \right) dE_{\lambda }  \label{V.b.e.2.12}
\end{equation}%
and the kernel $K_{n}\left( m,M,f;\cdot \right) $ has the representation 
\begin{align}
K_{n}\left( m,M,f;\lambda \right) & :=\left( M-\lambda \right) \left(
\int_{m}^{\lambda }\left( \lambda -t\right) ^{n}d\left( f^{\left( n\right)
}\left( t\right) \right) \right)  \label{V.b.e.2.12.a} \\
& +\left( -1\right) ^{n+1}\left( \lambda -m\right) \left( \int_{\lambda
}^{M}\left( t-\lambda \right) ^{n}d\left( f^{\left( n\right) }\left(
t\right) \right) \right)  \notag
\end{align}%
for $\lambda \in \left[ m,M\right] .$
\end{theorem}

\begin{proof}
Utilising Lemma \ref{V.b.l.2.1} we have the representation%
\begin{align}
f\left( \lambda \right) & =\frac{1}{M-m}\left[ \left( M-\lambda \right)
f\left( m\right) +\left( \lambda -m\right) f\left( M\right) \right]
\label{V.b.e.2.13} \\
& +\frac{\left( M-\lambda \right) \left( \lambda -m\right) }{M-m}  \notag \\
& \times \sum_{k=1}^{n}\frac{1}{k!}\left\{ \left( \lambda -m\right)
^{k-1}f^{\left( k\right) }\left( m\right) +\left( -1\right) ^{k}\left(
M-\lambda \right) ^{k-1}f^{\left( k\right) }\left( M\right) \right\}  \notag
\\
& +\frac{1}{\left( M-m\right) n!}\left[ \left( M-\lambda \right)
\int_{m}^{\lambda }\left( \lambda -t\right) ^{n}d\left( f^{\left( n\right)
}\left( t\right) \right) \right.  \notag \\
& \left. +\left( -1\right) ^{n+1}\left( \lambda -m\right) \int_{\lambda
}^{M}\left( t-\lambda \right) ^{n}d\left( f^{\left( n\right) }\left(
t\right) \right) \right] ,  \notag
\end{align}%
for any $\lambda \in \left[ m,M\right] .$

If we integrate (\ref{V.b.e.2.13}) in the Riemann-Stieltjes sense on the
interval $\left[ m,M\right] $ with the integrator $E_{\lambda },$ then we get%
\begin{align}
& \int_{m-0}^{M}f\left( \lambda \right) dE_{\lambda }  \label{V.b.e.2.14} \\
& =\frac{1}{M-m}\int_{m-0}^{M}\left[ \left( M-\lambda \right) f\left(
m\right) +\left( \lambda -m\right) f\left( M\right) \right] dE_{\lambda } 
\notag \\
& +\int_{m-0}^{M}\frac{\left( M-\lambda \right) \left( \lambda -m\right) }{%
M-m}\sum_{k=1}^{n}\frac{1}{k!}\left\{ \left( \lambda -m\right)
^{k-1}f^{\left( k\right) }\left( m\right) \right.  \notag \\
& \left. +\left( -1\right) ^{k}\left( M-\lambda \right) ^{k-1}f^{\left(
k\right) }\left( M\right) \right\} dE_{\lambda }+\frac{1}{\left( M-m\right)
n!}  \notag \\
& \times \left[ \int_{m-0}^{M}\left( M-\lambda \right) \left(
\int_{m}^{\lambda }\left( \lambda -t\right) ^{n}d\left( f^{\left( n\right)
}\left( t\right) \right) \right) dE_{\lambda }\right.  \notag \\
& \left. +\left( -1\right) ^{n+1}\int_{m-0}^{M}\left( \lambda -m\right)
\left( \int_{\lambda }^{M}\left( t-\lambda \right) ^{n}d\left( f^{\left(
n\right) }\left( t\right) \right) \right) dE_{\lambda }\right] .  \notag
\end{align}%
Now, on making use of the spectral representation theorem we deduce from (%
\ref{V.b.e.2.14}) the equality (\ref{V.b.e.2.1}) with the remainder
representation (\ref{V.b.e.2.2}).
\end{proof}

\begin{remark}
\label{V.b.r.2.2}Let $A$ be a selfadjoint operator in the Hilbert space $H$
with the spectrum $Sp\left( A\right) \subseteq \left[ m,M\right] $ for some
real numbers $m<M$, $\left\{ E_{\lambda }\right\} _{\lambda }$ be its 
\textit{spectral family. }In the case when the function $f$ is continuous
and of bounded variation on $\left[ m,M\right] $, then we get the
representation%
\begin{align}
f\left( A\right) & =\frac{1}{M-m}\left[ f\left( m\right) \left(
M1_{H}-A\right) +f\left( M\right) \left( A-m1_{H}\right) \right]
\label{V.b.e.2.15} \\
& +\frac{1}{\left( M-m\right) }\left[ \int_{m-0}^{M}\left( M-\lambda \right) %
\left[ f\left( \lambda \right) -f\left( m\right) \right] dE_{\lambda }\right.
\notag \\
& \left. -\int_{m-0}^{M}\left( \lambda -m\right) \left[ f\left( M\right)
-f\left( \lambda \right) \right] dE_{\lambda }\right] .  \notag
\end{align}%
Also, if the derivative $f^{\prime }$ is of bounded variation, then we have
the representation%
\begin{eqnarray}
f\left( A\right) &=&\frac{1}{M-m}\left[ f\left( m\right) \left(
M1_{H}-A\right) +f\left( M\right) \left( A-m1_{H}\right) \right]
\label{V.b.e.2.16} \\
&&+\frac{1}{\left( M-m\right) }\left[ \int_{m-0}^{M}\left( M-\lambda \right)
\left( \int_{m}^{\lambda }\left( \lambda -t\right) d\left( f^{\prime }\left(
t\right) \right) \right) dE_{\lambda }\right.  \notag \\
&&\left. +\int_{m-0}^{M}\left( \lambda -m\right) \left( \int_{\lambda
}^{M}\left( t-\lambda \right) d\left( f^{\prime }\left( t\right) \right)
\right) dE_{\lambda }\right] .  \notag
\end{eqnarray}
\end{remark}

\begin{example}
\label{V.b.ex.1.1}a. Let $A$ be a selfadjoint operator in the Hilbert space $%
H$ with the spectrum $Sp\left( A\right) \subseteq \left[ m,M\right] $ for
some real numbers $m<M$ and $\left\{ E_{\lambda }\right\} _{\lambda }$ be
its \textit{spectral family. }If we consider the exponential function, then
we get from (\ref{V.b.e.2.11}) and (\ref{V.b.e.2.12}) that 
\begin{align}
e^{A}& =\frac{1}{M-m}\left[ e^{m}\left( M1_{H}-A\right) +e^{M}\left(
A-m1_{H}\right) \right]  \label{V.b.e.2.17} \\
& +\frac{\left( M1_{H}-A\right) \left( A-m1_{H}\right) }{M-m}  \notag \\
& \times \sum_{k=1}^{n}\frac{1}{k!}\left\{ e^{m}\left( A-m1_{H}\right)
^{k-1}+\left( -1\right) ^{k}e^{M}\left( M1_{H}-A\right) ^{k-1}\right\} 
\notag \\
& +\frac{1}{\left( M-m\right) n!}\times \left[ \int_{m-0}^{M}\left(
M-\lambda \right) \left( \int_{m}^{\lambda }\left( \lambda -t\right)
^{n}e^{t}dt\right) dE_{\lambda }\right.  \notag \\
& \left. +\left( -1\right) ^{n+1}\int_{m-0}^{M}\left( \lambda -m\right)
\left( \int_{\lambda }^{M}\left( t-\lambda \right) ^{n}e^{t}dt\right)
dE_{\lambda }\right] .  \notag
\end{align}

b. If $A$ is a positive definite selfadjoint operator with the spectrum $%
Sp\left( A\right) \subseteq \left[ m,M\right] \subset \left( 0,\infty
\right) $ and $\left\{ E_{\lambda }\right\} _{\lambda }$ is its \textit{%
spectral family, then we have the representation}%
\begin{align}
& \ln A=\frac{1}{M-m}\left[ \left( M1_{H}-A\right) \ln m+\left(
A-m1_{H}\right) \ln M\right]  \label{V.b.e.3.21} \\
& +\frac{\left( M1_{H}-A\right) \left( A-m1_{H}\right) }{M-m}  \notag \\
& \times \sum_{k=1}^{n}\frac{1}{k}\left\{ \left( -1\right) ^{k-1}\frac{%
\left( A-m1_{H}\right) ^{k-1}}{m^{k}}-\frac{\left( M1_{H}-A\right) ^{k-1}}{%
M^{k}}\right\}  \notag \\
& +\frac{1}{\left( M-m\right) }\left[ \left( -1\right)
^{n}\int_{m-0}^{M}\left( M-\lambda \right) \left( \int_{m}^{\lambda }\frac{%
\left( \lambda -t\right) ^{n}}{t^{n+1}}dt\right) dE_{\lambda }\right.  \notag
\\
& \left. -\int_{m-0}^{M}\left( \lambda -m\right) \left( \int_{\lambda }^{M}%
\frac{\left( t-\lambda \right) ^{n}}{t^{n+1}}dt\right) dE_{\lambda }\right] .
\notag
\end{align}
\end{example}

The case of functions for which the $n$-th derivative $f^{\left( n\right) }$
is absolutely continuous is of interest for applications. In this case the
remainder can be represented as follows:

\begin{theorem}[Dragomir, 2010, \protect\cite{V.b.SSD8}]
\label{V.b.t.2.2}Let $A$ be a selfadjoint operator in the Hilbert space $H$
with the spectrum $Sp\left( A\right) \subseteq \left[ m,M\right] $ for some
real numbers $m<M$, $\left\{ E_{\lambda }\right\} _{\lambda }$ be its 
\textit{spectral family,} $I$ be a closed subinterval on $\mathbb{R}$ with $%
\left[ m,M\right] \subset \mathring{I}$ and let $n$ be an integer with $%
n\geq 1.$ If $f:I\rightarrow \mathbb{C}$ is such that the $n$-th derivative $%
f^{\left( n\right) }$ is absolutely continuous on the interval $\left[ m,M%
\right] $, then we have the representation (\ref{V.b.e.2.11}) where the
remainder $T_{n}\left( f,m,M\right) $ is given by%
\begin{equation}
T_{n}\left( f,m,M\right) :=\frac{1}{\left( M-m\right) n!}\int_{m-0}^{M}W_{n}%
\left( m,M,f;\lambda \right) E_{\lambda }d\lambda  \label{V.b.e.2.22}
\end{equation}%
and the kernel $W_{n}\left( m,M,f;\cdot \right) $ has the representation 
\begin{align}
W_{n}\left( m,M,f;\lambda \right) & :=\left( -1\right) ^{n}\int_{m}^{\lambda
}\left( \lambda -t\right) ^{n-1}\left[ nM+t-\left( n+1\right) \lambda \right]
f^{\left( n+1\right) }\left( t\right) dt  \label{V.b.e.2.23} \\
& -\int_{\lambda }^{M}\left( t-\lambda \right) ^{n-1}\left[ t+nm-\left(
n+1\right) \lambda \right] f^{\left( n+1\right) }\left( t\right) dt  \notag
\end{align}%
for $\lambda \in \left[ m,M\right] .$
\end{theorem}

\begin{proof}
Observe that, by Leibnitz's rule for differentiation under the integral
sign, we have%
\begin{align}
& \frac{d}{d\lambda }\left[ \left( M-\lambda \right) \left(
\int_{m}^{\lambda }\left( \lambda -t\right) ^{n}f^{\left( n+1\right) }\left(
t\right) dt\right) \right]  \label{V.b.e.2.24} \\
& =-\int_{m}^{\lambda }\left( \lambda -t\right) ^{n}f^{\left( n+1\right)
}\left( t\right) dt+\left( M-\lambda \right) \frac{d}{d\lambda }\left(
\int_{m}^{\lambda }\left( \lambda -t\right) ^{n}f^{\left( n+1\right) }\left(
t\right) dt\right)  \notag \\
& =-\int_{m}^{\lambda }\left( \lambda -t\right) ^{n}f^{\left( n+1\right)
}\left( t\right) dt+n\left( M-\lambda \right) \int_{m}^{\lambda }\left(
\lambda -t\right) ^{n-1}f^{\left( n+1\right) }\left( t\right) dt  \notag \\
& =\int_{m}^{\lambda }\left( \lambda -t\right) ^{n-1}\left[ nM+t-\left(
n+1\right) \lambda \right] f^{\left( n+1\right) }\left( t\right) dt  \notag
\end{align}%
for any $\lambda \in \left[ m,M\right] .$

Integrating by parts in the Riemann-Stieltjes integral we have%
\begin{align}
& \int_{m-0}^{M}\left( M-\lambda \right) \left( \int_{m}^{\lambda }\left(
\lambda -t\right) ^{n}f^{\left( n+1\right) }\left( t\right) dt\right)
dE_{\lambda }  \label{V.b.e.2.25} \\
& =\left. \left( M-\lambda \right) \left( \int_{m}^{\lambda }\left( \lambda
-t\right) ^{n}d\left( f^{\left( n\right) }\left( t\right) \right) \right)
E_{\lambda }\right\vert _{m-0}^{M}  \notag \\
& -\int_{m-0}^{M}\left( \int_{m}^{\lambda }\left( \lambda -t\right) ^{n-1} 
\left[ nM+t-\left( n+1\right) \lambda \right] f^{\left( n+1\right) }\left(
t\right) dt\right) E_{\lambda }d\lambda  \notag \\
& =-\int_{m-0}^{M}\left( \int_{m}^{\lambda }\left( \lambda -t\right) ^{n-1} 
\left[ nM+t-\left( n+1\right) \lambda \right] f^{\left( n+1\right) }\left(
t\right) dt\right) E_{\lambda }d\lambda .  \notag
\end{align}

By Leibnitz's rule we also have%
\begin{align}
& \frac{d}{d\lambda }\left[ \left( \lambda -m\right) \left( \int_{\lambda
}^{M}\left( t-\lambda \right) ^{n}f^{\left( n+1\right) }\left( t\right)
dt\right) \right]  \label{V.b.e.2.26} \\
& =\int_{\lambda }^{M}\left( t-\lambda \right) ^{n}f^{\left( n+1\right)
}\left( t\right) dt+\left( \lambda -m\right) \frac{d}{d\lambda }\left(
\int_{\lambda }^{M}\left( t-\lambda \right) ^{n}f^{\left( n+1\right) }\left(
t\right) dt\right)  \notag \\
& =\int_{\lambda }^{M}\left( t-\lambda \right) ^{n}f^{\left( n+1\right)
}\left( t\right) dt-n\left( \lambda -m\right) \int_{\lambda }^{M}\left(
t-\lambda \right) ^{n-1}f^{\left( n+1\right) }\left( t\right) dt  \notag \\
& =\int_{\lambda }^{M}\left( t-\lambda \right) ^{n-1}\left[ t+nm-\left(
n+1\right) \lambda \right] f^{\left( n+1\right) }\left( t\right) dt  \notag
\end{align}%
for any $\lambda \in \left[ m,M\right] .$

Utilising the integration by parts and (\ref{V.b.e.2.27}) we get 
\begin{align}
& \int_{m-0}^{M}\left( \lambda -m\right) \left( \int_{\lambda }^{M}\left(
t-\lambda \right) ^{n}f^{\left( n+1\right) }\left( t\right) dt\right)
dE_{\lambda }  \label{V.b.e.2.27} \\
& =\left. \left( \lambda -m\right) \left( \int_{\lambda }^{M}\left(
t-\lambda \right) ^{n}f^{\left( n+1\right) }\left( t\right) dt\right)
E_{\lambda }\right\vert _{m-0}^{M}  \notag \\
& -\int_{m-0}^{M}\left( \int_{\lambda }^{M}\left( t-\lambda \right) ^{n-1} 
\left[ t+nm-\left( n+1\right) \lambda \right] f^{\left( n+1\right) }\left(
t\right) dt\right) E_{\lambda }d\lambda  \notag \\
& =-\int_{m-0}^{M}\left( \int_{\lambda }^{M}\left( t-\lambda \right) ^{n-1} 
\left[ t+nm-\left( n+1\right) \lambda \right] f^{\left( n+1\right) }\left(
t\right) dt\right) E_{\lambda }d\lambda .  \notag
\end{align}%
Finally, on utilizing the representation (\ref{V.b.e.2.12}) for the
remainder $T_{n}\left( f,m,M\right) $ and the equalities (\ref{V.b.e.2.25})
and (\ref{V.b.e.2.27}) we deduce (\ref{V.b.e.2.22}). The details are omitted.
\end{proof}

\begin{remark}
\label{V.b.r.2.3}The case $n=1$ provides the following equality%
\begin{align}
f\left( A\right) & =\frac{1}{M-m}\left[ f\left( m\right) \left(
M1_{H}-A\right) +f\left( M\right) \left( A-m1_{H}\right) \right]
\label{V.b.e.2.28} \\
& +\frac{1}{\left( M-m\right) }\int_{m-0}^{M}W_{1}\left( m,M,f;\lambda
\right) E_{\lambda }d\lambda ,  \notag
\end{align}%
where%
\begin{equation}
W_{1}\left( m,M,f;\lambda \right) :=\int_{m}^{\lambda }\left( 2\lambda
-M-t\right) f^{\prime \prime }\left( t\right) dt+\int_{\lambda }^{M}\left(
2\lambda -t-m\right) f^{\prime \prime }\left( t\right) dt  \label{V.b.e.2.29}
\end{equation}%
for $\lambda \in \left[ m,M\right] .$
\end{remark}

\subsection{Error Bounds for $f^{\left( n\right) }$ of Bonded Variation}

The following result that provides bounds for the absolute value of the
kernel $K_{n}\left( m,M,f;\cdot \right) $ holds:

\begin{lemma}[Dragomir, 2010, \protect\cite{V.b.SSD8}]
\label{V.b.l.3.1}Let $I$ be a closed subinterval on $\mathbb{R}$ with $\left[
m,M\right] \subset \mathring{I}$, let $n$ be an integer with $n\geq 1$ and
assume that $f:I\rightarrow \mathbb{C}$ is such that the $n$-th derivative $%
f^{\left( n\right) }$ exists on the interval $\left[ m,M\right] $.

1. If $f^{\left( n\right) }$ is of bounded variation on $\left[ m,M\right] ,$
then 
\begin{align}
& \left\vert K_{n}\left( m,M,f;\lambda \right) \right\vert  \label{V.b.e.3.1}
\\
& \leq \left( M-\lambda \right) \left( \lambda -m\right)
^{n}\dbigvee\limits_{m}^{\lambda }\left( f^{\left( n\right) }\right) +\left(
\lambda -m\right) \left( M-\lambda \right) ^{n}\dbigvee\limits_{\lambda
}^{M}\left( f^{\left( n\right) }\right)  \notag \\
& \leq \frac{1}{4}\left( M-m\right) ^{2}\left[ \left( \lambda -m\right)
^{n-1}\dbigvee\limits_{m}^{\lambda }\left( f^{\left( n\right) }\right)
+\left( M-\lambda \right) ^{n-1}\dbigvee\limits_{\lambda }^{M}\left(
f^{\left( n\right) }\right) \right]  \notag \\
& \leq \frac{1}{4}\left( M-m\right) ^{2}J_{n}\left( m,M;\lambda \right) 
\notag
\end{align}%
where%
\begin{align}
& J_{n}\left( m,M;\lambda \right)  \label{V.b.e.3.1.a} \\
& :=\left\{ 
\begin{array}{l}
\left[ \frac{1}{2}\left( M-m\right) +\left\vert \lambda -\frac{m+M}{2}%
\right\vert \right] ^{n-1}\dbigvee\limits_{m}^{M}\left( f^{\left( n\right)
}\right) ; \\ 
\\ 
\left[ \left( \lambda -m\right) ^{p\left( n-1\right) }+\left( M-\lambda
\right) ^{p\left( n-1\right) }\right] ^{1/p} \\ 
\times \left[ \left( \dbigvee\limits_{m}^{\lambda }\left( f^{\left( n\right)
}\right) \right) ^{q}+\left( \dbigvee\limits_{\lambda }^{M}\left( f^{\left(
n\right) }\right) \right) ^{q}\right] ^{1/q} \\ 
\text{if }p>1,\frac{1}{p}+\frac{1}{q}=1; \\ 
\left[ \frac{1}{2}\dbigvee\limits_{m}^{M}\left( f^{\left( n\right) }\right) +%
\frac{1}{2}\left\vert \dbigvee\limits_{m}^{\lambda }\left( f^{\left(
n\right) }\right) -\dbigvee\limits_{\lambda }^{M}\left( f^{\left( n\right)
}\right) \right\vert \right] \\ 
\times \left[ \left( \lambda -m\right) ^{n-1}+\left( M-\lambda \right) ^{n-1}%
\right]%
\end{array}%
\right.  \notag
\end{align}%
and $\lambda \in \left[ m,M\right] .$

2. If $\lambda \in \left( m,M\right) $ and $f^{\left( n\right) }$ is $%
L_{n,1,\lambda }$-Lipschitzian on $\left[ m,\lambda \right] $ and $%
L_{n,2,\lambda }$-Lipschitzian on $\left[ \lambda ,M\right] ,$ then%
\begin{align}
& \left\vert K_{n}\left( m,M,f;\lambda \right) \right\vert  \label{V.b.e.3.2}
\\
& \leq \frac{1}{n+1}\left[ L_{n,1,\lambda }\left( M-\lambda \right) \left(
\lambda -m\right) ^{n+1}+L_{n,2,\lambda }\left( \lambda -m\right) \left(
M-\lambda \right) ^{n+1}\right]  \notag \\
& \leq \frac{1}{4\left( n+1\right) }\left[ L_{n,1,\lambda }\left( \lambda
-m\right) ^{n}+L_{n,2,\lambda }\left( M-\lambda \right) ^{n}\right]  \notag
\\
& \leq \frac{1}{4\left( n+1\right) }  \notag \\
& \times \left\{ 
\begin{array}{l}
\left[ \left( \lambda -m\right) ^{n}+\left( M-\lambda \right) ^{n}\right]
\max \left\{ L_{n,1,\lambda },L_{n,2,\lambda }\right\} \\ 
\\ 
\left[ \left( \lambda -m\right) ^{pn}+\left( M-\lambda \right) ^{pn}\right]
^{1/p}\left( L_{n,1,\lambda }^{q}+L_{n,2,\lambda }^{q}\right) ^{1/q} \\ 
\text{if }p>1,\frac{1}{p}+\frac{1}{q}=1; \\ 
\left[ \frac{1}{2}\left( M-m\right) +\left\vert \lambda -\frac{m+M}{2}%
\right\vert \right] ^{n}\left( L_{n,1,\lambda }+L_{n,2,\lambda }\right)%
\end{array}%
\right.  \notag
\end{align}%
and $\lambda \in \left[ m,M\right] .$

In particular, if $f^{\left( n\right) }$ is $L_{n}$-Lipschitzian on $\left[
m,M\right] ,$ then%
\begin{align}
& \left\vert K_{n}\left( m,M,f;\lambda \right) \right\vert
\label{V.b.e.3.2.a} \\
& \leq \frac{L_{n}}{n+1}\left[ \left( M-\lambda \right) \left( \lambda
-m\right) ^{n+1}+\left( \lambda -m\right) \left( M-\lambda \right) ^{n+1}%
\right]  \notag \\
& \leq \frac{L_{n}\left( M-m\right) ^{2}}{4\left( n+1\right) }\left[ \left(
\lambda -m\right) ^{n}+\left( M-\lambda \right) ^{n}\right]  \notag
\end{align}%
and $\lambda \in \left[ m,M\right] .$

3. If the function $f^{\left( n\right) }$ is monotonic nondecreasing on $%
\left[ m,M\right] ,$ then 
\begin{align}
& \left\vert K_{n}\left( m,M,f;\lambda \right) \right\vert  \label{V.b.e.3.3}
\\
& \leq \left( M-\lambda \right) \left[ n\int_{m}^{\lambda }\left( \lambda
-t\right) ^{n-1}f^{\left( n\right) }\left( t\right) dt-\left( \lambda
-m\right) ^{n}f^{\left( n\right) }\left( m\right) \right]  \notag \\
& +\left( \lambda -m\right) \left[ \left( M-\lambda \right) ^{n}f^{\left(
n\right) }\left( M\right) -n\int_{\lambda }^{M}\left( t-\lambda \right)
^{n-1}f^{\left( n\right) }\left( t\right) dt\right]  \notag \\
& \leq \left( M-\lambda \right) \left( \lambda -m\right)  \notag \\
& \times \left[ \left( \lambda -m\right) ^{n-1}\left[ f^{\left( n\right)
}\left( \lambda \right) -f^{\left( n\right) }\left( m\right) \right] +\left(
M-\lambda \right) ^{n-1}\left[ f^{\left( n\right) }\left( M\right)
-f^{\left( n\right) }\left( \lambda \right) \right] \right]  \notag \\
& \leq \frac{1}{4}\left( M-m\right) ^{2}  \notag \\
& \times \left[ \left( \lambda -m\right) ^{n-1}\left[ f^{\left( n\right)
}\left( \lambda \right) -f^{\left( n\right) }\left( m\right) \right] +\left(
M-\lambda \right) ^{n-1}\left[ f^{\left( n\right) }\left( M\right)
-f^{\left( n\right) }\left( \lambda \right) \right] \right]  \notag \\
& \leq \frac{1}{4}\left( M-m\right) ^{2}T_{n}\left( m,M;\lambda \right) 
\notag
\end{align}%
where%
\begin{align}
& T_{n}\left( m,M;\lambda \right)  \label{V.b.e.3.3.a} \\
& :=\left\{ 
\begin{array}{l}
\left[ \frac{1}{2}\left( M-m\right) +\left\vert \lambda -\frac{m+M}{2}%
\right\vert \right] ^{n-1}\left[ f^{\left( n\right) }\left( M\right)
-f^{\left( n\right) }\left( m\right) \right] ; \\ 
\\ 
\left[ \left( \lambda -m\right) ^{p\left( n-1\right) }+\left( M-\lambda
\right) ^{p\left( n-1\right) }\right] ^{1/p} \\ 
\times \left[ \left( f^{\left( n\right) }\left( M\right) -f^{\left( n\right)
}\left( \lambda \right) \right) ^{q}+\left( f^{\left( n\right) }\left(
\lambda \right) -f^{\left( n\right) }\left( m\right) \right) ^{q}\right]
^{1/q}\text{ } \\ 
\text{if }p>1,\frac{1}{p}+\frac{1}{q}=1; \\ 
\\ 
\left[ \frac{1}{2}\left[ f^{\left( n\right) }\left( M\right) -f^{\left(
n\right) }\left( m\right) \right] +\left\vert f^{\left( n\right) }\left(
\lambda \right) -\frac{f^{\left( n\right) }\left( M\right) +f^{\left(
n\right) }\left( m\right) }{2}\right\vert \right] \\ 
\times \left[ \left( \lambda -m\right) ^{n-1}+\left( M-\lambda \right) ^{n-1}%
\right] .%
\end{array}%
\right.  \notag
\end{align}
\end{lemma}

\begin{proof}
1. It is well known that if $p:\left[ a,b\right] \rightarrow \mathbb{C}$ is
a continuous function, $v:\left[ a,b\right] \rightarrow \mathbb{C}$ is of
bounded variation then the Riemann-Stieltjes integral $\int_{a}^{b}p\left(
t\right) dv\left( t\right) $ exists and the following inequality holds%
\begin{equation}
\left\vert \int_{a}^{b}p\left( t\right) dv\left( t\right) \right\vert \leq
\max_{t\in \left[ a,b\right] }\left\vert p\left( t\right) \right\vert
\dbigvee\limits_{a}^{b}\left( v\right) ,  \label{V.b.e.3.4}
\end{equation}%
where $\dbigvee\limits_{a}^{b}\left( v\right) $ denotes the total variation
of $v$ on $\left[ a,b\right] .$

Utilising the representation (\ref{V.b.e.2.12.a}) and the property (\ref%
{V.b.e.3.4}) we have successively%
\begin{align}
& \left\vert K_{n}\left( m,M,f;\lambda \right) \right\vert  \label{V.b.e.3.5}
\\
& \leq \left( M-\lambda \right) \left\vert \int_{m}^{\lambda }\left( \lambda
-t\right) ^{n}d\left( f^{\left( n\right) }\left( t\right) \right)
\right\vert +\left( \lambda -m\right) \left\vert \int_{\lambda }^{M}\left(
t-\lambda \right) ^{n}d\left( f^{\left( n\right) }\left( t\right) \right)
\right\vert  \notag \\
& \leq \left( M-\lambda \right) \left( \lambda -m\right)
^{n}\dbigvee\limits_{m}^{\lambda }\left( f^{\left( n\right) }\right) +\left(
\lambda -m\right) \left( M-\lambda \right) ^{n}\dbigvee\limits_{\lambda
}^{M}\left( f^{\left( n\right) }\right)  \notag \\
& =\left( M-\lambda \right) \left( \lambda -m\right) \left[ \left( \lambda
-m\right) ^{n-1}\dbigvee\limits_{m}^{\lambda }\left( f^{\left( n\right)
}\right) +\left( M-\lambda \right) ^{n-1}\dbigvee\limits_{\lambda
}^{M}\left( f^{\left( n\right) }\right) \right]  \notag \\
& \leq \frac{1}{4}\left( M-m\right) ^{2}\left[ \left( \lambda -m\right)
^{n-1}\dbigvee\limits_{m}^{\lambda }\left( f^{\left( n\right) }\right)
+\left( M-\lambda \right) ^{n-1}\dbigvee\limits_{\lambda }^{M}\left(
f^{\left( n\right) }\right) \right]  \notag \\
& \leq \frac{1}{4}\left( M-m\right) ^{2}I_{n}\left( m,M;\lambda \right) 
\notag
\end{align}%
for any $\lambda \in \left[ m,M\right] .$

By H\"{o}lder's inequality we also have%
\begin{align}
& I_{n}\left( m,M;\lambda \right)  \label{V.b.e.3.6} \\
& \leq \left\{ 
\begin{array}{l}
\left[ \frac{1}{2}\left( M-m\right) +\left\vert \lambda -\frac{m+M}{2}%
\right\vert \right] ^{n-1}\dbigvee\limits_{m}^{M}\left( f^{\left( n\right)
}\right) ; \\ 
\\ 
\left[ \left( \lambda -m\right) ^{p\left( n-1\right) }+\left( M-\lambda
\right) ^{p\left( n-1\right) }\right] ^{1/p} \\ 
\times \left[ \left( \dbigvee\limits_{m}^{\lambda }\left( f^{\left( n\right)
}\right) \right) ^{q}+\left( \dbigvee\limits_{\lambda }^{M}\left( f^{\left(
n\right) }\right) \right) ^{q}\right] ^{1/q} \\ 
\text{if }p>1,\frac{1}{p}+\frac{1}{q}=1; \\ 
\left[ \frac{1}{2}\dbigvee\limits_{m}^{M}\left( f^{\left( n\right) }\right) +%
\frac{1}{2}\left\vert \dbigvee\limits_{m}^{\lambda }\left( f^{\left(
n\right) }\right) -\dbigvee\limits_{\lambda }^{M}\left( f^{\left( n\right)
}\right) \right\vert \right] \\ 
\times \left[ \left( \lambda -m\right) ^{n-1}+\left( M-\lambda \right) ^{n-1}%
\right] .%
\end{array}%
\right.  \notag
\end{align}%
for any $\lambda \in \left[ m,M\right] .$

On making use of (\ref{V.b.e.3.5})\ and (\ref{V.b.e.3.6}) we deduce (\ref%
{V.b.e.3.1}).

2. We recall that if $p:\left[ a,b\right] \rightarrow \mathbb{C}$ is a
Riemann integrable function and $v:\left[ a,b\right] \rightarrow \mathbb{C}$
is Lipschitzian with the constant $L>0$, i.e.,%
\begin{equation*}
\left\vert f\left( s\right) -f\left( t\right) \right\vert \leq L\left\vert
s-t\right\vert \text{ for any }t,s\in \left[ a,b\right] ,
\end{equation*}%
then the Riemann-Stieltjes integral $\int_{a}^{b}p\left( t\right) dv\left(
t\right) $ exists and the following inequality holds%
\begin{equation*}
\left\vert \int_{a}^{b}p\left( t\right) dv\left( t\right) \right\vert \leq
L\int_{a}^{b}\left\vert p\left( t\right) \right\vert dt.
\end{equation*}

Now, on applying this property of the Riemann-Stieltjes integral we have%
\begin{align}
& \left\vert K_{n}\left( m,M,f;\lambda \right) \right\vert  \label{V.b.e.3.7}
\\
& \leq \left( M-\lambda \right) \left\vert \int_{m}^{\lambda }\left( \lambda
-t\right) ^{n}d\left( f^{\left( n\right) }\left( t\right) \right)
\right\vert +\left( \lambda -m\right) \left\vert \int_{\lambda }^{M}\left(
t-\lambda \right) ^{n}d\left( f^{\left( n\right) }\left( t\right) \right)
\right\vert  \notag \\
& \leq \frac{1}{n+1}\left[ L_{n,1,\lambda }\left( M-\lambda \right) \left(
\lambda -m\right) ^{n+1}+L_{n,2,\lambda }\left( \lambda -m\right) \left(
M-\lambda \right) ^{n+1}\right]  \notag \\
& =\frac{\left( M-\lambda \right) \left( \lambda -m\right) }{n+1}\left[
L_{n,1,\lambda }\left( \lambda -m\right) ^{n}+L_{n,2,\lambda }\left(
M-\lambda \right) ^{n}\right]  \notag \\
& \leq \frac{\left( M-m\right) ^{2}}{4\left( n+1\right) }\left[
L_{n,1,\lambda }\left( \lambda -m\right) ^{n}+L_{n,2,\lambda }\left(
M-\lambda \right) ^{n}\right]  \notag \\
& \leq \frac{\left( M-m\right) ^{2}}{4\left( n+1\right) }  \notag \\
& \times \left\{ 
\begin{array}{l}
\left[ \left( \lambda -m\right) ^{n}+\left( M-\lambda \right) ^{n}\right]
\max \left\{ L_{n,1,\lambda },L_{n,2,\lambda }\right\} \\ 
\\ 
\left[ \left( \lambda -m\right) ^{pn}+\left( M-\lambda \right) ^{pn}\right]
^{1/p}\left( L_{n,1,\lambda }^{q}+L_{n,2,\lambda }^{q}\right) ^{1/q} \\ 
\text{if }p>1,\frac{1}{p}+\frac{1}{q}=1; \\ 
\left[ \frac{1}{2}\left( M-m\right) +\left\vert \lambda -\frac{m+M}{2}%
\right\vert \right] ^{n}\left( L_{n,1,\lambda }+L_{n,2,\lambda }\right)%
\end{array}%
\right.  \notag
\end{align}%
which prove the desired result (\ref{V.b.e.3.2.a}).

3. From the theory of Riemann-Stieltjes integral is well known that if $p:%
\left[ a,b\right] \rightarrow \mathbb{C}$ is continuous and $v:\left[ a,b%
\right] \rightarrow \mathbb{R}$ is monotonic nondecreasing, then the
Riemann-Stieltjes integrals $\int_{a}^{b}p\left( t\right) dv\left( t\right) $
and $\int_{a}^{b}\left\vert p\left( t\right) \right\vert dv\left( t\right) $
exist and%
\begin{equation}
\left\vert \int_{a}^{b}p\left( t\right) dv\left( t\right) \right\vert \leq
\int_{a}^{b}\left\vert p\left( t\right) \right\vert dv\left( t\right) \leq
\max_{t\in \left[ a,b\right] }\left\vert p\left( t\right) \right\vert \left[
v\left( b\right) -v\left( a\right) \right] .  \label{V.b.e.3.8}
\end{equation}%
By utilizing this property, we have%
\begin{align}
& \left\vert K_{n}\left( m,M,f;\lambda \right) \right\vert  \label{V.b.e.3.9}
\\
& \leq \left( M-\lambda \right) \left\vert \int_{m}^{\lambda }\left( \lambda
-t\right) ^{n}d\left( f^{\left( n\right) }\left( t\right) \right)
\right\vert +\left( \lambda -m\right) \left\vert \int_{\lambda }^{M}\left(
t-\lambda \right) ^{n}d\left( f^{\left( n\right) }\left( t\right) \right)
\right\vert  \notag \\
& \leq \left( M-\lambda \right) \int_{m}^{\lambda }\left( \lambda -t\right)
^{n}d\left( f^{\left( n\right) }\left( t\right) \right) +\left( \lambda
-m\right) \int_{\lambda }^{M}\left( t-\lambda \right) ^{n}d\left( f^{\left(
n\right) }\left( t\right) \right)  \notag \\
& =H_{n}\left( m,M;\lambda \right)  \notag
\end{align}%
By the second part of (\ref{V.b.e.3.8}) we also have that%
\begin{align}
& H_{n}\left( m,M;\lambda \right)  \label{V.b.e.3.10} \\
& \leq \left( M-\lambda \right) \left( \lambda -m\right) ^{n}\left[
f^{\left( n\right) }\left( \lambda \right) -f^{\left( n\right) }\left(
m\right) \right]  \notag \\
& +\left( \lambda -m\right) \left( M-\lambda \right) ^{n}\left[ f^{\left(
n\right) }\left( M\right) -f^{\left( n\right) }\left( \lambda \right) \right]
\notag \\
& =\left( M-\lambda \right) \left( \lambda -m\right)  \notag \\
& \times \left[ \left( \lambda -m\right) ^{n-1}\left[ f^{\left( n\right)
}\left( \lambda \right) -f^{\left( n\right) }\left( m\right) \right] +\left(
M-\lambda \right) ^{n-1}\left[ f^{\left( n\right) }\left( M\right)
-f^{\left( n\right) }\left( \lambda \right) \right] \right]  \notag \\
& \leq \frac{1}{4}\left( M-m\right) ^{2}  \notag \\
& \times \left[ \left( \lambda -m\right) ^{n-1}\left[ f^{\left( n\right)
}\left( \lambda \right) -f^{\left( n\right) }\left( m\right) \right] +\left(
M-\lambda \right) ^{n-1}\left[ f^{\left( n\right) }\left( M\right)
-f^{\left( n\right) }\left( \lambda \right) \right] \right]  \notag \\
& =\frac{1}{4}\left( M-m\right) ^{2}L_{n}\left( m,M;\lambda \right)  \notag
\end{align}%
with%
\begin{align}
& L_{n}\left( m,M;\lambda \right)  \label{V.b.e.3.11} \\
& \leq \left\{ 
\begin{array}{l}
\left[ \frac{1}{2}\left( M-m\right) +\left\vert \lambda -\frac{m+M}{2}%
\right\vert \right] ^{n-1}\left[ f^{\left( n\right) }\left( M\right)
-f^{\left( n\right) }\left( m\right) \right] ; \\ 
\\ 
\left[ \left( \lambda -m\right) ^{p\left( n-1\right) }+\left( M-\lambda
\right) ^{p\left( n-1\right) }\right] ^{1/p} \\ 
\times \left[ \left( f^{\left( n\right) }\left( M\right) -f^{\left( n\right)
}\left( \lambda \right) \right) ^{q}+\left( f^{\left( n\right) }\left(
\lambda \right) -f^{\left( n\right) }\left( m\right) \right) ^{q}\right]
^{1/q}\text{ } \\ 
\text{if }p>1,\frac{1}{p}+\frac{1}{q}=1; \\ 
\\ 
\left[ \frac{1}{2}\left[ f^{\left( n\right) }\left( M\right) -f^{\left(
n\right) }\left( m\right) \right] +\left\vert f^{\left( n\right) }\left(
\lambda \right) -\frac{f^{\left( n\right) }\left( M\right) +f^{\left(
n\right) }\left( m\right) }{2}\right\vert \right] \\ 
\times \left[ \left( \lambda -m\right) ^{n-1}+\left( M-\lambda \right) ^{n-1}%
\right] .%
\end{array}%
\right.  \notag
\end{align}%
Integrating by parts we have%
\begin{align}
& H_{n}\left( m,M;\lambda \right)  \label{V.b.e.3.12} \\
& =\left( M-\lambda \right) \int_{m}^{\lambda }\left( \lambda -t\right)
^{n}d\left( f^{\left( n\right) }\left( t\right) \right) +\left( \lambda
-m\right) \int_{\lambda }^{M}\left( t-\lambda \right) ^{n}d\left( f^{\left(
n\right) }\left( t\right) \right)  \notag \\
& =\left( M-\lambda \right) \left[ \left. \left( \lambda -t\right)
^{n}f^{\left( n\right) }\left( t\right) \right\vert _{m}^{\lambda
}+n\int_{m}^{\lambda }\left( \lambda -t\right) ^{n-1}f^{\left( n\right)
}\left( t\right) dt\right]  \notag \\
& +\left( \lambda -m\right) \left[ \left. \left( t-\lambda \right)
^{n}f^{\left( n\right) }\left( t\right) \right\vert _{\lambda
}^{M}-n\int_{\lambda }^{M}\left( t-\lambda \right) ^{n-1}f^{\left( n\right)
}\left( t\right) dt\right]  \notag \\
& =\left( M-\lambda \right) \left[ n\int_{m}^{\lambda }\left( \lambda
-t\right) ^{n-1}f^{\left( n\right) }\left( t\right) dt-\left( \lambda
-m\right) ^{n}f^{\left( n\right) }\left( m\right) \right]  \notag \\
& +\left( \lambda -m\right) \left[ \left( M-\lambda \right) ^{n}f^{\left(
n\right) }\left( M\right) -n\int_{\lambda }^{M}\left( t-\lambda \right)
^{n-1}f^{\left( n\right) }\left( t\right) dt\right] .  \notag
\end{align}%
On making use of (\ref{V.b.e.3.9})-(\ref{V.b.e.3.12}) we deduce the desired
result (\ref{V.b.e.3.3}).
\end{proof}

On making use of the bounds for the kernel $K_{n}\left( m,M,f;\cdot \right) $
provided above, we can establish the following error estimates for the
remainder $T_{n}\left( f,m,M\right) $ in the representation formula (\ref%
{V.b.e.2.11}).

\begin{theorem}[Dragomir, 2010, \protect\cite{V.b.SSD8}]
\label{V.b.t.3.1}Let $A$ be a selfadjoint operator in the Hilbert space $H$
with the spectrum $Sp\left( A\right) \subseteq \left[ m,M\right] $ for some
real numbers $m<M$, $\left\{ E_{\lambda }\right\} _{\lambda }$ be its 
\textit{spectral family,} $I$ be a closed subinterval on $\mathbb{R}$ with $%
\left[ m,M\right] \subset \mathring{I}$ and let $n$ be an integer with $%
n\geq 1.$ If $f:I\rightarrow \mathbb{C}$ is such that the $n$-th derivative $%
f^{\left( n\right) }$ is of bounded variation on the interval $\left[ m,M%
\right] $, then we have the representation%
\begin{align}
\left\langle f\left( A\right) x,y\right\rangle & =\frac{1}{M-m}\left[
f\left( m\right) \left\langle \left( M1_{H}-A\right) x,y\right\rangle
+f\left( M\right) \left\langle \left( A-m1_{H}\right) x,y\right\rangle %
\right]  \label{V.b.e.3.13} \\
& +\frac{1}{M-m}  \notag \\
& \times \left\{ \sum_{k=1}^{n}\frac{1}{k!}f^{\left( k\right) }\left(
m\right) \left\langle \left( M1_{H}-A\right) \left( A-m1_{H}\right)
^{k}x,y\right\rangle \right.  \notag \\
& \left. +\sum_{k=1}^{n}\frac{1}{k!}\left( -1\right) ^{k}f^{\left( k\right)
}\left( M\right) \left\langle \left( A-m1_{H}\right) \left( M1_{H}-A\right)
^{k}x,y\right\rangle \right\}  \notag \\
& +T_{n}\left( f,m,M;x,y\right) ,  \notag
\end{align}%
where the remainder $T_{n}\left( f,m,M;x,y\right) $ is given by%
\begin{equation}
T_{n}\left( f,m,M;x,y\right) :=\frac{1}{\left( M-m\right) n!}%
\int_{m-0}^{M}K_{n}\left( m,M,f;\lambda \right) d\left\langle E_{\lambda
}x,y\right\rangle  \label{V.b.e.3.14}
\end{equation}%
and the kernel $K_{n}\left( m,M,f;\cdot \right) $ has the representation (%
\ref{V.b.e.2.12.a}).

Moreover, we have the error estimate%
\begin{align}
& \left\vert T_{n}\left( f,m,M;x,y\right) \right\vert   \label{V.b.e.3.15} \\
& \leq \frac{1}{4n!}\left( M-m\right) \dbigvee\limits_{m-0}^{M}\left(
\left\langle E_{\left( \cdot \right) }x,y\right\rangle \right)   \notag \\
& \times \max_{\lambda \in \left[ m,M\right] }\left[ \left( \lambda
-m\right) ^{n-1}\dbigvee\limits_{m}^{\lambda }\left( f^{\left( n\right)
}\right) +\left( M-\lambda \right) ^{n-1}\dbigvee\limits_{\lambda
}^{M}\left( f^{\left( n\right) }\right) \right]   \notag \\
& \leq \frac{1}{4n!}\left( M-m\right) ^{n}\dbigvee\limits_{m}^{M}\left(
f^{\left( n\right) }\right) \dbigvee\limits_{m-0}^{M}\left( \left\langle
E_{\left( \cdot \right) }x,y\right\rangle \right)   \notag \\
& \leq \frac{1}{4n!}\left( M-m\right) ^{n}\dbigvee\limits_{m}^{M}\left(
f^{\left( n\right) }\right) \left\Vert x\right\Vert \left\Vert y\right\Vert 
\notag
\end{align}%
for any $x,y\in H.$
\end{theorem}

\begin{proof}
The identity (\ref{V.b.e.3.13}) with the remainder representation (\ref%
{V.b.e.3.14}) follows from (\ref{V.b.e.2.11}) and (\ref{V.b.e.2.12}).

Now, on utilizing the property (\ref{V.b.e.3.4}) for the Riemann-Stieltjes
integral we deduce from (\ref{V.b.e.3.14}) that%
\begin{equation}
\left\vert T_{n}\left( f,m,M;x,y\right) \right\vert \leq \frac{1}{\left(
M-m\right) n!}\max_{\lambda \in \left[ m,M\right] }\left\vert K_{n}\left(
m,M,f;\lambda \right) \right\vert \dbigvee\limits_{m-0}^{M}\left(
\left\langle E_{\left( \cdot \right) }x,y\right\rangle \right) 
\label{V.b.e.3.16}
\end{equation}%
for any $x,y\in H.$

Further, by (\ref{V.b.e.3.1}) and (\ref{V.b.e.3.1.a}) we have the bounds%
\begin{align}
& \left\vert K_{n}\left( m,M,f;\lambda \right) \right\vert
\label{V.b.e.3.17} \\
& \leq \frac{1}{4}\left( M-m\right) ^{2}\left[ \left( \lambda -m\right)
^{n-1}\dbigvee\limits_{m}^{\lambda }\left( f^{\left( n\right) }\right)
+\left( M-\lambda \right) ^{n-1}\dbigvee\limits_{\lambda }^{M}\left(
f^{\left( n\right) }\right) \right]  \notag \\
& \leq \frac{1}{4}\left( M-m\right) ^{2}\left[ \frac{1}{2}\left( M-m\right)
+\left\vert \lambda -\frac{m+M}{2}\right\vert \right] ^{n-1}\dbigvee%
\limits_{m}^{M}\left( f^{\left( n\right) }\right) ;  \notag
\end{align}%
for any $\lambda \in \left[ m,M\right] .$

Taking the maximum over $\lambda \in \left[ m,M\right] $ in (\ref{V.b.e.3.17}%
) we deduce the first and the second inequalities in (\ref{V.b.e.3.15}).

The last part follows by the Total Variation Schwarz's inequality and we
omit the details.
\end{proof}

\begin{corollary}[Dragomir, 2010, \protect\cite{V.b.SSD8}]
\label{V.b.c.3.1}With the assumptions from Theorem \ref{V.b.t.3.1} and if $%
f^{\left( n\right) }$ is $L_{n}$-Lipschitzian on $\left[ m,M\right] ,$ then%
\begin{align}
& \left\vert T_{n}\left( f,m,M;x,y\right) \right\vert   \label{V.b.e.3.19} \\
& \leq \frac{1}{\left( n+1\right) !\left( M-m\right) }L_{n}\dbigvee%
\limits_{m-0}^{M}\left( \left\langle E_{\left( \cdot \right)
}x,y\right\rangle \right)   \notag \\
& \times \max_{\lambda \in \left[ m,M\right] }\left[ \left( M-\lambda
\right) \left( \lambda -m\right) ^{n+1}+\left( \lambda -m\right) \left(
M-\lambda \right) ^{n+1}\right]   \notag \\
& \leq \frac{1}{4\left( n+1\right) !}\left( M-m\right)
^{n+1}L_{n}\dbigvee\limits_{m-0}^{M}\left( \left\langle E_{\left( \cdot
\right) }x,y\right\rangle \right)   \notag \\
& \leq \frac{1}{4\left( n+1\right) !}\left( M-m\right) ^{n+1}L_{n}\left\Vert
x\right\Vert \left\Vert y\right\Vert   \notag
\end{align}%
for any $x,y\in H.$
\end{corollary}

\subsection{Error Bounds for $f^{\left( n\right) }$ Absolutely Continuous}

The following result that provides bounds for the absolute value of the
kernel $W_{n}\left( m,M,f;\cdot \right) $ holds:

\begin{lemma}[Dragomir, 2010, \protect\cite{V.b.SSD8}]
\label{V.b.l.4.1}Let $I$ be a closed subinterval on $\mathbb{R}$ with $\left[
m,M\right] \subset \mathring{I}$, let $n$ be an integer with $n\geq 1$ and
assume that $f:I\rightarrow \mathbb{C}$ is such that the $n$-th derivative $%
f^{\left( n\right) }$ is absolutely continuous on the interval $\left[ m,M%
\right] $. Then we have the bound%
\begin{equation}
\left\vert W_{n}\left( m,M,f;\lambda \right) \right\vert \leq
\sum_{i=1}^{4}B_{n}^{\left( i\right) }\left( m,M,f;\lambda \right)
\label{V.b.e.4.1}
\end{equation}%
where%
\begin{align}
& B_{n}^{\left( 1\right) }\left( m,M,f;\lambda \right)  \label{V.b.e.4.2} \\
& :=n\left( M-\lambda \right) \int_{m}^{\lambda }\left( \lambda -t\right)
^{n-1}\left\vert f^{\left( n+1\right) }\left( t\right) \right\vert dt\leq
n\left( M-\lambda \right)  \notag \\
& \times \left\{ 
\begin{array}{l}
\frac{1}{n}\left( \lambda -m\right) ^{n}\left\Vert f^{\left( n+1\right)
}\right\Vert _{\left[ m,\lambda \right] ,\infty }\text{ if }f^{\left(
n+1\right) }\in L_{\infty }\left[ m,\lambda \right] ; \\ 
\\ 
\frac{1}{\left[ \left( n-1\right) p_{1}+1\right] ^{1/p_{1}}}\left( \lambda
-m\right) ^{n-1+1/p_{1}}\left\Vert f^{\left( n+1\right) }\right\Vert _{\left[
m,\lambda \right] ,q_{1}}\text{ } \\ 
\text{if }f^{\left( n+1\right) }\in L_{q_{1}}\left[ m,\lambda \right]
,p_{1}>1,\frac{1}{p_{1}}+\frac{1}{q_{1}}=1; \\ 
\\ 
\left( \lambda -m\right) ^{n-1}\left\Vert f^{\left( n+1\right) }\right\Vert
_{\left[ m,\lambda \right] ,1};%
\end{array}%
\right.  \notag
\end{align}%
\begin{align}
& B_{n}^{\left( 2\right) }\left( m,M,f;\lambda \right)  \label{V.b.e.4.3} \\
& :=\int_{m}^{\lambda }\left( \lambda -t\right) ^{n}\left\vert f^{\left(
n+1\right) }\left( t\right) \right\vert dt  \notag \\
& \leq \left\{ 
\begin{array}{l}
\frac{1}{n+1}\left( \lambda -m\right) ^{n+1}\left\Vert f^{\left( n+1\right)
}\right\Vert _{\left[ m,\lambda \right] ,\infty }\text{ if }f^{\left(
n+1\right) }\in L_{\infty }\left[ m,\lambda \right] ; \\ 
\\ 
\frac{1}{\left( np_{2}+1\right) ^{1/p_{2}}}\left( \lambda -m\right)
^{n+1/p_{2}}\left\Vert f^{\left( n+1\right) }\right\Vert _{\left[ m,\lambda %
\right] ,q_{2}}\text{ if }f^{\left( n+1\right) }\in L_{q_{2}}\left[
m,\lambda \right] , \\ 
p_{2}>1,\frac{1}{p_{2}}+\frac{1}{q_{2}}=1; \\ 
\\ 
\left( \lambda -m\right) ^{n}\left\Vert f^{\left( n+1\right) }\right\Vert _{ 
\left[ m,\lambda \right] ,1}%
\end{array}%
\right.  \notag
\end{align}%
\begin{align}
& B_{n}^{\left( 3\right) }\left( m,M,f;\lambda \right)  \label{V.b.e.4.4} \\
& :=\int_{\lambda }^{M}\left( t-\lambda \right) ^{n}\left\vert f^{\left(
n+1\right) }\left( t\right) \right\vert dt  \notag \\
& \leq \left\{ 
\begin{array}{l}
\frac{1}{n+1}\left( M-\lambda \right) ^{n+1}\left\Vert f^{\left( n+1\right)
}\right\Vert _{\left[ \lambda ,M\right] ,\infty }\text{ if }f^{\left(
n+1\right) }\in L_{\infty }\left[ \lambda ,M\right] ; \\ 
\\ 
\frac{1}{\left( np_{3}+1\right) ^{1/p_{3}}}\left( M-\lambda \right)
^{n+1/p_{3}}\left\Vert f^{\left( n+1\right) }\right\Vert _{\left[ \lambda ,M%
\right] ,q_{3}}\text{ if }f^{\left( n+1\right) }\in L_{q_{3}}\left[ \lambda
,M\right] , \\ 
p_{3}>1,\frac{1}{p_{3}}+\frac{1}{q_{3}}=1; \\ 
\\ 
\left( M-\lambda \right) ^{n}\left\Vert f^{\left( n+1\right) }\right\Vert _{ 
\left[ \lambda ,M\right] ,1}%
\end{array}%
\right.  \notag
\end{align}%
and%
\begin{align}
& B_{n}^{\left( 4\right) }\left( m,M,f;\lambda \right)  \label{V.b.e.4.5} \\
& :=n\left( \lambda -m\right) \int_{\lambda }^{M}\left( t-\lambda \right)
^{n-1}\left\vert f^{\left( n+1\right) }\left( t\right) \right\vert dt\leq
n\left( \lambda -m\right)  \notag \\
& \times \left\{ 
\begin{array}{l}
\frac{1}{n}\left( M-\lambda \right) ^{n}\left\Vert f^{\left( n+1\right)
}\right\Vert _{\left[ \lambda ,M\right] ,\infty }\text{ if }f^{\left(
n+1\right) }\in L_{\infty }\left[ \lambda ,M\right] ; \\ 
\\ 
\frac{1}{\left[ \left( n-1\right) p_{4}+1\right] ^{1/p_{4}}}\left( M-\lambda
\right) ^{n-1+1/p_{4}}\left\Vert f^{\left( n+1\right) }\right\Vert _{\left[
\lambda ,M\right] ,q_{4}}\text{ } \\ 
\text{if }f^{\left( n+1\right) }\in L_{q_{1}}\left[ \lambda ,M\right]
,p_{4}>1,\frac{1}{p4}+\frac{1}{q_{4}}=1; \\ 
\\ 
\left( M-\lambda \right) ^{n-1}\left\Vert f^{\left( n+1\right) }\right\Vert
_{\left[ m,\lambda \right] ,1};%
\end{array}%
\right.  \notag
\end{align}%
for any $\lambda \in \left[ m,M\right] ,$ where the Lebesgue norms $%
\left\Vert \cdot \right\Vert _{\left[ a,b\right] ,p}$ are defined by%
\begin{equation*}
\left\Vert g\right\Vert _{\left[ a,b\right] ,p}:=\left\{ 
\begin{array}{l}
\left( \int_{a}^{b}\left\vert g\left( t\right) \right\vert ^{p}dt\right)
^{1/p}\text{ if }g\in L_{p}\left[ a,b\right] ,p\geq 1 \\ 
\\ 
ess\sup_{t\in \left[ a,b\right] }\left\vert g\left( t\right) \right\vert 
\text{ if }g\in L_{\infty }\left[ a,b\right] .%
\end{array}%
\right.
\end{equation*}
\end{lemma}

\begin{proof}
From (\ref{V.b.e.2.23}) we have 
\begin{align}
& \left\vert W_{n}\left( m,M,f;\lambda \right) \right\vert  \label{V.b.e.4.6}
\\
& \leq \left\vert \int_{m}^{\lambda }\left( \lambda -t\right) ^{n-1}\left[
nM+t-\left( n+1\right) \lambda \right] f^{\left( n+1\right) }\left( t\right)
dt\right\vert  \notag \\
& +\left\vert \int_{\lambda }^{M}\left( t-\lambda \right) ^{n-1}\left[
t+nm-\left( n+1\right) \lambda \right] f^{\left( n+1\right) }\left( t\right)
dt\right\vert  \notag \\
& \leq \int_{m}^{\lambda }\left( \lambda -t\right) ^{n-1}\left\vert
nM+t-\left( n+1\right) \lambda \right\vert \left\vert f^{\left( n+1\right)
}\left( t\right) \right\vert dt  \notag \\
& +\int_{\lambda }^{M}\left( t-\lambda \right) ^{n-1}\left\vert t+nm-\left(
n+1\right) \lambda \right\vert \left\vert f^{\left( n+1\right) }\left(
t\right) \right\vert dt  \notag \\
& \leq \int_{m}^{\lambda }\left( \lambda -t\right) ^{n-1}\left[ n\left(
M-\lambda \right) +\left( \lambda -t\right) \right] \left\vert f^{\left(
n+1\right) }\left( t\right) \right\vert dt  \notag \\
& +\int_{\lambda }^{M}\left( t-\lambda \right) ^{n-1}\left[ \left( t-\lambda
\right) +n\left( \lambda -m\right) \right] \left\vert f^{\left( n+1\right)
}\left( t\right) \right\vert dt  \notag \\
& =\sum_{i=1}^{4}B_{n}^{\left( i\right) }\left( m,M,f;\lambda \right)  \notag
\end{align}%
for any $\lambda \in \left[ m,M\right] ,$ which proves (\ref{V.b.e.4.1}).

The other bounds follows by H\"{o}lder's integral inequality and the details
are omitted.
\end{proof}

\begin{remark}
\label{V.b.r.4.1}It is obvious that the inequalities (\ref{V.b.e.4.1})-(\ref%
{V.b.e.4.5}) can produce 12 different bounds for $\left\vert W_{n}\left(
m,M,f;\lambda \right) \right\vert .$ However, we mention here only the case
when $f^{\left( n+1\right) }\in L_{\infty }\left[ \lambda ,M\right] ,$ namely%
\begin{align}
& \left\vert W_{n}\left( m,M,f;\lambda \right) \right\vert  \label{V.b.e.4.7}
\\
& \leq \left( M-\lambda \right) \left( \lambda -m\right) ^{n}\left\Vert
f^{\left( n+1\right) }\right\Vert _{\left[ m,\lambda \right] ,\infty }+\frac{%
1}{n+1}\left( \lambda -m\right) ^{n+1}\left\Vert f^{\left( n+1\right)
}\right\Vert _{\left[ m,\lambda \right] ,\infty }  \notag \\
& +\frac{1}{n+1}\left( M-\lambda \right) ^{n+1}\left\Vert f^{\left(
n+1\right) }\right\Vert _{\left[ \lambda ,M\right] ,\infty }+\left( \lambda
-m\right) \left( M-\lambda \right) ^{n}\left\Vert f^{\left( n+1\right)
}\right\Vert _{\left[ \lambda ,M\right] ,\infty }  \notag \\
& \leq \left[ \left( M-\lambda \right) \left( \lambda -m\right) ^{n}+\left(
\lambda -m\right) \left( M-\lambda \right) ^{n}\right.  \notag \\
& \left. +\frac{1}{n+1}\left( \lambda -m\right) ^{n+1}+\frac{1}{n+1}\left(
M-\lambda \right) ^{n+1}\right] \left\Vert f^{\left( n+1\right) }\right\Vert
_{\left[ m,M\right] ,\infty }  \notag
\end{align}%
for any $\lambda \in \left[ m,M\right] .$
\end{remark}

Finally, we can state the following result as well:

\begin{theorem}[Dragomir, 2010, \protect\cite{V.b.SSD8}]
\label{V.b.t.4.1}Let $A$ be a selfadjoint operator in the Hilbert space $H$
with the spectrum $Sp\left( A\right) \subseteq \left[ m,M\right] $ for some
real numbers $m<M$, $\left\{ E_{\lambda }\right\} _{\lambda }$ be its 
\textit{spectral family,} $I$ be a closed subinterval on $\mathbb{R}$ with $%
\left[ m,M\right] \subset \mathring{I}$ and let $n$ be an integer with $%
n\geq 1.$ If $f:I\rightarrow \mathbb{C}$ is such that the $n$-th derivative $%
f^{\left( n\right) }$ is absolutely continuous on the interval $\left[ m,M%
\right] $, then we have the representation (\ref{V.b.e.3.13}) where the
remainder $T_{n}\left( f,m,M;x,y\right) $ is given by%
\begin{equation}
T_{n}\left( f,m,M;x,y\right) :=\frac{1}{\left( M-m\right) n!}%
\int_{m-0}^{M}W_{n}\left( m,M,f;\lambda \right) \left\langle E_{\lambda
}x,y\right\rangle d\lambda  \label{V.b.e.4.8}
\end{equation}%
and the kernel $W_{n}\left( m,M,f;\cdot \right) $ has the representation (%
\ref{V.b.e.2.23}).

We also have the error bounds%
\begin{align}
& \left\vert T_{n}\left( f,m,M;x,y\right) \right\vert  \label{V.b.e.4.9} \\
& \leq \frac{1}{\left( M-m\right) n!}\int_{m-0}^{M}\left\vert W_{n}\left(
m,M,f;\lambda \right) \right\vert \left\vert \left\langle E_{\lambda
}x,y\right\rangle \right\vert d\lambda  \notag \\
& \leq \frac{1}{\left( M-m\right) n!}\int_{m-0}^{M}\left\vert W_{n}\left(
m,M,f;\lambda \right) \right\vert \left\langle E_{\lambda }x,x\right\rangle
^{1/2}\left\langle E_{\lambda }y,y\right\rangle ^{1/2}d\lambda  \notag \\
& \leq \frac{1}{\left( M-m\right) n!}\left\Vert x\right\Vert \left\Vert
y\right\Vert \int_{m}^{M}\left\vert W_{n}\left( m,M,f;\lambda \right)
\right\vert d\lambda  \notag
\end{align}%
for any $x,y\in H.$
\end{theorem}

\begin{remark}
\label{V.b.r.4.2}On making use of Lemma \ref{V.b.l.4.1} one can produce
further bounds. However, the details are left to the interested reader.
\end{remark}

\end{document}